\documentclass[11pt,reqno]{amsart}
  \usepackage[margin=0.98in]{geometry}

 \usepackage{amsmath,amssymb}
 \usepackage{algpseudocode}
 \usepackage{algorithm}
 \usepackage{algorithmicx}
 \usepackage{caption}
 \usepackage{subcaption}
 \usepackage{amsthm}
 \numberwithin{equation}{section}
 \usepackage{amsfonts}
 \usepackage{graphicx}
  \usepackage{mathrsfs}
 \usepackage{hyperref}
 \usepackage{makecell}  
 \usepackage{longtable}
 \usepackage{comment}
 \usepackage{enumitem}

\usepackage{tikz}
 
 \hypersetup{bookmarksdepth=2}

\usepackage{amsmath}
\usepackage{amssymb}
\usepackage{amsthm}
\usepackage{amsmath}
\usepackage{setspace}
\usepackage{xcolor}
\usepackage{fancyhdr}
\usepackage{amssymb}
\usepackage{amsthm}
\usepackage{listings}
    \usepackage{cite}
    \usepackage{tabularx}
    \usepackage{graphicx}
    \usepackage{epstopdf}    
    \usepackage{epsfig}
    \usepackage{float}
    \usepackage{ listings}
    \usepackage{appendix}
 \theoremstyle{plain}

 \newtheorem{thm}{Theorem}[section]
 
\newtheorem*{conjecture*}{Conjecture}

 \newtheorem{prop}[thm]{Proposition}
 \newtheorem{lem}[thm]{Lemma}
 \newtheorem{cor}[thm]{Corollary}
 \theoremstyle{definition}

 \newtheorem{remark}[thm]{Remark}
  \newtheorem{ques}{Question}

\newcommand\red[1]{{\color{red}{#1}}}
\newcommand\ang[1]{ { \la {#1}\ra } }

\newcommand \nnra[1]{ {\| #1 \|}_{\wwra} }
\newcommand \nnza[1]{ {\| #1 \|}_{\wwza} }
\newcommand \nna[1]{ {\| #1 \|}_{\linfza} }
\newcommand \nnla[1]{ {\| #1 \|}_{\linfa} }
\newcommand \nlinf[1]{ {\| #1 \|}_{\linf} }
\newcommand \ncXc[1]{ {\| #1 \|}_{\cXc} }

\newcommand \nnr[1]{ {\| #1 \|}_{\wwra} }
\newcommand \nnrr[1]{ {\| #1 \|}_{\wwrr} }

\newcommand \nnz[1]{ {\| #1 \|}_{\wwza} }
\newcommand \nnn[1]{ {\| #1 \|}_{\wonea} }
\newcommand \cZZ[1]{ {\| #1 \|}_{\cZ} }

\def \beps {\bar \e}
\def \epa { \e_2}

\let \hyr = \hyperref
\let \its = \itshape

\newcommand \nnab[1]{ {\| #1 \|}_{  L_{\mw{1D}, \rhoo \rhoc}^{\infty}   } }
\newcommand \nnrb[1]{ {\| #1 \|}_{ \dot \cW^{1,\infty}_{r, \rhoo \rhoc } } }
\newcommand \nnlb[1]{ {\| #1 \|}_{  L^{\infty}_{ \rhoo \rhoc } } }

\newcommand \nnzb[1]{ {\| #1 \|}_{  \dot \cW^{1,\infty}_{z, \rhoo \rhoc } } }
\newcommand \llz[1] { { E_{z, R}( #1 ) } }
\newcommand \llr[1]{ { E_{r, R}( #1) } } 
\newcommand \llw[1]{ { E_{r, R}( #1) } }

\newcommand \lla[1]{ { E_{ 0, R}( #1) } } 
\newcommand \lln[1]{ { E_{\cW^{1,\infty}, R}( #1) } }

\def \KRF {K_{\R^5}}

\def \ffar {\mw{far}}
\def \nne {\mw{ne}}

\def \cXc {\cX_3}
\def \cZc {\cZ_{\mathbb{C}}}
\def \cZs {\cZ_{\mathsf{s}}}
\def \cZu {\cZ_{\mathsf{u}}}
\def \ssf {\mathsf{s}}

\def \sss {\mathsf{g}}
\def \sfu {\mathsf{u}}

\def \cFc {\cF_{\mw{3D}}}
\def \cLa {\cL_{\al}}

\def \cRa { \td \cR_{\al}}

\def \cLc {\cL_{\mw{bd}}}
\def \cNc {\cN_{\mw{3D}}}
\def \wwd {\om_{\mw{1D}}}
\def \tcRd {\td \cR_{\mw{1D}}}

\def \vpa {\vp_\sst}

\def \vpc {\vp_c}
\def \vpcc {\vp_c}

\def \vpk {\vp_{\mw{one}}}

\def \hk { {  \f{1+\kp}{2} } }
\def \mhk { {  \f{1-\kp}{2} } }

\def \sst { \mf{ne}}

\def \muc { \mu_{\mf{ctr}} }

\def \alb{\bar \al}
\def \vmu{ \vec \mu }

\def \he {\hat \e_{\b}}

\def \hal {\hat \al_l}
\def \hau {\hat \al_u}
\def \als { \al_* }
\def \aln {\al_{\phi}}

\def \ze {\zeta}
\def \NNN {\mathsf{N}}

\def \CCs { {\mathscr{C}} }

\def \CCa { {\mathscr{C}_1} }
\def \CCj {\mathscr{C}_J}

\def \EEs {\mathscr{E}}

 \def\tw{{\wt W}}
 \def\bw{ \bar W}

\def \Lpsi {\cL_{\psi}}
\def \BS {\mathsf{BS}}
\def \BSa {\mathsf{BS}_{\mf{0}}}
\def \FFF {\cF_{\mathrm{\R^3}}}

\def \linf {L^{\infty}}

\def \bbb { \kp_\mf{O}}

\def \bbu{ \kp_{\G} }
\def \wa{ W_{\b} }

\def \waa{ \bar{W}_{\al} }
\def \vaa {\bar V_{\al}}
\def \NNN { \mathsf{N} }

\def\omthgg { \Om_{*,\g}^{\th}}
\def\omthss { \Om_{*,\al}^{\th}}

\def\omsth {\bar \Om_{\al}^{\th}}
\def \delin { \d_{\init}}

\def \omth {\Om^{\th}}

\def \wwwa {\overline W_{\alb}}
\def \vvva {\bar V_{\alb}}

\def \ws {\bar \Om_{\al}}

\def \wwb {\bar \Om}
\def \www { \Om }
\def \QQb {\bar \QQ}
\def \ww {\om}
\def \psio {\mathring{\psi}}
\def \psiox {\mathring{\psi}_x}

\def \psid {\psi_{\mathsf{1D}}}
\def \psiad {\psi_{\al, \mathsf{1D}}}
\def \psiod { \mr{\psi}_{ \mathsf{1D}}}
\def \psioad {\mathring{\psi}_{\al ,\mathsf{1D}}}

\def \psis {\bar \Psi_{\al}}
\def \QQs {\bar \QQ_{\al}}
\def \QQQs {\bar Q_{\al}}

\def \bcB{ \bar \cB}
\def \tcB{ \td \cB}
\def \cBs {\bar \cB_{\al}}
\def \upse {\upsilon_{\e}}
\def \phim {\phi_{\mathsf{M}}}
\def \phia {\phi_{\mathsf{A}}}

\def \TTa {\mf{T}_{\al}}
\def \AAa {\mf{A}_{\al}}
\def \AAaa {\mf{A}_{\alb}}

\def \CCl {\mf{C}_{ \mf{l}} } 
\def \CCw {\mf{C}_{ \mf{\om, \th}} }

\def \bcomth { \bar {c}_{\om,\th} }
\def \comth { c_{\om,\th} }
\def \cxx { \mf{c_{ x} }}
\def \cxs { \mf{c_{ x, \al} }}
\def \cxg { \mf{c_{ x, \gam} }}

\def \cls {\bar c_{l, \al}}
\def \cws  {\bar c_{\om, \al}}
\def \dds {\bar d_{\al}}

\def \barcp {  \mathsf{\bar \kp}_{\psi} }
\def \cpsia { \mathsf{\kp}_{\psi} }
\def \cpsi { \mathsf{\kp}_{\psi, 1} }

\def \cpsib { \mathsf{\kp}_{\psi,2} }
\def \wwd {\ww_{\mathsf{1D}}}

\def \Kin{ K_{\mathsf{in}}}

\def \tz {{\td z}}
\def \tr {{\td r}}

\def \bcl{ \bar c_l }
\def \bcw {\bar c_{\om}}

\def \cw {c_{\om}}
\def \tcl{ \td c_l }
\def \tcw {\td c_{\om}}

\def \mfL {\mf{L}}
\def \mfM  {\mf{P}}

 \def \cJa{ \cJ_{\al} }
 
  \def \cJab{ \cJ_{\al, \b} }
 \def \cJaa{ \la \cJ_{\al, \b} \ra }

 \def \cJak { |\cJaa(\xx)|^{\kp} }
 \def \cJakkk { |\cJaa(\xx)|^{\kp + 1} }

\def \JJ {\cJ_{\mw{3D}}}

\def \lamcL {\lam_{\cL}}

\def \cff {C_{\mathsf{far}}}

\def \cmr {\cM^r}
\def \cmz {\cM^z}

\def \vpr{ \vp^r } %\vp^{(r)}}
\def \vpz{ \vp^z }

\def \bpsi {\bar \Psi}

\def \qqr {q^r}
\def \qqz {q^z}

\def \rhoo { \G }
\def \rhor{ {\G^r }}
\def \rhoz{  {\G^z} }
\def \rkeya {\G_{30}}
\def \rkey {\G_3}
\def \rhoc {\G_{\e}}

\def \rag{ \G_{\mw{ag}}}

\def \lgp { \log_+ }

\def \com {c_{\om}}

\def \init  {\mathsf{in}}

\def \ccb {\mf{C}_1}

\def \kag{\kp_{\mw{ag}} }

\def \Rreg {\bar \cR_{\mw{reg}}}
\def \Rsin {\bar \cR_{\mw{sin}}}

\def \Lreg{  \cL_{\mw{reg}} }
\def \Lsin{  \cL_{\mw{sin}} }

\def \Qreg { Q_{\mw{reg}} }
\def \Qsin { Q_{\mw{sin}} }

\def \Rrem { \bar \cR_{\mw{rem}} }

\def \linff {L^{\infty}}

\def \wwrr { {\cW}^{1,\infty}_{r,\rhoo}}
\def \wwra {\dot{\cW}^{1,\infty}_{r,\rhoo}}
\def \wwza {\dot{\cW}^{1,\infty}_{z, \rhoo}}
\def \linfza {L_{\mw{1D}, \rhoo}^{\infty}}

\def \linfa {L_{\rhoo}^{\infty}}
\def \wonea {\cW^{1,\infty}_{\rhoo}}

\let \tts = \textstyle
\let \ssk = \substack
 \let \mw = \mathrm
 \let \we = \wedge
\let \mf = \mathsf
\let  \rsa =\rightsquigarrow

\let \texps = \texorpdfstring

\let \tf = \tfrac

\let \nolim = \nolimits

\let \co = \circ
 \let\pa=\partial
 \let\al=\alpha
 \let\b=\beta
 \let\d=\delta
 \let\g=\gamma
 \let \gam = \gamma 
 \let\e=\varepsilon

 \let\z=\zeta
 \let \kp = \kappa
 \let\lam=\lambda
 
 \let\s=\sigma
 \let\f=\frac
 \let \iin = \infty

 \let \les = \lesssim
  \let \gtr = \gtrsim
 \let\om=\omega
 
 \let \th = \theta
  
 \let \pr = \prime
 \let \vp = \varphi
 \let\G= \Gamma
 
\let\B = \Big
 \let\D=\Delta
 \let\Lam=\Lambda
 \let\S=\Sigma
 \let \Ups = \Upsilon
 \let\Om=\Omega
 \let\td = \tilde

 \let \mr = \mathring

 \let\teq \triangleq
 
 \let\pa=\partial
 \let \bsh = \backslash
 \let \vs = \vspace
 \let \cc = \circ

 \def\cA{{\mathcal A}}
 \def\cB{{\mathcal B}}
 \def\cC{{\mathcal C}}
 
 \def\cE{{\mathcal E}}
 \def\cF{{\mathcal F}}
 
 \def\cH{{\mathcal H}}
 
 \def\cJ{{\mathcal J}}
 
 \def\cL{{\mathcal L}}
 \def\cM{{\mathcal M}}
 \def\cN{{\mathcal N}}

 \def\cR{{\mathcal R}}
 \def\cS{{\mathcal S}}
 \def\cT{{\mathcal T}}
 \def\cU{{\mathcal U}}
 
 \def\cW{{\mathcal W}}
 \def\cX{{\mathcal X}}
 \def\cY{{\mathcal Y}}
 \def\cZ{{\mathcal Z}}

 \def\tw{\widetilde{w}}

 \def\cM{{\mathcal M}}
 \def\cN{{\mathcal N}}

 \def\na{\nabla}
 \def\la{\langle}
 \def\ra{\rangle}

\def\one{\mathbf{1}}
\def\udb{\underbrace}

  \def\Re{\mathrm{Re}}

 \newcommand{\bit}{\begin{itemize}}
 \newcommand{\eit}{\end{itemize}}

 \newcommand{\bseq}{\begin{subequations}}
 \newcommand{\eseq}{\end{subequations}}

 \newcommand{\beq}{\begin{equation}}
 \newcommand{\eeq}{\end{equation}}
  \newcommand{\bal}{\begin{aligned} }
  \newcommand{\eal}{\end{aligned}}
    \newcommand{\bga}{ \begin{gathered} }
  \newcommand{\ega}{ \end{gathered} }
 \newcommand{\ben}{\begin{eqnarray}}
 \newcommand{\een}{\end{eqnarray}}
 \newcommand{\beno}{\begin{eqnarray*}}
 \newcommand{\eeno}{\end{eqnarray*}}
  \newcommand{\bb}{\mathbf{b}}
 \newcommand{\ee}{\mathbf{e}}
 \newcommand{\hh}{\mathbf{h}}

 \newcommand{\nn}{\mathbf{n}}
  \newcommand{\pp}{\mathbf{p}}
  \newcommand{\qq}{\mathbf{q}}

 \newcommand{\uu}{\mathbf{u}}

 \newcommand{\xx}{\mathbf{x}}
 \newcommand{\yy}{\mathbf{y}}

  \newcommand{\XXb}{ \bar X }

\newcommand{\MXX}{ \mathscr{X} }

\newcommand{\MSS}{ \mathscr{S} }

 \newcommand{\QQ}{\mathbf{Q}}

 \newcommand{\R}{\mathbb{R}}
\newcommand{\C}{\mathbb{C}}
 
 \newcommand{\UU}{\mathbf{U}}

 \newcommand{\Id}{\mathrm{Id}}

 \newcommand{\Tr}{\mathrm{Tr}}

 \newcommand{\supp}{\mathrm{supp}}

 \author{Jiajie Chen}
 \date{ \today}

\address{Department of Mathematics, University of Chicago, Chicago, IL 60637.}
\email{\href{jiajiechen@uchicago.edu}{jiajiechen@uchicago.edu}}

\title[3D Profiles, Blowup, and Limiting behavior]{
Asymptotically Self-Similar Blowup for 3D Incompressible Euler with $C^{1, 1/3-}$ Velocity II: 3D Profiles, Blowup, and Limiting behavior}

\begin{document}

\begin{abstract} 

For any $\alpha \in (0,\tfrac13)$, we construct exact $C^{\alpha}$ self-similar blowup profiles for 
the vorticity of the 3D incompressible Euler equation without swirl, and build on them to prove asymptotically self-similar blowup from $C_c^\alpha$ initial vorticity and $C^{1,\alpha}\cap L^2$ initial velocity.
Moreover, we provide a complete characterization of the limiting behavior of 
the $C^{\alpha}$ vorticity profiles and the associated blowup solutions as $\alpha\to(\tfrac13)^-$. Specifically, as $\alpha \to(\tfrac13)^-$, the spatial blowup rate $\mathsf{c}_{\mathsf{x},\alpha}$ diverges to $\infty$, while the $C^{\alpha}$ vorticity profile $\Omega_{*,\alpha}^{\theta}$ asymptotically factorizes and converges strongly in a weighted $L^\infty$ norm to 
a nonzero constant multiple of $r^{1/3}\overline W_{\frac{1}{3}}(z)$, where $\overline W_{\frac{1}{3}}$ is a $C^\infty$ 1D blowup profile.
Our construction is inspired by the Hou--Zhang blowup scenario.
Using a fixed-point argument, we lift the $C^\infty$ blowup profiles for a 1D model constructed in the companion work \cite{chen2026eulerI} to exact 3D blowup profiles. To overcome the lack of $r$-directional decay in the approximate profile and capture the anisotropic structure, we develop a family of anisotropic weighted estimates and introduce a crucial integration-by-parts method along trajectories that exploits the equation twice.
We then develop a finite codimension stability argument in a low-regularity setting to prove stability of the 3D profiles and establish asymptotically self-similar blowup.
This blowup result is sharp in view of the global regularity theory for axisymmetric Euler without swirl with $ C_c^{\alpha}$ initial vorticity for all $\alpha \geq \frac13$. To the best of our knowledge, our results provide the first example in which a singularity from a 1D nonlocal fluid model is lifted to construct blowup for incompressible fluid equations in $\mathbb{R}^2$ or $\mathbb{R}^3$.

\end{abstract}

 \maketitle
 \tableofcontents

\section{Introduction}

Whether the 3D incompressible Euler equations can develop a finite-time singularity from smooth initial data with finite energy is one of the major open questions in mathematical fluid mechanics \cite{constantin2007euler,hou2009blow,majda2002vorticity}. In the 
 vorticity formulation, the 3D Euler equations read 
\begin{equation}\label{euler0}
   \om_{t} + \uu \cdot \nabla \om = \om \cdot \nabla \uu,  \quad \uu = \na \times (-\D)^{-1} \om,
\end{equation}
where $\uu$ is the velocity vector and $\om = \na \times \uu$ is the vorticity vector. 
Two main difficulties arise in analyzing \eqref{euler0}. The first is the nonlocal Biot-Savart law $ \uu = \na \times (-\D)^{-1} \om$, which significantly complicates the analysis.  The second  is the competition between the stabilizing effect of advection $\uu \cdot \na \om$ \cite{hou2006dynamic,hou2008dynamic,lei2009stabilizing} and the destabilizing effect of vortex stretching $\om \cdot (\na \uu)$. 

Owing to this strongly nonlinear and nonlocal structure, the question of finite-time blowup has long been studied through numerical simulations dating back to the 1990s \cite{grauer1991numerical,pumir1992development,kerr1993evidence}, and more recently through neural-network-based methods \cite{Tristan2022,wang2025discovery}. An influential scenario was proposed by Hou--Luo \cite{luo2013potentially-2}, who presented convincing numerical evidence for a potential finite-time singularity in the 3D axisymmetric Euler equations with smooth initial data and boundary. 
This scenario inspired several subsequent works \cite{kiselev2013small,choi2015finite,choi2014on}.

In the remarkable work \cite{elgindi2019finite}, Elgindi established singularity formation for 3D axisymmetric Euler without swirl with \(C^{1,\al}\) velocity for small \(\al\).
The result was later extended to allow small swirl \cite{elgindi2019stability}. 
With Hou \cite{chen2019finite2}, we established finite-time blowup for 2D Boussinesq and 3D axisymmetric Euler with \(C^{1,\al}\) velocity, large swirl, and boundary.
In these works, the parameter $\alpha$ plays an essential role as a small parameter 
to weaken the advection and simplify the equations.

In joint works with Hou \cite{ChenHou2023a,ChenHou2023b}, inspired by the Hou-Luo scenario \cite{luo2013potentially-2}, we established stable nearly self-similar blowup for 2D Boussinesq and 3D axisymmetric Euler with smooth boundary from smooth initial data. 
This provides the first construction of a 3D Euler singularity from smooth data in a smooth domain. 
See also the review paper \cite{chen2025singularity}. The proof combines singularly weighted energy estimates, sharp functional inequalities, quantitative finite rank perturbation, and computer-assisted rigorous numerics based on interval arithmetic.

Subsequently, in \cite{cordoba2023finite},  Cordoba, Martinez-Zoroa, and Zheng developed 
a new framework based on iterative layer constructions and established blowup for 3D axisymmetric Euler with $C^{1,\alpha} \cap C^{\infty}(\R^3 \backslash \{0 \})$ velocity and small $\al$. 
Using this framework and by introducing a solution-dependent bounded force $f \in C^{1,1/2-}$, Cordoba and Martinez-Zoroa \cite{cordoba2023blow} established blowup for the forced 3D Euler with a smooth velocity. This line of work was later extended to the forced Boussinesq \cite{cordoba2025finite}
and IPM \cite{cordoba2024finite}. 
Inspired by \cite{cordoba2023finite},  in \cite{chen2024remarks}, we refined the construction in \cite{elgindi2019finite,chen2019finite2} to obtain a self-similar blowup solution with $C^{\infty}$ regularity except the blowup point. 
In  \cite{elgindi2023instability},  Elgindi-Pasqualotto considered a new blowup scenario 
to construct finite-time blowup for Boussinesq equations in $\R^2$ and for 3D axisymmetric Euler  away from the axis with $C^{1,\al}$ velocity and small $\al$.

\subsection{Motivation}

After the works mentioned above, especially \cite{elgindi2019finite},\cite{ChenHou2023a,ChenHou2023b}, 
the following questions naturally arise for \eqref{euler0}.

\begin{ques}\label{ques:blowup}
Establish finite-time blowup for Euler equations \eqref{euler0} in $\R^3$ with $C^{1,\al}$ velocity 
without requiring $\al$ to be small. 
In the axisymmetric Euler without swirl setting, establish finite-time blowup with $C^{1,\al} \cap L^2$ velocity for any $\al \in (0, \f13)$.
\end{ques}

In the axisymmetric Euler without swirl setting, Danchin \cite{danchin2007axisymmetric} and
Shao--Wei--Zhang  \cite{shao2026global} proved global regularity for solutions with
\(C_c^\alpha\) initial vorticity when \(\alpha\geq \frac13\) (with 
initial velocity in $C^{1,\alpha}\cap L^2$);
\footnote{
The global regularity result for $C_c^{1/3}$ initial vorticity $\om_0$ was independently obtained 
by the author and Shao-Wei-Zhang \cite{shao2026global} using the critical 
conserved quantity $ \| r^{-1} \om^{\th} \|_{L^{3, \infty}}$ (private communication with Dongyi Wei, July 2024). The compact support assumption can also be weakened. See \cite{shao2026global} and \cite{danchin2007axisymmetric}.}
see also the work of Saint-Raymond \cite{saint1994remarks}. This suggests the threshold \(\alpha=\frac13\) in Question~\ref{ques:blowup}. 

One of the main contributions of this work is to answer Question~\ref{ques:blowup} below the threshold \(\alpha=\frac13\) by a self-similar approach. More precisely, for \emph{any} \(\alpha\in(0,\frac13)\), we establish asymptotically self-similar blowup for 3D axisymmetric Euler without swirl with \(C_c^\alpha\) initial vorticity and \(C^{1,\alpha}\cap L^2\) initial velocity.

We also note that, in a very recent and remarkable work, Shkoller \cite{shkoller2026incompressible} 
developed a new Lagrangian framework to prove finite-time blowup for 3D axisymmetric Euler without swirl, with $C^{1,\alpha} \cap L^2$  initial velocity for the \emph{entire range} $\alpha \in (0, \f13)$. The two works were carried out independently and use different methods:
a Lagrangian clock-and-strain framework in \cite{shkoller2026incompressible}, 
and a self-similar approach in the Eulerian framework in the present work. 
Moreover, we develop a new lifting method to construct \(C^\alpha\) self-similar blowup profiles 
for the vorticity and provide a complete characterization of their limiting behavior as \(\alpha\to(\frac13)^-\).

More broadly,  self-similar methods have been effectively applied to construct finite-time singularities in fluids and related equations, including 3D
incompressible Euler \cite{elgindi2019finite,ChenHou2023a,ChenHou2023b},  
compressible fluids  \cite{merle2022implosion2,buckmaster2022smooth,chen2024Euler,chen2024vorticity,chen2026implosion}, and, more recently, the Landau equation in kinetic theory \cite{bedrossian2026finite}.

\vs{0.05in}

\paragraph{\bf Continuation of $C^{\al}$ vorticity profiles}
Motivated by the self-similar singularity  with  $C^{\al}$  vorticity for $\al \ll 1$  in \cite{elgindi2019finite}, 
a natural approach to Question \ref{ques:blowup} is to use $\al$ as a continuation parameter to construct $C^{\al}$ self-similar profiles for 
the vorticity in \eqref{euler0} with larger $\al$. See discussion in \cite[Problem 12]{Elg22}. 
Such an approach has been investigated numerically by Hou-Zhang \cite{hou2024potential} 
for the axisymmetric Euler without swirl. In this setting, \eqref{euler0} reduces to 
\beq\label{eq:Euler}
 \pa_t \om^{\th} + ( u^r \pa_r + u^z \pa_z ) \om^{\th} =  \f{u^r}{r} \om^{\th} ,
\eeq
where $\om^{\th}$ is the angular vorticity and $(u^r, u^z)$ are the associated velocity components.

To impose the $C^{\al}$-regularity, in \cite{hou2024potential}, Hou-Zhang considered a solution $\om^{\th}$ to \eqref{eq:Euler} of the form $\om^{\th} = r^{\al} \om(r, z)$ and performed simulations on the resulting equations \eqref{eq:Euler2} for $\om(r, z)$ instead of $\om^{\th}$. Starting from 
\emph{smooth initial data} for $\om(r, z)$, the authors provided strong numerical evidence 
on potential self-similar singularity formation for $\al \leq 0.3$ and computed numerically 
the self-similar profiles when $\al = 0.05 i, 0\leq i\leq 6$. 
The solution in  \cite{hou2024potential}
is smooth near the $z$-axis, which is different from that in \cite{elgindi2019finite}. Due to the global regularity result \cite{danchin2007axisymmetric}, this leads to the following questions:

\begin{ques}\label{ques:SS}
Can one construct self-similar profiles and blowup solutions for \eqref{eq:Euler} 
with $C^{\al}$ vorticity for all $\al \in (0, \f13)$? What are the blowup rates and 
limiting behaviors as $\al \to \f13$? 
 What is the mechanism preventing the construction at $\al = \f13$ or beyond?

\end{ques}

The study of Question~\ref{ques:SS} for \eqref{eq:Euler} may also shed light on the possible construction of vorticity profiles for the full 3D Euler \eqref{euler0} with regularity beyond $C^{1/3}$. Due to the inclusion $C^{ \b} \subset C^{ \al}$ for $\b >  \al$, the key for Question \ref{ques:SS} is the range close to the endpoint $\f13$. Such a range cannot be reached numerically in \cite{hou2024potential}. Nevertheless, the authors reported an important feature that as $\al$ increases to $\f13$, the profile $\om(r, z)$ becomes very \emph{flat}, changes slowly in the $r$-direction, and is close to $\om(0, z)$.

\vs{0.05in}

\paragraph{\bf Lifting 1D singularities}

 A common approach to studying potential singularity formation 
 for \eqref{euler0} is through simplified 1D \emph{nonlocal} models.  This line of work has a long history, going back to the Constantin--Lax--Majda model \cite{CLM85}. %
Subsequently, many models have been proposed, including De Gregorio \cite{DG90}, CCF \cite{cordoba2005formation}, gCLM \cite{OSW08}, CKY \cite{choi2015finite}, and HL \cite{choi2014on} models. In several cases, the model is derived \emph{exactly} by restricting the fluid equations to the axis with a closure ansatz; see Section~\ref{sec:intro_further}.

Although many 1D nonlocal fluid models develop singularities from \emph{smooth data}, including 
all the above-mentioned models on the real line, these singularities have not led to singularities for incompressible fluid equations in a smooth domain. \footnote{
See the works \cite{elgindi2017finite,elgindi2018finite} for singularity formation in domains with corners.
}
We also note that in the perturbative regime \(\al\ll 1\), 3D Euler with \(C^{1,\al}\) velocity can be approximated by a 1D model or an ODE system with essentially \emph{local} Riccati-type dynamics \cite{elgindi2019finite,cordoba2023finite},
\footnote{
  For $\al \ll 1$, after several highly nontrivial transforms, \eqref{euler0} can be approximated by 
 a 1D nonlocal PDE, which is further reduced to a Riccati PDE $\pa_t f(t, x) = f(t, x)^2$ by introducing a suitable variable  \cite{elgindi2019finite}. The approximation error is of size \(O(\al)\).
 }
 which drive the blowup from \(C^{1,\al}\) data.

 In view of the highly nonlocal nature of 3D Euler, this naturally leads to the following questions:

\begin{ques}\label{ques:1D_model}
Can one lift a singularity from a \emph{genuinely nonlocal} 1D fluid model to an incompressible fluid equation in $\R^2$ or $\R^3$ with \(C^{1,\al}\cap L^2\) velocity---for example, Euler, Boussinesq, SQG, or IPM? More importantly, can this be achieved without requiring small $\al$? 
\end{ques}

There are two fundamental obstructions to such a lifting: 
many known connections between 1D models and fluid equations produce \emph{infinite-energy} solutions, and the velocity depends globally on the solution, whereas a 1D model captures only \emph{lower-dimensional} information.

\subsection{Main result}

Motivated by the above questions and the work \cite{hou2024potential}, we study the axisymmetric Euler \eqref{eq:Euler} without swirl. 
Our results give affirmative answers to Questions \ref{ques:SS}--\ref{ques:1D_model} and to Question \ref{ques:blowup} for any \(\al\in(0,\tfrac13)\), and provide a rigorous justification of the Hou-Zhang scenario in several respects. All arguments in this paper are purely analytic.

We first recall the parameters \(\kp\) and \(\kp_1\) from \cite{chen2026eulerI} 
\beq\label{def:kp} 
 \cff \teq 2^{\f43 } -  \f23 , \ \
   \kp \teq \f{\cff + 2}{4} < 1, \ \  
    \kp_1 \teq \f{1-\kp}{1000} ,  
  \quad  \kag = \f{1}{1000}, 
    \quad   \epa =  \kp_1 \e  - c_{\wwwa} \e^{2-\kp} ,
\eeq
with $\cff \approx 1.85$, where \(c_{\wwwa}\) is an absolute constant, independent of \(\e\), chosen in \cite{chen2026eulerI}.

Our first main result constructs \(C^{\al}\) self-similar blowup profiles 
for the vorticity and establishes their connection to $C^{\infty}$ 1D profiles.

\begin{thm}\label{thm:self-similar}

For any \(\al \in (0,\tfrac13)\), there exists a nontrivial \(C^{\al}\) self-similar blowup profile \(\omthss \not \equiv 0\) for the vorticity of the axisymmetric Euler without swirl \eqref{eq:Euler}. Moreover, there exists \(\beps \in (0, \f13) \) sufficiently small such that, for any \(\al \in (\tfrac13-\beps,\tfrac13)\), the following results hold.
\begin{enumerate}[label=(\roman*),leftmargin=1.5em,font=\normalfont]

\item \emph{\bf Blowup scaling}.
Let $\e = \f13 - \al$. The Euler equation \eqref{eq:Euler} admits a self-similar blowup solution
\bseq\label{eq:thm_SS}
\beq
  \ww_{\al}^{\th}(t, \xx) = \f{1}{1 -t} \omthss \big( \f{\xx}{ (1-t)^{ \cxs }} \big) , 
\eeq
with scaling exponent
\beq
\cxs \asymp \e^{-1},  \quad \cxs =  \f89 \e^{-1} + O( \e^{ -\hk}),  
\eeq
where $\kp \in (0,1)$ is defined in \eqref{def:kp}.

\eseq

\vs{0.05in}

  \item \emph{\bf Anisotropic regularity and decay}. 
 The vorticity profile admits the decomposition \\ $\omthss(\xx) = \AAa r^{\al } \ws(r, z) \in C^{\al}(\R^3)$
 for some constant $\AAa \asymp 1$, where
 $\ws$ satisfies higher regularity $\ws \in C^{1,\al}(\bar B_R)$ for any $R>0$
 with $\bar B_R =\{ (r,z):|(r,z) | \leq R\}$, and the following estimates hold
  \beq\label{eq:thm_profi_decay:a}
  |\AAa - \AAaa |  \les \e^{ \mhk}, \quad 
    |\bar \Om_{\al}(\xx)| \les_{\e}  \ang \xx^{-\als}, 
 \quad 
   |\na \bar \Om_{\al}(\xx)| \les_{\e}  \ang \xx^{-\als-1} .
  \eeq
Here $\alb = \f13,  \AAaa \asymp 1$ is an $\alpha$-independent constant defined in \eqref{def:AAa}. Moreover, the decay rate $\als$ and the scaling $\cxs$ in \eqref{eq:thm_SS} satisfies
  \beq\label{eq:thm_profi_decay:b}
\als = \f{1}{3} + \f{\e}{8} + O( \e^{1 + \mhk}), \quad \als > \f13, \quad \cxs = \f{1}{ \als - \al}.
\eeq

\vs{0.05in}
\item \emph{\bf Limiting behavior and asymptotic factorization}.
Let $\wwwa$ be the $C^{\infty}$ $\f13$-profile for the 1D profile equation 
\eqref{eq:1D_dyn} constructed in \cite{chen2026eulerI}; see Theorem \ref{thm:reg_alb}. 
The 3D profile $\ws$ is close to the 1D profile $\wwwa$ and satisfies
\beq\label{eq:thm_global_conv}
 \|  \ang z^{ \f13  + \kag - C \e^{2-\kp}}  \ang {|\xx|} ^{-\kag - \epa} | \log (2 + |\xx|) |^{-\kp}  (\ang z^{ -\b(\al) }  \ws - \wwwa) \|_{L^{\infty}} \les \e^{ 1 - 5 \kp_1} ,
\eeq
for some  $ \b(\al) = - \f{\e}{8} + O( \e^{2-\kp})$ and absolute constant $C>0$, where $|\xx| = |(r, z)|$, and $\kp_1, \kp, \epa , \kag$ are defined in \eqref{def:kp}.  In particular, as $\al \to (\f13)^-$, 
$ \ang z^{ - \b(\al)} \ws $ converges strongly to $\wwwa$ in the above norm;
 $ \ws$ converges to $\wwwa$, and $\omthss$ asymptotically factorizes as a tensor product:
\footnote{
The factor $\AAaa$ defined in \eqref{def:AAa}  comes from some normalization conditions in constructing the 1D and 3D profiles. See \eqref{eq:profi_3D_rescale}-\eqref{eq:profi_3D0}. 
To simplify some presentations for the 1D profiles, we keep this factor.
}
\beq\label{eq:thm_est_conv}
 \lim_{\al \to (1/3)^- } \ws(r, z) =  \wwwa(z), 
   \quad    \lim_{\al \to (1/3)^- } \omthss(r, z) = \AAaa r^{1/3} \wwwa(z),
   \quad  \alb = \tf13,
\eeq
where both convergences are uniform on compact subsets.

\end{enumerate}

\end{thm}

Result (ii) shows that $\omthss$ is more regular and decays faster in $z$ than
in $r$: locally, for fixed $r> 0$ one has
$\omthss(r,\cdot)\in C^{1,\al}$, while for fixed $z$ one only has
$\omthss(\cdot,z)\in C^\al$.

Due to slow spatial decay, the profile $\omthss$ does not have finite energy, i.e. $\| \uu \|_{L^2} = \infty$. Our second main result establishes asymptotically self-similar blowup
from finite-energy initial data.

\begin{thm}\label{thm:main_blowup}

There exists $ \al_0 < \f13$ such that the following holds for any $\al \in (0, \f13 )$, any $ \delin >0$, and $\g=\max(\al,\al_0)$. There exists compactly supported initial vorticity $\ww_0^{\th} = r^{ \g }g \in C_c^{\g}(\R^3) \subset C_c^{\al}(\R^3)$, with $ g \in C_c^{\infty}(\R^3)$ and the associated initial velocity $\uu_0 \in C^{1,\al}(\R^3) \cap L^2(\R^3)$, such that
\beq\label{eq:init_ass}
 \nlinf{ \om_0^{\th} - \omthgg } \les_{\g}  \delin ,
\eeq
 and the corresponding local solution to \eqref{eq:Euler} develops an asymptotically self-similar blowup in finite time $T < \infty$: 
\footnote{
The exponent \(0.01\) in \((T-t)^{0.01}\) can be tracked and easily improved, although we do not attempt to optimize it.
}
\beq\label{eq:main_SS_blowup}
  \| (T - t) \cdot \, \om^{\th} \big( t, \, \cE(t) \cdot (T-t  )^{\cxg} \xx \big) - 
     \omthgg(\xx)  \|_{ L^{\infty} }  \les_{\g} (T-t)^{0.01} \delin ,  \ \forall  \, t \in [0, T),
\eeq
where $\cxg,  \omthgg \not \equiv 0 $ are the scaling and profile in Theorem \ref{thm:self-similar}, and  $\cE(t)$ is some function depending on the initial data and satisfies $ \cE(t) \in C^1([ 0, T)) , \cE(t) \asymp_{\g} 1$. In addition, we have 
\beq\label{eq:thm_blowup_time}
    |T - 1 | \les_{\g}   \delin, \quad  \nlinf{ \cE - 1 } \les_{\g}   \delin. 
\eeq
In particular, for fixed $\al$ and hence fixed \(\gamma\), by taking $\delin$ small enough, the initial data $\om_0^{\th} \in C_c^{\g} \subset C_c^{\al}$, the blowup time $T$, and $\cE(t)$ can be made arbitrarily close to $ \omthgg$, $1$, and $1$, respectively.

\end{thm}

Since $C^{\al_2} \subset C^{\al_1}$ for any $0< \al_1 < \al_2 <1$, 
it suffices to consider $\al \in [\al_0, \f13)$ with $\al_0$ close to $\f13$ in Theorem \ref{thm:main_blowup}. 
For $\al < \al_0$, we take $\gamma = \alpha_0$
and the inclusion $C^{\gamma}_c \subset C^{\alpha}_c$ ensures $\omega_0^\theta \in C^\alpha_c$.

Below, we first record a few remarks on the main theorems, and then discuss several aspects of the results and their proofs.

\begin{remark}[\bf Regularity of initial data]\label{rem:reg_init}
 The low $C^{\al}$ regularity of the initial data $\om^{\th}_0$ in  Theorem \ref{thm:main_blowup} 
with $\al < \f13 $ comes only from the vanishing order near $r=0$;
away from \(r=0\), \(\om_0^{\th}\) is \(C^\infty\).  In view of the global regularity result \cite{danchin2007axisymmetric} for initial vorticity $\om_0^{\th} \in L^{\infty}$ satisfying $ \f{1}{r} \om_0^{\th} \in L^{3, 1}$, 
the decisive  condition  for global regularity \emph{is not the regularity} itself, but the \emph{vanishing order} near $r=0$ ensuring $ \f{1}{r} \om_0^{\th} \in L^{3, 1}$. The initial data $\om_0^{\th}$ in Theorem \ref{thm:main_blowup} violate this condition and attain the sharp vanishing order 
near $r=0$ for blowup. In particular, $\om_0^{\th}$ is $C^{\infty}$ near the $z$-axis for $r>0$.

\end{remark}

\begin{remark}[\bf Set of initial data]\label{rem:set_init}

There exists an open set $X_1$ (a ball in a weighted $C^{1,\al}$ space) 
and a finite-dimensional subspace $X_2 \subset C_c^{\infty}$ such that 
the initial data $g$ in Theorem \ref{thm:main_blowup} may be taken as $ g = g_1 + g_2$, 
where \(g_i\in X_i\) for $i=1,2$ and \(g_1, g_2\) have compact support; see  Remark \ref{rem:init_data2}. 
One may prescribe $g_1\in X_1$ arbitrarily within the symmetry class,  including \(C_c^\infty\) functions,  while $g_2$ is determined by a fixed point argument to avoid potentially unstable modes of the profile. 
We also note that \cite[Section 5.1]{hou2024potential} presents some numerical evidence 
suggesting that the blowup profile for $r^{-\al} \om^{\th}$ becomes less stable  as \(\al\) increases. Thus unstable modes may appear as $\al\to(\f13)^-$.

\end{remark}

\begin{remark}[\bf Convergence of 3D profiles in stronger norm]

We estimate the weighted $C^1$ norm for $\ws - \waa$ in Theorem \ref{thm:3D_solu} 
 and the weighted $C^k$ norm for $\waa , \wwwa$ in Theorems \ref{thm:reg_alb}, \ref{thm:1D_profile_prop}  with \emph{$\e$-independent} constants. 
 Interpolating these estimates and \eqref{eq:thm_global_conv} in Theorem 
\ref{thm:self-similar}, we obtain convergence of $\ws \to \wwwa$ in a weighted $C^{\g}$ norm for any $\g \in (0,1)$.

\end{remark}

\paragraph{\bf Complete characterization of limiting behavior as $\al \to (\f13)^-$}
\label{para:complete_char}

As $\al \to \f13$, by taking $\delin \ll_{\al} 1$, the asymptotically self-similar blowup solution constructed in Theorem \ref{thm:main_blowup}, after rescaling, remains sufficiently close to the blowup profile $\omthss$, which asymptotically factorizes and converges strongly to the tensor product $ \AAaa r^{1/3} \wwwa(z)$ in some weighted $L^{\infty}$-norm by Theorem \ref{thm:self-similar}.
\footnote{
 Since $|\b(\al)| \les \e$ by item (iii) in Theorem \ref{thm:self-similar},
 using $ |\ang z^{ -\b(\al) } - 1| \les \e \log(2 + \ang z) \ang z^{|\b(\al)|} $, $|\wwwa(z)| \les \ang z^{-1/3}$ from \eqref{eq:bw_smooth} and estimate \eqref{eq:thm_global_conv}, and imposing a fast decaying weight $\rho$, we obtain $\nlinf{ (\omthss - \AAaa \wwwa) \rho } \to 0$ as $\al \to \f13$. We present \eqref{eq:thm_global_conv} since it captures the sharp convergence estimates.
}
The 1D limiting profile \(\wwwa\) solves the 1D profile equation 
\eqref{eq:1D_dyn} with $\al = \f13$ 
and is constructed in \cite{chen2026eulerI}. See Theorem \ref{thm:reg_alb}. Theorem \ref{thm:self-similar} precisely characterizes the leading order behavior of the blowup rate $\cxs$ as $\al \to \f13$.

\vs{0.1in}

\paragraph{\bf Hou-Zhang scenario}

Our results rigorously justify the Hou-Zhang scenario \cite{hou2024potential} in several respects. 
First, we justify this scenario  by constructing blowup solutions to \eqref{eq:Euler2} of the form $\om^{\th} = r^{\al} g$, as proposed in \cite{hou2024potential}, with smooth initial data for $g$.
 Second, in Proposition \ref{prop:LWP}, by analyzing the 
 self-similar equation for $g \teq r^{-\al}\om^{\th}$, we show that, for general initial data $g_0(r, z) \in C^{1,\al}(\R_+\times \R)$, the solution $g$ propagates $C^{1,\al}$ regularity, which is \emph{one order} higher than that of $\om^{\th}$. This makes high-accuracy numerical computation for $g$  feasible. 
Third, 
as $\al \to \f13$, Theorems \ref{thm:self-similar} and \ref{thm:main_blowup} justify the formation of a profile with a 1D structure in  \cite{hou2024potential}, and show that $g$ eventually becomes a function of $z$ only, thereby proving a conjecture in \cite[Section 8]{hou2024potential} with 
$\al^{*} = \f13$. 
Fourth, Theorem \ref{thm:self-similar} shows that the blowup rate $\cxs \to \infty$ as $\al \to \f13$, proving a conjecture made in 
\cite[Section 5.1]{hou2024potential}.

\subsubsection{\bf Lifting 1D singularities}

Following \cite{hou2024potential}, we derive a gCLM-type 1D model  \eqref{eq:1D_dyn}
in  $z \in \R$ along the axis $r=0$. 
For $\al< \f13$ sufficiently close to $\f13$, we extend the blowup profiles $ \waa$ for the 1D model constructed in \cite{chen2026eulerI} 
constantly in $r$ as an approximate 3D profile $\wwb(r, z) \equiv \waa(z)$ and 
construct 3D profiles $\ws$ in Theorem \ref{thm:self-similar} around $\wwb$.
\footnote{
Since $C^{\g} \subset C^{\al}$ for any $\g \in [0, \al]$, Theorem \ref{thm:self-similar} constructs $C^{\al}$ self-similar profiles for the \emph{entire range} $\al \in (0, \f13)$.
}
We then use these 3D profiles to construct blowup solutions in Theorem \ref{thm:main_blowup}. See Section \ref{sec:idea} for further discussion.
Together with the results of limiting behavior established above,
these results give an affirmative answer to 
 Question \ref{ques:1D_model} 
and Question \ref{ques:SS} on the existence of self-similar profiles 
and blowup solutions with $C^{\al}$ vorticity,  and their limiting behaviors as $\al\to (\f13)^-$.
 The blowup construction is illustrated in Figure \ref{fig:lift}.
 \footnote{
 \label{foot:direct_Wa}
A direct construction of $\ws$ from the $\f13$-profile $\overline W_{1/3}$ is not available: the spatial blowup rate of $\overline W_{1/3}$, analogous to $\cxs$ in \eqref{eq:thm_SS}, is \emph{infinite}, whereas $\ws$ has a \emph{finite} spatial blowup rate. 
The intermediate 1D profile $\waa$ plays an important role in capturing 
the spatial blowup rate and far-field decay  in $z$-direction of the 3D profile $\ws$.
 }

To the best of our knowledge, our results provide the first example in which a singularity from a 1D nonlocal fluid model is lifted to construct blowup for incompressible fluid equations in \(\R^2\) or \(\R^3\). Moreover, we \emph{do not} require \(\al\) to be small. More broadly, our results provide a new approach to construct blowup by lifting singularities from 1D models.

\vs{0.05in}
\paragraph{\bf Mechanism: anisotropic flow induced by incompressibility}

The key mechanism enabling the lifting is the anisotropic 
structure of the self-similar flow associated with the 1D $\al$-profile $\waa$:
\beq\label{eq:intro_SS_flow}
\bar Q^r \teq  \bar c_{l, \al} r + \bar U^r(r, z)  , \quad  \bar Q^z \teq \bar c_{l, \al} z + \bar U^z(r, z)  , \quad  \tf{1}{r} \bar Q^r, \ \bar c_{l, \al}  \asymp ( \tf13 - \al)^{-1},
\quad  \tf{1}{z} \bar Q^z \asymp 1 , 
\eeq
in the main part of the domain. 
Here, $\bar \UU$ is the self-similar velocity and $\bar c_{l, \al}$ is the spatial scaling exponent associated with $\waa$. See Section \ref{sec:appr_profi}.  In the \(z\)-direction, 
for $z>0$, \(\bar U^z<0\) almost cancels the scaling field \(\bar c_{l,\al}z\). In the \(r\)-direction, due to incompressibility, \(\bar c_{l,\al} r\) and \(\bar U^r\) have the same sign and comparable size, producing a flow with \( \f{1}{r} \bar Q^r\) much larger than $ \f{1}{z} \bar Q^z$ as $\al \to \f13$.

\vs{0.1in}
\paragraph{\bf Key ingredients}

Note that the 1D profile $\waa$ neither decays in the \(r\)-direction nor solves 
the 3D profile equation for \(r>0\). The lifting construction therefore relies on three key ingredients:

\begin{itemize}[leftmargin=1em]

\item Fast transport of the perturbation to spatial infinity by the \(r\)-advection;

\item  
Angular averaging of the Biot--Savart law
in \((r,z)\)-variables,  transferring \(z\)-decay to \(|(r,z)|\)-decay;

\item Integration by parts along trajectories, exchanging  an \emph{unbounded} \(r\pa_r\) derivative for a \(z\pa_z\) derivative by exploiting the profile equation \emph{twice}, and overcoming
the key unbounded factor $\frac{\ang{r,z}}{\ang{z}}$ for $r\gg|z|$.

\end{itemize}

We develop a family of anisotropic weighted estimates to capture the first two  mechanisms,
and prove  Lemma \ref{lem:IBP} for the third ingredient. 
See Section~\ref{sec:idea}  for the proof ideas for Theorem \ref{thm:self-similar}.

\begin{figure}[t]
\centering

\captionsetup{width=0.9\textwidth}

\begin{tikzpicture}[>=stealth, line width=0.7pt]
\coordinate (O) at (0,0);
\coordinate (A) at (1.6,  0.5);
\coordinate (B) at (8.6, 0.5);
\coordinate (C) at (1.6, 2.8);
\coordinate (D) at (8.6,2.8);

\draw[->] (O) -- (10.0,0) node[right] {$\alpha$};
\draw[->] (O) -- (0, 3.2) node[above] {};

\node[rotate=90] at (-0.8, 1.6) {Dimen.};

\node[left=0pt] at (0, 0.5) {1D};
\node[left=0pt] at (0, 2.8) {3D};

\draw (1.6,-0.08) -- (1.6 ,0.08);
\node[below=0pt] at (1.6 , 0) {$\frac13-\varepsilon$};

\draw (8.6,-0.08) -- (8.6,0.08);
\node[below=0pt] at (8.6,0) {$\frac13$};

\node (Omega) at (C) {$\overline{\Omega}_{\alpha}(r, z) \in C^{1,\al}$};
\node (Wtop)  at (D) {$\wwwa(z) \in C^{\infty}$};

\node (Wleft) at (A) {$\overline{W}_{\alpha}(z) \in C^{\infty}$};
\node (Wbot)  at (B) {$\wwwa(z) \in C^{\infty} $};

\draw[->, shorten >=10pt, shorten <=10pt]
  (Omega) -- (Wtop)
  node[midway, above= 1pt] {Limit, $\alpha\to (\f13)^{-}$};

\draw[->, shorten >=4pt, shorten <=4pt]
  (Wleft) -- (Omega)
  node[midway, right=6pt] {Lift};

\draw[->, shorten >=10pt, shorten <=10pt]
  (Wbot) -- (Wleft) 
  node[midway, above= 1pt] {perturb};

\node at (8.6, 1.6) {\Large $\Vert$};

\end{tikzpicture}

\vspace{-1.0em}
\caption{Lifting construction and limiting behavior as $\alpha\to (\frac13)^-$.
$\ws$ is the blowup profile for $r^{-\al} \om^{\th}$. $\waa$ and $\wwwa$ are the blowup profiles for the 1D model, and $\alb =\tf13$.
}

\vspace{-1.2em}
\label{fig:lift}

\end{figure}

\vs{0.05in}

\paragraph{\bf Breakdown at $\al=\frac13$}

The construction of 3D self-similar profiles breaks down at $\al = \f13$. As $\al \to (\f{1}{3})^-$, from \eqref{eq:thm_global_conv}, \eqref{eq:thm_est_conv}, the profile $\omthss   \to \AAaa r^{\alb} \wwwa(z)$ develops a heavy tail in $r$ and no longer decays in $r$.
 Consequently, the associated
  velocity with correction near $0$ develops logarithmic growth and $\na \UU(0)$ diverges 
\beq\label{eq:vel_mild}
  | \na \UU(\ws) - \na \UU(\ws)(0) | \asymp \min( \log ( |(r, z)| + 2), \e^{-1}),
  \quad  |\na \UU(0)| \asymp  \e^{-1}, \quad \e = \tf13 - \al.
\eeq

The constant in the decay estimates \eqref{eq:thm_profi_decay:a} for $\omthss$ depends on $\e = \f13 - \al$. In fact, as $\al \to \f{1}{3}$, the profile  begins to decay in $r$ only for $r  \gtr c_{\e}$,  where $c_{\e} \to \infty$.   The scaling exponent $\cxs$ in \eqref{eq:thm_SS} also blows up. 
These singular structures prevent the construction of $C^{ 1/3 }$ profile for \eqref{eq:Euler} 
with Lipschitz velocity, thereby explaining the mechanism in Question \ref{ques:SS}.

\vs{0.05in}

The methods developed in this work may help guide the rigorous analysis of 
 recent numerical findings on 3D Euler \cite{huang2026novel,HouHuang2022},
 where simulations suggest that 3D Euler may develop self-similar blowup solutions with singular structure or anisotropic scaling from \emph{smooth} initial data.

\subsubsection{\bf Criticality, regularity, and singularity}

The blowup results in Theorem \ref{thm:main_blowup} are fundamentally connected to the 
criticality of 3D axisymmetric Euler without swirl \eqref{eq:Euler}. Due to the conservation of $ \| \tf{1}{r} \om^{\th}(s) \|_{L^{p, q}}$ by 
\eqref{eq:Euler} 
for all $p \in [1, \infty), q \in [1, \infty]$, the Lorentz spaces $L^{p, q}$
 can be  classified according to the criticality of this conserved quantity under the spatial scaling $\om_{0, \lam}^{\th}(\xx) \teq \om_0^{\th}( \lam \xx)$:\footnote{
3D Euler \eqref{eq:Euler} satisfies the spatial scaling symmetry: if \( \om\) solves 3D Euler, so does $\om_{\lam}(t,\xx) \teq \om(t,\lam \xx)$ for any \(\lam>0\).
}
 \[
  \mw{subcritical}:   p > 3;   \quad \mw{critical} :  p =  3 ; \quad  \mw{supercritical} :  p \in [1, 3) ,
\]
with any $q \in [1, \infty]$. 
Under some mild assumption, e.g. $\om_0^{\th} \in C_c^{\d}$ with some $\d >0$, the result \cite{danchin2007axisymmetric,shao2026global}  implies global regularity in 
every critical and subcritical space $L^{p, q}, p \geq 3, q \in [1, \infty]$.

\vs{0.1in}

\paragraph{\bf Self-similar blowup in any supercritical $L^{p, q}$ spaces}

Theorem \ref{thm:main_blowup} shows that for \emph{every} supercritical \(L^{p,q}\) space with \(p\in[1,3)\)
and any $\al \in (0, \f13)$, there exists initial data with $\om_0^{\th} \in C_c^{\al},\f{1}{r} \om_0^{\th} \in L^{p, q}$ 
\footnote{
We note that $\dot C^{1/3}$ is not invariant under the scaling 
$ f(\xx) \to f_{\lam}(\xx) = f(\lam \xx)$ for an axisymmetric function.
}
that leads to  asymptotically self-similar blowup of \eqref{eq:Euler}.
 Therefore, the blowup result in Theorem \ref{thm:main_blowup} is sharp both in the H\"older regularity of $\om_0^{\th}$
 and in the $L^{p, q}$ regularity of $\f{1}{r}\om_0^{\th}$.

A key idea for obtaining such a sharp self-similar blowup result is to perturb
from a \emph{non-blowup profile}; see Figure~\ref{fig:lift}. Here this profile
is the 1D limiting profile $\wwwa$. This approach was developed in our earlier
works on the gCLM, De Gregorio, and Landau models
\cite{chen2020slightly,chen2021regularity,chen2023nearly}.

\vs{0.1in}
\paragraph{\bf Criticality characterizes the regularity threshold.}

Together with \cite{shkoller2026incompressible,danchin2007axisymmetric,shao2026global}, this gives a \emph{sharp} classification of the regularity regimes under the mild assumption \(\om_0^{\th}\in C_c^{\d}\), \(\d>0\). 
The axisymmetric Euler without swirl \eqref{eq:Euler} is globally regular when \(r^{-1}\om_0^{\th}\) belongs to a critical or subcritical \(L^{p,q}\) space, and in particular when \(\om_0^{\th}\in C^{\al}\) with \(\al\ge \tfrac13\). 
In contrast, for any \(p\in[1,3)\) and \(\al<\f13\), it admits blowup from some initial data \(\om_0^{\th}\) such that \(r^{-1}\om_0^{\th}\) belongs to a supercritical \(L^{p,q}\) space and \(\om_0^{\th}\in C^{\al}\).
Thus, the criticality of the conserved quantity \(\|r^{-1}\om^{\th}(s)\|_{L^{p,q}}\) 
\emph{precisely} characterizes the regularity threshold for blowup versus global regularity.

\subsubsection{\bf Outgoing self-similar flow \(\Rightarrow\) finite codimension stability}

Let $\QQs$ be the self-similar flow associated with the 3D profile 
\eqref{eq:intro_SS_flow}. The profile satisfies a crucial global outgoing property:
\beq\label{eq:intro_outgo0}
\QQs \cdot \xx \geq \lam_Q |\xx|^2,
\quad 
\forall \ \xx = (r, z) .
\eeq
In the linearized self-similar equation, this condition implies that the perturbation is transported 
from the origin to spatial infinity, generating a strong stability mechanism. 
This stability mechanism has played a fundamental role in recent stability analysis for 3D incompressible Euler \cite{ChenHou2023a} and compressible Euler \cite{chen2026implosion,chen2024Euler,chen2024vorticity}.

In \cite{ChenHou2023a}, by exploiting this mechanism using singularly weighted estimates, Chen-Hou proved stability of the linearized operator up to a finite rank perturbation. 
They further developed a crucial \emph{quantitative} finite rank perturbation method with computer assistance to prove \emph{full} stability. Since no smooth \emph{exact} 3D Euler profile is available, in  \cite{ChenHou2023a}, computer assistance plays a key role in rigorously constructing approximate solutions, which provides crucial small parameters for the blowup analysis. This naturally leads to the following fundamental question:

\vs{0.05in}

\emph{Under the \emph{outgoing} condition \eqref{eq:intro_outgo0}, is an \emph{exact} self-similar blowup profile for 3D incompressible Euler, with suitable regularity and decay, nonlinearly stable up to a finite-dimensional instability?}

\vs{0.05in}

One main contribution of this work is a new analytic framework for finite codimension stability  without imposing high regularity or fast decay on the Euler profile.
\footnote{
The 3D profile $\omthss$ constructed in Theorem \ref{thm:self-similar} 
only has $C^{ \al}$ regularity but not $C^{1/3}$ 
and $\ws$ has $C^{1,\al}$ regularity but not 
$C^{1,1/3}$. Moreover, $\omthss$ and the gradient of the velocity profile 
decay only at the rate $\ang \xx^{- c \cdot \e }$ with $\e \ll 1$. 
}
The key observation is that the nonlocal operator in 3D Euler is \emph{essentially} local up to a \emph{compact operator}, since the Biot--Savart kernel is smooth away from the diagonal \(y=x\). 
\footnote{
In \cite{ChenHou2023a,chen2025singularity}, this property has been exploited by showing that the commutator between the weight and nonlocal operator is more regular.
}
By exploiting the outgoing property and singularly weighted estimates, we establish nonlinear finite codimension stability of the 3D profile constructed in Theorem \ref{thm:self-similar} and prove the blowup results in Theorem~\ref{thm:main_blowup}. 
The \emph{exact}  profile allows us to use several \emph{qualitative} functional analytic arguments to control the compact perturbation. See Section~\ref{sec:idea_main_blowup} for further discussion.
In the approximate-profile setting of \cite{ChenHou2023a},  
the corresponding step necessarily requires \emph{quantitative} estimates with \emph{computable} constants, which is achieved through the \emph{quantitative} finite rank perturbation method \cite{ChenHou2023a}.

Finally, we remark that constructing an \emph{exact} Euler self-similar blowup profile in $\R^3$ with higher regularity remains a fundamental open problem.

\subsection{Further questions}\label{sec:intro_further}
Following the present work, several related questions remain open. 
First, one may extend the blowup construction to other fluid equations by 
lifting 1D singularities; see the discussion below for further details.
Second, one may extend the construction of \(C^{\al}\) vorticity profiles beyond \(\al=\tfrac13\) by introducing swirl. It is also of interest to investigate whether \(C^{1/3+}\)  
vorticity  profiles with swirl can be obtained by perturbing the profiles in Theorem \ref{thm:self-similar}.

\vs{0.05in}

\paragraph{\bf{1D models and fluid equations}}

Suppose that \(\om(\xx,t)\) solves a 2D or 3D fluid equation with axisymmetry, and let \(\xx_h=x_1\) for \(d=2\), or \(\xx_h=(x_1,x_2)\) for \(d=3\). A canonical way to derive a 1D model along \(\xx_h=0\) is to consider solutions of the form $\om(\xx,t)=|\xx_h|^a f(x_n),$
and then compute the nonlocal velocity \(\uu(\xx)\) along \(\xx_h=0\), which depends only on \(f(x_n)\). The 1D model \eqref{eq:1D_dyn} derived from \eqref{eq:Euler2} falls into this class. There are many such examples, and we list a few below.

\begin{enumerate}[label=(\roman*), leftmargin=1.2em]

\item \textsl{Hou-Luo model \cite{choi2014on} and Euler/Boussinesq} : 
Using the ansatz  $\om(x, y, t) =f(x, y, t)$, $\th(x, y, t) = g(x, t)$ and 
 the Boussinesq in $\R^2_+$ for $(\om, \th)$, one derives the 1D HL model for $(f, g)$ \emph{exactly} on the boundary $y= 0$. 
A simplification of the HL model leads to the CKY model \cite{choi2015finite}. 

 \item \textsl{CCF model \cite{cordoba2005formation} and IPM}.
Using the ansatz $\th(x, y, t) = f(y, t)$ and the incompressible porous media (IPM) equation in $\R^2_+$ for $\th$, one derives the CCF model  $f(y, t)$ on the boundary $y= 0$. 

\item \textsl{De Gregorio model \cite{DG90} and SQG}. 
Using the ansatz  $\om(x, y, t) = x  f(y, t)$ and the SQG equation in $\R^2$, 
one derives the De Gregorio model on $y= 0$.

\end{enumerate}

These 1D models all develop finite-time self-similar singularities from smooth
initial data \cite{chen2019finite,chen2021HL,huang2024self,huang2025exact}. See also
related works on gCLM models \cite{chen2020singularity,Elg17}.

\subsection{Organization}

The rest of the paper is organized as follows.
In Section \ref{sec:setup}, we introduce the setup, the fixed point map for the 3D profile, 
and outline the proof ideas in Section \ref{sec:idea}. In Section~\ref{sec:nonlocal}, we develop  anisotropic weighted estimates for the Biot--Savart law.
In Section~\ref{sec:est_LQR}, we use these estimates to control the nonlocal and error terms in the fixed-point map, and obtain \emph{qualitative regularity estimates}.
In Section~\ref{sec:3D_EE}, we perform energy estimates for the fixed-point map \(\FFF\). 
In Section~\ref{sec:3D_solu}, we prove \emph{qualitative} estimates for \(\FFF\), construct the 3D profiles, estimate their decay and regularity, and prove Theorem~\ref{thm:self-similar}.
In Section~\ref{sec:blowup}, we prove finite codimension stability of the 3D profiles and prove Theorem~\ref{thm:main_blowup}. Finally, we prove basic functional inequalities and estimates in Appendix~\ref{app:basic}, and establish the local well-posedness theory in Appendix~\ref{app:LWP}.

\vs{0.05in}

\paragraph{\bf Suggested reading order}

While the paper is organized in logical order, we note that Section~\ref{sec:nonlocal}
(except for Section~\ref{sec:1D_3D}), Section~\ref{sec:qual_reg}, and Section~\ref{sec:3D_profile_prop} contain technical estimates supporting the main construction.
The absolute constants in the estimates of Sections~\ref{sec:BS_est_ker}--\ref{sec:iso_nonlocal}
are not needed explicitly in the later analysis, while the estimates in Sections~\ref{sec:qual_reg}
and~\ref{sec:3D_profile_prop} are largely \emph{qualitative}, except for the scaling terms
and outgoing conditions in Section~\ref{sec:3D_profile_prop}.
We therefore suggest that readers may first focus on the main construction and return to these sections as needed.

The core section is Section~\ref{sec:3D_EE}, which mainly relies on the quantitative
estimates in Sections~\ref{sec:est_LQR}--\ref{sec:err}.
We emphasize that Section~\ref{sec:nonlocal}, Sections~\ref{sec:est_LQR}--\ref{sec:3D_EE},
and Sections~\ref{sec:3D_solu}--\ref{sec:blowup} can be read mostly independently of one another.

\section{Fixed-point problem for 3D profiles and proof ideas}\label{sec:setup}

In this section, we recall the 1D profiles from \cite{chen2026eulerI}, formulate
the fixed-point problem for the 3D profile, and outline the main ideas of the proof 
in Section \ref{sec:idea}.

\subsection{Axisymmetric Euler and Biot-Savart law}
In the axisymmetric setting, given the angular vorticity $\om^{\th}$ in \eqref{eq:Euler}, we obtain the angular stream function 
$\psi^{\th}(r, z)$ and velocity $(u^r, u^z)$ via
\beq\label{eq:Euler_b}
\bga
 - (\pa_{rr} + \f{1}{r} \pa_r + \pa_z^2 + \f{1}{r^2} \pa_{\th}^2  )( \psi^{\th} \cos \th) = 
 \om^{\th} \cos \th  , \quad
 u^r = - \f{1}{r} \pa_z (r \psi^{\th}),
 \quad u^z = \f{1}{r} \pa_r( r \psi^{\th} ) .
\ega
\eeq
 Following \cite{hou2024potential}, for any $\al \in (0, \f13)$, we introduce 
\beq\label{eq:Euler_change}
 \psi = \f{ \psi^{\th} }{r}, \quad  \om = \f{\om^{\th}}{r^\al} .
\eeq

Using the vorticity equation \eqref{eq:Euler} and \eqref{eq:Euler_b}, we derive the following equation on $\om$
\bseq\label{eq:Euler2}
\begin{gather}
 \pa_t \om + (u^r \pa_r + u^z \pa_z) \om = -(1-\al) \psi_z \om, \label{eq:Euler2:a} \\
 u^r(\psi) \teq - r \psi_z, \quad  u^z(\psi) \teq r^{-1} \pa_r(r^2 \psi) = 2 \psi + r \pa_r \psi, \label{eq:Euler2:c}
\end{gather}
\eseq
where $\psi$ solves the elliptic equation 
\bseq\label{eq:Euler2_psi} 
\beq\label{eq:Euler2_psi:BS0} 
   - (\pa_{rr} + \pa_{zz} )(r \psi)  -  \pa_r \psi =  \om r^{\al} , \quad   \psi \teq \BSa(\om) .
\eeq
We use $\psi = \BSa(\om)$ to denote the above \emph{original} Biot-Savart law.  For the absolute constant $\cpsia$ chosen in \eqref{def:cpsi}, we introduce the following Biot-Savart law 
\beq\label{eq:Euler2_psi:BS}
  \BS(\ww) \teq \cpsia \BSa(\ww). 
\eeq
\eseq

In this paper, we consider $\om, \psi$ odd in $z$, which is preserved by \eqref{eq:Euler}.
We use superscripts for the \(r,z\) components of vectors, use \(\th\) for variables before dividing out the \(r\)-power, and omit it afterward as in \eqref{eq:Euler_change}.
We refer to \(\psi\) as the stream function, since \(\psi^\th\) is rarely used.

\subsection{3D self-similar profile equation}
For some scaling exponent $\cxx > 0$ to be determined,  plugging the self-similar blowup ansatz 
\beq\label{eq:SS_ansatz_1}
\om^{\th}(t, x) = \f{1}{1 - t} \omth_*( \f{x}{ (1-t)^{ \cxx } } ) 
\eeq
to \eqref{eq:Euler}, we obtain the 3D self-similar profile equation
\beq\label{eq:profi_3D_omth0}
  ( \cxx r + U_*^r , \ \cxx z + U_*^z) \cdot \na_{r, z} \omth_* = ( - 1 - \pa_z \Psi_* ) \omth_* ,
  \quad \Psi_* = \BSa( r^{-\al} \omth_*), 
\eeq
where $\Psi_*, U_*^r, U_*^z $ are the stream function and velocity associated with $r^{-\al}\omth_* $ \eqref{eq:Euler_change}, \eqref{eq:Euler2}, \eqref{eq:Euler2_psi}. 
To impose normalization conditions, we introduce parameters \(c_{\om},c_l,\cpsia\) into \eqref{eq:profi_3D_omth0} and consider
\beq\label{eq:profi_3D_omth}
\bga
   ( c_l r + U^r , c_l z + U^z ) \cdot \na_{r, z } \omth = 
 ( \cw + \al c_l - \Psi_z ) \omth,  \\
   \Psi = \BS( r^{-\al} \omth) 
 = \cpsia \BSa( r^{-\al} \omth), \quad  
 \UU = \tf{1}{r} (-\pa_z, \pa_r)( r^2 \Psi) ,
 \ega
\eeq
with $\BS \neq \BSa$ defined in  \eqref{eq:Euler2_psi}. 
Given a solution to \eqref{eq:profi_3D_omth} with $\cw + \al c_l \neq 0$, using the relation \footnote{
We multiply $\Om^{\th}$  by the extra constant $\cpsia$ to recover the original Biot-Savart law $\BSa$
in \eqref{eq:profi_3D_omth0}. 
}
\beq\label{eq:profi_3D_rescale}
( \cxx, \  \omth_*, \ \Psi_*, \ \UU_* ) 
\teq  - ( c_{\om} + \al c_l )^{-1} \cdot ( c_l, \ \cpsia \omth, \  \Psi, \  \UU),
\eeq
and multiplying \eqref{eq:profi_3D_omth} by 
$ ( c_{\om} + \al c_l )^{-2}  \cpsia$, 
we obtain a solution to \eqref{eq:profi_3D_omth0}. 
 To distinguish between solutions of \eqref{eq:profi_3D_omth0} and \eqref{eq:profi_3D_omth}, we use the star notation $\omth_*$ to denote a solution of \eqref{eq:profi_3D_omth0}. 

Dividing 
\eqref{eq:profi_3D_omth}  by $r^{\al}$, 
we obtain the profile equation for $\Om = r^{-\al} \omth$:
\beq\label{eq:profi_3D0}
 ( c_l r + U^r , c_l z + U^z ) \cdot \na_{r, z } \www = 
 ( \cw - (1-\al) \Psi_z  ) \www, \quad \Psi = \BS(\www) , \quad  \UU = \tf{1}{r} (-\pa_z, \pa_r)( r^2 \Psi),
\eeq

While equation \eqref{eq:profi_3D0} and \eqref{eq:profi_3D_omth0} are equivalent, 
formally, 
we have two modulation parameters $c_l, c_{\om}$ in \eqref{eq:profi_3D0}, which make it convenient to impose vanishing conditions. 
Moreover, we introduce an extra free parameter  \(\cpsia\) in the 3D Biot--Savart law.
\footnote{
We have this degree of freedom since if we fix $c_l, \UU, c_{\om},\Psi$ in 
\eqref{eq:profi_3D_omth}, equation \eqref{eq:profi_3D_omth} is homogeneous in $\Om^{\th}$.
}
We call \(\BS\) in \eqref{eq:Euler2_psi:BS} the Biot--Savart law and \(\BSa\) in 
\eqref{eq:Euler2_psi:BS0} the \emph{original} Biot--Savart law; in the rest of the paper, we mainly use \(\BS\).

\vs{0.05in}
\paragraph{\bf Constant $\cpsia$ for $\Psi$}

To match the 3D and 1D stream function considered in 
\cite{chen2026eulerI} (see \eqref{def:psi_1D} 
and Lemma \ref{lem:vel_bc}), we choose  constants $\cpsi, \cpsib, \cpsia$, with $\cpsia$ used in \eqref{eq:Euler2_psi}, as follows 
\bseq\label{def:cpsi}
\beq
\bal
 \cpsib(\al) \teq \int_0^{\infty} \f{s^{2+\al}}{ (s^2+1)^{5/2 } } d s,
 \quad \cpsi(\al) \teq \f{\al}{3 \cpsib(\al)}, 
\quad \cpsia(\al) =  8 \cpsi(\al),
  \quad \barcp \teq \cpsia( \tf13 ) .
\eal
\eeq
For any $\al \in (0, 1/3]$ and $s >0$, using $|s^{\al} - s^{ 1/3 }| \les | (\al - \f13) \log s |  ( s^{\al} + s^{ 1/3 })$, 
it is easy to obtain 
\beq
  \cpsia(\al) \asymp 1, \quad   | \cpsia(\al) -  \barcp | \les |\al - \tf13| .
\eeq
\eseq
We suppress the \(\al\)-dependence when no confusion can arise.

Below, we consider $\psi= \BS(\ww)$ \eqref{eq:Euler2_psi:BS}, Using the fundamental solution $\f{1}{4\pi |\xx|}$ in $\R^3$, \eqref{eq:Euler_b}, we get 
\beq\label{eq:psi_reg:3D}
\bal 
 r \psi(r, z)  & = \psi^{\th} \cos \th |_{\th = 0}
 = \cpsia (-\D_{3D})^{-1} (  \om^{\th} \cos \th)  |_{\th = 0} 
 \\
 & = \f{\cpsia }{4\pi}
   \lim_{L \to \infty}  \int_{\R_+ \times \R} \int_0^{2 \pi} 
  \one_{  |(\tr, \tz)| < L }
  \f{ \td r^{1+\al} \om( \td r ,\td z) \cos \th }{  (r^2 + \td r^2 - 2 r \td r \cos \th + |z- \td z|^2 )^{1/2} }
  d \th d \td r d \td z , 
\eal 
\eeq
for $\cpsia > 0$ chosen in \eqref{def:cpsi}. 
For $\om(\tr, \tz)$ odd in $z$ with suitable decay, the integrand is integrable after symmetrization in $\tr , \tz$, and the above formula is well-defined. 

\vspace{0.1in}
\paragraph{\bf Fundamental solution in $\R^5$}

The Poisson equation \eqref{eq:Euler2_psi} can be viewed as a Laplacian equation in $\R^5$ 
for axisymmetric functions \cite{luo2013potentially-2}. In fact, 
for $\psi =\BS(\ww)$ \eqref{eq:Euler2_psi}, since $\cos \th d \th = d \sin \th$ 
and $\f{\cpsia}{4 \pi} = \f{2 \cpsi}{ \pi }$ by \eqref{def:cpsi}, applying integration by parts in $\th$ to \eqref{eq:psi_reg:3D}, we obtain
\beq\label{eq:psi_reg:5D}
\bal
 r \psi(r, z) 
 & =  \f{ 2 \cpsi}{\pi}  \lim_{L \to \infty}  \int_{\R_+ \times \R} \int_0^{\pi} 
  \one_{  |(\tr, \tz)| < L } 
  \f{   r  \td r^{2  + \al } \om(\tr, \tz) \cdot \sin^2 \th  }{ ( r^2 + \td r^2 - 2 r \td r \cos \th + |z- \tz|^2)^{3/2} }  d \tr d \tz d \th,
 \eal
\eeq
where the boundary term at $\th = 0, \pi$ vanishes since $\sin \th = 0$.

For $\yy \in \R^5$, we denote $\yy_h = (y_1, y_2,y_3, y_4)$. We abuse notation to write
\[
  \psi(\yy) = \psi( |\yy_h|, y_5 ),
  \quad   \om(\yy) = \om(|\yy_h|, y_5), \quad \tr = |\yy_h|, \quad  \tz = y_5.
\]

Then the functions $\ww$ and $\psi$ are defined in $\R^5$ with axisymmetry with respect to the axis $\yy_h = 0$.

We parametrize $\yy \in \R^5$ using spherical coordinates for $\yy_h$:
\[
\bal
  \yy &= ( \tr \cos \th , \, \tr \sin \th \cos \th_2, \,
\tr \sin \th \sin \th_2 \cos \th_3,  \,
\tr \sin \th \sin \th_2 \sin \th_3, \,  \tz  ) ,  \\
d \yy &= \tr^3 \sin^2 \th \sin \th_2 d \tr d \th d \th_2 d \th_3 d \tz ,
\qquad \tr \geq 0 , \quad \th, \th_2 \in [0, \pi], \  \th_3 \in [0, 2\pi] .
\eal
\]

For any function $f(\yy) = f(|\yy_h|, \tz) $ and domain $\Ups$ only depending on $\tr = |\yy_h|, \tz = y_5$, 
(thus independent of $\th_2, \th_3$),
using the spherical coordinates, 
$\int_0^{\pi} \int_0^{2 \pi} \sin \th_2 d \th_2 d \th_3 = 4 \pi$, 
and integrating over $\th_2, \th_3$, we obtain
\beq\label{eq:5D_int_id}
\bal
 \int_{\R^5} \one_{(\tr, \tz) \in \Ups}  |\xx-\yy|^k f( \yy ) d \yy
  &= 4 \pi \int_{(\tr, \tz) \in \Ups, \th \times [0,\pi]} A( r, z-\tz , \tz, \th)^k  \sin^2 \th \cdot \tr^3 f(\tr, \tz) d \tr d \tz d \th , \\
        A(r , s,  \tr, \th) & \teq ( r^2 + \td r^2 - 2 r \td r \cos \th + s^2 )^{1/2} . \\
\eal
\eeq

Dividing \eqref{eq:psi_reg:5D} by $r$ and using the above identity, we derive 
\bseq\label{eq:5D_BS}
\beq\label{eq:5D_BS:a}
\bal
  \psi( |\yy_h|, y_5)  =   \int K_{\R^5}( \yy - \td \yy ) |\td \yy_h|^{\al-1} \ww( \td \yy) d \td \yy =  \KRF \ast ( \tr^{\al-1} \ww(\tr, \tz) ) , 
  \quad  \KRF( s ) = \f{ \cpsi}{ 2\pi^2  } |s|^{-3}.
\eal
\eeq
Using the axisymmetry in $\R^5$, we rewrite the Biot-Savart law $\BS$ \eqref{eq:Euler2_psi} equivalently as 
\beq
\psi = \BS(\ww) \iff
   - \D_{\R^5} \psi =\cpsia \om r^{\al - 1} ,
\eeq
with solution formula in \eqref{eq:5D_BS:a}.

\eseq

Dividing $r$ in \eqref{eq:psi_reg:5D}, taking $r =0$, and then using
$\int_0^{\pi} \sin^2 \th d \th = \pi / 2$, we obtain
\beq\label{eq:psi_cpsi}
  \psi(\om)(0, z)     =  \cpsi  \int_{ \R_+^2 \times \R}    \f{ \td r^{2 + \al} \om(\tr, \td z)}{ 
 ( \td r^2 + (z- \td z)^2 )^{3/2}  } d \td r d \tz .
 \eeq

Denote $\xx=  (r, z)$. For $\om$ odd in $\tz$, we define the main term for the stream function
  \bseq\label{def:JJ_3D}
  \beq
    \JJ(\ww)(\xx) =   \f{3 \cpsi }{ 2 \al } \int_{  |\tz| \leq |\xx|, \tr \geq 0 } \f{ \tr^{2 + \al} \tz }{ (\tr^2 + \tz^2)^{5/2} } \ww(\tr, \tz) d \tr d \tz .
\eeq  

Using \eqref{eq:5D_int_id} with $k=0, f =\f{ \tr^{\al-1} \tz }{ (\tr^2 + \tz^2)^{5/2} } \ww(\tr, \tz)  $ and using $\int_0^{\pi} \sin^2 \th = \pi/2$, we rewrite $\JJ$ as 
\beq\label{def:JJ_3D:b}
\bal
  \JJ(\ww)(\xx) &= \f{3 \cpsi }{\pi \al} \int_{  |\tz| \leq |\xx| ,\tr \geq 0} \int_0^{\pi/2} \f{\tr^{2 + \al} \tz}{ |(\tr, \tz)|^5 } 
  \om(\tr, \tz) \sin^2 \th d \tr d \tz  \\
   &=  \f{3 \cpsi}{4 \pi^2  \al } \int_{ |\tz| \leq |\xx| } \f{ \tz }{|\yy|^5 } \ww(\yy)\tr^{\al-1} d \yy 
= \f{1}{2 \al} \int_{|\tz| \leq |\xx|} (\pa_5 \KRF)(-\yy) \cdot \ww(\yy)\tr^{\al-1} d \yy  ,
\eal
\eeq
\eseq
where in the last identity, we have used $y_5 = \tz, |\yy| = |(\tr,\tz)|$ and $\f{1}{2\al} (\pa_5 \KRF)(-\yy)=  \f{1}{2 \al} \cdot \f{3 \cpsi}{2 \pi^2} \f{ y_5}{|\yy|^5} $ due to \eqref{eq:5D_BS:a}.
By definition of \eqref{eq:psi_cpsi} and \eqref{def:JJ_3D}, we have 
\beq\label{eq:psi_JJ}
  \pa_z \psi(\om)(0) = 2 \al \JJ(\om)(\infty).
\eeq

\subsection{1D Profiles Equation and Properties }\label{sec:1D_profile}
We choose  $\cpsia$ in \eqref{eq:Euler2_psi} (with formula in \eqref{def:cpsi}) so that for  $\ww(r, z) = \wwd(z)$ odd in $z$, constant in $r$ with an appropriate decay in $z$, the 3D stream function $\psi = \BS(\ww)$ \eqref{eq:Euler2_psi}  coincides  with the 1D stream function on $r = 0$  with \emph{prefactor $1$} as in \cite{chen2026eulerI}:
\beq\label{def:psi_1D}
 \psi( \wwd(\cdot) )(0, z) = \psid(\wwd)(z), \quad \psiad( \wwd )(x) \teq \int_0^{\infty} (|x+y|^{\al} -  |x-y|^{\al} )  \wwd(y) d y .
\eeq
We prove the above identity in Lemma \ref{lem:vel_bc}. In this paper, we use $f_{\mw{1D}}$ with a subscript $\mw{1D}$ to indicate that  $f_{\mw{1D}}$ is defined along $r=0$.
We omit this subscript when no confusion may arise.

Using $U^r(0, z) = 0$, the 3D profile equation \eqref{eq:profi_3D0}, and the above relation \eqref{def:psi_1D}, we derive the 1D profile equation for \eqref{eq:profi_3D0} along $r=0$:
\beq\label{eq:1D_dyn_profile}
  (c_l z + 2  \psiad )  \pa_z \wwd = ( \com - (1-\al) \pa_z \psiad) \wwd .
\eeq

Similarly, one can derive the evolution PDE associated with \eqref{eq:1D_dyn_profile} along $r=0$:
\footnote{
Model \eqref{eq:1D} was essentially derived in \cite[Section 8]{hou2024potential} 
by restricting 3D Euler \eqref{eq:Euler2} on $r=0$ and assuming that $\om(r, z)$ is constant in $r$.
Our main observation is that the nonlocal operator can be significantly simplified as in \eqref{def:psi_1D}. 
We include \eqref{eq:1D} only to provide background and context; it is
\emph{not used in any estimates} in this work.
}
\beq\label{eq:1D}
\pa_t \wwd +  2  \psiad  \pa_z \wwd =  - (1-\al) \pa_z \psiad \wwd .
\eeq

We introduce the self-similar velocity $V$ and the modified stream function $\psioad$
\bseq\label{eq:1D_normal}
\beq
 V = \f12  c_l x  + \psiad  , 
 \quad  \psioad \teq \psiad - \pa_x \psiad(0) x 
  =    \int_0^{\infty} (|x+y|^{\al} -  |x-y|^{\al}  - 2 \al x y^{\al-1} ) \wwd(y) d y .
\eeq

We impose normalization conditions on $V$ at $x=0$ and choose the scaling parameters 
$c_l, c_{\om}$ as:
\beq\label{eq:1D_normal:c}
V_x(0) = 1,  \quad c_l + 2 \pa_x \psiad (0) = 2 = c_{\om} - (1-\al) \pa_x \psiad(0) .
\eeq

\eseq

Then the equation \eqref{eq:1D_dyn_profile} for the 1D self-similar profile $W$  reduces to
\beq\label{eq:1D_dyn}
\bga 
 2 V \pa_x W = (3- \al - (1-\al) V_x) W ,
 \quad \quad V = x + ( V - V_x(0) x) =  x + \psioad(W)(x) .
\ega 
\eeq

We introduce the main term $\cJa$ for the 1D stream function $\psiad$:
 \beq\label{eq:Jw}
 \cJ_{\al}( w )(x) = \int_0^x y^{\al - 1} w(y) d y.
 \eeq

For $w$ with a sufficiently fast decay, from the definition of $\psi$ in \eqref{def:psi_1D}, we obtain
\beq\label{eq:iden_psiz_J}
 \pa_x \psi_{\al, \mw{1D}}(w)(0) 
 = 2 \al \int_0^{\infty} y^{\al-1} w(y) d y = 2 \al \cJ_{\al}(w)(\infty).
\eeq

To simplify notation, we suppress the dependence of $\psi, \psio$ on $\al$ when no confusion arises.
Denote 
\beq\label{def:al_ep}
\bga
\alb \teq  \tf{1}{3}, \quad \e \teq \tf13 - \al .
\ega
\eeq

We introduce the function $\lgp(x)$, which satisfies  $\lgp(x) \gtr 1$ for any $x \geq 0$
\beq\label{def:lgp}
  \lgp(x) =  \log(x + 2)  .
\eeq

\subsubsection{Existence of $\al$-profile}

We call a nontrivial, odd, smooth solution \(W\) of \eqref{eq:1D_dyn} an \(\alpha\)-profile. We recall the following result from \cite[Theorem 5.5]{chen2026eulerI}.

\begin{thm}[\bf $\f13$-profile]\label{thm:reg_alb}

There exists an odd $C^{\infty}$ $\f13$-profile $(\wwwa, \vvva)$ to \eqref{eq:1D_dyn} with $\al = \f13, \vvva = x +  \psio_{\alb}( \wwwa)$. The profile satisfies  
$|\pa_x \wwwa(0) + 1 | \leq 10^{-6}$,
\footnote{
In \cite{chen2026eulerI}, we construct the approximate profile $\bw$ ($W$ with a ``bar") numerically 
and impose $\pa_x \bw(0)=-1$ in the construction. Due to numerical error,
$\pa_x \bw(0)=-1$ is not preserved and we get $| \pa_x \bw(0) + 1| \leq 10^{-6}$.
All the 1D $\al$-profiles satisfy the normalization $\pa_x \wwwa(0)= \pa_x \waa(0) = \pa_x \bw(0)$.
Thus, we get the estimate of $\pa_x \wwwa(0), \pa_x \waa(0)$ in Theorems \ref{thm:reg_alb}, \ref{thm:1D_profile_prop}. Since $\lgp x \gtr 1$ for $x \geq 0$ \eqref{def:lgp}, the last two estimates in  \eqref{eq:wwwa_basic} are trivial for 
$|x| \les 1$.
}
\beq\label{eq:wwwa_basic}
\wwwa \leq 0 ,  \quad | \wwwa|   \asymp \min( |x|, \ang x^{-\alb}), \quad \vvva   \asymp x \lgp x , 
\quad      | \wwwa + 6 x^{-1/3} |  \les x^{-1/3} | \lgp x|^{-1/3} , 
\eeq
for any $x \geq 0$. Moreover, for any $k \geq 0$, we have
\beq\label{eq:bw_smooth}
\bal
  |\pa_x^k \wwwa | & \les_k \ang x^{-k -\alb},  & \quad |\pa_x^{k+2} \psio_{\alb}| & \les  \ang x^{-k-1}  ,
  \quad \forall \, k\geq 0 .
\eal
\eeq
\end{thm}

Some properties of $(\wwwa,\vvva)$ in Theorem~\ref{thm:reg_alb} are not used later;
we include them for completeness.

 Recall the $\cJa$-operator from \eqref{eq:Jw}.
For $\b \in [-\e, 0] $, we introduce the $\cJaa$-function
\beq\label{eq:wa}
     \cJab(x) := \cJa( \ang y^{\b} \wwwa(y) )(x) ,
\quad \cJaa = (1 + \cJab^2)^{1/2}. 
\eeq
Since $\wwwa(x)<0$ for $x>0$ by Theorem \ref{thm:reg_alb}, we obtain that $\cJaa(x)$ is increasing in $x$.

We recall the following estimates: for $\b(\al)$ from
\cite[Theorem 6.1]{chen2026eulerI}, for $\cJaa$ from
\cite[Theorem 6.9]{chen2026eulerI}, and for the $\al$-profile from
\cite[Theorem 7.10]{chen2026eulerI}.

\begin{thm}[\bf Regularity and decay of $\al$-profile]\label{thm:1D_profile_prop}

Let $\kp, \kp_1 < 1$ be the parameters in \eqref{def:kp},
and $(\wwwa, \vvva)$ be the $\f13$-profile in Theorem \ref{thm:reg_alb}. There exists an absolute constant $\beps_4 > 0$ 
\footnote{ We keep the notation \(\beps_4\) the same as in \cite{chen2026eulerI} for consistency.
}
such that for any $\e = \f13 - \al \in (0, \beps_4)$, the following results hold. There exists a $C^{\infty}$ $\al$-profile $(\waa, \vaa)$ solving \eqref{eq:1D_dyn}  and an exponent $\b = \b(\al)$ with 
\beq\label{eq:beta_est}
\b  = \b(\al) =  - \tf{1}{8} \e + O(\e^{4/3}) \in [-\tf{1}{2} \e, \ 0], 
\eeq
such that 
\bseq
\begin{gather}
 \pa_x \waa(0) = \pa_x \wwwa(0)  \in [-1 - 10^{-6}, -1 + 10^{-6}], \quad   \waa(x) < 0,  \quad  \forall \, x > 0,  \label{eq:waa_sign}  \\
   \nlinf{  (\waa - \ang x^{\b(\al) } \wwwa) (   |x|^{-\bbb} + |x|^{\alb - \b - \kp_1 \e } \cJaa^{-\kp}  ) }  \les \e , \quad \bbb = 1.2  \ . \label{eq:waa_error}   
\end{gather}
\eseq
Here,  $\cJaa$ is defined in \eqref{eq:wa} and satisfies %the following estimates
\beq\label{eq:Ja_hat}
\bal
\cJaa(x) \asymp \min( \lgp x, \f{1}{\e} )  , \  \cJaa =   6 \f{1 - \ang x^{-\he}}{\he} +  O( \min( \lgp x, \f{1}{\e} )^{\f23} ) ,
\   \he \teq \e - \b(\al) .
\eal
\eeq

\begin{enumerate}[label=(\roman*), leftmargin=1.5em]

\item \textsl{(Regularity and Asymptotics)}. We have the following estimates for $\waa$
\beq\label{eq:waa_reg}
  |\pa_x^k \waa |   \les_k \ang x^{- k  -\hal },    \quad  \forall  \ k \geq 0,  
 \quad  | \pa_x^{k+1} \vaa |   \les_k \ang x^{-k + \al -\hal} , \quad  \forall \  k \geq 1 ,
\eeq
with implicit constants independent of $\al, \e$, and 
 \beq\label{eq:wa_upper_lower}
 \bga
     |\xx| \we |\xx|^{ -\hau } \les |\waa(\xx)| \les |\xx| \we |\xx|^{-\hal}, \\
           \hau =  \f13 + \f18 \e + c_{\wwwa} \e^{2-\kp}, 
     \quad \hal = \f13 + \f18 \e - c_{\wwwa} \e^{2 - \kp} ,
     \ega
 \eeq 
where $c_{\wwwa}$ is an absolute constant independent of $ \e$.

\item \textsl{(Scaling)}. Recall $\pa_z \psi(f)(0) = 2 \al \cJa(f)(\infty)$ from \eqref{eq:iden_psiz_J}. Denote 
\footnote{
The scaling exponents $(\bar c_{l,\al,\mw{1D}}, \bar c_{\om,\al,\mw{1D}})$
are the same as $(\bar c_{l,\al}, \bar c_{\om,\al})$ in
\cite{chen2026eulerI}. We add the subscript \emph{1D} to distinguish these exponents from
the 3D analogues used in the present work.
}
\bseq\label{eq:wa_decay_1D}
\beq\label{eq:wa_decay_1D:a}
\bar c_{l,\al, \mw{1D}} = 2 - 4\al \cJa( \waa)(\infty), \quad  \bar c_{\om, \al, \mw{1D}}  = 2 + (1 - \al) 2 \al \cJa(\waa)(\infty) .
\eeq
We have 
\beq
\cJa(\waa)(\infty) = - \tf{16}{3 } \e^{-1} + O(\e^{2-\kp}), 
   \quad \bar c_{l,\al, \mw{1D}}  , \ - \bar c_{\om, \al, \mw{1D}} \asymp \e^{-1},
\eeq
and 
\beq\label{eq:wa_decay_1D:c}
      \f{ \bar c_{\om, \al, \mw{ 1D} }}{ \bar c_{l, \al, \mw{ 1D} } }= -   \f{1}{3} - \f18 \e 
      +  O (\e^{2 - \kp}), \quad
 \bar c_{l,\al, \mw{1D}} 
       = \f{64}{9} \e^{-1} + O(\e^{-\kp}) , 
    \quad \bar c_{\om, \al, \mw{1D} } = - \f{64}{27} \e^{-1} +  O(\e^{-\kp}) \, .
\eeq

\eseq

\item \textsl{(Log-improvement)}. We have the following estimates on $\waa$
\begin{align}
   |\f{x \pa_x \waa}{\waa} + \f{1}{3} | & \les \e  + |\lgp x|^{\kp-2}
     , \qquad  |  \f{ \pa_x \waa}{\waa} |  \les  |x|^{-1} . 
\label{eq:1D_prof_est2:c}
\end{align}
and 
\beq\label{eq:1D_prof_dxx} 
   | \pa_x  ( \f{x \pa_x \waa}{\waa} )  |  \les 
   \ang x^{-1} \cJaa^{-1} , \quad 
   | x\pa_x  ( \f{\pa_x \waa}{\waa} )  | \les |x|^{-1}  \, .
\eeq

In particular, there exists $C_{\wwwa}$ only depending on $\wwwa$, such that %for any $ \e \leq 1$, we have
\beq\label{eq:1D_prof:low}
   x  \f{ \pa_x \waa}{\waa} + \f{1-\al}{2} \geq - C_{\wwwa} \f{|\xx|}{\ang \xx } ( \e^{2-\kp} + |\lgp x|^{\kp - 2}  ) \, .
\eeq

\item \textsl{Estimates for $\vaa$}. Let $\he = \e - \b(\al)$.
For any $x \geq 0$, we have the following estimates for $\vaa$ 
\begin{align}
   |\pa_x \vaa - \tf{1}{x} \vaa - 4 \ang x^{-\he} | & \les 
   \ang x^{-\he}  |\lgp x |^{-1/3}  +  \e^{1- \kp} , \label{eq:1D_prof:va} \\
 \vaa(x) & \asymp \cJaa(x) \cdot |x| ,\qquad  \vaa  \geq 0.
   \label{eq:vaa_low}
\end{align}
All the implicit constants are \emph{independent} of $\al, \e$.
\end{enumerate}

\end{thm}

In the rest of the paper, we simplify \(\b(\al)\) as \(\b\) whenever no confusion can arise.

\begin{remark}[\bf $ \waa$ is close to $ \wwwa$]

Estimates in items (i), (iii), (iv) show that $\waa$ is sufficiently close to $\wwwa \ang x^{\b}$ in weighted space,
\footnote{
Using \eqref{eq:waa_error} and \eqref{eq:1D_dyn}, the proof of Theorem \ref{thm:1D_profile_prop} is given in
\cite[Section 7.9]{chen2026eulerI}, and is mostly self-contained.
}
and \(\waa\) has far-field asymptotics close to 
 that of $\wwwa$: $\wwwa = x^{-1/3} ( -6 + \ccb |\log x|^{-1/3} )$.
Formally, we get smallness $\e$ or faster decay $|\lgp x|^{\kp-2}$ in the first bound in \eqref{eq:1D_prof_est2:c} since $\waa$ has a decay rate close to $x^{-1/3}$;
we gain an extra $\cJaa^{-1}$ factor in the first bound in \eqref{eq:1D_prof_dxx}
since $ \f{x \pa_x \waa}{\waa} = C + l.o.t.$ for some  $C$ close to $-\f13$ and the derivative of the lower order term
decays faster than $\f{1}{x}$.
Estimate \eqref{eq:1D_prof:va} comes from perturbing the \emph{exact} identity 
\[ 
( \pa_x \psio_{\alb, \mw{1D}} - \tf{1}{x} \psio_{\alb, \mw{1D}})( f)(x)
\equiv 4  , \quad  f(x) = -6 \cdot \mw{sign}(x) \cdot |x|^{-1/3}  \, , 
\] 
where $ \pa_x \psio_{\alb, \mw{1D}} - \tf{1}{x} \psio_{\alb, \mw{1D}}$ is defined via the 
associated integral formula \eqref{def:psi_1D}. See \cite[Section 2]{chen2026eulerI}.
 The function $\cJaa$ behaves like $\lgp x$ 
 and we gain log-improved estimates in \eqref{eq:1D_prof_est2:c} and \eqref{eq:1D_prof_dxx}, 
 which are crucial for overcoming logarithmic growth in later fixed point estimates. 
\end{remark}

\subsubsection{Contraction estimates in 1D and norms}
We recall the weights $\vpa, \vpcc, \vpk$ and norms:
\footnote{
We combine the power weights $\ang x^{-\kp_1 \e}$ and $\ang x^{c_{\wwwa}\e^{2-\kp}}$ 
used in \cite[Section 7.10]{chen2026eulerI} into the weight $\ang x^{-\epa}$ to simplify notation. 
To reduce the number of parameters, we write \(\e_1\) explicitly as \(\kp_1\e\), rather than introducing it as in \cite{chen2026eulerI}.
}
\bseq\label{norm:Xc}
\beq\label{def:vpc}
\bal
  \vpa(x) & \teq \max_{i} ( \mu_{\sst, i}^{-1} |x|^{\bb_{\sst, i}} ),  \ 
\vpcc(x)   \teq |\waa|^{-1} \ang x^{- \epa } \cJaa^{- \kp} \vpk,  \ 
     \vpk(x) \teq \exp \big( 9 \kp_1^{-1}  \cdot \ang x^{- \kp_1 \e}  \big),     
\eal
\eeq
and the $\cXc$-norms:
\beq
\bal
    \| w \|_{\cXc} = \| w  \cdot \max( \muc \vpc , \vpa ) \|_{ \linff} , 
\eal
\eeq
where $\vmu_\sst, \bb_\sst, \epa$ are parameters given by 
\footnote{
The parameters $\vmu_\sst, \bb_\sst$ except for $\bbb$ are not used in this paper. 
We present them for completeness.
}
\beq
   \epa =  \kp_1 \e  - c_{\wwwa} \e^{2-\kp} , \quad 
   \vmu_{\sst} =  ( 1, 4, 30 ) , 
    \quad 
   \bb_\sst = (-1.2, -0.5, 0),  \quad \bbb \teq -\bb_{\sst, 1} = 1.2 ,
\eeq
\eseq
 $\kp_1, \kp$ are defined in \eqref{def:kp},
 and  $c_{\wwwa}$ is an absolute constant. 
 Recall the parameters $\b, \hal,\hau$ from Theorem \ref{thm:1D_profile_prop} and \eqref{eq:wa_upper_lower} and $\e= \f13 - \al$.
 By requiring $\e$ small enough, we impose 
 \bseq\label{ran:ep_all}
\beq\label{ran:ep}
 \epa = \kp_1 \e  - c_{\wwwa} \e^{2-\kp} , \quad \epa -  \e^2 \in [ \tf12 \kp_1 \e, \, 2 \kp_1 \e] \subset [0, \tf{1}{500} 
 \e ] , \quad \e \in [0, 10^{-4}] , \quad \e^2 < \tf14 \epa ,
\eeq
and 
\beq\label{ran:ep2}
 \he = \e -  \b     ,  \,  \hal - \al   , \,  \hau - \al  \in [ \f{9}{8} \e - C \e^{ 2 -\kp }, \, \f{9}{8} \e + C \e^{ 2 -\kp } ]
     \subset [ \f{8.9}{8} \e,  \ \f{9.1}{8} \e ] ,
     \quad \hal - \f13 \in [ \f{0.9}{8} \e,  \ \f{1.1}{8} \e ] . 
\eeq
\eseq

Using \eqref{eq:wa_upper_lower} and the definition of the weight, we obtain 
\beq\label{eq:wg_equiv}
\bga
 \vpa(x)  \asymp   |\xx|^{-\bbb} + 1 , \quad \vpk(x) \asymp 1 , 
 \quad \max(\vpa, \vpcc) |\waa| 
 \asymp |x|^{ -\bbb+1} +  |x|^{ -\epa}   \cJaa^{-\kp} .
  \ega
\eeq

We recall the following key contraction estimate for the $\al$-profile from \cite[Theorem 7.12]{chen2026eulerI}.

\begin{thm}[\bf Contraction estimates]\label{thm:contra_lin_exact}

Let \(\e=\f13-\al\) and let $ \beps_4$ be as in Theorem \ref{thm:1D_profile_prop},
and let $\muc$ be the parameter in the $\cXc$-norm \eqref{norm:Xc}. 
Define the linear operator  $\cLa(w)$ 
associated with the $\al$-profile $(\waa, \vaa)$ in  Theorem \ref{thm:1D_profile_prop}:
\beq\label{eq:lin_1D_alpha}
  \cLa (w) \teq \waa  \int_0^{x}  \f{  \cRa(w)}{2 \vaa }  ,  
  \quad  \cRa( w )   = - (1-\al) \psiox (w) - 2 \psio  \f{\pa_x \waa}{\waa}  .
\eeq
 There exists $\beps_5 \in (0, \beps_4)$ such that for any $\e \in (0, \beps_5 )$ and any $x\geq 0$,  we have
\[
 |\max(  \vpa, \muc \vpcc ) \cLa(w)| \leq  
  \max(  \vpa, \muc \vpcc ) \B| \waa  \int_0^{x}  \f{  \cRa(w)  }{2 \vaa }  \B|
\leq  \min( \lamcL  , C \e^{-1} \ang x^{-\kp_1 \e /2} )
\| w \|_{\cXc} ,
\]
where $\lamcL= \tf34 + \tf14 \max( \f{\cff}{2 \kp} ,  0.95 ) \in (0, 1) $,
and $\cff, \kp, \kp_1$ are the constants defined in \eqref{def:kp}.
\end{thm}

\subsection{Approximate profile and fixed point map}\label{sec:appr_profi}
Let $\waa$ be the 1D profile constructed in Theorem \ref{thm:1D_profile_prop}. We construct 
the scaling parameters using \eqref{eq:1D_normal} and \eqref{eq:wa_decay_1D:a}
\beq\label{eq:normal_bar}
  \bcl \teq  \bar c_{l,\al, \mw{1D}} 
   = 2 - 4 \al \cJa(\waa)(\infty),  \quad \bcw \teq  \bar c_{\om,\al, \mw{1D}}  =  2 + 2 \al (1-\al) \cJa(\waa)(\infty).
\eeq

We extend $\waa$ constantly in $r$ as a profile $ \wwb$ in $\R^3$ :
\beq\label{def:3d_Om}
   \bar \Om_{\al}(r, z) = \waa(z) ,\quad \forall z \in \R. 
\eeq

Since we fix $\alpha$ in most derivations, we simplify $\bar{\Omega}_\alpha$ as $\wwb$.
 Moreover, we use ``bar" notation $\bar \cdot$ to denote the variable generated from $\wwb$. We decompose the solution $(\Psi, \www, \UU, c_l ,c_{\ww})$ to \eqref{eq:profi_3D0} into the profile 
$(\bar \Psi, \wwb, \bar \UU, \bar c_l ,\bar c_{\ww})$ and perturbation $(\psi, \ww, \uu, \tcl, \tcw)$
\beq\label{eq:pertb_decom}
  \Psi = \bar \Psi + \psi, \quad \www = \wwb + \ww ,
  \quad \UU = \bar \UU +  \uu, 
  \quad 
    c_l = \bcl + \tcl, \quad c_{\om} = \bcw + \tcw,
\eeq

To determine $c_l, c_{\om}$ in  \eqref{eq:profi_3D0}, we impose the 
normalization conditions as in the 1D case \eqref{eq:1D_normal:c}:
\bseq\label{eq:normal_cond}
\beq
 c_l + 2 \Psi_z(0) = c_{\om} - (1 - \al) \Psi_z(0) \equiv 2.
\eeq

These conditions and \eqref{eq:normal_bar} imply the following conditions for the perturbation 
\beq
  \tcw = (1 - \al)  \psi_z(0),  \quad 
\tcl = - 2   \psi_z(0). 
\eeq

\eseq

We define the stream function and velocity associated to $\wwb, \www, \ww$ via \eqref{eq:Euler2_psi} and \eqref{eq:Euler2:c}
\beq\label{def:3d_vel}
\bal
\bar \Psi &= \BS( \wwb ), \quad & \bar U^r = - r \pa_z \bar \Psi, \quad & \bar U^z = 2 \bar \Psi + r \pa_r \bar \Psi , \\
\psi & = \BS(\ww), \quad  & u^r = - r \pa_z \psi,   \quad & u^z  = 2 \psi + r \pa_r \psi , \\
 \Psi & = \BS(\www), \quad & U^r = - r \pa_z \Psi,  \quad & U^z = 2 \Psi  + r \pa_r \Psi . 
 \quad
\eal
\eeq

We introduce the stream function with correction near $0$
\beq\label{def:psio}
  \psio(r, z) \teq \psi(r, z) - \pa_z \psi(0, 0) \cdot z .
\eeq
Since $\psi$ is odd in $z$, we have $\pa_r \psi(0) = 0$.
By definition, we have $\pa_r\psio = \pa_r \psi$ and $\pa_z^2 \psio = \pa_z^2 \psi$.

We denote the coefficients of the transport term in the self-similar variable as
\bseq\label{def:Q}
\beq
\bal
  \QQb  & \teq ( \bar Q^r, \bar Q^z) = (\bar c_l r + \bar U^r, \bar c_l z + \bar U^z) 
  = ( \bar c_l r  - r \pa_z \bar \Psi, \bar c_l z +  2 \bar \Psi + r \pa_r \bar \Psi ) \, , \\
   \td \QQ(\ww) & \teq (\td Q^r, \td Q^z)= ( \tcl r + u^r, \tcl z + u^z )  
   = ( \tcl r - r \pa_z \psi , \tcl z + 2 \psi + r \pa_r \psi  )  \,  , \\
  \QQ & \teq (Q^r, Q^z) 
  =  (c_l r + U^r, c_l z + U^z) 
= (c_l r - r \pa_z \Psi, c_l z + 2 \Psi  + r \pa_r \Psi )  \, .
\eal 
\eeq
By definition, we have 
\beq
 \QQ = \QQb + \td \QQ(\ww) \, .
\eeq
\eseq

Using the notation \eqref{def:psio} and the normalization \eqref{eq:normal_cond}, we obtain
\beq\label{eq:psio_res1}
 \td Q^z = \tcl z + 2 \psi + r \pa_r \psi =  2 \psio + r \pa_r \psi , 
   \quad \tcw - (1-\al) \psi_z =  - (1-\al) \psio .
\eeq

Using these notations, we  rewrite the profile equation \eqref{eq:profi_3D0} as 
\beq\label{eq:profi_3D}
  \QQ \cdot \na \www  = (c_{\om} - (1-\al) \Psi_z) \www \ . 
\eeq

Since $\wwb$ only depends on $z$, linearizing the above equation around the approximate profile $(\wwb, \bcl, \bcw)$, we obtain 
\bseq\label{eq:lin_2D}
\beq\label{eq:lin_2Da}
\bal
  \QQ \cdot \na \ww  =  (c_{\om} - (1-\al) \Psi_z) (\ww + \wwb) - Q^z \pa_z \wwb \, .
  \eal 
\eeq

To simplify the notation, we introduce 
\beq\label{eq:lin_2Db}
\bal
 \bcB & \teq \bar c_{\om} - (1-\al) \bar \Psi_z,  
\quad \tcB( \ww ) = \td c_{\om} -(1-\al) \psi_z ,  \\
   \bar \cR  & \teq \bcB -  \bar Q^z \f{\pa_z \wwb}{\wwb} =  \bar c_{\om} - (1 - \al) \bar \Psi_z -  (\bar c_l z +  2 \bar \Psi + r \pa_r \bar \Psi  ) 
\f{ \pa_z \wwb }{\wwb} , \\
\cL_{\psi}(\ww)  & \teq \tcB(\ww) - \td Q^z  \f{ \pa_z \wwb }{\wwb}  = \tcw - (1 -\al) \psi_z - (\tcl z + 2 \psi + r \pa_r \psi) \f{ \pa_z \wwb }{\wwb} .
\eal
\eeq

The term $\bar \cR$ measures the residue error of the profile $\wwb$ in the profile equation
\eqref{eq:profi_3D0}.  Using the $\psio$ notation \eqref{def:psio}, \eqref{eq:psio_res1}, and the normalization condition \eqref{eq:normal_cond}, we  rewrite $\Lpsi, \tcB$ as 
\beq
 \Lpsi = - (1 - \al) \psio_z  -  (2 \psio + r \pa_r \psi)  \f{ \pa_z \wwb }{\wwb} ,
 \quad \tcB(\ww) = - (1- \al) \psio_z.
 \eeq
Clearly, $\Lpsi, \tcB$ are linear operators of $\ww$. Then we can rewrite \eqref{eq:lin_2Da} as 
\beq\label{eq:lin_2Daa}
  \QQ \cdot \na \ww = ( \bcB+ \tcB(\ww ))  \ww 
  +  (\bar \cR + \Lpsi) \wwb,
  \quad \QQ = \bar \QQ + \td \QQ(\ww).
\eeq

\eseq

\subsubsection{Fixed point map}

We reformulate the profile equation \eqref{eq:lin_2D} as a fixed point problem. 
Given $\ww$, we consider the following linear equations in $\eta$ 
\bseq\label{eq:fix_pt_2D}
\beq\label{eq:fix_pt_2D:a}
\bal 
\QQ \cdot \na \eta = ( \bcB + \tcB(\ww)) \eta + (\bar \cR + \Lpsi(\ww)) \wwb ,
\quad 
 \QQ =  \bar \QQ + \td \QQ(\ww) .
\eal
\eeq

From \eqref{eq:lin_2Db}, we obtain $\bcB + \tcB(\ww) - Q^z \f{\pa_z \wwb}{\wwb} = \bar \cR + \Lpsi $.
Dividing both sides in \eqref{eq:fix_pt_2D:a} by $\wwb$, we rewrite 
\eqref{eq:fix_pt_2D} equivalently as 
\beq\label{eq:fix_pt_2DG}
     \QQ \cdot \na  G   = ( \bar \cR + \Lpsi) (G + 1) ,
 \quad G \teq \f{\eta}{\wwb} ,
 \quad  \QQ =  \bar \QQ + \td \QQ(\ww) .
\eeq

In Proposition \ref{prop:reg_solu}, 
for $\ww$ small enough in a weighted $W^{1,\infty}$ space $\wonea$ defined in \eqref{def:3d_norm}, we  construct a unique $C^{1,\al}$ solution $\eta$ to \eqref{eq:fix_pt_2D}. We encode the relation from $\ww$ to $\eta$  as a map 
\beq\label{eq:fix_pt_2Deta}
  \eta = \FFF( \ww ) .
\eeq
\eseq

Then we reformulate solving \eqref{eq:lin_2D} as a fixed point problem $   \ww = \FFF(\ww) $.

\subsection{Anisotropic weights and norms}

Below, we introduce the anisotropic weights and norms, and will motivate them in 
Section \ref{sec:idea}. Denote 
\bseq\label{def:3d_wg}
\beq\label{def:rag}
 \xx \teq (r, z), 
  \quad  \rag = \big( \f{\ang z}{ r + \ang z  } \big)^{\kag}, 
  \quad \kag = \f{1}{1000} .
\eeq
where $ \la a \ra  =( |a|^2 + 1)^{1/2}$ denotes the Japanese bracket, and $\textbf{ag}$ is short for angular.

Recall $\epa, \bbb$ from the 1D $\cXc$-norm \eqref{norm:Xc} 
and $\kp$ from \eqref{def:kp}. We design 
a radial weight $\rhoo$ as
\beq\label{def:3d_wg:a}
\rhoo(\xx) = (|\xx|^{-\bbb + 1} + 1)  \ang \xx^{ - \epa  } |\cJaa(\xx)|^{-\kp},
\quad \cJaa(\xx) \teq \cJaa(|\xx|).
\eeq
Note that $\rhoo$ depends only on the radial variable $R$, while $\rag$ depends primarily on the angular variable $\xi$, in the polar coordinates $(R,\xi)$ for $(r,z)$.
Given $\rhoo, \rag$, we design the weights 
\beq\label{def:3d_wg:rz}
\bal
  \rhor(\xx)  &= \rhoo(\xx)  \rag(\xx) \cdot \la \xx \ra ,
     \quad \rhoz(\xx) =   \rhoo(\xx) \rag(\xx) \cdot |z| .
\eal
\eeq
for weighted $W^{1,\infty}$ estimates in Section \ref{sec:3D_EE}. 
We design $\rhoc$ asymptotically stronger than $1$ for $|\xx| \gg 1$
\beq\label{def:rhoc}
\rhoc \teq \ang x^{\e^2}.
\eeq
\eseq

 The exponent $\e^2$ can be replaced by any other \emph{sufficiently small, positive} number, and it is treated perturbatively. See \hyr[sec:idea_comp]{\its Compactness of $\FFF$} for motivation.

Let $\hal$ be the exponent in \eqref{eq:wa_upper_lower}.  Given $\rhoo, \rag$, we define various norms associated with $\rhoo, \rag$
\begin{align}\label{def:3d_norm}
  \| f \|_{\linfza } & \teq \nlinf{ ( |z|^{-1} + \ang z^{\hal})  (\rhoo f)(0, z) } ,  
 &   \quad 
    \| f \|_{\linfa } & \teq \nlinf{ 
     \rhoo \rag  \cdot  \max( |\xx|^{-1} , 1) \ang z^{\hal} f  } ,  \notag \\ 
    \| f \|_{\wwra} & \teq  \|  \rhoo \rag \cdot \ang z^{\hal}   \la \xx \ra \cdot \pa_r f (r, z)  \|_{L^{\infty}},  
    & 
       \nnrr{f}  & \teq  \max( \nnla{f} , \, \nnra{f}  ) , \\
  \| f \|_{\wwza} &\teq    || \rhoo    \rag  \cdot \ang z^{\hal + 1} \cdot \pa_z f (r, z) ||_{L^{\infty}}  , %\\
  & 
\| f \|_{\wonea} & \teq \max( \nnla{ f} ,  \  \nnr{ f } ,   \ \nnz{ f}^{1-\al} \nnrr{f}^{\al}  ) .
\notag
\end{align}

We use $|z|$ in the weight $\rhoz$ \eqref{def:3d_wg:rz} for energy estimates in Section \ref{sec:3D_EE}, and the \emph{Japanese bracket} $\ang z$ in $\wwza$ \eqref{def:3d_norm}  for nonlocal estimates in Section \ref{sec:nonlocal}.
These two weights are for \emph{different purposes}.

From the definition of $\cXc$-norm \eqref{norm:Xc}  and $\rhoo$ in \eqref{def:3d_wg}, 
using \eqref{eq:wg_equiv} and \eqref{eq:wa_upper_lower}, we obtain $  \rhoo (|z|^{-1} + \ang z^{\hal}) \les  \rhoo |\waa|^{-1} \asymp  \max(\vpa, \vpc) $ and 
\beq\label{norm:comp_bd}
  \nlinf{ \rhoo |\waa|^{-1} f(0,z)} \asymp \| f \|_{\cXc},
  \quad 
 \nna{f} =  \|  (|z|^{-1} + \ang z^{\hal} ) (\rhoo f)(0, z)  \|_{L^{\infty}} 
 \les  \| f \|_{\cX_3} .
\eeq

\begin{remark}[\bf Radial weight $\rhoo$ ]\label{rem:3d_wg}
We choose the radial weight $\rhoo$ in \eqref{def:3d_wg} so that it is  comparable to the weight 
 in the 1D $\cXc$-norm \eqref{norm:Xc}. See estimates \eqref{norm:comp_bd}. Moreover, 
from \eqref{eq:wa_upper_lower},   \( |\waa|^{-1} \)  
and the weight \( |z|^{-1} + \ang z^{\hal} \) in the $\linfza{\cdot}$-norm \eqref{def:3d_wg} are essentially comparable. 
\end{remark}

In Sections \ref{sec:3D_boundary}-\ref{sec:EE_C1r}, we perform weighted estimates for $\eta$ along $r=0$, $\pa_r \eta$, and  $\pa_z \eta$, respectively, to control 
energies stronger than
$\linfza, \wwra, \wwza$-norms.
 We use the ``dot" notation to distinguish the homogeneous $\nnr{\cdot}$ 
and non-homogeneous $\nnrr{\cdot}$ weighted norms for $ f$.

If  $f$ satisfies $\nnr{f}, \nnz{f} < \infty$, we obtain that $\pa_r f$ decays much faster than $\pa_z f$: 
\[
 |\rhoo \rag \pa_r f | \les \ang \xx^{-1} \ang z^{-\hal},
 \quad  |\rhoo \rag \pa_z f | \les \ang z^{-\hal-1} \, .
\]

A typical function with the above properties 
is the approximate profile $\wwb(z)$. We have $\pa_r \wwb = 0$, while $|\pa_z \wwb|= |\pa_z \waa| \les \ang z^{-1-\hal}$ by \eqref{def:3d_Om}, \eqref{eq:waa_reg} \emph{only} decays in $z$. 
We expect that the perturbation around the profile 
has similar decay properties, and use anisotropic weighted derivatives  $ \ang \xx \pa_r f, \ang z \pa_z f$ in the norms \eqref{def:3d_norm}  to capture the specific decay structures.

In later estimates of the fixed point equation \eqref{eq:fix_pt_2D}, we do not estimate $\nnla{\eta}$ directly.
Instead, we estimate $\nna{ \eta}$ along $r=0$ building on Theorem \ref{thm:contra_lin_exact}, and 
estimate $\nnr{\eta}$ for the bulk part $\pa_r \eta$ in the whole space.
We control global $\linfa$-norm using $\linfza$, $\wwra$-norms and Lemma \ref{lem:linf_R3}.

\subsection{Outline of Proofs}\label{sec:idea}

In this section, we outline the proof of Theorems \ref{thm:self-similar}, \ref{thm:main_blowup}.

For completeness, we briefly review the construction of 
1D profiles in \hyr[sec:idea_step1]{\its Step 1},\, \hyr[sec:idea_step2]{\its 2} established in \cite{chen2026eulerI}, and summarize the main results in 
\cite{chen2026eulerI} in Theorems \ref{thm:reg_alb}, \ref{thm:1D_profile_prop}, \ref{thm:contra_lin_exact}. 
We then outline the constructions of 3D profiles in \hyr[sec:idea_step3]{\its Step 3}. We discuss the ideas of proof for Theorem \ref{thm:main_blowup} in Section \ref{sec:idea_main_blowup}. 

\vs{-0.05in}

\subsubsection*{\bf Step 1. Construct $\f13$-profile in 1D}\label{sec:idea_step1} 
First, we construct a profile $(\wwwa, \mw{\psi}_{\alb})$ to \eqref{eq:1D_dyn} with $\al = \alb$. 
The profile has an \emph{$\infty$} blowup rate $c_l = \infty, \com = -\infty$, and a slow decay rate $\wwwa \sim |x|^{-1/3}$. Hence, 
the integrand in \eqref{def:psi_1D} only has a non-integrable decay rate $y^{-1}$ and its stream function is not well-defined by the formula \eqref{def:psi_1D}.  
To capture these unbounded quantities, we impose the normalization in \eqref{eq:1D_normal:c} 
\emph{a priori} so that the $+\infty$ from $c_l$ and  $-\infty$ 
from $\pa_x \psiad (0)$ are \emph{cancelled}, leading to a normalized value $2$. 
The same cancellation applies to $\com, \pa_x \psiad (0)$. 
To overcome the slow decay of the profile, we introduce  the modified stream function in \eqref{eq:1D_normal}, which is well-defined for $\wwd$ with a slow decay rate due to the improved decay rate of the kernel in $y$ \eqref{eq:1D_normal}. 

After imposing normalization \eqref{eq:1D_normal:c} and introducing $\psioad$,  a key difficulty remains:  the self-similar velocity $ V $ has a 
\emph{superlinear} growth $ V \asymp x \lgp x$,
\footnote{
This stands in sharp contrast to existing works on (nearly) self-similar blowup in incompressible fluids \cite{elgindi2019finite,elgindi2019stability},\cite{chen2019finite2,ChenHou2023a,ChenHou2023b}, where the velocity has sublinear growth.
}
and $\pa_x V$ is not uniformly bounded. 
To overcome this difficulty, we derive the leading order asymptotics \emph{a priori} and construct an approximate steady state $\bw$  to the dynamic equation of \eqref{eq:1D_dyn}  with $\bw =  |x|^{- 1/3 }  ( -6 + \ccb \log |x|^{-1/3 }  + l.o.t. ), \ccb \approx 1.5745$  for $x \gg 1$. 
To obtain the exact profile, we linearize \eqref{eq:1D_dyn} around $\bw$ and 
reformulate the equation of perturbation $\td \ww = \wwd - \bw$, as a fixed-point problem along the characteristic. 
By exploiting a crucial \emph{local-nonlocal} cancellation along the characteristic, we construct 
an exact profile $\wwwa$ with asymptotics 
$\wwwa = |x|^{-\f13}( -6 + O( |\log x|^{-1/3} ))$ for large $x$, which is crucial for later constructions around $\wwwa$. We further establish a \emph{near-field contraction estimate} for a linear operator in the fixed point map 
around the exact profile $\wwwa$.
Computer assistance is used to construct the approximate profile $\bw$
and to verify some inequalities for  the fixed-point argument.

\vs{0.05in}
The following steps are purely analytic (pen-and-paper). We construct 1D $\al$-profiles and  3D $C^{1,\al}$-profiles with $\al = \f13 - \e$ 
around the $\f13$-profile $\wwwa$ by exploiting the smallness of $\e \ll 1$.

\vs{-0.05in}
\subsubsection*{\bf Step 2. Construct $\al$ -profile in 1D}\label{sec:idea_step2} 

A key difficulty is that the $\f13$-profile $\wwwa$ cannot be used directly as an approximate profile for \eqref{eq:1D_dyn} with $ \al < \f13$, since it has an \emph{$O(1)$ relative residual error} for large $x$ and it does not capture the correct decay 
rate of $(\f13-\e)$-profile, which \emph{determines} the \emph{finite} blowup rate $(c_l, c_{\om})$ of the profile. A key step is to construct an approximate profile of the form
$\wa = \ang x^{\b} \wwwa$ for \eqref{eq:1D_dyn}.
By carefully expanding the residual error for $\xx \gg 1$, we derive the main term, which determines $\b = - \f{\e}{8} + l.o.t.$, and obtain a small error with improved decay rate.

We then perform a fixed point argument around $\wa$ 
and perturb the  \emph{near-field contraction estimate} of $\cL_{\wwwa}$ in 
\hyr[sec:idea_step1]{\its Step 1} to obtain a similar contraction estimate for a similar operator $\cL_{ \wa}$.
Combining this property and a far-field contraction estimate for  $\cL_{ \wa}$, which requires $\kp < 1$ in \eqref{def:kp}, we construct the exact $\al$ -profile $\waa$. We further perturb the contraction estimate for $\cL_{ \wa}$ to obtain that for $\cL_{\waa}$. See Theorem \ref{thm:contra_lin_exact}.  For \(\al<\alb\), the \(\e\)-improved decay of the integrand in \eqref{def:psi_1D}, namely \(O(y^{-1-C\e})\), allows us to define the 1D stream function \(\psiad(\waa)\) through \eqref{def:psi_1D}, with 
uniformly bounded \(\pa_x\psiad(\waa), \tf{1}{x} \psiad(\waa)  \), in sharp contrast to the logarithmic growth at \(\al=\tfrac13\).  Then using normalization \eqref{eq:1D_normal:c}, we determine finite blowup rate 
$\bar c_{\om, \al, \mw{ 1D} }, \bar c_{l, \al, \mw{ 1D} }$ 
in \eqref{eq:wa_decay_1D} for the 1D model. These results are summarized in Theorem \ref{thm:1D_profile_prop}.

In this paper, we construct 3D profiles and prove Theorems \ref{thm:self-similar},
\ref{thm:main_blowup}.

\vs{-0.05in}
\subsubsection*{\bf Step 3. 3D profiles}\label{sec:idea_step3} 
We use the 1D $\al$-profile $\waa$ in 
\hyr[sec:idea_step2]{\its Step 2} as an approximate profile 
for the 3D equation \eqref{eq:profi_3D0},
\footnote{
We cannot use $\wwwa$ or the 1D approximate profile $\wa$ as the approximate profile for 3D equation \eqref{eq:profi_3D0} since the residual error is not small in weighted $W^{1,\infty}$-norm. 
See discussion in \cite[Section 6]{chen2026eulerI}.
See also Footnote \ref{foot:direct_Wa}.
}
and formulate the fixed point problem for the perturbation \eqref{eq:fix_pt_2D}
$\FFF: \om \to \eta$. 

\vs{0.1in}
\paragraph{\bf Compactness of $\FFF$}\label{sec:idea_comp}
Given $\ww$, we solve $\eta$ along the characteristics in Proposition \ref{prop:reg_solu}.
Since the nonlocal parts $\tcB(\ww), \Lpsi(\ww)$ are more regular than $\ww$ 
in \emph{any bounded domain}, the map $\FFF: \om \to \eta$
  gains regularity,  which we use to prove compactness of $\FFF$. 
To obtain compactness of $\FFF$ for very large $|\xx|$, 
we  prove stronger far-field estimates for $\rhoc \f{\eta}{\wwb}$ than for $ \f{ \ww }{\wwb}$; see Section \ref{sec:3D_comp}, where $\rhoc$ is defined in \eqref{def:rhoc}.

\vs{0.1in}

To prove that $\FFF$ maps a ball $Y$ in some weighted $W^{1,\infty}$ space into itself, we proceed as follows.
\subsubsection{\bf Contraction for $\eta$ along $r=0$}\label{sec:idea_r0}
We build on the key 1D contraction estimate in Theorem \ref{thm:contra_lin_exact}. 
Let $\psi$ be the \emph{3D stream function} associated with $\ww$.
Along $r = 0$, \eqref{eq:fix_pt_2D}  reduces to
\beq\label{eq:idea_step3_1}
    ( 2 \vaa(z) + \psio(0,z) )  \pa_z  G = 
\big( - (1 - \al) \psio_z  -  2 \psio   \tf{ \pa_z \wwb }{\wwb} \big) (G+1) , 
   \quad G = \tf{\eta}{\wwb} , \quad \psi = \BS(\ww).
\eeq
Our key observation is that after dropping nonlinear terms, the only difference between 
\eqref{eq:idea_step3_1} and the 1D equation $\eta = \cLa(\ww)$ is the difference 
between two stream functions  $I =\psi(\ww) - \psid(\wwd)$ with $\wwd(z)= \ww(0, z)$,  where $\cLa$ is the 1D operator in Theorem \ref{thm:contra_lin_exact}. We bound error $I$ by purely \emph{bulk} norm $\nnr{\ww}$ for $\pa_r \ww$ using Lemma \ref{lem:vel_bc}.
Applying  Theorem \ref{thm:contra_lin_exact}, we establish
\beq\label{eq:idea_bc}
  \ncXc{ \rhoc \eta } \leq \tf{1 + \lamcL}{2} \ncXc{\ww} + C \nnr{\ww} + l.o.t. , \quad 
  \tf{1 + \lamcL}{2} < 1 ,
\eeq
in Proposition \ref{prop:3D_bd}, where $C$ is an absolute coupling constant and it may be large.

To estimate the bulk part of $\eta$, we have the following  structures for 
stream function and flow.
\paragraph{\bf Nonlocal estimates and $\cJaa$-growth}\label{sec:idea_nonlocal}
Due to the slow decay of the profile or perturbation, the main part for $\psi = \BS(\ww)$ 
is  formally  given by 
\[
\bal
 \psio & = - 2 \al z \JJ(\ww) + O( |z \cdot \ww| ), 
 & \quad  \psio_z & = - 2 \al \JJ(\ww) + O( |\ww|), \\  
\eal
\]
where the ring-variable is defined as $\psio = \psi -\pa_z \psi(0) z$ \eqref{def:psio}, 
$\JJ$ is defined in \eqref{def:JJ_3D}, and the lower order term is measured in a suitable norm. %  for the lower order part. 
We have logarithmic growth for the main term $\JJ(\ww)$: 
\[
1 - \JJ(\wwb) \gtr \cJaa, \quad  | \JJ(\ww) | \les |\cJaa|^{\kp+1}  \nnla{\ww} .
\]
The $\cJaa$ function \eqref{eq:Ja_hat} can be seen as a $\e^{-1}$-regularization of $\lgp x$. Other first and second order mixed derivatives of $\psio$ have $O(1)$  size 
relative to $\ww$. See  Proposition \ref{prop:vel_est} and Corollary \ref{cor:u_bar}. To compensate the $\cJaa$-growth in nonlocal estimates, we 
exploit the far-field cancellations for the 1D profile in Theorem 
\ref{thm:1D_profile_prop}, e.g. \eqref{eq:1D_prof_est2:c}, \eqref{eq:1D_prof_dxx}, \eqref{eq:1D_prof:va}. 

\begin{remark}[\bf Perturbative estimates for $\psi$]
Using the contraction estimate for the 1D operator in Theorem~\ref{thm:contra_lin_exact}
and exploiting the small parameter \(\e\), we treat the nonlocal term %perturbation term 
\(\psio\)
in the estimates for $\pa_r \eta, \pa_z \eta$ in Sections \ref{sec:idea_drW}, \ref{sec:idea_dzW}
\emph{perturbatively}, without tracking $\e$-independent absolute constants explicitly. In this sense, the 1D profiles \(\wwwa,\waa\) are \emph{strongly stable} in the fixed-point formulation.
\end{remark}

\paragraph{\bf Anisotropic flow}

In Proposition \ref{prop:Q}, for perturbation $\ww$ small in the weighted norm, we show that the self-similar flow $ \QQ = \bar \QQ + \td \QQ(\ww)$ \eqref{def:Q}  satisfies 
\beq\label{eq:idea_aniso}
 Q^r  \gtr \e^{-1} r, \quad  Q^z \gtr z  \cJaa(\xx) + l.o.t.,
\quad \tf{1}{r}  Q^r - \tf{1}{z}  Q^z \gtr \e^{-1} \ang \xx^{\al-\hal} + l.o.t. ,
\quad \cJaa \asymp \min( \lgp \xx, \e^{-1}) ,
\eeq
where $\al-\hal \asymp -\e$ is some $\e$-size exponent \eqref{ran:ep2}. We also establish comparable upper bounds. Since $\e^{-1} \gtr \cJaa$ when $\lgp \xx \les \e^{-1}$, the flow is strongly anisotropic.

\subsubsection{\bf A toy model for anisotropic estimate}\label{sec:model_aniso}

To motivate the following  estimates for $\pa_r \eta, \pa_z \eta$, we consider  a  model problem
with anisotropic flow capturing \eqref{eq:idea_aniso} for $\lgp \xx \ll \e^{-1}$
\beq\label{eq:model_idea}
 \pa_t  f + ( \e^{-1} r \pa_r + z \pa_z ) f = c  \cdot f  ,
\eeq
with some $\e$-independent constant $|c| \les 1$.  We aim to develop robust stability estimates so that we can further control $O(1)$ additional terms.
Firstly, taking $\pa_r f$, we obtain
\[
 \pa_t \pa_r f +  ( \e^{-1} r \pa_r + z \pa_z ) \pa_r f = ( c - \e^{-1} ) \cdot \pa_r f  .
\]

Since $\e \ll 1$, we obtain a strong damping term $ (c - \e^{-1}) \pa_r f$, leading to exponential decay for $\pa_r f$. 
Similarly, to obtain spatial decay estimates for $\pa_r f$, we estimate 
 $\ang \xx^{a} \pa_r f$ with $a < 1$ and obtain a strong damping term 
 $( c - \f{1-a}{\e}+l.o.t. ) \ang \xx^{a} \pa_r f$ for $  \ang \xx^{a} \pa_r f$. %with a size $ \f{1-a}{\e}$.
  In contrast, for $\pa_z f$, we only obtain 
\[
 \pa_t \pa_z f +  ( \e^{-1} r \pa_r + z \pa_z ) \pa_z f = ( c - 1) \cdot \pa_z f  ,
\]
with a weak damping term $(c-1) \pa_z f$. If we estimate $ \ang \xx^a \pa_z f$ with some $a>0$, 
we obtain a strong \emph{anti-damping term} from the advection $ ( c-1 +  a \ang \xx^{-2} ( \e^{-1} r^2  + z^2  ) ) \cdot \ang  \xx^a  \pa_z f $ with a positive coefficient when $r \gtr |z|$. 
On the other hand, we can estimate $z \pa_z f$ (or $\ang z \pa_z f$) without generating 
an $O(\e^{-1})$ \emph{anti-damping term} from the advection.

Model \eqref{eq:model_idea} motivates anisotropic weighted estimates for
$\ang \xx\,\pa_r f$ and $z\pa_z f$ in the time-independent problem.
The stability mechanism is driven mainly by \emph{outgoing transport}: the
$r$-advection rapidly transports perturbations toward infinity, producing strong stability.

\vs{-0.05in}
\subsubsection{\bf Contraction for $\pa_r \eta$ }\label{sec:idea_drW}

The key step is to prove that the bulk part $\pa_r \eta$ is small globally. Motivated by the above model, to capture the sharp \emph{global} decay for $\pa_r G$, which is crucial to close the estimate, 
we estimate the weighted quantity $ \ang \xx \rag \pa_r G$ using \eqref{eq:fix_pt_2DG}: 
\beq\label{eq:idea_step3_3}
\bal
  \QQ \cdot \na (  \ang \xx \rag  \pa_r G ) & =  \B(\udb{ - \pa_r  Q^r + \f{  \QQ \cdot \na \ang \xx}{\ang \xx} }_{ II_r} +  \udb{ \f{ \QQ \cdot \na \rag}{\rag} }_{  I_{\mw{ag}}   }  + 
  \bar \cR  + \Lpsi  \B) 
 \cdot \ang \xx \rag \pa_r G \\
 & \quad + \ang \xx \rag ( - \pa_r  Q^z \pa_z G + ( \pa_r \Lpsi + \pa_r \bar \cR  )  \cdot (G+1)) ,
 \eal 
\eeq
where $\Lpsi$ and $\QQ$ depend on the input $\ww$ nonlocally. 
The scaling fields 
$c_l \xx$ in $ \QQ$  are \emph{exactly cancelled} in $II_r$.
\footnote{
 This is expected: $\ang \xx\,\pa_r G$ has a decay rate similar to that of $G$,
while the $G$-equation \eqref{eq:fix_pt_2DG} does not contain a large damping
term from the scaling field $c_l\xx$.
}
For small perturbation $\ww$, we have $\QQ \approx \bar \QQ$, 
 and only get a partial $\e^{-1}$ damping 
\beq\label{eq:idea_damp_IIr}
II_r \leq  - C \e^{-1} \f{\ang z^2}{\ang \xx^2} \cdot \ang \xx^{\al-\hal}, \quad \al-\hal \asymp - \e.
\eeq

In any bounded region $B_R$ with $R$ independent of $\e$, we have $O(\e^{-1})$ strong damping.

\vs{0.1in}
\paragraph{\bf Angular weight}

To obtain a full $\e^{-1}$ damping, we use the crucial angular weight 
$\rag$ to capture the anisotropic flow \eqref{eq:idea_aniso}. 
The damping coefficient $I_{\mw{ag}}$ reads
\beq\label{eq:damp_ag_iden}
I_{\mw{ag}} =  \f{  \QQ \cdot \na \rag}{\rag}
  =  - \kag \cdot \f{r}{r + \ang z} \cdot (  \f{  Q^r}{r} -  \f{Q^z}{z} +  \f{Q^z}{z} \cdot \f{1}{\ang z^2} ) .
  \eeq

Using the crucial anisotropic structure 
\eqref{eq:idea_aniso} that  $\f{1}{r}  Q^r$ is much stronger than $\f{1}{z} Q^z$,
we obtain a large damping in the $r$-direction:  
\beq\label{eq:damp_ag_iden2}
I_{ag} \leq - C \e^{-1} \kag  \f{r}{\ang \xx} \cdot \ang \xx^{\al-\hal} + l.o.t.
\eeq
similar to \eqref{eq:idea_damp_IIr}. Combining the angular damping $I_{\mw{ag}}$ with $II_{r}$ in \eqref{eq:idea_damp_IIr}, we obtain
a \emph{full} strong $\e^{-1}$ damping term in \eqref{eq:idea_step3_3} for $\ang \xx \rag \pa_r G$.

We show that the coupling coefficient $\pa_r Q^z$ except the \hyr[sec:idea_singular]{\its singular term} in $ \pa_r  Q^z \pa_z G$ has $O(1)$ size.
We also show that  $\bar \cR = O(1)$, and  $\rhoo \ang \xx \pa_r( \Lpsi + \bar \cR) \cdot (G+1)$ 
 except the \hyr[sec:idea_singular]{\its singular term}  has a  $O(  \nnr{ \ww} )$ size, where $\rhoo$ is the weight chosen in \eqref{def:3d_wg}. 
Using the $\e^{-1}$ strong damping, we estimate the non-singular terms in \eqref{eq:idea_step3_3} 
and the $O(1)$-coupling term $C \nnr{\ww}$ in \eqref{eq:idea_bc} perturbatively.

We choose $\rhoo$ in \eqref{def:3d_wg} so that it is comparable with the weight in 1D. See Remark \ref{rem:3d_wg}. This leads to the weight $\rhoo \cdot \rag \ang \xx$ and norm  $\wwra$ in \eqref{def:3d_wg} and \eqref{def:3d_norm}.

For any $ k_1 \geq \ell_1 >0$ and $k, \ell>0$, 
since $\cJaa(y)^{-1} \asymp |\lgp y|^{-1} + \e$ by \eqref{eq:Ja_hat}, we obtain 
\beq\label{eq:err_ag_decay}
\bal
  \ang z^{-\ell_1} (\f{\ang z}{\ang \xx})^{k_1} \leq \ang \xx^{- \ell_1} ,  
\quad   \cJaa(z)^{-\ell}  (\f{\ang z}{\ang \xx})^{k}  \les_{\ell, k}  \cJaa(\xx)^{-\ell} ,
\quad \xx = (r, z) .
\eal
\eeq
Using estimate \eqref{eq:err_ag_decay} with an angular weight, we transfer $z$-decay estimates to
$|\xx|$-decay estimates.

\begin{remark}[\bf Improved regularity for $\psi$]\label{rem:psi_improved}
While using the angular weight $\rag \les 1$ weakens the control of $ \ww, \pa_r \ww, \pa_z \ww$
in the norms \eqref{def:3d_norm}, 
the Biot--Savart law \eqref{eq:Euler2_psi}, \eqref{eq:5D_BS} gains regularity from $\ww \to \na \psi, \psi$. As a result, we control the weighted quantities \(\rhoo\,\na^i\psi\), \(i\le 1\) with the 
stronger radial weight \(\rhoo\), rather than the weaker weight \(\rhoo\rag\). 
See the discussion below Remark \ref{rem:no_decay_r}.

\end{remark}

\vs{-0.05in}
\subsubsection{\bf Contraction for $\pa_z \ww$ }\label{sec:idea_dzW}

Variable $\pa_z \ww$ satisfies an equation similar to \eqref{eq:idea_step3_3}. 
In contrast to the estimate for $\pa_r \ww$, in $z$-direction the scaling field $c_l z$ 
and velocity $ U^z$ \emph{strongly cancel}, resulting in a much weaker outgoing effect for $ Q^z$ \eqref{eq:idea_aniso}, without an $\e^{-1}$ factor.
 Motivated by the  model \eqref{eq:model_idea}, we can only obtain decay estimates in $z$ direction.
 \footnote{
This is expected since the profile $\pa_z \wwb(r, z) = \pa_z \waa(z)$ \emph{only} decays in $z$.
 }
 We estimate the weighted term $ \ang z \rag \pa_z \ww$, which has 
 damping coefficient $- C + I_{\mw{ag}}$ with $C = O(1)$. 
 Away from $r=0$, we have a much stronger damping term due to \eqref{eq:damp_ag_iden}.
Since we have strong contraction estimates 
 for $\ww(0, z), \pa_r \ww $ and the nonlocal terms $\na^i \psio, i \leq 2$ are more regular than 
 $\na \ww$, we use interpolation to control $\na^i \psio,i\leq 2$ by 
 \[
     \nu^{\al} \nnza{ \ww} + C \nu^{\al-1} \nnrr{\ww} , \quad \  \forall  \, \nu \leq \tf12 ,
 \]
where the norms are defined in \eqref{def:3d_norm}. Taking $\nu$ suitable small, we control the nonlocal terms in $\pa_z G$ equation \emph{perturbatively}.
 Moreover, there are no singular terms as \eqref{eq:idea_singu} in  $\pa_z G$-equation. Hence, we can prove a contraction estimate for $\pa_z G$ using only a \emph{weak} global damping 
 of size $O(1)$.

\subsubsection{\bf Integration by parts along trajectories for singular terms}\label{sec:idea_singular}

The most difficult term $S$ arises from the coupling term $ \pa_r Q^z \cdot \pa_z G$ in the weighted estimate for $   \rhoo \rag \ang \xx    \pa_r G$
\eqref{eq:idea_step3_3} :
\beq\label{eq:idea_singu}
 S =   \rhoo \rag \ang \xx    \pa_r ( r \pa_r \Psi ) \cdot \pa_z G 
 = \tf{\ang \xx}{\ang z } \cdot \pa_r ( r \pa_r \Psi ) \cdot  ( \rhoo \rag \ang z \pa_z G )
 \teq P_1 \cdot P_2.
\eeq
Since $ r \pa_r^2 \Psi$ has only weak decay $\ang \xx^{-O(\e)}$ (see Proposition \ref{prop:vel_est}) 
 and $P_2$ is controlled only by the energy, without any additional decay,
\footnote{
The energy for estimating the fixed point map is defined in \eqref{def:local_norm}, 
which is slightly stronger than those in \eqref{def:3d_norm}. 
}
the key obstruction is the unbounded factor $\tf{\ang \xx}{\ang z }$ when $r \gg |z|$. This is a \emph{structural difficulty}: the approximate profile $\wwb, \pa_z \wwb$ do not decay in $r$, so one can only obtain partial decay of $\pa_z G$ in $z$-direction.

To overcome the above difficulty,  along the characteristic, we apply integration by parts to $S$
and estimate $S_1, S_2$ and lower order terms instead:
\[
S \to  S_1, \, S_2, \, l.o.t. , \quad S_1 = 
 \rhoo \rag \ang \xx    \pa_r ( z \pa_z \Psi ) \cdot \pa_z G 
 , \quad S_2 =   \rhoo \rag \ang \xx    \pa_r  \Psi  \cdot (\QQ \cdot \na)( \pa_z G ).
\]
See Lemma \ref{lem:IBP}. Although $\na\pa_z G$ is not controlled by the energy, applying
 equation \eqref{eq:fix_pt_2D} \emph{again} obtains that $\QQ\cdot\na\pa_z G$ has the \emph{same
derivative order} as $\na G$.
\footnote{
In later estimates, we work along the flow map \(X(s,\xx)\) associated with \(\QQ\), and use
$\tf{d}{ds}(\pa_z G)\cc X(s,\xx)$ in place of \((\QQ\cdot\na)\pa_z G\); see \eqref{eq:traj_smooth}.
}
 Moreover, the coefficient $\pa_r\Psi$ in $S_2$ is
more regular than $\pa_r(r\pa_r\Psi)$ in $S$.
Thus, up to lower order or more regular terms, we exchange the $r \pa_r$-derivative from $\Psi$ 
in $S$ \eqref{eq:idea_singu} to $z \pa_z$-derivative of $\Psi$ in $S_1$.
Since the \emph{mixed derivative} $\ang \xx\,\pa_{rz}\Psi$ enjoys \emph{much better} decay estimates (see Proposition \ref{prop:vel_est} and Corollary \ref{cor:u_bar}), 
we bound $S_1$ using such an estimate and $ |\ang z \rhoo \rag  \cdot \pa_z G| \leq \nnz{G}$  \eqref{def:3d_norm}.
We apply the same integration by parts to estimate other singular terms 
in $ ( \pa_r \Lpsi + \pa_r \bar \cR  )  \cdot (G+1)$ in \eqref{eq:idea_step3_3}. 
 The above argument shows that the average effect of $S$ along the trajectory
is much smaller.

\vs{0.1in}
\paragraph{\bf Error terms}

Since we use the exact 1D profile $\wwb = \waa$, the residual error $\bar \cR$ \eqref{eq:lin_2Db} \emph{vanishes} on $r=0$. It has a $O(1)$ size away from $r=0$. 
Up to log factors and a factor of $\ang \xx^{\al-\hal}$, we prove
\beq\label{eq:idea_cR}
  | z \pa_z \bar \cR| \les  r \ang z^{-1}, 
  \quad  | \ang \xx \pa_r \bar \cR_{\mw{reg}} | \les 1 ,
\eeq
in Proposition~\ref{prop:3D_error}, where \(\bar \cR_{\mw{reg}}\) denotes the regular part of \(\bar \cR\), with the singular part treated using methods in Section \ref{sec:idea_singular}.

By comparing the $\e^{-1}$ damping term for $\pa_r G$ in \eqref{eq:idea_step3_3} and the above $O(1)$ 
error, the error $|\pa_r \bar \cR_{\mw{reg}}|$ only contributes $O(\e)$ term in the estimate for $\pa_r G$. Using the crucial vanishing of $z\pa_z \bar \cR$ on $r=0$ 
\eqref{eq:idea_cR} and $\e^{-1} \f{r}{\ang \xx}$ angular damping \eqref{eq:damp_ag_iden2},  the error only contributes $O(\e)$ term in the estimate for $\pa_z G$. The effect of the error in these estimates may be captured by a damped ODE 
\beq
   \tf{d}{d t} a(t) = -\e^{-1} a(t) + 1, \quad a(0 ) = 0, 
\eeq
where $1$ models the $O(1)$ error. It is easy to obtain $|a(t)| \les \e$ using a bootstrap argument.

Using the above methods, we prove in Theorem \ref{thm:into} that  $\FFF$ maps some weighted $W^{1,\infty}$ ball with radius  $\e^{1-6\kp_1}$ into itself, where $\kp_1$ 
is a small parameter  in \eqref{def:kp}.

\vs{0.05in}
\paragraph{\bf Continuity of $\FFF$ and existence of 3D profile}

Proving that the map $ \FFF(\ww)$ is continuous in $\ww$ in the weighted $W^{1,\infty}$ space
is non-trivial since the equation \eqref{eq:fix_pt_2D} for $\eta=\FFF(\ww)$ is \emph{quasilinear} 
and $\eta$ only has $C^{1,\al}$ regularity. By contrast, it is much easier to prove that the nonlinear map $\FFF$ is closed using the uniqueness of the solution.
Combining \hyr[sec:idea_comp]{\its Compactness} and closedness of $\FFF$, we prove 
continuity. Using a Schauder fixed-point argument, we construct the 3D profile.

\subsubsection{\bf Proof ideas for Theorem~\ref{thm:main_blowup}}
\label{sec:idea_main_blowup}

To prove Theorem~\ref{thm:main_blowup}, we first upgrade the estimates for $\ws$ obtained from the fixed-point argument. By iterating trajectory estimates and elliptic estimates in Section~\ref{sec:3D_profile_prop}, we obtain higher regularity and decay estimates; see  Proposition \ref{prop:3D_profile}. Note that establishing the decay estimates is highly nontrivial: the \emph{approximate} profile $\wwb$ is constant in $r$, while the \emph{exact} profile gains about $C_{\e}\ang \xx^{-1/3}$ decay in the $r$-direction. We then establish nonlinear finite codimension stability of the profile $\ws$ using the equation for $\ww=r^{-\al}\ww^\th$ in a low-regularity regime. Compared with the weighted $H^k$, $k\gg1$, framework in \cite{chen2024Euler,chen2024vorticity,bedrossian2026finite}, a key difficulty here is that $\ws$ has only $C^{1,\al}$ regularity and weighted $W^{2,\infty}$ regularity in a suitable weighted space.

For linear stability estimates, we decompose the linear operator $\cA$ (see \eqref{eq:dyn_lin2}) as 
\beq\label{eq:lin_intro}
\cA = \cT + \cU , \quad 
    \cT \ww \teq - \QQs \cdot \na \ww + \cBs \cdot \ww,
    \quad \cU \ww \teq  \tcB(\ww) \cdot \ws - \td \QQ(\ww) \cdot \na \ws .
\eeq
We use the crucial  outgoing property \eqref{eq:intro_outgo0} and 
singularly weighted estimates to establish full stability estimates for the \emph{local} parts $\cT$, and construct the semigroup $e^{\cT s}$ by perturbing the \emph{transport} semigroup 
generated by $ - \QQs \cdot \na \ww$. Using the regularity and decay of the profile in Proposition \ref{prop:3D_profile}, 
we prove that the nonlocal part $\cU$ is a compact perturbation to $\cT$.
Using semigroup method, we obtain  decay estimates for $ e^{\cA s} =  e^{ (\cU + \cT) s}$. 

For the nonlinear estimates, we build on the splitting method and Duhamel representation in \cite{ChenHou2023a} and the fixed-point argument in \cite{chen2024Euler,chen2024vorticity}. We reformulate constructing a global solution in the self-similar equation as a fixed point problem with a map $\cH$. Due to the low regularity of the profile $\ws$, we perform stability estimates in some
weighted $L^{\infty}$ spaces, and prove that $\cH$ maps a ball $B(\d, Y)$ in some space-time weighted $L^{\infty}$ space $Y$ to itself. The ball $B(\d, Y)$ serves as a closed convex set in the Schauder fixed-point argument.  By choosing weighted $C^{1,\al}$ initial data and using a BKM-type continuation criterion, the \emph{quantitative} weighted $L^{\infty}$ stability estimates provide 
\emph{qualitative} weighted $W^{1,\infty}$ estimates \emph{for free}, which is sufficient 
for proving compactness of $\cH$. 

We further prove closedness of the nonlinear map $\cH$ using uniqueness of the solution, and obtain a global solution using a Schauder fixed-point argument.
By regularizing unstable directions, we allow the potentially unstable part of initial data, i.e. $X_2$ in Remark \ref{rem:set_init}, to lie in $C_c^\infty$. The asymptotically self-similar blowup results then follow from the nonlinear finite codimension stability of the profile.

\subsection{Notations}

We collect here the main notation used throughout the paper. 
For each variable or operator, we list right of it the equation or result in which it is first defined or determined.

We reserve the letters $ \alb, \al, \b, \e ,\hal, \hau $ for the parameters related to the regularity index $\al$ or $\e$:
\[
\alb \teq  \tf{1}{3}, \quad \e \teq \tf13 - \al, \quad   
\b: \, \mw{Thm} \, \ref{thm:1D_profile_prop}, 
    \quad \he \teq  \e - \b : \eqref{ran:ep2}, 
    \quad  \epa :  \eqref{norm:Xc}, 
  \quad \hal, \hau:  \eqref{eq:wa_upper_lower}
\]

We use parameters $\kp, \kp_1, \epa$ in the weight \eqref{norm:Xc}, $\kag$ in \eqref{def:3d_wg}, 
\[
 \kp , \, \kp_1:\eqref{def:kp}, 
 \quad  \kag :  \eqref{def:rag}, 
 \quad \bbb : \eqref{norm:Xc}, \eqref{def:3d_wg}
\]

We use $\lam$-notations for contraction or stability estimates 
\[
\lamcL: \mw{Thm} \, \ref{thm:contra_lin_exact},
\quad \lam_Q : \mw{Prop} \, \ref{prop:reg_solu}, 
 \quad  \bar \lam_1 : \mw{Lem} \, \ref{lem:3D_damp}, 
\quad \lam : \mw{Lem} \, \ref{lem:tran_wg}. 
\]

We use the notation $\vp$ and $\rhoo$ for weights and $\cJaa$ in the weight: 
\[
\vpa, \, \vpcc , \, \vpk : \eqref{norm:Xc},  \quad \rhoo, \, \rag, \, \rhor , \, \rhoz, \ \rhoc : \eqref{def:3d_wg},
\quad \cJaa : \eqref{eq:Ja_hat}, \eqref{eq:wa}
\]

We use the notation $\cW$ for weighted $W^{1,\infty}$ norms and notation $L^{\infty}$ for weighted $L^{\infty}$
\[
 \wwra, \, \wwza, \, \wwrr, \, \wonea : \eqref{def:3d_norm},
 \quad   \linfza, \, \linfa : \eqref{def:3d_norm}.
\]

We use calligraphic letters to denote operators or their values on the profile:
\[
\cJa: \eqref{eq:Jw}, \  \JJ : \eqref{def:JJ_3D:b}, \quad 
 \Lpsi , 
  \bar \cB , \, \td \cB, \, \bar \cR : \eqref{eq:lin_2D}, \ \cLa  : \eqref{eq:lin_1D_alpha}, \  
 \FFF \, \eqref{eq:fix_pt_2Deta},  
 \   \cT, \cU, \cA : \, \eqref{def:cT}.
\]

We use ``bar"-notation for the exact or approximate profile, or related variables 
e.g. $\waa, \wwb, \QQb, \bar \cR, \bcl, \bcw $; we use ``tilde"-notation for the perturbation, e.g. $\td \QQ, \tcB$. We use lowercase variable for perturbation, e.g. $\psi, \ww$, and uppercase variable to denote the background solution, e.g. $\Psi, \Om, \bpsi, \wwb$.

We use $A \teq B$ to say that quantity $A$ is defined to be equal to quantity/expression $B$.

\section{Anisotropic weighted estimates for the Biot-Savart law}\label{sec:nonlocal}

In this section, we establish several anisotropic weighted estimates of the Biot-Savart law 
associated with the anisotropic norms  \eqref{def:3d_norm}
in Proposition \ref{prop:vel_est} and Corollary \ref{cor:u_bar},
and isotropic weighted estimates in Proposition \ref{prop:iso_est}. In Lemma \ref{lem:vel_bc}, we estimate the difference between the 3D stream function and 1D stream function. Since these estimates are quite technical, the reader may first skip their proof on a first reading and return to it later.

Our main nonlocal anisotropic weighted estimates for the nonlocal terms are the following.
\begin{prop}
[\bf Transfer of $z$-decay in $\ww$ to global $|\xx|$-decay in $\psi$ ]\label{prop:vel_est}
Let $\e = \f13 - \al \in [0, 10^{-3}]$, $\psi(\ww) = \BS(\ww)$ be the stream function associated with $\ww$ \eqref{eq:Euler2_psi}, \eqref{eq:5D_BS}. 
Suppose that $\ww$ is odd in $z$. Denote $\xx= (r, z)$.  Consider the following weights $\rhoo , \rag$  
\bseq\label{def:u_wg}
\beq
\bga
  \rhoo(\xx) = (|\xx|^{- \bbu + 1} + 1)  \ang \xx^{ - \epa  }  \cJaa^{-\kp},  \quad
   \rag = \big( \f{ \la z \ra }{ r + \ang z } \big)^{\kag},
  \ega
\eeq
with parameters satisfying 
\beq\label{def:u_wg:b}
\kp \in [0, 1], \quad  \kag \in [0, \tf{1}{1000}], \quad  \bbu  \in [1, 1.3 ],  \quad \epa \in [0, \e].
 \eeq
\eseq

Let $\linfa, \wwra, \wwza, \wwrr, \wonea $ be the norms  in  \eqref{def:3d_norm} associated with the above $\rhoo, \rag$. 
For any $\nu \in (0, \f12 ]$,  we have 
\footnote{
We only get an estimate involving $\rag^{-1}$ for $\psi_{zz}$ in the upper bound in \eqref{eq:vel_zz}.
 The upper bound is still bounded. We use the weight $\rag^{-1}$ to simplify the upper bound.
}
\bseq\label{eq:vel_est}
\begin{align}
 \rhoo  r | \psi_{rr} | & \les \min( |\xx|^{\al+1} , |\xx|^{ \al - \hal  } ) \nnrr{\ww} ,\label{eq:vel_rr} \\
  \rhoo |\psi_{r z} | & \les \min( |\xx|^{\al} , |\xx|^{ \al - \hal - 1 } )
  \cdot ( \nu^{\al} \nnza{ \ww} + \nu^{\al-1} \nnrr{\ww}  ) , 
  \label{eq:vel_rz}  \\
  \rhoo  |\psi_{zz}| & \les \min( |\xx|^{\al}, \rag^{-1} \ang \xx^{\al-1} \ang z^{-\hal}  )
   ( \nu^{\al} \nnza{ \ww} + \nu^{\al-1} \nnrr{\ww} ) ,
    \label{eq:vel_zz} 
\end{align}
where $\hal$ is defined in \eqref{eq:wa_upper_lower} with  $\hal - \al \asymp \e$ \eqref{ran:ep2}. By choosing $\nu = \f{ \nnrr{ \ww } }{ 2 \max(\nnz{\ww}, \nnrr{\ww}) } \leq \tf12$ and using the norm $\wonea$ \eqref{def:3d_norm}, we optimize the above estimates and obtain
\beq\label{eq:vel_bd_opt}
   \nu^{\al} \nnza{ \ww } + \nu^{\al-1} \nnrr{\ww}
   \les  \nnz{\ww}^{1-\al} \nnrr{\ww}^{ \al} 
   +  \nnrr{\ww} \les \nnn{\ww}.
\eeq

Recall $\psio = \psi - \psi_z(0)z$ \eqref{def:psio} and $\JJ$ from \eqref{def:JJ_3D}. We have the folowing estimates for $\psi, \na \psi$
\begin{align}
  \rhoo | \psi_r | & \les   |z| \min( |\xx|^{\al},  | \xx |^{\al - \hal -1 } )  \nnn{\ww},
    \label{eq:vel_r}  \\
  \rhoo | \psi_r | & \les  |\xx| \min( |\xx|^{\al},  | \xx |^{\al - \hal -1 } )  \nnla{\ww},
    \label{eq:vel_r1}  \\
   \rhoo | \psio_z + 2 \al \JJ(\om) | +   \rhoo  | \tf{1}{z} \psio   + 2 \al  \JJ(\ww)| 
   +   \rhoo |\psi_z  - \tf{1}{z} \psi |  & \les  \min( |\xx|^{\al+1}, | \xx|^{\al -\hal} ) \nnla{\ww}  \label{eq:vel_psi}  .
\end{align}
\eseq

Let $\cJaa$ be as in \eqref{eq:wa}, \eqref{eq:Ja_hat}.
For $\JJ$ \eqref{def:JJ_3D},   we have $\psi_z(0) =  2 \al \JJ(\infty)$ and  
\bseq \label{eq:vel_J}
\begin{align}
      |\JJ(\xx)| & \les \min( |\xx|^{1+\al}, \cJaa^{\kp+1} ) \nnla{\ww} ,  \\
      |\JJ(\xx) - \JJ(\infty)| & \les \e^{-1- \kp}  \ang \xx^{\al -\hal + \epa}  \nnla{\ww} .
  \end{align}
\eseq
All the implicit constants are \emph{independent} of $\al, \hal, \hau, \e$ and parameters in \eqref{def:u_wg:b}.

\end{prop}

Estimate \eqref{eq:vel_r} is stronger than \eqref{eq:vel_r1}, yet we need to use the stronger norm $\nnn{\ww}$ to obtain  \eqref{eq:vel_r}. We apply \eqref{eq:vel_r} to obtain the vanishing factor $|z|$ when $|z| \ll |\xx|$.

\begin{remark}[\bf Notation of weights]

The weights in \eqref{def:3d_wg} satisfy the assumptions in \eqref{def:u_wg}. We formulate the estimates for a broader class of weights, so that they apply later to the nonlocal terms associated with the profile and perturbation. We keep the same notation \(\rhoo,\rag\) in \eqref{def:u_wg} and \eqref{def:3d_wg} 
 to make the later nonlocal estimates easier to compare with the above results.

\end{remark}

\begin{remark}[\bf Angular averaging for global decay]\label{rem:no_decay_r}

As motivated in Sections \ref{sec:model_aniso}, \ref{sec:idea_drW}, and \ref{sec:idea_dzW}, 
we can \emph{only} obtain decay estimates for \(\pa_z\ww, \ww\) in the $z$-direction. By exploiting 
the property that the Biot--Savart law \eqref{eq:Euler2_psi}, \eqref{eq:5D_BS} gains regularity 
from $\ww \to \na \psio, \psio$, 
we transfer the $z$-decay of $\ww$ and $\pa_z \ww$
to the  decay of $\psio$ in $(r, z)$ in Proposition \ref{prop:vel_est};
see e.g. \eqref{eq:vel_rz}-\eqref{eq:vel_psi}.
Moreover, we control weighted estimates for $\psi$ with weights 
$\rhoo$  (except \eqref{eq:vel_zz}) 
stronger than the those for $\ww$ with weights $\rhoo \rag$. 
See the definition of norms in \eqref{def:3d_norm}.

\end{remark}

In estimate \eqref{eq:vel_est}, the radial weight $\rhoo$ essentially commutes with the nonlocal operators; the reader may therefore treat $\rhoo$ as if it were $1$ in the proof of \eqref{eq:vel_est}.

We have the following estimates for the approximate profile.

\begin{cor}\label{cor:u_bar}
Let $\bar \Psi = \BS( \wwb )$ \eqref{def:3d_vel} be the stream function associated with $\wwb$,
$\mr{ \bpsi } = \bpsi - \bpsi(0)z $ \eqref{def:psio}, $\JJ$ be as in \eqref{def:JJ_3D}, 
$\hal, \hau$ be the parameters in \eqref{eq:wa_upper_lower}, and $\xx= (r, z)$. We have 
\beq\label{eq:u_bar}
\bal
   | \bar \Psi_{rz} | & \les  \min( |\xx|^{\al}, |\xx|^{\al - \hal - 1} ) ,  \quad && | \bar \Psi_{zz} |  \les   \min( |\xx|^{\al},  \ang \xx^{\al-1}  \ang z^{ - \hal  } ) , \\ 
  |\bar \Psi_r| & \les |z| \min( |\xx|^{\al},  | \xx |^{\al - \hal -1 } ) , \quad  &&  | \bar \Psi_z -
  \tf{1}{z}  \bar \Psi  | \les  \min(|\xx|^{\al+1}, |\xx|^{\al - \hal}), \\
  |r \bar \Psi_{rr}| & \les  \min(|\xx|^{\al+1}, |\xx|^{\al - \hal})  .
\eal
\eeq

  For $|z| \geq r$, we have
  \beq\label{eq:u_bar_rzz}
     |\pa_{r zz} \bar \Psi | \les \min( |\xx|^{\al-1}, \ang \xx^{\al-\hal -2} ) \, ,
  \eeq
and 
\beq\label{eq:u_bar:J}
  \bal
 |  \mr{\bpsi}_z + 2 \al \cJ_{3D}(\wwb)| 
 +  | \tf{1}{z} \mr{\bpsi} + 2 \al \cJ_{3D}(\wwb)|  & \les  \min(|\xx|^{\al+1}, |\xx|^{\al-\hal} ), \\
 |\cJ_{3D}(\wwb)(\xx)| & \les \min(|\xx|^{1+\al}, \cJaa(\xx)).
  \eal
\eeq

Let $\xx = (r, z)$. For $\JJ$ \eqref{def:JJ_3D} and $z \geq 0$, we have the following estimates 
\bseq\label{eq:JJ_bound}
\begin{align}
  1- \JJ( \wwb )(\xx) &\gtr \cJaa(\xx),  \label{eq:JJ_bound:ne} \\
  \quad  \JJ(\wwb)(\xx) - \JJ(\wwb)(\infty)  & \gtr \e^{-1} \ang \xx^{\al - \hau}
  \label{eq:JJ_bound:far}  .
\end{align}
\eseq
All the implicit constants are \emph{independent} of $\al,\hal, \hau, \e$.

\end{cor}

\begin{proof}

We apply Proposition \ref{prop:vel_est} with $\om = \wwb$ and  weights $\rhoo, \rag$ admitting parameters:
$ \bbu = 1, \kag = 0 , \kp = 0$, and $\epa = 0$. Using \eqref{eq:waa_reg} and  \eqref{eq:wa_upper_lower},  and $\ww = \wwb = \waa$, we yield 
\[
\bga
  \rhoo \equiv  1, \quad \rag \equiv 1, \quad 
  (|\xx|^{-1} + 1) \ang z^{\hal} |\om |
  \les    (|\xx|^{-1} + 1)   \min( |z|, 1 )  \les 1, 
  \quad \pa_r \ww \equiv 0 ,  \quad  | \ang z^{\hal + 1} \pa_z \ww| 
  \les 1.  
  \ega
\]

From the definition of norms \eqref{def:3d_norm},  we obtain 
\[
 \nnla{\ww} \les 1, \quad 
 \nna{\ww} \les 1, \quad  \nnr{\ww} = 0,  \  \quad \nnz{\ww} \les 1, \
  \quad \nnn{\ww} \les 1.
\]

Taking $\nu = \f12$ in Proposition \ref{prop:vel_est} and using \eqref{eq:vel_est}, \eqref{eq:vel_psi},
we prove \eqref{eq:u_bar} and \eqref{eq:u_bar:J}.

Next, we prove \eqref{eq:JJ_bound}. Let $\xx = (r, z)$.  Since 
$ -\wwb= -\waa \gtr \min( z, z^{-\hau} )$ for $z \geq 0$ by \eqref{eq:waa_sign}
and $\wwb$ is odd in $z$, using definition \eqref{def:JJ_3D} 
and the polar coordinate  $ (\tr, \tz) = (R \cos \b, R \sin \b)$,   we obtain 
\[
\bal
  1 -   \JJ(\wwb)(\xx)
 & \geq C + C \int_{ 0\leq \tz \leq |\xx|, \tr \geq 0 }  \f{ \tr^{2 + \al} \tz }{ |(\tr, \tz)|^5 } 
\min( |\tz| , |\tz|^{-\hau} ) d \tr d \tz  \\
& \geq C + C \int_{R \leq |\xx|} \int_{\pi/8}^{\pi/4} R^{\al -1} (\cos \b)^{2+\al} \min( R (\sin \b)^2, 
(R \sin \b)^{1-\hau} ) d R d \b \, .
  \eal
\]

Since $\al - \hau \in [-2\e, -\e/2]$ by \eqref{ran:ep2},  we obtain
\[
\bal
    1 -   \JJ(\wwb)(\xx)&  \gtr  
    C + C \int_{0}^{  |\xx|} R^{\al -1} \min( R ,
R^{-\hau} ) d R  \gtr C + C \one_{|\xx| \geq 2} \cdot \f{1 - |\xx|^{\al -\hau}}{ |\al - \hau|} \\
& \gtr C + \one_{|\xx| \geq 2} \min(\e^{-1} , \lgp x ) \, .
\eal
\]

Similarly, we have 
\[
   II\teq \JJ(\wwb)(\xx) - \JJ(\wwb)(\infty) 
 \geq  C \int_{  |\xx| \leq  \tz , \tr \geq 0   }  \f{ \tr^{2 + \al} \tz }{ |(\tr, \tz)|^5 } 
 \min( |\tz| , |\tz|^{-\hau} ) d \tr d \tz .
\]
Since $\al -\hau \asymp  -\e$ by \eqref{ran:ep2} and $ |\xx| \leq  \tz $ implies $|\xx|  \leq R = (\tr, \tz)$, passing to the polar coordinate and using estimates similar to the above, we prove 
\[
    II 
     \geq C \int_{ |\xx|}^{\infty}R^{\al -1} \min( R ,
R^{-\hau} ) d R
\geq C \cdot \f{ \min( |\xx|^{\al-\hau}, 1) }{\e}
\geq C \f{\ang \xx^{\al-\hau}}{\e}.
\]
We prove \eqref{eq:JJ_bound}.

For \eqref{eq:u_bar_rzz}, we defer the proof to Section \ref{sec:psi_rzz}.
\end{proof}

\subsection{Estimate kernels and functions}\label{sec:BS_est_ker}

Recall the stream function associated with $\ww$ \eqref{eq:5D_BS}:
\bseq\label{eq:5D_BS:rec}
\beq
  \psi(r, z) = (-\D_{5D})^{-1} ( \cpsia \ww r^{\al-1} )
  = \f{ \cpsi}{2\pi^2} \int \f{1}{ | (r, 0,0,0, z) - \yy|^3 } W(\yy) d \yy
  =     ( \KRF \ast  W)(\xx) \, , 
\eeq
for $\cpsi$ chosen in \eqref{def:cpsi}, where $\KRF$ denotes 
a multiple of the fundamental solution for $-\D$ in $\R^5$ 
\beq\label{eq:5D_BS:ker}
  K_{\R^5}(s) = \f{ \cpsi}{2 \pi^2} |s|^{-3}, 
\eeq
\eseq
$W(\yy)$ is radially symmetric in $y_1, y_2, y_4, y_4$ and is given by
\beq\label{eq:psi_W}
W(\yy) = \ww(\tr, \tz) \tr^{\al-1} , \quad 
\yy \in \R^5, \quad 
 \yy_h = (y_1, y_2, y_3, y_4), \quad  \tr = |\yy_h|, \ \tz = y_5 , 
\eeq
and we abuse the $\xx$-notation to denote 
 \bseq\label{def:x_y}
\beq\label{def:5D_x}
  \xx = (r, 0,0,0 , z),  \quad x_1 = r, \quad 
\xx_h = (r,0,0,0),  \quad x_5 = z .
\eeq

We use $\tr,  \tz \in \R$ and $\yy \in \R^5$ to denote the variables in the integral. 
Since $|\yy| = |(\tr, \tz)|$, $|\xx| = |(r, z)|$, 
for any function $f$ only depending on $\tr, \tz$, we also abuse notation to write 
\beq\label{eq:y_abuse}
|\yy| = |(\tr, \tz)|, \quad |\xx| = |(r, z)|,
\quad  f(\yy) = f(\tr, \tz) .
\eeq

\eseq

Below, we use ``boldface" font to denote vectors, e.g. $\xx, \yy$. 
We introduce
\bseq\label{eq:psi_A}
\begin{align}
      A(r , s,  \tr, \th) & =  ( r^2 + \td r^2 - 2 r \td r \cos \th + s^2 )^{1/2}, \\
    D_{\pm}(r, s, \tr) &= (| r \pm \tr|^2 + s^2)^{1/2} . \label{def:Dpm}
\end{align}

Using the definition of $\xx$ and $\yy$, we obtain
\beq
  |\xx-\yy| = ( (r - \tr \cos \th)^2 + |\tr \sin \th|^2 + (z -\tz)^2  )^{1/2}
  = A(r, z - \tz, \tr , \th ).
\eeq

\eseq

\subsubsection{ Estimate of $W$ }

Using a direct calculation on $W$ \eqref{eq:psi_W}, 
$\tr = |\yy_h|$, and $ \pa_{y_1} \tr  = \f{y_1}{\tr}$, we yield
\begin{align}\label{eq:ker_W_id}
W(y)  & = \tr^{\al-1}  \om(\tr, \tz),  \notag \\
  \pa_\tz W(y) & = \td r^{\al-1}  \pa_{\tz} \om(\td r, \td z), \\
   \pa_{y_1} W(y) & =    \pa_{y_1} \tr \cdot \pa_\tr W
  =  \f{y_1}{\tr} \pa_\tr ( \td r^{\al-1} \om(\td r, \td z))    = \f{y_1}{\tr} ( (\al-1) \tr^{\al-2}  \om(\td r, \td z)
  + \tr^{\al-1} \pa_{\tr} \om  ). \notag
\end{align}

Recall $|\yy| = |(\tr, \tz)|$ \eqref{def:x_y}. Using the norms $\linfa, \wwrr, \wwra, \wwza$ in \eqref{def:3d_norm},  we get 
\bseq\label{eq:ker_wwlinf}
\begin{align}
  |\om(\tr, \tz)|  &\les
\min( |\yy|, 1 ) \ang \tz^{-\hal}
   \cdot \rhoo^{-1} \rag^{-1} \nnla{\ww} \, , \label{eq:ker_wwlinf:a} \\
   |\pa_r \om(\tr, \tz) | & \les  \ang \tz^{-\hal}   \ang  \yy^{-1} \rhoo^{-1} \rag^{-1}  \nnr{\ww} 
   \les  \ang \tz^{-\hal}  \ang  \yy^{-1} \rhoo^{-1} \rag^{-1}  \nnrr{\ww} ,  \\
    |\pa_{\tz} \om(\tr, \tz)|  & \les \ang z^{-1 - \hal}  \cdot \rhoo^{-1} \rag^{-1} 
    \nnz{\ww}.
  \end{align}
\eseq

Since $|y_1| \leq \tr , \tr \ang \yy^{-1}\les \min(|\yy|, 1)$ and 
$W(\yy) = \om^{\th} \cos \th =  \td r^{\al-1} y_1 \om( \td r, \td z )$, 
combining the above two estimates, we bound
\beq\label{eq:ker_Wlinf}
\bal
  | W(\yy)| & \les  \tr^{\al-1} |\om| 
   \les \min( |\yy|, 1 ) \tr^{ \al-1 } \ang \tz^{-\hal} \rhoo^{-1} \rag^{-1} \nnla{\ww}, \\
  |\pa_{y_1} W(\yy)|&  \les \tr^{\al-1} |\om |   + \tr^{\al} | \pa_r \om|
  \leq ( \min( |\yy|, 1 )  \tr^{\al-2} 
  +  \tr^{\al-1} \ang \yy^{-1}  ) \ang \tz^{-\hal} \rhoo^{-1} \rag^{-1} \nnrr{\ww} \\
  & \les ( \min( |\yy|, 1 )  \cdot  \tr^{\al-2}  
  \ang \tz^{-\hal} \rhoo^{-1} \rag^{-1} \nnrr{\ww} .
\eal
\eeq

We further introduce the polar coordinate for $(\tr, \tz)$
\beq\label{eq:polar_trz}
  \tr = R \cos \xi, \quad  \tz = R \sin \xi ,
  \quad R = |(\tr, \tz) | = |\yy|, \quad \xi \in [-\pi/2, \pi/2].
\eeq

By definition of $\rag$ \eqref{def:rag}, we obtain
\bseq\label{eq:wg_est}
\beq\label{eq:wg_est:ag}
  \rag(\tr, \tz)^{-1}
  \asymp (\f{ \la \tr, \tz \ra}{  \ang \tz} )^{\kag} 
  = (\f{R^2 + 1}{R^2 \sin^2 \xi + 1})^{\kag/2}
   \leq |\sin \xi|^{-\kag} = \f{|(\tr,\tz)|^{\kag}}{|\tz|^{\kag}}.
\eeq

Since $\rhoo(\tr, \tz)$ \eqref{def:u_wg} only depends on $|(\tr, \tz)|$, we abuse notation to write 
$\rhoo(\yy) = \rhoo(\tr, \tz)$. Since $\bbu \in [1,2)$ 
and $\cJaa \asymp \min( \lgp x, \, \e^{-1} )$, for any $m > 0$, we obtain
\beq\label{eq:wg_est:rb}  
\bal    
       \rhoo( \pp ) \asymp_m \rhoo(\qq)  , \quad 
       \forall \ m^{-1}|\qq|  \leq |\pp| \leq m |\qq| .
  \eal
\eeq

Since $\rhoo^{-1}(\pp)$ is increasing in $|\pp|$, if $|\yy| \leq m |\xx|$, we obtain
\beq\label{eq:wg_est:rc}  
\rhoo^{-1}(\yy) 
\leq \rhoo^{-1}( m \xx )
\les_m \rhoo^{-1}(\xx).
 \eeq

\eseq

For $ |\yy| = |(\tr, \tz)| \leq 4 |\xx|$, using \eqref{eq:wg_est:rc} and 
 \eqref{eq:wg_est:ag}, we bound 
\[
\bal
\rhoo^{-1}(\yy) \rag(\yy)^{-1} \ang \tz^{-\hal} 
& \les \rhoo^{-1}(\xx)
\cdot \f{\ang \yy^{\kag}}{ \ang \tz^{\kag}} \ang \tz^{-\hal}
\les \rhoo^{-1}(\xx) \cdot \ang \xx^{\kag} \ang \tz^{-\hal - \kag} .
\eal
\]

Thus, for $|\yy| = |(\tr, \tz)| \leq 4 |\xx|$, from the upper bound in \eqref{eq:ker_wwlinf},
 we further obtain 
\bseq\label{eq:ker_wwlinf_ne}
\begin{align}
|\om(\tr, \tz)| & \les P(\xx) \cdot   \min( |\xx|, 1 )  \ang \tz^{-\hal -\kag }  
\nnla{\ww}, \\
|\pa_r \om(\tr, \tz)| & \les P(\xx) \cdot  \ang \yy^{-1}   \ang \tz^{-\hal -\kag }    \nnrr{\ww}, \\
|\pa_z \om(\tr, \tz)| & \les 
P(\xx) \cdot
 \ang \tz^{-\hal -\kag -1}     \nnz{\ww} \, ,
\end{align}
where we denote 
\beq
  P(\xx) = \ang \xx^{\kag} \rhoo^{-1}(\xx) .
\eeq

From the assumption \eqref{def:u_wg:b} and \eqref{ran:ep2}, we obtain the following range 
\beq\label{para:range}
\bbu \in [1, 1.3 ], \quad \kag \in [1, \tf{1}{1000}], \quad \hal , \al \in [0.32, 0.34], \quad  \hal + \kag \in [0.32, 0.34]. 
\eeq

\eseq

\subsection{Symmetrized estimates of the kernel \texps{$\KRF$}{ K R5 } }\label{sec:sym_ker}

In this section, we estimate the kernel $\KRF$ for $-\D_{\R^5}$ and its symmetrization.
For any $s \in \R^n$ and $f: \R^n \to \R$, we denote 
\beq\label{def:df}
 (\pa_i f)(s) = \pa_{s_i} f(s).
\eeq

We recall the notation from \eqref{eq:psi_W} and \eqref{eq:5D_BS:ker}
\beq\label{eq:xy_recall}
\KRF(s) = \f{ \cpsi}{2 \pi^2} |s|^{-3}, 
\quad  \yy_h = (y_1, y_2, y_3, y_4), \quad  \tz = y_5,  \quad   z = x_5.
\eeq

We have the following basic estimates of the kernel
\bseq\label{eq:K5D_est}
\beq\label{eq:K5D_est:dK}
  |\na^k \KRF(s)| \asymp_k |s|^{-k-3}  \, .
\eeq
For any $\xx, \yy$  with $|\yy| \geq 2 |\xx|$, we obtain $|\yy-\xx| \asymp |\yy|$ and 
\beq\label{eq:K5D_est:far}
  |\na_{\xx, \yy}^k  \KRF(\xx- \yy)| \asymp_k |\xx-\yy|^{-k-3}  \asymp_k |\yy|^{-k-3}.
\eeq

\eseq

Recall $z = x_5$ from \eqref{def:5D_x}. Consider any $|\yy| \geq 4 |\xx|$.
Since $ \pa_{s_5} \KRF(s)$ is odd in $s_5$ and even in $s_i, i\leq 4$,  using \eqref{eq:K5D_est:far} 
and  the second order Taylor expansion of $\pa_5 \KRF$ at $-\yy$, we obtain
\[
\bal
  K_{\R^5, z}^{\mw{sym}}(\xx, \yy) \teq &      \pa_{z} \KRF(\xx- \yy) - \pa_z \KRF(  \xx_h-\yy_h , x_5 + y_5)  %\\
 + \pa_{z} \KRF( \xx_h + \yy_h, x_5 - y_5) 
- \pa_{z} \KRF( \xx + \yy)  \\
=  &  (\pa_5 \KRF)(\xx- \yy) + (\pa_5 \KRF)  ( \xx_h - \yy_h,  -x_5 -y_5) \\
& \quad  + (\pa_5 \KRF)( -\xx_h -\yy_h, x_5 - y_5)   + (\pa_5 \KRF) ( -\xx-\yy ) \\
= & 4 (\pa_5 \KRF)(-\yy)
+ \na \pa_5 \KRF(-\yy) \cdot ( \xx + ( \xx_h , -x_5 ) + (-\xx_h, x_5) + (-\xx) ) %\\
+ O( |\xx|^2 |\yy|^{-6} )  \\
  =  &  4 (\pa_5 \KRF)(-\yy) + O( |\xx|^2 |\yy|^{-6} ) .
\eal
\]

Since $\pa_{s_1} \KRF(s)$ is even in $s_5$ and odd in $s_h$: $\pa_{s_1} \KRF(s) = - \pa_{s_1} \KRF(-s_h, s_5)$, for $|\yy| \geq 4 |\xx|$, using \eqref{eq:K5D_est:far} and the second order Taylor expansion of $\pa_5 \KRF$ at $\xx=0$, we obtain
\[
\bal
K_{\R^5, r}^{\mw{sym}} 
& \teq  \pa_{x_1} \KRF(\xx-\yy) - \pa_{x_1 } \KRF( \xx_h - \yy_h, x_5  + y_5 )
 + \pa_{x_1} \KRF( \xx_h + \yy_h, x_5 - y_5  )
 - \pa_{x_1} \KRF( \xx + \yy  )   \\
& =  (\pa_1 \KRF)( \xx-\yy ) -  
 ( \pa_{ 1 } \KRF)(  \xx_h - \yy_h, - x_5 - y_5 ) \\
& \quad  - ( \pa_1 \KRF)( -\xx_h - \yy_h, x_5 - y_5  )
 + ( \pa_1 \KRF)( -\xx-\yy ) \\
& = \pa_1 \KRF(-\yy) \cdot (1-1-1+1) 
+ \na \pa_1 \KRF(-\yy) \cdot ( \xx - (\xx_h, -x_5) - (-\xx_h, x_5) + (-\xx) ) 
+ O(|\xx|^2 |\yy|^{-6}) \\
& =  O(|\xx|^2 |\yy|^{-6} ).
\eal
\]

For any $|\yy| \geq 4 |\xx|$, 
since $\KRF(s)$ is even in $s_i$, 
for fixed $\xx_h, \yy_h$, using third order Taylor expansion at $x_5=0$ and \eqref{eq:K5D_est:far}, we obtain 
\[
  \bal
 &  \KRF( \xx_h \pm \yy_h, x_5 - y_5) - \KRF(   \xx_h \pm \yy_h,  x_5 + y_5 )   = \KRF( \xx_h \pm \yy_h, x_5 - y_5) - \KRF(   \xx_h \pm \yy_h,- x_5 -y_5 )  \\
 & \quad  \qquad =   2 x_5 (\pa_5 \KRF)  ( \xx_h \pm \yy_h, - y_5)  
+ O( |x_5|^3 \sup\nolimits_{|\xi| \leq x_5}  |\pa_5^3 \KRF( \xx_h \pm \yy_h, \xi-y_5) | )  \\
& \quad \qquad = 2 x_5 (\pa_5 \KRF)  ( \pm \xx_h + \yy_h, - y_5)   + O( |x_5|^3 |\yy|^{-6} ) . 
\eal
\]

Using Taylor expansion at $\xx_h =0$ and the above estimate
, we obtain the following symmetrized estimates 
\[
\bal
  K_{\R^5}^{\mw{sym}} 
 &  \teq \KRF( \xx_h - \yy_h, x_5 - y_5) - \KRF(   \xx_h - \yy_h,  x_5 + y_5 )
   +    \KRF( \xx_h + \yy_h, x_5 - y_5) - \KRF(   \xx_h + \yy_h,  x_5 + y_5 ) \\
  & = 4 x_5 (\pa_5 \KRF)( \yy_h , - y_5 )
  + 2 x_5  (\na \pa_5 \KRF)(\yy_h, -y_5) \cdot ( (\xx_h, 0) + (-\xx_h, 0) )
  + O( |x_5| \cdot |\xx|^2 |\yy|^{-6} ) \\
& = 4 x_5 (\pa_5 \KRF)( \yy_h , - y_5 ) + O( |x_5| \cdot |\xx|^2 |\yy|^{-6} ) 
=  4 x_5 (\pa_5 \KRF)( -\yy ) + O( |x_5| \cdot |\xx|^2 |\yy|^{-6} ) .
\eal 
\]

Recall $z =x_5$. Thus, for any $|\yy| \geq 4 |\xx|$, we prove the following symmetrized estimates 
\bseq
\begin{align}
   | K_{\R^5, z}^{\mw{sym}}(\xx,\yy) -  4 (\pa_5 \KRF)(-\yy) | & \les |\xx|^2 |\yy|^{-6}, \label{eq:K5D_est:symz} \\ 
 |K_{\R^5, r}^{\mw{sym}}(\xx, \yy) | & \les  |\xx|^2 |\yy|^{-6} ,  \label{eq:K5D_est:symr}  \\
      | K_{\R^5}^{\mw{sym}}(\xx, \yy) -  4 z (\pa_5 \KRF)(-\yy) | &  \les |z| \cdot |\xx|^2 |\yy|^{-6} 
      \label{eq:K5D_est:sym0}  .
\end{align}
\eseq

\subsection{Angular averaging in integrals}

In this section, we establish integral estimates that capture angular averaging effects in the far field and near the singularity. We first establish the far-field estimates, where angular averaging makes the decay $\ang{\tz}^{-\hal}$ comparable to the global decay $\ang{\yy}^{-\hal}$ and yields $\rag\asymp 1$.

\begin{lem}[\bf Transfer of $z$-decay to global decay]\label{lem:int_W}
Let $\rhoo, \rag$ be defined in \eqref{def:u_wg}.
Denote $\rhoo(\yy) = \rhoo(\tr, \tz), \rag(\yy) = \rag(\tr, \tz), \tr = |\yy_h|, \tz = y_5$. 
For any $ \ell \geq 4 $, we have
\beq\label{eq:int_W_far}
 \int_{| \yy | \geq 4 |\xx| , \R^5} |\yy|^{-\ell} 
    \min(|\yy|, 1) \ang \tz^{-\hal} 
 \rhoo^{-1} \rag^{-1}( \tr, \tz )   \tr^{\al -2} d \yy
 \les 
 \rhoo(\xx)^{-1} 
  \min( |\xx|^{\al + 4 -\ell} + |\xx|^{-\al},  |\xx|^{\al - \hal + 3 - \ell} ) .
\eeq
with an implicit constant independent of $\ell$. For any $ |\xx| \geq 0$, we have 
\beq\label{eq:int_W_near}
 \int_{  \tz \leq |\xx|  } 
|\tz| \ang \yy^{-1} |\yy|^{\al - 2}
 \rhoo^{-1} \rag^{-1}(\tr, \tz)
 \cdot \ang \tz^{-\hal}
d \tr d \tz 
\les  \min( |\xx|^{1 + \al}, \cJaa(\xx)^{1 + \kp} ).
\eeq
For any $|\xx| \geq 0$ and any constant $c_1, c_2 > 0$ independent of $\xx$, we have 
\beq\label{eq:int_W_near2}
 \int_{  \tz \leq c_1 |\xx| , \tr \geq c_2 |\xx| } 
 \ang \yy^{-1} |\yy|^{\al - 2}
 \rhoo^{-1} \rag^{-1}(\tr, \tz)
 \ang \tz^{-\hal}
d \tr d \tz 
\les_{c_1, c_2}    \rhoo^{-1}( \xx ) \min( |\xx|^{\al}, \ang \xx^{\al-\hal-1} ).
\eeq

\end{lem}

\begin{proof}

Since the integrand only depends on $\tr, \tz$, using \eqref{eq:5D_int_id} with $k=0$ 
and $\tr \leq |\yy|$, we obtain
\[
  \mathrm{LHS}_{\eqref{eq:int_W_far} } =  \int_{| (\tr, \tz) | \geq 4 |\xx|} |\yy|^{-\ell} 
    \min(|\yy|, 1) \ang \tz^{-\hal} 
 \rhoo^{-1}( \tr, \tz ) \rag^{-1}( \tr, \tz )   \tr^{\al + 1} d \tr d \tz  ,
\]
where we abuse the notation to write $|\yy| = |(\tr, \tz)|$. Below, for any $\ell \geq 4$, we prove 
\beq\label{eq:int_W_far1}
\bal
  M_2 & = \int_{| (\tr, \tz) | \geq 4 |\xx|} |\yy|^{-\ell+ \al + 2} 
\ang \yy^{-1} \ang \tz^{-\hal} 
 \rhoo^{-1} \rag^{-1}( \tr, \tz )   d \tr d \tz  \\
& \les |\cJaa(x)|^{\kp}    \min(1 + |\xx|^{\al + \bbu + 3 - \ell}, |\xx|^{ \al - \hal + 3 -\ell + \epa }) .
 \eal
\eeq

Clearly, we have $  \mathrm{LHS}_{\eqref{eq:int_W_far} }  \leq M_2$. 
Passing to the polar coordinate $(R, \xi)$ of $(\tr, \tz)$, $d \tr d \tz = R d R d \xi$, 
and using $\rag^{-1} \les |\sin \xi|^{-\kag}$ \eqref{eq:wg_est:ag}, we bound 
\[
  \bal
M_2
   \les  \int_{  R \geq 4 |\xx|,  |\xi| \leq \pi/2 }
\rhoo(R)^{-1}\
\ang R^{-1} R^{-\ell + \al + 3}
   \min(1,   | R \sin \xi| ^{-\hal}  )
   |\sin \xi|^{- \kag}  
  d R d \xi .
   \eal 
\]

Since $ |\hal - \f{1}{3}| ,|\al - \f13| , \kag< \f{1}{10}$ by \eqref{para:range}
 it is easy to obtain that the above integrand is integrable in $\xi$. Thus, we obtain 
\[
M_2 \les 
      \int_{ R \geq 4 |\xx|}
\rhoo(R)^{-1} \cdot 
\ang R^{-1 -\hal} R^{-\ell+\al + 3 }  d R .
\]

Using $\rhoo^{-1} \les \min( R^{\bbu-1}, 1) \ang R^{\epa} \cJaa(R)^{\kp}$ 
with $\bbu \in [1, 3/2]$, and 
\beq\label{eq:cJaa_rat}
  \cJaa(R)^{\kp} \les_m  \ang \xx^{-0.1} \cJak \ang R^{0.1} ,\quad \forall  \ R \geq m |\xx| \, ,
\eeq
by \eqref{eq:Ja_hat}, we bound  
\[
\bal
M_2
& \les  
\ang \xx^{-0.1} \cJak
\int_{R \geq 4|\xx|} 
\min( R^{\bbu-1}, 1) \ang R^{  \epa + 0.1 -1 -\hal} R^{ -\ell + \al + 3  }  d R   \\
& \les \ang \xx^{-0.1}  \cJak   \int_{R \geq 4|\xx|} \one_{R \leq 1} R^{\bbu -\ell+ \al+ 2 }
+ \one_{R > 1} R^{-\ell + \al + 2   -\hal + \epa + 0.1} d R .  \\
\eal
\]
Since $\ell \geq 4, \bbu + \al + 2 < \f32 + \f13 + 2 < \ell, \al + 2 - \ell -\hal + \epa + 0.1 < \f52 - \ell \leq -\f32 $, we obtain
\[
 M_2  \les  \cJak   \min(1 + |\xx|^{\al + \bbu + 3 - \ell}, |\xx|^{ \al - \hal + 3 -\ell + \epa }) .
 \]
 with implicit constant independent of $\ell$. Thus, we prove \eqref{eq:int_W_far1}.

Using $\rhoo(\xx) = (1 + |\xx|^{-\bbu+1}) \ang \xx^{-\epa} \cJaa^{-\kp}$ \eqref{def:u_wg}, estimating the case $|\xx| \leq 1$ and $|\xx| > 1$, and using $\bbu < 1 + \al$ by \eqref{def:u_wg:b}, we further obtain
\[
\rhoo M_2 \les    \min( |x|^{\al + 4 -\ell} + |\xx|^{-\bbu+1},  |x|^{\al - \hal + 3 - \ell} )
\les   \min( |x|^{\al + 4 -\ell} + |\xx|^{-\al},  |x|^{\al - \hal + 3 - \ell} ) \, .
\]
Since $M \leq M_2$, rearranging $\rhoo$ in the above estimate, we prove 
\eqref{eq:int_W_far} in Lemma \ref{lem:int_W}.

\vspace{0.1in}
\paragraph{\bf Proof of \eqref{eq:int_W_near2} }

Below, the implicit constants allow to depend on $c_1, c_2$ in \eqref{eq:int_W_near2}. 
Recall $\rag, \rhoo$ from \eqref{def:3d_wg}.  In the domain of the integral, since $\tr \gtr |\xx| \gtr |\tz|, |\yy| \asymp \tr$, using \eqref{eq:cJaa_rat}, we obtain
\[
\bal
  I & = \ang \yy^{-1} |\yy|^{\al - 2}
 \rhoo^{-1} \rag^{-1}(\tr, \tz)
 \ang \tz^{-\hal} 
 \les  \ang \tz^{-\hal -\kag}  \min(|\yy|^{\bbu-1}, 1) |\yy|^{\al-2} \ang \yy^{-1 + \kag + \epa} |\cJaa(\yy)|^{\kp} \\
 & \les \ang \tz^{-\hal -\kag} 
 \min( \tr^{\bbu-1}, 1)
  |\tr|^{\al-2} \ang \tr^{-1 + \kag + \epa}  \cJaa(\xx)^{\kp} (  \f{\ang \tr}{ \ang \xx} )^{0.1}.
 \eal
\]

Integrating over $|\tz| \leq c_1 |\xx|, \tr \geq c_2 |\xx|$, 
 using $\al -3 + \kag + \epa + 0.1 < -2, \bbu \in [1,3/2]$ by \eqref{para:range}, \eqref{def:u_wg:b}, 
 and estimates similar to the above $M_2$, we prove
\[
\bal
  \int_{ \ssk{ |\tz| \leq c_1 |\xx|,\\ \tr \geq c_2 |\xx| } } I d \tz d \tr
 & \les \ang \xx^{-0.1} \cJaa(\xx)^{\kp} \int_{ |\tz| \leq c_1 |\xx| } \ang \tz^{-\hal -\kag} d \tz 
  \int_{ \tr \geq c_2 |\xx| } |\tr|^{\al-2}  \min( \tr^{\bbu-1}, 1)   \ang \tr^{-1 + \kag + \epa + 0.1}    \\
& \les \ang \xx^{-0.1} \cJaa(\xx)^{\kp} 
\min(|\xx|, \ang \xx^{1 -\hal -\kag}) \cdot ( |\xx|^{ \al + \bbu-2 } \one_{|\xx| < 1} + \one_{|\xx|\geq 1} |\xx|^{\al-2+\kag + \epa + 0.1} )  \\
& \les |\cJaa(\xx)|^{\kp} ( \one_{|\xx| < 1} |\xx|^{\al + \bbu -1} 
+ \one_{|\xx| \geq 1} \ang \xx^{\al-\hal + \epa-1}  )
\les \rhoo^{-1}(\xx) \min( |\xx|^{\al}, \ang \xx^{\al-\hal-1} ).
\eal
\]

\vspace{0.1in}
\paragraph{\bf Proof of \eqref{eq:int_W_near} }

We decompose the integral in $\mathrm{LHS}_{\eqref{eq:int_W_near} }$ and  use $|\tz| \les |\xx|$ and $|\tz| \les |\yy|$ to obtain
\[
\mathrm{LHS}_{\eqref{eq:int_W_near} } 
\les 
 \int_{  \tz \leq |\xx|  } 
 ( |\yy|^{\al - 1} \one_{\tr \leq |\xx|}
 + |\xx| \cdot |\yy|^{\al-2} \one_{\tr \geq |\xx|}  )
 \ang \yy^{-1}  
 \rhoo^{-1} \rag^{-1}(\tr, \tz)
 \ang \tz^{-\hal}
d \tr d \tz \teq I_1 + I_2.
\]

For the integral $I_2$ associated wit the second integrand, since $\rhoo^{-1} \les \ang \xx^{\epa} \cJaa^{\kp}$ and $\al-\hal +\epa<0$ by \eqref{ran:ep_all}, using \eqref{eq:int_W_near2}, we bound $I_2$ by 
$\mathrm{RHS}_{\eqref{eq:int_W_near} }$. Thus, it remains to bound $I_1$. 

Let $R, \xi$ be the polar coordinates for $(\tr, \tz)$. Recall $\rhoo , \rag$ from \eqref{def:3d_wg}. Using $\rag(\tr, \tz)^{-1}  \les |\sin \xi|^{-\kag}$ by \eqref{eq:wg_est:ag} and
$\rhoo^{-1} \les \cJak \ang \xx^{\epa}$ by \eqref{def:3d_wg}, and passing to the polar coordinate of $(\tr, \tz) = (R \cos \xi, R \sin \xi)$,  we estimate
\[
\bal
I_1 & \les |\cJaa(x)|^{\kp} \int_{ 0  }^{2 |\xx|} \int_0^{\pi/2}
 \ang R^{-1 + \epa} R^{\al }
 |\sin \xi|^{-\kag}
 \min( 1, ( R \sin \xi)^{-\hal} ) d R d \xi  .
\eal
\]
Since $\hal + \kag \in [0, 0.34]$ \eqref{para:range}, 
using Lemma \ref{lem:log_ineq_J} with $k=0$ 
and $\al - \hal + \epa \asymp -\e$ by \eqref{ran:ep_all}, we prove
\[
   I_1 \les |\cJaa(x)|^{\kp} \int_{ 0 }^{2 |\xx|} \ang R^{-1 -\hal +\epa} R^{\al } d R
   \les \min( |\xx|^{\al+1} , | \cJaa(\xx)|^{1 + \kp} ).
\]
Thus, $I_1$ is bounded by $ \mathrm{RHS}_{\eqref{eq:int_W_near} } $. Combining the 
estimates of $I_1, I_2$, we prove \eqref{eq:int_W_near}.
\end{proof}

The weighted estimates of the integral near the singularity are much more involved. 
We have the following estimates for the integral over the azimuthal angle.

\begin{lem}[\bf Averaging over the azimuthal angle $\th$]\label{lem:int_singu}
Let $p , q >0$ be some parameters. For $k > \ell+ 1.01, \ell \in [0,10]$, we get 
\beq\label{eq:int_singu0}
\int_0^{2} \f{t^\ell}{ (p^2 + q^2  t^2 )^{k/2} } d t
\les %\f{ k -\ell}{k-\ell-1} 
p^{-k} \min(  \f{p}{q}, 1 )^{\ell+1}.
\eeq
with absolute constants independent of $p,q, k, \ell$. 
Recall $A, D_{\pm}$ from  \eqref{eq:psi_A}, \eqref{def:Dpm}:
\[
   A(r , s,  \tr, \th) 
   =  ( r^2 + \td r^2 - 2 r \td r \cos \th + s^2 )^{1/2}, \quad 
    D_{\pm}(r, s, \tr) = (| r \pm \tr|^2 + s^2)^{1/2} .
\]

Recall the notation $\tr = |\yy_h|, y_5 = \tz $ for $ \yy \in \R^5$. For any $k \geq 4$, 
 $\xx = (r,0,0,0,z)$ \eqref{def:5D_x}, domain $\Ups \subset \R_+ \times \R$,  and function $f(\tr, \tz)$ only depending on $\tr, \tz$, we have 
\beq\label{eq:int_singu}
   \int_{\R^5} \one_{(\tr, \tz) \in \Ups} \f{ |f(\tr, \tz)|  }{|\xx-\yy|^k} d \yy 
\les \int_{ (\tr, \tz) \in \Ups } \f{  |f(\tr, \tz)|  \tr^3  }{ D_-(r, z -\tz ,\tr)^{k-3} D_+(r, z - \tz , \tr)^3 }d \tr d \tz.
\eeq

\end{lem}

We impose $k \geq 4$ in the assumption for \eqref{eq:int_singu} as we only use \eqref{eq:int_singu} with $k \geq 4$. 

\begin{proof}
A direct calculation yields 
\[
\mw{LHS}_{ \eqref{eq:int_singu0} }
\leq \int_0^{2} (\one_{t \leq \f{p}{q}} + \one_{t \geq \f{p}{q}}) \f{t^\ell}{ (p^2 + q^2  t^2 )^{k/2} } d t
\les \int_0^2 \one_{t \leq \f{p}{q}} \f{t^\ell}{p^{k}}
+ \one_{t > \f{p}{q}}  q^{- k} t^{\ell- k} .
\]
Since $ k > \ell +1.01, \ell \in [0, 10]$,  we prove \eqref{eq:int_singu0}:
\[
  \mw{LHS}_{ \eqref{eq:int_singu0} } \les p^{- k} \min( \f{p}{q}, 2 )^{\ell+1}
  +  \one_{ \f{p}{q} < 2 }  q^{- k} ( \f{p}{q} )^{\ell+1- k}
  \les  p^{- k} \min( \f{p}{q}, 2 )^{\ell+1}
  \les p^{- k} \min( \f{p}{q}, 1 )^{\ell+1}.
\]

For \eqref{eq:int_singu}, since $\one_{(\tr, \tz) \in \Ups}  |f(\tr, \tz)|$ only depends on $\tr, \tz$, using identity \eqref{eq:5D_int_id}, we obtain
\[
  I \teq
     \int_{\R^5} \one_{(\tr, \tz) \in \Ups} \f{ |f(\tr, \tz)|   }{|\xx-\yy|^k} d \yy 
    \les \int_{ (\tr, \tz) \in \Ups } \int_0^{\pi} \f{  |f(\tr, \tz)| \tr^3 \sin^2 \th   }{ A(r, z -\tz, \tr)^k }  d \tr d \tz d \th.
\]

Next, we integrate $\th$. From \eqref{eq:psi_A}, for $\th \in [0, \pi]$, we obtain
\[
  A(r, z-\tz, \tr)^2 = r^2 + \tr^2 - 2 \cos \th r \tr + (z-\tz)^2 
  =(r-\tr)^2 + (z-\tz)^2 + 4 \sin^2 \tf{\th}{2} \cdot \tr r
  \asymp D_-^2 + \tr r \th^2 \, .
\]
Below, we simplify $D_{\pm}(r, z -\tz, \tr)$ as $D_{\pm}$. Since $|\sin \th | \asymp |\th|$, applying \eqref{eq:int_singu0} with $p = D_-, q = (\tr r)^{1/2}$, $\ell=2$, and $k\geq 4 > \ell + 1.01$, we prove 
\[
I \les   \int_{ (\tr, \tz) \in \Ups } \f{  |f(\tr, \tz)| \tr^3    }{ D_-^{k}  } 
\min( \f{D_-}{ (r \tr)^{1/2} } , 1)^{3} d \tr d \tz   .
\]

Recall $D_{\pm}$ from \eqref{def:Dpm}. We have
\[
  \min( \f{D_-}{ (r \tr)^{1/2} } , 1)
  \asymp  \f{ D_- }{\max( D_-, (r \tr)^{1/2} )}
\asymp \f{ D_- }{  (D_-^2 +  r \tr)^{1/2} } \asymp \f{D_-}{D_+}.
\]
Combining the above two estimates, we prove \eqref{eq:int_singu}.
\end{proof}

We have the following weighted estimate near the singularity $(r, z)$.

\begin{lem}[\bf Average decay near the singularity]\label{lem:CS_singu}
For any $k /2 + \ell  < 0.99, k+ m < 1.99, n \in [0, 1]$ and $k, m ,\ell\geq 0$, we have
\beq\label{eq:CS_singu}
\bal 
   \int_{ |(\tr, \tz)| \leq 8 |\xx|}  \f{1}{ |(r -\tr, z- \tz)|^k  
|(r +\tr, z- \tz)|^m} \la \td z \ra^{- \ell} \ang {\tr, \tz}^{- n}
  d \tr d \tz  & \les |\xx|^{2- k-m} \ang \xx^{- \ell - n } . 
\eal 
\eeq
where $ |\xx| =| (r, z)|$, and the implicit constants are independent of $k,m, n, \ell$.
\end{lem}

The above bound shows that on average, the decay $\ang \tz^{-\ell} \ang {\tr, \tz}^{- n} $ is comparable to $\ang \xx^{-\ell-n}$, and the singular part $ |(r \pm \tr, z- \tz)|$ is comparable to $|\xx|$.

\begin{proof}

\textbf{ Case $m=0, n =0$}.
Below, we use notation \eqref{def:Dpm} and simplify $D_{\pm} = ( |r \pm \tr|^2 + |z-\tz|^2 )^{1/2}$. We first prove the estimate for $m=0$.   Since $ k /2 + \ell < 1$, we can choose $p, q > 1$ as
\[
p = (k/2 + \ell) (k/2)^{-1}, \quad  q = (k/2 + \ell) \ell^{-1}, 
\quad 
  p^{-1} + q^{-1 } = 1, \quad \ell q < 0.99,  \quad k p < 1.99.
\]

Using H\"older inequality, 
and $|(\tr, \tz) -(r,z)| \leq 8 |\xx|  +   |\xx| \leq 9 |\xx|$, we yield 
\[
  \mw{LHS} \leq  (\int_{ |(\tr,  \tz) - (r, z)| \leq 9 |\xx| } |(\tr, \tz)- (r, z)|^{- k p})^{1/p} ( \int_{|(\tr, \tz)| \les  |\xx|} \la z \ra^{- \ell q} d \tr d \tz)^{1/q} . 
\]
Since $k p < 1.99$ and $\ell q < 0.99$, we prove 
\[
  \mw{LHS} \les  |\xx|^{  \f{2-k p}{p} }  (  |\xx|  \min( |\xx|,   |\xx|^{1-\ell q}) )^{ \f{1}{q} }
  \les   |\xx|^{2 - k} \min(1 ,  |\xx|^{-\ell})
  \les  |\xx|^{2 -k} \ang \xx^{-\ell}.
\]

\vspace{0.1in}
\paragraph{\bf Case $m >0, n =0$}
Next, we consider  $n = 0, m >0$. 
 Using the definition of $D_+$, we obtain
\[
  \max( D_+, |\tz| ) \gtr r + \tr + |z-\tz| + |\tz| \gtr r + |z| \gtr |\xx|.
\]

Using $D_+ \gtr |\xx|$ or $|\tz| \gtr |\xx|$, and $D_+ \geq D_-$,  we obtain
\[
  \f{1}{D_-^k D_+^m} \ang \tz^{-\ell}
  \les |\xx|^{-m} \f{1}{D_-^k } \ang \tz^{-\ell}
  +  \f{1}{D_-^k D_+^m} \ang \xx^{-\ell}
  \les |\xx|^{-m} \f{1}{D_-^k } \ang \tz^{-\ell} + \ang \xx^{-\ell} \f{1}{D_-^{k+m}}.
\]
For the integral of the first term, we apply \eqref{eq:CS_singu} when $m=0$. 
For the integral of the second term, we estimate it directly. 
Since $k + m \leq 1.99$, we prove
\[
  \int_{|(\tr, \tz)|\leq 8 |\xx|}   \f{ \ang \tz^{-\ell}}{D_-^k D_+^m}
\les |\xx|^{-m} |\xx|^{2-k} \ang \xx^{-\ell}
+ \ang \xx^{-\ell} \int_{ |(r -\tr,   z - \tz | \leq 9 |\xx| }  |(r -\tr,   z - \tz |^{-k-m} d \tr d \tz
\les |\xx|^{2-k-m} \ang \xx^{-\ell}.
\]

\paragraph{\bf Case $m>0, n\in [0,1]$}
In the region $|(\tr, \tz)| \geq |\xx|/2$, 
using \eqref{eq:CS_singu} with $n =0$, we obtain 
\[
  I_1 =  \int_{ |\xx|/2 \leq  |(\tr, \tz) \leq 8 |\xx| }   \ang {\tr ,\tz}^{-n} \f{1}{D_-^k D_+^m} \ang \tz^{-\ell}
\les \ang \xx^{-n} 
 \int_{  |(\tr, \tz) \leq 8 |\xx| }  \f{1}{D_-^k D_+^m} \ang \tz^{-\ell}
 \les \ang \xx^{-n - \ell} |\xx|^{2-k-m}.
\] 

In the region $ |(\tr, \tz)| \leq |\xx| /2$, we obtain $|(\tr \pm r, \tz - z)|  \geq |\xx| - |\xx|/2 \gtr |\xx|$ and 
\[
\bal
 I_2 =  \int_{ |(\tr, \tz) | \leq |\xx|/2 }   \ang {\tr ,\tz}^{-n} \f{1}{D_-^k D_+^m} \ang \tz^{-\ell}
\les |\xx|^{-k-m}  \int_{ |(\tr, \tz) | \leq |\xx|/2  } \ang {\tr ,\tz}^{-n} \ang \tz^{-\ell} \, . \\
\eal
\]

Since $\ell \leq 0.99$ and $n \leq 1$, using the polar coordinate of $\tr = R \cos \xi, \tz = R \sin \xi$, we bound
\[
\bal
  I_2 & \les |\xx|^{-k-m} \int_0^{|\xx|/2} \int_0^{\pi} \ang R^{-n} \ang {R \sin \xi}^{-\ell} R d R d \xi
\les  |\xx|^{-k-m} \int_0^{|\xx|/2} \ang R^{-n} \ang R^{-\ell} R d R  \\
& \les |\xx|^{ - k -m}  \min( |\xx|^2, 1 + |\xx|^{2-\ell-n} )
\les |\xx|^{2 - k -m} \ang \xx^{-\ell-m} .
\eal
\]
Combining the estimate of $I_1, I_2$, we complete the proof.
\end{proof}

\subsection{Estimate of \texps{$ \psi_{rz}$ }{ psi,rz}, \texps{$r \psi_{rr}$}{r psi, r r} }\label{sec:psi_rz}

In the following sections \ref{sec:psi_rz}-\ref{sec:d_psi}, we estimate $\psi_{rz}, r \psi_{rr}, 
 \psi_{zz},  \psi_r, \psi_z, \psi$ in order. Below, we assume $z > 0, r > 0$ without 
loss of generality. 
For $\nnrr{\ww}, \nnz{\ww} < \infty$, we show in Lemma \ref{lem:psi_reg} that $\psi_z \in C^{1,\al}, r \psi \in C^{2,\al}$. Thus, the case $z =0$ or $r =0$ can be obtained using continuity.   We decompose the integral in two regions: 
(a) far-field $|\yy| \geq 4 |\xx|$ and (b) near the singularity $|\yy| \leq 4 |\xx|$.
We symmetrize the kernels in the estimate for $|\yy| \geq 4 |\xx|$.

The estimate of $\psi_{r z}$ with $r > z$ and $r \psi_{rr}$ are very similar, so we estimate them together.

Recall $\xx = (r, 0,0,0,z)$. Since $\pa_r \psi(r, z) = \pa_{ x_1 } \psi( r, z ) |_{y_1 = r} $, using \eqref{eq:psi_W} and formulas of $W$ in \eqref{eq:ker_W_id}, we obtain
\beq\label{id:psi_rz}
\bal
    \pa_{r z} \psi(r, z) & =  \pa_z \pa_{x_1} (\KRF \ast W )(\xx)
    =  P.V. (\pa_z \KRF \ast \pa_{y_1} W)(\xx)   \\
   & = P.V. \B( \int_{|\yy | \leq 4 |\xx|} + \int_{|\yy| \geq 4 |\xx|} \B) \pa_z \KRF(\xx-\yy) \cdot \pa_{y_1} W(\yy) d \yy
   \teq I_{  r z }^{\nne} + I_{r z}^{\ffar} ,
\eal
\eeq
where $P.V$ denotes the principle value of the integral near  $\xx$. Similarly, for $\pa_{rr} \psi$, we obtain
\beq\label{id:psi_rr}
\bal
   \psi_{rr}(r, z) & = \pa_{x_1}^2 (\KRF \ast W)(\xx)
  =  P.V. ( \pa_{x_1} \KRF \ast \pa_{y_1} W )(\xx)   \\
   & = P.V. \B( \int_{|\yy | \leq 4 |\xx|} + \int_{|\yy| \geq 4 |\xx|} \B) \pa_{x_1} \KRF(\xx-\yy) \cdot \pa_{y_1} W(\yy) d \yy \teq I_{  r r }^{\nne} + I_{r r}^{\ffar}  \, .
\eal
\eeq

\subsubsection{Estimate in the far-field for $\psi_{rz}, \psi_{rr}$}
For $I_{r r}^{\ffar}, I_{r z}^{\ffar}$, the integral is not singular. 
Since $\pa_{y_1} W(\yy)$ is odd in $\yy_h$ from \eqref{eq:ker_W_id} 
and $ \pa_5 \KRF(s)$ is even in $s_h$, 
symmetrizing the integral in $\yy_h$, we rewrite $  I_{rz}^{\ffar} $ in \eqref{id:psi_rz} as
\[
  I_{rz}^{\ffar} = \f{1}{2} \int_{|\yy| \geq 4 |\xx|}  ( \pa_z \KRF(\xx- \yy) 
 - \pa_z \KRF( - \xx_h - \yy_h, x_5 - y_5 ) ) \pa_{y_1} W(\yy) d \yy .
\]
Similarly, for the kernel in \eqref{id:psi_rr}, since $\pa_{y_1} W(\yy)$ is odd in $y_5$
and $\pa_{ s_1} \KRF(s)$ is even in $s_5$,  symmetrizing the integral in $y_5$, we obtain
\[
    I_{rr}^{\ffar} = \f{1}{2} \int_{|\yy| \geq 4 |\xx|}  ( \pa_{x_1} \KRF(\xx- \yy) 
 - \pa_{x_1} \KRF(  \xx_h - \yy_h, - x_5 - y_5 ) ) \pa_{y_1} W(\yy) d \yy .
\]

Using \eqref{eq:K5D_est:far} we obtain 
\[
\bal
     | \pa_z \KRF(\xx- \yy) 
 - \pa_z \KRF(  -\xx_h - \yy_h, x_5 - y_5 ) |  
 & \les |\xx|  \sup\nolimits_{|\xi| \leq |\xx|} |\na^2 \KRF(\xi - \yy)| \les |\xx| \cdot |\yy|^{-5} , \\
 |  \pa_{x_1} \KRF(\xx- \yy)  - \pa_{x_1} \KRF(  \xx_h - \yy_h, - x_5 - y_5 ) |
 & \les |\xx|  \sup\nolimits_{|\xi| \leq |\xx|} |\na^2 \KRF(\xi - \yy)| \les |\xx| \cdot |\yy|^{-5} . 
 \eal
\]

Applying the above estimate, \eqref{eq:ker_Wlinf} on $\pa_{y_1} W$,
 we bound $I_{rz}^{\ffar}, I_{rr}^{\ffar}$ as 
\[
  |I_{rz}^{\ffar} | + I_{rr}^{\ffar}|
  \les |\xx| \int_{|\yy| \geq 4 |\xx|} |\yy|^{-5} \min(|\yy|, 1) \ang \tz^{-\hal} 
  \tr^{\al-2} \rhoo^{-1}( \yy ) \rag^{-1}(\yy) d \yy \nnrr{\ww}.
\]
Applying Lemma \ref{lem:int_W} with $\ell = 5$, $\al \in [0.32, 0.34]$ \eqref{para:range}, 
and $\al- \hal \in [-2\e,0]$ \eqref{ran:ep2}, we bound 
\beq\label{eq:psi_Irz:far}
\bal
    |I_{rz}^{\ffar}(\xx) | 
+     |I_{rr}^{\ffar}(\xx) | 
   &  \les \rhoo^{-1} 
    |\xx| \min( |\xx|^{\al-1} + |\xx|^{-\al} , |\xx|^{\al-\hal -2} ) \nnrr{\ww} \\
   & \les \rhoo^{-1} \min( |\xx|^{\al}, \ang \xx^{\al-\hal -1} ) \nnrr{\ww}.
\eal
\eeq

\subsubsection{Estimate of $\psi_{rz}$ for $r \gtr |z|$ and $\psi_{rr}$ near the singularity}

 We recall $D_{\pm}, A$ from \eqref{eq:psi_A}
\beq\label{eq:Dpm:recall1}
    D_{\pm}(r, z -\tz, \tr) = (| r \pm \tr|^2 + |z-\tz|^2)^{1/2} .
\eeq
Below, to simplify notation, we simplify $  D_{\pm}(r, z -\tz, \tr)$ as $D_{\pm}$.

For $|\yy| \leq 4 |\xx|$, using \eqref{eq:ker_W_id} and \eqref{eq:ker_wwlinf_ne}, 
and $|y_1| \leq \tr \les |\xx|, \tr \leq \ang \yy$,  we bound 
\[
\bal
  |\pa_{y_1} W|  \les |\pa_\tr W | & \les  \tr^{\al-2} |\ww| + \tr^{\al-1} |\pa_{\tr} \ww|
  \les P(\xx) ( \min(|\xx|, 1)  + \tr \ang \yy^{-1} ) \tr^{\al-2} \ang \tz^{-\al-\kag}
  \nnrr{\ww} \\
  & \les P(\xx) \min( |\xx|, 1 ) \tr^{\al-2} \ang \tz^{-\al-\kag}
  \nnrr{\ww} .
\eal
\]

Since $ \pa_{\tr} W(\yy)$ only depends on $\tr, \tz$ by definition \eqref{eq:psi_W}
and $|\yy| = |(\tr, \tz)| \leq 4 |\xx|$,  using $|\na \KRF( s ) | \les |s|^{-4}$ by \eqref{eq:K5D_est:dK},
and then estimate \eqref{eq:int_singu} with $k = 4$,
we bound  $I_{rz}^{\nne}$ in \eqref{id:psi_rz} 
and $I_{rr}^{\nne}$ in \eqref{id:psi_rr}  as 
\[
\bal
  |I_{rz}^{\nne} | 
  +   |I_{rr}^{\nne} | 
 &  \les  \int_{|\yy| \leq 4 |\xx|} |\xx-\yy|^{-4} 
 |\pa_\tr W |
 d \yy   \les  
  \int_{ |(\tr, \tz) | \leq 4 |\xx|}  \f{  |\pa_\tr W |  \tr^3 }{D_- D_+^3} 
 d \tr d \tz .
\eal
\]

Using the above bound of $\pa_\tr W$, we obtain
\[
     |I_{rz}^{\nne} | +      |I_{rr}^{\nne} |  \les 
P(\xx) 
\min(|\xx| , 1) \nnrr{\ww}  \int_{ |(\tr, \tz) | \leq 4 |\xx|}  \f{  \tr^{\al+1} }{D_- D_+^3} 
    \ang \tz^{-\hal -\kag} d \tr d \tz.
\]

Below, we estimate $   |I_{rz}^{\nne} | $ for $r \gtr |z|$ and $  | I_{rr}^{\nne} | $ separately.

\vs{0.1in}
\paragraph{\bf Estimate of $ I_{rr}^{\nne}$}

Recall $D_{\pm}$ from \eqref{eq:Dpm:recall1}. For $I_{rr}^{\nne}$, using $\tr , r \leq D_+$, 
and then Lemma \ref{lem:CS_singu} with $k=1, m = 1 - \al \in [0.66, 0.68], \ell = \hal + \kag \in [0.32,0.34]$ and $n=0$, we bound 
\[
\bal
  |I_{rr}^{\nne}| & \les P(\xx) 
\min(|\xx| , 1) \nnrr{\ww} \cdot \f{1}{r}   \int_{ |(\tr, \tz) | \leq 4 |\xx|}   \f{  1 }{D_- D_+^{1-\al} }     \ang \tz^{-\hal -\kag} d \tr d \tz \\
& \les P(\xx) 
\min(|\xx| , 1) \nnrr{\ww} \cdot \f{1}{r}  |\xx|^{\al} \ang \xx^{-\hal - \kag} \, .
\eal
\]
Recall $P(\xx) = \ang \xx^{\kag} \rhoo^{-1}(\xx)$ from \eqref{eq:ker_wwlinf_ne}. We prove 
\[
  r    |I_{rr}^{\nne}|
\les  \min(|\xx|, 1) \rhoo^{-1}(\xx) |\xx|^{\al} \ang \xx^{-\hal} \nnrr{\ww}
   \les \rhoo^{-1} \min( |\xx|^{\al + 1}, \ang \xx^{\al-\hal} ) \nnrr{\ww}.
\]
Since $r \leq |\xx|$, combining the above estimate and \eqref{eq:psi_Irz:far},  we prove \eqref{eq:vel_rr}.

\vs{0.1in}
\paragraph{\bf Estimate of $\psi_{rz}$ for $r \gtr |z|$}
Since $r \gtr |z|$, we obtain $r \gtr |\xx|$.  Using $D_+ \geq r +\tr \gtr |\xx| + \tr$, 
$|\yy| \les |\xx|$, and then Lemma \ref{lem:CS_singu} with $k = 1, m=0, n=0,\ell = \hal + \kag \in [0.32, 0.34]$, we bound 
\[
\bal
      |I_{rz}^{\nne} | 
     & \les \min(|\xx|, 1) P(\xx) \nnrr{\ww} |\xx|^{\al-2} 
      \int_{ |(\tr, \tz)| \leq 4 |\xx| } \f{1}{D_-} \ang \tz^{-\hal - \kag} d \tr d \tz \\
   &  \les \min(|\xx|, 1) P(\xx) \nnrr{\ww} |\xx|^{\al-2}  
   \cdot  |\xx|  \cdot \ang \xx^{- \hal -\kag}.
  \eal
\]
Recall $P(\xx) = \ang \xx^{\kag} \rhoo^{-1}(\xx)$ from \eqref{eq:ker_wwlinf_ne}. We prove 
\[
   |I_{rz}^{\nne} | 
   \les \min(|\xx|, 1) \rhoo^{-1}(\xx) |\xx|^{\al-1} \ang \xx^{-\hal} \nnrr{\ww}
   \les \rhoo^{-1} \min( |\xx|^{\al}, \ang \xx^{\al-\hal-1} ) \nnrr{\ww}.
\]
Combining the above estimate and \eqref{eq:psi_Irz:far}, we prove 
\eqref{eq:vel_rz} for $\nu \leq \f{1}{2}$ in the case$r \geq z$.

\subsubsection{Estimate of $\psi_{rz}$ near the singularity: case $r \leq z$}\label{sec:psi_rz:z_large}

In the case $r \leq z$, we estimate $\psi_{rz}$ using 
 a different decomposition. Using \eqref{eq:ker_W_id}, we first decompose $I_{rz}^{\nne}$ in \eqref{id:psi_rz} as follows 
\begin{align}\label{id:psi_rz2}
  I_{rz}^{\nne} & = 
   P.V.  \int_{|\yy | \leq 4 |\xx|}  \pa_z \KRF(\xx-\yy) \cdot 
  y_1  \tr^{\al-2} \pa_\tr \tw 
  +   P.V.  \int_{|\yy | \leq 4 |\xx|}  \pa_z \KRF(\xx-\yy) \cdot 
   (\al-1)  y_1 \tr^{\al-3}  \ww  \notag \\
   & \teq I_{rz}^{\nne}(\pa_\tr \ww) + I_{rz}^{\nne}(  \ww ).
\end{align}

\paragraph{\bf Estimate of $I_{rz}^{\nne}( \pa_{\tr} \ww)$}

We first estimate $I_{rz}^{\nne}(\pa_\tr \ww) $. 
Since $\pa_{\tr} \ww$ only depends on $\tr, \tz$ 
and $|\yy| = |\tr, \tz|$, using $|\na \KRF(s) |\les |s|^{-4}$,
$|y_1| \leq \tr$,  
and then \eqref{eq:int_singu} with $k =4$,  we bound  $I_{rz}^{\nne}(\pa_\tr \ww)$ in \eqref{id:psi_rz} as 
\[
\bal
  |I_{rz}^{\nne}(\pa_\tr \ww) | 
 &  \les  \int_{|\yy| \leq 4 |\xx|} |\xx-\yy|^{-4} 
 |\tr^{\al-1} \pa_\tr \ww |  d \yy   \les 
       \int_{ |(\tr, \tz)| \leq 4 |\xx| } \f{ |\tr^{\al-1 }  \tr^3  \pa_\tr \tw| }{D_- D_+^3} d \tr d \tz  \, .
\eal
\]

Applying $\tr \leq D_+$ and  the bound of $\pa_\tr \ww$ in \eqref{eq:ker_wwlinf_ne}, we obtain
\[
    |I_{rz}^{\nne}(\pa_\tr \ww) | 
    \les 
     \int_{ |(\tr, \tz)| \leq 4 |\xx| } \f{ |  \pa_\tr \tw| }{D_- D_+^{1-\al}} d \tr d \tz 
  \les P(\xx) \nnrr{\ww} \int_{ |(\tr, \tz)| \leq 4 |\xx| } 
  \f{ 1  }{D_- D_+^{1-\al}} \ang {\tr, \tz}^{-1} \ang \tz^{-\hal -\kag} d \tr d \tz .
\]

Applying Lemma \ref{lem:CS_singu} with $k = 1, m=1-\al, n =1, \ell = \hal + \kag \in [0.32, 0.34]$, we bound 
\[
\bal
      |I_{rz}^{\nne}(\pa_\tr \ww) | 
   & \les   P(\xx) \nnrr{\ww} 
  |\xx|^{2 - (2 -\al)} \ang \xx^{-\hal -\kag-1}
\les  P(\xx) \nnrr{\ww}  |\xx|^{\al} \ang \xx^{-\hal -\kag-1}.
  \eal
\]
Recall $P(\xx) = \ang \xx^{\kag} \rhoo^{-1}(\xx)$ from \eqref{eq:ker_wwlinf_ne}. We bound 
\beq\label{eq:psi_rz:Iwr}
        |I_{rz}^{\nne}(\pa_\tr \ww) | 
        \les \rhoo^{-1}(\xx) |\xx|^{\al} \ang \xx^{-\hal}\nnrr{\ww}
      \les \rhoo^{-1}(\xx) \min (|\xx|^{\al}, \ang \xx^{\al-\hal}) \nnrr{\ww}.
\eeq

\vspace{0.1in}
\paragraph{\bf Estimate of $I_{rz}^{\nne}(\ww)$}

We fix $\nu \in [0, 1/2]$ from Proposition \ref{prop:vel_est}
and recall $D_-(r, z -\tz, \tr)= |(\tr- r, z-\tz)|$ from \eqref{def:Dpm}. 
We simplify it as $D_-$. We expand the principle integral in $I_{rz}^{\nne}$ from \eqref{id:psi_rz2}
and perform the decomposition
\begin{align}\label{eq:psi_rz:Iw1}
I_{rz}^{\nne}(  \ww ) 
& = \lim_{a\to 0} \int_{ \ssk{a\leq |\yy- \xx| , \\ |\yy | \leq 4 |\xx| } } 
 \pa_z \KRF(\xx-\yy) \cdot (\al-1)  y_1 \tr^{\al-3}  \ww   d \yy \\
& = \lim_{a\to 0}  \B( \int_{ \ssk{a \leq |\yy- \xx| , D_- \leq \nu z, \\ |\yy | \leq 4 |\xx| } } 
+ \int_{ \ssk{ D_- > \nu z, \\ |\yy | \leq 4 |\xx| } } \B)
 \pa_z \KRF(\xx-\yy) \cdot (\al-1)  y_1 \tr^{\al-3}  \ww   d \yy
 \teq  \lim_{a\to 0}  I_{a ,\nu} + I_{\geq \nu}. \notag
\end{align}

For $I_{a, \nu}$, the domain of the integral in $\tz$ is symmetric with respect to $z$.  Since $\pa_z \KRF(\xx-\yy) = C \f{z - \tz}{ |\xx-\yy|^5 }$ 
is odd in $z -\tz$ by \eqref{eq:psi_W}, we  symmetrize the integrand in $\tz$ as 
\[
  I_{a,\nu}  = \f{1}{2} \int_{ \ssk{a \leq |\yy- \xx| , \, D_-\leq \nu z, \, |\yy | \leq 4 |\xx| } } 
 \pa_z \KRF(\xx-\yy) \cdot (\al-1)  y_1 \tr^{\al-3}  ( \ww(\tr , \tz) 
 - \ww(\tr, 2 z - \tz) d \yy .
\]

In the domain, since $|z -\tz| \leq \nu z \leq z/2$, %we obtain $\tz \in [z /2 , 2 z ]$ and $|\tz| \asymp |\xx|$. 
applying bound of $\pa_\tz \ww $ from \eqref{eq:ker_wwlinf_ne}, we bound 
\[
  | \ww(\tr , \tz) 
 - \ww(\tr, 2 z - \tz)|
 \les  |z-\tz| \sup\nolimits_{ |\xi- z| \leq z/2 } |\pa_\tz \ww(\tr, \xi)|
 \les |z-\tz| P(\xx)  \ang z^{-\hal -\kag - 1} \nnz{\ww} .
\]

Since $|z| \geq r$, we get $|z| \asymp |\xx|$. Applying \eqref{eq:K5D_est:dK} on $\KRF$, the above bound, $|y_1| \leq \tr$, and $|z| \asymp |\xx|$, we estimate
 \beq\label{eq:psi_rz:Iw_est1}
   I_{a, \nu}
  \les 
  P(\xx) \ang \xx^{-\hal -\kag - 1} \nnz{\ww} \udb{ \int_{ \ssk{  D_- \leq \nu z } } 
  |\xx- \yy|^{-4} \tr^{\al-2} \cdot |z -\tz| d \yy }_{:=M_1} .
 \eeq

Since the domain $\{ D_-\leq \nu z \}$ 
and $\tr^{\al-2} |z-\tz|$ only depend on $\tr, \tz$, applying \eqref{eq:int_singu} with $k=4$,
and then use $ D_+  \gtr \tr , |z-\tz | , D_- $,
 we bound 
\[
  M_1 \les  \int_{ D_- \leq \nu z } \f{\tr^{\al+1} |z -\tz|  }{ D_- D_+^3 } d \tr d \tz 
\les \int_{ D_- \leq \nu z } \f{ 1  }{ D_-^{2-\al} } d \tr d \tz 
\les  \int_{ |z-\tz | \leq \nu z } \f{ 1  }{ D_-^{2-\al} } d \tr d \tz  \, .
\]

Recall $D_{\pm} = |(\tr \pm r, z -\tz)|$ \eqref{def:Dpm}.  
Integrating in $r$ and then in $z$, and using $|z| \leq |\xx|$, we bound 
\[
  M_1 \les  |\nu z|^{\al} \les \nu^{\al }|\xx|^{\al} .
\]

Recall $P(\xx) = \ang \xx^{\kag} \rhoo^{-1}(\xx)$ from \eqref{eq:ker_wwlinf_ne}. Applying the above estimate to \eqref{eq:psi_rz:Iw_est1}, we bound 
\beq\label{eq:psi_rz:Iw_est2}
  |   I_{a, \nu}| \les \nu^{\al} |\xx|^{\al} \ang \xx^{\kag} \rhoo^{-1}(\xx)
  \ang \xx^{-\kag - \hal - 1} \nnz{\ww}
  \les \nu^{\al} \rhoo^{-1}(\xx) \min( |\xx|^{\al}, \ang \xx^{ \al- \hal - 1} )
  \nnz{\ww}  \, ,
\eeq
uniformly in $a$.

\vspace{0.1in}

\paragraph{\bf Estimate of $I_{ \geq \nu }$}

Recall $I_{\geq \nu}$ from \eqref{eq:psi_rz:Iw1}. 
Since $\ww(\tr, \tz)$ and  the domain $ \{ D_- \geq \nu z, |(\tr, \tz)| \leq 4 |\xx| \}$
only depend on $\tr, \tz$, using \eqref{eq:K5D_est:dK},
$|y_1| \leq |\tr|$, and then \eqref{eq:int_singu} with $k=4$,
we bound 
\[
I_{\geq \nu}
\les
 \int_{ \ssk{ D_- > \nu z, \\ |\yy | \leq 4 |\xx| } } 
  \f{1}{|\xx-\yy|^4}  \tr^{\al-2} |\ww(\yy)| d \yy
  \les 
  \int_{ \ssk{ D_- > \nu z, \\ |(\tr, \tz) | \leq 4 |\xx| } } \f{  \tr^{\al- 2} \tr^3 
|\ww(\tr, \tz)|
  }{ D_- D_+^3 }  d \tr d \tz .
\]

Using $D_+ \geq \tr, D_-$,
and the bound of $\ww$ in \eqref{eq:ker_wwlinf_ne}, we bound 
\bseq\label{eq:psi_rz:I_geq}
\beq\label{eq:psi_rz:I_geq:a}
  I_{\geq \nu}
  \les  \int_{ \ssk{ D_- > \nu z, \\ |(\tr, \tz) | \leq 4 |\xx| } } \f{   
|\ww(\tr, \tz)|
  }{ D_- \cdot D_-^{2-\al}  }  d \tr d \tz
  \les  P(\xx) \min(|\xx|,1) \nnrr{\ww} \udb{ \int_{ \ssk{ D_- > \nu z, \\ |(\tr, \tz) | \leq 4 |\xx| } } \f{1    }{ D_-^{3-\al} }  \ang \tz^{-\hal -\kag} d \tr d \tz}_{:=M_2}.
\eeq

Recall $D_- = |(\tr - r, z -\tz)|$ and $z \geq r, z \asymp |\xx|$. Since $\max(D_-, |\tz|)\gtr |z-\tz| + |\tz| \gtr |\xx|$, we obtain
\[
   \f{ 1    }{ D_-^{3-\al} }  \ang \tz^{-\hal -\kag}
   \les |\xx|^{\al-3} \ang \tz^{-\hal -\kag}
   + \f{1}{D_-^{3-\al}} \ang \xx^{-\hal -\kag}.
\]

For the integral of the first term, we apply Lemma \ref{lem:CS_singu}
with $k =0, \ell = \hal + \kag \in[0.32, 0.34], m=n=0$. For the second term, we integrate it directly. We then yield 
\beq
\bal
  M_2 
  & \les \int_{|\tr ,\tz| \leq 4 |\xx|} |\xx|^{\al-3} \ang \tz^{-\hal - \kag}
  + \ang \xx^{-\hal -\kag} \int_{ |(\tr, \tz) - (r, z)| \geq \nu z } \f{1}{ |(\tr, \tz) - (r, z)|^{3-\al} } d \tr d \tz  \\
& \les |\xx|^{\al-3 + 2} \ang \xx^{-\hal -\kag} + (\nu z)^{\al-1}  \ang \xx^{-\hal -\kag}
\les \nu^{\al-1} |\xx|^{\al-1} \ang \xx^{-\hal -\kag}.
\eal
\eeq
where in the last inequality, we have used $|z| \asymp |\xx|, \nu \leq 1/2$.

Recall $P(\xx) = \ang \xx^{\kag} \rhoo^{-1}(\xx)$ from \eqref{eq:ker_wwlinf_ne}. We prove 
\beq
  I_{\geq \nu} \les P(\xx) \min(|\xx|,1) \nnrr{\ww} \nu^{\al-1} |\xx|^{\al-1} \ang \xx^{-\hal -\kag}
  \les \rhoo^{-1} \nu^{\al-1} \nnrr{\ww}
  \min( |\xx|^{\al}, \ang \xx^{\al-\hal -1} ).
\eeq
\eseq

Combining the above estimate, \eqref{eq:psi_rz:Iw_est2} on $I_{a, \nu}$,
 \eqref{eq:psi_rz:Iwr} on $I_{rz}^{\nne}(\pa_\tr \ww)$, and decomposition
\eqref{eq:psi_rz:Iw1} on $I_{rz}^{\nne}$, we prove 
\[
\bal
  |I_{rz}^{\nne}| & \leq   |I_{rz}^{\nne}( \pa_{\tr} \ww )|
  +  |I_{rz}^{\nne}(  \ww )|  \leq |I_{rz}^{\nne}( \pa_{\tr} \ww )|
  +  \lim_{a\to 0} |I_{a,\nu}| + I_{\geq \nu} \\
 &  \les \rhoo^{-1}( \nu^{\al-1} \nnrr{\ww} + \nu^{\al} \nnz{\ww} )
  \min( |\xx|^{\al}, \ang \xx^{\al-\hal -1} ).
\eal
\]
in the case $z \geq r$. Combining the above estimate and \eqref{eq:psi_Irz:far}, 
we prove \eqref{eq:vel_rz} in the case $z \geq r$.

\subsection{Estimate of \texps{$\psi_{zz}$}{psi,zz} }\label{sec:psi_zz}

Recall $\psi$ from \eqref{eq:5D_BS:rec}. 
We consider $z \neq 0$. The estimate for $z =0$ is obtained using continuity, as discussed in 
Section \ref{sec:psi_rz}.
Using the convolution in $z$ and then performing integration by parts in $D_- = |(\tr - r, z -\tz)| \geq \nu z$, we obtain 
\begin{align}\label{id:psi_zz}
  \pa_{zz} \psi & =   \int  (\pa_{z} \KRF)( \xx -\yy )  \cdot \pa_{\tz} W(\yy)  d \yy    
 = 
 (   \int_{ D_- \leq \nu z} +    \int_{ D_- > \nu z} ) (\pa_{z} \KRF)( \xx -\yy )  \cdot \pa_{\tz} W(\yy)  d \yy  \notag \\
 & =   \int_{ D_- \leq \nu z} 
 (\pa_{z} \KRF)( \xx -\yy )  \cdot \pa_{\tz} W(\yy)  d \yy 
+ \int_{  D_- = \nu z  } \nn_{\tz} \pa_z \KRF(\xx-\yy) \cdot W(\yy) \notag \\
& \quad 
 + \int_{ D_- \geq \nu z , \, |\yy| \leq 4 |\xx|  } \pa_{zz} 
 \KRF (\xx- \yy) \cdot  W(\yy) 
 + \int_{   |\yy| \geq 4 |\xx|  } \pa_{zz}  \KRF (\xx- \yy) \cdot  W(\yy) 
  \notag \\
&  \teq I_{zz}(\pa_\tz \ww)  + I_{zz, \nu } +I_{zz, \geq \nu}
+ I_{zz}^{\ffar} \, ,
\end{align}
where $\nu\leq 1/2$ and  $\nn_{\tz}$ denotes the unit outer normal in the $\tz$ direction. 
Since $ |\yy| \geq 4 |\xx| $ implies $D_-= |(\tr, \tz) - (r, z)| \geq |\yy|-   |\xx| \geq |\xx|\geq \nu z$, 
we have simplified the domain in $I_{zz}^{\ffar}$.

\subsubsection{Estimate far-field integral: $I_{zz}^{\ffar}$}

Recall $z = x_5$ from \eqref{def:5D_x}. We introduce the symmetrized kernel in $z$:
\[
   K_{\R^5, zz}^{\mw{sym}}(\xx, \yy) \teq \pa_{zz}  \KRF (\xx- \yy) - \pa_{zz} \KRF( \xx_h - \yy_h, x_5 + y_5  ).
\]

Since $\pa_{5}^2 \KRF(s)$ is even in $s_5$, 
for $|\yy| \geq 4 |\xx|$, using mean value theorem and \eqref{eq:K5D_est:dK}, we bound 
\[
\bal
 |    K_{\R^5, zz}^{\mw{sym}}(\xx, \yy)  |  & = |\pa_5^2 \KRF(  \xx_h -\yy_h,  x_5 - y_5 )
 - \pa_5^2 \KRF(  \xx_h -\yy_h,  -x_5 - y_5 )| \\
 & \les |\xx|  \sup\nolim_{ |\xi |\leq |\xx| }| \na^3 \KRF(\xi-\yy) |
 \les |\xx| \cdot |\yy|^{-6} .
 \eal
\]

Since $W(\yy)$ is odd in $\tz$, symmetrizing the integral in $I_{zz}^{\ffar}$ \eqref{id:psi_zz} and using the above estimate and bound of $W$ from \eqref{eq:ker_Wlinf}, we bound 
\[
\bal
   |I_{zz}^{\ffar} |
   & = 
   \f{1}{2} \B| \int_{|\yy| \geq 4 |\xx|}  K_{\R^5, zz}^{\mw{sym}}(\xx,\yy)   \, W(\yy) d \yy \B| 
    \les |\xx| \int_{|\yy| \geq 4 |\xx|} %|\yy|^{-6} 
   \f{ \min( |\yy|, 1 )}{|\yy|^6} \tr^{\al-1} \ang \tz^{-\hal} \rhoo^{-1}(\yy) \rag^{-1}(\yy) d \yy
   \nnrr{\ww} .
   \eal
\]
Since $|\yy|^{-6} \tr^{\al-1} \les |\yy|^{-5} \tr^{\al-2}$, applying 
Lemma \ref{lem:int_W} with $\ell=5$, we bound 
\beq\label{eq:psi_zz:Ifar}
\bal
     |I_{zz}^{\ffar} |
     & \les |\xx| \rhoo^{-1}(\xx) \min( |\xx|^{\al-1} + |\xx|^{-\al}, |\xx|^{\al-\hal -2} ) \nnrr{\ww} \\
     & \les \rhoo^{-1}(\xx) \min( |\xx|^{\al} , \ang \xx^{\al- \hal -1} ) \nnrr{\ww}. 
    \eal
\eeq

\subsubsection{Estimate near the singularity: $I_{zz}^{\nne}$}

Below, we estimate $I_{zz}^{\nne}( \pa_{\tz} \ww ), I_{zz, \nu}$ and $I_{zz,\geq \nu}$ in \eqref{id:psi_zz}.

\vspace{0.1in}
\paragraph{\bf Estimate of $I_{zz}^{\nne}( \pa_{\tz} \ww )$} 

Recall $I_{zz}^{\nne}( \pa_{\tz} \ww )$ from \eqref{id:psi_zz}. Since $\pa_{\tz} W(\yy)$ and the domain $\{ D_- \leq \nu z \}$ only depend on $\tr, \tz$, applying estimate of $\KRF$ \eqref{eq:K5D_est:dK}, and \eqref{eq:int_singu} with $k=4$, we bound 
\[
  |I_{zz}(\pa_\tz \ww)| 
  = \B| \int_{ D_- \leq \nu z} 
 (\pa_{z} \KRF)( \xx -\yy )  \cdot \pa_{\tz} W(\yy)  d \yy \B|
 \les \int_{ D_- \leq \nu z} \f{|\pa_{\tz} W(\yy)|}{|\xx-\yy|^4}   d \yy 
 \les \int_{ D_- \leq \nu z} \f{\tr^3  |\pa_{\tz} W(\yy)|   }{D_- D_+^3}  d \yy.
\]

Recall $D_{\pm} = |(r\pm \tr, z -\tz)|$ \eqref{eq:Dpm:recall1}. Since $\nu |z| \leq |\xx|$, in the domain of the integral, we obtain 
$|(\tr, \tz)| \leq D_- + |\xx| \leq 4 |\xx| $. 
Moreover, since $\nu \leq 1/2$, $|z-\tz| \leq D_- \leq \nu z$, 
and $z \neq 0$,  we obtain $\tz \asymp |z|$. Therefore, using the bound of $\pa_\tz W$ from \eqref{eq:ker_W_id} and \eqref{eq:ker_wwlinf_ne}, we obtain
\[
  \f{ \tr^3 |\pa_\tz W(\tr, \tz)|}{D_- D_+^3}
  =  \f{ \tr^3 \tr^{\al-1} | \pa_\tz \ww |  }
  {D_- D_+^3}
  \les P(\xx)  \f{\tr^{\al+2}  \ang \tz^{ -\hal -\kag-1 } }{D_- D_+^3} \nnz{\ww}
  \les P(\xx) \nnz{\ww} \cdot \f{\tr^{\al+2}}{D_- D_+^3 } \ang z^{-\hal-\kag-1}  \, .
\]

Since $D_+ + |z| \gtr r + |z| \gtr |\xx|$, we have $D_+ \gtr |\xx|$ or $z \gtr |\xx|$. Thus, 
using $\tr \leq |\xx|$ and $\tr, D_- \leq D_+$, for any $ \ell \geq 0$, we estimate (the decay factor is in $\ang z$ rather than  the ``tilde" $\ang \tz$)
\beq\label{eq:psi_zz:int_est1}
\bal
  \f{\tr^{\al+2}}{D_- D_+^3 } \ang z^{-\hal-\kag- \ell} 
&  \les \one_{D_+ \gtr |\xx|} \f{\tr^{\al+2}}{ D_- |\xx|^3 }  
  \ang z^{-\hal-\kag- \ell}  + \one_{|z| \gtr |\xx|} 
  \ang \xx^{-\hal-\kag - \ell} \f{\tr^{\al+2}}{ D_- D_+^3 } \\
  & \les |\xx|^{\al-1} \f{1}{D_-}  \ang z^{-\hal-\kag- \ell}  
  +    \ang \xx^{-\hal-\kag - \ell} \f{ 1 }{ D_-^{2-\al} } \one_{z \gtr |\xx|} .
\eal
\eeq

Integrating the above bound with $\ell=1$ in $ \{ D_- \leq \nu z \}$, we establish
\beq\label{eq:psi_zz:est1}
\bal
    |I_{zz}(\pa_\tz \ww)| 
 &  \les P(\xx) \nnz{\ww}  \int_{ D_-\leq \nu z } |\xx|^{\al-1} \f{1}{D_-}  \ang z^{-\hal-\kag-1}  
  +     \ang \xx^{-\hal-\kag - 1} \f{ 1 }{ D_-^{2-\al} } d \tr d \tz \\
  & \les P(\xx) \nnz{\ww} \big( |\xx|^{\al-1} \nu |z|  \ang z^{-\hal-\kag-1}   
  + \ang \xx^{-\hal -\kag - 1} \nu^{\al} |z|^{\al}  \big) .
\eal
 \eeq

 Since $\nu \leq 1/2$, using $|z| \leq |\xx|$ and by considering $|z|\leq 1$ and $|z| > 1$,
we further obtain
 \beq\label{eq:psi_zz:Izz_wz}
       |I_{zz}(\pa_\tz \ww)| 
       \les P(\xx) \nnz{\ww} \nu^{\al} \min (|\xx|^{\al}, \ang \xx^{\al-1} \ang z^{-\kag - \hal}) .
 \eeq

\vspace{0.1in}

\paragraph{\bf Estimate of $I_{zz, \nu}$}
Recall $I_{zz, \nu}^{\nne}$ from \eqref{id:psi_zz}
\[
  I_{zz, \nu} = \int_{  D_- = \nu z  } \nn_{\tz} \pa_z \KRF(\xx-\yy) \cdot W(\yy)  \, .
\]

Since $W(\yy)$ and the domain $\{ D_-= \nu z \}$ only depend on $\tr, \tz$, using bound 
of $\KRF$ \eqref{eq:K5D_est:dK} and \eqref{eq:int_singu} with $k=4$, 
and bound of $W$ from \eqref{eq:ker_W_id} and \eqref{eq:ker_wwlinf_ne}, 
we estimate 
\[
\bal
  |I_{zz, \nu}|
 &  \les \int_{  D_- = \nu z  } \f{  |W(\yy) |  }{|\xx-\yy^4|} 
  \les \int_{D_- = \nu z} \f{ |W(\yy) \tr^{3}|  }{D_- D_+^3} d \tr d \tz \\
  & \les  P(\xx) \min(|\xx|,1) \nnrr{\ww}  
   \int_{D_- = \nu z} \f{ |\tr^{\al-1} \tr^{3}|  }{D_- D_+^3} \ang \tz^{-\hal -\kag} d \tr d \tz \, .
   \eal
\] 

For $D_- = \nu z$, we obtain $|\tz - z| \leq z/2$, which 
implies 
  $|\tz| \asymp |z|$. 
 Applying estimate of $W = \tr^{\al-1} \ww$  from  \eqref{eq:ker_W_id} and \eqref{eq:ker_wwlinf_ne}, 
 $|\tz| \asymp |z|$ and \eqref{eq:psi_zz:int_est1} with $\ell=0$, we bound 
 \[
 \bal
 \f{ \tr^{\al+2} }{D_- D_+^3} \ang \tz^{-\hal - \kag } 
 & \les   \f{ \tr^{\al+2} }{D_- D_+^3} \ang z^{-\hal - \kag } 
 \les  |\xx|^{\al-1} \f{1}{D_-}  \ang z^{-\hal-\kag }  
  +    \ang \xx^{-\hal-\kag } \f{ 1 }{ D_-^{2-\al} } \one_{|z|\gtr |\xx|} .
  \eal
 \]

Since the surface measure  of domain $\{ D_- = \nu z \} $ 
in $(r,z)$-coordinate is $C \nu z$, using the above estimate, we bound
\begin{align}\label{eq:psi_zz:est2_nu}
  |I_{zz, \nu}|
& \les   P(\xx) \min(|\xx|,1) \nnrr{\ww}  \cdot \nu |z| 
\cdot \big( |\xx|^{\al-1} \f{1}{ \nu |z| }  \ang z^{-\hal-\kag }  
  +    \ang \xx^{-\hal-\kag } \f{ 1 }{ (\nu z)^{2-\al} } \one_{|z|\gtr |\xx|} \big) \notag \\
& \les P(\xx) \min(|\xx|,1) \nnrr{\ww} 
( |\xx|^{\al-1} \ang z^{-\hal - \kag} 
+  \ang \xx^{-\hal - \kag}  \nu^{\al-1} |\xx|^{\al-1} ) \notag \\
& \les 
\nu^{\al-1} P(\xx)  \nnrr{\ww} 
\min( |\xx|^{\al}, \ang \xx^{\al-1} \ang z^{-\hal - \kag}  ).
\end{align}

\vspace{0.1in}
\paragraph{\bf Estimate of $I_{zz, \geq \nu}$}
Recall $I_{zz, \geq \nu}$ from \eqref{id:psi_zz}
\[
 I_{zz, \geq \nu} =  \int_{ D_- \geq \nu z , |\yy| \leq 4 |\xx|  } \pa_{zz} 
 \KRF (\xx- \yy) \cdot  W(\yy) 
\]

Note that $|\yy| =| (\tr, \tz)|$. Since the domain of the integral 
and $W$ do not depend on $W$, applying bound of $\KRF$ \eqref{eq:K5D_est:dK},
 \eqref{eq:int_singu} with $k =5$, and $\tr \leq D_+$, we obtain
\beq\label{eq:psi_zz:I_lg_est1}
\bal
    |I_{zz, \geq \nu}|
  & \les \int_{ \ssk{ D_- \geq \nu z , \\ |\yy| \leq 4 |\xx|  }  } \f{  |W(\yy) |  }{|\xx-\yy|^5} d \yy
   \les \int_{ \ssk{ D_- \geq \nu z , \\ |(\tr, \tz)| \leq 4 |\xx| }  } \f{  \tr^{\al-1 + 3}  |\ww |  }
   {D_-^2 D_+^3} d \tr d \tz 
\les \int_{ \ssk{ D_- \geq \nu z , \\ |(\tr, \tz)| \leq 4 |\xx| }  } \f{    |\ww |  }
   {D_-^2 D_+^{1-\al}} d \tr d \tz .
  \eal
\eeq

\paragraph{\bf Case: $|z| \geq r$}

When $|z| \geq r$, since $  D_-\leq D_+$, we obtain 
\[
    |I_{zz, \geq \nu}| \les 
    \int_{  D_- \geq \nu z ,  |(\tr, \tz)| \leq 4 |\xx| }   \f{   |\ww |  }{D_-^{3-\al}  }.
    d \tr d \tz  \, .
\]
The above upper bound is the same as that of $I_{\geq \nu}$ in \eqref{eq:psi_rz:I_geq:a}
where we assume $|z| \geq r$ in Section \eqref{sec:psi_rz:z_large}.
Thus using estimate of $I_{\geq \nu}$ in \eqref{eq:psi_rz:I_geq}, we bound 
\beq\label{eq:psi_zz:I_lg_ez}
      |I_{zz, \geq \nu}| \les  \rhoo^{-1} \nu^{\al-1} \nnrr{\ww}
  \min( |\xx|^{\al}, \ang \xx^{\al-\hal -1} ).
\eeq

\paragraph{\bf Case: $r  \geq |z|$}

Below, we consider $|(\tr, \tz)| \leq 4 |\xx|$. For $r \geq |z|$, we obtain $r \gtr |\xx|$ and $D_+ \gtr |\xx|$. 
To bound $I_{zz, \geq \nu}$, since $\ww$ is odd in $z$, %we  use vanishing condition of $\ww$ near $z = 0$. Since $\ww(\tr,0) = 0$, 
using $\ww(\tr,0) = 0$, \eqref{eq:ker_wwlinf_ne}, for some $ |\xi| \leq |\tz|$, we obtain 
\[
  |\ww(\tr, \tz)| \les |\tz| \cdot |\pa_\tz \ww(\tr, \xi)|  \les P(\xx) |\tz| \nnz{\ww}. 
\]
Minimizing the above bound and the bound in \eqref{eq:ker_wwlinf_ne}, we yield 
\beq\label{eq:psi_zz:I_lg_est2}
    |\ww(\tr, \tz)|  \les P(\xx) \cdot B(\xx, \tz), \quad B(\xx, \tz) \teq \min\B( \min( |\xx|, 1 ) \ang \tz^{-\hal-\kag} \nnrr{\ww} , \, |\tz|
    \cdot  \nnz{\ww} \B) .%\teq P(\xx) \cdot B(\xx, \tz).
\eeq

Using the above bound and $D_+ \gtr |\xx|$, we estimate $  I_{ zz, \geq \nu }$ in \eqref{eq:psi_zz:I_lg_est1} as 
\beq\label{eq:psi_zz:I_lg_est3}
    I_{ zz, \geq \nu } 
    \les |\xx|^{\al-1}  \int_{  D_- \geq \nu z ,  |(\tr, \tz)| \leq 4 |\xx|   } \f{   |\ww|  }
   {D_-^2  } d \tr d \tz 
   \les P(\xx) |\xx|^{\al-1} \udb{ \int_{  D_- \geq \nu z , |\tz| \leq 4 |\xx|   } \f{  B(\xx, \tz) }
   {D_-^2  } d \tr d \tz }_{\teq M_3} \, .
\eeq

Recall $D_- \asymp |z-\tz| + |r -\tr|$ \eqref{eq:Dpm:recall1}. Denote $s = r -\tr$. Integrating $M_3$ over $\tr$, we obtain
\[
\bal
  M_3  & \les  
   \int_{|\tz|\leq 4 |\xx|} \B( \int_{|z -\tz| + |s| \geq \nu z} \f{  B(\xx,\tz)  }{ (|z -\tz| + |s|)^2  } d s  \B)
   d \tz   \les   \int_{|\tz|\leq 4 |\xx|}  \f{  B(\xx,\tz)   }{ \max( |z-\tz| , \nu z ) }    d \tz 
 \les  \int_{|\tz|\leq 4 |\xx|}  \f{ B(\xx,\tz)  }{  |z-\tz| + \nu z  }      d \tz  \\
 &  =(  \int_{|\tz|\leq |z|/2 } 
+ \int_{ |z|/2 \leq |\tz| \leq 2 |z|}
+ \int_{ 2|z| <  |\tz | < 4 |\xx| } )  \f{ B(\xx,\tz)  }{  |z-\tz| + \nu z  }      d \tz  
\teq M_{31} + M_{32} + M_{33}. 
\eal
\]

In the domain of $M_{31}$, we have $|z-\tz| \asymp |z|$. 
For $M_{31}$ and $M_{32}$, we apply the  first bound in \eqref{eq:psi_zz:I_lg_est2} for $B$. 
Using a direct calculation and a change of variable $\tz = s |z|$, we obtain 
\[
\bal
\int_{|\tz| \leq |z|/2} \f{\ang \tz^{-\kag -\hal}}{ |z-\tz| + \nu z  }  d \tz 
& \les |z|^{-1} \int_{|\tz| \leq |z|/2} \ang \tz^{-\kag -\hal}  d \tz 
\les |z|^{-1} \min( |z|, |z|^{1-\kag -\hal } ) \les \ang z^{-\kag - \hal}, \\
\int_{ |z|/2 \leq |\tz| \leq 2 |z|}\f{\ang \tz^{-\kag -\hal}}{ |z-\tz| + \nu z  }  d \tz 
& \les  \ang z^{-\kag -\hal}  
\int_{ 1/2 \leq |s| \leq 2 }  \f{ 1 }{ |s-1| + \nu } d s
\les \ang z^{-\kag -\hal}  \log \nu^{-1}.
\eal
\]

Since $\nu \leq 1/2$, using the above estimates and the first bound in \eqref{eq:psi_zz:I_lg_est2} for $B$, we obtain
\[
  M_{31} + M_{32} \les \min(|\xx|, 1)  \nnrr{\ww} 
  \ang z^{-\kag -\hal}  \log \nu^{-1}.
\]

In the domain of $M_{33}$, we have $|z-\tz| \gtr |\tz|$ and obtain
\[
    M_{33} 
    \les \int_{ 2 |z| \leq |\tz| \leq 4 |\xx|} \f{B(\xx, \tz)}{ |\tz| } d \tz
    \les \int_{ 2 |z| \leq |\tz| \leq 4 |\xx|} \f{B(\xx, \tz)}{ |\tz| }  
\one_{|\tz| \leq \nu \min(1, |\xx|)} 
+ \one_{|\tz| > \nu \min(1, |\xx|)} 
    d \tz \teq M_{33,1} + M_{33,2}.
\]

If $M_{33,1} \neq 0$, we have $|z| \les 1$ and $\ang z^{-\kag -\hal} \asymp 1$.  We apply the second bound of $B$ \eqref{eq:psi_zz:I_lg_est2} to obtain
\[
  M_{33,1} \les \ang z^{-\kag -\hal} \int_0^{\nu \min(1, |\xx|)} 
  \f{ |\tz| \cdot \nnz{\ww} }{|\tz|} d \tz 
  \les  \nu \ang z^{-\kag -\hal}  \min(  |\xx|, 1  ) \nnz{\ww}. 
\]

For the second part $M_{33,2}$, we use the first bound of $B$ \eqref{eq:psi_zz:I_lg_est2}.
If $|z| > 1$, we obtain
\[
  M_{33,2} \les \min(|\xx|, 1) \nnrr{\ww}  \int_{2 |z|}^{\infty} \f{1}{ |\tz| } \ang \tz^{-\kag - \hal} d \tz
  \les  \ang z^{-\kag -\hal} \min(|\xx|, 1) \nnrr{\ww} .
\]

If $|z| \leq 1$, since $\kag + \hal \in [0.32, 0.34]$ by \eqref{para:range}, by discussing $|\xx|\leq 1$ and $|\xx| > 1$, we obtain
\[
    M_{33,2} \les \min(|\xx|, 1) \nnrr{\ww} \int_{\nu \min(1, |\xx|)}^{ 4 |\xx| } \f{1}{ |\tz| } \ang \tz^{-\kag - \hal} d \tz
\les \ang z^{-\kag -\hal} \min(|\xx|, 1) \nnrr{\ww}  \log \nu^{-1}  .
\]

Applying the above estimates of $M_{3j}, M_{33,i}$ and $\log \nu^{-1} 
\les \nu^{\al-1}$ to \eqref{eq:psi_zz:I_lg_est2}, \eqref{eq:psi_zz:I_lg_est3}, we prove 
\beq\label{eq:psi_zz:I_lg_est4}
\bal
I_{zz, \geq \nu}
& \les P(\xx) |\xx|^{\al-1} M_{3} 
\les P(\xx) |\xx|^{\al-1} \ang z^{-\kag -\hal}  \min(  |\xx|, 1  ) ( \nu \nnz{\ww} + \nu^{\al-1} \nnrr{\ww} ) \\
& \les P(\xx) \min( |\xx|^{\al}, \ang \xx^{\al-1} \ang z^{-\kag-\hal} ) 
( \nu^{\al} \nnz{\ww} + \nu^{\al-1} \nnrr{\ww} ) .
\eal
\eeq

Therefore, combining estimate \eqref{eq:psi_zz:I_lg_est4},  \eqref{eq:psi_zz:I_lg_ez} on $I_{zz,\geq \nu}$, \eqref{eq:psi_zz:est2_nu} on $I_{zz, \nu}$, \eqref{eq:psi_zz:Izz_wz} on $I_{zz}(\pa_{\tz} \ww)$, 
 \eqref{eq:psi_zz:Ifar} on $I_{zz}^{\ffar}$
 and using the decomposition \eqref{id:psi_zz}, we prove 
\[
\bal
  |\pa_{zz} \psi| 
 & \les |I_{zz,\geq \nu}|
  + |I_{zz, \nu}|
  + |I_{zz}(\pa_\tz \ww)| + |I_{zz}^{\ffar}| \\
 & \les P(\xx) \min( |\xx|^{\al}, \ang \xx^{\al-1} \ang z^{-\kag-\hal} ) 
( \nu^{\al} \nnz{\ww} + \nu^{\al-1} \nnrr{\ww} )  \\
& \quad + \nu^{\al-1} \rhoo^{-1} \min( |\xx|^{\al}, \ang \xx^{\al-\hal-1} ) \nnrr{\ww}.
\eal
\]

Recall $P(\xx) = \ang \xx^{\kag} \rhoo(\xx)^{-1},| \xx| = |r, z|$, 
and $\rag  \asymp \f{\ang z^{\kag}}{\ang \xx^{\kag}} \les  1$ by \eqref{def:3d_wg}. Since 
\[
  P(\xx) \ang z^{-\kag} =\rhoo^{-1}(\xx) 
  \cdot  \ang \xx^{\kag} \ang z^{-\kag}
  \asymp \rhoo^{-1} \rag^{-1},
  \quad \ang \xx^{-\hal}
\leq \ang z^{-\hal}
\leq \ang z^{-\hal} \rag^{-1},
\]
combining the above two estimates, for any $\nu \in (0,1/2]$, we prove \eqref{eq:vel_zz}.

\subsection{Estimate of  \texps{$ \psi_r, \psi_z, \psi$}{ psi,r \  psi,z  \ psi} }\label{sec:d_psi}

The estimates of $\na_{r, z}$ are much simpler. Recall the formula of $\psi$ from \eqref{eq:5D_BS:rec}, $\JJ$ from \eqref{def:JJ_3D:b}, and $\pa_z \psi(0) = 2 \al \JJ(\infty) $. 
For $ \cT   \in \{ \Id, \pa_r, \pa_z \}$, we decompose 
\beq\label{eq:dpsi_decomp}
\bal
& \cT (\psi - z \psi_z(0) ) +  2 \al (\cT z )  \cdot \JJ(\xx) =
\cT \psi - 2 \al (\cT z) \cdot ( \JJ(\infty) - \JJ(\xx) ) \\
& =  \int \cT \KRF(\xx-\yy) W(\yy) d \yy  - 
 2 \al ( \cT z ) \cdot  ( \JJ(\infty) - \JJ(\xx) )  \\
& = \int_{|\yy| \leq 4 |\xx|} \cT \KRF(\xx-\yy) W(\yy) d \yy
+ \int_{|\yy| > 4 |\xx|} ( \cT \KRF(\xx-\yy) -  ( \cT z ) \cdot  (\pa_5 \KRF)(-\yy) ) W( \yy ) d \yy \\
& \quad - ( \cT z ) \cdot \int %_{ \tr \vee \tz \geq r \vee |z|, \, |\yy| \leq 4 |\xx| } 
( \one_{|\yy| \leq 4 |\xx|} - \one_{\tz \leq |\xx|}  )
(\pa_5 \KRF)(-\yy) ) W( \yy ) d \yy  \\
& \teq I_{\cT}^{\nne} 
+ I_{\cT}^{\ffar}  
+ ( \cT z ) \cdot I_{\cJ}^{\nne}.
\eal
\eeq

\subsubsection{Estimate of the far-field integral}

For $I_\cT^{\ffar}$, since $W(\yy) = |\yy_h|^{\al-1} \om( |\yy_h|, y_5 )$ is even in $\yy_h$ 
and odd in $y_5$, and the domain $|\yy|> 4 |\xx|$ is symmetric in $\yy$, 
symmetrizing the integral in $I_{\cT}^{\ffar} $, we obtain
\[
\bal
    I_{\cT}^{\ffar} 
   & = \f{1}{4} \int_{|\yy| > 4 |\xx|} \B( (\cT \KRF)(\xx-\yy) 
    + (\cT \KRF)( \xx_h + \yy_h, x_5-y_5 )  \\
   & \quad  - (\cT \KRF)( \xx_h - \yy_h, x_5 + y_5 ) 
    -   (\cT \KRF)( \xx_h + \yy_h, x_5 + y_5 )  - 4 (\cT z) \cdot (\pa_5 \KRF)(-\yy) \B) W( \yy ) d \yy.
\eal 
\]

Using the symmetrized estimates in \eqref{eq:K5D_est:symz}, \eqref{eq:K5D_est:symr}, \eqref{eq:K5D_est:sym0} and \eqref{eq:ker_Wlinf} on $W(\yy)$, we bound 
\[
    | z^{-1} I_{\Id}^{\ffar} |  +  |I_{\pa_r }^{\ffar} | +  |I_{\pa_z }^{\ffar} | 
     \les |\xx|^2 M,
\]
where 
\[
  M =   \int_{|\yy| \geq 4 |\xx|} |\yy|^{-6}  \min(|\yy|, 1) \tr^{\al-1}  \ang \tz^{-\hal} \rhoo^{-1} \rag^{-1} d \yy \nnla{\ww} .
\]

Since $|\yy| = |(\tr, \tz)|$ and  $|\yy|^{-6} \tr^{\al-1} \leq |\yy|^{-5} \tr^{\al-2} $, 
applying Lemma \ref{lem:int_W} with $\ell = 5$ and using $\al-1 < -\al$, we bound 
\[
  M \les \rhoo^{-1} \min( |\xx|^{\al-1} + |\xx|^{-\al}, \ |\xx|^{\al-\hal - 2} ) \nnla{\ww}
  \les \rhoo^{-1} \min( |\xx|^{\al-1} , \ |\xx|^{\al-\hal - 2} ) \nnla{\ww}.
\]
It follows 
\beq\label{eq:dpsi_far}
      |z^{-1} I_{\Id}^{\ffar} | 
+  |I_{\pa_r }^{\ffar} | +   |I_{\pa_z }^{\ffar} |  \les 
    \rhoo^{-1} \min( |\xx|^{\al+1} , \ |\xx|^{\al-\hal } ) \nnla{\ww} .
\eeq

\subsubsection{Estimate near the singularity $I_{\pa_r}^{\nne}, I_{\pa_z}^{\nne}$}

Recall $I_{\cT}^{\nne}$ from \eqref{eq:dpsi_decomp}. For $\cT = \pa_r, \pa_z$, 
since $W(\yy)$ and the domain $\{ |\yy| \leq 4 |\xx| \}$ only depends on $\tr, \tz$, applying \eqref{eq:K5D_est:dK} and then \eqref{eq:int_singu} with $k = 4$ and $W(\yy) = \tr^{\al-1} \om(\tr, \tz)$ by 
\eqref{eq:ker_W_id}, for $ i \in \{ r, z \}$, we bound 
\bseq\label{eq:dpsi:singu}
\beq
  |I_{ \pa_i }^{\nne}|
  \les \int_{|\yy| \leq 4 |\xx|} \f{1}{|\xx-\yy|^4} |W(\yy)| d \yy
  \les \int_{ |(\tr, \tz)| \leq 4 |\xx| }  \f{ \tr^3 \tr^{\al-1}  |\om(\tr, \tz)| }{ D_-  D_+^3 } 
  d \tr d \tz,
\eeq
where $D_{\pm} = |(r \pm \tr, z - \tz)|$. Applying \eqref{eq:ker_wwlinf_ne} on $\om$, $\tr \leq D_+$, 
and Lemma \ref{lem:CS_singu} with $k = 1, m = 1-\al \in [0.66, 0.68], 
\ell = \hal + \kag \in [0.32,0.34]$, 
and using $P(\xx) = \ang \xx^{\kag} \rhoo^{-1}$ from \eqref{eq:ker_wwlinf_ne}, we obtain 
\beq
\bal
    |I_{ \pa_i }^{\nne}|
   &  \les  P(\xx) \min( |\xx|, 1 ) \nnla{\ww} \int_{ |(\tr, \tz)| \leq 4 |\xx|  } \f{1}{D_- D_+^{1-\al}}
  \ang \tz^{-\hal -\kag} d \tr d \tz \\
  & \les P(\xx) \min( |\xx|, 1 ) \nnla{\ww} |\xx|^{\al} \ang \xx^{-\hal -\kag}
\les \rhoo^{-1} \min( |\xx|^{\al+1}, \ang \xx^{\al-\hal} ) \nnla{\ww}.
\eal 
\eeq
\eseq

\subsubsection{Estimate near the singularity $I_{\Id}^{\nne}$ }

To estimate $I_{\Id}^{\nne}$ in \eqref{eq:dpsi_decomp}, we assume $z \geq 0$ without loss of generality. Denote $B_{\pm} = | (\xx_h - \yy_h, z \pm \tz)|$. We estimate the symmetrized kernel in $\tz$:
\[
  M \teq  \f{1}{B_-^3} - \f{1}{B_+^3}
= \f{ B_+^3 - B_-^3 }{ B_+^3 B_-^3 }
= \f{ ( B_+^2 + B_- B_+ +  B_-^2)(B_+^2 - B_-^2)  }{ B_+^3 B_-^3 (B_+ + B_-) }.
\]
Clearly, for $\tz \geq 0$, we obtain $B_- \leq B_+, |z| + |\tz| \leq B_+$.  By definition, we obtain
\[
  |B_+^2 - B_-^2| = |(z+\tz)^2 - (z-\tz)^2 | =| 2 z \tz| \les |z| B_+,
  \quad \Rightarrow \quad  |M| \les \f{B_+^2 \cdot |z| B_+}{ B_-^3 B_+^4  } \les \f{|z|}{B_-^3 B_+}
  \les \f{|z|}{|\xx-\yy|^4}.
\]
Since $W(\yy)$ is odd in $y_5 =\tz$, using the above estimate and the integral estimate in \eqref{eq:dpsi:singu}, we bound 
\beq\label{eq:psi:singu}
\bal
   |\f{1}{z} I_{\Id}^{\nne}| 
  & = \f{1}{2 z} \B| \int_{|\yy|\leq 4|\xx|} \big( \f{1}{|\xx-\yy|^3} -\f{1}{ |(\xx_h - \yy_h, z + \tz)|^3 } \big) W(\yy) \B|
  \les  \int_{|\yy|\leq 4|\xx|} \f{1}{|\xx-\yy|^4} |W(\yy) | d \yy \\
& \les \rhoo^{-1} \min( |\xx|^{\al+1}, \ang \xx^{\al-\hal} ) \nnla{\ww}.
  \eal
\eeq

\subsubsection{Estimates of $I_{\cJ}^{\nne}$}

Recall $I_{\cJ}^{\nne}$ from \eqref{eq:dpsi_decomp} and $|\yy| = |(\tr, \tz)|$. Firstly, we have 
\[
  | \one_{|\yy| \leq 4 |\xx|} - \one_{|\tz| \leq |\xx|} |
  \leq | \one_{|\yy| \leq 4 |\xx|} - \one_{|\tz| \leq |\xx|, \tr \leq |\xx|}  |
  +  \one_{|\tz| \leq |\xx|, \tr \geq |\xx|} 
\leq \one_{ |\xx| \leq |\yy| \leq 4 |\xx| }   +  \one_{|\tz| \leq |\xx| \leq  \tr }  \, . 
\]

Since   $|\pa_5 \KRF(\yy) | \les \f{ |y_5| }{|\yy|^5}$,  using \eqref{eq:5D_int_id} and the above estimate, we bound 
\[
\bal
  | I_{\cJ}^{\nne}| & \les  
  \int ( \one_{|\tz| \leq |\xx|\leq \tr } +  \one_{ |\xx| \leq |\yy| \leq 4 |\xx| }      )
 \f{\tr^{\al-1} |\ww(\tr, \tz) \tr^3 \tz|}{ |(\tr, \tz)|^5  } d \tr d \tz 
 \teq I_1 + I_2. \\
\eal
\]
where $a\vee b = \max(a, b)$. In the domain of the integral, we have $|\xx|/4 \leq |(\tr, \tz)| \leq 4 |\xx|$.

For $I_1$,  using $|\tz|\leq |\xx|$, \eqref{eq:ker_wwlinf} for $\ww$, 
$ \min(|\yy|,1) |\yy|^{-1} \asymp \ang \yy^{-1}$, and then \eqref{eq:int_W_near2}, we bound
\[
  I_1 \les  |\xx| \int_{ |\tz| \leq |\xx|\leq \tr  }  |\yy|^{\al-2} \ang \yy^{-1}
\rag^{-1} \rhoo^{-1} \ang \tz^{-\hal } d \tr d \tz \nnla{\ww} 
\les \rhoo(\xx)^{-1} \min( |\xx|^{\al+1}, |\xx|^{\al-\hal} ).
\]

For $I_2$, since $|\yy| = |(\tr, \tz)| \asymp |\xx|$ in the domain, 
using \eqref{eq:ker_wwlinf_ne},  $\tr, |\tz| \leq |(\tr, \tz)|$, and then Lemma \ref{lem:CS_singu} with $k=m=n=0, \ell=\hal + \kag \in [0.32, 0.34]$ \eqref{para:range}, we bound 
\[
\bal
|I_2| &\les 
 \ang \xx^{\kag} \rhoo^{-1}  \min(|\xx|,1) \nnla{\ww}
 \int_{ \ssk{|\xx| \leq |(\tr, \tz)|\leq 4|\xx| } }  |(\tr, \tz)|^{\al-2}   \ang \tz^{-\hal -\kag} d \tr d \tz  \\
 & \les \ang \xx^{\kag} \rhoo^{-1}  \min(|\xx|,1) \nnla{\ww} |\xx|^{\al-2} \cdot |\xx|^2 \ang \xx^{-\hal-\kag}.
 \eal
\]
Combining the above estimates for $I_1, I_2$, we prove 
  \beq\label{eq:Ine_JJ_goal}
  | I_{\cJ}^{\nne}|  \les  \rhoo(\xx)^{-1}    \min( |\xx|^{\al+1 }, |\xx|^{\al - \hal} ) \nnla{\ww}.
\eeq

\vspace{0.1in}
\paragraph{\bf Proof of \eqref{eq:vel_r}, \eqref{eq:vel_psi}
on $\psi_r, \psi_z, \psi$} \

Combining \eqref{eq:psi:singu},  \eqref{eq:Ine_JJ_goal}, and \eqref{eq:dpsi_far}, we prove 
the estimate \eqref{eq:vel_psi} for $\psio$.

Combining \eqref{eq:dpsi:singu}, \eqref{eq:dpsi_far},  \eqref{eq:Ine_JJ_goal}, and \eqref{eq:dpsi_decomp}, we prove
\[
      \rhoo | \psi_r | 
    + \rhoo | \psio_z + 2 \al \JJ(\om) |  \les  \min( |\xx|^{\al+1}, | \xx|^{\al -\hal} ) \nnla{\ww} .
\]
It follows the estimate \eqref{eq:vel_psi} for $\psi_z$ and \eqref{eq:vel_r1} for $\psi_r$.

Using \eqref{def:psio}, we obtain $ \psi_z - \f{1}{z} \psi = \psio_z - \f{1}{z} \psio$. Using 
the estimates in \eqref{eq:vel_psi} for $\tf1z \psio + 2 \al \JJ, \pa_z \psio + 2 \al \JJ $, which we have just proved, and the triangle inequality, we prove estimate \eqref{eq:vel_psi} for $ \psi_z - \f{1}{z} \psi$.

When $|z| \gtr r$, since $|z| \asymp |\xx|$, estimate \eqref{eq:vel_r} for $\psi_r$ follows from the  estimate \eqref{eq:vel_r1} for $\psi_r$.

When $|z| \leq r$, we obtain $r \asymp |\xx|$.   Since $\psi(r, z), \ww(r, z)$ are odd in $z$, we obtain $\psi_r(r, 0) = 0$. 
 Using \eqref{eq:vel_rz} proved in Section \ref{sec:psi_rz} and \eqref{eq:vel_bd_opt}, for some $|\xi| \leq |z| $, we establish  \eqref{eq:vel_r}:
\[
   |\psi_r(r,z) | =| \psi_r(r, z) -\psi_r(r, 0)|
   =|z \cdot \pa_{rz}  \psi (r, \xi)| \les |z| \min( |\xx|^{\al}, |\xx|^{\al-\hal-1} )
   \nnn{\ww}.
\]

\subsubsection{ Estimate of $\JJ$}

We use the notation $\yy = (\tr, \tz)$. 
Using \eqref{eq:ker_wwlinf:a}, \eqref{eq:int_W_near}, 
and \eqref{eq:int_W_near2}, we bound 
  \[
\bal
   |\JJ(\om)(\xx)| & \les   \int_{ 0\leq \tz \leq |\xx|  } 
 \f{\tr^{2 + \al} \tz}{ |(\tr, \tz)|^5 } |\om| 
  d \tr d \tz  
  \les  \int_{ 0\leq \tz \leq |\xx|  }  
  |\tz|
  \ang \tz^{-\hal} |\yy|^{\al-3} \min(|\yy|, 1) \rhoo^{-1} \rag^{-1}(\tr, \tz) d \tr d \tz 
  \cdot \nnla{\ww} \\
  & \les \min( |\xx|^{\al+1},  |\cJaa(\xx)|^{1 + \kp} )  \nnla{\ww}.
\eal
\]

Note that $|\xx| \leq   |\tz | $ implies $|\xx|  \leq | (\tr, \tz)|$. 
Applying \eqref{eq:ker_wwlinf}, \eqref{eq:wg_est:ag} for $\rag$, using $\rhoo(y)^{-1} \les \ang \yy^{\epa} |\cJaa(\yy)|^{\kp}
\les \e^{-\kp}\ang \yy^{\epa}  $, the polar coordinate  $(\tr, \tz)
= (R \cos \xi, R \sin \xi)$, we estimate 
\[
\bal
  |\JJ(\om)(\xx) - \JJ(\om)(\infty)|    
 & \les  \int_{ |\xx|
  \leq  |(\tr , \tz )|     } 
   \f{\tr^{2 + \al} |\tz | }{ |(\tr, \tz)|^5 } |\om| 
   \les \int_{ |\xx|
 \leq  |(\tr , \tz )|     } 
 \f{\tr^{2 + \al} |\tz | }{ |(\tr, \tz)|^5 } |\sin \xi|^{-\kag} 
 \rhoo^{-1} \ang \tz^{-\hal} \min(|\yy|,1) \nnla{\ww} \\
& \les \e^{-\kp}
 \int_{ |\xx|    }^{\infty} \int_0^{\pi/2}
 |R|^{\al } \ang R^{-1+\epa}  |\sin \xi|^{-\kag}  \min( 1, (R \sin \xi)^{-\hal} ) d R d \xi \cdot \nnla{\ww}.
\eal
\]
Since $\al - \hal + \epa \asymp -\e$ by \eqref{ran:ep_all}, and $\hal + \kag \in [0,0.4]$ \eqref{para:range}, integrating $\xi$, we prove \eqref{eq:vel_J}:
\[
\bal
  |\JJ(\om)(\xx) - \JJ(\om)(\infty)|   
& 
   \les \e^{-\kp} \int_{|\xx| }^{\infty} 
   R^{\al} \ang R^{-1+ \epa -\hal} d R  \nnla{\ww} 
   \les \e^{-\kp - 1} \ang \xx^{\al -\hal + \epa}  \nnla{\ww} .
   \eal
\]

\subsection{Estimate of \texps{$\bar \Psi_{rzz}$}{ bar Psi,rzz } with \texps{$|z| \geq r$}{ |z| > r} }\label{sec:psi_rzz}

In this section, we prove \eqref{eq:u_bar_rzz} with $|z| \geq r$. We assume that $r>0$. 
The case $r =0$ can be obtained using continuity. 

 We abuse notation $\xx = (r, 0,0,0,z)$ \eqref{def:5D_x}. Recall the profile  $\wwb=\waa \in C^{3}$ from \eqref{def:3d_Om}, \eqref{eq:waa_sign}.  We introduce $\bar W(\yy) = \tr^{\al-1} \wwb$ as in \eqref{eq:psi_W}.  Since $\pa_r \bpsi(r, z) = \pa_{ x_1 } \bpsi( \xx ) |_{x_1 = r} $, using \eqref{eq:5D_BS:rec}, we obtain
\beq\label{id:psi_rzz}
\bal
    \pa_{r z z} \bpsi(r, z) & =  \pa_z^2 \pa_{x_1} (\KRF \ast \bar W )(\xx)
    =  P.V. (\pa_z \KRF \ast \pa_{y_1} \pa_{y_5} \bar W)(\xx)   \\
   & = P.V. \B( \int_{ |z - \tz| \leq |z| / 2 } + \int_{ |z - \tz| \geq |z|/2 } \B) \pa_z \KRF(\xx-\yy) \cdot \pa_{y_1} \pa_{\tz} \bar W(\yy) d \yy
   \teq I_{  r zz }^{\nne} + I_{r zz}^{\ffar} ,
\eal
\eeq
where $z = x_5$ and $P.V$ denotes the principle value of the integral near the singularity $\xx$.

Since $\bar W(\tr, \tz) = \tr^{\al-1} \wwb(\tz) =\tr^{\al-1} \waa(\tz)$, 
using estimates \eqref{eq:wa_upper_lower}, \eqref{eq:waa_reg}, and formulas of $\bar W$ in \eqref{eq:ker_W_id}, 
for $i\leq 2$, we obtain
\beq\label{eq:W5D_est}
\bga
 |\wwb| \les \min(|\tz|, \ang \tz^{-\hal}), \quad 
 \     |\pa^i_{\tz} \wwb| \les   \ang \tz^{-\hal-i},\quad
 | \pa_{y_1} \pa_\tz^i \bar W| \les \tr^{\al-2} |\pa_\tz^i \wwb| .
\ega
\eeq

\paragraph{\bf Estimate of $I_{rz z}^{\ffar}$}
For $I_{rz z}^{\ffar}$, since the kernel is not singular, using integration by parts in $\tz$, we obtain
\[
  I_{rzz}^{\ffar} = \int_{ |z -\tz| \geq |z| / 2 } 
  \pa_{zz} \KRF(\xx-\yy) \pa_{y_1} \bar W(\yy) d \yy
  + \int_{ |z -\tz| = |z| / 2  } \nn_{\tz} \cdot  \pa_z \KRF( \xx- \yy ) \pa_{y_1} \bar W(\yy) d \yy \, ,
\]
where $ \nn_{\tz}$ is the unit outer normal. 

For $I_1$, since the bounds of $|\pa_{y_1} \bar W|, |\pa_{y_1} \pa_{\tz} \bar W|$ in 
\eqref{eq:W5D_est} and the domains of the above integral only depends on $(\tr, \tz)$, using estimate \eqref{eq:K5D_est:dK} and \eqref{eq:int_singu} with $k = 4, 5$, we estimate 
\[
\bal
  |  I_{rzz}^{\ffar} |
& \les \int_{ |z -\tz| \geq |z| / 2 } 
  |\xx-\yy|^{-5} \tr^{\al-2} |\wwb| d \yy
  + \int_{ |z -\tz| = |z| / 2  } |\xx-\yy|^{-4} \tr^{\al-2} |\wwb| d \yy \\
& \les \int_{ |z -\tz| \geq |z| / 2 } 
  \f{\tr^{\al+1}}{ D_-^2 D_+^3 } |\wwb| d \tr d \tz
  + \int_{ |z -\tz| = |z| / 2  } \f{\tr^{\al+1} }{D_- D_+^3}  |\wwb| d \tr  .
\eal
\]

Recall $D_{\pm} = |( r \pm \tr, z - \tz )|$ from \eqref{def:Dpm}. 
Since $\wwb$ is constant in $r$, using $ \tr^{\al+1} / D_+^3 \leq D_-^{\al-2}$ and integrating over $\tr$, we bound 
\[
    |  I_{rzz}^{\ffar} | \les 
    \int_{ |z -\tz| \geq |z| / 2 } 
  \f{ |\wwb(0,\tz)| }{ D_-^{4 -\al} }   d \tr d \tz
  + \int_{ |z -\tz| = |z| / 2  }  \f{ |\wwb(0,\tz)| }{ D_-^{3 -\al} }  d \tr  
\les \int_{ |z -\tz| \geq |z| / 2 } 
  \f{ |\wwb(0,\tz)| }{  |z- \tz|^{3-\al}  }   d \tz
  +   \f{ |\wwb(0,\tz)| }{ |  z- \tz |^{2 -\al} }   \B|_{|z-\tz| = |z|/2}.
\]

Using the bound \eqref{eq:W5D_est} and integrating the above bounds, we prove 
\begin{align}\label{eq:psi_rzz:far}
      |  I_{rzz}^{\ffar} | & \les \min( |z|, \ang z^{-\hal} )  |z|^{\al-2}
      +  \int_{ |z|/2 \leq |\tz-z| \leq 4 |z|} |z|^{\al-3} \min( |\tz|, \ang \tz^{-\hal} ) 
      + \int_{ |\tz| \geq 3 |z| } |\tz|^{\al-3}  \min( |\tz|, \ang \tz^{-\hal} )  \notag  \\
      & \les  \min( |z|, \ang z^{-\hal} )  |z|^{\al-2}
    \les \min( |z|^{\al-1}, \ang z^{\al-\hal-2} ). 
\end{align}

\vs{0.1in}
\paragraph{\bf Estimate of $I_{rz z}^{\nne}$}
Next, we estimate $ I_{rzz}^{\nne}$ in \eqref{id:psi_rzz}. Since $\pa_z \KRF(\xx-\yy) = c \f{z - \tz}{|\xx-\yy|^5}$ is odd in $z -\tz$ and the domain is symmetric in $z-\tz$, we obtain
\[
  I_{rzz}^{\nne} = 
  \f{1}{2} P.V. \int_{ |z - \tz| \leq |z| / 2 }  \pa_z \KRF(\xx-\yy) \cdot ( \pa_{y_1} \pa_{\tz} \bar W(\yy_h, 
\tz) - \pa_{y_1} \pa_{\tz} \bar W(\yy_h, 2 z - \tz) ) d \yy.
\]

For $|z-\tz| \leq |z|/2$, we have $|z| \asymp |\tz|$. Using \eqref{eq:W5D_est} for $\bar W$,
and \eqref{eq:K5D_est:dK} for $\KRF$, we bound 
\[
    I_{rzz}^{\nne} \les \int_{|z-\tz| \leq |z|/2} \f{  |z - \tz| }{  |\xx-\yy|^4} \sup_{ |\xi - z| \leq |z|/2 
    } \tr^{\al-2} | \pa_z^2 \wwb| d \yy
  \les  |z|^{-1} \ang z^{-\hal-1}  \int_{|z-\tz| \leq |z|/2} \f{  |z - \tz| }{|\xx-\yy|^4} \tr^{\al-2} d \yy.
\]

Using \eqref{eq:int_singu} with $k=4$, and then use $\tr , D_-, \leq D_+$ and $|z-\tz|\leq D_-$, we bound 
\[
      I_{rzz}^{\nne} 
      \les  |z|^{-1} \ang z^{-\hal-1} \int_{|z-\tz| \leq |z|/2} \f{  |z - \tz| }{ D_- D_+^3 } \tr^{\al+1} d \tr d \tz
    \les  |z|^{-1} \ang z^{-\hal-1} \int_{|z-\tz| \leq |z|/2} \f{ 1}{ D_-^{2-\al} } d  \tr d \tz.
\]
By first integrating $\tr$ and then $\tz$, we prove 
\[
       | I_{rzz}^{\nne} |  \les  |z|^{-1} \ang z^{-\hal-1} 
         \int_{|z-\tz| \leq |z|/2} \f{ 1}{ |z-\tz|^{1-\al} } d  \tr d \tz
        \les |z|^{\al-1}  \ang z^{-\hal-1}  .
\]
Combining the above estimate and \eqref{eq:psi_rzz:far}, and using $|z| \geq r$, $|z|\gtr |\xx|$,
we prove \eqref{eq:u_bar_rzz}.

\subsection{Isotropic nonlocal estimates}\label{sec:iso_nonlocal}

We develop the following nonlocal estimates with radial weights and use them for nonlinear stability analysis in Section \ref{sec:blowup}.

\begin{prop}\label{prop:iso_est}
Suppose that $\psi = \BS(\ww)$ \eqref{eq:Euler2_psi}. Recall $\psio =\psi - \psi_z(0) z$ from \eqref{def:psio} and $\psi_z(0) = 2 \al \JJ(\infty)$. For any 
\beq\label{eq:iso_est:wg1}
\phi = |\xx|^{-\bbu} + |\xx|^{\g} ,
\quad  \g \in [\hal - \e, \hal ] ,\quad   \bbu \in [1, 3/2] ,
\eeq
we have the following estimates
\bseq\label{eq:vel_iso_est}
\begin{align}
   \phi \cdot \big( | \tf{1}{z} \psio + 2 \al \JJ|  + | \pa_z \psio + 2 \al \JJ| + |\pa_r \psi| \big) &  \les  |\xx|^{\al }  
   \nlinf{   \ww \phi } , \label{eq:vel_iso:1st}  \\  
   |\JJ(\xx) | & \les   \min( |\xx|^{\al + 1}, \e^{-1}  )    \nlinf{  \phi \ww } , \label{eq:vel_iso:J}   \\  
  |2 \al \JJ(\xx)- \psi_z(0)|  +   |\JJ(\xx) - \JJ(\infty)| & \les \e^{-1} \ang \xx^{\al - \g}   \nlinf{  \phi \ww } ,
      \label{eq:vel_iso:J2}
  \end{align}
  and 
  \begin{align}
  |\xx| \phi   \cdot \big(  | \pa_{zz} \psi| + |\pa_{r z} \psi| + \f{r}{|\xx|} | \pa_{rr} \psi | \big) & \les 
     |\xx|^{\al }  
    \nlinf{ \phi ( |\ww| + |\xx| \cdot |\na \ww| ) } .
 \label{eq:vel_iso:2nd}  
\end{align}
\eseq

\end{prop}

\begin{proof}

Given $\g \in [\hal - \e, \hal]$, we apply Proposition \ref{prop:vel_est} with $\kag = 0, \epa = \hal - \g \in [0, \e] , \kp = 0$ and 
\beq\label{eq:iso_est:wg2}
  \rhoo(\xx) = ( |\xx|^{-\bbu+1} + 1 ) \ang \xx^{- \hal + \g}, \quad \rag(\xx) \equiv 1.
\eeq
By definitions of $\nnla{\cdot}, \nnn{\cdot}, \nnrr{\cdot}$ norms \eqref{def:3d_norm} associated with $\rhoo, \rag \equiv 1$, we obtain 
\[
\bal
& \nnla{\ww} = 
\nlinf{ (|\xx|^{-1} + 1) \ang z^{\hal} \rhoo  \ww } 
 \les \nlinf{ (|\xx|^{-\bbu} + 1) \ang \xx^{\hal} \ang \xx^{-\hal + \g} \ww } 
\les  \nlinf{\phi \ww}, \\
& \nnr{\ww} + \nnz{\ww}
  \les 
   \nlinf{ \rhoo \cdot \ang \xx^{1 + \hal} ( | \pa_r \ww| + |\pa_z \ww|) }
  \les  \nlinf{ ( \ang \xx^{\g} + |\xx|^{-\bbu}) ( |\ww| + |\xx| \cdot |\na \ww| ) }
\teq M(\ww) .
\eal
\]

Estimate  \eqref{eq:vel_iso:J}, \eqref{eq:vel_iso:J2}  follows from  \eqref{eq:vel_J} and $\cJaa \les \e^{-1}$ by \eqref{eq:Ja_hat}.

Using \eqref{eq:vel_r1}, \eqref{eq:vel_psi}, Proposition \ref{prop:vel_est}, and the above bound, we prove 
\bseq\label{eq:iso_est:pf1}
\begin{align}
  \rhoo \big( |\xx|^{-1} r |\pa_{rr} \psi|
  + |\pa_{r z} \psi| \big) & \les  \min( |\xx|^{\al}, \ang \xx^{\al-\hal-1} ) M(\ww), 
\label{eq:iso_est:pf1:a}
  \\
\rhoo \cdot \big( | \tf{1}{z} \psio + 2 \al \JJ|  + | \pa_z \psio + 2 \al \JJ| + |\pa_r \psi| \big)   & \les  \min( |\xx|^{\al+1}, \ang \xx^{\al-\hal } ) \nlinf{ \phi \ww}.
\label{eq:iso_est:pf1:b}
\end{align}
\eseq

For $|z| > r$, using Proposition \ref{prop:vel_est} and $|z| \gtr |\xx|$, we obtain the following estimates for $\pa_{zz} \psi$:
\bseq\label{eq:iso_est:pf2}
\beq
\rhoo  |\pa_{zz} \psi| \les  \min( |\xx|^{\al},  \ang \xx^{\al-1} \ang z^{-\hal} ) \nnn{\ww}
  \les \min( |\xx|^{\al}, \ang \xx^{\al-1 -\hal }  ) M(\ww) \, .
\eeq

For $r > |z|$, we obtain $r \gtr |\xx|$. Using the elliptic equation \eqref{eq:Euler2_psi}, \eqref{eq:vel_iso_est} 
on $\pa_r \psi, \pa_{rr} \psi$, $\bbu \geq 1$,
and $|\ww| \les \min( |\xx|, |\xx|^{-\g}) M(\ww)$, we bound 
\beq
  \rhoo |\pa_{zz} \psi | \les \rhoo \tf{1}{r} ( r |\pa_{rr} \psi| + |\pa_r \psi| + |\ww r^{\al}| )
  \les \min( |\xx|^{\al}, \ang \xx^{\al-\hal-1} ) M(\ww).
\eeq
\eseq

From $\rhoo$ \eqref{eq:iso_est:wg2} and $\phi$ \eqref{eq:iso_est:wg1}, we have $\rhoo  \asymp \ang \xx^{- \hal - 1} |\xx| \phi$. 
Thus, multiplying \eqref{eq:iso_est:pf1:a}, \eqref{eq:iso_est:pf2} by 
$\ang \xx^{\hal +1}$, we prove \eqref{eq:vel_iso:2nd}. 
Multiplying \eqref{eq:iso_est:pf1:b} by $\ang \xx^{\hal + 1} |\xx|^{-1}$,
 we prove  prove \eqref{eq:vel_iso:1st}. 
\end{proof}

\subsection{Relation between the 1D and 3D velocity}\label{sec:1D_3D}

For a function $f$ defined on $\R$, we recall the 1D stream function $\psi$ from
\eqref{def:psi_1D}, $\cJa$ from \eqref{eq:Jw} 
\[
\psid( f )(z) = \int_{z \geq 0} ( |z+y|^{\al} - |z- y|^{\al})  f(y) d y , 
\quad \cJa(f)(z) = \int_0^z f(y) y^{\al-1} d y .
\]

We have the following relation between the 1D and 3D velocity.

\begin{lem}[\bf Error between the 1D and 3D stream function]\label{lem:vel_bc}
Suppose that $\al \in [0,\f13]$ and $\ww$ is odd in $z$. Let $\psi=  \BS(\ww)$ 
 be the 3D stream function associated with $\ww$ \eqref{eq:Euler2_psi}. Denote 
\beq\label{def:Kin}
\bga
\wwd(z)  = \ww(0, z) ,% \quad \psid( \wwd )(z) = \int_{\R} ( |z+y|^{\al} - |z- y|^{\al})  \wwd(y) d y , \ 
\quad  \psid = \psid(\wwd), \quad 
\psiod = \psid - \pa_z \psid(0) z, 
\quad \JJ(\xx) = \JJ(\ww)(\xx), 
\\
\Kin(r, z, \tz)  \teq   \int_{\td r \geq r}   \f{ \td r^{2 + \al}  }{ 
 ( \td r^2 + (z-\tz)^2 )^{3/2}  } - \f{ \td r^{2 + \al}  }{ 
 ( \td r^2 + (z+\tz)^2 )^{3/2}  }  d \td r . \\
 \ega
\eeq
Let $\cpsi$ be the constant defined in \eqref{def:cpsi}. We have 
\beq\label{eq:iden_bc:1}
  \psi(\ww)(0, z) = \psid (  \wwd  )(z) +
  \cpsi \int_{\R_+ \times \R_+} \Kin( \td r, z , \tz ) \pa_{\tr} \ww( \tr, \tz ) d \tr d \tz .
\eeq
As a result, if $\ww$ is constant in $r$, we obtain $\psi(\ww)(0, z) = \psid(\wwd)(z)$.

Let $\rhoo, \nnr{\cdot}$ be the weight and norm in \eqref{def:3d_wg} \eqref{def:3d_norm}. 
We have the following estimates 
\begin{align}\label{eq:comp_vel_bc}
  | (  \tf1z \psio(0,z)  + 2 \al   \JJ(0,z)  - \tf1z \psio_{\mw{1D}}(z) - 2 \al  \cJa(\wwd)(z) |  &\les
 \min( |z|,   \ang z^{ \al-\hal + \epa } \cJaa^{\kp} ) \nnra{ \ww},\notag  \\
|   \pa_z \psio(0,z)  +  2 \al \JJ(0, z)-  
\pa_z \psio_{\mw{1D}}(z) - 2 \al  \cJa(\wwd)(z)| &\les 
 \min( |z|,   \ang z^{ \al-\hal + \epa } \cJaa^{\kp} )   \nnra{ \ww}, \notag  \\
 |\JJ(0, z) - \cJa(\wwd)(0,z) | & \les \min( |z|,  \cJaa^{\kp + 1 } )  \nnra{\ww} \, ,  
\end{align}
where $\cJaa$ evaluates on $(0, z)$.

\end{lem}

\begin{proof}

We assume that $\ww$ is in the Schwarz class, and the general case follows from approximation. 
We consider $z \geq 0$. Recall the formula of $\psi$ from \eqref{eq:psi_cpsi}. Using the definition of $\Kin$ and integration by parts in $r$ , we obtain 
\beq\label{eq:deri_Kin0}
\bal
 & \psi(\ww)(0, z)  = \cpsi \int_{\R_+^2} - \pa_r \Kin(\tr, z - \tz) \ww(\tr, \tz) d \tr d \tz \\
 &  = \cpsi \B(  \int_{\R^+} ( - \Kin( \infty, z , \tz ) \ww( \infty, \tz) 
+ \Kin(0, z, \tz) \ww(0, \tz) ) d \tz 
 + \int_{\R_+^2}  \Kin(\tr, z,  \tz) \pa_r \ww(\tr, \tz) d \tr d \tz \B) .
\eal
\eeq
For $z, \tz \geq 0$, it is easy to obtain $|\Kin(r, z , \tz)| \les z \tz r^{\al - 2}$.  
Thus, we get $\Kin(\infty, z, \tz) \ww(\infty,\tz) = 0$. 
Next, we compute $\Kin(0, z , \tz )$. Using a change a variable $\tr \to y |z + \tz|$, %\quad  a = \tf{z - \tz}{z + \tz},
we obtain
\bseq\label{eq:deri_Kin1}
\beq
  \Kin(0, z, \tz) = \int_{\tr \geq 0}
  \f{ \td r^{2 + \al}  }{ 
 ( \td r^2 + (z-\tz)^2 )^{3/2}  } - \f{ \td r^{2 + \al}  }{ 
 ( \td r^2 + (z+\tz)^2 )^{3/2}  }  d \td r
 = |z + \tz|^{\al} H( \f{z - \tz}{z + \tz} ) \, ,
\eeq
where 
\beq
   H(a) := \int_0^{\infty} h(a, y) - h(1, y) d y, \quad h(a, y) =     \f{y^{2 + \al}}{ (y^2 + a^2)^{3/2} } \, .
\eeq
\eseq
Since the integrand is in $L^1$ for fixed $a$, $H(a)$ is well-defined.  Moreover, $H(a)$ is even in $a$ and $H(1) = 0$. Consider $a > 0$. For small $|\e| < a/2$, since 
$
  |\f{h(a + \e) - h(a)}{\e} | \leq \sup_{\xi \in [a/2, 2a]}
  | \pa_{\xi} \f{y^{2+\al}}{ (y^2 + \xi^2)^{3/2} } |% \les \sup_{\xi \in [a/2, 2a]}  (y^2 + \xi^2)^{\al/2-2} 
  \les (y^2 + a^2)^{\al/2-1}, 
$
which is $L^1$ integrable in $y$, using dominated convergence theorem and change of variable  $y \to a s$, we compute 
\beq\label{eq:iden_bc:pf0}
  H^{\pr}(a) = - 3 a \int_0^{\infty} \f{ y^{2+\al}}{ (y^2 + a^2)^{5/2} } d y
 = - 3 a^{ 4 + \al -5 } \int_0^{\infty} \f{s^{2+\al}}{ (s^2 + 1)^{5/2} }  d s
= - 3 a^{\al-1} \cpsib.
\eeq
where we use $\cpsib$ from \eqref{def:cpsi}. Thus, for $a>0$, $H(a) \in C^1$.
Since $H(1) = 0$, $\lim_{a \to 1^-} H(a) = H(1) = 0 $, and $H(a)$ is even in $a$, integrating $H^{\pr}(a)$, we obtain
\beq\label{eq:iden_bc:Ha}
  H(a) =   3 \al^{-1} \cpsib ( 1-  |a|^{\al} ) = \cpsi^{-1} ( 1-  |a|^{\al} ),
\eeq
where we have used $\cpsi = \f{\al}{ 3 \cpsib}$ from  \eqref{def:cpsi}. 
Using \eqref{eq:deri_Kin1} and \eqref{eq:iden_bc:Ha} with $ a = |\f{z - \tz}{z + \tz}|$, 
we prove 
\beq\label{eq:deri_Kin2}
  \Kin(0, z, \tz) =  \cpsi^{-1} |z + \tz|^{\al} ( 1 - |a|^{\al}   )
  = \cpsi^{-1} ( |z + \tz|^{\al} - |z - \tz|^{\al} ).
\eeq
Combining \eqref{eq:deri_Kin0},  \eqref{eq:deri_Kin1}, 
and \eqref{eq:deri_Kin2}, we prove \eqref{eq:iden_bc:1}.

\vspace{0.1in}
\paragraph{\bf Difference between $\JJ$ and $\cJa$}

Denote 
\beq\label{def:Kin_H2}
  f(r, z) = \f{r^{2+\al} }{ (r^2 + z^2)^{3/2} }, \quad 
   H_2( r, \tz ) = \int_r^{\infty}    \f{ \tr^{2 + \al} \tz }{ (\tr^2 + \tz^2)^{5/2} } d \tr
   = - \f{1}{3} \int_r^{\infty}  \pa_\tz f(\tr, \tz ) d \tr .
\eeq
For $r = 0$, using definition of $\JJ$ \eqref{def:JJ_3D}, 
 symmetrizing the integral in $\tz$, and using integration by parts, we obtain
\[
\bal
  \JJ(\ww)(0,z)  
  &=   \f{3 \cpsi}{   \al }  \int_{ \tr \geq 0, 0\leq \tz \leq |z| } 
   \f{ \tr^{2 + \al} \tz }{ (\tr^2 + \tz^2)^{5/2} } \ww(\tr, \tz) d \tr d \tz 
=  - \f{3 \cpsi}{   \al }  \int_{\tr \geq 0, 0\leq \tz \leq |z| } \pa_{\tr} H_2(\tr, \tz) \ww(\tr, \tz) d \tr d \tz  \\
& =  \f{3 \cpsi}{   \al }  \int_{ 0\leq \tz \leq |z|} H_2( 0, \tz) \ww(0, \tz) d \tr d \tz
+ \f{3 \cpsi}{   \al }  \int_{\tr \geq 0, 0\leq \tz \leq |z| } H_2 \pa_{\tr} \ww(\tr, \tz) d \tr d \tz
\teq I_{1} + I_2.
\eal
\]

Using the identity \eqref{eq:iden_bc:pf0}  divided by $-3 a$ with $a=\tz$, we obtain $H_2(0, \tz) =\tz \cdot   \cpsib  |\tz|^{\al-2}$. 
Since $\cpsib \cdot  \f{3 \cpsi}{   \al }=1$ by \eqref{def:cpsi},  using the operator \eqref{eq:Jw}, we obtain
\beq\label{eq:iden_bc:2}
I_1=   \int_{ 0}^{|z|}  \tz^{\al-1} \ww(0, \tz) d \tz
  = \cJa(\wwd)(z) ,
  \quad \JJ(\ww)(0, z)- \cJa(\wwd)(0, z) = I_2(z).
\eeq

Denote $\yy = (\tr, \tz)$. For any $\ell \geq 5$, a direct calculation yields 
\bseq\label{eq:iden_bc:pf1}
\beq
  P_{\ell} \teq \int_{\tr}^{\infty}  \f{ s^{2 + \al}  }{ ( s^2 + \tz^2)^{\ell/2} } d s 
  \les \int_{\tr}^{\infty} ( s + |\tz|)^{\al+ 2 - \ell} d s
  \les  (|\tr| + |\tz|)^{\al+3 -\ell} \les |\yy|^{\al+ 3 -\ell}.
\eeq
Using the above estimate with $\ell=5$, we bound $H_2$ in \eqref{def:Kin_H2} as 
\beq
 |H_2(\tr, \tz)| \les |\tz| \cdot |\yy|^{\al-2}, \quad |\yy| = |(\tr, \tz)|.
\eeq
\eseq

Using the definition of $\nnr{\cdot}$ from \eqref{def:3d_norm}, we obtain
\beq\label{eq:iden_bc:ww}
  |\pa_{\tr} \ww(\tr, \tz )| \les  
  \rag^{-1} \rhoo^{-1} \ang \yy^{-1} \ang \tz^{-\hal}  \nnr{\ww} .
\eeq

Using \eqref{eq:iden_bc:pf1} for $H_2$, \eqref{eq:iden_bc:ww} for $\pa_{\tr} \ww$,
and \eqref{eq:int_W_near} with $|\xx| = |\z|$, we prove the third estimate in \eqref{eq:comp_vel_bc}: %for $\cJ$-operator:
\[
\bal
   |I_2| & \les  \int_{\tr \geq 0, \tz \leq |z| } |H_2 \pa_{\tr} \ww(\tr, \tz)| d \tr d \tz 
  \les \int_{ 0\leq \tz \leq |z| }
|\tz| \cdot |\yy|^{\al-2}
  \rag^{-1} \rhoo^{-1} \ang \yy^{-1} \ang \tz^{-\hal}  d \tr d \tz  \nnr{\ww}  \\
& \les \min( |z|^{1 + \al}, |\cJaa(z)|^{1 + \kp} ) \nnr{\ww}.
\eal
\]

\vs{0.1in}
\paragraph{\bf Proof of \eqref{eq:comp_vel_bc}}

We focus on the estimate for $\pa_z \psi$.
Since $\pa_{z} \psi(0) = 2 \al \JJ(\infty)$, using identities \eqref{eq:iden_bc:1} and \eqref{eq:iden_bc:2}, we obtain
\begin{align}\label{eq:iden_bc:pf2}
I_{\D} &\teq   (\pa_z \psi(\ww) - \pa_z \psi(\ww)(0) ) + 2 \al \JJ(\ww)(0, z)
  -   (\pa_z \psi(\wwd) - \pa_z \psi(\wwd)(0) ) - 2 \al \cJa(\wwd)(0, z)  \notag \\
& =  \pa_z \psi(\ww) 
  -   \pa_z \psi(\wwd)  + 2 \al ( I_2(z) -  I_2(\infty)) 
   = \int_{\R_+ \times \R_+} 
M(\tr, \tz)
\cdot  \pa_{\tr}  \ww(\tr, \tz) d \tr d \tz \, ,
\end{align}
where $M(\tr, \tz)$ is given by 
\[
M(\tr, \tz) = \cpsi \pa_z  \Kin(\tr, z, \tz)
- 6 \cpsi  \one_{\tz \geq |z|} H_2(\tr, \tz) 
=  \cpsi \B( \int_{\tr}^{\infty} \pa_z ( f(s, z -\tz) -  f(s, z + \tz) )
+  2\one_{\tz \geq |z|} \pa_{\tz} f(s, \tz) d s \B)  \, ,
\]
and $f$ is defined in \eqref{def:Kin_H2}.

Since $\pa_z f$ is odd in $z$, for $|(s, \tz)| \geq 4 |z|$, using Taylor expansion at $z = 0$, we obtain
\[
\bal
 |\pa_z ( f(s, z - \tz) - f(s, z + \tz) ) + 2   \pa_\tz f(s, \tz)|
& =| - ( \pa_z f)(s, \tz - z) - (\pa_z f)(s, z + \tz) + 2  ( \pa_z f)( s, \tz) |  \\
& \les \sup\nolim_{|\xi| \leq z} | \pa_z^3 f( s, \tz + \xi ) | \cdot |z|^2 
\les |z|^2 s^{ 2 + \al } |(s, \tz)|^{-6} . \\
\eal
\]

\paragraph{\bf Estimate of $M$ for $|\tz| \geq 4 |z|$}.
Denote $\yy = (\tr, \tz)$. For $|\tz| \geq 4 |z|$, 
using the above estimates and \eqref{eq:iden_bc:pf1} with $\ell=6$, we bound 
\beq\label{eq:iden_bc:ker1}
|M(\tr, \tz)|
 \les  \B| \int_{s \geq \tr} 
\pa_z (f(s, z - \tz)  - f(s, z +\tz) ) + 2 \pa_z f(s, \tz)  d s \B|  \les \int_{s \geq \tr} \f{z^2 s^{2 + \al}}{ |(s , \tz)|^6 } d s
\les z^2  |\yy|^{\al-3}.
\eeq

\paragraph{\bf Estimate of $M$ for $|\tz| \leq 4 |z|$}. For $|\tz| \leq 4 |z|$, we bound $\Kin, H_2$ separately. If $s \geq 4 |z|$, 
by definition \eqref{def:Kin_H2}, we obtain
\[
 f_{\D}(s, z) \teq |\pa_z ( f(s, z - \tz) - f(s, z + \tz) )| 
  \les |\tz| \sup_{ |\xi| \leq 5 z } |\pa_\xi^2 f(s, \xi)|
\les |\tz| s^{2 + \al -5 } \les |\tz| s^{\al-3}.
\]

Since $|\pa_z f (s, z \pm \tz)| \les (| s |+ |z \pm \tz|)^{\al-2}$, for $|\tz| \leq 4 |z|$, we bound 
\[
\bal
  |\pa_z \Kin(\tr, z , \tz) | &\les \int_{\tr}^{\infty} f_{\D}(s, z)  d s
   \les \one_{\tr \leq 4 |z|} \int_{\tr}^{4 |z|} ( s + |z-\tz| )^{\al-2}  ds + \int_{\max(4 |z|, |\tr| )}^{\infty}
|\tz| s^{\al-3} d s \\
& \les \one_{\tr \leq 4 |z|}   ( \tr + |z-\tz| )^{\al-1} + |\tz| \cdot ( |z| + \tr)^{\al-2} .
\eal 
\]

For $|\tz| \leq 4 |z|$, 
using the above bound and 
\eqref{eq:iden_bc:pf1}  for $H_2$, $|\yy| \les |z| + |\tr|$ for $\yy = (\tr, \tz)$, 
and $|z-\tz| \les |z|$, we bound 
\[
|M(\tr, \tz)| \les 
    |\pa_z \Kin(\tr, z , \tz) | 
+  \one_{\tz \geq |z|} | H_2(\tr, \tz) | 
  \les  (\tr + |z-\tz| )^{\al-1} \one_{\tr \leq  4 |z|} + |\tz| \cdot (|z| + \tr)^{\al-2} 
 + \one_{\tz \geq |z|} |\tz| \cdot |\yy|^{\al-2} \, . \\
\]
When $\tr \leq |z|$, since $|z-\tz| \les |z|$, the last two terms can be bounded by the first term. We prove 
\beq\label{eq:iden_bc:ker2}
    |M(\tr, \tz)|
  \les  (\tr + |z-\tz| )^{\al-1} \one_{\tr \leq  4 |z|} 
 + \one_{\tr \geq |z|} |\tz| \cdot |\yy|^{\al-2} , \quad \forall |\tz| \leq 4 |z|.
 \eeq

Applying estimate \eqref{eq:iden_bc:ww} for $\pa_\tr \ww$, 
estimates \eqref{eq:iden_bc:ker1},  \eqref{eq:iden_bc:ker2} for $M(\tr, \tz)$, we estimate \eqref{eq:iden_bc:pf2} as 
\[
\bal
  |I_{\D}|  \les \int \B(  (\tr + |z-\tz|)^{\al-1} & \one_{\tr, |\tz| \leq  4 |z|} +
 \one_{\tr \geq |z|, |\tz| \leq 4 |z|} |\tz| \cdot |\yy|^{\al-2} 
  +  \one_{|\tz| \geq 4 |z|} |z|^2 |\yy|^{\al-3}  \B) \\
  & \qquad \cdot \rag^{-1} \rhoo^{-1} \ang \yy^{-1} \ang \tz^{-\hal}   d \tr d \tz  \nnr{\ww} \teq ( I_{\D, 1} + I_{\D, 2} + I_{\D, 3} ) \cdot \nnr{\ww} ,
\eal
\]
where $I_{\D, i}$ denotes the integral for the $i$-th summand.

For $I_{\D, 2}$, applying $|\tz| \leq |z|$, \eqref{eq:int_W_near2} with $|\xx| = |z|$, 
and $\rhoo^{-1}(z) \les \ang z^{\epa} \cJaa(z)^{\kp}$ \eqref{def:3d_wg}, we obtain
\[
  I_{\D, 2} \les |z| \rhoo^{-1}(z) \min( |z|^{\al}, \ang z^{\al-\hal} )
\les   \min( |z|^{1 + \al} , \cJaa^{\kp} \ang z^{\al - \hal + \epa} ) .
\]

For $I_{\D, 3}$, since $\one_{|\tz| \geq 4 |z|} \leq \one_{ |(\tr, \tz)|  \geq 4 |z| }$, applying \eqref{eq:int_W_far1} with $\ell=5$ and using $\bbu \geq 1$, we bound 
\[
  I_{\D, 3} \les |z|^2
   |\cJaa(z)|^{\kp}    \min(1 + |z|^{\al + \bbu + 3 - 5}, |z|^{ \al - \hal + 3 -5 + \epa })   
   \les \min( |z|^{1 + \al}, \ang z^{\al-\hal + \epa}    |\cJaa(z)|^{\kp} ) .
\]

For $I_{\D, 1}$, using \eqref{eq:wg_est:rc} , $\rag^{-1} \les \ang \yy^{\kag} \ang \tz^{-\kag}$, 
and then Lemma \ref{lem:CS_singu} with
$ (r, z) = (0, z)$, $k = 1 -\al, m = 0, \ell = \hal + \kag \in [0.32, 0.34], n = 1 - \kag < 1$, we bound 
\[
  I_{\D, 1} \les \rhoo(|z|)^{-1}  \int_{ |(\tr, \tz)| \leq 8 |z| } (\tr + |z-\tz|)^{\al-1} 
  \ang \yy^{\kag - 1} \ang \tz^{-\hal - \kag } d \tr d \tz 
\les \rhoo(|z|)^{-1} |z|^{2 - (1-\al) } \ang z^{ -1 -\hal} .
\]

Since $\rhoo(z)^{-1} \les \ang z^{\epa} \cJaa(z)^{\kp}$ by \eqref{def:3d_wg}, combining the above estimates, we prove estimate \eqref{eq:comp_vel_bc} for $I_{\D}$ defined in \eqref{eq:iden_bc:pf2}
\[
  |I_{\D}| \les \min( |z|, \cJaa^{\kp} \ang z^{\al-\hal + \epa} ) \nnr{\ww}. 
\]

The first estimate in \eqref{eq:comp_vel_bc}  for $\psio$  is proved similarly and it thus omitted.
\end{proof}

\section{Estimates of nonlocal operators and qualitative regularity estimates}\label{sec:est_LQR}

In this section, we first use Proposition \ref{prop:vel_est} 
and Corollary \ref{cor:u_bar} to estimate the linear nonlocal 
operator $\Lpsi$ \eqref{eq:fix_recall}, the transport coefficient $\QQ$ \eqref{def:Q},  and the error $\bar \cR$ \eqref{eq:fix_recall}. Then we establish qualitative regularity estimates for solution of the fixed point map \eqref{eq:fix_pt_2D}.  

Recall  $\bar \cR$ and $\Lpsi(\ww)$ from \eqref{eq:lin_2D} 
\beq\label{eq:fix_recall}
\bal
\bar \cR & =   
    \bar c_{\om} - (1 - \al) \bar \Psi_z -  (\bar c_l z +  2 \bar \Psi + r \pa_r \bar \Psi  )  \f{ \pa_z \wwb }{\wwb}
     , \\
 \Lpsi(\ww) & = - (1 - \al) \psio_z  -  (2 \psio + r \pa_r \psi)  \f{ \pa_z \wwb }{\wwb} .
\eal
\eeq

We extract the singular part $r \pa_r \psi, r \pa_r \bpsi, r \pa_r \Psi$ in these terms and 
decompose them as follows 
\beq\label{def:LR_decomp}
\bal
\bar \cR & = \Rreg + \Rsin, \quad 
&& \Rreg \teq  
   \bar c_{\om} - (1 - \al) \bar \Psi_z -  (\bar c_l z +  2 \bar \Psi   )  \tf{ \pa_z \wwb }{\wwb}, 
\quad  && \Rsin \teq -  r \pa_r  \bar \Psi \cdot \tf{\pa_z \wwb}{\wwb}  , \\
\Lpsi & = \Lreg + \Lsin, \quad  && \Lreg \teq 
- (1 - \al) \psio_z  -  2 \psio  \cdot  \tf{ \pa_z \wwb }{\wwb} , 
\quad && \Lsin \teq -  r \pa_r  \psi \cdot \tf{\pa_z \wwb}{\wwb}  , \\
Q^z & = \Qreg^z + \Qsin^z, \quad && \Qreg^z \teq c_l z + 2 \Psi, 
\quad && \Qsin^z \teq r \pa_r \Psi. 
\eal
\eeq

We introduce similar decompositions for $\bar \QQ, \td \QQ$, 
but we do not need to decompose $Q^r$.  We estimate the regular parts $\Rreg, \Lreg, \Qreg^z$ using  Proposition \ref{prop:vel_est} and Corollary \ref{cor:u_bar}.  In the following estimates, all weighted terms are bounded for $|\xx|$ not very large. The main difficulty is to overcome the factor $\cJaa \asymp \min(\lgp x, \e^{-1})$ for large $|\xx|$. 
The reader may focus on $\psio, \pa_z \psio$ terms involving $\JJ$. 
Other mixed derivatives of $\psio$ are \emph{much smaller}. See also the discussion 
in \hyr[sec:idea_nonlocal]{\its Nonlocal Estimates}.

We cannot bound the singular parts $\Rsin, \Lsin, \Qsin^z$ directly using the energy norms.
Instead, we  estimate them using the crucial integration by parts in Lemma \ref{lem:IBP}
and following Section \ref{sec:idea_singular}.

\subsection{Estimate linear nonlocal operators $\Lpsi$}

We have the following results for $\Lpsi(\ww)$ in \eqref{eq:fix_recall}.
\begin{prop}\label{prop:Lpsi}
Recall the weights $\rhoo, \rag$ from \eqref{def:3d_wg}, parameters $\kp, \hal, \epa$ from \eqref{def:kp}, \eqref{ran:ep2}, \eqref{norm:Xc}, and function $\cJaa$ from \eqref{eq:Ja_hat}.
Let $\psi, \bar \Psi$ be the stream functions associated with $\om, \wwb$.
Denote $\xx = (r, z)$. Then $\Lpsi(\ww)$ satisfies
\bseq\label{eq:Lpsi_linf_est} 
\begin{align}
  |\Lpsi(\ww)|  & \les 
   \big(1 +  \one_{ r> z }  |\lgp z|^{ \kp - 2 }  \cdot |\cJaa(\xx) |    \big)  |\cJaa(\xx)|^{\kp} \cdot  \nnn{\ww}, \label{eq:Lpsi_linf_est:a} \\
  |\Lpsi(\ww)| & \les |\xx |^{\al+1}  \nnn{\ww}   \label{eq:Lpsi_linf_est:b}.
\end{align}
\eseq
The $z$-derivative  satisfies 
\bseq\label{eq:Lpsi_Cz_est}
\begin{align}
  \rag \rhoo | z \pa_z \cL_{\psi}(\ww) | & \les 
  \min\big( |\xx|^{1+\al},     \ang \xx^{- \epa} \big) \nnn{ \ww }, \label{eq:Lpsi_Cz_est:wg} \\
| z \pa_z \cL_{\psi}(\ww) | & \les ( \ang \xx^{\al-\hal + \epa} + \cJaa(\xx) \cJaa(z)^{-1} ) 
\cdot |\cJaa(\xx)|^{\kp}
\nnn{\ww} \label{eq:Lpsi_Cz_est:nowg} .
\end{align}
\eseq
Finally,  the $r$-derivative of the regular part $\Lreg$ \eqref{def:LR_decomp} satisfies
\beq\label{eq:Lreg_Cr_est}
\bal
 |\G \pa_r \Lreg(\ww)| & \les  \min( |\xx|^{\al }, |\xx|^{\al- \hal - 1} ) \cdot \nnn{\ww} , \\
 | \pa_r \Lreg(\ww)|  
 & \les    \ang \xx^{-1 + \al -\hal + \epa}  \cJak  \nnn{\ww}. \\
 \eal
 \eeq

All the implicit constants are \emph{independent} of $ \epa, \hal, \al, \e$.
\end{prop}

\begin{proof}

\textbf{Estimate of $\Lpsi$}
We rewrite $\Lpsi(\ww)$ from \eqref{eq:fix_recall} as 
\[
\bal
  \Lpsi(\ww)  & = - (1-\al) ( \pa_z \psio - \f{\psio}{z} ) - \psio \B( \f{2 \pa_z \wwb}{\wwb} + \f{ 1-\al }{z}  \B)
  - r \pa_r \psi \f{\pa_z \wwb}{\wwb}
  \teq P_1 + P_2 + P_3 .
  \eal
\]

To simplify notation, we write $\Lpsi$ for $\Lpsi(\ww)$. We recall $\wwb = \waa$ \eqref{def:3d_Om}.

Using \eqref{def:psio}, we obtain $ \pa_z \psio - \f{\psio}{z} = \pa_z \psi - \f{\psi}{z}$. 
For $P_1, P_3$, applying Proposition \ref{prop:vel_est}, \eqref{eq:vel_bd_opt} with weights $\G$ defined in \eqref{def:3d_wg} and  using  $ \rhoo^{-1} \les \ang \xx^{\epa} \cJak  $, 
$\epa + \al -\hal < 0$ by \eqref{ran:ep_all}, 
$r\leq |\xx|$, and $|\f{\pa_z \wwb}{\wwb}| \les \f{1}{|z|}$ by \eqref{eq:1D_prof_est2:c}, we bound  
\[
\bal
  |P_1| & \les \min( |\xx|^{1+\al},    \ang \xx^{ \epa+ \al -\hal} ) \cJak \nnrr{\ww}
  \les  \min(  |\xx|^{1+\al}, 1) \cdot \cJak  \nnrr{\ww} , \\
  |P_3| & \les | \f{r \pa_r \psi}{z} | \les \min(|\xx|^{1+\al},  \ang \xx^{ \epa + \al - \hal } )   \cJak \nnn{\ww}
  \les \min(  |\xx|^{1+\al}, 1) \cdot \cJak \nnn{\ww}. 
\eal 
\]

Next, we estimate the coefficient in $P_2$.  Since $ 1 - \al = \f{2}{3} + \e$, using \eqref{eq:1D_prof_est2:c}, we obtain 
\[
  \B| \f{2 \pa_z \wwb}{\wwb} + \f{1-\al}{z} \B| 
  \les   \B| \f{2 \pa_z \wwb}{\wwb} + \f{2}{3 z} \B| + \f{\e}{z}
  \les \e |z|^{-1} + |\lgp z|^{\kp-2} |z|^{-1} .
\]

For $P_2$, applying \eqref{eq:vel_psi} and \eqref{eq:vel_J}, 
and $\rhoo^{-1} \les \ang \xx^{\epa} \cJak$ \eqref{def:3d_wg}, we bound 
\[
\bal
  |P_2| &  \les (\e + |\lgp z|^{ \kp - 2 } )|\f{\psio}{z}| 
  \les  (\e + |\lgp z|^{ \kp - 2 } )  \big( |\f{ \psio}{z} + 2\al \JJ|  
  + |2 \al \JJ |\big) \\
  & \les (\e + |\lgp z|^{ \kp - 2 } ) \cdot \min\B(|\xx|^{1+\al},  \ang \xx^{ \al- \hal + \epa} 
  + \cJaa(\xx)  \B) \cJak \nnrr{\ww}.
\eal
\]

For $r \leq |z|$, we obtain $\cJaa(\xx) \les \cJaa(z) \les \lgp z$ \eqref{eq:Ja_hat}. 
 For $r > |z|$, we bound $\cJaa(\xx) \les \e^{-1}$. 
 Since $\al -\hal + \epa< 0$ by \eqref{ran:ep_all}, $\kp<1$ \eqref{def:kp}, we bound
 \[
\bal
   |P_2| & \les (\e + |\lgp z|^{ \kp - 2 } ) 
  \cdot \min( |\xx|^{1+\al}, \cJaa(\xx) ) \cJak \nnrr{\ww} \\
  & \les  \min\B( |\xx|^{1+\al},  \e \cdot \e^{-1} + 1
+\one_{ r> z }  |\lgp z|^{\kp-2} 
  \cdot \cJaa(\xx)   \B) \cJak \nnrr{\ww} .
 \eal
 \]

Summarizing the above estimates for $P_i$, and using $\nnrr{\ww} \leq \nnn{\ww}$, we establish \eqref{eq:Lpsi_linf_est}:
\[
\bal
    |\Lpsi(\ww) | \les |\xx|^{1+\al} \nnn{\ww}, \quad 
  |\Lpsi(\ww) | &\les  
  \one_{ r> z }  |\lgp z|^{ \kp - 2 }  \cdot \cJakkk  \nnn{\ww}   + \cJak \nnn{\ww} .
  \eal
\]

\paragraph{\bf Estimate of $\pa_z \Lpsi$}

Recall $\Lpsi$ from \eqref{eq:fix_recall} and $\psio, \pa_r \psio =\pa_r \psi, \pa_z^2 \psio =\pa_z^2 \psi$ from \eqref{def:psio}.  We decompose 
\[
z \pa_z \Lpsi = - (1-\al) z \pa_{z z} \psi -
  2 z \pa_z ( \f{ \psio}{z} \cdot \f{ z \pa_z  \wwb}{\wwb} ) 
  -  z \pa_z ( r \pa_r \psi \cdot \f{ \pa_z \wwb}{ \wwb }  ) .
\]

Using \eqref{eq:1D_prof_est2:c}, \eqref{eq:1D_prof_dxx}, and $\wwb(r, z) = \waa(z)$, we bound 
\[
   |\f{z\pa_z \wwb}{\wwb} | \les 1,  \quad  | z \pa_z (\f{z \pa_z \wwb}{\wwb} ) | \les  \cJaa(z)^{-1} ,
   \quad | z \pa_z (\f{ \pa_z \wwb}{\wwb} ) | \les z^{-1}.
\]

By definition of $\psio$ \eqref{def:psio}, we get $\pa_z \f{\psio}{z} = \pa_z \f{\psi}{z}$. Using the above estimate,
we bound  $z \pa_z \Lpsi$ as 
\[
\bal 
 \rhoo | z \pa_z \cL_{\psi} |
&  \les \rhoo \B( | z\pa_{zz} \psi| + 
|z \pa_z ( \f{\psi}{z} )  \cdot\f{z \pa_z \wwb}{\wwb} |
+ \B|  \f{\psio}{z} \cdot z \pa_z ( \f{ z \pa_z \wwb}{\wwb}) \B| 
+ | r \psi_{r z}  \cdot \f{ z \pa_z \wwb}{ \wwb } |
+ | r \psi_r \cdot z \pa_z ( \f{\pa_z \wwb}{ \wwb } ) | \\
& \les \rhoo \B( | z\pa_{zz} \psi| + 
| \pa_z \psi - \f{\psi}{z}  |
+ \B|  \f{\psio}{z} \B| \cJaa(z)^{-1} 
+ | r \psi_{r z}   |
+ | r \psi_r \cdot z^{-1} | \B) .
\eal  
\]

Using Proposition \ref{prop:vel_est}, the bound \eqref{eq:vel_bd_opt}, 
and $|(r, z)| = |\xx|$, we obtain %the near-field estimate 
\bseq\label{eq:nloc_Cz_est1}
\beq\label{eq:nloc_Cz_est1:a}
 \rhoo | z \pa_z \cL_{\psi} | \les |\xx|^{\al+1} \nnn{\ww},
\eeq
and %the far-field estimate 
\beq
 \rhoo | z \pa_z \cL_{\psi} |   \les 
 ( \rag^{-1} \ang \xx^{\al-1} |z| \ang z^{-\hal} 
 + |\xx|^{\al-\hal}  )  \nnn{\ww}
+  \rhoo \B| \f{\psio}{z} \B| | \cJaa(z)^{-1}  .
\eeq

Since $\kag < \f{1}{100}$ by \eqref{def:rag} and $1-\hal, 1-\al > \f{1}{2}$ by \eqref{ran:ep2}, we get 
\beq\label{eq:rag_wg}
\rag^{-1} \ang \xx^{\al-1} \ang z \cdot \ang z^{-\hal} \les \rag^{-1} \ang \xx^{\al-1} \ang z^{1-\hal}\les \ang \xx^{\al-1 + \kag} \ang z^{1-\hal -\kag} \les \ang \xx^{\al -\hal}.
\eeq
It follows 
\beq\label{eq:nloc_Cz_est1:c}
   \rhoo | z \pa_z \cL_{\psi} |   \les 
  \ang \xx ^{\al-\hal}    \nnn{\ww}
+  \rhoo | \tf{\psio}{z} |  \cJaa(z)^{-1}  .
\eeq

\eseq

\paragraph{\bf Anisotropic weight and decay in $\xx$}

Recall $\rhoo$ from \eqref{def:3d_wg}. Using Proposition \ref{prop:vel_est}, we bound
\beq\label{eq:nloc_Cz_est2}
\bal
 II & \teq \rhoo \B| \f{\psio}{z} \B| | \cJaa(z)^{-1}
  \les \B( \ang \xx^{\al - \hal}  
  +  \rhoo \cdot \min( |\xx|^{1 + \al} , \cJaa(\xx)^{1 + \kp} )  \B) \cdot \nnrr{\ww} \cJaa(z)^{-1}  .\\
\eal
\eeq

Since $ \al - \hal < - \epa$ by \eqref{ran:ep_all} and $|z|\leq |\xx|$, we further bound 
\[
  II \les \ang \xx^{-\epa} (1 + \cJaa(\xx) \cJaa(z)^{-1}) \nnrr{\ww}
  \les \ang \xx^{-\epa}  \cJaa(\xx) \cJaa(z)^{-1} \nnrr{\ww}.
\]

For $\xx \gg z$ and large $\xx$, the coefficients in the above bound become very large and has size $\e^{-1}$. This difficulty arises since the profile $\wwb$ \textit{only} decays in $z$. 
To overcome this difficulty, 
we use crucially the angular weight $\rag$ to obtain global decay estimates in $\xx$
from $\cJaa(z)^{-1}$ in \eqref{eq:nloc_Cz_est2}.

Recall  $\rag \asymp (\f{\ang z }{\ang \xx})^{\kag}$. 
Using \eqref{eq:err_ag_decay} and $\kag  \gtr 1$ by \eqref{def:rag}, we obtain 
\[
\bal
\rag II 
 & \les 
 \cJaa(\xx)   \nnrr{\ww} \cdot \ang \xx^{- \epa} \cJaa(\xx)^{-1}  \les \ang \xx^{- \epa}  \nnrr{\ww}.
\eal
\]

Thus, combining the above estimate, \eqref{eq:nloc_Cz_est1}, using the norm \eqref{def:3d_norm}, 
and $\al -\hal < -\epa$ by \eqref{ran:ep_all},  we prove \eqref{eq:Lpsi_Cz_est:wg}:
\[
  \rag \rhoo | z \pa_z \cL_{\psi} | \les 
  \min( |\xx|^{1+\al},   \ang \xx^{\al-\hal}  +  \ang \xx^{- \epa}  ) \nnn{ \ww }
  \les   \min( |\xx|^{1+\al},   \ang \xx^{- \epa}  ) \nnn{ \ww } .
\]

Since $\rhoo^{-1} \les \ang \xx^{\epa} \cJaa $, from  \eqref{eq:nloc_Cz_est2} and \eqref{eq:nloc_Cz_est1:c},  we obtain
\[
\bal
 | z \pa_z \cL_{\psi} | & \les \ang \xx^{\al - \hal + \epa} \cJak \nnn{\ww}
  +   |  \tf{\psio}{z} |  \cJaa(z)^{-1}  , \\
  | \tf{\psio}{z} |\cdot | \cJaa(z)^{-1}
 & \les  ( \ang \xx^{\al-\hal + \epa} 
+  \cJaa(\xx) \cJaa(z)^{-1} ) \cJak  \nnrr{\ww} .
\eal
\]

Since $ \nnrr{\ww} \leq \nnn{\ww}$, combining the above two estimates, We prove \eqref{eq:Lpsi_Cz_est:nowg}.

\vs{0.05in}

\paragraph{\bf Estimate of $\pa_r \Lreg$}

Recall $\Lreg$ from \eqref{eq:fix_recall} and \eqref{def:LR_decomp}:
\[
  \pa_r \Lreg  =  - (1-\al) ( \pa_{r z} \psi - \f{\pa_r \psi}{z} ) - \f{ \pa_r \psi }{z} \cdot \B( \f{2 z\pa_z \wwb}{\wwb} + (1-\al)  \B). 
  \]

Using Proposition \ref{prop:vel_est} and \eqref{eq:1D_prof_est2:c},  $\G^{-1} \les \ang \xx^{\epa} \cJak $ by \eqref{def:3d_wg}, and $ \epa + \al -\hal < 0$ by \eqref{ran:ep_all}, we obtain
\[
\bal
\rhoo(\xx) | \pa_r \Lreg | & \les  \min( |\xx|^{\al }, |\xx|^{\al- \hal - 1} ) \cdot   \nnn{\ww}  , \\
    |\pa_r \Lreg| 
 & \les    \ang \xx^{\al- \hal - 1 + \epa }  \cdot \cJak \nnn{\ww} ,
\eal 
\]
and \eqref{eq:Lreg_Cr_est}. We complete the proof.
\end{proof}

\subsection{Estimate transport coefficient}

In this section, we estimate $\QQ = \bar \QQ + \td \QQ(\ww)$ from \eqref{def:Q}. Recall $\UU$ from \eqref{def:3d_vel}.  We introduce the following notations
\beq\label{def:qq}
\bal
\qqz & \teq  \f{Q^z}{z}  = c_l + \f{ U^z }{z}
= c_l + \f{2 \Psi}{z} + \f{ r \pa_r \Psi }{z} ,
\quad \qqr \teq  \f{ Q^r }{r} = c_l +  \f{ U^r}{r} = c_l - \pa_z \Psi \, . \\
\eal
\eeq

Recall $\rhoo$ from \eqref{def:3d_wg}. To track different bounds, we introduce 
 \beq\label{def:CC_bounds}
 \bal
 \CCa(\ww, \xx) &= \min( |\xx|^{\al}, |\xx|^{\al- \hal - 1} ) \cdot  (1 + \ang \xx^{\epa} 
 \cJak \nnn{\ww}) , \\
\CCj(\ww, \xx) & = \min( |\xx|^{\al+1}, \cJaa ) \cdot (1 +  \cJak \nnn{\ww}).
\eal
 \eeq

If the perturbation $\nnn{\ww}$ is small enough (see \eqref{eq:3D_boot_q:1}), $\CCa$ has a size of $O(1)$.
We use the subscript \textit{J} to indicate that $\CCj$ contains a large factor $\cJaa$ compared to $\CCa$. 

Since $\al-\hal + \epa < 0$ by \eqref{ran:ep_all} and $\cJaa(\xx) \gtr 1$, by definition, we obtain
\beq\label{eq:CC_bounds}
\bal
 |\xx| \CCa(\ww, \xx) & \les 
 \min( |\xx|^{\al + 1}, \ang \xx^{\al- \hal} ) \cdot ( 1 + \ang \xx^{\epa} \cJaa^{\kp} \nnn{\ww} )
 \les  \CCj.
\eal
\eeq

In the next proposition, we estimate the transport coefficients $\QQ$.
\begin{prop}\label{prop:Q}
 Recall  $\Qreg^z =  c_l z + 2 \Psi, Q^r =  c_l r - r \pa_z \Psi$
from \eqref{def:Q}, \eqref{def:3d_vel}, \eqref{def:LR_decomp}. We have the following estimates 
\bseq\label{eq:Q_mix_est}
\begin{align}
  |\pa_z  Q^r | & \les r \ang z^{-1} \ang \xx \CCa(\ww, \xx), \label{eq:Q_mix_est:Qr} \\
   |\pa_r \Qreg^z | & \les |z| \CCa(\ww, \xx)  \label{eq:Q_mix_est:Qz} ,
\end{align}
\eseq
and 
\bseq\label{eq:dq_bound}
\begin{align}
  |r \pa_r (Q^r / r) | & \les r \CCa(\ww, \xx),  \label{eq:dq_bound:Qr}
   \\
  |z \pa_z (Q^z/z) | &\les |\xx| \CCa(\ww, \xx),  \label{eq:dq_bound:Qz}
\end{align}
\eseq
and 
\bseq\label{eq:Q_bound}
\begin{align}
 |q^z | & \les  1 +  \CCj(\ww, \xx) ,
 \label{eq:Q_bound:a} \\
  |\pa_z Q^z| & \les 1 +  \CCj(\ww, \xx) .
   \label{eq:Q_bound:b}
\end{align}
\eseq

Recall $q^r = \f{1}{r} Q^r, q^z = \f{1}{z} Q^z$ from \eqref{def:qq}. We have the following lower bounds 
\bseq\label{eq:Q_lower}
\begin{align}
  q^z &\geq  ( 2 \bar C_Q - C_1 \cJak \nnrr{\ww} ) \cJaa (\xx) -  C r \CCa(\ww, \xx) , 
 \label{eq:Qz_lower} \\
  q^r & \geq ( 2 \bar C_Q - C_1 \e^{-\kp} \nnrr{\ww} ) \e^{-1}  - C, \label{eq:Qr_lower} \\
  q^r -  q^z & \geq 2  \bar C_Q \e^{-1} \ang \xx^{\al - \hau}
  - C \e^{-1 -\kp} \ang \xx^{\al-\hal + \epa} \nnn{\ww}
    - C \ang \xx^{\al-\hal} ,
\end{align}
\eseq
for some absolute constant $\bar C_{Q}, C, C_1> 0$ and upper bound 
\beq\label{eq:q_upper}
   |q^r| % + |q^z| 
   \les \e^{-1} (1 + \e^{-\kp} \nnn{\ww}).
\eeq

\end{prop}

\begin{proof}

\textbf{\bf Estimates of $\Psi$.}
Recall $\QQ$ from \eqref{def:Q} and $\Qreg^z$ from \eqref{def:LR_decomp}
\beq\label{eq:recall_Q}
 Q^r = c_l r - r \pa_z \Psi, \quad Q^z = c_l z + 2 \Psi  + r \pa_r \Psi  , 
 \quad \Qreg^z \teq c_l z + 2 \Psi. 
\eeq
and $\CCa, \CCj$ from \eqref{def:CC_bounds}. Firstly, using $\Psi =\psi + \bpsi, \Om = \ww + \wwb$, Proposition 
\ref{prop:vel_est} and Corollary \ref{cor:u_bar}, $\nnrr{f} \leq \nnn{f}$ by \eqref{def:3d_norm}, 
and $\rhoo^{-1} \les \ang x^{\epa} \cJak$,  we bound 
\beq\label{eq:tot_u_est1}
\bal
   | \pa_r \tf{1}{z} \Psi| 
   +    |\pa_{r z} \Psi|
   & \les  \min( |\xx|^{\al}, \ang \xx^{\al-\hal -1} ) ( \rhoo(\xx)^{-1}  \nnn{\ww} + 1 )
   \les  \CCa(\ww, \xx) \, . \\
  \eal
\eeq
Recall the notation \eqref{def:psio}. Since $\f{1}{z} \Psi - \pa_z \Psi = \f{1}{z}\mr{\Psi} - \pa_z \mr{\Psi}$, using Proposition \ref{prop:vel_est} and Corollary \ref{cor:u_bar}, we bound 
\beq\label{eq:tot_u_estJ}
\bal
   |\tf{1}{z} ( \mr{\Psi}   + 2 \al & \JJ(\Om) )  | 
  + |\pa_z ( \mr{\Psi} + 2 \al \JJ(\Om) )| + |\tf{1}{z} \Psi - \pa_z \Psi  |  \\
 \qquad & \les \min( |\xx|^{\al+1}, \ang \xx^{\al- \hal} ) ( \rhoo^{-1}(\xx) \nnrr{\ww} + 1 )
  \les |\xx| \CCa(\ww, \xx), \\
   |\JJ(\Om) | & \les \min( |\xx|^{\al+1}, \cJaa(\xx) ) ( \cJak \nnrr{\ww} + 1)
\les \CCj(\ww, \xx).
\eal
\eeq

\paragraph{\bf Decomposition of $\QQ, \qq, \Qreg^z$}

Using the formula \eqref{def:qq}, $\pa_z \Psi(0) = 2 \al \JJ( \Om)(\infty)$ \eqref{eq:psi_JJ}, 
the normalization conditions \eqref{eq:normal_cond},  and $\Om = \om + \wwb$, we extract the main $\JJ$-terms 
\beq\label{eq:q_recall}
\bal
   q^z  &= \tf{1}{z} Q^z = c_l + \tf{2}{z}  \Psi + r \pa_r \tf{\Psi}{z}
  = 2 + \tf{2}{z}  \mr{\Psi} + r \pa_r \tf{\Psi}{z},   \\
  q^r &= c_l -  \Psi_z
    = 2 - 2 \Psi_z(0) - (\Psi_z - \Psi_z(0) ) - \Psi_z(0)  = 2 - 6 \al \JJ(\Om)(\infty)  -  \mr{\Psi}_z \\
    & = 2 - 4 \al \JJ(\Om)(\infty) 
  - ( \mr{ \Psi  }_z + 2 \al \JJ(\Om)(\xx))
  + 2\al ( \JJ(\Om)(\xx) - \JJ(\Om)(\infty)  ) ,
  \\
  q^r - q^z &=  - 6 \al \JJ(\Om)(\infty) -  \mr{ \Psi  }_z
  - \tf{2}{z}  \mr{\Psi} - r \pa_r \tf{\Psi}{z} \\
  & = 6 \al ( \JJ(\Om)(\xx) -
   \JJ(\Om)(\infty)    )
   - ( \mr{\Psi}_z + 2 \al \JJ(\Om)(\xx) )
   - ( \tf{2}{z}  \mr{\Psi} + 4 \al \JJ(\Om)(\xx) ) - r \pa_r \tf{\Psi}{z} .
\eal
\eeq

\paragraph{\bf Proof of \eqref{eq:Q_mix_est} and \eqref{eq:dq_bound} }

Recall  $Q^r$ from \eqref{eq:recall_Q}. Applying \eqref{eq:vel_est}, \eqref{eq:vel_bd_opt}, and using \eqref{eq:rag_wg}, $ \rag^{-1} \les \ang \xx^{\kag} \ang z^{-\kag}$
and $\kag + \hal<1$ \eqref{para:range}, we prove \eqref{eq:Q_mix_est:Qr}:
\[
\bal
    |\pa_z  Q^r | & \les r | \pa_{zz} \Psi| 
   \les  r \min( |\xx|^{\al},  \rag^{-1} \ang \xx^{\al-1} \ang z^{-\hal} ) 
( \rhoo^{-1}   \nnn{\ww} + 1) \\
 & \les r  \min( |\xx|^{\al}, \ang z^{-1}  \ang \xx^{\al - \hal } ) ( \rhoo^{-1} \nnn{\ww}+1)  \\
 & \les r \ang z^{-1}\ang \xx \min( |\xx|^{\al},  \ang \xx^{\al - \hal-1 } ) ( \rhoo^{-1} \nnn{\ww}+1) \les  r \ang z^{-1} \ang \xx  \CCa(\ww, \xx).
\eal
\]
For $\pa_r \Qreg^z$, using 
$\Qreg^z = c_l z + 2 \Psi$ \eqref{eq:recall_Q} and \eqref{eq:tot_u_est1},  we prove \eqref{eq:Q_mix_est:Qz}:
\[
   |\pa_r \Qreg^z|  \les |\pa_r \Psi| \les |z| \CCa(\ww, \xx).
\]

Recall $Q^z / z = c_l  + 2 \Psi / z + \f{r \pa_r \Psi}{z}$ 
and $Q^r /r= c_l - \pa_z \Psi $ from \eqref{def:qq}. 
Using \eqref{eq:tot_u_est1}, we prove \eqref{eq:dq_bound:Qr}
\[
\bal
    |r \pa_r (Q^r / r) | & \les |r \pa_{r z} \Psi| \les r \CCa(\ww, \xx).
\eal
\]

Using \eqref{eq:tot_u_est1}, \eqref{eq:tot_u_estJ}, $\nnrr{f} \leq \nnn{f}$ by \eqref{def:3d_norm}, 
and $r \leq |\xx|$,  we prove \eqref{eq:dq_bound:Qz}
\[
\bal
     |z \pa_z ( Q^z / z)| & \les  |z \pa_z \f{\Psi}{z} | + | r z \pa_z ( \f{\pa_r \Psi}{z} ) |
   \les |  \pa_z \Psi - \f{1}{z} \Psi |
  + | r \pa_{r z} \Psi| + |r \pa_r \f{\Psi}{z}| \\
   & \les |\xx| \CCa(\ww, \xx) + r \CCa(\ww, \xx) \les |\xx| \CCa(\ww, \xx).  \\
  \eal
\]

\vspace{0.05in}

\paragraph{\bf Proof of \eqref{eq:Q_bound} on $Q^z, q^z$}

Using  \eqref{eq:tot_u_estJ}, and \eqref{eq:CC_bounds}, we bound 
\beq\label{eq:Q_bound:pf}
\bal
  |    2 + \tf{2}{z}  \mr{\Psi} | 
  & \les 1 +  |\xx| \CCa(\ww, \xx) + \CCj(\ww, \xx)  \les 1 + \CCj(\ww, \xx).
\eal
\eeq

Recall $q^z$ from \eqref{eq:q_recall}. Using the above estimate, \eqref{eq:tot_u_est1} on $  \pa_r \f{\Psi}{z}$, 
and  \eqref{eq:CC_bounds}, we prove \eqref{eq:Q_bound:a}:
\[
  |q^z| \leq   |    2 + \tf{2}{z}  \mr{\Psi} | +| \tf{1}{z} r \pa_r \Psi|
  \les  1 + \CCj(\ww, \xx) + r \CCa(\ww, \xx)
  \les 1 + \CCj(\ww, \xx)  . 
\]

 Combining the above estimate on $q^z = Q^z / z$ and \eqref{eq:dq_bound}, and using \eqref{eq:CC_bounds}, we prove the estimate on $\pa_z Q^z$ in \eqref{eq:Q_bound:b}:
 \[
   |\pa_z Q^z | \leq 
   | z \pa_z (Q^z/z)|  + | Q^z / z | 
   \les 1 + \CCj(\ww, \xx) + |\xx| \CCa(\ww, \xx) \les 1 + \CCj(\ww, \xx) .
 \]

\vspace{0.1in}

\paragraph{\bf Lower bound of $q^z$ }

Next, we prove the lower bound of $q^r, q^z, q^r - q^z$ in \eqref{eq:Q_lower}.  

 We first estimate $q^z$. Recall $q^z$ from  \eqref{eq:q_recall}. We have 
\beq\label{eq:Q_lower:decomp1}
 q^z  =   2 + \f{2}{z}  \mr{\Psi} +\f{1}{z} r \pa_r \Psi,
 \quad   2 + \f{2}{z}  \mr{\Psi} = 2  - 4 \al \JJ(\Om) 
  + 2 ( \f{1}{z} \mr{\Psi}(\Om) + 2 \al \JJ(\Om) ) .
  \eeq

If $r \geq z$, using \eqref{eq:tot_u_estJ} and $\Om = \om + \wwb$,
 we obtain
\[
\bal
  2 + \f{2}{z}  \mr{\Psi} &    \geq 2 - 4 \al \JJ(\wwb + \om)  
  - C  \CCa(\ww,\xx)\cdot |\xx| . \\ 
\eal
\]

 Using  $r \asymp |\xx|$, \eqref{eq:vel_J} and \eqref{eq:JJ_bound:ne}, we bound 
\[
\bal
  2 + \f{2}{z}  \mr{\Psi}
    &    \geq 2- 4 \al \JJ(\wwb) -  C \cJakkk   \nnrr{\ww} 
  - C  r \CCa(\ww,\xx) \\
 & \geq \bar C_1 \cJaa  - C \cJakkk   \nnrr{\ww}  - C r  \CCa(\ww,\xx).
\eal
\]

Combining the above estimates and estimate \eqref{eq:tot_u_est1} on $r \cdot \f{1}{z}  \pa_r \Psi$, we prove \eqref{eq:Qz_lower} when $r > z$.

If $r < z$, using \eqref{eq:tot_u_est1} and $|z| \asymp |(r, z)| = |\xx|$,  for some $\xi \in [0, r]$,  we obtain
\[
    \f{2}{z}  \mr{\Psi}(r, z)
 =     \f{2}{z}  ( \mr{\Psi}(r, z) - \mr{\Psi}(0, z) )
 + \f{2}{z}\mr{\Psi}(0, z) 
 \geq \f{2}{z}\mr{\Psi}(0, z)  - C  \f{r}{|z|} |\pa_r \Psi(\xi, z)| 
 \geq \f{2}{z}\mr{\Psi}(0, z) - C r \CCa(\ww,\xx).
\]

Since $\mr{\Psi} = \mr{\bpsi} + \mr{\psi}$, using Proposition \ref{prop:vel_est}, $\cJaa \geq 1$, 
and $\rhoo^{-1} \les \ang \xx^{\epa} \cJak$, we obtain 
\[
  |\f{2}{z} \psio(0, z)| \les 
  ( \ang z^{\al - \hal + \epa} + \cJaa(z) ) |\cJaa(z)|^{\kp} \nnrr{\ww}
  \les   |\cJaa(z)|^{\kp+1}  \nnrr{\ww}.
\]

Using Lemma \ref{lem:vel_bc} and $\wwb = \waa$, we obtain $\f{1}{z} \mr{\bpsi}(0, z) = \mr{\psi}_{1D}(\waa)(0, z)$, which along with \eqref{eq:vaa_low} and $|z|\asymp |\xx|$ imply
\[
   2 +  \f{2}{z} \mr{\bpsi}(0, z)
  = 2 +\f{2}{z} \mr{\psi}_{1D}(\waa)(0, z) 
  = \f{1}{z} \vaa \gtr \cJaa(z) 
  \gtr \cJaa(\xx).
\]

Combining the above estimates and estimate \eqref{eq:tot_u_est1} on $r \cdot \f{1}{z} \pa_r \Psi$,  we prove \eqref{eq:Qz_lower} when $r < z$.

\vspace{0.1in}

\paragraph{\bf Lower bound of $q^r$}
For $q^r$ in \eqref{eq:q_recall}, since $\cJaa \gtr 1$, using 
\eqref{eq:tot_u_estJ}, $\al - \hal + \epa<0$, 
and $|\xx| \CCa(\ww, \xx) \les \e^{-\kp} \nnn{\ww}  + 1 $ 
by \eqref{def:CC_bounds} and $\cJaa \les \e^{-1}$, we obtain
\[
\bal
  q^r & \geq 2 - 4 \al \JJ(\Om)(\infty) 
  + 2\al ( \JJ(\Om)(\xx) - \JJ(\Om)(\infty)  ) - |\xx| \CCa(\ww, \xx) \\
& \geq  2 - 4 \al \JJ(\om + \wwb)(\infty) 
  + 2\al ( \JJ(\om + \wwb)(\xx) - \JJ(\om + \wwb)(\infty)  ) - C (\e^{-\kp}\nnn{\ww} + 1) .
\eal 
\]
Using  \eqref{eq:vel_J}, \eqref{eq:JJ_bound}, $- \JJ(\wwb)(\infty) \gtr \e^{-1},
\JJ(\wwb)(\xx) - \JJ(\wwb)(\infty) \geq 0$ by \eqref{eq:JJ_bound}, 
$ \cJaa \les \e^{-1}$, we bound
\[
\bal
  q^r 
  & \geq C_2 \e^{-1} - C  \e^{-1- \kp}   \nnn{\ww}
   -  C ( \e^{-\kp} \nnn{\ww} + 1) 
   \geq C_2 \e^{-1 } -  C \e^{-1 - \kp}   \nnn{\ww}- C \, ,
\eal
\]
and prove \eqref{eq:Qr_lower}.

\vspace{0.1in}

\paragraph{\bf Lower bound of $q^r - q^z$}
For $q^r - q^z$ in \eqref{eq:q_recall}, similarly, 
using \eqref{eq:tot_u_est1}, \eqref{eq:tot_u_estJ}, and the definition of $\CCa(\ww,\xx)$ \eqref{def:CC_bounds},
 we estimate 
\[
\bal
  q^r - q^z  & \geq 
  6 \al ( \JJ( \Om )(\xx) -
   \JJ(\Om)(\infty)    ) -C  |\xx| \CCa(\ww, \xx) \\
& \geq 
  6 \al ( \JJ( \ww + \wwb )(\xx) -
   \JJ( \ww + \wwb )(\infty)    ) 
   -C  |\xx| \CCa(\ww, \xx). \\
\eal
\]
Using \eqref{eq:vel_J}, and \eqref{eq:JJ_bound}, we yield 
\[
\bal
  q^r - q^z   \geq C_3 \e^{-1} \ang \xx^{\al -\hau}
- C \e^{-1 - \kp} \ang \xx^{\al-\hal + \epa} \nnrr{\ww}
- C r \CCa(\ww, \xx) .
\eal
\]

Using the upper bound of $\CCa$ from \eqref{eq:CC_bounds}, 
$\cJaa \les \e^{-1}$, $r \leq |\xx|$, and $\nnrr{\ww} \leq \nnn{\ww}$, we simplify the above bound as 
\[
    q^r - q^z   \geq C_3 \e^{-1  } \ang \xx^{\al -\hau}
    - C \e^{-1 -\kp} \ang \xx^{\al-\hal + \epa} \nnn{\ww}
    - C \ang \xx^{\al-\hal}.
\]
We prove \eqref{eq:Q_lower}.

\paragraph{\bf Upper bound of $|q^r|$}

Using \eqref{eq:q_recall} on $q^r$,  estimates \eqref{eq:tot_u_estJ} on $\Psi$,
 \eqref{eq:vel_J}, \eqref{eq:u_bar:J} on $\JJ$, 
and $\cJaa \les \e^{-1}$, we bound 
\[
  |q^r|  \les 1 +  |\JJ(\Om)(\infty) |
+ |  \mr{ \Psi  }_z + 2 \al \JJ(\Om)(\xx) |
  + | \JJ(\Om)(\xx) - \JJ(\Om)(\infty) |
\les \e^{-1} (1 + \e^{-\kp} \nnn{\ww}).
\]
We prove \eqref{eq:q_upper} and complete the proof.
\end{proof}

\subsection{Estimate of residual error }\label{sec:err}

Recall the residual error from \eqref{eq:fix_recall}.
Since the approximate profile  $\wwb$ \eqref{def:3d_Om} is constant in $r$, using $\nnr{\wwb}= 0$,
Lemma \ref{lem:vel_bc}, the definitions of $Q^z$ \eqref{def:Q}, $V$ \eqref{eq:1D_normal}, and $\vaa$ 
in Theorem \ref{thm:1D_profile_prop}, we obtain 
\bseq\label{eq:iden_r0}
\beq\label{eq:iden_r0:a}
   \bar \Psi(0, z) = \psid(\waa)(z) , 
   \quad \bar Q^z(0, z) = \bar c_l z + \bar U^z = \bar c_l z + 2 \bar \Psi
   =  \vaa .
\eeq

Since $(\waa, \psid(\waa), \bcl, \bcw)$ solve the 1D-profile equation \eqref{eq:1D_dyn} \emph{exactly}, we obtain %and using the above identities, we obtain 
\beq%\label{eq:iden_r0:b}
  \bar \cR(0, z) = \bcw - (1-\al) \pa_z \psid(\waa) - (\bcl z  + 2 \psid) \tf{\pa_z \waa}{\waa}
  = 0.
\eeq
\eseq
Using the above identities, $\wwb = \waa$,  and the formula of $ \bar U^z$ \eqref{def:3d_vel}, we rewrite $\bar \cR$ \eqref{eq:fix_recall} as 
\beq\label{eq:err_recall}
\bal
  \bar \cR   =   \bar \cR - \bar \cR(0, z) &= - (1-\al) (\bar \Psi_z - \bar \Psi_z(0, z)) 
 - ( \bar U^z(r, z) - \bar U^z(0, z) ) \tf{\pa_z \wwb}{\wwb}  \\
 & =  - (1-\al) (\bar \Psi_z - \bar \Psi_z(0, z)) 
 - ( 2 ( \bar \Psi(r, z) - \bar \Psi (0, z) ) 
+ r \pa_r \bar \Psi(r, z) ) \tf{\pa_z \wwb}{\wwb}  .
 \eal
\eeq

In the following estimates for $\bar\cR$, we show that each term either has a good sign, has size $O(1)$, or decays faster in $z$. Since $\bar \cR$ is odd in $z$, we focus on $z \geq 0$.

\begin{prop}\label{prop:3D_error}
 For $z \geq 0$, we have the following upper bound for the residual error $\bar \cR$ 
\footnote{
We \emph{do not} take absolute values in the estimate \eqref{eq:3D_err:linf}. This one-sided estimate is useful in controlling the damping coefficients; see Lemma \ref{lem:3D_damp}.
}
\bseq\label{eq:3D_err_est}
\beq\label{eq:3D_err:linf}
\bal
  \bar \cR 
  &
 \leq \bar C_{\cR} \min\B( |\xx|^{\al+1}, 
 \one_{r > z} \ang z^{\al-\hal} 
 + r \ang \xx^{\al-\hal - 1}
+  \e^{1-\kp} + \f{|z|}{\ang z} \cJaa(\xx) |\lgp z|^{\kp-2} \B)   ,
  \eal
\eeq
with some absolute constant $\bar C_{\cR} = C(\wwwa) $, and the following estimates of the \emph{absolute value} of $\bar \cR$ 
\beq\label{eq:3D_err:abs}
  |\bar \cR| \les |\xx|^{1+\al}.
\eeq

The $z$-derivative satisfies 
\beq\label{eq:3D_err:Cz}
    | z \pa_z \bar \cR |  \les \f{r}{\ang \xx} \min\B( |\xx|^{\al}, 
\ang z^{\al - \hal } + 
 \min\B(  \ang z^{\al- \hal}   \log \f{\ang \xx}{\ang z} , \cJaa(\xx) \B) \cdot \cJaa(z)^{-1}  \B)  .
\eeq

Finally,  the $r$-derivative of the regular part $\Rreg$ \eqref{def:LR_decomp} satisfies
\beq\label{eq:3D_err:Cr}
|\pa_r \Rreg| \les \min( |\xx|^{\al}, \la \xx \ra^{-1 + \al - \hal}) .
\eeq
\eseq
\end{prop}

\begin{remark}[\bf Estimate of $\bar \cR$]

In the following proof of (one-sided) upper bound \eqref{eq:3D_err:linf}, we use the sign condition $
\wwb(r, z) = \waa(z) \leq 0$ for $z \geq 0$ \eqref{eq:waa_sign}
and \eqref{eq:1D_prof:low} to estimate the $I_{22}$ term. The estimates of other terms in the proof 
of Proposition \ref{prop:3D_error} are perturbative. 
The main difficulty in estimating $\bar \cR$ lies in the region $r \gg |z|$, since the profile $\wwb$ does not decay in $r$. 
In this region, the error has a logarithmic growth  $\cJaa(\xx) |\lgp z|^{\kp-2} $. %for $|\xx| \gg |z|$.  
We decompose the upper bound of $\bar \cR$ carefully  and use the damping term generated from the weight and the transport terms to control it in later energy estimates. 
See Lemma \ref{lem:3D_damp}.
\end{remark}

\begin{proof}

Below, we consider $z \geq 0$. 

Firstly, we rewrite $\bar \cR$ \eqref{eq:err_recall} as follows 
\[
\bal
  \bar \cR  & =  - (1-\al) \B( (\bar \Psi_z(r,z)  - \tf{1}{z} \bpsi(r,z) ) - (\bar \Psi_z(0, z) 
  - \tf{1}{z} \bpsi(0,z) ) \B) 
 - \f{1}{z} (  \bar \Psi(r, z) - \bar \Psi (0, z) ) \B(  \f{ 2  z \pa_z \wwb}{\wwb} +  1 -\al \B) \\
 & \quad  - r \pa_r \bar \Psi(r, z)  \tf{\pa_z \wwb}{\wwb}   \teq I_1 + I_2 + I_3.
\eal
\]

Applying $\wwb = \waa$ \eqref{def:3d_Om},  estimates \eqref{eq:1D_prof_est2:c}, 
$\e = \f13 - \al$,  and  $| \f{1 - \al}{2} - \f{1}{3}  |\les \e$, we bound 
\beq\label{eq:3D_err:pf1}
  |\f{z \pa_z \wwb }{\wwb}| \les 1, 
  \quad |\f{ z \pa_z \wwb }{\wwb} +  \f{ 1 - \al }{2} | \les \e + |\lgp z|^{ \kp - 2 }.
  \quad 
\eeq

For $r > |z|$, using Corollary \ref{cor:u_bar}, we obtain
  \[
  |I_1|  \les \min( |\xx|^{\al+1}, \ang z^{\al-\hal}  ) .
  \]

For $r \leq z$, applying the bound for $\pa_{rz} \bpsi, \pa_r \bpsi$ in  Corollary \ref{cor:u_bar},
 the mean-value theorem, and using $|z| \asymp |\xx|$, for some $\xi_1, \xi_2 \in [0, r]$, , we obtain
\[
  |I_1 |  \les |r \pa_{rz} \bpsi(\xi_1, z) |+ |\f{r}{z} \pa_r \bpsi(\xi_2, z)|
  \les r \min( |\xx|^{\al}, \ang \xx^{\al-\hal-1} ).
\]

Combining the above estimates of $I_1$, we obtain
\[
\bal
  |I_1|& \les  \min( |\xx|^{\al+1}, r \ang \xx^{\al-\hal-1} + \one_{r > z} \ang z^{\al-\hal} ) .
\eal
\]

Applying  Corollary \ref{cor:u_bar},  we obtain 
\[
\bal
  |I_3| & \les r |z| \min( |\xx|^{\al}, \ang \xx^{\al-\hal - 1} ) \cdot |z|^{-1}
  \les r \min( |\xx|^{\al}, \ang \xx^{\al-\hal - 1 } ).
\eal
\]

Using the notation $\psio = \psi -\pa_z \psi(0) z $ \eqref{def:psio}, $\cJ_{3D}$, and $\xx = (r, z)$,
we further decompose $I_2$ as 
  \[
  \bal
      I_2 &= - \f{2}{z} (  \mr{\bar \Psi}(r, z) - \mr{ \bar \Psi} (0, z) ) \B(  \f{  z \pa_z \wwb}{\wwb} +  \f{1 -\al }{2}\B) \\
      & =  - 2 \B(  ( \f{1}{z}\mr{\bar \Psi}(\xx) + 2 \al \JJ(\wwb)(\xx) ) - ( \f{1}{z}\mr{ \bar \Psi} (0, z) + 2 \al \JJ(\wwb)(0, z) ) \B) \B(  \f{  z \pa_z \wwb}{\wwb} +  \f{1 -\al }{2}\B)    \\
      & \quad + 4 \al ( \JJ(\wwb)(\xx) - \JJ(\wwb)(0, z) ) \cdot \B(  \f{  z \pa_z \wwb}{\wwb} +  \f{1 -\al }{2}\B) \teq I_{21} + I_{22}.
       \\
    \eal
  \]

If $r \leq |z|$, using $\pa_r \mr{\bpsi}(s, z) = \pa_r \bpsi(s, z)$, Corollary \ref{cor:u_bar}, $|z| \asymp |\xx|$, and \eqref{eq:3D_err:pf1} 
we bound 
\[
  |I_2| \les  \sup_{s \in [0, r]} \f{r}{|z| } |\pa_r \bpsi(s, z) |  
  \les  r \min( |\xx|^{\al}, \ang \xx^{\al-\hal-1} ) .
\]

For $r > |z|$, we bound $I_{21}$ using  \eqref{eq:u_bar}, 
and \eqref{eq:3D_err:pf1} 
\[
      |I_{21}|  \les \one_{r > z} \min( |\xx|^{1+\al}, \ang z^{\al-\hal} ) .  
\]

For $I_{22}$, we shall use the sign of $\JJ$. Recall the definition of $\JJ$ from \eqref{def:JJ_3D}
and $\xx = (r, z)$. Since $\JJ(\wwb)(\xx)$ is odd in $z$, we assume $z \geq 0$. 
Since $|\xx| \geq |z|$, we have
\[
  \JJ(\wwb)(\xx) - \JJ(\wwb)(0, z)
  = C \int_{ z  \leq \, \tz \, \leq |\xx| , \R^2_{++}} \f{\tr^{2 + \al} \tz }{ |(\tr, \tz)|^{5} } 
  \wwb(\tz) d \tr d \tz \,  ,
\]
for some absolute constant $C >0$. Since $\wwb(\tz) \leq 0$ for $\tz \geq 0$, from the above formula, we obtain 
\beq\label{eq:JJ_ran}
  \JJ(\wwb)(\xx)\leq  \JJ(\wwb)(\xx) - \JJ(\wwb)(0, z) \leq 0.
\eeq

Using the lower bound of $\f{z \pa_z \bar \Om}{\bar \Om} + \f{1-\al}{2}$ in \eqref{eq:1D_prof:low} 
and the above sign,  we obtain
\[
  I_{22} \leq - 4 \al C_{\wwwa} (    \JJ(\wwb)(\xx) - \JJ(\wwb)(0, z)) \f{|z|}{\ang z} ( \e^{2 -\kp} + |\lgp z|^{\kp-2} ).
\]

Using \eqref{eq:JJ_ran} and then Corollary \ref{eq:u_bar:J} and \eqref{eq:Ja_hat},  we obtain
\[
  |\JJ(\wwb)(\xx) - \JJ(\wwb)(0, z)| 
  \leq  |\JJ(\wwb)(\xx) | \les \min( |\xx|^{1 + \al}, \cJaa(\xx) ).
\]

Since $\cJaa \les \e^{-1}$, it follows 
\[
\bal
  I_{22} & \leq C_{\wwwa, 2} \min( |\xx|^{1+\al}, \tf{|z|}{\ang z} \cJaa(\xx) ) ( \e^{2-\kp} + |\lgp z|^{\kp-2} ) \\
  & \leq C_{\wwwa, 3} \min\B( |\xx|^{1+\al}, \, \f{|z|}{\ang z} \cJaa(\xx)   |\lgp z|^{\kp-2} + \e^{1-\kp} \B),
  \eal
\]
for some absolute constant depending only on the profile $\wwwa$. 
Combining the above estimates, we prove \eqref{eq:3D_err:linf}. Moreover, using \eqref{eq:u_bar} and \eqref{eq:3D_err:pf1}, we obtain
\[
  |I_{22}| \les  |\xx|^{1+\al}.
\]
Combining the above estimates of $I_{21}, I_1,I_2, I_3$, we prove \eqref{eq:3D_err:abs}.

\vs{0.1in}

\paragraph{\bf Estimate of $z \pa_z \bar \cR$}

Taking $z \pa_z$ on \eqref{eq:err_recall}, we obtain
\[
\bal
  z \pa_z \bar \cR & =  
   - (1-\al) z \pa_z (\bar \Psi_z - \bar \Psi_z(0, z)) 
   - z r \pa_r \pa_z \bpsi \cdot \tf{\pa_z \wwb}{\wwb}
   - r \pa_r \bpsi \cdot z \pa_z \tf{\pa_z \wwb}{\wwb} \\
  & \quad - 2 z \pa_z (\f{ \bpsi(r, z ) - \bpsi(0,z) }{z}  ) \cdot \tf{z \pa_z \wwb}{\wwb}
  - 2  (\f{ \bpsi(r, z ) - \bpsi(0,z) }{z}  ) \cdot z \pa_z ( \f{z \pa_z \wwb}{\wwb} ) .
\eal 
\]

Applying \eqref{eq:1D_prof_est2:c} and 
\eqref{eq:1D_prof_dxx} for $z \pa_z ( \f{z \pa_z \wwb}{\wwb}), z \pa_z ( \f{\pa_z \wwb}{\wwb})$, we obtain
\beq\label{eq:3D_err:pf3}
    | z \pa_z \bar \cR |
    \les |z \pa_z (\bar \Psi_z - \bar \Psi_z(0, z)) |
    + |r \pa_{rz} \bpsi| + |\f{r}{z} \pa_r \bpsi|
    + \B| z \pa_z (\f{ \bpsi(r, z ) - \bpsi(0,z) }{z}  ) \B|
    + |\f{\bpsi(r, z ) - \bpsi(0,z) }{z} | \cdot \cJaa(z)^{-1} .
\eeq

\vs{0.1in}
\paragraph{\bf Case $|z| \geq r$}

For $|z| > r$, we get $|z| \asymp |\xx|$. Using mean-value theorem and \eqref{eq:u_bar}, 
\eqref{eq:u_bar_rzz}, we bound 
\[
\bal
  | z \pa_z \bar \cR| 
  & \les \sup_{\tr \in [0, r]} |r z \pa_{r zz} \bpsi(\tr, z)|
   + \sup_{\tr \in [0, r]} r |\pa_{r z} \bpsi(\tr,  z) | 
   + \sup_{\tr \in [0, r]} \tf{r}{z} | \pa_{r } \bpsi(\tr,  z) |  \\
  & \les  
  z r \min( |\xx|^{\al-1}, \ang \xx^{\al-\hal-2} ) 
  + r \min( |\xx|^{\al}, \ang \xx^{\al-\hal-1} ) 
\les \f{r}{\ang \xx} \min( |\xx|^{\al}, \ang \xx^{\al-\hal} ) .
\eal
\]
Since $|z| \asymp |\xx|$, we prove \eqref{eq:3D_err:Cz} when $|z| > r$.

\vs{0.1in}
\paragraph{\bf Case $|z| < r$}

In this case, we obtain $r \asymp |\xx|$. Applying Corollary \ref{cor:u_bar} and $|z| \leq |\xx|$, we obtain
\begin{align}\label{eq:3D_err:pf4}
 \B| z \pa_z (\f{ \bpsi(r, z ) - \bpsi(0,z) }{z}  ) \B|
 &\les | \pa_z \bpsi(r,z) - \tf{1}{z} \bpsi(r,z) |
 + | \pa_z \bpsi(0,z) - \tf{1}{z} \bpsi(0,z) | 
 \les \min( |\xx|^{\al+1} , \ang z^{\al-\hal}  ), \notag \\
|z\pa_z (\bar \Psi_z - \bar \Psi_z(0, z)) | 
& \les z ( \min( |\xx|^{\al} , \ang \xx^{\al-1} \ang z^{-\hal} )  
+\min( |z|^{\al} , \ang z^{\al-\hal-1} )    ) 
\les  \min( |\xx|^{\al+1}, \ang z^{\al-\hal} ), \notag \\
     |r \pa_{rz} \bpsi| + |\tf{r}{z} \pa_r \bpsi| & 
    \les r \min( |\xx|^{\al}, |\xx|^{\al-\hal-1} )
    \les \min( |\xx|^{\al+1}, \ang \xx^{\al-\hal} ).
\end{align}

Applying \eqref{eq:u_bar:J} in Corollary \ref{cor:u_bar} and $\bpsi(r, z) -\bpsi(0, z) = \mr{\bpsi}(r, z) -
\mr{ \bpsi}(0, z)$, we bound 
\[
    |\f{\bpsi(r, z ) - \bpsi(0,z) }{z} | 
  \les \min( |\xx|^{1+\al}, \ang z^{\al-\hal} 
  +  |\JJ(\wwb)(r, z) - \JJ(\wwb)(0,z)|  ).
\]

Using the definition of $\JJ$ \eqref{def:JJ_3D},
$|\wwb| \les \min( |z|, |z|^{-\hal} )$, and integrating in $\tr$,  we bound 
\[
\bal
  I & \teq |\JJ(\wwb)(r, z) - \JJ(\wwb)(0,z)| \\
  & \les \int_{ |z| \leq  \tz \leq |\xx|  } \f{ \tr^{2 + \al} \tz }{ (\tr^2 +\tz^2)^{5/2} } 
  \min( \tz , \tz^{-\hal} ) d \tz d \tr  \les  \int_{|z|}^{ |\xx| } \tz^{\al-1} \min( \tz , \tz^{-\hal}) d \tz   \, .
  \eal
\]

Since the integrand is bounded by $1$ and $| |a| - \ang a | \les \ang a^{-1}$ for any $a \geq 0$, we further bound  
\[
  I 
  \les \ang z^{-1} + \int_{\ang z}^{ \ang \xx }   \tz^{\al-\hal-1} d R
  \les \ang z^{-1} + \ang z^{\al -\hal} \log \f{\ang \xx}{\ang z}.
\]

Since $\cJaa(\xx) \geq \cJaa(z)$, using Corollary \ref{cor:u_bar}, we obtain two estimates 
\[
      |\f{\bpsi(r, z ) - \bpsi(0,z) }{z} | 
      \les \min( |\xx|^{\al+1}, 
      \ang z^{\al -\hal} + \min( \cJaa(\xx), \ang z^{\al-\hal} \log \f{\ang \xx}{\ang z} ) ).
    \]

For $ r > |z|$, since $\f{r}{\ang \xx} \asymp \f{|\xx|}{\ang \xx}$, combining the above estimates, \eqref{eq:3D_err:pf3}, \eqref{eq:3D_err:pf4}, $\cJaa(z)^{-1} \les 1$,  we prove \eqref{eq:3D_err:Cz}
when $r > |z|$.

\vs{0.1in}

\paragraph{\bf Estimate of $\pa_r \bar \cR$}

Applying \eqref{eq:err_recall} and \eqref{def:LR_decomp}, we obtain
\[
\pa_r \Rreg = \pa_r \big(  \bar c_{\om} - (1 - \al) \bar \Psi_z -  (\bar c_l z +  2 \bar \Psi   )  \tf{ \pa_z \wwb }{\wwb}   \big)
= - (1-\al) \pa_{r z} \bar \Psi
 -  2  \pa_r \bar  \Psi \cdot  \tf{\pa_z \wwb}{\wwb} .
\]

Applying Corollary \ref{cor:u_bar} and \eqref{eq:1D_prof_est2:c}, we bound 
\[
  | \pa_r \Rreg| \les |\pa_{rz} \bpsi| + | \tf{1}{z} \pa_r \bpsi| 
  \les \min( |\xx|^{\al}, |\xx|^{\al-\hal -1} )
  \les \min( |\xx|^{\al}, \ang \xx^{\al-\hal -1} ). 
\]
We complete the proof.
\end{proof}

\subsection{Qualitative regularity and smallness assumption}\label{sec:qual_reg}

In the next two subsections, we derive \emph{qualitative} regularity estimates for stream function and solution to the fixed point equation \eqref{eq:fix_pt_2D}. We will use these properties to justify the energy estimates and prove compactness in Section \ref{sec:3D_solu}.

Firstly, we impose the following smallness on the input of the map $\eta = \FFF(\ww)$ \eqref{eq:fix_pt_2Deta}:
\bseq\label{eq:3D_size:ww}
\beq
  \nnla{\ww} + \nnr{\ww} + \nnz{\ww} \leq \e^{ \f{1 + \kp}{2} },
\eeq
to simplify later nonlinear estimates. From the definition of norm in \eqref{def:3d_norm}, we obtain
\beq
 \nnrr{\ww} \leq  \nnn{\ww} = \max(\nnla{\ww}, \nnr{\ww}, \nnr{\ww}^{\al} \nnz{\ww}^{1-\al}) \leq \e^{ \f{1 + \kp}{2} } .
\eeq
\eseq

Under the assumption \eqref{eq:3D_size:ww}, we can simplify some estimates in Proposition \ref{prop:Q}. Using $\ang \xx^{\al-\hal + \epa} \les 1$ by \eqref{ran:ep_all} and $\cJaa \les \e^{-1}$, we bound $\CCa, \CCj$ \eqref{def:CC_bounds} as 
\beq\label{eq:3D_size:C}
\bal
 \CCa(\ww, \xx) &\les \min( |\xx|^{\al}, |\xx|^{\al- \hal - 1} ) 
 \cdot ( 1 +  \ang \xx^{\epa} \e^{ - \kp } \cdot \e^{ \hk } ) \\
& \les 
\min( |\xx|^{\al},    \ang \xx^{\al-\hal -1 }(1 + \e^{\mhk} \ang \xx^{\epa} ) )
\les  \ang \xx^{\al-\hal -1} +  \e^{\mhk} \ang \xx^{-1}   , \\
\CCj(\ww, \xx) & \les \min( |\xx|^{\al+1}, \cJaa ) \cdot (1 + \e^{-\kp} \cdot \e^{\kp})
\les \min( |\xx|^{\al+1}, \cJaa ).
\eal
\eeq

Recall $\epa$ from \eqref{ran:ep_all} and $\bbb$ from \eqref{norm:Xc}.  From \eqref{ran:ep_all}, we have 
\beq\label{eq:ep_para1}
  \bbb - 1  < \al, \quad \al - \hal + \epa + \e^2 < - \tf12 \e , 
  \quad \e^2  - \epa \leq - \tf12 \kp_1 \e  < 0.
\eeq

Using the definitions of $\rhoo$ \eqref{def:3d_wg}, $\rhoc$ \eqref{def:rhoc}, 
and \eqref{eq:ep_para1}, we obtain
\beq\label{eq:ep_ineq1}
\bal
\rhoo \rhoc \min( |\xx|^{\al} , \ang \xx^{\al-\hal} ) & \les  \ang \xx^{\al -\hal + \e^2 - \epa} 
\les \ang \xx^{\al-\hal} \les 1, \\
\rhoc \min( |\xx|^{\al} , \ang \xx^{\al-\hal} ) & \les  \ang \xx^{\al -\hal + \e^2  } \les 1 \, . \\
  \eal
\eeq

From \eqref{eq:3D_size:C} and the above estimates, we obtain
\beq\label{eq:3D_size:C2}
 \rhoc \rhoo \cdot \ang \xx \CCa(\ww, \xx)   \les 
\ang \xx^{\e^2 - \epa}
\cdot \ang \xx \cdot  
\ang \xx^{\al -\hal - 1  +\epa}\les \ang \xx^{\al-\hal +  \e^2}.
\eeq

Let $\beps_5$ be as in Theorem \ref{thm:contra_lin_exact}.  We choose 
\beq\label{def:beps_6}
\beps_6 \in (0, \beps_5] \, ,
\eeq
 small enough, so that for any $\e \leq \beps_6 $, using Proposition \ref{prop:Q}, 
\eqref{eq:Q_lower} on the lower bound,  \eqref{eq:q_upper},\eqref{eq:Q_bound:a} on the upper bound, 
and $|\xx| \CCa\les 1 \les \cJaa \les \e^{-1}$,  we obtain \eqref{ran:ep_all} and
\beq\label{eq:3D_boot_q:1}
\bal
 q^z  & \geq  \tf32 \bar C_Q \cJaa(\xx) - C r \min( |\xx|^{\al} , \ang \xx^{\al-\hal-1} ) , \\
   q^r & \geq   \bar C_Q \e^{-1} , \\
  |q^z| & \les \cJaa \les \e^{-1}, \quad  q^r \asymp \e^{-1} .
  \eal
\eeq

\begin{comment}

\subsubsection{Smoothness in Axisymmetric coordinate and in $\R^n$}\label{sec:smooth_coordinate}

Consider a scalar axisymmetric function $G(\xx)$ $\in \R^n$ with $n\geq 3$, namely $G(\xx)= g( |\xx_h|, x_n)$,
 for some $g(r, z)$ defined on $\R_+ \times \R$, where $\xx_h = (x_1,..,x_{d-1})$. 
 By calculus, we have the following simple relation.

 \begin{lem}\label{lem:smooth_coordinate}
 Consider an axisymmetric domain $D_G = \{ \xx: |\xx_h| < a, |x_n| < b \}$ and $D_g = \{ (r, z): r \in [0, a), |z| < b \}$ for some $a, b> 0$. We have $G \in C^{k, \g}(D_G)$ for some integer $k\geq 0$ and $\g \in (0, 1)$ if and only if $ g \in C^{k, \g}(D_g)$ and it satisfies the compatibility conditions $\pa_r^i g(0) = 0$ for all odd $i\leq k$. 
\end{lem}

See also \cite{liu2009characterization}. In the rest of the paper, unless otherwise specified, the notation \(f\in C^{k,\g}\) means that \(f(r,z)\) is a \(C^{k,\g}\) function of \((r,z)\in D\subset \R_+\times\R\). 
To obtain a \(C^{k,\g}\) function in \(\R^n\), we will additionally verify the compatibility conditions.

\end{comment}

\subsubsection{Regularity of the stream function}\label{sec:qual_reg_stream}

In the rest of the paper, unless otherwise specified, \(f\in C^{k,\g}(D)\) means that \(f(r,z)\) is a \(C^{k,\g}\) function of $(r,z)$ in domain $D\subset \R_+ \times\R $, with derivatives $\na_{r, z}^i f, i\leq k$ continuous \emph{up to} $r=0$.  
It is well-known that a $C^{k ,\g}(D)$ function $f$ in $(r, z)$ can be extended as a $C^{k, \g}$ 
function $F$ in $\R^3$ if it satisfies the 
compatibility conditions $\pa_r^i f(0, z) = 0$ for all odd $i\leq k$.

We first establish the following \textit{qualitative} estimates of the stream function.

\begin{lem}\label{lem:psi_reg}
Let $\Om = \wwb + \ww$ and $\Psi = \BS(\Om)$ be the stream function associated with $\Om$. 
Denote $B_R = \{ (r, z):| (r, z)| < R \}$, 
$\QQ = \bar \QQ + \td \QQ(\ww)$ as in \eqref{def:Q}. For any $R > 1$ and $|\xx| < R$, we have : 
\bseq\label{eq:psi_reg}
\begin{align}
    | \Psi(\xx) | + |\na \Psi(\xx)| + | \na_{r, z} \pa_z \Psi(\xx)| 
    + 
       \|  r \Psi  \|_{ C^{2, \alpha}(B_R)} 
       +          \|  \pa_z  \Psi \|_{C^{1, \al}(B_R)} 
    & \les_{R, \e} 1 + \nnn{\ww}  \label{eq:psi_reg:a} , \\ 
   \| \QQ \|_{C^{1, \al}(B_R)} & \les_{R, \e} 1 + \nnn{\ww}. \label{eq:psi_reg:d}
\end{align}
\eseq
\begin{comment}
Moreover, we have 
\beq\label{eq:psi_r0_van}
\pa_r f(0, z) = 0, \quad  r \pa^2  f(0, z) = 0, \quad \mbox{for \ } f= \bar \Psi, \, \psi(\ww), \, \Psi = \bar \Psi + \psi(\ww). 
\eeq
\end{comment}

\end{lem}

The estimates are \emph{qualitative} since they do not quantify 
the dependence on $R, \e$.

\begin{proof}

Using estimate \eqref{eq:CC_bounds}, \eqref{eq:tot_u_estJ}, and \eqref{eq:3D_size:ww}, we obtain 
\[
|\Psi_z(0)| + |\JJ(\Om) | \les \CCj(\ww, \xx) \les_{\e} 1 + \nnn{\ww} .
\]

Using the above estimate, $\Psi =\bpsi + \td \psi$, $ \mr{\Psi} = \Psi - \Psi_z(0) z$, Proposition \ref{prop:vel_est} and Corollary \ref{cor:u_bar}, and $\rhoo^{-1} \rag^{-1} \les_{R, \e} 1$ for $\xx \in \bar B_R$ by \eqref{def:3d_wg},  we prove 
\beq\label{eq:psi_reg:pf1}
  |\Psi|  + |\na \Psi| + | \na \Psi_z| 
  \les_{R, \e} 1 + \nnn{\ww}  , \quad \forall x \in \bar B_R .
\eeq

Since $\Om =\wwb + \ww$, from the definition of the norm \eqref{def:3d_norm} and Lemma \ref{lem:linf_R3}, we obtain
\beq\label{eq:psi_reg:Om}
\bal
   |\Om(\xx)| + |\na \Om(\xx) | & \les C(|\xx|, \e) (1 + \nnn{\ww}) ,
   \quad  C(|\xx|, \e) < \infty, \\
  |\Om(\xx) | 
  & \les_{\e} |\wwb|  (\rag^{-1} \rhoo^{-1} \nnn{\ww} + 1)
  \les_{\e} \ang \xx^{\kag + \epa} (1 + \nnn{\ww}).
\eal 
\eeq

\paragraph{\bf $C^{2,\al}$ estimates for $r \Psi$}

Let $(r, z ,\th)$ be the cylindrical coordinates in $\R^3$. Recall the relation  $ \Phi^{\th}  = \Psi r  ,  \Om^{\th} = \Om r^{\al}  $ from \eqref{eq:Euler_change} and the elliptic equation \eqref{eq:Euler_b}:
\beq\label{eq:elli_R3_recall}
  - \D_{\R^3} ( \Psi r \cos \th  )  = \cpsia \Om r^{\al} \cos \th.
\eeq

Let $\td B_R(\R^3) = \{ \yy \in \R^3: |\yy|  < R\} $ be the ball in $\R^3$
and $B_R = \{ (r, z) \in \R_+ \times \R, | (r, z)| < R \}$. 
Denote $\yy =(r \cos \th, r \sin \th, z)$. 
We can view $\Psi r \cos \th, \Om r^{\al} \cos \th $ as  functions in $\R^3$. Applying Schauder's  estimates \cite[Theorem 4.8]{gilbarg1998elliptic} to $\Psi r \cos \th $, \eqref{eq:psi_reg:a} to $\Psi$,
and \eqref{eq:psi_reg:Om} to $\Om$, we obtain 
\beq\label{eq:psi_reg:schau1}
\| r\Psi  \cos \th \|_{ C^{2, \al}( \td B_{R}(\R^3) } \les_{R, \e}
\| \Psi    \|_{L^{\infty}(\td B_{2R}(R^3))}
+ \| \Om r^{\al} \cos \th  \|_{C^{\al}(\td B_{2R}(R^3))}
\les_{R, \e}
  1 + \nnn{\ww} .
\eeq

For any function $ g(r ,z) r \cos \b \in C^{k, \al}(\td B_R)$, a direct calculation yields 
\[
\pa_r^i \pa_z^j (g \cdot r)  = \pa_{y_1}^i \pa_z^j ( g \cdot r \cos \th ) |_{\th = 0},
\quad \forall i + j \leq k,
\quad r \neq 0 .
\]
Taking $\yy = (r, 0, z)$ with $r \to 0$ in the above identity, 
and using the continuity of the right hand side, we obtain that $ g(r , z) \cdot r \in C^k$.  Thus, we prove: 
\beq\label{eq:psi_reg:pf2}
  \| r \Psi(r, z) \|_{C^{2, \al}(B_R) } \les \| r \Psi \cos \th \|_{ C^{2, \al}( \td B_R(\R^3) )} \les C_{R, \e} (  1 + \nnn{\ww} ).
\eeq

\vspace{0.1in}
\paragraph{\bf $C^{1,\al}$ estimate for $\pa_z \Psi$}

Recall from \eqref{eq:5D_BS} the formula of $\Psi$:
\[
  \Psi(r, z) =  (\KRF \ast W)( r, 0,0,0,z ), 
  \quad W(\yy) \teq |\yy_h|^{\al-1} \Om( |\yy_h|, y_5 ),
  \quad  K_{\R^5}( s ) = \tf{\cpsi}{ 2\pi^2 } |s|^{-3}.
\]

Let $\chi_1(x) \geq 0: \R \to \R_+$ be a $C^{\infty}$-smooth cutoff function with 
\beq\label{def:chi1}
  \chi_1(x) = 1 , \quad |x| \leq 1 , \quad  \chi_1(x) = 0,\quad \forall \ |x| \geq 2 .
\eeq

We fix $R > 1$ and decompose 
\[
   \Om  = \Om \chi_1( \f{\xx}{2 R} ) +  \Om  (1 - \chi_1( \f{\xx}{2 R} ))
   \teq \Om_1 + \Om_2,
   \quad W_i(\yy) = \Om_i( |\yy_h|, y_5)  |\yy_h|^{\al-1} ,
   \quad \Psi_i =  K_{\R^5} \ast W_i.
  \]

From \eqref{eq:5D_BS}, we have  $-\D_{\R^5} \Psi_2 = \cpsia \Om_2(\yy) |\yy_h|^{\al-1} =0 $ for $|\yy| \leq 2 R$. Applying Schauder's  estimates \cite[Theorem 4.8]{gilbarg1998elliptic} and 
then using \eqref{eq:vel_psi}, \eqref{eq:vel_J}, and the norm \eqref{def:3d_norm}, we obtain
 \[
    \| \Psi_2 \|_{ C^{3, \al}(\bar B_R)} \les_R \| \Psi_2 \|_{ L^{\infty}(\bar B_{2 R}) }
    \les_R \| \Om  (1 - \chi_1( \f{\xx}{2 R} )) \|_{\linfa}
    \les_R \| \Om \|_{\linfa}
    \les_R 1 + \| \ww \|_{\linfa}.
 \]

For $\Psi_1$, since $ \pa_z \Psi_1 = \KRF \ast \pa_\tz W_1
= \KRF \ast ( \pa_\tz \Om_1 |\yy_h|^{\al-1} )$ and $\Om_1$ has a compact support, applying  Lemma \ref{lem:hol_singu}
and then \eqref{eq:psi_reg:Om}, we prove 
\[
  \|\na_{r, z} \pa_z \Psi_1(r, z) \|_{\dot C^{\al}}
  \les \sup\nolimits_{|\yy| < 4 R} | \pa_{\tz} \Om_1(\yy) |
  \les_{R, \e}   1 + \nnn{\ww} .
\]
  Combining the above estimates for $\Psi_1, \Psi_2$, and using $\Psi_1 + \Psi_2 = \Psi$, we prove 
  $ \na \pa_z \Psi \in \dot C^{\al}(B_R)$. From \eqref{eq:psi_reg:pf1}, \eqref{eq:psi_reg:pf2}, we obtain $ |\na^2 \Psi(1,0)| \les 1  $. Thus, we prove 
  \beq\label{eq:psi_reg:pf3}
     \| \na \pa_z \Psi \|_{ C^{\al}(\bar B_R)} \les_{R, \e} 1 + \nnn{\ww}.
  \eeq
Combining estimates \eqref{eq:psi_reg:pf1}, \eqref{eq:psi_reg:pf2}, 
\eqref{eq:psi_reg:pf3}, we prove \eqref{eq:psi_reg:a}.

Estimate \eqref{eq:psi_reg:d} follows from \eqref{eq:psi_reg:a}.
\begin{comment}
Denote $\xx = (r, 0,0,0,z)$.
Since $\pa_1 \KRF(s) = C s_1 |s|^{-5}$ is odd in $ s_1$ and $W(\yy)$ is even in $y_1$, we obtain
\[
 \pa_r \Psi(0, z) = \pa_{x_1} (\KRF \ast W)(\xx) |_{x_1 = 0}
 = C \cdot ( \tf{y_1}{|y|^5} \ast W)(\xx) |_{x_1 = 0} = 0.
\]

Using \eqref{eq:Euler2_psi} with $r = 0$, the boundedness of $\om, \pa_{zz} \Psi$, and the above identity, we prove 
\[
 |r \pa_r^2 \Psi(0,z) | 
  \les |\pa_r \Psi(0, z) | + 0 \cdot |\pa_{zz} \Psi(0, z) | + 0 \cdot |\Om(0, z)| = 0.
\]
We prove the identities \eqref{eq:psi_r0_van} for $\Psi$. The identities for $\bar \Psi, \psi(\ww)$ are proved in the same way.
\end{comment}
 \end{proof}

\vs{0.1in}

Recall $\QQ = \bar \QQ + \td \QQ(\ww)$ from \eqref{def:Q}. We define the flow map $X = (X^r, X^z)$ associated with $\QQ$
\beq\label{eq:traj}
    \tf{d}{d s} X(s, \xx )  = \QQ \circ  X(s , \xx ),  \quad  X(0, \xx) = \xx .  \\
\eeq

With Lemma \ref{lem:psi_reg}, we establish the \emph{qualitative} regularity estimate 
of solution $\eta$ to \eqref{eq:fix_pt_2D}.

\begin{prop}[\bf Regularity of the solution]\label{prop:reg_solu}

Let $\beps_6$ be as in \eqref{def:beps_6}. Suppose that $\ww$ satisfies assumption \eqref{eq:3D_size:ww}. There exists $ \beps_7 \in (0, \beps_6]$ such that for any $ \e \leq \beps_7$, the following statement holds.

\begin{enumerate}[label=(\roman*),leftmargin=1.5em]

\item Let $\xx = (r, z)$ and $\QQ$ be as in \eqref{def:Q}. There exists some absolute constants $\lam_Q, C > 0$ such that 
\beq\label{eq:traj_outgo}
\QQ  \cdot \xx =  q^r r^2 + q^z z^2 \geq  C \cJaa(\xx) |\xx|^2 \geq \lam_Q |\xx|^2 \, .
\eeq
For any $\xx \neq 0$, $| X(s,  \xx )|$ of the trajectory $X(s, \xx)$ associated to $\QQ$ 
\eqref{eq:traj} is strictly increasing in $s$. Moreover, for any $s \leq 0$ and any $\xx$, we have
\beq\label{eq:traj_pf1}
    \quad |X(s, \xx )| \leq |\xx| e^{ \lam_Q  s} .
\eeq

For any $s \leq 0, z \neq 0$, and $\xx \in B_R$,  we have $X^z(\xx, s) \neq 0$, and 
\beq\label{eq:traj_ratio}
 \f{X^r(s, \xx)}{|X^z(s, \xx )|} \les_{R, \e} \f{r}{|z|},
 \quad 
   |X^z(s, \xx )|  \asymp_{R, \e} e^{2 s } |z| .
\eeq

\item
There exists a unique continuous solution $\eta(r, z)$ to the equation \eqref{eq:fix_pt_2D}, which is odd in $z$ and satisfies the vanishing order near $\xx = 0$
\beq\label{eq:3D_ODE:init}
 \lim_{|\xx| \to 0 } \f{\eta}{|\xx|} = 0.
\eeq

The solution satisfies $\eta \in C^{1, \al}(B_R)$ 
for the ball $B_R\subset \R_+ \times \R$ with any $R>0$.
 Let $\rhoo, \rag, \rhoc$ be the weights from \eqref{def:3d_wg}.
For any $R > 0$ and $\xx \in B_R$,  $\eta$ satisfies the following estimates 
\bseq
\begin{gather}
 |\eta(\xx) |  \les_{R , \e } |z| \cdot |\xx|^{ \al }, 
 \quad  |\na \eta(\xx) |  \les_{R , \e } |\xx|^{ \al }, \quad 
  \| \eta \|_{C^{1, \al}(B_R)}  \les_{R, \e} 1 ,  \label{eq:local_reg}  \\
   |\rhoc \rhoo \rag   \f{\eta}{\wwb} | +  | \f{z}{\ang z} \rhoc \rhoo \rag \ang \xx \pa_r  \f{\eta}{\wwb} |
  + | \rhoc \rhoo \rag z \pa_z  \f{\eta}{\wwb} | \les_{R , \e} 1  
  \label{eq:local_reg:nnn}. 
\end{gather}

Denote $G = \f{\eta}{\wwb}$. 
For fixed $ z \neq 0$ 
and $i  \in \{ r, z \}$, the function $(\pa_i G ) \cc X(s, \xx)$  is locally $C^1$ in $s$ and satisfies
\beq\label{eq:traj_smooth}
\tf{d}{ds} (\pa_i G) \cc X(s, \xx)
= \big[ \pa_i ( ( \bar \cR + \Lpsi) (G + 1) ) 
  -  (\pa_i  \QQ) \cdot \na G \big ] \cc X(s, \xx)  .
\eeq
\eseq

\end{enumerate}

\end{prop}

If $G$ is locally $C^2$, \eqref{eq:traj_smooth} follows trivially from the chain rule. We will perform \textit{quantitative} estimates on $\eta, \na_{r, z} \eta$ in Section \ref{sec:3D_EE}.

\begin{remark}[\bf Improved regularity]\label{rem:high_reg}

In Proposition \ref{prop:reg_solu}, 
 $\eta$ has weighted $C^{1,\al}$ regularity and is more regular than the input $\ww$, which belongs to a weighted $W^{1.,\infty}$ space. 
We  use this crucial improved regularity to prove that the fixed point map 
$ \FFF$ from $\ww$ to $\eta$ defined in \eqref{eq:fix_pt_2D} 
is compact. This is an important step in applying the Schauder fixed point argument.
See Section \ref{sec:3D_comp}.

\end{remark}

In the rest of this Section \ref{sec:qual_reg}, we  prove Proposition \ref{prop:reg_solu}. 
We first estimate the regularity of the coefficients of the equations 
\eqref{eq:fix_pt_2D}, including $\QQ, \Lpsi, \bar \cR$, in Sections \ref{sec:reg:eqn}, 
\ref{sec:reg:traj_C1a}. We then perform $C^{1,\al}$ estimates on $\eta$ in Sections 
\ref{sec:traj_solu_C1a}, \ref{sec:reg:eta_C1a}.

\subsubsection{Proof of \eqref{eq:traj_outgo} and Item (a) }

Recall the notation $Q^r = q^r r, Q^z = q^z z$. Using \eqref{eq:3D_boot_q:1}, Cauchy-Schwarz inequality, 
 $1 \les \cJaa(\xx) \les \e^{-1}$,  and choosing $\e$ small, we obtain
\[
\bal
  q^r r^2 + q^z z^2 
& \geq \bar C_Q \e^{-1} r^2 + \bar C_Q \cJaa(\xx) z^2 
- C r z^2 \min( |\xx|^{\al}, |\xx|^{\al-\hal - 1} )  \\
&\geq \tf{1}{2} ( \bar C_Q \e^{-1} r^2 + \bar C_Q \cJaa(\xx) z^2 )
+  C_1  r |z| ( \e^{-1/2}- C z \ang \xx^{-1})  \\
& \geq \tf{1}{2} ( \bar C_Q \e^{-1} r^2 + \bar C_Q \cJaa(\xx) z^2 )
\geq C \cJaa(\xx) |\xx|^2 \geq \lam_Q |\xx|^2,  
\eal
\]
for some absolute constant $\lam_Q > 0$.

Using the ODE \eqref{eq:traj} and \eqref{eq:traj_outgo}, for any $\xx$, we obtain
\[
  \tf{1}{2} \cdot \tf{d}{d s} |X(s, \xx)|^2 = \QQ( X(s, \xx) ) \cdot X(s, \xx) \geq \lam_Q |X(s, \xx)|^2 ,  \   \Rightarrow  \ \tf{d}{d s} |X(s, \xx)| \geq \lam_Q |X(s, \xx)|.
\]
which implies \eqref{eq:traj_pf1}.
Thus, for $\xx \neq 0$,  $|X(s, \xx)|$ increases strictly in $s$.

\subsubsection{Equation of $\eta$ }\label{sec:reg:eqn}

Recall the equation of $\eta$ from \eqref{eq:fix_pt_2D} 
\bseq\label{eq:reg:eqn_w}
\beq
    \QQ \cdot \na \eta =    \cB  \cdot \eta + H(\ww) , \quad  \QQ = \bar \QQ + \td \QQ(\ww) .
\eeq
From \eqref{eq:lin_2D},  \eqref{eq:fix_pt_2D} , 
and \eqref{eq:normal_cond}, $\cB, H, \com $ are given by
\beq
\bal
\cB    & =  \bcB + \tcB( \ww) =   c_{\om} - (1 -\al) \Psi_z  ,  
\quad c_{\om} = 2 + (1 -\al) \Psi_z(0) ,  \\
   H(\ww) & = \wwb ( \bar \cR + \Lpsi(\ww) ) 
= ( \com - (1 -\al ) \Psi_z) \wwb - Q^z \pa_z \wwb.
\eal
\eeq
\eseq

Below, we simplify $\Lpsi(\ww)$ as $\Lpsi$. We first estimate the regularity of the coefficients $\wwb (\bar \cR + \Lpsi),\cB $.  Using Lemma \ref{lem:psi_reg},  and $\nnn{\ww} \leq 1$ by \eqref{eq:3D_size:ww}, for any $R >1$,  we obtain
\bseq\label{eq:reg:est_coe}
\beq\label{eq:reg:est_coe:a} 
\bal
   \| \QQ \|_{C^{1,\al}(B_R)}  +
  & \| \wwb (\bar \cR + \Lpsi) \|_{C^{1,\al}(B_R)} 
   + \| \cB  \|_{C^{1,\al}(B_R)}    \\
   &  \les_{\e} 1 + \| \pa_z \Psi \|_{C^{1,\al}(B_R)}  + \| \QQ  \|_{C^{1,\al}(B_R)}  
     \les_{R, \e} 1 + \nnn{\ww}  \les_{R, \e} 1 .
 \eal
\eeq

Using \eqref{eq:Lpsi_linf_est} on $\Lpsi$ and \eqref{eq:3D_err:linf} on $\bar \cR$, we obtain 
\beq\label{eq:reg:est_coe:cR}
 |\wwb( \Lpsi + \bar \cR )|
  \les_{R,  \e } |z| \cdot |\xx|^{1 + \al}(1 + \nnn{\ww}) 
  \les_{R, \e} |z| \cdot |\xx|^{1 + \al}, \quad \forall \  |\xx| \leq R \, .
\eeq

Using \eqref{eq:3D_boot_q:1}, we introduce $Q_{r0}$ as follows and it satisfies 
\beq\label{def:Qr0}
 Q_{r0} \teq \pa_r Q^r(0) = c_l - \pa_z \Psi(0),
 \quad Q_{r0} \asymp \e^{-1}.
\eeq

Using \eqref{eq:normal_cond} and \eqref{def:Q},  we obtain
\[
\bal
\cB   &= 2 + (1-\al)\Psi_z(0) - (1-\al) \Psi_z = 2 - (1-\al) \mr{\Psi}_z,  \\ 
 Q^z  & = ( 2 - 2 \Psi_z(0) ) z + 2 \Psi + r \pa_r \Psi  = 2 z + \mr{\Psi} + \pa_r ( r \mw{\Psi} ), \\
 Q^r & = (c_l - \pa_z \Psi(0) ) r - r \pa_z \mr{\Psi} 
 = Q_{r0} r - r \pa_z \mr{\Psi} . 
 \eal
\]

Using \eqref{eq:tot_u_est1} and \eqref{eq:tot_u_estJ} for $\mr{\Psi}_z$,  \eqref{def:CC_bounds} 
for $\CCa, \CCj$, for $|\xx| \leq R$, we obtain
\beq\label{eq:reg:est_coe:p}
\bga
  |\cB -2|  \les_{\e, R} |\xx|^{1 + \al} (1 + \nnn{\ww})
  \les_{\e, R} |\xx|^{1+\al}    , \\
    |Q^z - 2 z| \les_{\e, R} |z| \cdot |\xx|^{1 + \al},   \quad |Q^r -  Q_{r0} r|  \les_{\e, R} r \cdot |\xx|^{\al + 1}  .
\ega
\eeq

Since $\cB , \QQ \in C^{1,\al}$, above estimates for $\cB, \QQ = (Q^r, Q^z)$ implies 
\beq\label{eq:reg:est_coe:q}
\bga
|\na_{r, z} \cB |  \les_{\e, R} |\xx|^{\al} , 
\quad  | \na_{r, z} \QQ - \mw{diag}( Q_{r0}, 2  ) |   \les_{\e, R} |\xx|^{\al} .
  \ega
\eeq

\eseq

Next, we establish the $C^{1,\al}$-regularity for the flow map $X(s, \xx)$ defined in \eqref{eq:traj}.

\subsubsection{$C^{1,\al}$ estimates of the flow map}\label{sec:reg:traj_C1a}

To solve $\eta$, we refine the estimate of the flow map $X(s, \xx)$. 
Recall the ODE of $X(s, \xx)$ from \eqref{eq:traj}. Using the lower bound of $Q^r / r$  in \eqref{eq:3D_boot_q:1}, 
and by requiring $\e$ small so that $\bar C_Q \e^{-1} \geq 2$, for any $s \leq 0$, we obtain
\bseq\label{eq:traj_est:Xrz}
\beq\label{eq:traj_est:Xr}
\tf{d}{ds} X^r(s, \xx) = Q^r \cc X(s, \xx) 
\geq \bar C_Q \e^{-1} X^r (s, \xx ) \geq 2 X^r(s, \xx) ,
\quad \Rightarrow  \  X^r(s, \xx) \leq r e^{ 2 s }. 
\eeq

Using the  ODE of $X^z$ in \eqref{eq:traj}, we obtain 
\[
  \f{d}{d s} X^z(s, \xx) e^{ -2 s} = X^z(s, \xx) e^{-2 s} \cdot ( \f{Q^z}{z} - 2  ) \cc X(s, \xx) .
\]

Using $ |Q^z - 2 z| \les_{R, \e} |z| \cdot |\xx|^{1 + \al} $ by \eqref{eq:reg:est_coe:p}, and \eqref{eq:traj_pf1}, for $|\xx| < R$, we derive 
\[
   \int_{-\infty}^0  \B|( \f{Q^z}{z} - 2) \cc X(\tau, \xx ) d \tau \B|
   \les |\xx|^{1  + \al} \int_{-\infty}^0  e^{ C s } d s \les |\xx|^{1 + \al} \les_R 1.
\]
for some absolute constant $C>0$. Integrating the ODE of $X^z(\xx, s) e^{ -2 s}$ from $s$ to $0$, we obtain
\beq\label{eq:traj_est:Xz}
  |X^z(s, \xx ) e^{-2 s}| \asymp_{R, \e} |X^z(0, \xx ) | \asymp_{R, \e} |z|,  \quad \forall  \ |\xx| < R.
\eeq
\eseq
with constant independent of $s$. Combining the above estimates of $X^r, X^z$, we prove 
\beq\label{eq:traj_est:linf}
 |X(s, \xx )| \les_{R, \e} e^{2 s} |\xx|,
 \quad  \forall \ |\xx| < R.
\eeq

Moreover, for $z \neq 0$, we obtain $X^z(\xx, s) \neq 0$ for any $s$. For any $s \leq 0$, we prove \eqref{eq:traj_ratio}:
\[
   \B|\f{X^r(s, \xx)}{X^z(s, \xx)} \B| \les_{R, \e} \f{r}{|z|}.
\]

Next, we perform $C^{1,\al}$ estimate of $X(s, \xx)$. For $i, j \in \{ r, z \}$
using the ODE \eqref{eq:traj}, we obtain
\bseq\label{eq:3D_ODE:X_C1}
\beq\label{eq:3D_ODE:X_C1:a}
  \tf{d}{d s} \pa_j X^i(s, \xx ) = \pa_j ( Q^i( X(t, \xx) ) )
  = \la (\pa_j X^r, \pa_j X^z),  (\na_{r, z} Q^i) \cc X(s, \xx) \ra \, ,
\eeq
where $\la a, b \ra = a_1 b_1 + a_2 b_2$ denotes the inner product in $\R^2$. Below, 
we simplify $(\na \QQ ) \cc X(s, \xx) $ as $\na \QQ$. We introduce the following $2\times 2 $ matrix and matrix norm 
\[
   (\na X)_{ij} = \pa_j X^i, \quad  (\na Q)_{ij} = \pa_j Q^i,
   \quad |M|^2 =  \Tr( M^{\top}, M ),\quad M\in \R^{2 \times 2}.
  \]
We can rewrite the above ODE as 
\beq\label{eq:3D_ODE:X_C1:b}
  \tf{d}{d s} \na X = \na \QQ \cdot \na X.
\eeq
\eseq

For $|\xx| < R$ and $s \leq 0$, using \eqref{eq:reg:est_coe:p} and \eqref{eq:reg:est_coe:q} on $\QQ$, we obtain
\beq\label{eq:dQ:exp}
  (\na \QQ ) \cc X(s, \xx) = O_{R, \e}( |X(s, \xx)|^{\al} ) + \mw{diag}( Q_{r0}, 2 ) .
\eeq
Since $Q_{r0} \asymp \e^{-1}$, by requiring $\e$ small so that $Q_{r0} \geq 2$, we obtain
\[
  \tf{1}{2}\tf{d}{ds} \Tr( (\na X)^T , \na X )(s, \xx )
  = \Tr( (\na X)^T ,  \na \QQ \cdot \na X )(s, \xx)
\geq (2 - C(R, \e) |X|^{\al}) \Tr( (\na X)^T , \na X )(s, \xx) .
\]

Using the estimate \eqref{eq:traj_pf1} and $|\xx| < R$, for any $s\leq 0$, we estimate $|\na X|^2$ as 
\[
  \tf{1}{2} \tf{d}{ds} |\na X|^2(s, \xx)
\geq ( 2 - C(R,\e)  e^{- C s} )  |\na X|^2(s, \xx).
\]

Using $\na X(s, \xx)|_{s =0} = \Id$, 
for any $s \leq 0$ and $|\xx|\leq R$,  we obtain
\beq\label{eq:traj_est:Lip}
  |\na_{r, z} X(s, \xx)|
 \leq \exp( - \int_s^0 ( 2 - C(R,\e)  e^{ C \tau } ) d \tau ) 
 \les_{R , \e} e^{  2 s}.
\eeq

Next, we further estimate the $C^{\al}$ norm of $\na X$. Fix $\xx_1 \neq \xx_2 \in B_R$, 
we have $ X(s, \xx_i) \in B_R$ for any $s \leq 0$.   Denote $\yy_i = X(s,\xx_i)$. We expand
\beq\label{eq:hol_exp1}
\bal
  & (\na \QQ)(\yy_1)  \cdot (\na X)(s, \xx_1)
  - (\na \QQ)(\yy_2) \cdot (\na X)(s, \xx_2) \\
  = & ( (\na \QQ)(\yy_1) - (\na \QQ )(\yy_2) ) \cdot (\na X)(s, \xx_1)
  +  (\na \QQ )(\yy_2) \cdot ((\na X)(s, \xx_1) - (\na X)(s, \xx_2)) \\
   \teq &  I_1 + I_2.
  \eal
\eeq
For $I_1$, using $\QQ \in C^{1,\al}$  \eqref{eq:reg:est_coe:a}, and \eqref{eq:traj_est:Lip}, we bound 
\[
  |I_1| \les_{R, \e} |\yy_1 - \yy_2|^{\al} \| \QQ \|_{C^{1,\al}(  B_R )}
= |X(s, \xx_1) - X(s, \xx_2)|^{\al}  \| \QQ \|_{C^{1,\al}(  B_R )}
\les_{R, \e} e^{2 \al s} |\xx_1 - \xx_2|^{\al} \| \QQ \|_{C^{1,\al}(  B_R )}.
\]

Denote 
\[
  \Ups(s) = (\na X)(s, \xx_1 ) - (\na  X) (s, \xx_2) \in \R^{2 \times 2 }.
  \]

For $I_2$, using \eqref{eq:dQ:exp} and $Q_{r0} \geq 2$, we obtain
\[
  \Tr( \Ups^{\top} , I_2 ) =
    \Tr( \Ups^{\top} , (\na \QQ)(y_2) \Ups ) 
 \geq  ( 2 - C(R, \e) |\yy_2|^{\al} )  \Tr( \Ups^{\top}, \Ups ).
\]

Using \eqref{eq:3D_ODE:X_C1:b}, for any $s \leq 0$, we estimate $|\Ups|^2$ as follows 
\[
  \tf{1}{2} \tf{d}{ds}
  \Tr( \Ups^{\top}, \Ups )
  = \Tr(\Ups^{\top}, I_1  +I_2)
  \geq ( 2 - C(R, \e) |\yy_2|^{\al} )  \Tr( \Ups^{\top}, \Ups )
  - C(R, \e) e^{2 \al s} |\xx_1 - \xx_2|^{\al} |\Ups| \cdot  \| \QQ \|_{C^{1,\al}(  B_R )} .
\]

Since  $|\yy_2| =| X(s, \xx_2 )| \les_{R}  e^{ C s}$ by \eqref{eq:traj_pf1}, 
and $|\Ups|^2=\Tr(\Ups^{\top}, \Ups)$, we derive 
\[
  \tf{d}{ds} |\Ups| \geq (2 - C(R, \e) e^{ C s } ) |\Ups| - C(R, \e) e^{2 \al s} |\xx_1 - \xx_2|^{\al}
  \| \QQ \|_{C^{1,\al}(  B_R )}. 
\]

We introduce the integrating factor $\ell_1(s) = \exp( \int_s^0 (2 - C(R, \e) e^{C \tau}) d \tau )$. 
Since $\Ups(0) =0$, for $s \leq 0$, we obtain
\[
  0 = |\Ups(0)| \geq |\Ups(s)| \ell_1(s) - \int_s^0 \ell_1(\tau) C(R, \e) e^{2 \al \tau} d \tau |\xx_1 - \xx_2|^{\al} \| \QQ \|_{C^{1,\al}(  B_R )}.\| \QQ \|_{C^{1,\al}(  B_R )}.
\]
Since $\ell_1(s) \asymp_{R, \e} e^{- 2 s}$ for any $s \leq 0$ and $\al \in (0, 1)$,
from the above estimate, we prove 
\beq\label{eq:traj_est:C1a}
\bal
 |(\na X)(s, \xx_1) - (\na X)(s, \xx_2 ) | & =   |\Ups(s)| 
 \les_{R, \e} |\xx_1 - \xx_2|^{\al} \int_s^0 \ell_1(\tau) \ell_1(s)^{-1} e^{2 \al \tau} d \tau 
\| \QQ \|_{C^{1,\al}(  B_R )} 
 \\
 & \les_{R, \e} e^{2 \al s } |\xx_1 -\xx_2|^{\al} \| \QQ \|_{C^{1,\al}(  B_R )}. 
 \eal
\eeq

\subsubsection{$C^{k,\g}$ estimate of the solution}\label{sec:traj_solu_C1a}

We have the following estimates for the solution. 

\begin{lem}\label{lem:traj_reg}

Suppose that $\ww$ satisfies \eqref{eq:3D_size:ww} 
and $\QQ = \bar \QQ + \td \QQ(\ww)$ be the transport coefficient defined in \eqref{def:Q}.
Suppose that $a(\xx), f(\xx) \in C^{k, \g}$ with 
$k=0, \g \in (0, 1]$ or $k=1, \g \in (0, \al]$, and
\beq\label{eq:eqn_traj_reg}
   \QQ \cdot \na \zeta = a(\xx) \zeta + f(\xx) , \quad  f(0) = 0, \quad \na a(0) = 0, \quad 
   \lim_{\xx \to 0} \zeta(\xx) =  \zeta(0) .
\eeq
In addition, we consider additional assumptions in three cases 
\beq\label{eq:traj_reg:case}
\bal
& (1) \, a(0) = 0, \  k=0 ; \qquad 
 (2)  \,  a(0) = 0, \  k=  1 ; \\ % \quad  \mbox{or } \\
& (3) \, a(0) \leq 2 ,   \quad   \na f(0) = 0, \quad k = 1,
\quad \lim_{|\xx| \to 0} \f{\ze}{|\xx|} = 0 .
  \eal
\eeq

For any $R >0$, we have 
\beq\label{eq:traj_reg:est1}
     \|   \zeta  \|_{ C^{k, \g}(\bar B_R)}    \les_{R, \e}    \exp( C_R \| a \|_{ C^{k, \g}(\bar B_R)}  ) \cdot ( |\ze(0)| +   \| f \|_{C^{k, \g}( \bar B_R) } ).
\eeq

For the case (3) $a(0)\leq 2, \na f(0) = 0, k = 1$, we have $\ze(0) = 0, \na \ze(0)=0$ and for any $\xx \in \bar B_R$:
\beq\label{eq:traj_reg:est2}
  | \zeta(\xx) | \les_{R, \e} |\xx|^{1 + \g}   \exp( C_R \| a \|_{ C^{1, \g}(\bar B_R)}  ) \cdot     \| f \|_{C^{1, \g}( \bar B_R) } .
\eeq

\end{lem}

\begin{proof}

Firstly, from \eqref{eq:psi_reg:d}, we obtain $\QQ \in C^{1,\al}$ with $\| \QQ \|_{C^{1,\al}(\bar B_R)} \les_{\e, R} 1.$ Moreover, the associated flow map $X(s, \xx)$ satisfies the estimates in Section \ref{sec:reg:traj_C1a}.

To estimate three cases together, we modify the equation \eqref{eq:eqn_traj_reg} and $\ze$. For parameter $\hh = (h_r, h_z) $ to be determined, we consider 
$\td \ze = \ze - \ze(0) - \hh \cdot \xx = \ze - \ze(0) - h_r r - h_z z$. Using the equation \eqref{eq:eqn_traj_reg}, we obtain
\beq\label{eq:traj_reg:eqn1}
  \QQ \cdot \na (\td \ze + \ze(0) + \mathbf{h} \cdot \xx)
  = a(\xx) (\td \ze + \ze(0) +  \mathbf{h}  \cdot \xx) + f(\xx) ,
  \quad \td \ze = \ze - \ze(0) - \hh \cdot \xx.
\eeq

Using $\QQ = (q^r r, q^z z)$, we rewrite the above equation as 
\beq\label{eq:traj_reg:eqn2}
  \QQ \cdot \na \td \ze = a(\xx) \td \ze + \td f, 
  \quad \td f(\xx) \teq f(\xx) + a(\xx) \ze(0) + ( a(\xx) - q^r ) h_r r +
  ( a(\xx) - q^z ) h_z z.
\eeq

For case (1) $a(0)=0, k=0$ or the case (3) $k=1$ $a(0) \neq 0 , \na f(0) = 0$ \eqref{eq:traj_reg:case}, we choose $\hh = 0$. 
For the case (2) where $a(0) = 0$ and $k = 1$, we obtain $f, a \in C^{1, \g}$ with $\g \in (0, \al ]$. 
Since $q^r(0) = Q_{r0} > 0, q^z(0) = 2$ by \eqref{def:Qr0}, we choose  
\beq\label{eq:traj_reg:h}
  h_{i} = - (q^i(0))^{-1}  \cdot  \pa_i ( f(\xx)  + a(\xx) \ze(0) )|_{\xx=0} , \quad i \in \{ r, z \}.
\eeq
By definition, we obtain $\td f(0) = 0$ and  $ \na \td f(0 ) = 0$ when $k = 1$. Using \eqref{eq:3D_boot_q:1}, we bound 
\beq\label{eq:traj_reg:h_bd}
  |h_i| \les |\na f(0) |+ |\na a(0) | \cdot |\ze(0)|.
\eeq

Let $X(s, \xx)$ be the flow map  \eqref{eq:traj} associated with $\QQ$.  Since $\td \ze(0 ) = 0$, solving $\td \zeta$ \eqref{eq:traj_reg:eqn2}  along the characteristics, we obtain
\beq\label{eq:traj_reg:ze1}
     \tf{d}{d s } \td \ze \cc X(s, \xx  ) = ( a(\xx) \cdot \td \ze ) \cc X(s, \xx) + \td f  \cc X(s, \xx) .
\eeq

Below, we estimate the $C^{k, \g}$ bound of $a \cc X(\cdot , s), \, \td f \cc X(\cdot ,s)$.

\vs{0.1in}
\paragraph{\bf Composition bound}

For any $\xx \in B_R$, from \eqref{eq:traj_pf1}, we have $X(s, \xx) \in B_R$ for $s \leq 0$.  We consider $F \in C^{k,\g}(B_R)$ with $F(0) =0$ and an extra assumption $\na F(0) = 0$ when $k = 1$. For any $\xx, \xx_i \in B_R$, $s \leq 0$, using \eqref{eq:traj_est:Lip}, we have 
\beq\label{eq:compo:fX0}
\bal
  \| F \cc X(s, \xx) \|_{ \dot C^{\g}(\bar B_R) } & \les  \| F \|_{ \dot C^{\g}(\bar B_R)} 
   \| \na X(\cdot , s ) \|_{ L^{\infty}(B_R) }^{\g} 
  \les e^{ 2 \g s } \| F \|_{\dot C^{\g}(\bar B_R)} ,  \\
  |F \cc X(s, \xx)| & \les |X(s, \xx)|^{\g} \| F \|_{ \dot C^{ \g}(\bar B_R)} 
  \les e^{ 2 \g s } |\xx|^{\g}\| F \|_{ C^{\dot \g}(\bar B_R)}
  \les_R e^{ 2 \g s } \| F \|_{ C^{\dot \g}(\bar B_R)}  \, .  \\
  \eal
\eeq
When $k =1$, $\na F(0) = 0$ and $\g \leq \al$, we estimate 
\beq\label{eq:dF_compos}
\bal
\pa_i ( F \cc X(s, \xx) )
& = \la  \pa_i X(s, \xx), (\na F) \cc X( s, \xx ) \ra, \\
|(\na F ) \cc X(s, \xx_1 ) - (\na F ) \cc X(s, \xx_2)|
&\les %\| (\na F ) \cc X( \cdot, s  ) \|_{\dot C^{\g}(B_R)} 
\| \na F  \|_{\dot C^{\g}(B_R)} 
|X(s, \xx_1) - X(s, \xx_2)|^{\g} \\
&\les \| F \|_{C^{1,\g}(B_R) }  \| (\na X)(\cdot ,s ) \|_{L^{\infty}(B_R)}^{\g}
|\xx_1 - \xx_2|^{\g}.
\eal
\eeq

Applying the $L^{\infty}$ and the 
$\dot C^{\g}$ estimate to the products and using \eqref{eq:traj_est:Lip},
for $\g \leq \al$,  we obtain
\[
\bal
  | \na ( F \cc  X(s, \xx ) )| & \les \| \na X(\cdot, s) \|_{L^{\infty}(B_R)} 
  \| (\na F ) \cc X(\cdot , s ) \|_{L^{\infty}(B_R)} 
  \les_{R, \e} e^{ 2 s } \| (\na F ) \cc X(\cdot , s ) \|_{L^{\infty}(B_R)} , \\
  \| \na ( F \cc X(\cdot , s) )  \|_{\dot C^{\g}(B_R)} 
 &  \les \| \na X(\cdot , s) \|_{\dot C^{\g}}    \| (\na F)\cc X(\cdot ,s) \|_{L^{\infty}(B_R)} \\
 & \quad + \| \na X(\cdot , s) \|_{L^{\infty}(B_R)} 
  \|  F \|_{ C^{1, \g}(B_R)}  \| \na X(\cdot , s) \|^{\g}_{L^{\infty}(B_R)}  \\
  & \les_{R, \e }  e^{2 \g s} \| (\na F)\cc X(\cdot ,s) \|_{L^{\infty}(B_R)}
  + e^{2 (1 + \g) s }   \|  F \|_{ C^{1, \g}(B_R)} .
\eal
\]

Suppose that $\na F(0) = 0, F(0)=0$. Using $ F \in C^{1, \g}$, \eqref{eq:traj_est:Lip}, 
and \eqref{eq:traj_est:linf}, for any $|\xx| < R$ and $ s\leq 0$, we obtain
\beq\label{eq:compo:fX2}
\bal
  |(\na F ) \cc X(s, \xx) | & \les \| F \|_{ C^{1,\g}(B_R) } |X(s, \xx)|^{\g}
  \les_{R, \e} e^{ 2 \g s } \| F \|_{ C^{1,\g}(B_R) }  , \\
      | F \cc X(s, \xx)| & \les \| F \|_{ C^{1, \g}(B_R) } |X(s, \xx)|^{1 + \g}
        \les_{R, \e} e^{ 2 (1 + \g) s } |\xx|^{1 + \g} \cdot \| F \|_{ C^{1,\g}(B_R) } .
\eal 
\eeq

Combining the above estimates, % \eqref{eq:compo:fX1} and \eqref{eq:compo:fX2}, 
for any $f  \in C^{k, \g}$ with $f(0) = 0$ and an extra assumption $\na f (0) = 0$ when $k = 1, \g \in (0, \al] $, and for any $s \leq 0$, we establish 
\beq\label{eq:compo:fX3}
  \| F \cc X(\cdot ,s) \|_{C^{k,\g}(B_R)} \les_{R, \e} e^{ 2 (\g+ k) s } \| f \|_{ C^{k,\g}(B_R) }.
\eeq

From the assumption of $\td f(\xx), a(\xx)$ from \eqref{eq:eqn_traj_reg} 
and \eqref{eq:traj_reg:h},  $F(\xx) = a(\xx) - a(0) \in C^{1,\al}$ and $f(\xx) = \td f(\xx)$ satisfies $f(0)$, and $\na f(0)=0, f \in C^{1, \g}$ when $k=1$.

Integrating \eqref{eq:traj_reg:ze1}, we obtain
\beq\label{eq:traj_reg:ze}
  \td \ze(\xx) = 
  \lim_{s \to -\infty} e^{ \cS(s, \xx) } \td \ze( X(s, \xx) )+
   \int_{s}^0 e^{ \cS(\tau, \xx) } \td f( X(\tau, \xx) ) d \tau,
  \quad \cS(s, \xx) = \int_s^0 a( X(\tau, \xx ) ) d \tau.
\eeq

 For any $s \leq 0$, using \eqref{eq:compo:fX3}, we obtain
\[
  \| \cS(s, \cdot ) + a(0) s  \|_{ C^{k, \g}(\bar B_R)}
  \les \int_s^0 \|  a(X(\tau, \xx)) - a(0) ) \|_{ C^{k, \g}(\bar B_R)}
   \les_R \int_s^0 e^{2 (k+ \g) s}  d \tau \cdot \| a - a(0) \|_{ C^{k, \g}(\bar B_R)} 
  \les_R  \| a  \|_{ C^{k, \g}(\bar B_R)}  .
\]
Since $\na e^{\cS} = \na \cS \cdot e^{\cS(s, \xx)}$, for $s \leq 0$, we obtain
\[
\bal
    \| e^{S(s, \xx) } \|_{ C^{\g}(\bar B_R)} & \les   \| a \|_{ C^{\g}(\bar B_R)}  e^{ - a(0) s + C_R \| a \|_{\dot C^{\g}(\bar B_R)}  } 
    \les      e^{ -a(0) s + C_R \| a \|_{\dot C^{\g}(\bar B_R)}  } , \\
    \|  e^{S(s, \xx) } \|_{ C^{1, \g}(\bar B_R)} 
  & \les ( \| \na \cS(s, \cdot) \|_{C^{\g}(\bar B_R)}   + 1 )     \| e^{S(s, \xx) } \|_{ C^{\g}(\bar B_R)}
  \les   e^{ - a(0) s + C_R \| a \|_{ C^{1, \g}(\bar B_R)}  } ,
\eal 
\]
where we have used $  \| a \|_{ C^{k, \g}(\bar B_R)} \les  \exp(  \| a \|_{ C^{k, \g}(\bar B_R)} )$.  

In case (1), (2) in \eqref{eq:traj_reg:case}, we have $a(0) = 0$, $X(s, \xx) \to 0$  as $s \to -\infty$, and $\td \ze(\xx) \to 0$ as $|\xx| \to 0$. 
In case (3) in  \eqref{eq:traj_reg:case}, we have 
 $a(0) \leq 2$, $\ze(0) = h_i= 0, \td \ze(\xx) = \ze(\xx) 
= o_{|\xx|\to 0}(|\xx|)$ by \eqref{eq:traj_reg:case}  and \eqref{eq:traj_reg:h}. In all three cases, using the above estimates and \eqref{eq:traj_est:linf}, we have 
\beq\label{eq:traj_reg:ze_init}
  \limsup_{s \to -\infty} | e^{ \cS(s, \xx) } \td \ze( X(s, \xx) ) |
  \les_{R, \e}  e^{  C_R \| a \|_{ C^{1, \g}(\bar B_R)}  }  
    \limsup_{s \to -\infty} \f{ |\td \ze( X(s, \xx) ) | }{ |X(s, \xx)| }
     e^{- a(0) s} \cdot e^{2 s} |\xx| = 0.
\eeq
Thus, we get $\lim_{s \to -\infty} e^{ \cS(s, \xx) } \td \ze( X(s, \xx) ) = 0 $.

Therefore, applying the above estimates,
 \eqref{eq:compo:fX2},  \eqref{eq:compo:fX3} to $\td f, a(\xx) -a(0)$, and product rule, for any 
$k= 0, \g \in (0, 1]$ or $k=1, \g \in (0, \al]$,  we bound $\td \ze$ in \eqref{eq:traj_reg:ze} as
\beq\label{eq:3D_traj_est2}
\bal
   \| \td  \zeta  \|_{ C^{k, \g}(\bar B_R)}   &  \les   \int_{-\infty}^0
  \| e^{ \cS( \cdot, \tau) } \|_{  C^{k, \g}(\bar B_R)}  \| \td f( X(\cdot, \tau) ) \|_{C^{k, \g}(\bar B_R)}
 d \tau  \\
 & \les  \int_{-\infty}^0 e^{ - a(0) \tau + 2 (k + \g) \tau } d \tau \cdot
 e^{ C_R \| a \|_{ C^{k, \g}(\bar B_R)}  }  \| f \|_{C^{k, \g}( \bar B_R) } , \\
  | \td  \zeta(\xx) |  & 
    \les   \int_{-\infty}^0  e^{ - a(0) \tau + 2 (1 + \g) \tau } d \tau \cdot 
    |\xx|^{ 1 + \g }  e^{  C_R \| a \|_{\dot C^{\g}(\bar B_R)}  } 
     \| \td f \|_{C^{1, \g}(\bar B_R)} .
  \eal
\eeq

Recall three cases from \eqref{eq:traj_reg:case}. We have $a(0) = 0$ for case (1), (2), or $a(0) \leq 2$ and $k=1$ for case (3). In all cases, we obtain $ -a(0) \tau + 2 (k + \g) \tau \leq 2 \g \tau$ 
for $\tau \leq 0$. Thus, the above integrals are bounded. Combining \eqref{eq:3D_traj_est2}  for $\td \ze$, \eqref{eq:traj_reg:h_bd}, 
and using \eqref{eq:traj_reg:eqn1}, and triangle inequality, we prove
\eqref{eq:traj_reg:est1}. 

For case (3) in \eqref{eq:traj_reg:case}, since we choose $h_i =0$ \eqref{eq:traj_reg:h} and $\zeta(0) = 0$ by \eqref{eq:traj_reg:case}, we have $\td \zeta = \zeta $ and obtain $\na \zeta(0) = 0$ from \eqref{eq:3D_traj_est2}. We prove \eqref{eq:traj_reg:est2}. 
\end{proof}

\subsubsection{$C^{1,\al}$ estimates of $\eta$}\label{sec:reg:eta_C1a}

Recall the equation of $\eta$ from \eqref{eq:fix_pt_2D} 
\[
    \QQ  \cdot \na \eta =    \cB  \cdot \eta + H(\ww) .
\]
and recall $H(\ww), \cB$ from \eqref{eq:reg:eqn_w}. Using estimates \eqref{eq:reg:est_coe:a}, \eqref{eq:reg:est_coe:cR}  on $\wwb (\bar \cR + \Lpsi)$, we obtain
\bseq\label{eq:reg:C1a_H}
\begin{align}
 |H(\ww)|  & \leq |\wwb (\bar \cR + \Lpsi) | 
 \les_{R,\e } |z| \cdot |\xx|^{1 + \al} ,
     \quad \| H(\ww ) \|_{C^{1,\al}(B_R)}
    \les_{R , \e } 1 . 
   \label{eq:reg:C1a_H:c} 
\end{align}
\eseq
For $\cB$, using \eqref{eq:reg:est_coe:p} and \eqref{eq:reg:est_coe:a}, we obtain 
\[
  |\cB  - 2| \les_{\e, R} |\xx|^{1+\al}, \quad  \| \cB \|_{C^{1,\al}(\bar B_R)} \les_{\e, R} 1 .
\]
The above regularity and vanishing orders near $\xx=0$ imply
 $\na H(\ww)(0) = 0,  \na \cB(0)  = 0$.

Note that $\QQ \in C^{1, \al}$ and $\ww$ satisfies the assumptions in Lemma \ref{lem:traj_reg}. 
Moreover, the above estimates and \eqref{eq:3D_ODE:init} for $\eta$ imply that 
$a = \cB , f = H(\ww), \eta $ satisfy the assumptions in \eqref{eq:eqn_traj_reg} and in \eqref{eq:traj_reg:case} for the case (3).  Since we choose $\hh = 0$ for case (3) \eqref{eq:traj_reg:h},
applying Lemma \ref{lem:traj_reg} and the solution formula \eqref{eq:traj_reg:ze}, 
\eqref{eq:traj_reg:ze_init} for $\td \zeta = \eta$, we obtain 
\beq\label{eq:traj_reg:eta_form}
  \eta(\xx) =  \int_{s}^0 e^{ \cS(\tau, \xx) } H(\ww)( X(\tau, \xx) ) d \tau,
  \quad \cS(s, \xx) = \int_s^0 \cB( X(\tau, \xx ) ) d \tau.
\eeq

Applying \eqref{eq:traj_reg:est1},\eqref{eq:traj_reg:est2} in Lemma \ref{lem:traj_reg} with $(k, \g) = (1,\al)$, we prove estimates \eqref{eq:traj_reg:est1}, \eqref{eq:traj_reg:est2} for $\eta$
\[
 \| \eta \|_{C^{1,\al}(B_R)} \les_{R, \e }   1,
 \quad \na \eta(0) = 0,
 \quad |\na^i \eta(\xx) | \les_{R, \e} |\xx|^{ \al + 1 -i }.
\]

Since $\eta$ is odd in $z$, from $|\na \eta(\xx) | \les_{R, \e } |\xx|^{\al}$, we prove $ |\eta(\xx) | \les_{R, \e } |z| \cdot |\xx|^{\al}$. We prove \eqref{eq:local_reg}.

\begin{comment}

\vs{0.05in}

\paragraph{\bf Proof of \eqref{eq:local_reg:van}}

Recall  $H(\ww)(\xx)$ and $\cB(\xx)$ from \eqref{eq:reg:eqn_w} and \eqref{eq:lin_2Db}. 
From \eqref{eq:psi_r0_van} in Lemma \ref{lem:psi_reg}, we obtain
\[
 \pa_r \Psi(0, z) = r \pa_r^2 \Psi(0, z) = 0,
\, \Rightarrow  \pa_r Q^z(0, z) = 0, \quad
\pa_r f(0, z)  = 0,  \quad \mbox{for} \ f = H(\ww), \cB.
\]

Using \eqref{eq:traj_est:Xr}, we obtain $X^r(s, 0, z) \equiv 0$. 
Since $\na X(s, \xx)|_{s =0} = \Id$, we obtain $\pa_r X^z(0, 0, z) |_{s = 0} = 0$.
Using the ODE \eqref{eq:3D_ODE:X_C1} for $\pa_r X^z$ and $\pa_r Q^z(0, z) = 0$, we yield  
\[
  \tf{d}{ds} \pa_r X^z(s, 0, z) = \pa_r X^r  \cdot \pa_r Q^z + \pa_r X^z \cdot \pa_z Q^z \B|_{r=0}
  = \pa_r X^z \cdot \pa_z Q^z |_{r=0} .
\]
Thus, we obtain $\pa_r X^z(s, 0, z) \equiv 0$. Using \eqref{eq:dF_compos}, 
for $ f = H(\ww), \cB$, we obtain
\[
 \pa_r f(X(s, \xx)) |_{r = 0} = \pa_r X(s, 0, z) \cdot \na f (X(s, 0, z))
 =  \pa_r X^r(s,0,z) \pa_r f( 0, X^z(s, 0, z) ) = 0.
\]

Using the formula \eqref{eq:traj_reg:eta_form}, we prove  $\pa_r \eta(0, z ) = 0$ and 
\eqref{eq:local_reg:van}.

\end{comment}

\vs{0.05in}

\paragraph{\bf Proof of \eqref{eq:local_reg:nnn}}
Finally, we prove \eqref{eq:local_reg:nnn}. By definitions of \eqref{def:3d_wg} 
and \eqref{def:rhoc}, for $\xx \in B_R$, 
we have $\rag , \rhoc \les_R 1, \rhoo   \les |\xx|^{1 - \bbb } 
\les_R |\xx|^{-\al }$. 
Recall $\wwb = \waa$ from \eqref{def:3d_Om}. Using \eqref{eq:wa_upper_lower} and \eqref{eq:1D_prof_est2:c} :
\[
 |\wwb| = |\waa| \asymp_R |z| ,\quad |\pa_z \wwb^{-1}| = |\f{\pa_z \wwb}{\wwb^2}|
 \les_R |z|^{-2},\quad  \forall |z| \leq R \, ,
 \]
estimate \eqref{eq:local_reg}, and $G = \f{\eta}{\wwb}$,  we obtain
\[
\bal
  |\rhoc \rhoo \rag G | 
  +   |\f{z}{\ang z} \rhoc \rhoo \rag \ang \xx \cdot \pa_r G | 
  +|  \rhoc \rhoo \rag |z|\cdot \pa_z G | 
& \les_R |\xx|^{-\al} ( \f{|\eta|}{|\wwb|} + |z| \cdot \B| \f{ \pa_r \eta}{\wwb}  \B|
+ |z| \cdot \B| \pa_z \f{\eta}{\wwb}\B| ) \\
& \les_R |\xx|^{-\al } ( |z|^{-1} |\eta| + |\pa_r \eta|
+ |\pa_z \eta| )  \les_{R, \e } 1,
\eal
\]
for any $\xx \in B_R$.  We prove  \eqref{eq:local_reg:nnn}.

\vs{0.05in}
\paragraph{\bf Proof of \eqref{eq:traj_smooth}}

Fix $ z \neq 0$ and $i \in \{ r, z \}$. We introduce the finite difference $\D_h f(\xx) \teq \f{1}{h} ( f(\xx  + h \cdot \ee) - f(\xx) )$  
with $\ee = (1, 0)$ for $i=r$ and $\ee = (0, 1)$ for $i = z$. 
Using  \eqref{eq:fix_pt_2D}, we obtain
\[
 \QQ(\xx) \cdot \na \D_h G(\xx) = \D_h F(\xx) - \D_h \QQ( \xx) \cdot \na G(\xx + h)  , \quad F \teq (\bar \cR + \Lpsi) (G+1), 
\]
for $ 0 < h < |z|/2 $. 
Since $h >0$, we obtain $ \xx + h \cdot \ee, \xx \in \R_+ \times \R$. Since $z \neq 0$, from \eqref{eq:traj_est:Xrz}, 
we obtain $X^z(s, \xx) \neq 0$. 
Integrating the above 
equation along trajectory, we obtain
\[
  (\D_h G) \cc  X( s+ \tau, \xx) = (\D_h G) \cc  X(s, \xx) +  \int_{s}^{s + \tau}
  \big[ \D_h F( \yy ) - \D_h \QQ( \yy ) \cdot \na G( \yy + h \cdot \ee ) \big] \B|_{ \yy =   X( \ze, \xx) } d \ze.
\]

Since $X( s , \xx)$ is uniformly bounded away from $z = 0$, we obtain 
$|\wwb(X(\ze, \xx))| \gtr 1$ for $\ze \in [s , s + \tau]$. 
Since $\eta, F, \QQ $ are $C^{1,\al}$ locally and away from $z = 0$
by  \eqref{eq:reg:est_coe},\eqref{eq:local_reg}, 
taking $h\to 0$, we obtain 
\[
  (\pa_i G) \cc X( s+ \tau, \xx)  = (\pa_i G) \cc X(s, \xx) +  \int_{s}^{s + \tau}
  \big[ \pa_i F( \yy ) - \pa_i \QQ( \yy ) \cdot \na G( \yy ) \big] \B|_{ \yy = X( \ze, \xx) } d \ze.
\]
Taking $ \f{d}{d \tau} $ in the above identities in $\tau$ and evaluating at $\tau =0$, we prove \eqref{eq:traj_smooth}.

 We complete the proof of Proposition \ref{prop:reg_solu}.

\section{Energy estimates for the fixed-point map}\label{sec:3D_EE}

In this section, we perform energy estimates on the fixed point equation \eqref{eq:fix_pt_2DG} 
for $\FFF$ \eqref{eq:fix_pt_2D}:
\beq\label{eq:fix_recall:G}
     \QQ \cdot \na  G   = ( \bar \cR + \Lpsi(\ww)) (G + 1) , 
 \quad G \teq \tf{\eta}{\wwb} ,
 \quad \QQ = \bar \QQ + \td \QQ(\ww) .
\eeq

Our goal is to prove Theorem \ref{thm:into}, which shows that \(\FFF\) maps a ball in a weighted $W^{1,\infty}$ space into itself.

\subsection{Trajectory estimates and integration by parts along trajectories} 

To perform estimates along the trajectory, we use the following estimates.

\begin{lem}[\bf Weighted trajectory estimates]\label{lem:traj_est}

Suppose that $\ww$ satisfies \eqref{eq:3D_size:ww} and $\QQ$ is the 
transport coefficient associated with $\ww$. Suppose that $m, \vp \geq 0$ are some weights satisfying  $\QQ \cdot \na ( m \vp^{-1} )$ does not change sign, $m(\xx) , \vp(\xx)> 0$ for any $\xx \neq 0$, and $ \lim_{(r,z) \to 0} \vp m^{-1}(0) =0 $. If  
\beq\label{eq:traj_pt}
   m(\xx) | H(\xx) | \leq M |\f{\QQ \cdot \na ( m \vp^{-1}) }{m \vp^{-1} }| \, , 
\eeq
holds for any $\xx$, then we have 
\beq\label{eq:traj_est}
m(\xx)  \vp(\xx)^{-1} \int_{-\infty}^0   | (\vp H) \circ X( s, \xx )| d s 
   \leq M .
\eeq
\end{lem}

Let's motivate Lemma \ref{lem:traj_est}.
Consider the PDE: $ \QQ \cdot \na f = c(\xx) f + H(\xx)$ with some coefficient $c(\xx)$ and  a forcing $H$. We introduce an integrating factor : $\QQ \cdot \na \vp = - c(\xx) \vp$. 
If $(f \vp)(\xx)\to 0$ as $\xx \to 0$,  we obtain $ f(\xx) = \vp(\xx)^{-1} \int_{-\infty}^0 (\vp H) \cc X( s, \xx ) d s$.
To estimate $m f$, we rewrite the PDE as 
\[
 \QQ \cdot \na ( f m ) =  d(m, \vp)
 \cdot f m +  m H(\xx) , \quad d(m,\vp) \teq \tf{ \QQ \cdot ( m \vp^{-1}) }{m \vp^{-1} } ,
\]
  We use Lemma \ref{lem:traj_est} with $Q \cdot \na (m \vp^{-1}) < 0$ and treat 
$ d(m,\vp) < 0$ as a damping coefficient. 
By comparing $mH$ and the ``amount of damping" $ d(m,\vp)$, we 
estimate $m f$ without solving $f$ explicitly. 
\footnote{
Consider a simple ODE: $ \tf{d }{dt} a(t) = - a(t) + c$ with $a(0) = 0, c>0$. 
We obtain $ |a(t)| < c$ for any $t >0$ using a simple bootstrap argument. The above comparison generalizes this simple estimate without solving the PDE.
}
Moreover, 
we can design the weight $m$ to obtain tighter estimates, i.e. smaller $M$ in \eqref{eq:traj_est}.

\begin{proof}[Proof of Lemma \ref{lem:traj_est}]

Using the assumptions \eqref{eq:traj_pt} and those on the sign in the statement, we bound 
\[
  I:= \int_{-\infty}^0  | (\vp H) \circ X( s, \xx )| d s
  \leq \B| \int_{-\infty}^0 \big(  M \vp m^{-1} \cdot \f{\QQ \cdot \na ( m \vp^{-1}  ) }{m \vp^{-1} } 
  \big) \circ X( s, \xx  ) d s \B| \, .
\]

Using the characteristic \eqref{eq:traj}, $X(\xx, 0) = \xx$, and $X(\xx,-\infty) = 0$,  we obtain
\[
  I 
\leq    M \B|   \int_{-\infty}^0   \f{\QQ \cdot \na ( m \vp^{-1} ) }{ ( m \vp^{-1} )^2 } \circ X( s, \xx  ) d s \B|
  \leq M \B|  \int_{-\infty}^0 \f{d}{d s }  \big( \f{1}{m \vp^{-1}} \circ X( s, \xx )  \big) d s \B| 
  = M  \B| \f{ 1 }{ (m \vp^{-1})(\xx)} \, -\,  \f{ 1 }{ (m \vp^{-1}) (0)} \B|  .
\]
Since $m\vp^{-1}(0) = 0$, rewriting the above estimate, we prove the lemma.
\end{proof}

\subsubsection{Integration by parts along trajectories}\label{sec:IBP}

To exploit the anisotropic structure and overcome the unbounded estimates in the region with $z \ll r$,
we develop a fundamental estimate.

\begin{lem}\label{lem:IBP}
Let  $X(s, \xx)$ be the flow map associated with $\QQ$ \eqref{eq:traj}.  Suppose that $\ww$ satisfies \eqref{eq:3D_size:ww} and the assumption in Proposition \ref{prop:reg_solu}, and
\beq\label{eq:IBP_ass}
 g \in C^1(D), \quad  f, \  \tf{1}{r} g , \   r^{-1}  z \pa_z g \in C^0(D), 
 \quad 
  \tf{d}{ds} f(X(s, \xx)) \in C^0( [-T, T]\times D) , \quad g(0, z) = 0, 
\eeq
for any compact domain $D \subset (\R_+ \times \R ) \bsh \{ z=0 \}$ and $T>0$. 
\footnote{
\label{foot:no_z0}
We exclude $\{z=0\}$ since our estimates are formulated for $\eta/\wwb$ 
using its equation \eqref{eq:fix_pt_2DG}, which loses $C^{1,\alpha}$ regularity there. This is only to simplify the
presentation. The main result, Theorem~\ref{thm:into}, can also be proved by 
using equation \eqref{eq:fix_pt_2D:a} for the more regular variable $\eta$, at the cost of lengthier expressions.
}
We consider $\xx = (r, z)$ with $z\neq 0$ and denote $X_s = X(s, \xx)$. 
For any $\tau > 0$, we have
\beq\label{eq:IBP}
  \int_{-\tau}^{0} (\pa_r g \cdot f) \co X_s d s
   = 
    - \int_{-\tau}^{ 0 }  \B[ (   \f{Q^z}{Q^r} \pa_z g \cdot f  
   )  \co  X_s
   +   g \cc X_s \cdot  \f{d}{d s}  \big( \f{f}{Q^r}  \cc X_s \big) \B] d s 
+   \f{fg}{Q^r} \cc X_s \B|_{s = -\tau}^{s =0}. %( \xx ).
\eeq

\end{lem}

The crucial feature of \eqref{eq:IBP} is that it exploits the cancellation in $\pa_r g$ along the characteristics, and distributes the $r$-derivative from $g$ to $\pa_z g$ and to $ \f{d}{ds} f \cc X_s$.  See \hyr[sec:idea_singular]{\its Integration by parts} for motivation.

\begin{proof}

Below, we simplify $X_s = X(s, \xx)$ and denote $\f{d}{ds} h \cc X_s 
= \f{d}{ds} h(X(s, \xx))$.  For any locally $C^1$ function $h$, 
since $ \f{d}{ds} X_s = \QQ \cc X_s$ \eqref{eq:traj}, 
we have
\beq\label{eq:IBP_pf_iden}
       \tf{d}{d s} h \cc X_s = (\QQ \cdot \na h ) \co  X_s, 
   \quad \QQ = (Q^r, Q^z).
\eeq

Since $g$ is locally $C^1$, using Leibniz rule, and \eqref{eq:IBP_pf_iden}, we obtain
\bseq\label{eq:IBP_pf}
\beq\label{eq:IBP_pf:a}
 \f{d}{ds} ( \f{f g}{Q^r} ) \cc X_s
    = (\pa_r g \cdot f) \cc X_s + \big( \f{Q^z}{Q^r} \pa_z g \cdot f \big) \cc X_s +  g \cc X_s \cdot  \f{d}{d s} 
    ( \f{f}{Q^r} ) \cc X_s
    \teq I_1 + I_2 + I_3.
\eeq

We verify that $I_i(s, \xx) \in C^0( D_{T}), D_{T} = [-T, 0] \times D$  for any compact 
domain $D \subset (\R_+ \times \R ) \bsh \{ z=0 \}$ and $T>0$. 
From estimate \eqref{eq:traj_ratio}, we obtain $|X^z(s, \xx)| \geq C_{D, T} > 0$ 
for any $(s, \xx) \in D_T$. Since $\ww$ satisfies \eqref{eq:3D_size:ww}, using \eqref{eq:3D_boot_q:1}, we obtain   $Q^r = r q^r \gtr r, |Q^z| \les \e^{-1} z$.  Using Lemma \ref{lem:psi_reg}, $Q^z(r,0) = 0, Q^r(0, z) =0$, 
and assumption \eqref{eq:IBP_ass}, we obtain $I_1, I_2 \in C^0$. 
Since $Q^r = r ( c_l - \pa_z \Psi ) $ by \eqref{def:Q}  and $\pa_z \Psi, \QQ \in C^{1, \al}$ by Lemma \ref{lem:psi_reg}, 
using assumption \eqref{eq:IBP_ass}, we obtain $\f{1}{r}Q^r \in C^0$ and
\[
  g \cc X_s \cdot  \f{d}{ds} \f{ f}{Q^r} \cc X_s=   \f{g}{Q^r} \cc X_s  \cdot \f{d}{d s} f \cc X_s-  \B( \f{ f g}{ (Q^r)^2}
 ( Q^r \pa_r Q^r + r Q^z \pa_z ( c_l - \pa_z \Psi) ) \B) \cc X_s \in C^0(D_T).
\]
Thus, $I_3 \in C^0(D_T)$.  Integrating \eqref{eq:IBP_pf:a}, we prove 
\beq
\int_{-\tau}^0 \B[ \B( \pa_r g \cdot f + \f{Q^z}{Q^r} \pa_z g \cdot f \B) \cc X_s 
+  g \cc X_s  \cdot  \f{d}{ d s} ( \f{f}{Q^r} ) \cc X_s \B] d s
= \f{ fg }{Q^r} \cc X(s, \xx) \B|_{s=-\tau}^{s=0} .
\eeq
\eseq
Since $I_i \in C^0(D_T)$, the integral is well-defined.  Rearranging the terms,  we complete the proof.
\end{proof}

\subsection{$L^{\infty}$ Estimate along $r = 0$}\label{sec:3D_boundary}

In this section, we perform the weighted $L^{\infty}$ estimate of $G$ along $r = 0$ 
\eqref{eq:fix_recall:G} using the contraction estimate in Theorem \ref{thm:contra_lin_exact}.

Recall the fixed point equation from \eqref{eq:fix_recall:G}, $\Lpsi$ from \eqref{eq:lin_2Db}, 
and $\QQ = \bar \QQ + \td \QQ(\ww)$ from \eqref{def:Q}. Restricting \eqref{eq:fix_recall:G} to $r = 0$ and using 
$\bar \cR(0, z) = 0$ \eqref{eq:iden_r0} and $Q^r(0, z) = 0$ \eqref{def:Q}, \eqref{def:3d_vel}
we obtain 
\beq\label{eq:3D_1D_eqn}
   Q^z \pa_z G =  \Lpsi(\ww) (G+1) , 
   \quad G = \f{\eta}{\wwb} ,  \quad Q^z(0, z) = c_l z + 2 \psi(0,z).
\eeq

Since $G(0) = 0$ by \eqref{eq:local_reg}, dividing $Q^z(0, z)$ and integrating 
\eqref{eq:3D_1D_eqn} from $0$ to $z$, we obtain
\beq\label{eq:lin_bd_3D0}
  G  = \cFc(\ww, \eta) \teq \int_0^z \f{ \Lpsi }{Q^z} (G+1) .
\eeq

Recall $G = \f{\eta}{\waa}, \waa = \wwb$ \eqref{def:3d_Om}. To compare $\cFc$ with the map in 1D 
$\cLa$ \eqref{eq:lin_1D_alpha}, we decompose 
\beq\label{eq:lin_bd_3D}
\bal
  \cFc(\ww, \eta) & =  \int_0^z \f{ \Lpsi }{ 2 \vaa } 
  + \int_0^z \Lpsi \B( \f{G+1}{Q^z} - \f{1}{2\vaa} \B)  
 \teq \cLc(\ww) + \cNc( \ww, \eta) , 
\eal
\eeq
where $\vaa$ is the velocity of the 1D $\al$-profile in Theorem \ref{thm:1D_profile_prop}. We use subscript \textit{bd} to indicate that $\cLc$ is an operator on $r=0$, which is 
similar to a boundary.  We remark that in the above operators, $\ww$ is a function in $\R^3$ rather than an 1D function in $z$.

Recall the weight $\rhoc$ from \eqref{def:rhoc}. Our goal is to prove the following estimates. 

\begin{prop}\label{prop:3D_bd}

Let $\beps_7$ be as in Proposition \ref{prop:reg_solu}. There exists some absolute constant $\beps_8 \in ( 0, \beps_7]$, such that for any $\e \leq \beps_8$ and any $\ww$ 
satisfying \eqref{eq:3D_size:ww}, the solution $\eta$ to \eqref{eq:fix_recall} satisfies 
$ \rhoc \eta \in \cXc$ and 
\[
       \| \rhoc \eta(0, \cdot) \|_{\cXc} 
       \leq \f{1 + \lamcL } {2} || \ww ||_{\cXc}  +  C \nnr{\ww}   +   C \e^{1 - \kp} \nnrr{\ww}  \, ,
\]
where $\lamcL < 1$ is the parameter from Theorem \ref{thm:contra_lin_exact}.
\end{prop}

It is crucial that the second term $ C \nnr{\ww}$ is independent of the boundary part of $\ww$,
and depends only on $\pa_r \ww$, which enjoys a much stronger contraction estimate. 
See Section \ref{sec:idea_drW}. 
Since $\rhoc(\xx) = \ang \xx^{\e^2} \geq 1 $, $\eta(0,z)$ satisfies stronger estimates than $\ww(0,z)$ for large $|z|$. We use this property to prove compactness of $\cF_{\R^3}$ in Section \ref{sec:3D_solu}.

To distinguish the boundary part $\ww(0, z)$ and the whole solution $\ww(r, z)$, 
as in \eqref{eq:1D_normal}, we denote
\bseq\label{def:w1d}
\beq
   \wwd(z) = \om(0, z), \quad   
\psiod(\wwd)(z) = \psid(\wwd) - \pa_z \psid(\wwd)(0) z.
\eeq
\eseq

Below, we estimate two parts in $\cFc$ \eqref{eq:lin_bd_3D} separately. We show that $\cLc$ is a perturbation to $\cL$ 
in Theorem \ref{thm:contra_lin_exact}, which satisfies a contraction estimate.
The estimates for $\cNc$ are purely perturbative.

\subsubsection{Estimate of $\cLc$}

Recall $ \waa = \wwb $ from \eqref{def:3d_Om} and  $\Lpsi$ from \eqref{eq:lin_2Db}. We compare $\Lpsi(\ww) |_{r=0}$ and $\td \cR_{\mw{1D}}(\ww)$ in 
\eqref{eq:lin_1D_alpha}
\beq\label{eq:Lpsi_1D_rec}
 \Lpsi(\ww)(0,z) = - (1 - \al) \psio_z  -  2 \psio   \f{ \pa_z \waa }{\waa} ,\quad  \tcRd(\wwd)(z) = - (1 - \al) \pa_z \psiod - 2 \psiod  \f{2 \pa_z \waa}{\waa}.
\eeq

We decompose $\cLc$ as 
\begin{align}
  \cLc(\ww)  & =   \int_0^x \f{\tcRd(\wwd)}{2 \vaa}  
  + \int_0^x \f{ \cL_{\psi} - \tcRd(\wwd) }{ 2 \vaa } 
   \teq I_{ \eqref{eq:cLc_decomp},1 } (\ww) + I_{ \eqref{eq:cLc_decomp}, 2} (\ww).
   \label{eq:cLc_decomp}
\end{align}

Below, we estimate each term $ I_{ \eqref{eq:cLc_decomp},i }$.

\vs{0.1in}
\paragraph{\bf Estimate of $I_{ \eqref{eq:cLc_decomp},1 }$}
Since $I_{ \eqref{eq:cLc_decomp},i }$ integrates along $r =0$, $I_{ \eqref{eq:cLc_decomp},i }(\ww)$ only depends on $\wwd$ 
and we have $I_{ \eqref{eq:cLc_decomp},i } = |\waa|^{-1} \cLa(\wwd)$. 
Applying Theorem \ref{thm:contra_lin_exact}, 
using the $\cXc$-norm \eqref{norm:Xc}, 
and $\rhoc$ from \eqref{def:rhoc}, we obtain
\beq
\rhoc \max( \vpa, \muc \vpc ) |\waa I_{ \eqref{eq:cLc_decomp},1 } (w) | 
=  \rhoc \max( \vpa, \muc \vpc ) | \cLa(w) | 
\leq \rhoc \min( \lamcL,  C \e^{-1} \ang x^{- \kp_1 \e / 2} ) \nna{\wwd}.
\label{eq:lin_bd_3D:I1_1}
\eeq

Since $\rhoc = \ang x^{\e^2}$ satisfies 
\[
\rhoc \leq  1 + C\e^{1/2} , \ \forall \, \e^{3/2} \lgp x \leq 1,
\quad \rhoc  \ang x^{- \kp_1 \e / 2}
\les e^{- \kp_1 \e \lgp x  / 4} \les e^{ - C \e^{-1/2}} 
\les \e^{4}, \ \forall \, \e^{3/2} \lgp x > 1,
\]
combining the above two estimates, and using $\lamcL < 1$ by Theorem \ref{thm:contra_lin_exact},  we prove 
\beq
\| \rhoc \waa I_{ \eqref{eq:cLc_decomp},1 } (\ww) \|_{\cXc}
= \nlinf{ \rhoc \max( \vpa, \muc \vpc ) \waa I_{ \eqref{eq:cLc_decomp},1 }  (\ww) }  \leq (\lam_{\cL}  + C \e^{1/2} ) \| \wwd \|_{\cXc} \, .
\label{eq:lin_bd_3D:I1} 
\eeq

\paragraph{\bf Estimate of $I_{ \eqref{eq:cLc_decomp}, 2 }$}
We rewrite $\cL_{\psi}- \tcRd$ in \eqref{eq:Lpsi_1D_rec} as 
\[
  \cL_{\psi} - \tcRd(\wwd)= - (1 -\al)\B( ( \pa_z \psio - \f{1}{z} \psio ) - ( \pa_z \psiod - \f1z \psiod) \B)
  - \B( \f{ \psio }{z} - \f{\psiod}{z}    \B) ( 1 -\al + \f{ 2 z \pa_z \waa }{\waa} ) .
\]

Since $\psi$ is the steam function in $\R^3$ associated with $\ww$, 
and since $\psid(\wwd)$ is the stream function in $\R$ associated with $\ww(0,z)$ from \eqref{def:w1d}, using Lemma \ref{lem:vel_bc} and  triangle inequality, we bound 
\[
\bal
 |( \pa_z \psio - \f{1}{z} \psio ) - ( \pa_z \psiod - \f{1}{z} \psiod  )|
 & \les \min(|z|,  \ang z^{ \al- \hal + \epa }  |\cJaa(z)|^{\kp} ) \nnr{\ww}, \\
  |\tf{1}{z} (\psio -  \psiod) |
  & \les   \min(|z|, 
  |\cJaa(z)|^{1+\kp} ) \nnr{\ww}.
\eal
\]

Using \eqref{eq:1D_prof_est2:c}, $1-\al = \f{2}{3} + \e$, we bound 
\[
   \B| 1 -\al + 2 \f{z \pa_z \waa}{\waa} \B| \les 
      \B| 2( \f{1}{3} +  \f{z \pa_z \waa}{\waa} ) \B|
+ \e 
\les  |\lgp z|^{\kp-2} + \e .
  \]

Combining the above two estimates, we obtain
\[
\bal
  |\cL_{\psi} - \tcRd(\wwd)|
& \les \min( |z|, \ang z^{\al-\hal + \epa}  \cJaa^{\kp} + \cJaa^{\kp+1} ( \e + |\lgp z|^{ \kp - 2} )   ) \nnr{\ww} .
\eal
\]

Since $\cJaa(z) \les \lgp z$  by \eqref{eq:Ja_hat}, $\kp < 1$ by \eqref{def:kp}, and $ \vaa \gtr |z|\cJaa$ \eqref{eq:vaa_low}, using the above estimate of 
$\cL_{\psi} - \td \cR(\wwd)$, we bound $I_{ \eqref{eq:cLc_decomp}, 2 }$ as 
\[
\bal
  & I_{ \eqref{eq:cLc_decomp}, 2 } \teq  \int_0^x \f{   |\cL_{\psi} - \td \cR(\wwd)| }{ 2 \vaa }(z) d z \\
 &  \les \B( \one_{x<2} x  + 
  \one_{x > 2} \int_{2}^x \f{1}{z \cJaa } \big( \ang z^{\al-\hal + \epa} \cJaa^{\kp} +  \e \cJaa^{\kp+1} +
  \cJaa^{\kp+1} |\lgp z|^{\kp-2} 
   \big) d z \B) \nnr{\ww} \\
    &  \les \B( \one_{x<2} x + 
  \one_{x > 2} \int_{2}^x  \ang z^{\al-\hal + \epa-1} \cJaa^{\kp-1} +  \f{\e \cJaa^{\kp} }{\ang z}
   + \f{ |\lgp z|^{2 \kp- 2} }{\ang z}
    d z \B) \nnr{\ww} \, . \\
\eal
\]

To estimate the first integrand, since $\al -\hal + \epa \leq - \e / 2, \kp > 0$ by \eqref{ran:ep_all}, \eqref{def:kp}, we use Lemma \ref{lem:log_ineq_J} with $ (c, k) \rsa ( \f12,  \kp - 1)$.
For the second integrand, we use $\cJaa(z) \leq \cJaa(x)$ for $z \leq x$. 
For the third integrand, we use $2 \kp - 2 > -1$ \eqref{def:kp}, a change of variable $z \to e^s, s \in [\log 2 , \log x]$, and estimate the integral directly.
Then, we  estimate $I_{ \eqref{eq:cLc_decomp}, 2 }$ as 
\[
  \bal
I_{ \eqref{eq:cLc_decomp}, 2 } & \les  \min\B( x,  \   \cJaa^{\kp}  + \e \lgp x \cdot \cJaa^{\kp} 
  + |\lgp x|^{2 \kp- 1}  \B) \nnr{\ww} .
\eal
\]

Using $\cJaa(x) \les \lgp x$, $2 \kp - 1 \in (0, \kp)$ by \eqref{def:kp}, and Lemma \ref{lem:lgx_pow} with $ \tf12 \epa  \gtr \e$ \eqref{ran:ep_all}, we prove
\[
  | I_{ \eqref{eq:cLc_decomp}, 2 } | \les \min( x, \ang x^{\epa/2} \cJaa^{\kp} ) \nnr{\ww} .
 \]

Recall $\rhoc = \ang x^{\e^2}$ from \eqref{def:rhoc}. 
Since $\e^2 < \epa/4$  by \eqref{norm:Xc}, using estimate \eqref{eq:wg_equiv} on the weights with $\bbb \in (1.1,2)$ \eqref{norm:Xc} and the above estimates,  we prove 
\begin{align}
   \rhoc \max( \vpa, \muc\vpcc) |\waa|  | I_{ \eqref{eq:cLc_decomp}, 2 } | 
 &  \les \rhoc ( |x|^{-\bbb+1} + \ang x^{-\epa} \cJaa^{-\kp})\min( x, \ang x^{\epa/2} \cJaa^{\kp} ) \nnr{\ww} \notag  \\
 & \les  \rhoc \ang x^{-\epa/2} \nnr{\ww } \les \ang x^{-\epa/4} \nnr{\ww}.
  \label{eq:lin_bd_3D:I2}
\end{align}

We emphasize that the norm $\nnrr{\ww}$ and $\nnr{\ww}$ are \emph{different}
and $\nnr{\ww}$ \emph{does} not depend on the value $\ww(0, z)$ on $r=0$.

\subsubsection{Estimate of $\cNc$}

We estimate the nonlinear term $\cNc$  perturbatively. 
Recall $\cNc$ from \eqref{eq:lin_bd_3D}
\[
\cNc(\ww)(x) =   \int_0^x \Lpsi \B( \f{G+1}{Q^z} - \f{1}{2\vaa} \B)(z) d z  ,
\quad \Lpsi = - (1 - \al) \psio_z  -  2 \psio   \f{ \pa_z \wwb }{\wwb} . 
\]

Using \eqref{eq:1D_prof_est2:c}, $\f{1-\al}{2} = \f{1}{3} + O(\e)$, and Proposition \ref{prop:vel_est}, we bound
\beq\label{eq:lin_bd_3D:non1}
\bal
  |\Lpsi(0,z)| & =\B| - (1 -\al) ( \pa_z \psio - \f{1}{z} \psio )
  - \f{ \psio }{z} ( 1 -\al + \f{ 2 z \pa_z \waa}{\waa} ) \B|(0,z) \\
  & \les  \min\B( |z|^{1+\al} ,  ( \ang z^{\al-\hal + \epa} \cJaa^{\kp} + |\cJaa(z ) |^{\kp+1} \cdot (\e + |\lgp z|^{\kp-2}) ) \B) \nnrr{\ww} , \\
\eal
\eeq
with $\kp < 1$ \eqref{def:kp}. Recall $Q^z$ from \eqref{def:Q}. 
From the normalization conditions \eqref{eq:normal_cond}, \eqref{def:psio}, 
and the identity \eqref{eq:iden_r0},  we obtain 
\[
  Q^z(0, z) = c_l z + 2 \Psi =  ( 2 -  2 \Psi_z(0) ) z   + 2 \Psi 
  = \bar Q^z(0, z) + \td Q^z(0, z) =  2 \vaa + 2 \psio  \, .
\]

Using Proposition \ref{prop:vel_est}, 
$\vaa \gtr |z| \cJaa$ by \eqref{eq:vaa_low},  $\cJaa \les \e^{-1}$ \eqref{eq:Ja_hat}, 
and the assumption \eqref{eq:3D_size:ww} on $\ww$,  for $\e$ small enough, we obtain
\[
  |Q^z - 2 \vaa| 
  \les |\psio|
  \les  |z| \cJaa^{\kp+1} \nnrr{\ww}
  \les   |z| \cJaa \e^{-\kp} \cdot \e^{ (1 + \kp)/2 } 
    \les   |z| \cJaa \cdot \e^{ \mhk } \les  \vaa \cdot \e^{\mhk} < \vaa \, ,
\]
which implies 
\[
 Q^z(0, z) \geq \vaa(z) \gtr |z| \cJaa(z) .
\]

Thus, we estimate
\beq\label{eQ:lin_bd_3D:non2}
  |\f{G+1}{Q^z} - \f{1}{2 \vaa}|
  = |\f{G}{Q^z} + \f{2 \vaa - Q^z}{2 \vaa Q^z}|
\les \f{|G | }{|z| \cJaa(z) } + \f{ \e^{ \mhk } }{ |z| \cJaa(z) }.
\eeq

From the definition of $\cXc$ in \eqref{norm:Xc}, 
$\rhoc(z)^{-1} \rhoo^{-1} \les \ang z^{\epa-\e^2} |\cJaa(z)|^{\kp}$ by \eqref{def:3d_wg}, \eqref{def:rhoc}, and $G = \f{\eta}{\waa}$,
we obtain
\beq\label{eq:lin_bd_3D:non3}
  |G(z)| = |\eta |\waa|^{-1}| \les  \| \rhoc \eta \|_{\cXc} 
  \ang z^{\epa- \e^2} |\cJaa(z)|^{\kp} .
\eeq

\paragraph{\bf Estimate integrals}
Combining  \eqref{eq:lin_bd_3D:non1}-\eqref{eQ:lin_bd_3D:non2}, we prove 
\beq
\bal
II & \teq |\Lpsi \B( \f{G+1}{Q^z} - \f{1}{2\vaa} \B)  |  \les \min( |z|^{1+\al} ,  h(z) )  \nnrr{\ww} \cdot  \f{ |G| + \e^{ \mhk } }{ |z| \cdot |\cJaa(z)| } , \\
h(z)  & \teq % \ang z^{\al-\hal + \epa}  + \e \cJaa(z ) + \cJaa |\lgp z|^{-4/3} .
 \ang z^{\al-\hal + \epa} \cJaa^{\kp} + |\cJaa(z ) |^{\kp+1} \cdot (\e + |\lgp z|^{\kp-2}) 
\teq h_{\eqref{eq:lin_bd_3D:non32}, 1 } + h_{\eqref{eq:lin_bd_3D:non32}, 2 } .
\label{eq:lin_bd_3D:non32}
\eal
\eeq

Using the bound of $G$ from \eqref{eq:lin_bd_3D:non3}, we obtain
\beq\label{eq:lin_bd_3D:non_II1}
  \int_0^{x \we 1 } II \les \int_0^{x \we 1} z^{\al}   \nnrr{\ww}   ( \|\rhoc \eta \|_{\cXc} + \e^{ \mhk } )
  \les (x \we 1)^{1+\al} \nnrr{\ww}   ( \| \rhoc \eta \|_{\cXc} +\e^{ \mhk } ).
\eeq

Next, we estimate the case with $x > 1$. We decompose $II$ into the constant part $II_C$ and the part $II_G$ involving $G$ :
\[
  \int_1^x II(z) d z  \les   (II_C(x) + II_G(x)) \nnrr{\ww},
    \  II_C(x)  =\int_1^x   \f{ h(z)  \e^{\mhk}  }{ |z| \cdot |\cJaa(z)|  } , 
  \  II_G(x) = \int_1^x    \f{ h(z) |G|  }{ |z| \cdot |\cJaa(z)|  } .
\]

For $II_C$, since $\cJaa(z) \les \lgp z, \e^{-1}$ \eqref{eq:Ja_hat}, we have $h(z) \les |\cJaa(z)|^{\kp}$ 
by \eqref{eq:lin_bd_3D:non32}. First using $\ang z \leq \ang x$ and then using 
Lemma \ref{lem:log_ineq_J} with $k\rsa \kp-1 > -1$, $\epa/ 10 \gtr \e$ and $\kp > 0$, we obtain
\beq\label{eq:lin_bd_3D:non_C}
 II_C(x)   \les \e^{ \mhk } \ang x^{ \f{\epa}{10} } \int_1^x  \ang z^{- \f{\epa}{10} -1} \cJaa(z)^{ \kp -1} d z 
  \les \e^{ \mhk } \ang x^{ \f{\epa}{10}} \cJaa^{\kp } .
\eeq

For $II_G$, using \eqref{eq:lin_bd_3D:non3} for $G$ and $1 \leq \ang x / \ang z $ we bound 
\beq
\bal
   II_G(x)
& \les  \| \rhoc \eta \|_{\cXc}
 \int_1^x  \f{  h(z) \cdot \cJaa^{\kp} }{ z \cJaa  } \ang z^{\epa-\e^2} 
 \B( \f{\ang x}{\ang z} \B)^{ \epa/2 } d z  \\
 & \les   \ang x^{\epa/2}  \| \rhoc \eta \|_{\cXc}
 \int_1^x    (  h_{\eqref{eq:lin_bd_3D:non32}, 1 }(z) + h_{\eqref{eq:lin_bd_3D:non32}, 2 }(z) )  \cdot  |\cJaa(z)|^{\kp-1}  \ang z^{\epa/2-1 - \e^2}  d z .
\eal
 \label{eq:lin_bd_3D:non_G1}
\eeq

Next, we estimate the integral of each part in $h(z)$.

Since $( \al -\hal + \epa) + \epa / 2 -\e^2 <   - \f{\e}{10}$ by \eqref{ran:ep_all}, using
Lemma \ref{lem:log_ineq_J} with $k \rsa \kp - 1 $, we bound 
\[
\bal
\int_1^x h_{\eqref{eq:lin_bd_3D:non32}, 1 }  \cdot  |\cJaa(z)|^{\kp-1}  \ang z^{\epa/2-1 - \e^2}  d z 
\les 
\int_1^x \ang z^{ ( \al-\hal + \epa) + \epa/2 -1 -\e^2} \cJaa^{ 2 \kp - 1 } d z
& \les  \cJaa^{2 \kp}(x).
\eal
\]

For the integral involving $h_{\eqref{eq:lin_bd_3D:non32}, 2 } $  in \eqref{eq:lin_bd_3D:non_G1},
since $\cJaa$ is increasing and $\epa -\e^2 \asymp \e$ by \eqref{ran:ep_all}, we estimate 
\[
\bal
  & \int_1^x  h_{\eqref{eq:lin_bd_3D:non32}, 2 } \cdot \cJaa^{\kp-1} \ang z^{\epa/2-1-\e^2} d z 
 = \int_1^x  ( \e  +  |\lgp z|^{ \kp - 2 } ) 
  |\cJaa(z )|^{\kp + 1} \cdot \cJaa^{\kp-1} \ang z^{\epa/2-1-\e^2} d z  \\
& \quad \les |\cJaa(x)|^{ 2 \kp } \B( \int_1^x  \e \ang z^{\epa/ 2 - 1 -\e^2} +   \ang x^{\epa/2-\e^2} \ang z^{-1} |\lgp z|^{ \kp-2}    d z \B)  \les \cJaa^{ 2 \kp } \ang x^{\epa/2 -\e^2} ,
\eal
\]
where the integral for $ z^{-1} |\lgp z|^{ \kp-2}$ is bounded by a constant since $\kp < 1$ \eqref{def:kp}.

Thus, plugging the above estimates of integrals in \eqref{eq:lin_bd_3D:non_G1}, we prove 
\beq\label{eq:lin_bd_3D:non_G}
\bal
   II_G(x)  & \les \ang x^{\epa/2 }  \| \rhoc \eta \|_{\cXc} \cdot \cJaa^{2 \kp} \ang x^{\epa/2 -\e^2} 
    \les \ang x^{\epa -\e^2} \cJaa^{2\kp}  \| \rhoc \eta \|_{\cXc} .
  \eal
\eeq

Since $\cJaa \gtr 1, \ang x \gtr 1$, summarizing \eqref{eq:lin_bd_3D:non_II1}, \eqref{eq:lin_bd_3D:non_C}, and \eqref{eq:lin_bd_3D:non_G}, we prove 
\[
\bal
  |\cNc(x)| & \leq \int_0^x II(z) d z   \les (II_C(x) + II_G(x) )  \nnrr{\ww}  \\
&  \les  \min( x , 1 )^{ \al + 1}
  \B(  \ang x^{\epa-\e^2}  \|\rhoc \eta \|_{\cXc} 
  \cJaa^{\kp}
  + \e^{ \mhk }  \ang x^{\epa/10}   \B) 
   \cJaa^{\kp}    \nnrr{\ww}   .
\eal
\]
 Recall $\rhoc = \ang x^{\e^2}$ from \eqref{def:rhoc}. 
 Using estimate \eqref{eq:wg_equiv} on the weights with $\bbb \in (1.1,2)$ \eqref{norm:Xc}, 
 $\e^2 < \epa / 4$ \eqref{norm:Xc}, $\cJaa \les \e^{-1}$, 
  assumption \eqref{eq:3D_size:ww}, 
and $\nnrr{w} =\max( \nnr{w} , \nnla{w} )$ and $\nnla{w} \les \| w\|_{\cXc}$ by \eqref{norm:comp_bd}, we prove
 \beq
 \bal
&   \rhoc \max( \vpa, \muc\vpcc) |\waa \cNc| 
 \les \rhoc ( |x|^{-\bbb+1} + \ang x^{-\epa} \cJaa^{-\kp}) |\cNc| \\
& \qquad \les \rhoc \ang x^{-\epa}  \B(  \ang x^{\epa-\e^2}  \|\rhoc \eta \|_{\cXc} 
  \cJaa^{\kp}
  + \e^{ \mhk }  \ang x^{\epa/10}   \B)   \nnrr{\ww} \\
  & \qquad \les ( \e^{-\kp}   \|\rhoc \eta \|_{\cXc}  \e^{ \hk}
  + \e^{\mhk} \nnrr{\ww}  ) \les
   \e^{\mhk} ( \|\rhoc \eta \|_{\cXc}   +  \nnr{\ww} + \| \ww \|_{\cXc}  ).
  \eal
  \label{eq:lin_bd_3D:non}
 \eeq
While the norm $\nnrr{\cdot}$ \eqref{def:3d_norm} involves the boundary part, 
due to the small parameter $\e^{\mhk}$, we treat it perturbatively.

\subsubsection{Proof of Proposition \ref{prop:3D_bd} }

From \eqref{eq:lin_bd_3D} and \eqref{eq:cLc_decomp}, we have 
$\cFc =\cLc + \cNc = I_{ \eqref{eq:cLc_decomp},1 } (\ww) + I_{ \eqref{eq:cLc_decomp}, 2 } (\ww)
+  \cNc $. Combining 
  \eqref{eq:lin_bd_3D:I1}, \eqref{eq:lin_bd_3D:I2}, 
    \eqref{eq:lin_bd_3D:non}, and using the $\cXc$-norm \eqref{norm:Xc}, we prove 
\beq
\bal
\|  \cFc(\ww, \eta) \rhoc  |\waa|   \|_{\cXc} 
& \leq 
(\lam_{\cL}  + C_{ \eqref{eq:lin_bd_3D:pertb1} } \e^{\mhk} ) \| \ww \|_{\cXc} + C_{ \eqref{eq:lin_bd_3D:pertb1} } ( \nnr{\ww}   + \e^{\mhk} \|  \rhoc \eta \|_{\cXc} ) \, ,
\eal
\label{eq:lin_bd_3D:pertb1}
\eeq
for some absolute constant $C_{ \eqref{eq:lin_bd_3D:pertb1} } >0$.

Since $G = \f{\eta}{\waa}$, from \eqref{eq:lin_bd_3D0} and \eqref{eq:lin_bd_3D},  we obtain
$  \eta =  \waa  \cFc(\ww, \eta)$.

Since we do not know a-priori $\rhoc \eta \in \cXc$, we cannot apply estimate \eqref{eq:lin_bd_3D:pertb1} directly to bound $\eta$. 
 Instead, we first obtain a local estimate of $\eta$ for $|x| \leq R$. 
Since the integration in \eqref{eq:lin_bd_3D} is from $0$ to $x$, 
and it only depends on $\eta(y)$ for $y \leq |x|$,  for any $R \geq 0$, we have 
\beq\label{eq:lin_bd_3D:eta2}
\bal
 & \one_{|x|\leq R}\eta(0, x)  = \one_{|x| \leq R} \cF(\ww, \eta)(x)
 =  \one_{|x| \leq R} \cFc(\ww,  \one_{|\cdot|\leq R} \eta(0, \cdot) )(x).
 \eal 
\eeq

Since Proposition \ref{prop:reg_solu} implies $|\eta(x) | \les |x|^{1 + \al}$, 
and $|x|^{1 + \al} |x|^{-\bbb} |x|^{-\e^2} \les C(R)$ for $|x| \leq R$, we obtain $\rhoc \eta(0, x) \one_{|x|\leq R} \in \cXc$ for any $R \geq 0$.  Since $\| f \one_{|x|\leq R} \|_{\cXc} \leq \| f \|_{\cXc}$ for any $f$, applying estimate \eqref{eq:lin_bd_3D:pertb1} with $\eta \rsa \eta \one_{|\cdot| \leq R}$ and using 
 identity \eqref{eq:lin_bd_3D:eta2}, 
 for any $R \geq 0$,  we prove
\beq
\bal
   \| \one_{ |x|\leq R } \rhoc \eta(0, \cdot)  \|_{\cXc} & \leq  
    (\lamcL + C_{ \eqref{eq:lin_bd_3D:pertb1} }  \e^{\mhk} ) \| \ww \|_{\cXc} +  
C_{ \eqref{eq:lin_bd_3D:pertb1} } (   \nnr{\ww}+
 \e^{ \mhk } \| \one_{ |x|\leq R } \rhoc \eta(0, \cdot) \|_{\cXc} ) .
\eal
\label{eq:lin_bd_3D:pertb}
\eeq

Since $\lamcL < 1$ from Theorem \ref{thm:contra_lin_exact} and $\kp < 1$ by \eqref{def:kp},
by choosing $\e$ small enough so that $\e^{ \mhk }$ is very small, and solving the inequality \eqref{eq:lin_bd_3D:pertb} on $\eta$, we  prove 
 \[
   \| \one_{ |x|\leq R } \rhoc \eta(0, \cdot)  \|_{\cXc}  
    \leq \tf{1}{2} (1 + \lamcL ) || \ww ||_{\cXc}  +  C \nnr{\ww}   ,
 \]
for some absolute constant $C$. 
Since the above estimate holds uniformly for any $R \geq 0$, taking $R \to \infty$, we prove  Proposition \ref{prop:3D_bd}.

\subsection{Damping terms in weighted $W^{1,\infty}$ estimates}\label{sec:C1_damp}

For $\xx = (r, z)$ with $z \neq 0$, using \eqref{eq:traj_smooth}, we obtain 
\bseq\label{eq:lin_dG}
\begin{align}
    \tf{d}{ds} ( \pa_r G ) \cc X(s, \xx) &=
  \big[ ( \bar \cR + \Lpsi - \pa_r  Q^r ) \cdot \pa_r G
 + \cmr \big] \cc  X(s, \xx)  \label{eq:lin_dG:Gr} ,
  \\
    \tf{d}{ds} ( \pa_z G ) \cc X(s, \xx) & = 
    \big[  ( \bar \cR + \Lpsi - \pa_z Q^z ) \cdot \pa_z G + \cmz  \big] \cc X(s, \xx) , \label{eq:lin_dG:Gz} 
\end{align}
$\cmr, \cmz$ denote the terms not involving $\pa_r G, \pa_z G$, respectively.
\begin{align}
  \cmr & =  - \pa_r Q^z \cdot \pa_z G 
  + \pa_r \bar \cR \cdot (G + 1) 
  + \pa_r \cL_{\psi} \cdot (G + 1) 
  \label{eq:lin_dG:Mr} \, , \\
  \cmz & =   - \pa_z Q^r \cdot \pa_r G 
  + \pa_z \bar \cR \cdot (G + 1)
  + \pa_z \cL_{\psi} \cdot (G+1) \, . 
    \label{eq:lin_dG:Mz}
\end{align}
\eseq
From Proposition \ref{prop:reg_solu}, the terms in \eqref{eq:lin_dG} are continuous in $\xx, s$ away from $z=0$.

\begin{remark}[\bf Nonlocal terms depend on $\om$]
The nonlocal terms $\QQ = \bar \QQ + \td \QQ(\ww)$ \eqref{def:Q} and $\Lpsi$  \eqref{eq:lin_2D} depend on the input $\om$. To simplify the notations, we omit this dependence. 
\end{remark}

We aim to apply the weighted trajectory estimate in Lemma \ref{lem:traj_est} 
to $\pa_r G, \pa_z G$ with the weight $m = \rhoc \rhor, \rhoc \rhoz$, which are stronger than 
$\rhor, \rhoz$. See the definitions of these weights in \eqref{def:3d_wg} 
and \eqref{def:rhoc}. We first estimate the damping terms, which correspond to 
$\f{\QQ \cdot \na ( m \vp^{-1}) }{m \vp^{-1} }$ in Lemma \ref{lem:traj_est}.

In the following key lemmas, we use the transport effects and design weights to 
generate strong damping terms, and will use them to control the nonlocal terms and error terms $\bar \cR$.

\begin{lem}\label{lem:3D_damp}
Recall the weights $\rhor = \rhoo(\xx)  \rag(\xx) \cdot \la \xx \ra , 
\rhoz = \rhoo(\xx) \rag \cdot |z|$ from \eqref{def:3d_wg}, $\rhoc$ from \eqref{def:rhoc}
and parameters $\kp, \kp_1$ from \eqref{def:kp}. Suppose that 
$\ww$ satisfies \eqref{eq:3D_size:ww}. 
Let $\beps_8$ be as in Proposition \ref{prop:3D_bd} and 
 $\beps_{\mw{pow}}$ be as in Lemma \ref{lem:interp1}.  There exists $\beps_9\in (0, \min( \beps_8 , \beps_{\mw{pow}})] $ small enough and a  Lipschitz weight  $\rkey(z) > 0$, which is even in $z$ and satisfies $\pa_z \rkey(z) \leq 0$, $\rkey(z) \asymp 1$ with constant independent of $\e$, 
such that the following statements hold. 
For any $\e \in (0, \beps_9]$, we have
\bseq
 \begin{align}
 A^r & \teq  \f{\QQ \cdot \na ( \rhoc \rhor \rkey )}{ \rhoc \rhor \rkey }  + \bar \cR 
 + \Lpsi - \pa_r Q^r 
 \leq - \bar \lam_1 \B( \f{1}{\e} \ang \xx^{\al-\hau}  + \cJaa(\xx) |\lgp z|^{ -\kp_1 -1} +1 \B), \label{eq:damp:Ar} \\
 A^z & \teq \f{\QQ \cdot \na ( \rhoc \rhoz \rkey )}{ \rhoc \rhoz \rkey }  + \bar \cR + \Lpsi- \pa_z Q^z  \leq -\bar \lam_1 \B(  \f{r}{\ang \xx} \big(  \f{1}{\e} \ang \xx^{\al-\hau} 
 +  \cJaa(\xx)   |\lgp z|^{-\kp_1 - 1} \big)   +1 \B) , \label{eq:damp:Az}
 \end{align}
 \eseq
for some constant $ \bar \lam_1 = \lam(\kag, \wwwa) > 0$ independent of $\e$.
\end{lem}

The coefficients $\bar \cR + \Lpsi - \pa_r Q^r , \bar \cR + \Lpsi- \pa_z Q^z$ arise from the damping coefficients in \eqref{eq:lin_dG}. The extra damping terms come from the weights $\rhor \rkey, \rhoz \rkey$.
The contribution from $\rhoc$ is treated perturbatively.

Since $\rkey$ only depends on $z$ 
we decompose $A^r, A^z$ as follows 
\beq\label{eq:decomp:damp}
\bal
 & \f{\QQ \cdot \na ( \rhoo \rag \rkey \rhoc) }{ \rhoo \rag \rkey \rhoc }
 = \f{\QQ \cdot \na \rhoo}{\rhoo }
  + \f{\QQ \cdot \na \rag}{\rag} + \f{Q^z \pa_z \rkey }{\rkey}
  + \f{\QQ \cdot \na \rhoc}{\rhoc} \teq I_R + I_{\mw{ag}} + I_{z} + I_{c}, \\
  A^r & = \f{\QQ \cdot \na ( \rhoo \rag \rkey \rhoc) }{ \rhoo \rag \rkey \rhoc }
+ \bar \cR + \Lpsi + (  \f{\QQ \cdot \na \ang \xx}{\ang \xx}  - \pa_r Q^r) \teq I_R + I_{\mw{ag}} + I_{z} + I_{c} + \bar \cR + \Lpsi + II_r, \\ 
  A^z &= \f{\QQ \cdot \na ( \rhoo \rag \rkey \rhoc) }{ \rhoo \rag \rkey \rhoc }
+ \bar \cR + \Lpsi +  (  \f{\QQ \cdot \na z}{z}  - \pa_z Q^z)  \teq I_R + I_{\mw{ag}} + I_{z} +I_c+ \bar \cR + \Lpsi+ II_z,
\eal
\eeq
where $I_R, I_{\mw{ag}}, I_{z}, I_c$ denotes the terms from the radial, angular, $z-$weights,
and $\rhoc$ respectively. Note that $A^r, A^z$ only differs in the last term $II_r, II_z$. We have estimated $\bar \cR$ 
in Proposition \ref{prop:3D_error}.

In the following estimates, we shall extract the main damping terms and treat all the lower order terms perturbatively to the main damping terms.

\subsubsection{Estimate of radial weights $I_R, I_c$}

Since $\rhoo(\xx) =
 (1 + |\xx|^{-\bbb+1}) \ang \xx^{-\epa} |\cJaa(\xx)|^{-\kp} $ \eqref{def:3d_wg}
 only depends on the radial variable $|\xx|$, using the definition \eqref{def:3d_wg}, 
$(Q^r, Q^z) = (q^r r, q^z z)$,   we obtain
\[
\bal
   I_R  & = \f{\QQ \cdot \na \rhoo}{\rhoo}
   = \f{ q^r r^2 + q^z z^2 }{|\xx|^2} \cdot  \f{ |\xx| \pa_{|\xx|} \big( ( 1 + |\xx|^{-\bbb+1}) \ang \xx^{-\epa} |\cJaa(\xx)|^{-\kp} \big) }{ (1 + |\xx|^{-\bbb+1}) \ang \xx^{-\epa} |\cJaa(\xx)|^{-\kp} }
   \\
 &  =  \f{(q^r r^2 + q^z z^2)}{|\xx|^2} \cdot 
 \B(  - \epa \f{|\xx|^2 }{\ang \xx^2}
  - (\bbb - 1) \f{  |\xx|^{1-\bbb} }{1 +|\xx|^{1-\bbb} } 
  + \f{ |\xx| \pa_{|\xx| } |\cJaa(\xx)|^{-\kp} }{|\cJaa(\xx)|^{-\kp}} 
   \B)   .
\eal
\]

Using \eqref{eq:traj_outgo}, we obtain
\[
    q^r r^2 + q^z z^2  = \QQ \cdot \xx \geq C \cJaa(\xx) |\xx|^2 >0.
\]

Since  $|\cJaa(\xx)|$ is increasing in $|\xx|$ from the discussion in \eqref{eq:wa}, using the above sign, we obtain
\[
  I_R \leq \f{(q^r r^2 + q^z z^2)}{|\xx|^2} \cdot 
 \B(  - \epa \f{|\xx|^2 }{\ang \xx^2}
  - (\bbb - 1) \f{  |\xx|^{1-\bbb} }{1 +|\xx|^{1-\bbb} } 
   \B)
  =  -\f{(q^r r^2 + q^z z^2)}{|\xx|^2} \cdot 
 \B(   \epa \f{|\xx|^2 }{\ang \xx^2}
  +  \f{  \bbb - 1 }{1 +|\xx|^{\bbb - 1} }  
   \B) .
\]

Since $\bbb - 1 > \f{1}{100}$ \eqref{norm:Xc} and $ \epa \f{|\xx|^2}{\ang \xx^2} + \f{ \bbb-1 }{ 1 + |\xx|^{\bbb-1} }
\gtr (\bbb-1) \ang \xx^{1-\bbb} + \epa$, and $\epa \gtr \e $ \eqref{norm:Xc}, we obtain
\beq\label{eq:3D_damp_IR}
  I_R \leq  - \bar C_R \cJaa(\xx) \cdot ( (\bbb-1) \ang \xx^{1-\bbb} + \e ).
\eeq
for some absolute constant $\bar C_R > 0$.

For $I_c$, recall $\rhoc(\xx) = \ang \xx^{\e^2} $ from \eqref{def:rhoc}. Using $\na \ang \xx^{\e^2}= \e^2  \xx \ang \xx^{\e^2 -2} $ and \eqref{eq:3D_boot_q:1}, we obtain
\beq\label{eq:3D_damp_Ic}
\bal
  |I_c| &=  |\f{\QQ \cdot \na \rhoc(\xx)}{ \rhoc( \xx ) }| 
    \leq  \e^2 \f{ |\QQ \cdot \xx| }{ \ang \xx^2 } 
  \les \e^2 (|q^r| + |q^z|)  
  \les \e^{-1} \cdot \e^2 \les \e . \\
\eal
\eeq

\subsubsection{Estimate of angular weights}\label{sec:damp_ag}

Recall $\QQ = (q^r r , q^z z)$ from \eqref{def:qq}. Using the definition \eqref{def:rag}, we obtain
\beq\label{eq:Iag_comp1}
\bal
 I_{\mw{ag}} &= \kag  \f{ (  r q^r \pa_r + z q^z \pa_z ) \f{\ang z}{ r + \ang z } }{ \f{\ang z}{ r + \ang z }  }
 = \kag \cdot  \f{   - q^r  \cdot \f{r \ang z}{ (r + \ang z)^2 }
 + z q^z \cdot \f{ \pa_z \ang z \cdot ( r + \ang z - \ang z ) }{  (r + \ang z)^2 } 
 }{ \f{\ang z}{ r + \ang z }  }   \\
 & = \kag \cdot \B( - q^r \f{r}{r + \ang z} + q^z  \cdot \f{r}{r + \ang z} \cdot \f{z^2}{\ang z^2}  \B)
 = - \kag \cdot \f{r}{r + \ang z} ( q^r - q^z + q^z \f{1}{\ang z^2} ) .
\eal 
\eeq

Using the lower bounds of $q$ \eqref{eq:Q_lower}, \eqref{eq:3D_boot_q:1},  $\nnn{\ww} \leq \e^{ \hk }$ by \eqref{eq:3D_size:ww},
and $|r| \leq |\xx|, \epa > 0$ \eqref{ran:ep}, we obtain
\[
\bal
q^r - q^z + q^z \f{1}{\ang z^2}
  & \geq    2 \bar C_Q \f{1}{\e} \ang \xx^{\al-\hau} -
     C \e^{-\hk} \ang \xx^{\al-\hal + \epa} 
  - C \ang \xx^{\al- \hal} \\
  & \qquad + \tf32 \ang z^{-2}  ( \bar C_Q \cJaa - C r \min( |\xx|^{\al-\hal} , |\xx|^{\al-\hal -1} ) )   \\
  & \geq  2 \bar C_Q \f{1}{\e} \ang \xx^{\al-\hau} 
  + \tf32 \bar C_Q  \ang z^{-2}   \cJaa(\xx) - C \e^{ -\hk } \ang \xx^{\al-\hal + \epa} .
\eal
\]

Since $\al < \alb$ and $\mhk > 10 \kp_1$ by \eqref{def:kp}, using \eqref{eq:interp:pow} in Lemma \ref{lem:interp1} with $ a =
 \e$ and requiring $\e \leq \beps_{\mw{pow}} $, we obtain
\beq\label{eq:switch_pow}
  \e^{ - \hk }  \ang \xx^{\al-\hal + \epa} 
\leq C \e^{- \hk - \kp_1} \ang \xx^{\al-\hau} + C\e^{1 -\hk} 
\leq \f{1}{2} \bar C_Q \e^{-1} \ang \xx^{\al-\hau} + C \e^{ \mhk }.
\eeq
Using the above  estimates, we obtain
\beq\label{eq:damp_ag:est1}
   q^r - q^z + q^z \f{1}{\ang z^2} \geq 
   \f{3}{2 \e } \bar C_Q  \ang \xx^{\al-\hau} 
  + \f{3}{2} \bar C_Q  \ang z^{-2}   \cJaa(\xx) - C \e^{\mhk}.
\eeq

Since 
\[
   \ang \xx^2 = r^2 + \ang z^2 \geq \f{1}{2} (r + \ang z)^2 >( \f{2}{3} (r + \ang z) )^2,  \ \Longrightarrow  \, 
    \f{r}{r + \ang z} \geq \f{2 r}{3 \ang \xx} \, , 
  \]
  combining the above bound and using \eqref{eq:Iag_comp1}, we prove 
\beq\label{eq:3D_damp:ag}
  I_{\mw{ag}} 
  \leq -\kag \f{r}{\ang \xx} \B(  \bar C_Q \f{1}{\e} \ang \xx^{\al-\hau} 
  +    \bar C_Q   \ang z^{-2} \cJaa(\xx)  \B) + C \e^{ \mhk} .
\eeq

Thus, when $r \gtr |z|$, the angular weight generates a large damping terms of order $\e^{-1} \ang \xx^{\al-\hau}$.
Moreover, if $|z| \les 1$ and $r \gtr |z|$, we gain a damping term $\cJaa(\xx) \sim \min( \lgp x, \e^{-1} )$.

\subsubsection{Estimate of $II_r, II_z$}

\textbf{ Estimate of $II_r$} Recall $II_r, II_z$ from \eqref{eq:decomp:damp}
\[
   II_r =    \f{\QQ \cdot \na \ang \xx}{\ang \xx}  - \pa_r Q^r,  \quad  II_z = \f{\QQ \cdot \na z}{z}  - \pa_z Q^z . 
\]

Using $Q^r = r q^r, Q^z = z q^z$, we estimate 
\[
  II_r =  \f{q^r  r^2 + q^z z^2}{ \ang \xx^2 } - r \pa_r q^r - q^r
  = \f{  - (q^r -q^z) ( z^2 + 1) - q^z   }{ \ang \xx^2 } - r \pa_r q^r
  = - \f{\ang z^2}{\ang \xx^2} (  q^r -q^z + \f{1}{\ang z^2} q^z )  - r \pa_r q^r .
\]

Applying \eqref{eq:dq_bound} and 
\eqref{eq:3D_size:C}, we obtain 
\[
  |r \pa_r q^r| \les r \CCa(\ww, \xx)
\les |\xx| ( \ang \xx^{\al-\hal-1} + \e^{\mhk}  \ang \xx^{-1} )
  \les \ang \xx^{\al - \hal } +  \e^{\mhk} .
\]
and using estimate \eqref{eq:damp_ag:est1}, we obtain
\beq\label{eq:damp_II:r}
  II_r 
  \leq - \f{\ang z^2}{\ang \xx^2} \B( \f{3}{2} \bar C_Q  \f{1}{\e} \ang \xx^{\al-\hau} 
  + \f{3}{2} \bar C_Q \ang z^{-2} \cJaa(\xx) \B) +  
  C (  \ang \xx^{\al - \hal } +  \e^{\mhk}  ).
\eeq

\vs{0.05in}
\paragraph{\bf Estimate of $II_z$}
For $II_z$, using $Q^z = z q^z$, we obtain 
\[
II_z = Q^z / z - \pa_z Q^z = - z \pa_z q^z.
\]

If $r \geq z$, using \eqref{eq:dq_bound}, 
and \eqref{eq:3D_size:C} for $\CCa$, and $\al - \hal + \epa < 0$,  we obtain
\beq\label{eq:damp_IIz_1}
  |II_z| \les |\xx| \CCa(\ww, \xx)
  \les |\xx|  ( \ang \xx^{\al- \hal  - 1}  + \e^{\mhk} \ang \xx^{-1} ) 
  \les  r  \ang \xx^{\al -\hal-1} + \e^{ \mhk } .
\eeq

If $r < z$, we compare $II_z(r, z)$ with $II_z(0, z)$. Using the definition of $q^z 
= \f{2 (\Psi(r,z) - \Psi(0,z)) }{z} + r \pa_r \f{1}{z} \Psi $ \eqref{def:qq}, 
and estimates \eqref{eq:tot_u_est1}, \eqref{eq:tot_u_estJ} for $\Psi$, 
 we obtain
\[
\bal
  | II_z(r,z) - II_z(0, z) |
  & \les | \f{2 (\Psi(r,z) - \Psi(0,z)) }{z} | +| r \pa_r \f{1}{z} \Psi| 
+ |\pa_z (\Psi(r, z) - \Psi(0, z ) )| + | r \pa_{r z} \Psi | \\
  &  \les \sup\nolim_{s \in [0, r]} | \tf{r}{z}\pa_r \Psi(s, z)|     
  + \sup\nolim_{ s\in[0,r] } |r \pa_{r z} \Psi(s, z)| 
  \les r \CCa(\ww, \xx) .
\eal
\]

Applying \eqref{eq:3D_size:C}  for $\CCa$ and $r \leq |\xx | $, we bound 
\beq\label{eq:damp_IIz_2}
    | II_z(r,z) - II_z(0, z) | 
    \les  r \ang \xx^{\al -\hal-1} + \e^{\mhk}.
\eeq

It remains to estimate $ II_z(0, z) = -z \pa_z q^z(0, z)$.
We decompose $Q^z(0, z) = \bar Q^z(0, z)
+ \td Q^z(0, z)$. Using $\bar Q^z = \vaa$ \eqref{eq:iden_r0:a} and \eqref{eq:1D_prof:va}, we obtain
\bseq\label{eq:damp_II_z:pf1}
\beq\label{eq:damp_II_z:pf1:a}
   - z \pa_z \bar q^z(0, z)
  =  - z \pa_z (\vaa / z) 
  = - \pa_z \vaa + \tf{\vaa}{z}
\leq - 4 \ang z^{ -\he }
+  C \ang z^{-\he} |\lgp z|^{-1/3} + C \e^{ \mhk }.
\eeq
Moreover, from \eqref{eq:dq_bound}, we obtain
\beq\label{eq:damp_II_z:pf1:b}
     |z \pa_z \bar q^z(0, z)| \les |z|^{\al+1}.
\eeq
\eseq

Combining the above two estimates, we obtain
\beq\label{eq:damp_IIz_3}
  - z \pa_z \bar q^z(0, z) \leq \f{|z|}{\ang z}( - 4 \ang z^{ -\he }
+  C \ang z^{-\he} |\lgp z|^{-1/3} ) + C \e^{ \mhk }.
\eeq
for some absolute constants. 
For $|z| \leq 1$, the estimate follows from 
\eqref{eq:damp_II_z:pf1:b} and choose $C$ large enough. For $|z| > 1$, the above estimate follows from 
\eqref{eq:damp_II_z:pf1:a} and $ -1 \leq -\f{|z|}{\ang x} $.

For $\td Q^z(0, z)$, we recall $\td Q^z(0,z) / z = \td c_l +  2 \f{1}{z}\psi(0,z) $ from \eqref{eq:psio_res1}. 
Using Proposition \ref{prop:vel_est}, we estimate 
\[
\bal
   | z \pa_z ( \td Q^z / z ) |
   \les | \psi_z - \tf{1}{z} \psi |
  \les \ang z^{\al-\hal } \rhoo(z)^{-1} \nnn{\ww} 
   \, .
  \eal 
\]

Using the bound of $\rhoo^{-1}$ \eqref{def:3d_wg}, $\al-\hal + \epa < 0$ by \eqref{ran:ep_all}, and $\nnn{\ww} \leq \e^{\hk}$  \eqref{eq:3D_size:ww}, we obtain
\beq\label{eq:damp_IIz_4}
     | z \pa_z ( \td Q^z / z ) | \les
           \ang z^{\al-\hal + \epa}  |\cJaa( z)|^{\kp}\cdot \e^{ \hk }
          \les  
           \e^{-\kp} \cdot \e^{\hk}
     \les \e^{\mhk}.     
\eeq

Combining \eqref{eq:damp_IIz_2}, \eqref{eq:damp_IIz_3}, and \eqref{eq:damp_IIz_4}, 
for $r < z$, we prove 
\[
\bal
  II_z(r, z)  & = (II_z(r, z) - II_z(0, z)  ) - z \pa_z \bar q^z(0, z)
  - z \pa_z ( \td Q^z(0,z) / z )  \\
 &  \leq  \f{|z|}{\ang z}( - 4 \ang z^{ -\he }
+  C \ang z^{-\he} |\lgp z|^{-1/3} ) + C \e^{ \mhk }
+ C r \ang \xx^{\al-\hal-1}.
\eal
\]

Since $ \al < \alb$, combining the above estimate for $r<z$ and \eqref{eq:damp_IIz_1} for 
$r \geq z$,  we prove 
\beq\label{eq:damp_II:z}
  II_z =  - \pa_z q^z \leq   - 4  \one_{r< z} |z|  \ang z^{-\he-1}  + \bar C_{II} |z|  \ang z^{-\he-1} |\lgp z |^{-1/3} 
  + C \e^{\mhk} + C r \ang \xx^{\al-\hal -1} , 
\eeq
for some absolute constant $ \bar C_{II}$.

\vs{0.05in}

\paragraph{\bf Summary of the estimates}
Now, we summarize the estimates of $I_R$ \eqref{eq:3D_damp_IR}, $I_{\mw{ag}}$ \eqref{eq:3D_damp:ag}, $II_r$ \eqref{eq:damp_II:r}, $II_z$ from \eqref{eq:damp_II:z}. Recall the estimate of $\bar \cR$ in \eqref{eq:3D_err:linf}. 
\[
  \bar \cR   
 \leq \bar C_{\cR} \min\B( |\xx|^{\al+1}, 
 \ang z^{\al-\hal}  \one_{r > z}
 + r \ang \xx^{\al-\hal - 1}
+  \e^{1 -\kp} + \f{|z|}{\ang z} \cJaa(\xx) |\lgp z|^{ \kp -2 } \B)  \, . 
\]

Using \eqref{eq:interp:pow2} in Lemma \ref{lem:interp1}, we bound 
\beq\label{eq:3D_interp:pow}
    \f{r}{\ang \xx} \ang \xx^{\al -\hal }  \les   \f{r}{\ang \xx} \ang \xx^{\al -\hau}  + 
    \e^2 .
\eeq

By requring $\e$ small,  we combine the lower order terms  $r \ang \xx^{\al-\hal - 1},
  \e^{1-\kp}, \e^{\mhk} $ in $\bar \cR$ and $I_{\mw{ag}}$ \eqref{eq:3D_damp:ag}  to obtain
\[
    I_{\mw{ag}} +  \bar \cR 
    \leq  - \kag \f{r}{\ang \xx} \bar C_Q \big( \f{7}{8}  \cdot \f{1}{\e} \ang \xx^{\al-\hau}  + \ang z^{-2} \cJaa(\xx)  \big) + C \e^{\mhk}  + \Rrem  \, , 
   \]
  where $\Rrem$ denotes the remaining terms in $\bar \cR_{\mw{rem}}$ 
\beq\label{eq:3D_damp:rem}
   \Rrem: = \bar C_{\cR}  \min \B( |\xx|^{\al+1}, 
\one_{r > z}\ang z^{\al-\hal}  + \f{|z|}{ \ang z } \cJaa(\xx) |\lgp z|^{ \kp - 2 } \B)  .
\eeq

Applying a similar estimate to $C r \ang \xx^{\al-\hal -1}$ in $II_z$ \eqref{eq:damp_II:z}
and adding $I_R$ \eqref{eq:3D_damp_IR}, $I_c$ \eqref{eq:3D_damp_Ic}
we yield
\bseq\label{eq:damp_sum1}
 \beq\label{eq:damp_sum1:z}
 \bal 
I_R + II_z + I_{\mw{ag}} & + I_c+ \bar \cR 
 \leq 
- \kag \f{r}{\ang \xx}   \bar C_Q  \B( \f{6}{8} \cdot \f{1}{\e} \ang \xx^{\al-\hau}  +  \ang z^{-2} \cJaa(\xx)  \B) + I_R   \\
& \quad 
 - 4  \one_{r< z} |z|  \ang z^{-\he-1}  + \bar C_{II} |z|  \ang z^{-\he-1} |\lgp z |^{-1/3} 
+ \Rrem 
  + C \e^{ \mhk } .
 \eal
\eeq

Recall $\ang y = \sqrt{y^2 + 1}$. Using $\kag < \f{1}{100}$ \eqref{def:rag}, $\f{r}{\ang \xx} + \f{\ang z^2}{\ang \xx^2} \geq 1$, we combine $II_r$ \eqref{eq:damp_II:r},
 $I_c$ \eqref{eq:3D_damp_Ic}, and $ I_{\mw{ag}} +  \bar \cR$ to obtain
\beq\label{eq:damp_sum1:r}
  I_R + I_{\mw{ag}} + II_r + I_c+ \bar \cR 
\leq 
- \kag  \bar C_Q  \B( \f{7}{8}\cdot \f{1}{\e} \ang \xx^{\al-\hau}  +  \ang z^{-2} \cJaa(\xx) \B)
+ C \e^{ \mhk } + I_R +   \Rrem  \, .
\eeq

\eseq

Recall $I_R$ from \eqref{eq:3D_damp_IR}:
\[
    I_R \leq  - \bar C_R \cJaa(\xx) \cdot ( (\bbb-1) \ang \xx^{1-\bbb} + \e ).
\]

We have not used $I_R$ in the above derivation.
Next, we further extract a uniform $O(1)$ damping terms. Since $\al - \hau , - \he \geq - 2 \e$ 
by \eqref{ran:ep2}, using \eqref{eq:interp:J},  we obtain
\[
\bal
   I_R - \kag \ang \xx^{\al-\hau} 
   \leq  - C ( \e \cJaa + \ang \xx^{\al-\hau} ) \leq  - C , 
  \eal
\]
for some absolute constant,
which along with \eqref{eq:damp_sum1:r} implies
\bseq\label{eq:damp_sum2}
\beq\label{eq:damp_sum2:r}
       I_R + II_r + I_{\mw{ag}} +I_c+  \bar \cR
       \leq 
       - \kag  \bar C_Q  \B(   \f{5}{8 \e} \ang \xx^{\al-\hau}  + \ang z^{-2} \cJaa(\xx) \B)
       - \bar C_1+ C \e^{\mhk} + \Rrem ,
\eeq
for some absolute constant $\bar C_1$.

For $r > z$ and $|\xx|\geq 1$, we obtain $\f{r}{\ang \xx} \ang \xx^{\al-\hau} \gtr \one_{r \geq z}  
\ang \xx^{\al-\hau} $; for $r < z$, and $|x| \geq 1$, we obtain 
\[
\one_{r < z}|z| \ang z^{ -\he -1} 
\gtr \one_{r < z} |\xx| \ang \xx^{ -\he -1}
\gtr \one_{r < z} \ang \xx^{ -\he }.
\]
Using \eqref{eq:interp:J}, 
and $-\he , \al - \hau \geq -2 \e$,  we obtain
\[
\bal
  &   -\one_{r < z}|z| \ang z^{ -\he -1} 
- \f{1}{8} \kag \f{r}{ \ang \xx } \ang \xx^{\al -\hau}
-  \bar C_R \cJaa ( (\bbb-1) \ang \xx^{\bbb-1} + \e ) \\
& \leq - C \B( \ang \xx^{\bbb-1} + 
\one_{|\xx| \geq 1} ( \cJaa \e  +  \one_{r < z} \ang \xx^{ -\he }
+ \one_{r > z} \ang \xx^{\al-\hau} ) \B)
\leq - C ( \ang \xx^{\bbb-1} +  \one_{|\xx| \geq 1}   )
\leq  - \bar C_2 \, ,
\eal 
\]
for some absolute constant $\bar C_2 >0$. 
Combining the above estimate and \eqref{eq:damp_sum1:z}, we prove
\beq\label{eq:damp_sum2:z}
\bal
       I_R + II_z + I_{\mw{ag}} + I_c+ \bar \cR
& \leq 
- \kag \f{r}{\ang \xx} \bar C_Q \B(   \f{5}{ 8 \e} \ang \xx^{\al-\hau}  + \ang z^{-2} \cJaa(\xx) 
    \B) 
     - \bar C_2 \\ 
    & \quad  
     + \bar C_{II} |z|  \ang z^{-\he-1} |\lgp z |^{-1/3} 
+ \Rrem + C \e^{ \mhk }   .
\eal
\eeq
\eseq

It remains to bound the positive terms in \eqref{eq:damp_sum2}, especially $\Rrem$ which can not be bounded by the good damping terms directly, when $r \gg |z| $  and $|z| \les 1$. 

\subsubsection{\bf From integrability to an $O(1)$ $z$-weight}

Lastly, we design a special $z$-weight $\rkeya$ to handle the 
remaining term $ \bar \cR_{\mw{rem}}$. Our key observation is that  $ \bar \cR_{\mw{rem}}$ 
\eqref{eq:3D_damp:rem} is integrable along the trajectory with $O(1)$ integral, and it can be controlled using a $O(1)$ $z$-weight.
Recall the parameter $\kp_1 \in (0, 1)$ from \eqref{def:kp}.  For $z \geq 0$, we construct $\rkeya$ as follows 
\bseq\label{def:rkey}
\beq\label{def:rkey:a}
  \rkeya(0) = 1,
  \quad  \f{ z \pa_z \rkeya }{\rkeya} =  - \f{z}{ \ang z  } | \lgp z|^{ -1 - \kp_1 } 
  - \e \f{z}{\ang z}  \ang z^{\al-\hal}. 
\eeq

Solving the ODE inequality, we obtain 
\beq
  \rkeya(z) = \exp\big( - \int_0^z \ang z^{-1} |\lgp z|^{ -1 - \kp_1  } - \e \ang z^{\al-\hal - 1} d z \big) \, .
\eeq
Since $ \al - \hal \leq - \f{8.9 \e}{8} < 0$ by \eqref{ran:ep2} and $\kp_1 \in (0, 1)$ \eqref{def:kp}, 
 $\rkeya(z)$ is decreasing in $z \geq 0$, we obtain that  $\ang z^{-1} |\lgp z|^{ -1- \kp_1  } \in L^1$
and $ \|  \e \ang z^{\al-\hal - 1} \|_{L^1} \les 1$. Hence, we obtain
\beq\label{def:rkey:c}
1\geq     \rkeya(z) \geq \rkeya(\infty) \gtr 1,
\eeq
with constant \emph{independent} of $\e$. 
For $\NNN \geq 1$ to be chosen and \emph{independent} of $\e$,  we choose 
\beq
\rkey = \rkeya^{\NNN} \, ,
\eeq
\eseq
Since $ \f{z \rkeya}{\rkeya} \leq 0$, and $ | \f{\ang z \rkeya}{\rkeya}| \les 1 $, using \eqref{eq:3D_boot_q:1}, for any $\NNN \geq 0$,  we obtain
\begin{align}\label{eq:damp_Iz0}
   I_z = \f{ Q^z \pa_z \rkey }{\rkey} 
   &= \f{ Q^z \pa_z \rkeya }{\rkeya} 
    \leq  \NNN \bar C_Q \cJaa( \xx )\cdot \f{z \pa_z \rkeya }{\rkeya} 
   + C  \NNN r \ang \xx^{\al -\hal -1}  \notag \\
  & = - \NNN \bar C_Q \cJaa( \xx )
   \cdot | \lgp z|^{-\kp_1-1}    \f{|z|}{\ang z}  
  - \NNN \bar C_Q \cJaa(\xx) \e  \ang z^{\al-\hal}  
  \f{|z|}{\ang z} 
     + C  \NNN r \ang \xx^{\al -\hal -1}  \notag \\
    & \teq \NNN  I_{z, 1} + \NNN I_{z, 2} +\NNN I_{z, 3} \, .
\end{align}

The last term $\NNN I_{z, 3}$ is treated as an error. Note that $I_{z, 1}, I_{z, 2} < 0$. We 
have a large factor $\cJaa(\xx)$ (rather than $\cJaa(z)$) in $I_{z, 1}, I_{z, 2} < 0$. 
Below, using the growth of $\cJaa(\xx)$ for large $|\xx|$ and choosing $\NNN$ large enough (independent of $\e$), 
we use $I_z$ to control the positive terms  in \eqref{eq:damp_sum2}.

\vspace{0.05in}
\paragraph{\bf Damping term $I_{z, 2}$}
Denote by $P_1$ part of the damping term in \eqref{eq:damp_sum2} 
\beq\label{eq:rkey_P1}
 P_1  \teq    \f{1}{8} \kag \f{r }{\ang \xx } \bar C_Q \B(   \f{1}{\e} \ang \xx^{\al-\hau}  + \ang z^{-2} \cJaa(\xx) \B) .
\eeq

Using Lemma \ref{lem:interp2} with $k =\hal, \ell = \e^{2/3}$,  and $\al-\hal, \al -\hau 
\in [- \f{9.1 \e}{8}, -\f{8.9\e}{8}]$ by \eqref{ran:ep2}, we combine $P_1$ and the damping term $I_{z, 2}$ in \eqref{eq:damp_Iz0} to obtain
\begin{align}\label{eq:rkey_P2}
 & P_1  +\tf12 \NNN | I_{z, 2} | 
    + \e^{1/3}  \notag \\
 &   \geq C(\kag, \bar C_Q) \min( \e^{-1/3}, \NNN ) 
     \cdot \B( \f{r }{\ang \xx  } (  \e^{-2/3}  \ang \xx^{\al-\hau}
  + \e \cJaa(\xx) \cdot \ang z^{-2}) +
  \e \cJaa(\xx) \ang z^{\al-\hal} \f{|z|}{\ang z} + \e^{2/3} \B) \notag \\
  & \geq C(\kag, \bar C_Q) \min( \e^{-1/3}, \NNN ) \cdot \min( |\xx|^{\al+1} , \one_{r > z} \ang z^{\al- \hal} ) .
\end{align}

The remarkable feature of the above estimate is that 
when $r \gg z$, $P_1$ is very small relative to $\ang z^{\al-\hal}$, but we can control the term 
$\ang z^{\al- \hal}$ using the new damping term $I_{z, 2}$.
  
\vspace{0.05in}
\paragraph{\bf Damping term $I_{z, 1}$}
Since $ \kp_1 < \f{1}{3}$ by \eqref{def:kp},  $\he \gtr \e$ \eqref{ran:ep_all}, and $ \ang z^{-\he} \lgp z \les \cJaa(z) \les \cJaa(\xx)$ by Lemma \ref{lem:lgx_pow},  we obtain
\[
 \bar C_{II}  |z| \ang z^{-\he -1 } |\lgp z|^{-1/3}
 \leq C  |z| \ang z^{-\he -1 } |\lgp z|^{ -\kp_1 }
 \leq \bar C_3 \f{|z|}{\ang z} \cdot |\lgp z|^{ -\kp_1-1} \cdot  \cJaa(\xx)
 = \f{\bar C_3}{\bar C_Q} | I_{z, 1} | \, ,
\]
for some absolute constant $\bar C_3>0$.  For the $\cJaa$ term in $\Rrem$ \eqref{eq:3D_damp:rem}, 
since $ \kp_1 < 1 -\kp$  \eqref{def:Kin} implies 
$\kp-2 < - \kp_1-1$, we obtain
\[
   \bar C_{\cR}  
\cdot  \f{|z|}{ \ang z } \cJaa(\xx) |\lgp z|^{ \kp-2 } 
\leq \bar C_4 \cdot  \f{|z|}{ \ang z } \cJaa(\xx) |\lgp z|^{ -\kp_1-1 } 
\leq  \f{ \bar C_4}{\bar C_Q} |I_{z, 1}| \,  .
\]
for some absolute constant $\bar C_4$. Now, we choose $\NNN$ as:
\beq\label{def:NNN}
   \NNN = 2 \big( \f{ \bar C_R  }{  C(\kag, \bar C_Q) } +  \f{ \bar C_3 }{\bar C_Q }  + \f{\bar C_4}{\bar C_Q} + \f{1}{\bar C_Q} \big) .
\eeq

Choosing  $\e$ with $\e^{-1/3} \geq \NNN$, using 
 $I_{z, 1}, I_{z, 2} \leq 0$, \eqref{eq:rkey_P1} and the above estimates, we obtain
\[
\bal
 & P_1 - \tf12 \NNN I_{z, 1} -\tf12 \NNN I_{z, 2} + \e^{1/3}
   = ( P_1 +\tf12 \NNN  |I_{z, 2}| + \e^{1/3} ) + \tf12 \NNN  |I_{z, 1} |   \\
& \qquad \geq 
C(\kag, \bar C_Q) \NNN \cdot \min( |\xx|^{\al+1} , \one_{r > z} \ang z^{\al- \hal} ) 
+ \tf{\bar C_3}{\bar C_Q} | I_{z, 1} | 
+  \tf{ \bar C_4}{\bar C_Q} |I_{z, 1}| \\
& \qquad \geq 
 \bar C_{II}  |z| \ang z^{-\he -1 } |\lgp z|^{-1/3}
 + \bar C_R \min( |\xx|^{\al+1} , \one_{r > z} \ang z^{\al- \hal} )  +
    \bar C_{\cR}   
\cdot  \f{|z|}{ \ang z } \cJaa(\xx) |\lgp z|^{ \kp - 2 }  \\
& \qquad \geq  C_{II}  |z| \ang z^{-\he -1 } |\lgp z|^{-1/3} +    \Rrem .
\eal
\]

Using the definitions of $\NNN$ \eqref{def:NNN} and $I_{z,i}$ \eqref{eq:damp_Iz0} and $-I_{z, 2} \geq 0$, we obtain
\[
\bal
  - \tf12 \NNN I_{z, 1} -\tf12 \NNN I_{z, 2} 
 - \NNN I_{z, 3}
&  \geq \cJaa(\xx)  | \lgp z|^{ -\kp_1- 1} \tf{|z|}{\ang z}
   - C \NNN r \ang \xx^{\al-\hal - 1}. 
  \eal
 \]

Adding the above two estimates and using  $I_z = \NNN ( I_{z, 1} + I_{z, 2} + I_{z, 3} ) $, we prove 
\beq\label{eq:damp_Iz}
\bal
P_1 - I_z \geq  C_{II}  |z| \ang z^{-\he -1 } |\lgp z|^{-1/3}
  + \Rrem + \cJaa(\xx)  | \lgp z|^{ -\kp_1 - 1}  
  \tf{|z|}{\ang z}   - C \e^{1/3} - C \NNN r \ang \xx^{\al-\hal - 1}. 
  \eal
\eeq

\vspace{0.05in}
\paragraph{\bf Combined damping terms}

Recall $P_1$ from \eqref{eq:rkey_P1}.  Combining \eqref{eq:damp_Iz} $\times (-1) $ and \eqref{eq:damp_sum2}, 
and using $\e^{1/3} \leq  \e^{ \mhk }$, we prove 
\beq\label{eq:damp_sum3}
  \bal
       I_R + II_r + I_{\mw{ag}}  + I_c + \bar \cR  + I_z
 &      \leq - \f{1}{2} \kag \bar C_Q \B(   \f{1}{\e} \ang \xx^{\al-\hau}  + \ang z^{-2} \cJaa(\xx)  \B) - \bar C_1         \\
& \quad - \cJaa(\xx)  | \lgp z|^{ -\kp_1 - 1 }     \tf{|z|}{\ang z}  
      + C \e^{ \mhk } + C \NNN r \ang \xx^{\al-\hal - 1}, \\ 
    I_R + II_z + I_{\mw{ag}} + I_c+  \bar \cR
    + I_z
& \leq 
- \f{1}{2} \kag \f{r}{\ang \xx} \bar C_Q \B(  \f{1}{\e} \ang \xx^{\al-\hau}  + \ang z^{-2} \cJaa(\xx) 
    \B) 
     - \bar C_2 \\
& \quad - \cJaa(\xx)  | \lgp z|^{ -\kp_1 - 1 }    \tf{|z|}{\ang z}    + C \e^{ \mhk }  + C \NNN r \ang \xx^{\al-\hal - 1}.
\eal
\eeq
We further simplify the bound. For the $\cJaa(\xx)$-terms, 
since $|z| \ang z^{-1} + \ang z^{-2 } \gtr 1$,  we have
\bseq\label{eq:damp_sum4}
\beq\label{eq:damp_sum4:a}
- \f{r}{\ang \xx} \cJaa(\xx) \ang z^{-2} - \cJaa(\xx) \f{|z|}{\ang z} |\lgp z|^{ -\kp_1 - 1 } \leq - 
C \f{r}{\ang \xx}  \cJaa(\xx) |\lgp z|^{ -\kp_1 - 1 }.
\eeq

Since $\NNN$ \eqref{def:NNN} is independent of $\e$, applying \eqref{eq:3D_interp:pow}, we bound the positive terms in \eqref{eq:damp_sum3} as
\beq
 C \NNN r \ang \xx^{\al-\hal - 1}
\les r \ang \xx^{\al-\hau - 1} 
+ \e^2 .
\eeq

\vspace{0.05in}
\paragraph{\bf Estimate of $\Lpsi$}

Using the estimate of $\Lpsi$ from \eqref{eq:Lpsi_linf_est:a}, the assumption \eqref{eq:3D_size:ww}
 for $\nnn{\ww}$, and $\cJak \les \e^{-\kp}$ by \eqref{eq:Ja_hat}, we obtain 
\[
\bal
  |\Lpsi| 
& \les   \big( \one_{ r> z }  |\lgp z|^{ \kp - 2 } 
  \cdot |\cJaa(\xx) |   + 1 \big) \cdot \e^{-\kp} \cdot \e^{\hk}  \, .
\eal
\]
Since 
\[
  |\lgp z|^{ \kp - 2 } 
  \cdot |\cJaa(\xx) |  
  \les 1, \  \forall \ |\xx| \leq 1, \quad \one_{ r> z } \les \tf{r}{\ang \xx} , \quad \forall \ |\xx| > 1,
\]
and $\kp- 2 < -\kp_1 - 1$ by \eqref{def:kp}, we obtain
\beq
   |\Lpsi|  \les \big( \tf{r}{\ang \xx}   |\lgp z|^{ -\kp_1 - 1 } 
  \cdot |\cJaa(\xx) |   + 1  \big) \e^{\mhk} .
\eeq
\eseq
Thus, $\Lpsi$ can be absorbed by the damping term in \eqref{eq:damp_sum4:a} 
and the constant term in \eqref{eq:damp_sum3} for $\e$ small enough. Combining estimates \eqref{eq:damp_sum3}, \eqref{eq:damp_sum4}, and taking $\e$ small enough, we prove
\[
  \bal
       I_R + II_r + I_{\mw{ag}} + I_c + \bar \cR  + \Lpsi +  I_r
 &      \leq - \f{1}{3} \kag \bar C_Q \B(   \f{1}{\e} \ang \xx^{\al-\hau}  +
 \bar C_3 |\lgp z|^{-\kp_1-1} \cJaa(\xx)  \B) - \f{1}{2} \bar C_1       , \\
    I_R + II_z + I_{\mw{ag}} + I_c+ \bar \cR + \Lpsi +  I_z
& \leq 
- \f{1}{3} \kag \f{r}{\ang \xx} \bar C_Q \B(   \f{1}{\e} \ang \xx^{\al-\hau}  + 
 \bar C_3 |\lgp z|^{ -\kp_1 - 1 }\cJaa(\xx) 
    \B) 
     - \f{1}{2} \bar C_2 \, ,
\eal
\]
for some $\bar C_4 >0$.  Upon renaming the absolute constant $\bar C_i, \bar C_Q$, we prove Lemma \ref{lem:3D_damp}.

\begin{remark}[\bf Integrability]

The above estimates exploit the key property that
$\Rrem$ \eqref{eq:3D_damp:rem} is integrable along the trajectory, and it 
has a faster decay in $z$. Moreover, we exploit the stronger $z$-outgoing effect 
$Q^z /z \asymp \cJaa(\xx)$ \eqref{eq:3D_boot_q:1}  for large $|\xx|$ and $r \les |z|$.
\end{remark}

\begin{remark}[\bf Absolute constants]
Since  $\NNN$ \eqref{def:NNN}  depends on $\bar \cR$ uniformly in $\e$, 
we treat $\NNN$ and  $\kag$ \eqref{def:rag}  as absolute constants. 
Since $\rkeya$ satisfies \eqref{def:rkey:c} with absolute constants 
and $\rkey = \rkeya^{\NNN}$ \eqref{def:rkey}, in the following estimates, the absolute constants may depend on $\rkey, \NNN,\kag$.
\end{remark}

\subsection{Weighted Estimates for $\pa_z G$}\label{sec:EE_C1z}

In this section and in Section \ref{sec:EE_C1r}, we perform weighted $\cW^{1,\infty}$ estimates on 
 $\pa_z G$ and $\pa_r G$ \eqref{eq:lin_dG}.  We first estimate $\pa_z G$, which is simpler. 

 Recall $\rhor =\rhoo \rag \ang \xx ,\rhoz = \rhoo \rag |z|$ from \eqref{def:3d_wg}. We introduce the local energy in $\bar B_R$ 
 \footnote{
For $\ww$ satisfying \eqref{eq:3D_size:ww}, from Proposition \ref{prop:reg_solu} and \eqref{eq:local_reg:nnn}, we obtain 
$\eta \in C^{1,\al}(B_R), \llr{G}|, \llz{G}, \lla{G} < \infty, \forall R >0$. 
Thus, the weighted terms in \eqref{def:local_norm} are continuous, the maximum is well-defined, 
and the local norms \eqref{def:local_norm} are bounded for any $R > 0$. 
 }
 \beq\label{def:local_norm}
 \bal
     \lla{G} & \teq   \max_{\xx \in \bar B_R } | \rhoc  \rhoo \rag  G(\xx)| ,   \\
   \llw{G} & \teq \max_{ \xx \in \bar B_R } |  \rhoc \rhor \f{z}{\ang z} \pa_r G(\xx)| 
   = \max_{ \xx \in \bar B_R } |  \rhoc \rhoo \rag \f{z}{\ang z}  \ang \xx \pa_r G(\xx)|  ,
   \\ 
      \llz{G} & \teq \max_{\xx \in \bar B_R } | \rhoc \rhoz \pa_z G(\xx)| 
      \teq \max_{\xx \in \bar B_R } | \rhoc \rhoo \rag |z| \pa_z G(\xx)| .
  \eal
  \eeq
for any $R > 1$. Note that we use $|z|$ rather than the Japanese bracket $\ang z$ in \eqref{def:local_norm}.

The above energy is closely related to the norms in \eqref{def:3d_norm} up to a factor $\wwb$. 
Using 
$  \f{z}{\wwb} \pa_z \eta = z \pa_z \f{\eta}{\wwb} + \f{z \pa_z \wwb}{\wwb} \eta $, $ |\f{z \pa_z \wwb}{\wwb} | \les 1 $ by \eqref{eq:1D_prof_est2:c}, 
and 
$|\wwb|^{-1} \gtr |z|^{-1} + \ang z^{\hal}$, $ |\f{z}{\wwb}  | \gtr \ang z^{\hal+1}$ by \eqref{eq:wa_upper_lower}, we have 
\bseq\label{eq:norm_comp2}
\begin{align}
   \nnrb{\eta} & \les E_{r, \infty}( \f{\eta}{\wwb} ), \quad  \qquad \qquad \nnlb{  \eta} \les E_{0, \infty}( \f{\eta}{\wwb} )  , \label{eq:norm_comp2:a} \\
  \nnzb{\eta}   
  & =  
  \nlinf{ \rhoc \rhoo \rag \ang z^{\hal + 1} \pa_z \eta}
  \les   \nlinf{ \rhoc \rhoo \rag \f{z}{\wwb} \pa_z \eta}
  \les    E_{z, \infty}( \f{\eta}{\wwb} ) + E_{ 0 ,\infty}( \f{ \eta}{\wwb} ). \label{eq:norm_comp2:b}
\end{align}
\eseq

The above norms $\nnlb{\cdot}, \nnrb{\cdot}, \nnzb{\cdot}$ relate to the weight $\rhoo \rhoc$, which is stronger than $\rhoo$ appearing in the norms $\nna{\cdot}, \nnr{\cdot}, \nnz{\cdot}$ in \eqref{def:3d_norm}. Since the power $\e^2$ in $\rhoc$ \eqref{def:rhoc}  is very small compared to $\e$, the reader can 
essentially treat $\rhoc$ as $1$ in the following estimates.

\begin{remark}[\bf Local estimate of $G$]
We do not bound $G$ using the global weighted $L^{\infty}$ norm 
since the boundedness of these norms does not follow from Proposition \ref{prop:reg_solu}. 
Instead, we estimate $G$ using \eqref{eq:lin_dG} in $\bar B_R$ (see also \eqref{eq:lin_Cz}, \eqref{eq:lin_Cr}) and obtain uniform bound in $R$. 
\end{remark}

Below, unless indicated specifically, all the implicit constants \textit{do not} depend on $\e$.

\subsubsection{Integrating factor $\vpz$}\label{sec:vpz}
We design $\vpz > 0$ to capture the local terms 
$ (\bar \cR - \pa_z Q^z ) \cdot \pa_z G$ in \eqref{eq:lin_dG:Gr},
which plays a role as an integrating factor.
We design $\vpz$ with $ \lim_{(r,z) \to 0} \f{\vpz}{z} = 1$ and
\beq\label{eq:lin_vpz}
  \QQ \cdot \na \vpz = -  ( \bar \cR + \Lpsi- \pa_z Q^z) \vpz , \iff 
  \QQ \cdot \na ( \f{\vpz}{z} ) = - ( \bar \cR + \Lpsi - \pa_z Q^z + \f{1}{z} Q^z )  \f{\vpz}{z}  .
\eeq

Using the characteristics $X(x, s)$ defined in \eqref{eq:traj}, we rewrite the equation as 
\[
   \f{d}{d s} ( \f{\vpz}{z} ) \circ  X(s, \xx)  = - \B( ( \bar \cR + \Lpsi- \pa_z  Q^z + \f{1}{z} Q^z ) \cdot ( \f{\vpz}{z}) \B) \cc  X(s, \xx) , \quad  \lim_{(r,z) \to 0} \f{\vpz}{z} = 1.
\]
We solve the ODE as follows 
\[
\bal
 &  \f{\vpz}{z} \circ X(s, \xx)   =  \exp\B( - \int_{-\infty}^s  ( \bar \cR  + \Lpsi- \pa_z Q^z + \f{Q^z}{z}  ) \cc X(\tau,  \xx) d \tau   \B) , \quad  \lim_{(r,z) \to 0} \f{\vpz}{z} = 1.
  \eal 
\]

From estimate \eqref{eq:Q_bound:b} and \eqref{eq:3D_err:abs},
we obtain $ |\bar \cR + \Lpsi- \pa_z Q^z + \f{Q^z}{z} | \les_{\e}  |\xx|^{1+\al} (1 + \nnn{\ww})$, which vanishes near $\xx=0$.  
 Since $|X(\xx, \tau)| \les e^{C \tau} |\xx|, \tau \leq 0$ by \eqref{eq:traj_pf1}, the above integrand is integrable. Since $X( \xx, 0 ) = \xx$, taking $s=0$, we obtain 
\beq\label{eq:lin_vpz2}
   \vpz( \xx )  = z \exp\B( - \int_{-\infty}^0  ( \bar \cR + \Lpsi - \pa_z Q^z + \f{Q^z}{z}  ) \cc X(\tau, \xx ) d \tau   \B) , \quad  \lim_{(r,z) \to 0} \f{\vpz}{z} = 1 , 
   \quad  |\f{\vpz(\xx)}{z}| \les C(\e, |\xx|) < \infty.
\eeq

Using $\vpz$ and \eqref{eq:lin_dG}, for any $z \neq 0$, we yield 
\bseq\label{eq:lin_Cz}
\beq
 \tf{d}{ds} ( \vpz \pa_z G) \cc X(s, \xx) = (\vpz \cmz ) \cc  X(s, \xx).
\eeq

Solving $\vpz \pa_z G$ along the flow $ X$, we obtain 
\beq
  (\vpz \cdot \pa_z G)(\xx) = \int_{-\infty}^0 \vpz \cmz \cc  X(s, \xx) d s.
\eeq
\eseq
Since $\vpz(r, z) \sim z$ near $0$ by \eqref{eq:lin_vpz2},
and $X(s, \xx ) \to 0$ as $s \to - \infty$, 
using $|\pa_z G| \les_{R, \e} 1$ by \eqref{eq:local_reg}, we obtain that  the boundary term at $s = -\infty$ vanishes:  $\lim_{s \to -\infty} \vpz \pa_z G \cc  X(s, \xx) = 0$.

By definition of $A^z$ \eqref{eq:damp:Az} and $\vpz$, we obtain 
\beq\label{eq:damp:Az2}
  A^z =    \f{\QQ \cdot \na ( \rhoc \rhoz \rkey )}{\rhoc \rhoz \rkey }  + \bar \cR + \Lpsi - \pa_z Q^z
  =  \f{\QQ \cdot \na ( \rhoc \rhoz \rkey (\vpz)^{-1} )}{ \rhoc \rhoz \rkey  (\vpz)^{-1}} .
\eeq

We aim to apply Lemma \ref{lem:traj_est} with weights $(\vp, m) \rsa ( \vpz, \rhoc \rhoz \rkey )$.  Recall  $\cmz$ from \eqref{eq:lin_dG}. We further decompose it as 
\beq\label{eq:decom_Cz}
\bal
    \cmz& = - \pa_z Q^r \cdot \pa_r G  + \pa_z \bar \cR \cdot (G+1)  +  \pa_z \cL_{\psi} \cdot (G+1)  .
  \eal
\eeq

From Proposition \ref{prop:reg_solu},  estimate \eqref{eq:reg:est_coe} for 
$ \wwb( \bar \cR + \Lpsi)$, and  $\f{1}{z} \wwb = \f{1}{z} \waa \in C^{2}$ locally by 
\eqref{eq:waa_reg}, we obtain that $z \cmz$ is locally $C^{\al}$.  We have the following estimates for $\eta$.

\begin{prop}\label{prop:Cz}
Suppose that $\ww$ satisfies \eqref{eq:3D_size:ww}. 
Let $\eta$ be the solution to 
\eqref{eq:fix_pt_2D}, $G =\f{\eta}{\wwb}$, and $\rhoo, \rag, \rhoc, \vpz,  \cmz$  be defined in \eqref{def:3d_wg}, \eqref{eq:lin_vpz}, \eqref{eq:lin_Cz}. 
For any $R > 1$ and  any $\xx \in \bar B_R$, we have the following pointwise estimates on $\cmz$ 
and $\pa_z G(x)$
\beq\label{eq:cmz_est}
   |\rhoc \rhoo \rag z \pa_z G(\xx)| 
    \les 
   \e^{ -\kp_1 } \lla{G}     + \e \llw{G}+  \e + \nnn{\ww} ,
\eeq
with implicit constants \textit{independent} of $R$.

\end{prop}

Each term in $\cmz$ \eqref{eq:decom_Cz} is estimated perturbatively based on
Propositions \ref{prop:Lpsi}-\ref{prop:3D_error}. 
In the following estimates, we fix $R > 1$ and assume $\xx \in \bar B_R, z \neq 0$. From Item (a) in Proposition \ref{prop:reg_solu}, we obtain
\beq\label{eq:support}
\xx=(r,z) \in \bar B_R, \  z \neq 0  \quad \Rightarrow \quad  X(s, \xx) \in \bar B_R, \quad 
X^z(s, \xx) \neq 0,  \forall \    s \leq 0 .
\eeq

In the following estimates, $G$ only evaluates in $B_R$, and we only need the local norms 
in \eqref{def:local_norm}.

\subsubsection{Estimate of $\cmz$}

Recall the weight $\rhor = \rhoo \rag \cdot \ang \xx, \rhoz = \rhoo \rag |z|$ from \eqref{def:3d_wg}.

Under the assumption \eqref{eq:3D_size:ww}, using \eqref{eq:Q_mix_est}
and the energies in \eqref{def:local_norm}, we obtain
\beq\label{eq:cmz_est:Gr}
\bal
  \rhoc \rhoo \rag |z|\cdot | \pa_z Q^r \cdot \pa_r G | 
  & \les \rhoc \rhoo \rag |z| \cdot r \ang z^{-1} \ang \xx \CCa(\ww, \xx) \cdot | \pa_r G| \\
  & \les  r \CCa(\ww, \xx) \cdot |z| \ang z^{-1} |\rhoc \rhoo \rag \ang \xx \pa_r G|  \\
    & \les r \CCa(\ww, \xx) \llw{G}.
\eal
\eeq

\paragraph{\bf Estimate of $\bar \cR$-terms }
For $\pa_z \bar \cR$, using \eqref{eq:3D_err:Cz}, we obtain
\[
      | z \pa_z \bar \cR |  \les 
      \f{r}{\ang \xx} \min\B( |\xx|^{\al}, 
\ang z^{\al - \hal } + 
 \min\B(  \ang z^{\al- \hal}   \log \f{\ang \xx}{\ang z} , \cJaa(\xx) \B) \cdot \cJaa(z)^{-1}  \B) .
\]

From \eqref{ran:ep_all} and \eqref{def:rag}, we obtain $|\al-\hal| \leq 2 \e  \ll \kag$ and $\kag \gtr 1$.  
Using $ \log t \les_a t^{a}$ for any $a \geq 0 , t \geq 1$, 
and  \eqref{eq:err_ag_decay}, we obtain
\[
  \rag \ang z^{\al- \hal}  
\les \ang \xx^{\al-\hal}, 
\quad  \rag  \ang z^{\al- \hal}   \log \tf{\ang \xx}{\ang z} \les  \ang \xx^{\al-\hal}  .
\]

Using the above estimates and \eqref{eq:ep_ineq1} on $\rhoo \rhoc$, we obtain
\[
  | \rhoc \rhoo \rag |z| \pa_z \bar \cR| 
  \les \f{r}{\ang \xx} \rhoc  \rhoo  \min( |\xx|^{\al}, \ang \xx^{\al-\hal} ) 
\les r \ang \xx^{\al-\hal - 1}.
\]

Next, we estimate the unweighted $z \pa_z \bar \cR$. 
Using $\cJaa(z) \les \e^{-1}$ \eqref{eq:Ja_hat}, 
and $|z| \leq |\xx|$, we obtain
\bseq\label{eq:interp:JxJz2}
\beq
\bal
P_J \teq 1 + \cJaa(\xx) \cJaa(z)^{-1} & 
  \les 
    \e^{ - \kp_1} \cJaa(\xx) \cJaa(z)^{-\kp_1 - 1} .
    \eal
\eeq

Moreover, since  $ \cJaa(\xx) \cJaa(z)^{-1} \les 1$  for $|z| > r$ or $|\xx| \les 1$ 
and since  $\cJaa(z)^{-1} \les |\lgp z|^{-1} + \e,\, \e \cJaa(\xx)\les 1$ by \eqref{eq:Ja_hat}, we yield 
\beq
P_J 
    \les 1 +  \e^{ - \kp_1}   \f{r}{\ang \xx} \cJaa(\xx) \cJaa(z)^{-\kp_1 - 1} 
   \les 1 +  \e^{-\kp_1}   \f{r}{\ang \xx}  \cJaa(\xx) |\lgp z|^{-\kp_1 -1} .
\eeq
\eseq

Using $\ang z^{\al-\hal} \les 1$, \eqref{eq:interp:JxJz2},  we obtain
\[
\bal
  |z \pa_z \bar \cR| & \les P_J \les  1 +  \e^{-\kp_1}   \f{r}{\ang \xx}  \cJaa(\xx) |\lgp z|^{-\kp_1 -1} .
\eal
\]

 The coefficient in the above upper bound relates to the damping terms $A^z$ in \eqref{eq:damp:Az}. 
 Compared to $A^z$, we lose a large factor $\e^{-\kp_1}$ but with a very small exponent. 

Combining the above estimate on $z \pa_z \bar \cR$ and using the norm \eqref{def:local_norm}, we establish 
\beq\label{eq:cmz_est:cR}
\bal
   \rhoc \rhoo \rag |z|\cdot | (\pa_z \bar \cR) \cdot (G+1)| 
  & \les |\rhoc \rhoo \rag z \pa_z \bar \cR|
   +  | z \pa_z \bar \cR | \cdot  | \rhoc \rhoo \rag G | \\
  & \les  r \ang \xx^{\al-\hal - 1} +   \big(1 +  \e^{-\kp_1} \f{r}{\ang \xx} \cdot \f{ \cJaa(\xx)}{  |\lgp z|^{\kp_1+1} } \big)  \lla{G} .
  \eal
\eeq

\paragraph{\bf Estimate of $\pa_z \Lpsi$-term}

For the $\pa_z \Lpsi$ term in \eqref{eq:decom_Cz}, using \eqref{eq:Lpsi_Cz_est}
and the definition of $\lla{\cdot}$ norm \eqref{def:local_norm}, 
we bound 
\[
\bal
  \rhoc \rhoo \rag |z| \cdot |\pa_z \Lpsi \cdot (G+1)|
& \les \rhoc \rhoo \rag |z \pa_z \Lpsi |
+ |z \pa_z \Lpsi| \cdot |\rhoc \rhoo \rag G|   \\
& \les \rhoc \rhoo \rag |z \pa_z \Lpsi |
+ |z \pa_z \Lpsi| \cdot \lla{G} \\
& \les \rhoc \min(|\xx|^{1+\al},  \ang \xx^{- \epa} ) \nnn{\ww}  \\
& \quad + ( \ang \xx^{\al- \hal + \epa} + \cJaa(\xx) \cJaa(z)^{-1}  ) \cJak \nnn{\ww} \lla{G}.
\eal
\]

Since $\rhoc = \ang \xx^{\e^2}$ \eqref{def:rhoc},  using  \eqref{eq:interp:JxJz2},  $\al -\hal + \epa < 0, \e^2 < \epa , \epa >0$ by \eqref{ran:ep_all}, 
and $\cJaa \les \e^{-1}$, we prove
\beq\label{eq:cmz_est:Lpsi_z}
  \rhoc  \rhoo \rag |z| \cdot |\pa_z \Lpsi  \cdot (G+1)|
  \les  \nnn{\ww} + \B( 1 + \e^{-\kp_1}   \f{r}{\ang \xx} \cdot \f{\cJaa(x)}{ | \lgp z |^{\kp_1+1} }  \B)
  \e^{-\kp} \nnn{\ww} \lla{G}.
\eeq

\subsubsection{Summary of estimates}

Recall the term $\cmz$ from \eqref{eq:decom_Cz}. Combining estimates \eqref{eq:cmz_est:Gr}-\eqref{eq:cmz_est:Lpsi_z} and similar terms, 
and then using $ \rkey \asymp 1 $ \eqref{def:rkey}, $ \e \leq 1$, for $\cmz$ and $z \neq 0$, we prove 
  \bseq\label{eq:cmz_est:tot}
  \beq
\bal
  \rkey \rhoc \rhoo \rag |z|  & \cdot |\cmz|  
 \les 
\f{r}{\ang \xx} \cdot \f{ \cJaa(\xx)}{  |\lgp z|^{\kp_1+1} } \B( \udb{ \e^{-\kp_1} \lla{G} }_{ \eqref{eq:cmz_est:cR} } +
\udb{ \e^{- \kp - \kp_1} \nnn{\ww} \lla{G} }_{ \eqref{eq:cmz_est:Lpsi_z} }   \B) \\
&  + \udb{ r \CCa(\ww, \xx) \llw{G} }_{\eqref{eq:cmz_est:Gr}  }  
+ \udb{  \lla{G} +  r \ang \xx^{\al-\hal-1} }_{  \eqref{eq:cmz_est:cR} }
 +  \nnn{\ww} ( \udb{ 1 + \e^{-\kp}   \lla{G} }_{ \eqref{eq:cmz_est:Lpsi_z} } ).
\eal
\eeq

Using $\kp + \kp_1 < \f{1+ \kp}{2}$ by \eqref{def:kp} and  $ \e^{-\kp-\kp_1}   \nnn{\ww} \les 1$ 
by \eqref{eq:3D_size:ww}, we simplify the bound as 
\beq
\bal
  \rkey \rhoc \rhoo \rag |z| \cdot |\cmz|  
  & \les  
  \f{r}{\ang \xx}  \cJaa |\lgp z|^{- \kp_1-1} \cdot \e^{- \kp_1}  \lla{G}   \\
 & \quad  +   r \CCa(\ww, \xx) \llw{G}
 + r \ang \xx^{\al-\hal-1} 
+   \lla{G} 
+ \nnn{\ww} .
 \eal
\eeq
\eseq

To apply Lemma \ref{lem:traj_est} with $(m, \vp) \rsa (\rhoc \rhoz \rkey, \vpz)$, we further compare the bound  \eqref{eq:cmz_est:tot} and the damping terms $A^z$ \eqref{eq:damp:Az}. Recall $A^z$ from 
\eqref{eq:damp:Az}, \eqref{eq:damp:Az2}
\[
  A^z  = \f{\QQ \cdot \na ( \rhoc \rhoz \rkey (\vpz)^{-1} )}{ \rhoc \rhoz \rkey (\vpz)^{-1} }  
 \leq 
 -\bar \lam_1 \B(  \f{r}{\ang \xx} \big(  \f{1}{\e} \ang \xx^{\al-\hau} 
 +  \cJaa(\xx)   |\lgp z|^{-\kp_1 - 1} \big)    +1 \B) < 0.
\]

Using \eqref{eq:interp:pow2} , we obtain
\beq\label{eq:damp:Az_comp}
\bal
\ang \xx^{ \al - \hal } \les \ang \xx^{ \al - \hau } + \e^2,  \ 
\Longrightarrow  \
  r \ang \xx^{ \al - \hal-1 }  \les   r  \ang \xx^{\al -\hau - 1} + \e^2  \les  \e |A^z| . \\
  \eal
\eeq

Using \eqref{eq:3D_size:C} and \eqref{eq:interp:pow} with $a = \e$, $ \kp_1 <  \f{1 - \kp}{2}$ by \eqref{def:kp}, and the above estimate, we bound 
\beq\label{eq:CC_bound2}
\bal
   \CCa(\ww, \xx) & \les \ang \xx^{\al-\hal - 1} + \ang \xx^{\al-\hal + \epa -1} \e^{\mhk}
   \les  \ang \xx^{\al-\hal - 1} 
   + \ang \xx^{-1} ( \e^{\mhk -\kp_1} \ang \xx^{\al-\hau} + \e   ) \\
  & \les  \ang \xx^{\al-\hal - 1}  + \e \ang \xx^{-1}
  \les \ang \xx^{\al-\hau -1} + \e \ang \xx^{-1} . 
  \eal
\eeq
It follows 
\[
   r  \CCa(\ww, \xx) \les \e |A^z|.
\]

Using \eqref{eq:cmz_est:tot}, the above bound on $A^z$, 
and $\e \leq 1$,  we have 
\beq\label{eq:IBP_cmz:sum1}
\bal
    \rhoc \rhoz \rkey |\cmz| 
    & \les |A^z| \cdot
    (    \e^{ -\kp_1 } \lla{G}     + \e \llw{G}+  \e + \nnn{\ww} ).
  \eal
\eeq

Recall from \eqref{eq:lin_Cz}
\[
      \vpz \cdot \pa_z G =
\tts{\int_{-\infty}^0}( \vpr  \cmr ) \cc  X(s, \xx) d s .
\]

Applying Lemma \ref{lem:traj_est} with $(m, \vp, H) \rsa (\rhoc \rhoz \rkey, \vpz, \cmz  )$ to the above integral, and combining the estimates \eqref{eq:IBP_cmz:sum1}, for $z\neq 0$, we prove 
\beq\label{eq:EE_Gz}
\bal
 |\rhoc \rhoz \rkey \pa_z G(\xx)|
  & = \rhoc \rhoz \rkey (\vpz)^{-1} \cdot | \vpz \pa_z G |  \les  \e^{ -\kp_1 } \lla{G}     + \e \llw{G}+  \e + \nnn{\ww} .
\eal
\eeq

Since $\rkey \asymp 1$ by \eqref{def:rkey}, we prove Proposition \ref{prop:Cz} 
for $z \neq 0$. Since $\eta \in C^{1,\al}$ locally and $\eta$ is odd in $z$ by Proposition \ref{prop:reg_solu}, 
using $z \pa_z G = z \pa_z \f{\eta}{\waa}
= \f{z}{\waa} \pa_z \eta - \f{z \pa_z \waa}{\waa} \cdot \f{\eta}{\waa}$, 
the vanishing near $\xx=0$ \eqref{eq:local_reg}, 
and continuity, we also prove Proposition \ref{prop:Cz} with $z=0$.

\subsection{Weighted Estimate of $\pa_r G$}\label{sec:EE_C1r}

In this section, we perform weighted $\cW^{1,\infty}$ estimates on \eqref{eq:lin_dG}.

\subsubsection{Integrating factor $\vpr$}

We design $\vpr > 0$ to capture the local terms 
$ (\bar \cR - \pa_r Q^r ) \cdot \pa_r G$ in \eqref{eq:lin_dG:Gr}, which plays a role as an integrating factor. 
Following \eqref{eq:lin_vpz}, we construct $\vpr$ 
\bseq\label{eq:lin_vpr}
\beq
  \QQ  \cdot \na \vpr
 = - ( \bar \cR + \Lpsi - \pa_r Q^r) \vpr , \quad  \lim_{(r,z) \to 0} \f{\vpr}{r} ( \xx ) = 1.
\eeq
and design $\vpr$ as follows 
\beq
   \vpr(\xx)  =  r \exp\big( - \int_{-\infty}^0  ( \bar \cR  + \Lpsi- \pa_r Q^r + \f{Q^r}{r}  ) \cc X( \tau, \xx ) d \tau   \big) , \quad  \lim_{(r,z) \to 0} \f{\vpr}{r} ( \xx ) = 1.
\eeq
\eseq

Similar to \eqref{eq:lin_Cz}, for $z\neq 0$, we obtain 
\beq\label{eq:lin_Cr}
\f{d}{d s} ( \vpr \pa_r G) \cc X(s, \xx) =  (\vpr \cmr ) \cc X(s, \xx), 
 \quad     (\vpr  \pa_r G)(\xx) =
 (\vpr  \pa_r G) \cc X(-\tau, \xx) + \int_{-\tau}^0 \vpr \cmr \cc  X(s, \xx ) d s .
\eeq

For any $\xx=(r, z)$ with $z\neq 0, |\xx|\leq R$,  and $ s \leq 0$, using \eqref{eq:traj_ratio}, we obtain
\[
    \f{X^r(s, \xx) }{ |X^z(s, \xx)| } \les_{R, \e} \f{r}{|z|},
    \quad |X^z(s, \xx)| \asymp_{R, \e} e^{ 2 s} |z| ,
    \quad  |X^z(s, \xx ) | \neq 0 .
\]

Since $\vpr(r, z) \sim r$ near $0$, $z \neq 0$, and $X(s, \xx) \to 0$ as $s \to - \infty$ by \eqref{eq:traj_pf1}, using \eqref{eq:local_reg}, the above estimates, 
and $|\wwb| = |\waa| \asymp |z|$ for $|z| \leq 1$ by \eqref{eq:wa_upper_lower}, we obtain
\beq\label{eq:lin_Cr2}
  \limsup_{s \to -\infty}|\vpr \pa_r G|  \cc  X( s,  \xx )
\les_{R, \e}     \limsup_{s \to -\infty} \f{X^r(s , \xx )}{|X^z(s, \xx)|} \cdot |z \pa_r G |   \cc  X(s, \xx) 
\les_{R,\e} \f{r}{|z|} \limsup_{s \to -\infty} | \pa_r \eta |   \cc  X(s, \xx) 
= 0.
\eeq
Thus, the boundary term $ (\vpr  \pa_r G) \cc X(-\tau, \xx)$ in \eqref{eq:lin_Cr} 
vanishes if we take $\tau \to -\infty$.

By definition of $A^r$ \eqref{eq:damp:Ar} and $\vpr$, we obtain 
\beq\label{eq:damp:Ar2}
  A^r =    \f{\QQ \cdot \na ( \rhoc \rhor \rkey )}{ \rhoc \rhor \rkey }  + \bar \cR + \Lpsi - \pa_r Q^r 
  =  \f{\QQ \cdot \na ( \rhoc \rhor \rkey (\vpr)^{-1} )}{ \rhoc \rhor \rkey  (\vpr)^{-1}} .
\eeq

We aim to apply Lemma \ref{lem:traj_est} with weights $(\vp, m) \rsa ( \vpr, \rhoc \rhor \rkey )$ and establish the following estimates.

\begin{prop}\label{prop:Cr}
Suppose that $\ww$ satisfies the bound \eqref{eq:3D_size:ww}. Let $\rhor = \rhoo(\xx)  \rag(\xx) \cdot \la \xx \ra , \vpr,  \cmr$ be defined in \eqref{def:3d_wg}, \eqref{eq:lin_vpr},\eqref{eq:lin_Cr}. 
For any $\xx \in \bar B_R$ and $z \neq 0$, we have 
\beq\label{eq:EE_cmr0}
\bal
    |\rhoc \rhor  \pa_r G|&  \les 
    \e \big( 1 + \llz{G}  +  \e^{-\kp_1} \lla{G} + \e \llw{G} \, \big) \, ,
  \eal
\eeq
with implicit constants independent of $\e, \xx, R$. 

\end{prop}

\begin{remark}[\bf  Constraint: $z \neq 0$]

Near $z =0$, the term $|\rhor \pa_r G| = |\rhoo  \rag \ang \xx \f{1}{\wwb} \pa_r \eta|$ 
in \eqref{eq:EE_cmr0} is more singular than the energy $\llr{\eta}$ 
\eqref{def:local_norm} and $ \rhoo \rag \ang \xx \f{|z|}{\ang z \wwb} \pa_r \eta$ related to the norm $\nnr{\eta}$ \eqref{def:3d_norm}.  We \emph{do not} use the estimate for $\rhor \pa_r G$ on $z=0$ later. See more discussion in 
Section \ref{sec:EE_sum}.

\end{remark}

To prove the above result, we estimate the forcing term $\cmr$ \eqref{eq:lin_dG}, and then apply Lemma \ref{lem:traj_est}. We extract the singular term $r \pa_r \psi$ in $\Lpsi$ and $r\pa_r \bpsi$ in $\bar \QQ$ 
and decompose $\cmr$ in \eqref{eq:lin_dG} using \eqref{def:LR_decomp}
\beq\label{eq:decom_Cr}
\bal
    \cmr& =  - \pa_r Q^z \cdot \pa_z G  + \pa_r \bar \cR \cdot (G+1)  +  \pa_r \cL_{\psi} \cdot (G+1)    \\
    & = \B( - \pa_r \Qreg^z \cdot \pa_z G + \pa_r \Rreg 
    \cdot (G+1)  +  \pa_r \Lreg \cdot (G+1)  \B)   \\
    & \quad  + \B(- \pa_r \Qsin^z \cdot \pa_z G \B) +  \B(- \pa_r ( \Rsin + \Lsin ) \cdot (G+1) \B) 
    \teq \cmr_{R} + \cmr_{S, 1} + \cmr_{S, 2} ,
  \eal
\eeq
where $\cmr_R, \cmr_S$ denotes the regular and singular part in $\cmr$. We estimate $\cmr_R$ perturbatively, and defer the estimate to Section \ref{sec:est:cmr_R}. 
See also the ideas of estimates in Sections \ref{sec:idea_drW}, \ref{sec:idea_singular}.

\subsubsection{Integration by parts for $\cmr_{S, i}$}

Recall the decomposition \ref{def:LR_decomp}. Using $\Psi = \bpsi + \psi$, we have 
\beq\label{eq:decom_Cr_term}
\bal
 \cmr_{S, 1} & =  -\pa_r \Qsin^z \cdot \pa_z G =  - \pa_r ( r \pa_r \Psi ) \cdot \pa_z G , \\
 \quad \cmr_{S, 2} & = - \pa_r ( \Rsin^z + \Lsin ) ( G+1 )
 = \pa_r( r \pa_r \Psi \f{\pa_z \wwb}{\wwb} ) \cdot (G+1).
 \eal
\eeq

For these singular terms, we \textit{cannot} apply the nonlocal estimates in Proposition \ref{prop:vel_est} and use the energy to bound it directly, as motivated in \hyr[sec:idea_singular]{\its Integration by parts}. Instead, we use the crucial integration by parts in Lemma \ref{lem:IBP} 
with $(f_1, g_1) = ( -\vpr \pa_z G, r \pa_r \Psi )$ and $(f_2, g_2) = ( -\vpr (G+1),  r \pa_r \Psi \f{ \pa_z \wwb }{\wwb} )$. 
In view of the proof of Lemma \ref{lem:IBP}, we first derive identities \eqref{eq:IBP_pf} for 
$(f, g) = (f_i, g_i)$. 

Recall the  flow map 
 $X(s, \xx)$ \eqref{eq:traj} associated with $\QQ$. 
Below, we fix $\xx$ with $z \neq 0$ and simplify 
\[
X_s \teq X(s, \xx), \quad   \tf{d}{ds} X_s = \QQ \cc X_s.
\]

\vs{0.1in}
\paragraph{\bf Formulas for $f_1, g_1$}
For $(f_1, g_1) = ( -\vpr \pa_z G, r \pa_r \Psi )$, using a direct calculation, we yield 
\bseq\label{eq:3D_IBP2}
\beq\label{eq:3D_IBP2:a}
  - \f{d}{ds} ( r \pa_r \Psi ) \cc X_s \cdot ( \f{\vpr \pa_z G}{Q^r} ) \cc X_s
  = ( r \pa_r \Psi ) \cc X_s \cdot \f{d}{d s} \big ( \f{\vpr \pa_z G}{Q^r} \big) \cc X_s
   + \f{d}{d s} ( - r \pa_r \Psi \cdot \f{\vpr \pa_z G}{Q^r}  ) \cc X_s  \, ,
\eeq
where we denote $\f{d}{ d s} f \cc X_s = \f{d}{d s} f( X(s, \xx))$. We use  equations  \eqref{eq:lin_dG}, \eqref{eq:lin_vpr} again to simplify the first term on the right hand side, which seems to involves $\na^2 G$.  Using \eqref{eq:lin_dG}, \eqref{eq:lin_vpr}, 
and $ \f{d}{d s} Q^r \cc X_s = (\QQ \cdot \na Q^r) \cc X_s = ( Q^z \pa_z Q^r + Q^r \pa_r Q^r ) \cc X_s$, we derive
\begin{align}
  \tf{d}{d s} ( \vpr \cdot \pa_z G ) \cc X_s
 & = \B( -   (\bar \cR + \Lpsi- \pa_r Q^r ) \vpr \cdot \pa_z G   
+  \vpr \cdot \big(  ( \bar \cR + \Lpsi - \pa_z Q^z )  \pa_z G
+ \cmz  \big) \B) \cc X_s  \notag  \\
 & =  \B( \vpr \cdot \big( ( \pa_r Q^r - \pa_z Q^z) \pa_z G  
+ \cmz
  \big) \B) \cc X_s , \notag  \\
 \f{d}{d s} \big( \f{ \vpr \cdot \pa_z G }{Q^r} \big) \cc X_s
  & = \B( \vpr \f{1}{Q^r}   \big( ( \pa_r Q^r - \pa_z Q^z) \pa_z G  
+ \cmz
  \big)   - \vpr \pa_z G \f{ \QQ \cdot \na Q^r }{ (Q^r)^2} \B) \cc X_s \notag \\
  & = \B(  \vpr \cdot \f{1}{Q^r}  \big( ( -\pa_z Q^z - \f{Q^z}{Q^r } \pa_z Q^r ) \pa_z G + \cmz  \big) 
  \B) \cc X_s.
  \end{align}
\eseq
Since $\pa_z G$ is not $C^1$, the chain rule  $ \f{d}{ds} (\pa_z G) \cc X_s = (\QQ \cdot \na \pa_z G) \cc X_s$ does not hold.

Using the above identities and $\f{d}{ds}( r \pa_r \Psi) \cc X_s 
= (\QQ \cdot \na ) ( r \pa_r \Psi) \cc X_s $, we obtain
\bseq\label{eq:3D_IBP1a}
\begin{align}
& (\vpr \cmr_{S, 1}) \cc X_s = \B( \vpr  \big( - \pa_r (r \pa_r \Psi) \cdot  ( \pa_z  G ) \big) \B)  \cc X_s  \label{eq:3D_IBP1a:a}
  \\
&   = \B( \vpr \cdot \udb{ \f{Q^z}{Q^r} \pa_z ( r \pa_r \Psi ) 
  \cdot  \pa_z G }_{  \teq \cmr_{S, 11}  }  +
    \vpr \cdot \udb{  \f{r \pa_r \Psi  }{Q^r}  \B( ( -\pa_z Q^z - \f{Q^z}{Q^r } \pa_z Q^r ) \pa_z G + \cmz  \B) }_{ \teq \cmr_{S, 12}} \B) \cc X_s 
  + 
   \f{d}{d s} \big( \vpr \cdot \udb{  \f{ - \pa_z G  r \pa_r \Psi  }{ Q^r} }_{ \teq \cmr_{S, 13} } \big) \cc X_s \, , \notag
\end{align}
where we define $\cmr_{S, 1i}, i=1,2,3$ to be the above terms \emph{without} the factor $\vpr$.
The second and third terms in \eqref{eq:3D_IBP1a:a} corresponds to $\mw{LHS}_{ \eqref{eq:3D_IBP2:a}}$,
and $\cmr_{S,12},\cmr_{S,13}$ corresponds to $\mw{RHS}_{ \eqref{eq:3D_IBP2:a}}$.

From Proposition \ref{prop:reg_solu}, \eqref{eq:lin_dG}, 
and \eqref{eq:lin_vpr}, $(f_1, g_1)=  ( -\vpr \pa_z G, r \pa_r \Psi )$ satisfies the regularity assumption \eqref{eq:IBP_ass}.
Using Lemma \ref{lem:IBP} on integration by parts 
with $(f_1, g_1) = ( -\vpr \pa_z G, r \pa_r \Psi )$ or integrating the above identities along trajectory, for any $\tau > 0$, we derive 
\beq
   \int_{- \tau }^0  \vpr \cmr_{S, 1} \cc X(s, \xx) d s
   =    \int_{- \tau}^0 \vpr ( \cmr_{S, 11} + \cmr_{S, 12} )  \cc X(s, \xx) d s
   + \vpr \cmr_{S, 13} \cc X(s, \xx) \B|_{s=-\tau}^{s=0} . %(\xx) .
\eeq
\eseq

\paragraph{\bf Formulas for $f_2, g_2$}

For $(f_2, g_2) = ( -\vpr (G+1),  r \pa_r \Psi \f{ \pa_z \wwb }{\wwb} )$, 
using a direct calculation, we yield 
\bseq\label{eq:3D_IBP3}
\beq\label{eq:3D_IBP3:a}
\bal
\f{d}{d s} \big(  r \pa_r \Psi \f{ \pa_z \wwb }{\wwb}   \big) \cc X_s
 \, \cdot \,  (\f{ \vpr (G+1) }{Q^r} ) \cc X_s
 & = (  - r \pa_r \Psi \f{ \pa_z \wwb }{\wwb} ) \cc X_s \, \cdot \, \f{d}{d s} ( \f{ \vpr (G+1) }{Q^r} ) \cc X_s \\
  & \qquad + \f{d}{ds} \B(  r \pa_r \Psi \f{ \pa_z \wwb }{\wwb}  
 \, \cdot \f{ \vpr (G+1) }{Q^r}  \B) \cc X_s.
 \eal
\eeq

We use  equations  \eqref{eq:fix_pt_2D},\eqref{eq:lin_vpr}  for $G, \vpr$ again to simplify the first term on the right hand side:
\begin{align}
& \f{d}{ds} ( \f{\vpr (G+1)}{Q^r} ) \cc  X_s  \notag \\
  & \qquad =\B( \f{ - (\bar \cR + \Lpsi - \pa_r Q^r) \vpr (G+1) 
+  ( \bar \cR + \Lpsi(\ww)) (G+1)     \vpr }{Q^r} 
 - \vpr (G+1)  \f{\QQ \cdot \na Q^r}{ (Q^r)^2} \B) \cc X_s \notag \\
   & \qquad =  \B( \vpr (G+1) \big( \f{ \pa_r Q^r }{Q^r} - \f{\QQ \cdot \na Q^r}{ (Q^r)^2 }  \big) \B) \cc X_s\ = - \big( \vpr \cdot    \f{ Q^z \pa_z Q^r}{ (Q^r)^2} (G+1)   \big) \cc X_s \, .
 \end{align}
 \eseq

Recall $\cmr_{S, 2}$ from \eqref{eq:decom_Cr}. Combining the above identities and using \eqref{eq:decom_Cr_term}, we derive 
\bseq\label{eq:3D_IBP1b}
\begin{align}
 & (\vpr \cdot \cmr_{S, 2} ) \cc X_s =  \B( \vpr 
 \cdot  \pa_r( r \pa_r \Psi \f{\pa_z \wwb}{\wwb} )
  \cdot (G+1) \B) \cc X_s  \label{eq:3D_IBP1b:a}
   \\
 & = \B( \vpr  \B( \udb{  - \f{Q^z}{Q^r} \pa_z ( r \pa_r \Psi \f{ \pa_z \wwb }{\wwb} )
  \cdot   (G+1)    }_{ \teq \cmr_{S, 21} }
  +  \udb{ r \pa_r \Psi \f{ \pa_z \wwb }{\wwb} 
\cdot      \f{ Q^z \pa_z Q^r}{ (Q^r)^2} (G+1)     }_{\teq \cmr_{S, 22}}   \B) \B) \cc X_s
+ \f{d}{ds} \big(  \vpr \cdot \udb{ \f{  (G+1)  r \pa_r \Psi \f{ \pa_z \wwb }{\wwb} }{Q^r} }_{\teq \cmr_{S, 23}}  \big) \cc X_s \notag \, .
\end{align}
The second and third terms in \eqref{eq:3D_IBP1b:a} corresponds to $\mw{LHS}_{ \eqref{eq:3D_IBP3:a}}$,
and $\cmr_{S,22},\cmr_{S,23}$ corresponds to $\mw{RHS}_{ \eqref{eq:3D_IBP3:a}}$. From Proposition \ref{prop:reg_solu} and \eqref{eq:lin_vpr}, $(f_2, g_2) = ( -\vpr (G+1),  r \pa_r \Psi \f{ \pa_z \wwb }{\wwb} )$ satisfies the regularity assumption \eqref{eq:IBP_ass}. Using Lemma \ref{lem:IBP}  with $(f_2, g_2) $ and the above identities (or integrating the above identities along the trajectory), 
for any $\tau>0$, we derive 
\beq
   \int_{- \tau }^0  \vpr \cmr_{S, 2} \cc X(s, \xx) d s
   =    \int_{- \tau}^0 \vpr ( \cmr_{S, 21} + \cmr_{S, 22} )  \cc X(s, \xx) d s
   + \vpr \cmr_{S, 23} \cc X(s, \xx) \B|_{s =-\tau}^{s=0} .%  (\xx) .
\eeq

Since $\lim_{\xx \to 0} \f{1}{r} \vpr(\xx) = 1$ 
and $Q^r \gtr r$ \eqref{eq:3D_boot_q:1}, using the regularity estimate \eqref{eq:local_reg:nnn} 
for $G$ and estimate \eqref{eq:waa_reg} for $\wwb = \waa$, we estimate the boundary term 
$\cmr_{S, 13 }, \cmr_{S, 23}$ in  \eqref{eq:3D_IBP1a} and \eqref{eq:3D_IBP1b}
\beq
| \vpr
\cmr_{S, 13} | \les |\vpr \pa_z G \pa_r \Psi | \les |r \pa_z G \cdot z |\xx|^{\al} |  ,
\  | \vpr \cmr_{S, 23} |
\les |\vpr (G+1) \pa_r \Psi \f{\pa_z \wwb}{\wwb}|
\les \B| \f{r}{z} (G+1) \pa_r \Psi \B| .
\eeq
\eseq

From Proposition \ref{prop:reg_solu} and  \eqref{def:LR_decomp}, \eqref{eq:decom_Cr}, $ \cmr_R, \cmr_{S, ij}$ are continuous in $\bar B_R$ for any $R > 0$, including on $z=0$. Using estimate \eqref{eq:EE_cmr:1} below and  Lemma \ref{lem:traj_est}, 
we show that  $(\vpr \cmr_R) \cc X(s, \xx), (\vpr \cmr_{S, ij} ) \cc X(s,\xx)$ are
$L^1$-integrable over $s \in (-\infty,0]$. Since $X(-\tau, \xx) \to 0$ as $\tau \to \infty$ by \eqref{eq:traj_pf1}, using identities and vanishing in \eqref{eq:3D_IBP1a}, \eqref{eq:3D_IBP1b}, \eqref{eq:lin_Cr}, \eqref{eq:decom_Cr}, \eqref{eq:lin_Cr2}, and taking $\tau \to \infty$, we obtain
\begin{align}\label{eq:lin_dG:Gr2}
      \vpr \cdot \pa_r G(\xx) &=
 \lim_{\tau \to \infty} \B( (\vpr  \pa_r G) \cc X(-\tau, \xx) +  \int_{- \tau }^0  \vpr (\cmr_R + \cmr_{S, 1} + \cmr_{S, 2} ) \cc  X(s, \xx) d s \B) \\
    & = \int_{-\infty}^0 \vpr ( \cmr_R + \cmr_{S, 11} + \cmr_{S, 12} + \cmr_{S, 21} + \cmr_{S, 22} )
     \cc X(s, \xx) d s + \cmr_{S, 13}(\xx) + \cmr_{S, 23}(\xx).  \notag 
\end{align}

The key feature of derivations  \eqref{eq:3D_IBP1a}, \eqref{eq:3D_IBP1b} for $\cmr_{S , \cdot}$  is that there is \emph{no $\pa_r (r \pa_r \Psi)$ term}. Moreover, we gain the small factor $\f{q^z}{q^r} $ in $\f{Q^z}{Q^r} = \f{z}{r } \cdot \f{q^z}{q^r} $ \eqref{def:qq}. Next, we estimate each term \emph{perturbatively}.

\subsubsection{Estimate of regular part $\cmr_R$}\label{sec:est:cmr_R}

Recall the weight $\rhor = \rhoo \rag \cdot \ang \xx, \rhoz = \rhoo \rag |z|$ from \eqref{def:3d_wg}. 

Below, we estimate $\cmr_R$. Using \eqref{eq:Q_mix_est} in Proposition \ref{prop:Lpsi} and $\epa>0$, we estimate
\beq\label{eq:cmr_est:Gz}
\bal
 |\rhoc \rhor \rkey \pa_r \Qreg^z \cdot \pa_z G|
 & \les \rhoc \rhoo \rag \cdot \ang \xx \cdot |z| \CCa(\ww, \xx) |\pa_z G|  \\
 & \les \ang \xx \CCa(\ww,\xx ) 
 \cdot \rhoc \rhoo \rag |z \pa_z G| 
\les \ang \xx \CCa(\ww,\xx) \llz{G}.
 \eal
\eeq

Using Proposition \ref{prop:3D_error} on $\bar \cR$, 
definition \eqref{def:local_norm} of $\lla{\cdot}$,
$\rag \les 1$, and \eqref{eq:ep_ineq1} on $\rhoc \rhoo$, we estimate 
\beq\label{eq:cmr_est:cR}
\bal
  \rhoc \rhor \rkey |\pa_r \Rreg \cdot (G+1)|
  & \les \rhoc \rhoo \rag \cdot \ang \xx \min( |\xx|^{\al},  |\xx|^{\al-\hal-1 } ) (G+1) \\
 & \les 
    \lla{G}  \ang \xx^{\al-\hal}  + \rhoc \rhoo \min( |\xx|^{\al},  |\xx|^{\al-\hal } ) \\
    &   \les  \lla{G}  \ang \xx^{\al-\hal}  +  \ang \xx^{\al-\hal} .
\eal
\eeq

Using \eqref{eq:Lreg_Cr_est} on $\pa_r \Lreg$, the norm $\lla{G}$ \eqref{def:local_norm},  $\rhoc = \ang \xx^{\e^2}$ \eqref{def:rhoc}, 
and $\cJaa \les \e^{-1}$, we bound 
\begin{align}\label{eq:cmr_est:Lreg_r}
\rhoc  \rhor  \rkey |\pa_r \Lreg \cdot (G+1) |
 & \les  |\ang \xx \pa_r \Lreg|  \cdot |\rhoc \rhoo \rag G|
 + \rhoc \rhoo \rag \cdot |\ang \xx \pa_r \Lreg| \notag  \\
  & \les  \cJak \ang \xx^{\al-\hal + \epa} \nnn{\ww}  \lla{G} 
  + \rhoc   \min( |\xx|^{\al}, |\xx|^{\al-\hal} ) \nnn{\ww}  \notag  \\
 & \les  \e^{-\kp} \ang \xx^{\al-\hal + \epa} \nnn{\ww}  \lla{G} 
 + \ang \xx^{\al- \hal + \e^2} \nnn{\ww}.
\end{align}

We will combine the above estimates for $\cmr_R$ in \eqref{eq:cmr_est:Mreg}.

\subsubsection{Estimate each parts in $\cmr_S$}\label{sec:EE_cmr_S}

 By imposing \eqref{eq:3D_size:ww}, and using \eqref{eq:3D_boot_q:1} on $q^r, q^z$,  we obtain
\beq\label{eq:qq_est:1}
Q^r = r q^r, \quad  Q^z = z q^z, \quad 
  |q^z| \les \cJaa(\xx),  \quad q^r \asymp \e^{-1} .
\eeq

We recall the following estimates for $\Psi$ from \eqref{eq:tot_u_est1} and $\CCa$ from \eqref{def:CC_bounds}
\beq\label{eq:tot_u:recall}
\bal
  & | \pa_r \tf{1}{z} \Psi|  +    |\pa_{r z} \Psi| \les  \CCa(\ww, \xx) .
  \eal
\eeq

\paragraph{\bf Estimate of $\cmr_{S,11}$}
For the integrand in $\cmr_{S,11}$ in \eqref{eq:3D_IBP1a}, using \eqref{eq:tot_u:recall} for $\pa_{r z} \Psi$ and \eqref{eq:qq_est:1}
 for $q^r, q^z$, and $ \cJaa \les \e^{-1}$ \eqref{eq:Ja_hat}, we obtain
\beq\label{eq:cmr_IBP:I11}
\bal
\rhoc \rhor \rkey  \cmr_{S,11}|  & = 
  \B| \rhoc \rhor \rkey \f{Q^z}{Q^r} \pa_z ( r \pa_r \Psi ) 
  \cdot  \pa_z G \B|
\les \ang \xx  | \f{q^z}{q^r} \pa_{r z} \Psi | \cdot  \rhoc \rhoo \rag |z \pa_z G| \\
& \les  \e \cJaa
\cdot
\ang \xx \CCa(\ww, \xx)  \llz{G}  \les \ang \xx \CCa(\ww, \xx)  \llz{G} .
\eal
\eeq

\vs{0.1in}

\paragraph{\bf Estimate of $\cmr_{S,12}$}
For $\cmr_{S,12}$, using \eqref{eq:3D_IBP1a} and \eqref{eq:3D_IBP2}, and $\wwb=\waa$,  we bound 
\[
\bal
  |\rhoc \rhor \rkey   \cmr_{S,12}|
    & \les  \rhoc \rhoo \rag \ang \xx |r \pa_r \Psi| 
  \cdot \f{1}{Q^r}  \cdot \B| ( -\pa_z Q^z - \f{Q^z}{Q^r } \pa_z Q^r ) \pa_z G + \cmz  \B|  \\
  & \les  \ang \xx |r \pa_r \f{ \Psi}{z}| \cdot 
  \rhoc \rhoo \rag |z| \f{1}{Q^r}  \cdot  \B| ( -\pa_z Q^z - \f{Q^z}{Q^r } \pa_z Q^r ) \pa_z G + \cmz  \B| .
  \eal
\]

Applying \eqref{eq:tot_u:recall} to $r \pa_r \f{\Psi}{z}$,  and \eqref{eq:qq_est:1}, we obtain
 \[
 P_1 = \ang \xx |r \pa_r \f{ \Psi}{z}| \cdot 
  \rhoc \rhoo \rag |z| \f{1}{Q^r}
  \les   \ang \xx \CCa(\ww, \xx)
\cdot \rhoc \rhoo \rag |z| \f{1}{q^r} 
\les \e \ang \xx \CCa(\ww, \xx) \cdot \rhoc  \rhoo \rag |z| .
 \]

Applying \eqref{eq:Q_mix_est} and \eqref{eq:Q_bound} for $\QQ$ and 
 and \eqref{eq:qq_est:1} for $q^r, q^z$ ,  we estimate 
\[
\bal
  P_2  & = |( -\pa_z Q^z - \f{Q^z}{Q^r } \pa_z Q^r )|
  \les %\cJaa( 1 + \nnn{\ww} ) + |\xx| \CCa(\ww, \xx)
  1 + \CCj(\ww, \xx)
  + \f{z}{r} \cdot \e \cJaa \f{r}{\ang z} \ang \xx \CCa(\ww, \xx) .
\eal
\]

Since $ \e \cJaa(\xx) \les 1, \cJaa \gtr 1$, using \eqref{eq:3D_size:ww} and \eqref{eq:3D_size:C},  we obtain
\[
  P_2 \les 1 + \CCj(\ww, \xx) + \ang \xx \CCa(\ww, \xx) 
  \les 1 + \cJaa + \ang \xx \cdot \ang \xx^{-1} \les  \cJaa.
\]

Combining the estimates of $P_1, P_2$, 
and using $\rhoz = \rhoo \rag |z|$, we bound 
\beq\label{eq:cmr_IBP:I12}
\bal
    |\rhoc \rhor \rkey   \cmr_{S,12}|
  & \les \e \ang \xx \CCa(\ww, \xx) \cdot \rhoc \rhoo \rag |z| \cdot (  \cJaa | \pa_z G | + |\cmz| )\\
  & \les \e  \ang \xx \CCa(\ww, \xx) \cdot ( \cJaa \llz{G} + \rhoc  \rhoz |\cmz| ) .
\eal
\eeq

We recall the estimate of $  \rhoc \rhoz |\cmz|$ from \eqref{eq:cmz_est:tot}.
\[
\bal
    \rhoc \rhoz \rkey |\cmz|  
  & \les  
  \f{r}{\ang \xx} \cJaa |\lgp z|^{- \kp_1-1} \cdot \e^{- \kp_1} \lla{G}   \\
 & \quad  +   r \CCa(\ww, \xx) \llw{G}
 + r \ang \xx^{\al-\hal-1} 
+   \lla{G} 
+ \nnn{\ww} .
\eal
\]

Since we have a small factor $\e$ in \eqref{eq:cmr_IBP:I12}, we derive a rough bound of 
$ \rhoc \rhoz |\cmz|$.  Using $\rkey \asymp 1$, $ \f{r}{\ang \xx}, |\lgp z|^{-1} \les 1$, $ \nnn{\ww }\les \e^{\hk } $ \eqref{eq:3D_size:ww},  $ r \CCa(\ww,\xx) \les 1$ by \eqref{eq:3D_size:C},
and $\cJaa \les \e^{-1}$, we bound 
\beq\label{eq:cmr_IBP:I12_2}
\bal
  \rhoc \rhoz |\cmz| &  \les \e^{-1} \cdot \e^{-\kp_1} \lla{G}
  + \llw{G} + 1 +  \lla{G} + \e^{\hk}  \\
  & \les  1 + \e^{-1 - \kp_1}   \lla{G} + \llw{G} .
  \eal
\eeq

Combining \eqref{eq:cmr_IBP:I12}, \eqref{eq:cmr_IBP:I12_2}, 
and using $\e \cJaa \les 1$, $\ang \xx \CCa(\ww,\xx) \les 1$ from \eqref{eq:3D_size:C}, we obtain
\begin{align}\label{eq:cmr_IBP:I12_tot}
      |\rhoc \rhor \rkey   \cmr_{S,12}| 
      & \les \e  \ang \xx \CCa(\ww, \xx) \cdot ( \cJaa \llz{G} +   \rhoc \rhoz |\cmz| 
       )
       \notag \\
     & \les \e  \ang \xx \CCa(\ww, \xx) \cdot ( \e^{-1} \llz{G} +  1 +  \e^{-1 -\kp_1} \lla{G} + 
     \llw{G}  )   \notag \\
     & \les   \ang \xx \CCa(\ww, \xx) \cdot (  \llz{G} +  \e +  \e^{-\kp_1} \lla{G} 
     + \e \llw{G} ).  
\end{align}

\vs{0.1in}

\paragraph{\bf Estimate of $\cmr_{S,13}$}
For $\cmr_{S,13}$ in \eqref{eq:3D_IBP1a}, using \eqref{eq:tot_u:recall} for $\pa_r \Psi$, we bound 
\beq\label{eq:cmr_IBP:I13}
\bal
  |\rhoc \rhor \rkey  \cmr_{S,13}| 
  & \les \rhoc \rhoo \rag \ang \xx \cdot  | \f{  \pa_z G  r \pa_r \Psi  }{ Q^r} |
  \les \rhoc \rhoo \rag  \ang \xx \cdot  \f{ |\pa_z G | }{q^r} \cdot  |z| \CCa(\ww, \xx) \\
  & \les  \e \ang \xx \CCa(\ww, \xx) \llz{G}.
\eal
\eeq

\vs{0.1in}

Next, we estimate $\cmr_{S, 2\cdot}$.
\paragraph{\bf Estimate of $\cmr_{S,21}$}

For $\cmr_{S,21}$ in \eqref{eq:3D_IBP1b}, using \eqref{eq:qq_est:1} for $Q^r, Q^z$, $\rhor =\rhoo \rag \ang \xx$, we have
\[
\bal
   |\rhoc \rhor \rkey  \cmr_{S,21} | & \les \rhoc \rhoo \rag \ang \xx  \B| \f{Q^z}{Q^r} \pa_z ( r \pa_r \Psi \f{ \pa_z \wwb }{\wwb} ) 
  \cdot (  G+1 ) \B|   \\
  &  \les \rhoc \rhoo \rag \ang \xx  \e \cJaa
  \cdot \f{|z|}{r}  \B| \pa_z ( r \pa_r \Psi \f{ \pa_z \wwb }{\wwb} ) 
  \cdot (  G+1 ) \B| .
  \eal
\]

Using $\wwb = \waa$, \eqref{eq:1D_prof_dxx}, and \eqref{eq:1D_prof_est2:c} 
\beq\label{eq:cmr_est:pf1}
  |z \pa_z \f{\pa_z \wwb}{\wwb} | \les |z|^{-1} ,\quad  |\f{\pa_z \wwb}{\wwb}| \les |z|^{-1}, 
\eeq
and \eqref{eq:tot_u:recall}, we obtain
\[
  \f{|z|}{r} \B| \pa_z ( r \pa_r \Psi \f{ \pa_z \wwb }{\wwb} ) \B|
  \les  |\pa_{r z} \Psi| + |\f{1}{z} \pa_r \Psi| \les \CCa(\ww, \xx).
\]

Using the above estimate and $\e \cJaa \les 1$ by \eqref{eq:Ja_hat}, we yield 
\[
   |\rhoc \rhor \rkey  \cmr_{S,21} | \les \rhoc \rhoo \rag \ang \xx \cdot \CCa(\ww, \xx)
  \cdot | G+1 |.
\]

Using the norm $\lla{\cdot}$ on $G$  \eqref{def:local_norm}, $\rag \les 1$ \eqref{def:3d_wg}, and 
\eqref{eq:3D_size:C2},  we prove
\beq\label{eq:cmr_IBP:I21}
\bal
      |\rhor \rkey  \cmr_{S,21} | \les  \rhoc \rhoo \rag  \ang \xx \CCa(\ww, \xx) | G + 1|
 & \les \ang \xx \CCa(\ww, \xx) \lla{G} 
 + \rhoc \rhoo \cdot \ang \xx \CCa(\ww,\xx)  \\
& \les \ang \xx \CCa(\ww, \xx) \lla{G}  + \ang \xx^{\al-\hal + \e^2}.
\eal
\eeq

\paragraph{\bf Estimate of $\cmr_{S,22}$}

For $\cmr_{S,22}$ in \eqref{eq:3D_IBP1b}, 
we bound 
\[
     |\rhoc \rhor \rkey  \cmr_{S,22} | 
  \les \rhoc \rhor  | r \pa_r \Psi \f{ \pa_z \wwb }{\wwb}  | 
  \cdot \f{1}{Q^r} \cdot    \B|  \f{Q^z}{ Q^r} \pa_z Q^r  \cdot  (G+1) \B|.
\]

Using \eqref{eq:cmr_est:pf1}, \eqref{eq:tot_u:recall} for $\pa_r \Psi$, and \eqref{eq:qq_est:1} for $q^r$ we bound 
\[
    P_{22,1} \teq | r \pa_r \Psi \f{ \pa_z \wwb }{\wwb}  |  \cdot  |\f{1}{Q^r}| 
     \les |\f{1}{z} \pa_r \Psi| \f{1}{q^r} \les  \e  \CCa(\ww, \xx).
\]

Using \eqref{eq:qq_est:1} for $q^r, q^z$, 
 \eqref{eq:Q_mix_est} for $\pa_z Q^r$, and $\e \cJaa \les 1$, we bound the remaining term
\[
\bal
  P_{22,2}  & \teq |  \f{Q^z}{Q^r} \pa_z Q^r | 
  \les 
   \e \cJaa \f{|z|}{r} \cdot  \f{r}{\ang z} \ang \xx \CCa(\ww, \xx) 
   \les \ang \xx \CCa(\ww, \xx) .
  \eal
  \]

  Using the above estimates,  \eqref{eq:cmr_IBP:I21} for $G+1$, and $\ang \xx \CCa(\ww, \xx) \les 1$ from \eqref{eq:3D_size:C}, we prove
\beq\label{eq:cmr_IBP:I22}
\bal
       |\rhoc \rhor \rkey  \cmr_{S,22} | 
    &   \les P_{22,1}   P_{22, 2} \cdot |\rhoc \rhor (G+1)| 
      \les \e \ang \xx  \CCa(\ww, \xx) \cdot  \CCa(\ww, \xx) \rhoc \rhoo \rag \ang \xx | G+1 | \\
    &   \les  \e \ang \xx \CCa(\ww,\xx) \cdot 
     (  \ang \xx \CCa(\ww, \xx)  \lla{G}+ \ang \xx^{\al-\hal + \e^2} )   \\
   & \les \e \ang \xx \CCa(\ww,\xx) \cdot  (  \lla{G}+ \ang \xx^{\al-\hal + \e^2} ).
\eal 
\eeq

\vs{0.1in}

\paragraph{\bf Estimate of $\cmr_{S,23}$}

Recall $\rhor = \rhoo \rag \ang \xx$, and $\cmr_{S,23}$ from \eqref{eq:3D_IBP1b}. Using \eqref{eq:cmr_est:pf1} for $\f{\pa_z \wwb}{\wwb}$, \eqref{eq:qq_est:1} for $q^r$, 
  we bound 
\[
\bal
  |\rhoc \rhor \rkey  \cmr_{S,23}|
  & \les \rhoc \rhoo \rag \ang \xx \cdot
  \B| \f{  (G+1)  r \pa_r \Psi \f{ \pa_z \wwb }{\wwb} }{Q^r} \B|
\les  \rhoc \rhoo \rag \ang \xx \cdot \e \B|   (G+1)   \pa_r \Psi \cdot \f{1}{z}  \B| \, . 
\eal
\]

Using \eqref{eq:tot_u:recall}  for $\pa_r \Psi$ and \eqref{eq:cmr_IBP:I21} for $G+1$, we further obtain
\beq\label{eq:cmr_IBP:I23}
\bal
  |\rhoc \rhor \rkey  \cmr_{S,23}| & \les \rhoc \rhoo \rag \ang \xx \cdot \e |G+1|  \CCa(\ww, \xx)
\les  \e   \ang \xx \CCa(\ww, \xx)  \lla{G}+ \e \ang \xx^{\al-\hal + \e^2} .
\eal
\eeq

\vs{0.1in}
\paragraph{\bf Summary of boundary terms}

Combining the estimates of $\cmr_{S,13}$ 
in \eqref{eq:cmr_IBP:I13} and $\cmr_{S,23}$ in \eqref{eq:cmr_IBP:I23},  using $\ang \xx^{\al-\hal + \e^2} \les 1$ by \eqref{ran:ep_all}, and $\ang \xx \CCa(\ww, \xx) \les 1$ by \eqref{eq:3D_size:C},   we estimate  boundary terms as
\beq\label{eq:EE_cmr:bd}
\bal
   \rhoc \rhor \rkey  ( |\cmr_{S,13} | +| \cmr_{S,23} | )
   & \les \e \ang  \xx \CCa(\ww, \xx) ( \llz{G} + \lla{G} )
   + \e \ang \xx^{\al-\hal + \e^2} \\
   & \les \e  ( \llz{G} + \lla{G} + 1 ).
\eal
\eeq

\subsubsection{Trajectory estimates}\label{sec:EE_Cr_sum}

To apply Lemma \ref{lem:traj_est} with $(m, \vp) \rsa (\rhor \rkey, \vpr)$, we shall further compare the bound of $\cmr_R, I_{ij}$ relative to the damping terms $A^r$ \eqref{eq:damp:Ar}, \eqref{eq:damp:Ar2}. We recall
\[
  A^r  =  \f{\QQ \cdot \na ( \rhoc \rhor \rkey  )}{ \rhoc \rhor \rkey  }  
 \leq - \bar \lam_1 \B( \e^{-1} \ang \xx^{\al-\hau}  + \cJaa(\xx) |\lgp z|^{ -\kp_1 -1 } +1 \B) < 0.
\]

Using \eqref{eq:CC_bound2}, we obtain
\bseq\label{eq:CC_bound:1}
\beq
\bal
  \ang \xx \CCa(\ww,\xx)
  \les \ang \xx ( \ang \xx^{\al-\hau - 1}  + \e \ang \xx^{-1} )
  \les   \ang \xx^{\al-\hau }  + \e \les \e |A^r|.
\eal
\eeq

Moreover, using \eqref{eq:interp:pow2},  we have 
\beq
\ang \xx^{-\e + \b + \e^2} + 
  \ang \xx^{\al-\hal} + \ang \xx^{\al-\hal + \e^2 }
\les \ang \xx^{\al-\hau} + \e^2
   \les \e |A^r| ,  \quad
   1 \les |A^r| .
\eeq

Using \eqref{eq:interp:pow} with $ a = \e$, we bound 
\beq
 \ang \xx^{\al -\hal + \epa} 
  \les \e^{-\kp_1} \ang \xx^{\al-\hau} + \e 
  \les \e^{1 -\kp_1 }  |A^r| .
\eeq

\eseq

\paragraph{\bf Relative bounds}

Recall $\cmr_R$ from \eqref{eq:decom_Cr}. Combining estimates \eqref{eq:cmr_est:Gz}, \eqref{eq:cmr_est:cR}, and \eqref{eq:cmr_est:Lreg_r},   we estimate $\cmr_R$ as
\beq\label{eq:cmr_est:Mreg}
\bal
 \rhoc  \rhor \rkey | \cmr_R| & \les  \udb{ \ang \xx \CCa(\ww, \xx)  \llz{G} }_{ \eqref{eq:cmr_est:Gz} } 
     +  \udb{ \lla{G}  \ang \xx^{\al-\hal}  + \ang \xx^{\al-\hal } }_{  \eqref{eq:cmr_est:cR} } \\
   & \qquad  + \udb{ \e^{ -\kp } \ang \xx^{\al-\hal + \epa} \nnn{\ww}  \lla{G}  
  +   \ang \xx^{\al- \hal + \e^2} \nnn{\ww}  }_{ \eqref{eq:cmr_est:Lreg_r} } .
  \eal
\eeq

Using the key estimates in \eqref{eq:CC_bound:1}, $\e \les 1$, 
and $\ang \xx^{\al-\hal + \epa}\les 1$,  we further  bound $\cmr_R$:
\[
\bal
 \rhoc \rhor \rkey | \cmr_R| 
     & \les |A^r| \big(  \e  \llz{G} 
     +  \e \lla{G} + \e 
+ \e^{1-\kp_1 - \kp} \nnn{\ww} \lla{G}
+ \e \nnn{\ww}       
    \big). 
\eal
\]

Since $\nnn{\ww}\leq \e^{\hk}<1$ from \eqref{eq:3D_size:ww} and  $\e^{1-\kp_1 -\kp} \nnn{\ww} \leq \e^{ \mhk + 1 -\kp_1 -\kp}
\les \e $ due to $ \f{1-\kp}{2} - \kp_1 > 0$ \eqref{def:kp},  we further obtain
\[
   \rhoc \rhor \rkey | \cmr_R| 
     \les \e |A^r| (    \llz{G} 
     +   \lla{G} + 1 
+  \nnn{\ww}       
    )
    \les \e |A^r| (    \llz{G} 
     +   \lla{G} + 1      ) .
\]

For $I_{ij}, j\leq 2$,  $\cmr_{S,11}$
\eqref{eq:cmr_IBP:I11}, $ \cmr_{S,12}$ \eqref{eq:cmr_IBP:I12_tot},  $\cmr_{S,21}$  \eqref{eq:cmr_IBP:I21},  similarly, we estimate 
\[
\bal
\rhoc \rhor \rkey  \cmr_{S,11}| &\les 
\ang \xx \CCa(\ww, \xx)  \llz{G} 
\les \e |A^r|  \llz{G} , \\
\rhoc \rhor \rkey  \cmr_{S,12}| &\les 
    \ang \xx \CCa(\ww, \xx) \cdot (  \llz{G} +  \e +  \e^{-\kp_1} \lla{G} 
     + \e \llw{G} )  \\
 &  \les  
\e |A^r| \cdot (  \llz{G} +  \e +  \e^{-\kp_1} \lla{G} 
     + \e \llw{G} ) , \\
|\rhoc \rhor \rkey  \cmr_{S,21} |     & \les   \ang \xx \CCa(\ww, \xx) \lla{G} + \ang \xx^{\al -\hal + \e^2}  \\
& \les \e |A^r| (   \lla{G} + 1  ) \, . \\
\eal
\]

For $\cmr_{S,22}$, using 
 \eqref{eq:cmr_IBP:I22} for $\cmr_{S,22}$, estimates \eqref{eq:CC_bound:1} for $A^r$, and $\ang \xx^{\al-\hal} \les 1$, we bound 
\[
\bal
    |\rhoc \rhor \rkey  \cmr_{S,22} | 
    &   \les  \e \ang \xx \CCa(\ww,\xx)
  \cdot  (     \lla{G}+  \ang \xx^{\al-\hal} ) 
   \les \e^2 |A^r| (        \lla{G}+ 1  )  .
  \eal
\]

Combining similar terms in the above estimates and using $\e \leq 1$, we bound 
\beq\label{eq:EE_cmr:1}
  \bal
 & \rhoc \rhor \rkey ( |\cmr|  + 
 |\cmr_{S,11}| + |\cmr_{S,12}| + |\cmr_{S,21}| + |\cmr_{S,22}|   ) \\ 
& \quad \les \e |A^r| \cdot  ( 1 + \llz{G}  +  \e^{-\kp_1} \lla{G} + \e \llw{G} ) .
\eal
\eeq

Applying Lemma \ref{lem:traj_est} with $(m, \vp, H) \rsa ( \rhoc \rhor \rkey, \vpr, \cmr_R + \cmr_{S,11} + \cmr_{S,12} + \cmr_{S,21} + \cmr_{S,22}) $ to equation \eqref{eq:lin_dG:Gr2}, using 
\eqref{eq:EE_cmr:1}, and the estimate of the boundary term \eqref{eq:EE_cmr:bd}, we prove 
\[
\bal
  |\rhoc \rhor \rkey \pa_r G|
  & = |\rhoc \rhor \rkey (\vpr)^{-1} \cdot \vpr \pa_r G|  \\
  & \les 
  \e  ( \llz{G} + \lla{G} + 1 )
  + \e \B( 1 + \llz{G}  +\nnn{\ww}  + \e^{-\kp_1} \lla{G} + \e \llw{G}\B) \\
  & \les \e ( 1 + \llz{G}    + \e^{-\kp_1} \lla{G} + \e \llw{G} ).
\eal
\]

Since $\rkey \asymp 1$ by \eqref{def:rkey},  we prove Proposition \ref{prop:Cr}.

\subsection{Proof of into estimates}

We use a best-constant argument to close the estimate, analogous to a bootstrap argument  for evolutionary nonlinear PDEs. We show that the perturbation has size \(O(\e^\g)\) with \(\g\) close to \(1\), which is determined by the size of error $\bar \cR$ in the fixed point estimate.

\vs{0.05in}

\paragraph{\bf Functional Space}

 Recall the parameter $\kp_1$ from \eqref{def:kp}. We introduce a functional space $\cY_{\e}$ for the fixed point argument and the $\cY$-norm as follows 
\bseq\label{def:cY} 
\begin{align} 
\| \ww \|_{\cY} & \teq \max( \| \ww \|_{\cXc} , \nnr{\ww}, \nnz{\ww} ) , \\
   \cY_{\e}  & \teq \B\{ \ww : \| \ww \|_{\cXc} \leq \e^{1- 2 \kp_1}, 
   \ 
   \nnr{\ww} \leq \e^{1 - \kp_1}, \   \nnz{\ww} \leq \e^{ 1 - 6 \kp_1 } 
      \B\}  \, . 
\end{align}
\eseq

 By definition \eqref{def:cY} and \eqref{norm:comp_bd}, 
 and Lemma \ref{lem:linf_R3}, we obtain
\beq\label{norm:comp_Y}
 \nnla{\ww} + \nnrr{\ww}  + \nnn{\ww} 
  \les \| \ww \|_{\cY} + \nna{\ww}
\les \| \ww \|_{\cY} + \| \ww \|_{\cXc} \les \| \ww \|_{\cY}.
\eeq

For $\ww \in \cY_{\e}$, since $\e \leq 1$ and $\al \in [0.32, 0.34]$, using  \eqref{norm:Xc},
the definition \eqref{def:3d_norm}, and Lemma \ref{lem:linf_R3},  we obtain
\beq\label{eq:ww_bound}
\bga
   \nna{ \ww } \les \| \ww \|_{\cXc} \les \e^{1 -2 \kp_1},
   \quad \nnrr{\ww} = \max(\nna{\ww}, \nnr{\ww}) \les \e^{1-2 \kp_1},  \\
   \quad \nnn{\ww} = \max(\nnrr{\ww},  \nnrr{\ww}^{\al} \nnz{\ww}^{1- \al} ) \les 
\e^{ (1 - 6 \kp_1) (1-\al) + \al (1 - 2 \kp_1) }
\les \e^{1 - 5\kp_1}.
\ega
\eeq

For $\e$ small enough, since $1 - 10 \kp_1 > 1 - \mhk = \f{1 + \kp}{2}$, 
the above estimates imply \eqref{eq:3D_size:ww}.

\subsubsection{Summary of the energy estimates}\label{sec:EE_sum}
For any $R > 1$, and $\xx \in \bar B_R$ and $z \neq 0$, under the assumption \eqref{eq:3D_size:ww} on $\ww$,  Proposition \ref{prop:Cr} proves
\beq\label{eq:EE_tot1}
\bga
      |\rhoc \rhor  \pa_r G|  \les  \cC_{r, R}(\ww, G), 
      \quad  \cC_{r, R} \teq 
     \e \B( 1 + \llz{G}  + \e^{-\kp_1} \lla{G} + \e \llw{G} \B) . 
  \ega
\eeq

Recall $G =\f{\eta}{\wwb}$. From estimate \eqref{eq:wa_upper_lower}, we obtain $ |\f{\wwb}{z}|  \gtr_R 1$  for any $|z| \leq R$.
Since $\eta(r, z)$ is odd in $z$, $\eta(r,0) =0$, using Item (b)  in Proposition \ref{prop:reg_solu}, 
definition of $\rhoo, \rhoc$ \eqref{def:3d_wg}, and $\al - \bbb  > 0$ by \eqref{norm:Xc}, \eqref{ran:ep_all}, for $\xx \in \bar B_R$, we obtain 
\[
\bal
I_r  &\teq \rhor \rkey \rhoc \tf{z}{\ang z} \pa_r G,   & \quad 
 | I_r | & \les_{\e, R} |\xx|^{-\bbb + 1 } |\pa_r \eta| \les_{\e, R} |z|^{\al} |\xx|^{-\bbb + 1 } \| \eta \|_{C^{1,\al}(B_{2R}) } , \\
I_0 & \teq \rhoc \rhoo \rag G,  & \quad  | I_0 | & \les_{\e, R}  |\xx|^{-\bbb+1} \tf{| \eta| }{|z|} 
  \les_{\e, R} |\xx|^{\al + 1 -\bbb } \| \eta \|_{C^{1,\al}(B_{2R}) }.
\eal
\]
The term $I_r$ vanishes both near $\xx=0$ and near $z=0$. Since the weights in $I_r, I_0$ are singular only at $\xx =0$ and $\eta \in C^{1,\al}$, we obtain that $I_r, I_0$ are continuous in $\R_+ \times \R$. 
Thus, estimate \eqref{eq:EE_tot1} and continuity of $I_r$ implies 
 \beq\label{eq:EE_tot1:Cr}
 \bal         
 | \tf{z}{\ang z} \rhoc \rhor  \pa_r G|  \leq \cC_{r, R}(\ww, G) , \quad \forall \xx \in \bar B_R.
    \eal
 \eeq

From the definition of $\rhoo , \rhoc$, the weight $\rhoo \rhoc$ satisfies the assumption in Lemma \ref{lem:linf_R3}. Thus, for $\xx \in \bar B_R$ with $z \neq 0$, using \eqref{eq:linf_R3:a} in Lemma \ref{lem:linf_R3} with $f \rsa G= \f{\eta}{\wwb} = \f{\eta}{\waa}$ and $\rhoo \rsa \rhoo \rhoc$, 
and then estimate \eqref{eq:EE_tot1} on $\pa_r G$, and estimate  \eqref{norm:comp_bd} with $f = \eta \rhoc$, we bound 
\beq\label{eq:EE_tot1:C0}
\bal
  |\rhoc \rhoo \rag G(r, z) |&  \les | (\rhoc \rhoo G)(0, z) | 
  + \max\nolimits_{\tr \leq r} | \rhoc \rhoo \rag \ang{\tr, z} \pa_r G(\tr, z) |  \\
& \les 
| (\rhoc \rhoo |\waa|^{-1} \eta )(0, z) |  + \cC_{r, R}(\ww, G).
\les  \| \rhoc \eta(0, z) \|_{\cXc} + \cC_{r, R}(\ww, G).
\eal 
\eeq

By continuity of $I_0$, \eqref{eq:EE_tot1:C0} holds  for \emph{any} $\xx \in \bar B_R$.

Using Propositions \ref{prop:3D_bd}, \ref{prop:Cz}, for any $\xx \in \bar B_R$,  we have
\begin{align}
       \| \rhoc \eta(0, \cdot) \|_{\cXc} & \leq \f{1 + \lamcL}{2} || \ww ||_{\cXc}  +  C \nnr{\ww}   +   C \e^{1 - \kp} \nnrr{\ww}  , \label{eq:EE_tot:a} \\ 
        |\rhoc \rhoz  \pa_z G(\xx)|  & \les 
           \e^{ -\kp_1 } \lla{G}     + \e \llw{G}+  \e + \nnn{\ww} . 
\label{eq:EE_tot:b}         
\end{align}

We recall the local norms $\llr{\cdot}, \llz{\cdot}, \lla{\cdot}$ from \eqref{def:local_norm}, the global norms  $ \nnn{\cdot}, \nnr{\cdot}$, $\nnz{\cdot}, \nna{\cdot}$ from \eqref{def:3d_norm}, 
and the 1D norm $\cXc$ from \eqref{norm:Xc}.

For $\ww$ satisfying \eqref{eq:3D_size:ww}, from Proposition \ref{prop:reg_solu}, we obtain $\llr{G}|, \llz{G}, \lla{G} < \infty$ for any $R >0$. Thus, taking maximum of the estimates \eqref{eq:EE_tot1:C0}, \eqref{eq:EE_tot1:Cr},
 \eqref{eq:EE_tot:b} over $\xx \in \bar B_R$, and using 
bound on $\nnn{\ww}$ \eqref{eq:ww_bound}, for $\ww \in \cY_{\e}$, we obtain
\beq\label{eq:EE_tot3}
\bal
  \lla{G} &\leq  C_* ( \| \rhoc \eta(0, z) \|_{\cXc} + \cC_{r, R}(\ww, G) ) , \\
  \llr{G} & \leq  C_* \cC_{r, R}(\ww, G) , \\
  \llz{G} & \leq  C_*
     (\e^{ -\kp_1 } \lla{G}     + \e \llw{G}+  \e + \e^{1 -5\kp_1} ) ,
\eal
\eeq
for some absolute constant $C_*$ independent of $\e$.

Below, for $\e$ small enough and $\ww \in \cY_{\e}$, we prove that  $\eta \in \cY_{\e}$. 

\vspace{0.1in}
\paragraph{\bf Estimate $\| \rhoc \eta \|_{ \cXc } $}

Using \eqref{eq:EE_tot:a}, $\lamcL < 1$,  $\rhoc \geq 1$ by \eqref{def:rhoc}, 
and $1-\kp > 10 \kp_1$ by \eqref{def:kp},  we obtain
\beq\label{eq:EE_tot4}
\bal
\| \eta(0, z) \|_{\cXc} & \leq \| \rhoc \eta(0, z) \|_{\cXc} 
\leq \tf{1 + \lamcL}{2} \e^{1-2 \kp_1} + C \e^{1 - \kp_1} 
+ C \e^{ 1 -\kp } \e^{1 - 2\kp_1} \\
& \leq 
\tf{1 + \lamcL}{2} \e^{1-2 \kp_1} + C \e^{1 - \kp_1} 
\leq  \tf{2 + \lamcL}{3} \e^{1 -2\kp_1} , 
\eal
\eeq
with $\f{2 + \lamcL}{3}< 1$, for $\e$ small enough. 

\vspace{0.1in}
\paragraph{\bf Estimate of $E_{\cdot, R}(G) $}

To solve the inequality \eqref{eq:EE_tot3}, we define 
\beq\label{def:lln}
  \lln{G} = \max(  \lla{G},  \llr{G}, \e^{ 2 \kp_1}  \llz{G} ).
\eeq

We bound $\cC_{r, R}$ in \eqref{eq:EE_tot1} in terms of $\lln{G}$:
 \bseq\label{eq:EE_tot5}
\beq
\bal
  \cC_{r, R}(\ww, G) 
 & \les     \e \B( 1 + \llz{G}  + \e^{-\kp_1} \lla{G} + \e \llw{G} \B)  
 \\
 & \les \e + \e ( \e^{-2\kp_1} + \e^{-\kp_1} + \e ) \lln{G}
 \les \e + \e^{1-2\kp_1} \lln{G}. 
\eal
\eeq

Applying \eqref{eq:EE_tot4} and \eqref{eq:EE_tot5}, we bound the upper bounds in \eqref{eq:EE_tot3} in terms of $\lln{G}$:
\beq
\bal
  \lla{G} &\les \e^{1 -2 \kp_1}
  + \e + \e^{1 - 2\kp_1} \lln{G} \les \e^{1- 2 \kp_1} + \e^{1 - 2\kp_1} \lln{G} ,\\
  \llr{G} &\les \e + \e^{ 1 - 2\kp_1} \lln{G}, \\
  E_{z, G} &
  \les ( \e^{-\kp_1} + \e  ) \lln{G} + \e^{1 - 5\kp_1} .
\eal
\eeq
\eseq

Since $ 10 \kp_1 < \f12$ by \eqref{def:kp}, taking the weighted maximum, we prove 
\[
\bal
 \lln{G}  & =  \max(\lla{G}, \llr{G}, \e^{2 \kp_1} \llz{G} ) \\
 & \leq \max(  C  \e^{1- 2 \kp_1} , C \e,  
 C_* \e^{1-3 \kp_1}  )
 + C \max(\e^{1 - 2 \kp_1}, \e^{ 1 - 2\kp_1 }, \e^{\kp_1} ) \lln{G}.
\eal
\]

Taking $\e$ small enough so that the coefficient of the last term is less than $\f12$, and then solving the inequality for $\lln{G}$, we prove 
\beq\label{eq:EE_tot6}
  \lln{G} \les  \e^{1- 3 \kp_1}.
\eeq

Plugging \eqref{eq:EE_tot6} into \eqref{eq:EE_tot5} and using $\kp_1 < \f{1}{100}$ by \eqref{def:kp},
we prove 
\beq\label{eq:EE_tot7}
\cC_{r, R}(\ww, G) \les \e,  \quad \lla{G} \les \e^{1-2 \kp_1},
 \quad  \llr{G} \les \e, \quad  \llz{G} \les \e^{1 - 5 \kp_1}.
\eeq

Since the constants in estimates \eqref{eq:EE_tot7} are independent of $R$, taking $R \to \infty$, using the relations 
\eqref{eq:norm_comp2}, $\rhoc \geq 1$ by \eqref{def:rhoc}, and the definition of norms \eqref{def:3d_norm}, we prove 
\[
\bal
 &\nnr{\eta} \leq  \nnrb{\eta} = E_{r, \infty}(G) \les \e, \\
  &  \nna{\eta} + \nnla{\eta} 
\leq \nnab{\eta}  + \nnlb{\eta} \les E_{ 0, \infty}(G) \les \e^{1 -2 \kp_1},  
\quad 
    E_{z,\infty}(G) \les 
  \e^{1 - 5\kp_1} .  \\
\eal
\]
Using \eqref{eq:norm_comp2:b}, the above estimates, and the definition of $\nnn{\cdot}$-norm 
in \eqref{def:3d_norm}, we prove 
\[
\bal
   \nnz{\eta} & \leq \nnzb{\eta} \leq E_{z,\infty}(G) + E_{ 0, \infty}(G)
  \les \e^{1 -5 \kp_1} , \\
     \nnn{\eta}  & \les \max( \nna{\eta},  \nnr{\eta}, \nnz{\eta} ) \les \e^{1 -5\kp_1}. 
\eal
\]

By requiring $\e$ small enough, using the above estimates and \eqref{eq:EE_tot4}, we prove  $\eta \in \cY_{\e}$ defined in \eqref{def:cY} and the following results.

\begin{thm}\label{thm:into}

Let $\cY_{\e}$ be the functional space defined in \eqref{def:cY}, 
and let $\beps_9$ be as in Lemma \ref{lem:3D_damp}. There exists $\beps_{10} \in (0, \beps_9]$ small enough 
such that for any $\e \leq \beps_{10}$ and any $\ww \in \cY_{ \e}$, 
we have $\nnn{\ww} \leq \e^{ \hk }$ and the solution $\eta = \FFF(\ww)$ defined via \eqref{eq:fix_pt_2D} satisfies $\eta \in \cY_{ \e}, \eta \in C^{1,\al}(\R_+ \times \R)$ and the following estimates 
\beq\label{eq:thm_into}
\bal
 \| \eta\|_{ \cXc } & \leq \| \rhoc \eta \|_{ \cXc } \leq \tf{2 + \lamcL}{3} \e^{1-2 \kp_1},   \quad  &&  \nnr{\eta} \leq \nnrb{\eta} \leq \bar C \e \leq \tf{1}{2} \e^{1- \kp_1} , \\
     \nnz{\eta} & \leq \nnzb{\eta} \leq C \e^{ 1 - 5 \kp_1 } \leq \tf{1}{2} \e^{1-6\kp_1} , \quad && \nnn{\eta} \leq C \e^{1- 5\kp_1} \ll \e^{\hk} .
  \eal
\eeq
Moreover, the solution $\eta$ satisfies the properties in Proposition \ref{prop:reg_solu}.

\end{thm}

\section{Construction of 3D profiles and their properties}\label{sec:3D_solu}

In this section, we first use Schauder fixed point theorem to construct the 3D profile to 
\eqref{eq:profi_3D0} and then prove Theorem \ref{thm:self-similar} on the $C^{1,\al}$ self-similar profiles and its properties in Section \ref{sec:proof_self_similar}.

Recall the map $\FFF$ from \eqref{eq:fix_pt_2D} and the ball $\cY_{\e}$ from \eqref{def:cY}. Clearly, $\cY_{\e}$ is a non-empty, closed, and convex set. 
In Theorem \ref{thm:into}, we have proved $\cF_{\R^3} : \cY_{\e} \to \cY_{\e}$. To apply the Schauder fixed point theorem, we further show that $\cF_{\R^3}$ is compact and continuous.

\begin{thm}\label{thm:3D_comp}

 The map $\cF_{\R^3}: \cY_{\e} \to \cY_{\e}$ is compact with respect to the $\cY$-norm defined in \eqref{def:cY}: 
for any sequence $\{ \ww_i\}_{i\geq 1} \in \cY_{\e}$, there is a subsequence of $ \eta_i =\cF_{\R_3}(\ww_i)$ that 
converges to some $\eta_{\infty} \in \cY_{\e}$ in $\cY$-norm.

\end{thm}

\begin{thm}\label{thm:3D_contin}
The map $\cF_{\R^3}: \cY_{\e} \to \cY_{\e}$ is continuous with respect to the $\cY$-norm.
\end{thm}

Proving the continuity of \(\cF_{\mathbb{R}^3}\) directly, or even by relying on the estimates in Section~\ref{sec:3D_EE}, is highly nontrivial, since the trajectory \(X(s, \xx)\) associated with \(\QQ(\ww)\) in \eqref{def:Q} depends on \(\ww\) in a strongly nonlinear manner. Instead, we will prove the compactness of $\cF_{\R^3}$ and the closedness of $\cF_{\R^3}$, and then use these properties to deduce continuity. We defer the proof of Theorem \ref{thm:3D_comp} to Section \ref{sec:3D_comp},
and Theorem \ref{thm:3D_contin} to Section \ref{sec:3D_contin}.  With the above two theorems, we 
construct the profile.

\begin{thm}\label{thm:3D_solu}
Let $\beps_{10}$ be as in Theorem \ref{thm:into}. For any $\e \leq \beps_{10}$, there exists a fixed point $\ww_{\al} = \cF_{\R^3}(\ww_{\al})$ such that $\ww_{\al} \in \cY_{\e}$
and $\ww_{\al}$ satisfies the properties in Theorem \ref{thm:into} and 
\beq\label{eq:thm_3D_est}
\bal
 \|  \ww_{\al} \|_{\cXc} & \les \e^{1-2\kp_1}, 
 \   \nnr{ \ww_{\al}} \les \e , \quad 
     \nnz{\ww_{\al}}   \les \e^{1 -5\kp_1}, 
    \  \nnn{\ww_{\al}} \leq C \e^{1 - 5 \kp_1} \ll \e^{\hk}.
\eal 
\eeq
Moreover, $\ww_{\al}$ is odd in $z$ and it satisfies the properties in Proposition \ref{prop:reg_solu}.

\end{thm}

\begin{proof}[Proof of Theorem \ref{thm:3D_solu}]

By definition \eqref{def:cY}, $\cY_{\e}$ is a non-empty, closed, and convex set in a 
Banach space. Theorems \ref{thm:into}, \ref{thm:3D_comp}, \ref{thm:3D_contin} establish that
 $\cF_{\R^3}: \cY_{\e} \to \cY_{\e}$ is continuous and compact. Using Schauder fixed point theorem, 
 we obtain a fixed point  $\ww_{\al} = \cF_{\R^3}(\ww_{\al})$ with $\ww_{\al} \in \cY_{\e}$. The estimates \eqref{eq:thm_3D_est} for $\ww_{\al}$ follows from 
\eqref{eq:thm_into} in Theorem \ref{thm:into} with $\eta =\ww_{\al}$. 
We prove Theorem \ref{thm:3D_solu}.
\end{proof}

\begin{comment}
Using Theorem \ref{thm:into}, we obtain that $\ww_{\al}$ satisfies properties in Theorem \ref{thm:into} 
and \eqref{eq:thm_into:impr} with $\eta =\ww_{\al}$:
\[
       \| \rhoc \ww_{\al}(0, \cdot) \|_{\cXc}  \leq \tf{1 + \lamcL}{2} || \ww_{\al} ||_{\cXc}  + C \e , 
       \quad 
      | \wwb^{-1} \rhoc \rhoo \rag \ww_{\al} |  \les          \| \rhoc \ww_{\al}(0, \cdot) \|_{\cXc} + \e .
\]

Since $\rhoc \geq 1$ by \eqref{def:rhoc} and $\tf{1 + \lamcL}{2} \in (0, 1)$  by Theorem \ref{thm:into}, from the first estimate, we prove 
\[
   \| \rhoc \ww_{\al}(0, \cdot) \|_{\cXc}  \leq \tf{1 + \lamcL}{2} || \rhoc \ww_{\al} ||_{\cXc}  + C \e \
   \Longrightarrow  \  \| \rhoc \ww_{\al}(0, \cdot) \|_{\cXc} \les \e.
\]

Combining the above two estimates, we prove $ |\wwb^{-1} \rhoc \rhoo \rag \ww_{\al} |  \les \e$.
Since $ |z|^{-1} + \ang z^{\hal} \les \wwb^{-1} $ by \eqref{eq:wa_upper_lower}, 
we prove the second estimate in \eqref{eq:thm_3D_est}. 
The last three estimates have been proved in Theorem \ref{thm:into}. We conclude the proof.
\end{comment}

\subsection{Compactness of $\cF_{\R^3}$}\label{sec:3D_comp}

Suppose that $\ww_i \in \cY_{\e}$ and we denote $\eta_i = \cF_{\R^3}(\ww_i)$ for $i\geq 1$. 
From Theorem \ref{thm:into}, we obtain $ \nnn{\ww_i} \leq \e^{ \hk}$.  Using Proposition \ref{prop:reg_solu} and estimates \eqref{eq:local_reg}, we obtain
\beq\label{eq:3D_comp:pf1}
\| \eta_i \|_{C^{1,\al}(\Bar B_n)} \les_{n, \e} 1 \, ,
\eeq
for any $i ,n \geq 0$. Using the Ascoli-Arzela theorem and a diagonal argument, we obtain a subsequence of 
$\{ \eta_i \}_{ i \geq 1 }$ that converges to $\bar \eta$ uniformly in $\bar B_m$ for any $m \geq 1$. 
To simplify the notation, we may denote the subsequence as $ \eta_i$. 
From \eqref{eq:3D_comp:pf1} and interpolation, for any $n\geq 1$, we obtain that 
\beq\label{eq:3D_comp:pf2}
\bal
\lim_{i\to \infty} \| \eta_i - \bar \eta \|_{C^{1,\b}(\bar B_n)} \to 0,  
& \quad  \forall \ \b \ \in (0, \al),   
 \eal
\eeq
and $\bar \eta \in C^{1,\al}(\bar B_n)$ for any $n\geq 1$. 
Using estimates \eqref{eq:local_reg}, and the above convergence, we obtain 
$ \bar \eta(0) = \na \bar \eta(0) = 0$ and
\beq\label{eq:3D_comp:vanish}
   \limsup_{ i \to \infty }  \max_{ |\xx| \leq 1 } \B( \f{ |\eta_i - \bar \eta| }{|z| \cdot |\xx|^{\b}}
   + \f{ |\na (\eta_i - \bar \eta) | }{ |\xx|^{\b}} \B) 
\les \lim_{i\to \infty} \| \eta_i - \bar \eta \|_{C^{1,\b}(\bar B_2)}
   = 0, \quad  \forall \ \b \ \in (0, \al) .
\eeq

Moreover, from the $C^{1,\b}$ convergence and estimates \eqref{eq:thm_into} on $\eta_i$, we derive 
\beq\label{eq:3D_comp:pf3}
   \nnzb{\bar \eta} \leq \f12 \e^{1 - 6 \kp_1}, \quad \nnrb{\bar \eta} %\leq \bar C \e 
   \leq \f{1}{2} \e^{1-\kp_1},
   \quad  \| \rhoc \bar \eta \|_{\cXc} \leq \f{2 + \lamcL }{3} \cdot \e^{1 - 2 \kp_1}. 
\eeq

Next, we prove that $\eta_i$ converges to $\bar \eta$ in $\cY$ norm \eqref{def:cY}.

Recall the cutoff function $\chi_1$ from \eqref{def:chi1}. 
We define a smooth cutoff function $\chi_{\d}$:
\[
   \chi_{\d}(\xx)  =\chi_1( \d \xx )  , 
   \quad \xx = (r, z) ,
   \quad \d < 1/2.
\]

By definition, we obtain 
\[
  \chi_{\d}(\xx) = 0,\quad \forall   \
  |\xx| \geq 2 \, \d^{-1}, 
  \quad \chi_{\d}(\xx) = 1,
  \quad 
  \forall \  |\xx| \leq \d^{-1} .
\]

Recall the definitions of norms and weights in \eqref{def:3d_norm}. Since the weights 
$\rhoo,\rhoc, \rag, \ang z, \ang \xx$ in the norms 
and those in $\cXc$ \eqref{norm:Xc} are weaker than $|\xx|^{-\bbb}$ 
with $\bbb < 1 + \al$ in any compact domain,
 using the convergence \eqref{eq:3D_comp:pf2} and vanishing near $0$ \eqref{eq:3D_comp:vanish}, 
for a fixed $\d < 1$, we obtain
\[
\bal
  & \limsup_{i \to \infty } \| \chi_{\d} (\eta_i - \bar \eta) \|_{\cXc}
  + \| \chi_{\d} \rhoc \rhoo \rag \ang z^{\hal} \ang \xx \pa_r (\eta_i - \bar \eta) \|_{ L^{\infty} }
  +  \| \chi_{\d} \rhoc \rhoo \rag \ang z^{\hal + 1} \pa_z (\eta_i - \bar \eta) \|_{ L^{\infty} }  \\
  & \qquad \les_{\d}   \limsup_{i \to \infty } \| \eta_i - \bar \eta \|_{C^{1,\al}(\bar B_{ 4 / \d })}  = 0.
\eal
\]

For $\xx$ away from the support of $\chi_{\d}$:  $|\xx| \gtr \d^{-1}$, using the definition of $\rhoc$ \eqref{def:rhoc}, we obtain $|\rhoc(\xx)|^{-1} \les \d^{ \e^2 }$.  Therefore, using \eqref{eq:3D_comp:pf3} and \eqref{eq:thm_into}, we prove 
\[
\bal
& \| (1-\chi_{\d}) (\eta_i - \bar \eta) \|_{\cXc}
  + \| (1-  \chi_{\d} )  \rhoo \rag  \ang z^{\hal} \ang \xx \pa_r (\eta_i - \bar \eta) \|_{ L^{\infty} }
  +  \| (1-\chi_{\d}) \rhoo \rag \ang z^{\hal+1} \pa_z (\eta_i - \bar \eta) \|_{ L^{\infty} }  \\
\leq & 
 \d^{\e^2} \B( \| \rhoc (\eta_i - \bar \eta) \|_{\cXc}
  + \| \rhoc \rhoo \rag \ang z^{\hal} \ang \xx \pa_r (\eta_i - \bar \eta) \|_{ L^{\infty} }
  +  \| \rhoc \rhoo \rag \ang z^{\hal+1} \pa_z (\eta_i - \bar \eta) \|_{ L^{\infty} } \B)
\eal
\]

Using the above estimates, convergence, the definition of $\cY$-norm \eqref{def:cY}, 
and taking $\d \to 0$, we prove  $\lim_{i\to 0} \| \bar \eta - \eta_i \|_{\cY} = 0.$
Thus, we prove Theorem \ref{thm:3D_comp}.

\subsection{Continuity of $\cF_{\R^3}$}\label{sec:3D_contin}

To prove continuity, we first show that $\cF_{\R^3}$ is closed.  

\begin{lem}\label{lem:closed}

Suppose that $\ww_i \in \cY_{\e}$ and $\eta_i = \cF_{\R^3}(\ww_i)$. 
If $\ww_i \to  \ww_{\infty}$ and $\eta_i \to \eta_{\infty}$ in the $\cY$-norm \eqref{def:cY} as $i\to \infty$, then we have $ \eta_{\infty} = \cF_{\R^3}(\ww_{\infty} )$.
\end{lem}

\begin{proof}
Clearly, the convergence implies $\eta_{\infty}, \ww_{\infty} \in \cY_{\e}$. Using the definition of $\cF_{\R^3}$ \eqref{eq:fix_pt_2D}, we obtain
\beq\label{eq:closed_eqn}
  (\bar \QQ + \td \QQ(\ww_{i} ) ) \cdot \na \eta_{ i } = ( \bcB +  \tcB(\ww_{i}) ) \eta_{i}
  + (\bar \cR + \Lpsi(\ww_{i})) \wwb .
\eeq

Since $\eta_i = \cF_{\R^3}(\ww_i)$ and $ \nnn{\ww} \leq \e^{\hk}$, $\eta_i$ satisfies 
the properties in Proposition \ref{prop:reg_solu} and estimates \eqref{eq:3D_comp:pf1}. 
From Lemma \ref{lem:linf_R3}, definition \eqref{def:3d_norm},  and \eqref{norm:comp_Y}, for any $\xx \in \bar B_R$,  we obtain
\beq\label{eq:closed:pf1}
\bal
  | \eta_i -  \eta_{\infty} | & \les \min(|\xx|, 1) (\rhoo \rag )^{-1} 
  \ang z^{-\hal}  \nnn{ \eta_i -  \eta_{\infty}}
  \les_R  |\xx|^{\bbb} 
  \| \eta_i - \eta_{\infty} \|_{\cY} .
\eal
\eeq

Thus, $\eta_i -  \eta_{\infty}$ converges to $0$ uniformly in $\bar B_n$ for any $n \geq 1$, which along with the boundedness \eqref{eq:3D_comp:pf1} imply the $C^{1,\b}$-convergence \eqref{eq:3D_comp:pf2}
\beq\label{eq:closed_reg}
   \lim_{i\to \infty} \|  \eta_i -  \eta_{\infty} \|_{C^{1,\b}( \bar B_R )} =0,
   \quad \forall \b \in (0, \al), 
   \quad  \eta_{\infty}  \in C^{1, \al}(\bar B_R), \quad  \forall R \geq 1.
\eeq

From definition \eqref{def:qq} and \eqref{eq:lin_2Db}, $\td \QQ(\ww), \tcB(\ww), \Lpsi(\ww)$ depends on $\ww$ linearly. Using the definition, the linearity, \eqref{norm:comp_Y}, and Proposition \ref{prop:vel_est}, for any $ \xx \in \bar B_R$, we have 
\[
\bal
 &  | \td \QQ(\ww_i) - \td \QQ( \ww_{\infty})|
  + \| \tcB(\ww_i) - \tcB( \ww_{\infty}) | 
  + | \wwb \Lpsi(\ww_i -  \ww_{\infty}) |  \\
& \quad  \les_{R, \e}  (1 + |\xx|) ( 1 + |\wwb|  )  ( | \psi_z(\ww_i- \ww_{\infty} )(0) | + 
|\psi(\ww_i-\ww_{\infty})| + | \na \psi(\ww_i- \ww_{\infty})|  ) \\
& \quad \les_{R, \e}  \nnn{ \ww_i - \ww_{\infty}} 
\les_{R, \e}  \|\ww_i -  \ww_{\infty} \|_{\cY}.
\eal
\]

Moreover, for any $R>0$, Lemma \ref{lem:psi_reg} implies that the coefficients $\td \QQ(\ww_i), \tcB(\ww_i), \Lpsi(\ww_i)$ are uniformly bounded in $\bar B_R$.  Thus, for any fixed $R > 1$, taking $i\to \infty$ in \eqref{eq:closed_eqn} and using the above uniform convergence in $\bar B_R$, we prove 
\[
   ( \bar \QQ +  \QQ(\ww_{\infty} ) ) \cdot \na \eta_{\infty} = ( \bcB + \tcB(\ww_{\infty} ) ) \eta_{\infty}
  + (\bar \cR + \Lpsi(\ww_{\infty} )) \wwb .
\]
for any $\xx$. Thus, we obtain that $ \eta_{\infty}$ solves the fixed point equation \eqref{eq:fix_pt_2D} 
with $C^{1,\al}$-regularity \eqref{eq:closed_reg}. Estimate \eqref{eq:closed:pf1}
and \eqref{eq:local_reg} on $\eta_i$ implies that $ |\eta_{\infty}(\xx) | \les_{\e,  \e} |\xx|^{\bbb}$ 
for $|\xx|\leq 1$ with $\bbb > 1$.

Since $\ww_{\infty} \in \cY_{\e}$ satisfies condition \eqref{eq:3D_size:ww}
and $\eta_{\infty}$ also satisfies \eqref{eq:3D_ODE:init}, by uniqueness in Proposition \ref{prop:reg_solu}, $\eta_{\infty}$ is the unique solution to \eqref{eq:fix_pt_2D}. We conclude $\eta_{\infty} = \FFF(\ww_{\infty})$.
\end{proof}

\vspace{0.1in}

Now, we are in a position to prove Theorem \ref{thm:3D_contin} on the continuity of $\FFF$.
\begin{proof}[Proof of  Theorem \ref{thm:3D_contin}]

Suppose that $ \ww_i,  \ww_{\infty} \in \cY_{\e}$ and $ \| \ww_i - \ww_{\infty} \|_{\cY} \to 0$ as $i \to \infty$.
Denote $\eta_i = \cF_{\R^3}(\ww_i)$. % \eta_{\infty} = \cF_{\R^3}(  \ww_{\infty} )$.
We shall prove $\eta_i \to \cF_{\R^3}(  \ww_{\infty} ) $ in $\cY$-norm. 

We argue by contradiction. If not, there exists $\d > 0$ and a subsequence $\eta_{n_i}$ 
of $\{ \eta_i \}_{i\geq 1}$ such that  $ \| \eta_{n_i} -  \eta_{\infty} \|_{\cY} > \d$ for any $i \geq 1$. As $\eta_{n_i} = \cF_{\R^3}(\ww_{n_i}), \ww_{n_i} \in \cY_{\e}$, using Theorem \ref{thm:3D_comp} on the compactness, we further extract a subsequence $\{  \eta_{ m_i } \}_{i\geq 1} \subset \{  \eta_{ n_i } \}_{i\geq 1}$, such that $\| \eta_{m_i} - \eta^{\pr} \|_{\cY} = 0$ as $i \to \infty$ for some $\eta^{\pr} \in \cY_{\e}$. Since $ \| \ww_{m_i} - \ww_{\infty} \|_{\cY} \to 0$ as $i\to \infty$,  using Lemma \ref{lem:closed} on the closedness, we obtain $\eta^{\pr} = \cF_{\R^3}(\ww_{\infty})$. 
As a result, we obtain $\| \eta_{m_i} -  \cF_{\R^3}(\ww_{\infty}) \|_{\cY}   = \| \eta_{m_i} - \eta^{\pr} \|_{\cY} \to 0$
as $i\to \infty$, which contradicts   $\| \eta_{m_i} - \eta_{\infty} \|_{\cY} > \d$ for any $i$. We conclude the proof.
\end{proof}

\subsection{Properties of the profile}\label{sec:3D_profile_prop}

Recall from  \eqref{eq:normal_bar}, \eqref{def:3d_Om},
and Section \ref{sec:appr_profi} the  approximate profile 
\[
\wwb(r, z) = \waa(z), \quad \bpsi = \BS( \wwb)  ,
\quad   \bcl 
   = 2 - 4 \al \cJa(\waa),  \quad \bcw =   2 + 2 \al (1-\al) \cJa(\waa).
\]
 From the definition of $\ww_{\al}$ and the fixed point map $ \cF_{\R^3}$, we obtain the profile 
to \eqref{eq:profi_3D0}
\bseq\label{eq:solu_al}
\beq
\bga
     \QQs \cdot \na \ws =  \cBs \cdot \ws, 
    \quad \ws \teq \ww_{\al} + \wwb, 
    \quad \psis \teq \BS(\ws) = \BS(\ww_{\al}) + \bpsi , \\
  \cls = \bcl + \td c_l , \quad \cws = \bcw + \td c_{\om}. 
\ega
\eeq
where we define $\QQs$ and $\cBs$ similar to $\bar \QQ$ \eqref{def:Q} and $ \bcB$ \eqref{eq:lin_2Db}:
\beq
\bal
        \QQs & \teq  ( \cls r - r \pa_z \psis, \  \cls z + 2 \psis + r \pa_r \psis  ),
   & \quad  \cBs &= \cws - (1 -\al) \pa_z \psis,  \\ 
      \cls & = 2 - 2 \pa_z \psis(0), & \quad  \cws & = 2 + (1 - \al) \pa_z \psis(0).
\eal
\eeq
\eseq

Moreover, from Theorem \ref{thm:3D_solu} and Theorem 
\ref{thm:1D_profile_prop}, since $\ww_{\al}$ satisfies Proposition \ref{prop:reg_solu} 
and $\pa_z \wwb(0) = \pa_z \waa(0) = \pa_z \wwwa(0)$, we obtain $\na \ww_{\al}(0) = 0$ and $\pa_z \ws(0) = \pa_z \wwb(0) = \pa_z \wwwa(0)$. 

\begin{remark}%[\bf $\ws \neq \waa$]

We denote the 3D profile for \eqref{eq:profi_3D0} by $\ws$ and the 1D profile for \eqref{eq:1D_dyn} by $\waa$.

\end{remark}

We establish several profile properties  and use them to construct blowup solution in Section \ref{sec:blowup}.

\begin{prop}\label{prop:3D_profile}

Let $\beps_{10}, \epa$ be as in Theorem \ref{thm:3D_solu} and  \eqref{ran:ep2}, respectively.  There exists $\beps_{ 11} \in (0, \beps_{10}]$ small enough, such that for any $\e = \f13 - \al \leq \beps_{ 11}$, the following results hold. 
\begin{enumerate}[label=(\roman*),leftmargin=2em]
\item \textbf{Scaling}. The scaling parameters satisfy
\bseq\label{eq:3D_solu_scale}
\begin{align}
 \cls
 & =   \tf{64}{9} \e^{-1} + O(\e^{- \hk})
  \asymp \e^{-1},  & \quad  & \cws = - \tf{64}{27} \e^{-1} + O(\e^{- \hk})   \asymp - \e^{-1},  \label{eq:3D_solu_scale:a} \\
  \als & \teq -\f{\cws}{\cls} = \f13 + \f{\e}{8} + O(\e^{1 + \mhk } ),
 &  \quad  &
 \max( \hal - \epa , \f13 + \f{\e}{9}) \leq \als \leq \f13 + \f{\e}{7}.
    \label{eq:3D_solu_scale:b}
\end{align}
\eseq

\item \textbf{Regularity and decay}.
The profile satisfies
\beq\label{eq:3D_solu_reg}
\bal
\| \ws \|_{ C^{1,\al}(\bar B_R) } + 
\| \QQs \|_{ C^{1,\al}(\bar B_R) } + \| \pa_z \psis \|_{ C^{1,\al}(\bar B_R) }  & \les_{R, \e} 1 , \\
\|  r \pa_r \ws \|_{ C^{1,\f18}(\bar B_R) }
+  \|  \pa_z \ws \|_{ C^{1,\f18}(\bar B_R) } & \les_{R, \e} 1. 
\eal
\eeq
as a function in $(r, z) \in \R_+ \times \R$ for any $R > 0$,  and 
\beq\label{eq:3D_solu_Q0}
 \pa_r \QQQs^r(0) \geq C \e^{-1} \geq 4, \quad  \pa_z \QQQs^z(0) = 2 .
\eeq
We have the decay and vanishing estimates 
\bseq\label{eq:3D_solu_decay}
\begin{align}
   |\ws| & \leq \upse \min( |z|,    \ang \xx^{-\als} ),  
   \quad  |\pa_r \ws| \leq \upse  |z|^{\al} ,  
\label{eq:3D_solu_decay:a} \\ 
  | \na \ws | + |r \pa_r \na \ws|  &  \leq \upse  \ang \xx^{-\als - 1}, 
  \quad  |\na \pa_z \ws|  \leq \upse \ang \xx^{-\als - 2}, 
  \label{eq:3D_solu_decay:b} \\
  |\tf{1}{z} \psis | +  |\na \psis | & \les \e^{-1} \ang \xx^{\al -\hal + \epa}, 
\label{eq:3D_solu_decay:c} \\
  |\tf{1}{z} \psis - \pa_z \psis(0) | +  | \pa_r \psis |
  + | \pa_z \psis - \pa_z \psis(0)| & \leq \upse |\xx|^{\al + 1}, 
\label{eq:3D_solu_decay:c2} \\
  r \ang \xx^{-1} |\pa_{rr} \psis| + |\pa_{r z} \psis|
  + |\pa_{ zz } \psis| & \leq \upse \min( |\xx|^{\al}, \ang  \xx^{\al-\hal + \epa - 1} ) . 
  \label{eq:3D_solu_decay:d}
\end{align}
\eseq
for some  constant $\upse > \e^{-1}$ depending on $\e$.

\item \textbf{Outgoing properties}. Let $\lam_Q$ be the constant in \eqref{eq:traj_outgo}.
For any $\xx$, we have 
\beq\label{eq:outgo}
  \bar \QQ_{\al}(\xx) \cdot \xx  \geq \lam_Q |\xx|^2 . 
\eeq

\end{enumerate}

\end{prop}

The main difficulty in proving Proposition \ref{prop:3D_profile} is the decay estimates 
\eqref{eq:3D_solu_decay:a}, \eqref{eq:3D_solu_decay:b}, which do not follow from the estimates 
in Theorem \ref{thm:3D_solu}. We exploit the crucial property that $\ws$ enjoys better estimates in $z$-direction than in $r$-direction. Below, we first derive basic regularity estimates from Theorem \ref{thm:3D_solu}, and then perform decay estimates by iterating the elliptic estimates for $\psis$ and trajectory estimates for $\ws$. 
The proof of \eqref{eq:3D_solu_decay} is \emph{qualitative}, in the sense that
the implicit constants may depend on \(\e\) and may be very large.

\vs{0.05in}
\paragraph{\bf Basic regularity estimates}
Using Lemma \ref{lem:psi_reg}
with $\ww = \ws$ and Proposition \ref{prop:reg_solu} with $\eta = \ww = \ws$, 
for any $R > 0$,  we obtain
\beq\label{eq:solu_basic_reg}
\| \ws \|_{ C^{1,\al}(\bar B_R) } + 
\| \QQs \|_{ C^{1,\al}(\bar B_R) } + \| \pa_z \psis \|_{ C^{1,\al}(\bar B_R) } \les_{R, \e} 1 ,
 \ \Longrightarrow \ \cBs , \ \psis 
\in C^{1,\al}(\bar B_R), 
\eeq
The regularity 
$\psis \in C^{1,\al}$ follows from  $\psis(r, 0) = 0, \psis(r, z) = \int_0^z \pa_z \psis(r, s) d s$, and $\pa_z \psis \in C^{1,\al}$.  

Since $\ws = \wwb + \ww_{\al}$ and $\ww_{\al}$ satisfies estimates 
\eqref{eq:thm_3D_est}, using \eqref{eq:3D_boot_q:1} 
with $\eta = \ww_{\al}$, we obtain
\beq\label{eq:qqs_low}
   q^r    \teq \tf{1}{r} \QQQs^r =  \cls - \pa_z \psis \gtr \e^{-1} ,
   \quad q^z(0, z) \teq \tf{1}{z} \QQQs^z(0, z) 
  = \cls + 2 \tf{1}{z} \psis \gtr 1 .
\eeq

Using $q^z(0) = 2$ by \eqref{eq:normal_cond} and the above estimates with $\xx=0$, we prove 
\eqref{eq:3D_solu_Q0}.

\subsubsection{Estimate of scaling parameters $\cls , \cws$}\label{sec:3D_solu:scale}

We decompose $\cls = \bcl + \td c_l , \cws = \bcw + \td c_{\om}$. 
Since $ \ws =\wwb + \ww_{\al} 
= \waa + \ww_{\al}$, using \eqref{eq:normal_bar}, \eqref{eq:normal_cond}, 
and the relation \eqref{eq:psi_JJ}, we obtain
\[
\bal
   & \bcl =  2 - 4 \al \cJa(\waa)(\infty),  \quad && \bcw = 2 + 2 \al (1-\al) \cJa(\waa)( \infty), \\
   & \td c_l =   - 2 \pa_z \psi(\ww_{\al})(0) = - 4 \al \JJ(\ww_{\al})(\infty),
   \quad && \td c_{\om} = (1-\al)  \pa_z \psi(\ww_{\al})(0) = (1-\al) 2 \al  \JJ(\ww_{\al})(\infty).
\eal
\]

Using estimates \eqref{eq:vel_J},  $\nnn{\ww_{\al}} \les \e^{1-5\kp_1}$ 
from \eqref{eq:thm_3D_est}, and $\cJaa \les \e^{-1}$ from \eqref{eq:Ja_hat}, we obtain
\[
|\tcl | + |\tcw| +
 | \JJ(\ww_{\al})(\infty)| \les \cJaa^{1 + \kp} \nnn{\ww_{\al}}
 \les \e^{-1-\kp  + 1 - 5\kp_1}
 = \e^{-\kp - 5\kp_1}.
\]

Since $\wwb= \waa$ is constant in $r$, using the estimate \eqref{eq:comp_vel_bc} and the estimates \eqref{eq:wa_decay_1D}, we obtain
\[
 \JJ(\wwb)(\infty) =\JJ(\waa)(\infty)=  \cJa(\waa)(\infty)
 =  - \tf{16}{3} \e^{-1} + O(\e^{-\kp}),
 \quad \bar c_l = \tf{64}{9} \e^{-1} + O(\e^{-\kp}).
\]

Since $\kp , \kp + 5 \kp_1 < \hk < 1$ by \eqref{def:kp},  combining the above two estimates,
using \eqref{eq:wa_decay_1D}, and by requiring $\e$ small enough, we obtain
\[
\JJ(\ws)(\infty) = - \tf{16}{3} \e^{-1} + O(\e^{- \hk }) ,
\quad \cls =  \tf{64}{9} \e^{-1} + O(\e^{- \hk}) \asymp \e^{-1}, 
\quad  -\cws \asymp \e^{-1}  .
\]

We prove \eqref{eq:3D_solu_scale:a}. Since $\kp , \kp + 5 \kp_1 < 1$ and $\e = \f13 - \al$, using the above estimates, we obtain
\[
\bal
3 \cws  + \cls  & =  3 \bar c_l + \bar c_{\om}  
+ 2 \al \JJ( \ww_{\al})(\infty) ( 3 (1-\al) - 2) 
= 3 \bar c_l + \bar c_{\om}   
+ 2 \al \cdot 3 \e \JJ( \ww_{\al} )(\infty) \\
& = 3 \bar c_l + \bar c_{\om}  + O( \e^{1- \hk}), \\
\cls & = \bar c_l + \td c_l 
= \bar c_l ( 1 + O( \e^{1 -\hk}) .\\
\eal 
\]

Thus, using estimate \eqref{eq:wa_decay_1D} and $\bar c_l \gtr \e^{-1} > 0$, 
for $\e$ small enough, we bound 
\[
 \f{3 \cws + \cls}{ 3 \cls}
 = \f{ 3 \bar c_l + \bar c_{\om}  + O( \e^{1-\hk}) }{ \bar c_l ( 1 + O( \e^{1 - \hk} ) ) } 
= O( \e^{2-\hk }) + \f{ 3 \bar c_l + \bar c_{\om} }{ \bar c_l} 
( 1 + O( \e^{1 - \hk } ) 
= - \f18 \e + O( \e^{1 + \mhk}).
\]

Recall $\epa,\hal$ from \eqref{ran:ep_all}. Now that $\als = - \cws / \cls$. Rewriting the above estimates and by requiring $\e$ small enough, we prove \eqref{eq:3D_solu_scale:b}:
\[
\als = \tf13 + \tf 18 \e + O( \e^{1 + \mhk}) 
\in [\tf 13 + \tf19 \e , \, \tf13 + \tf17 \e] , \quad  \als > \hal - \epa =  \tf13 + \tf 18 \e  - \epa+ O( \e^{1 + \mhk}) .
\]

\subsubsection{Regularity and decay of $\ws(0,z)$}

Restricting \eqref{eq:solu_al} on $r = 0$ and taking $\pa_z$,  we obtain
\beq\label{eq:3D_decay:pf0}
  \QQQs^z \pa_z \ws = \cBs \ws .
\eeq

Integrating from $1$ to $z > 1$, we obtain
\[
  \ws(0, z) = \ws(0, 1) \exp\big( \int_1^z \f{ \cBs }{ \QQQs^z }( y ) d y \big),
  \quad \ws(0, 1) \neq 0.
\]

Note that $\ws = \wwb + \ww_{\al}$ with $\nnn{\ww_{\al}} \leq \e^{ \hk}$.  Using Proposition
\ref{prop:vel_est} and \eqref{cor:u_bar}, $\Psi_z(0) =  2 \al \JJ(\infty)$,
and $\rhoo^{-1} \les_{\e} \ang \xx^{\epa}$ by \eqref{def:3d_wg}, we obtain
\beq\label{eq:3D_decay:bd:psi}
\bal
  |\pa_z \psis(\xx)| + \tf{1}{z} |\psis(\xx)| & \les | \pa_z \psis -  \pa_z \psis(0) + 2 \al \JJ(\ws) |
  + | \tf{1}{z} \psis - \pa_z \psis(0) + 2 \al \JJ(\ws) | \\
  & \quad + | 2 \al ( \JJ(\ws)(\xx) - \JJ(\ws)(\infty) )| \les_{\e} \ang \xx^{ \al -\hal + \epa }, \\
  |\pa_{zz} \psis(0, z) | &   \les_{\e}  ( \rhoo^{-1}(0, z ) + 1) \ang z^{\al-1-\hal} 
  \les_{\e} \ang z^{\al -\hal + \epa - 1}.
\eal
\eeq
Using the above estimates, the formulas of $\QQQs^z(0, z), \cBs$ from \eqref{eq:solu_al}, 
and \eqref{eq:qqs_low}, we obtain
\beq\label{eq:3D_decay:pf1}
  |\QQQs^z(0, z) - \cls z| \les_{\e} |z|^{\al -\hal + \epa + 1},
  \quad  | \cBs - \cws | \les_{\e}  |z|^{\al-\hal + \epa},
  \quad |\QQQs^z(0, z)| \gtr_{\e} |z| .
\eeq

Since $\al -\hal + \epa < 0$ by \eqref{ran:ep_all}, using $\QQQs^z(0,z) \gtr z, \forall z>0$  
by \eqref{eq:qqs_low}, for $z > 1$, we obtain
\[
  \f{\ws(0, z)}{ \ws(0,1) } = 
   \exp\big( \int_1^z \f{ \cBs }{ \QQQs^z }( y ) d y \big)
  = \exp\big( \int_1^z \f{\cws}{\cls} y^{-1} + O_{\e}( y^{-1 + \al -\hal + \epa} ) d y \big)
  \asymp_{\e}  z^{ \cws / \cls }.
\]

Using $\als = - \cws / \cls$, $\ws(0, 1) \neq 0$, $\pa_z \ws(0) = \pa_z \wwwa(0) < 0$, $\ws \in C^{1,\al}$ \eqref{eq:solu_basic_reg}, we prove
\beq\label{eq:3D_decay:bd:w1}
  |\ws(0,z)| \asymp_{\e} \min( |z|, \ang z^{ -\als } ),  \quad \als = -   \cws / \cls.
\eeq

 Taking $\pa_z$ on \eqref{eq:3D_decay:pf0},  we obtain
\beq\label{eq:3D_decay:bd:pf1}
    \QQQs^z \pa_{zz} \ws = (\cBs - \pa_z \QQQs^z ) \pa_z \ws + \pa_z \cBs \cdot \ws .  
\eeq

Using \eqref{eq:3D_decay:pf0}, \eqref{eq:3D_decay:pf1}, \eqref{eq:3D_decay:bd:w1}, we prove 
\bseq\label{eq:3D_decay:bd:w2}
\beq\label{eq:3D_decay:bd:w2:a}
  |\pa_z \ws| \les_{\e} |\QQQs^z|^{-1} \cdot | \cBs \ws  | 
  \les |z|^{-1} |\ws| \les \ang z^{ -\als -1 }.
\eeq

Since $ g(z) = \cBs - \pa_z \QQQs^z = \cws - (1-\al) \pa_z \psis - ( \cls + 2 \pa_z \psis)$ 
satisfies $g(0) = 0$ by \eqref{eq:solu_al}, using \eqref{eq:3D_decay:bd:psi}, we bound 
$|g(z)|\les_{\e} \min( |z|, 1)$. 
Thus, using the vanishing of $\ws(0,z), g(z)$ near $z=0$ and 
\eqref{eq:3D_decay:bd:pf1}, we obtain $\ws(0,z) \in C^{2}$. Using 
\eqref{eq:3D_decay:bd:psi}, \eqref{eq:3D_decay:bd:w1}, \eqref{eq:3D_decay:bd:w2:a}, 
and equation \eqref{eq:3D_decay:bd:pf1}, we prove 
\beq
 |\pa_{zz} \ws| \les_{\e} |z|^{-1} ( |g(z) \cdot \pa_z \ws | + |\pa_z \cBs \cdot \ws| )
 \les_{\e} \ang z^{-1} \cdot  \ang z^{-\als - 1} + \ang z^{-2} \cdot \ang z^{-\als} 
 \les_{\e}  \ang z^{-\als - 2} .
\eeq
\eseq

\subsubsection{Regularity of $\pa_z \ws$}

Next, we estimate higher order regularity of $\pa_z \ws$, which enjoys better regularity than $\pa_r\ws $. Taking $ \pa_z$ on \eqref{eq:solu_al}, we obtain
\[
  \QQs \cdot \na \pa_z \ws = ( \cBs - \pa_z \bar Q_{\al}^z ) \cdot \pa_z \ws 
  - \pa_z \QQQs^r \cdot \pa_r \ws + \pa_z \cBs \cdot \ws. 
\]

We aim to derive a transport-type equation of $\pa_z \ws$ with more regular source term. 
To this end, we divide $q^r = \QQQs^r / r$ on both sides to obtain
\[  
\QQs \cdot \na  \f{ \pa_z \ws}{q^r} = \big( \cBs - \pa_z \bar Q_{\al}^z -  \f{\QQs \cdot \na q^r}{q^r} \big) \cdot \f{ \pa_z \ws }{q^r}
  - \f{ \pa_z \QQQs^r \cdot \pa_r \ws}{q^r} +  \f{ \pa_z \cBs \cdot \ws }{q^r}. 
\]

We use the equation \eqref{eq:solu_al} and $\QQQs^r = r q^r$ to eliminate the term $\pa_r \ws$: 
\[
\bal
   - \f{ \pa_z \QQQs^r \cdot \pa_r \ws }{q^r} & =-  \f{ \pa_z \QQQs^r }{ q^r \cdot \QQQs^r} \cdot (\QQs \cdot \na \ws - \QQQs^z \pa_z \ws )
   = -\f{ \pa_z q^r  }{ (q^r)^2 } \cdot ( \cBs \ws - \QQQs^z \pa_z \ws )\\
  & =  \f{ \QQQs^z \pa_z q^r}{q^r} \cdot \f{\pa_z \ws}{q^r}
    -\f{ \pa_z q^r  }{ (q^r)^2 }  \cdot \cBs \ws .
\eal
\]

Combining the above derivations, we obtain
\bseq\label{eq:3D_tran:ws1}
\beq
    \QQs \cdot \na \Ups = a(\xx) \cdot \Ups + f(\xx), \quad 
    \Ups(r, z) \teq  \f{ \pa_z \ws }{q^r} ,
\eeq
  where $a, f$ are given by 
\beq
\bal
a(\xx) & =  \cBs - \pa_z \QQQs^z  +    \f{ \QQQs^z  \pa_z q^r  }{q^r} 
- \f{\QQs \cdot \na q^r}{q^r} \, ,
  \\
f(\xx) &= - \f{ \pa_z q^r  }{ (q^r)^2 } \cdot  \cBs \ws  + \f{ \pa_z \cBs \cdot \ws}{q^r}
= \pa_z (\f{\cBs }{q^r} ) \cdot \ws \, .
\eal
\eeq
Next, we simplify the formula. Using $\QQQs^z = \cls z + 2 \psis + r \pa_r \psis $ 
and $\cBs = \cws - (1-\al) \pa_z \psis$, we obtain
\beq
\bal
  a(\xx) & = \cBs
 -\pa_z \QQQs^z
  - r  \pa_r q^r =  \cws - (1-\al) \pa_z \psis - \cls - 2 \pa_z \psis - r \pa_{r z } \psis 
   + r \pa_{r z} \psis \\
  & = \cws - \cls - (3 -\al) \pa_z \psis . \\
\eal
\eeq

Using the definition \eqref{eq:solu_al}, we have 
 \[ 
\cBs= \cws - (1-\al) \cls + (1-\al) ( \cls - \pa_z \psis )
 = \cws - (1-\al) \cls + (1 -\al) q^r, 
\]
which along with $q^r = \tf{1}{r} \QQQs^r = \cls - \pa_z \psis$ (see \eqref{eq:solu_al}) implies 
\beq
  f(\xx) = ( \cws - (1-\al) \cls) \pa_z \f{1}{q^r} \cdot \ws 
  = \dds \cdot \ws \cdot  \f{\pa_{zz} \psis}{ (q^r)^2} , 
  \quad \dds \teq   ( \cws - (1-\al) \cls).
\eeq
\eseq

From the normalization conditions \eqref{eq:normal_cond} and the odd symmetry in $z$, we obtain 
\beq\label{eq:traj_est_van}
f(0) = 0, \quad
  a(0) = \cws - (1-\al) \pa_z \psis(0) - ( \cls + 2 \pa_z \psis(0) ) = 0, 
  \quad \Ups(0) = 1 / q^r(0) \asymp \e. 
\eeq

\paragraph{\bf Improved regularity}

We use the transport equation \eqref{eq:3D_tran:ws1}, the elliptic estimates in Lemma \ref{lem:schauder_decay} (to be proved), and trajectory estimates in Lemma \ref{lem:traj_reg}. We pick 
\[
\g_i = \tf{i}{8},  \quad \g_{n_{\g}} = \tf98,
\quad \g_{i+1} \leq  \tf{1}{8} + \tf{5}{4} \g_i , \quad  0\leq i\leq n_{\g} - 1 , \quad n_{\g} = 8.
\]

From \eqref{eq:solu_basic_reg},\eqref{eq:qqs_low}, we have $\ws, \QQs, q^r, a \in C^{1,\al}(\bar B_{2R}), f \in C^{\al}(\bar B_{2R})$ for any $R > 0$ and $q^r \gtr \e^{-1}$. 

Suppose that $\pa_z \ws \in C^{\g_i}, i \leq n_{\g} - 1$. 
We have the base case $\pa_z \ws \in C^{\g_0} = C^0$  due to \eqref{eq:solu_basic_reg}. Applying Lemma \ref{lem:schauder_decay}, 
for the above choices of $\g_i$,  we obtain $\pa_{zz} \psis \in C^{ \g_{i+1} }(\bar B_R)$ 
and thus $f(\xx) \in C^{ \g_{i+1} }(\bar B_R)$  for any $R >0$. Here, we abuse $C^{\s} = C^{1,\s-1}$ if $\s \in (1, 2)$. 

Applying Lemma \ref{lem:traj_reg} with case (1), (2) in \eqref{eq:traj_reg:case} to \eqref{eq:3D_tran:ws1}, using the above regularity estimates, and \eqref{eq:traj_est_van}, we obtain 
$\pa_z \ws \in C^{\g_{i+1}}(\bar B_R)$ for any $R > 0$. 
Repeating these steps for $i=0,.. , n_{\g}-1$, 
we prove
\bseq\label{eq:solu_high_reg1}
\beq
  \| \pa_{z} \ws \|_{C^{1, 1/8}(\bar B_R)} +  
    \| \pa_{zz} \psis \|_{C^{1, 1/8}(\bar B_R)} \les_{\e, R} 1 , \ \forall \, R > 0.
\eeq

Next, we estimate $\pa_r \ws$. Taking $\pa_r$ in \eqref{eq:solu_al}, we obtain
\[
\QQQs^r \cdot \pa_r^2 \ws + \QQQs^z \pa_{r z} \ws  = \pa_r ( \cBs \ws ) - \pa_r \QQs \cdot \na \ws.  
\]
Using \eqref{eq:solu_basic_reg}, \eqref{eq:solu_high_reg1}, $ \QQQs^r = r q^r
\gtr \e^{-1} r $ by \eqref{eq:qqs_low} ,
and $q^r \in C^{1,\al}$ by \eqref{eq:solu_basic_reg},  for any $R > 1$, we prove 
\beq
  \| q^r  r \cdot \pa_r^2 \ws \|_{C^{1, 1/6}(\bar B_R)}  
  \les_{\e, R } 1 ,   \ 
 \Longrightarrow \ 
\|  r \cdot \pa_r^2 \ws \|_{C^{1, 1/6}(\bar B_R)} \les_{\e, R } 1.  
\eeq
\eseq

\subsubsection{Elliptic estimates}

Recall the cutoff function $\chi_1$ from \eqref{def:chi1}. It has support $\supp(\chi_1) \subset [-2, 2]$. We decompose $\psis = \Psi_{1} + \Psi_{2}  $ into the parts from the boundary $r=0$, and from the interior :
\bseq\label{eq:elli_decomp}
\begin{align}
  - \D_{\R^5}  \Psi_{1}  &= \cpsia ( \ws(r, z) -\ws(0, z) \chi_1( \tf{8 r^2}{\ang z^2} ) )  r^{\al-1} \teq \cpsia F_1 ,
 && \Rightarrow \, - \D_{\R^5} \pa_z \Psi_1 = \cpsia \pa_z F_1.  \label{eq:elli_decomp:a}
  \\
   - \D_{\R^5}  \Psi_2  & = \cpsia \ws(0, z)  \chi_1( \tf{8 r^2}{\ang z^2} )  r^{\al-1} \teq \cpsia F_2 .
   \label{eq:elli_decomp:b}
\end{align}
\eseq

Using \eqref{eq:vel_r1}-\eqref{eq:vel_psi} and \eqref{eq:vel_J} in Proposition \ref{prop:vel_est}
with $\bbu = 1, \rhoo = \ang \xx^{-\epa} \cJaa^{-\kp}$ and $\rag = (\f{\ang z}{\ang \xx})^{\kag}$, we obtain
\beq\label{eq:elli_dec0}
\bal
   |\na \Psi_i| + |\tf{1}{z} \Psi_i |
   & \les_{\e} 
   \ang \xx^{\al-\hal + \epa} \B(   
    \nlinf{ ( |\xx|^{-1} + 1 ) \ang z^{\hal} \rhoo \rag | \ws(r, z) -\ws(0, z) \chi_1( \tf{8 r^2}{\ang z^2 }  ) | } \\
   & \qquad + \nlinf{ ( |\xx|^{-1} + 1 ) \ang z^{\hal} \rhoo \rag \ws(0, z) 
   \chi_1( \tf{8 r^2}{\ang z^2 }   ) } 
   \B)
   \les_{\e} \ang \xx^{\al-\hal + \epa}.
  \eal
\eeq

For any $a < b$, we define the annulus region
\beq\label{def:annu}
B_{a, b} \teq \{ \yy : |\yy| \in (a, b) \},
\quad \bar B_{a, b}  \teq \{ \yy : | \yy| \in [ a, b ] \}.
\eeq

We have the following estimates of $\Psi_i$. 

\begin{lem}\label{lem:schauder_decay}

Suppose that $\al \in (0.31, \f13)$. Let $\ws$ be the profile in Theorem \ref{thm:3D_solu} 
and $\psis = \BS(\ws)$.

(I) For any $\g \in [0, 1-\al ]$ and $ \g_* = \f{1}{8} + \f{5}{4} \g$,  we have 
\bseq\label{eq:schauder_reg1}
\beq
   \| \pa_{zz} \psis  \|_{ C^{ \g_*}(\bar B_{ R}) } 
   \les_{R, \e, \g} 1  
   + \| \pa_z \ws \|_{ C^{\g}(\bar B_{3 R}) } 
+ \| \pa_z^{\leq 2} \ws(0, z) \|_{ L^{\infty}( [-3 R, \, 3 R] )}  \, .
\eeq
For $\g \in (1-\al, 1)$ and $\gam_1 = \min( \al + \g-1, \tf18)$, we have 
\beq
     \| \pa_{zz} \psis  \|_{ C^{ 1, \gam_* }(\bar B_{ R}) } 
   \les_{R, \e, \g} 1 
    + \| \pa_z \ws \|_{ C^{\g}(\bar B_{3 R}) } 
+ \| \pa_z^{\leq 2} \ws(0, z) \|_{ L^{\infty}( [- 3 R, \, 3 R] )} .
\eeq
\eseq

(II)  Suppose that for some $\g \in [0, 1)$, we have 
\beq\label{eq:ass:schauder_decay}
  \|  \ws \|_{ L^{\infty}( B_{R, 2 R}) } \leq C_{\e} \ang R^{ - \als },  
  \quad 
  \| \pa_z \ws \|_{ L^{\infty}( B_{R, 2 R}) } \leq C_{\e} \ang R^{-1 - \als }, 
  \quad 
  \| \pa_z \ws \|_{ \dot C^{\g}( B_{R, 2 R}) } \leq C_{\e} \ang R^{-1 - \als - \g},
\eeq
for  any $R>0$, 
where $C_{\e}$ is some constant independent of $R$ and we write $C^{0} = L^{\infty}$.  

For $\g \leq 1-\al$ and $\g_* = \tf18 + \tf54 \g$, we have 
\bseq\label{eq:schauder_decay}
\begin{align}
R \| \pa_{zz} \psis \|_{ L^{\infty}{ (\bar B_{R, 2R} ) }}
+  R^{1 + \g_*} \|   \pa_{zz} \psis  \|_{ \dot C^{  \g_* }(  \bar B_{R, 2 R} )} 
 \les_{\e, \g} R^{\al -\hal + \epa} .
\end{align}
For $\g \in (1-\al, 1)$ and $\gam_* = \min( \al + \g-1, \tf18)$, we further obtain 
\beq
           R^2 \|  \na  \pa_{zz} \psis  \|_{ L^{\infty}( \bar B_{R, 2R} ) }
       + R^{ \gam_* + 2} 
       \|  \na  \pa_{zz} \psis  \|_{ \dot C^{ \gam_* }(  \bar B_{R, 2 R} )}
         \les_{\e, \g} R^{\al -\hal + \epa}.
\eeq
\eseq

\end{lem}

\begin{proof}

Denote 
\beq\label{eq:elli_pf_nota}
  M_R = \| \pa_z \ws \|_{ C^{\g}(\bar B_{3 R}) } 
  + \| \pa_z^{\leq 2} \ws(0, z) \|_{ L^{\infty}( [-3 R, 3 R] )} .
\eeq

We have 
\beq\label{eq:elli_dzF}
\bal
    \pa_z F_1 & =  \pa_z ( \ws(r, z) -\ws(0, z) \chi_1( \tf{ 8 r^2}{ \ang z^2 } ) )  r^{\al-1}  \\
 & =   \B( \big( \pa_z \ws(r, z) - \pa_z \ws(0, z) \chi_1\big( \f{ 8 r^2}{  \ang z^2} \big) \big)  
  -  \ws(0, z)  \cdot 8 \pa_z \f{r^2}{\ang z^2} \chi_1^{\pr}( \f{8 r^2}{\ang z^2} )    \B)  r^{\al-1}  .
\eal
\eeq

Since $\chi_1(8 r^2 / \ang z^2) = 1$ for $ 8r^2 \leq 2 \ang z^2$ and $\supp( \chi_1^{\pr}(8 r^2 /\ang z^2) ) \subset \{ \f{r}{ \ang z } \in [ \f14, \f12 ] \}$, for any $\xx \in \bar B_{2R}$, we obtain
\[
\bal 
|\pa_z F_1|  & \les_R \one_{r \les \ang z}  \| \pa_z \ws \|_{ C^{\g}(\bar B_{2 R}) }  r^{\al + \g - 1}
  + \one_{r \gtr  \ang z}  \| \pa_z \ws \|_{ L^{\infty}(\bar B_{2 R}) }  r^{\al-1}
   \les_R M_R ( 1 + r^{\al + \g - 1} ).
\eal
\]

If $\g > 1 -\al$, using Lemma \ref{lem:axi_reg} and the formula \eqref{eq:elli_dzF}, we obtain
\beq\label{eq:elli_dzF2}
   \| \pa_z F_1 \|_{ C^{\al + \g - 1}(\bar B_{2R}) }   \les_R   M_R.
\eeq

\paragraph{\bf Case $\g \leq 1-\al$}
We treat $\pa_z F_1(r, z)$ as an axisymmetric function in $\R^5$. 
For $\g < \f23$, we choose $p = \f{4}{ 0.7-\g}$ and estimate
\[
  \| \pa_z F_1 \|_{L^p(\bar B_{2R})} \les_R M_R ( 1 + 
  \B( \int_{|\yy| \leq 2 R} r^{p (\al + \g - 1)} d \yy \B)^{1/p} )
  \les_R M_R ( 1 +  ( \int_{ |(r, z)| \leq 2 R } r^{ p(\al + \g-1) + 3 } d r d z )^{1/p} ).
\]

Since $\g \leq 1-\al \leq 0.69$ and $\al > 0.31$, by definition, we get
\[
  p(\al + \g - 1) + 3 > p( \g - 0.7 ) + 3 = - 4  + 3 = -1, \, \quad  \forall  \, \g \in [0, 1-\al] .
\]

Thus, the above integral is finite and we obtain 
\beq\label{eq:elli_dzF3}
    \| \pa_z F_1 \|_{L^p(\bar B_{2R})}  \les_{R} M_R, \quad p = \f{4}{ 0.7-\g} 
    \in [ \f{4}{0.7}, \f{4}{0.01}]
\eeq

Applying the $W^{2,p}$ estimates in \cite[Theorem 9.11]{gilbarg1998elliptic} to \eqref{eq:elli_decomp:a}, \eqref{eq:elli_dec0}, and then using Morrey's inequality \cite[Theorem 7.17]{gilbarg1998elliptic} in dimension $d = 5$, we obtain
\beq\label{eq:schauder_reg:pf1}
      \|   \pa_z \Psi_1  \|_{ C^{1, \g_1 }(  \bar B_R )}
   \les  \|  \pa_z \Psi_1 \|_{ W^{2, p }( B_{3R/2}) }
   \les \| \pa_z F_1 \|_{L^{p  }( B_{2R} )} 
   +  \| \pa_z \Psi_1 \|_{L^{p }( B_{2R} )} 
   \les_R 1 + M_R ,
\eeq
where $\g_* =1 - \f{5}{ p }
= 1 - \f{5}{4}\cdot ( 0.7 - \g)
=  \f18 + \f54 \g$. 

\vs{0.1in}
\paragraph{\bf Case $\g > 1 -\al$}

In this case, the source term $\pa_z F$ is in $C^{\al+\g-1}$ by 
\eqref{eq:elli_dzF2}. Applying Schauder interior estimates \cite[Theorem 4.8]{gilbarg1998elliptic}, \eqref{eq:elli_dec0} and  \eqref{eq:elli_dzF2}, we obtain
\beq\label{eq:schauder_reg:pf2}
\bal
    \|   \pa_z \Psi_1  \|_{ C^{2, \al + \g -1}(  \bar B_R )}
   \les_{\g} \|  \pa_z \Psi_1 \|_{L^{\infty}( \bar B_{2R} )}
  +  \| \pa_z F_1  \|_{ C^{\al + \g -1}( \bar B_{2R} ) } 
  \les_R 1 + M_R.
\eal
\eeq

\paragraph{\bf Estimate of $\Psi_2$}

Next, we estimate $\Psi_2$ and $F_2$ from \eqref{eq:elli_decomp:b}. 
Applying derivations similar to the above for $\pa_z F_1$ with $\g =0$ and Leibniz rule, for 
$k=1,2$, we obtain
\[
\bal
  |\pa_z^k F_2(\xx)| \les_R ( r^{\al-1} + 1 ) \| \pa_z^{\leq k} \ws(0, z) \|_{L^{\infty}([-3 R, 3R])},
  \quad  \forall \ |\xx|  \leq 2 R .
\eal
\]
We treat $\pa_z^k F_2$ as an axisymmetric function in $\R^5$. Since $(\al-1) \f{40}{7} 
\geq - 0.69 \cdot \f{40}{7} > -4$, we obtain
\[
          \| \pa_z^k F_2 \|_{L^{ \f{40}{7} }(\bar B_{2R})} \les_R \| \pa_z^{\leq k} \ws(0, z) \|_{L^{\infty}([-3 R, 3R])}  .
\]

Applying the $W^{2,p}$ estimates in \cite[Theorem 9.11]{gilbarg1998elliptic} to \eqref{eq:elli_decomp:a}, \eqref{eq:elli_dec0}, and then using Morrey's inequality \cite[Theorem 7.17]{gilbarg1998elliptic} in dimension $d = 5$, for $1 - \f{5}{ 40/7} = \f18$, we obtain
\beq\label{eq:schauder_reg:pf3a}
\bal
        \|   \pa_z^k \Psi_2  \|_{ C^{1, \f{1}{8} }(  \bar B_R )}
  & \les \|  \pa_z^k \Psi_2 \|_{ W^{2, \f{40}{7} }( B_{3R/2}) }
   \les \| \pa_z^k F_2 \|_{L^{ \f{40}{7} }( B_{2R} )} 
   +  \| \pa_z^k \Psi_2 \|_{L^{ \f{40}{7} }( B_{2R} )}  \\
  & \les_R \| \pa_z^{\leq k} \ws(0, z) \|_{L^{\infty}([-3 R, 3R])} +   \| \pa_z^k \Psi_2 \|_{L^{ \f{40}{7} }( B_{2R} )}, 
  \eal
\eeq
for any $R > 0$. For $k = 0$, we bound $\|  \Psi_2 \|_{L^{40/7}( B_{2R})}$ in the upper bound using \eqref{eq:elli_dec0}. Applying the above estimates inductively to $k=1,2$, 
and adjust the radius $R$ if necessary, we prove 
\beq\label{eq:schauder_reg:pf3}
    \|   \pa_z^{\leq 2} \Psi_2  \|_{ C^{1, 1/ 8 }(  \bar B_R )}
    \les_{\e, R} 
    \|  \Psi_2 \|_{L^{40/7}( B_{2R} )} 
+ \| \pa_z^{\leq 2} \ws(0, z) \|_{L^{\infty}([-3 R, 3R])} 
\les_{\e, R} 1 + \| \pa_z^{\leq 2} \ws(0, z) \|_{L^{\infty}([-3 R, 3R])} .
\eeq
Combining \eqref{eq:schauder_reg:pf1}, \eqref{eq:schauder_reg:pf2},
 \eqref{eq:schauder_reg:pf3}, and using triangle inequality, we prove \eqref{eq:schauder_reg1}.

\vs{0.1in}
\paragraph{\bf (II) Decay estimates}

We fix $R > 1$ and define 
\beq\label{def:omR}
  \Om_R(\xx) = R^{\al} \ws(R \xx), 
  \quad \Psi_R(\xx) = R^{-1} \psis(R x).
\eeq
Then the elliptic equation for $\psis$ is equivalent to 
\[
- \D_{\R^5} \psis = \cpsia \ws r^{\al-1} \iff -\D_{\R^5}   \Psi_R(\xx)  = \cpsia  \Om_R(\xx) r^{\al-1}.
\]

We perform the decomposition as in \eqref{eq:elli_decomp}
\[
\bal
  - \D_{\R^5}  \Psi_{R,1}  &= \cpsia ( \Om_R(r, z) -\Om_R(0, z) \chi_1( \tf{8 r^2}{\ang z^2} ) )  r^{\al-1} \teq \cpsia F_{R, 1} ,
  \\
   - \D_{\R^5}  \Psi_{R, 2}  & = \cpsia \Om_{R}(0, z)  \chi_1( \tf{8 r^2}{\ang z^2} )  r^{\al-1} \teq \cpsia F_{R, 2} .
\eal
\]

In the annulus $B_{1,4}$, applying Lemma \ref{lem:axi_reg}, we have 
\beq\label{def:elli_dec_M1}
   \| \pa_z F_{R, 1} \|_{ C^{\al + \g - 1}(\bar B_{1.5, 3.5}) }   \les   
   \| \pa_z^{\leq 2} \Om_R(0, z) \|_{L^{\infty}( \{ 1\leq |z| \leq 4 \} )} 
    + \| \pa_z \Om_R \|_{ C^{\g}(\bar B_{1,4}) } 
   \teq M_1, 
\eeq
for $\g > 1 -\al$. For $\g \leq 1-\al$, following the proof of \eqref{eq:elli_dzF3}, we obtain
\[
         \| \pa_z F_{R, 1} \|_{L^p(\bar B_{1.5,3.5})}  \les M_1, \quad p = \tf{4}{ 0.7-\g} 
          \in [ \tf{4}{0.7},\, \tf{4}{0.01}]
\]

Applying the estimates and argument in the proof of 
\eqref{eq:schauder_reg:pf1}, \eqref{eq:schauder_reg:pf2}, \eqref{eq:schauder_reg:pf3a},  to $(\Psi_{R, 1}, \Psi_{R, 2}, \Om_{R})$ in the \emph{$R$-independent} annulus domain $B_{1,4}$, we obtain
\[
       \|   \pa_z \Psi_{R, 1}  \|_{ C^{1, \g_* }(  \bar B_{2,3} )}
   \les \| \pa_z F_{R, 1} \|_{L^{p  }( \bar B_{1.5, 3.5} )} 
   +  \| \pa_z \Psi_{R, 1} \|_{L^{p  }( \bar B_{1.5, 3.5} )} 
   \les M_1     +  \| \pa_z \Psi_{R, 1} \|_{L^{p  }( \bar B_{1.5, 3.5} )} , 
\]
for $\g \leq 1-\al, \g_* = \f{1}{8} + \f{5}{4} \g$;
\[
         \|   \pa_z \Psi_{R, 1}  \|_{ C^{2, \al+ \g - 1 }(  \bar B_{2,3} )}
         \les_{\g} \|  \pa_z \Psi_{R, 1} \|_{L^{\infty}( \bar B_{1.5, 3.5} )}
  +  \| \pa_z F_{R, 1}  \|_{ C^{\al + \g -1}( \bar B_{1.5,3.5} ) } 
  \les  M_1 + \|  \pa_z \Psi_{R, 1} \|_{L^{\infty}( \bar B_{1.5, 3.5} )},
\]
for $\g \in (1-\al , 1)$. Applying estimates similar to \eqref{eq:schauder_reg:pf3a}, \eqref{eq:schauder_reg:pf3} and by adjusting the size of the annulus, we obtain
\[
    \|   \pa_z^2 \Psi_{R, 2}  \|_{ C^{1,1/ 8 }(  \bar B_{2, 3} )}
   \les \| \pa_z^{\leq 2} \Om_R(0, z) \|_{L^{\infty}([ \{ 1\leq |z|\leq 4 \} ])} +   \|  \Psi_{R, 2} \|_{L^{40/7}( \{ \bar B_{1, 4} \} )} .
\]
All the above implicit constants are \emph{independent of $R$}.

Next, we apply \eqref{eq:elli_dec0} and \eqref{eq:ass:schauder_decay} to obtain the decay rates. 
Recall $M_1$ from \eqref{def:elli_dec_M1} and $\Om_R(\xx) = R^{\al} \ws(R \xx) $. 
Since $R > 1$, using assumption \eqref{eq:ass:schauder_decay} 
on $\ws$ and the decay estimates 
\eqref{eq:3D_decay:bd:w1}, \eqref{eq:3D_decay:bd:w2} for $\ws(0, z)$, we obtain
\[
\| \Om_R \|_{L^{\infty}(\bar B_{1,4} )}  + \| \pa_z \Om_R \|_{ C^{\g}(\bar B_{1,4}) } 
\les R^{\al-\als}, 
\quad \| \pa_z^{\leq 2} \Om_R(0, z) \|_{L^{\infty}([ \{ 1\leq |z|\leq 4 \} ])}
\les R^{\al-\als} \, ,
\quad M_1 \les  R^{\al-\als} .
\]

Recall $\Psi_{R, i} =R^{-1} \Psi_i(R \xx)$ from \eqref{def:omR}. Using \eqref{eq:elli_dec0}, we bound 
\[
\bal
 M_2 &\teq   \|  \Psi_{R, 2} \|_{L^{40/7}( \bar B_{1, 4} )} 
  +   \| \pa_z \Psi_{R, 1} \|_{L^p(  \bar B_{1, 4}  )} 
  + \|  \pa_z \Psi_{R, 1} \|_{L^{\infty}( \bar B_{1,4} )}  \\
& \les   \|  \Psi_{R, 2} \|_{L^{\infty}(  \bar B_{1, 4}  )} 
  +   \| \pa_z \Psi_{R, 1} \|_{L^{\infty}(  \bar B_{1, 4}  )} 
  + \|  \pa_z \Psi_{R, 1} \|_{L^{\infty}( \bar B_{1, 4} )} 
  \les R^{\al- \hal + \epa}.
\eal
\]

Since $-\hal + \epa > -\als$,, combining the above bounds, we prove $M_1 + M_2 \les R^{\al-\hal + \epa} $. Since $\Psi_R = \Psi_{R, 1} + \Psi_{R, 2}$, using triangle inequality and combining the above bounds on $\Psi_{R, i}$, we prove 
\[
       \one_{\g \leq 1-\al}  \|   \pa_{zz} \Psi_{R}  \|_{ C^{ \g_* }(  \bar B_{2,3} )}
       + \one_{\g > 1 -\al}  \|   \pa_{zz} \Psi_{R}  \|_{ C^{ 1, \min( \al + \g-1, \f18 ) } (  \bar B_{2,3} )}
         \les R^{\al -\hal + \epa}.
\]
Undoing the rescaling $\Psi_R(\xx)=  R^{-1} \psis(R \xx)$, and covering the 
dyadic annulus $B_{R, 2 R}$ by $B_{r_i, 3/2 r_i}$ with finite many $r_i \asymp  R$, we prove \eqref{eq:schauder_decay}.
\end{proof}

Let $X(s, \xx)$ be the flow map associated with $\QQs$
\beq\label{eq:solu_traj}
  \tf{d}{d s } X(s, \xx) = \QQs( X(s, \xx) ) , 
  \quad X(0, \xx) = \xx.
\eeq

Since $\ws = \wwb + \ww_{\al}$ satisfies the assumptions in Proposition \ref{prop:reg_solu}
and in Lemma \ref{lem:traj_reg}, the flow map $X(s, \xx)$ satisfies the estimates in 
Section \ref{sec:reg:traj_C1a}. 

Suppose that $f$ satisfies $|f(\xx)| \les \ang \xx^{-a}$ with some $a>0$ independent of $\xx$.
For any $|\xx|\geq 1,s \geq 0$, using $X(\tau, \xx) \geq e^{\bar c \tau} |\xx|
\geq e^{\bar c \tau} ,  \forall \, \tau \geq 0$ by \eqref{eq:traj_pf1}, we obtain
\beq\label{eq:basic_traj}
    | \int_0^s f( X(\tau, \xx) ) d \tau | \les \int_0^{\infty} \ang {X(\tau, \xx)}^{-a} d \tau
    \les \int_0^{\infty} e^{-a \bar  c \tau} d \tau \les_a 1.
\eeq

\subsubsection{Lower order decay estimates for $\ws$}

In this section, we prove the lower order decay bounds for the profile 
in \eqref{eq:3D_solu_decay:a}, \eqref{eq:3D_solu_decay:c}-\eqref{eq:3D_solu_decay:d}. Recall the weight $\rhoo^{-1} \les_{\e} \ang \xx^{\epa}$ from \eqref{def:3d_wg}. Using Proposition \ref{prop:vel_est} and Corollary \ref{cor:u_bar}, we obtain
\bseq\label{eq:solu_dec:pf1}
\begin{align}
  |\na \psis | + |\tf{1}{z} \psis| & \les \e^{-1} ( \rhoo^{-1}  \nnn{\ww_{\al}}
  + 1 )  \ang \xx^{\al-\hal} \les_{\e} \ang \xx^{-\d_{\e}}, \notag \\
 \tf{r}{\ang \xx} |  \pa_{rr} \psis|  + | \pa_{r z} \psis | 
+ ( \one_{ |z| > r } + \one_{|\xx|\leq 1})  |\pa_{zz} \psis |
 &  \les 
    ( \rhoo^{-1}  \nnn{\ww_{\al}} 
   +  1 )  \min(|\xx|^{\al},   \ang \xx^{\al-\hal-1} ) \notag \\
   &  \les  \min( |\xx|^{\al}, \ang \xx^{-\d_{\e}-1} ), 
\end{align}
where $\d_{\e}$ is defined as 
\beq
  \d_{\e}  \teq \hal - \al -\epa \gtr \e.
\eeq
\eseq
The inequality follows from \eqref{ran:ep_all}. Note that when $|z| > r$, we obtain $\ang z \asymp \ang \xx$. The estimate of $\pa_{zz} \psi$ in \eqref{eq:vel_zz}, \eqref{eq:u_bar} implies the above bound for $\pa_{zz}\psis$.

Using the formula of $\QQs$ from \eqref{eq:solu_al} and the above estimates, we obtain
\beq\label{eq:solu_dec:pf2}
  |\QQs - \cls \xx| \les |r \na \psis| + |\psis| 
  \les |\xx| \ang \xx^{ -\d_{\e} }. 
\eeq

Using \eqref{eq:solu_traj} and the derivation in \eqref{eq:3D_ODE:X_C1}, we obtain
\beq\label{eq:solu_traj_dX}
\bal
  \f{d}{d s}  | X(s, \xx) |  &=   \f{  \QQs( X(s, \xx)) 
  \cdot X(s, \xx) }{|X(s, \xx)|} , 
   \quad \f{d}{d s} \na X(s, \xx)
  = (\na \QQs)( X(s, \xx) ) \cdot \na X(s, \xx), 
\eal
\eeq
which implies
\[
\bal
 | X(s, \xx) |  &= e^{\cls s} |\xx| m(s, \xx),
  \quad m(s, \xx) = \exp\B( \int_0^s \f{ ( \QQs(X(\tau, \xx)) -\cls X(\tau, \xx) )
  \cdot X(\tau, \xx) }{| X(\tau, \xx)|^2}  d \tau   \B) . \\
  \eal
\]

Using estimate \ref{eq:basic_traj} with $f(\xx) =\f{ ( \QQs(\xx) -\cls \xx )
  \cdot \xx }{|\xx|^2} $ and \eqref{eq:solu_dec:pf2}, we bound 
  \[
    m(s, \xx) + m(s, \xx)^{-1} \les_{\e} 1. 
  \]

Therefore, we obtain
\beq\label{eq:traj_far}
|X(s, \xx)| \geq |\xx|, \quad 
   |X(s, \xx)| \asymp_{\e} e^{\cls s} |\xx|, \quad \forall \, |\xx| \geq 1 , \quad s \geq 0.
\eeq

\paragraph{\bf Decay estimates of $\ws$}

Recall the profile equation  for $\ws$  and $\cBs$ from \eqref{eq:solu_al} 
\[
     \QQs \cdot \na \ws =  \cBs \cdot \ws,
     \quad  \cBs= \cws -(1-\al) \pa_z \psis .
\]

Using the flow map $X(s, \xx)$, for any $s > 0$, we obtain
\beq\label{eq:traj_linf}
\ws(X(s, \xx)) =  e^{\cS_{\cB}(s, \xx)} \ws(\xx), 
 \quad \cS_{\cB}(s, \xx) = \int_0^s \cBs \cc X(\tau, \xx) d \tau .
\eeq

Using \eqref{eq:solu_dec:pf1} and estimate \ref{eq:basic_traj},  we obtain
\beq\label{eq:traj_linf:S1}
  \B| \int_0^s ( \cBs(X(\tau,\xx)) - \cws)  d \tau \B| \les_{\e } 1,
  \quad \forall \, |\xx| \geq 1 ,
  \quad \Rightarrow e^{\cS_{\cB}(s, \xx)} \asymp_{\e} e^{\cws s}.
\eeq

Using the above estimates and \eqref{eq:traj_far}, for $|\xx| \in [1, 2]$, we obtain
\[
  |X(s, \xx)|^{- \f{ \cws}{\cls} } |\ws(X(s, \xx))| \asymp_{\e} 
|\xx|^{-  \f{ \cws}{\cls}  } e^{ -\f{ \cws}{\cls}  \cdot \cls s  }
  \cdot e^{\cws s} | \ws(\xx) |  \asymp_{\e} 1,
\]
uniformly in $s$. 
Since $\ws \in C^{1,\al}(\bar B_1)$ \eqref{eq:solu_basic_reg}, $\ws(0)=0$,  and since for any $|\yy| \geq 1$, there exists $|\xx| \in [1,2]$ and $s \geq 0$ such that $X(s, \xx) = \yy$, we prove 
\beq\label{eq:ws_dec1}
   |\ws(\xx) | \les_{\e} \min( |\xx|,  \ang \xx^{\cws /\cls} ) = \min( |\xx|,  \ang \xx^{-\als} ),
   \quad \als = - \cws /\cls.
\eeq

For $r > |z|$ and $|\xx| \geq 1$, we obtain $r \asymp |\xx|$. Using the elliptic equation \eqref{eq:Euler2_psi}  
for $\psis$, and \eqref{eq:solu_dec:pf1} for $\pa_{rr} \psis, \pa_r \psis$, and \eqref{eq:ws_dec1} for $\ws$, we estimate 
\[
 |\pa_{zz} \psis| \les \f{1}{r} ( | \pa_{rr}(r \psis) |  + |\pa_r \psi| + |\ws r^{\al}| )
 \les_{\e} \ang \xx^{-\d_{\e}-1} 
 + \ang \xx^{-\als + \al-1} \les  \ang \xx^{-\d_{\e}-1} .
\]
where we have used $\al -\als < \al-\hal + \epa = -\d_{\e}$ 
by \eqref{eq:3D_solu_scale} in the last inequality. 

Combining the above estimate and  \eqref{eq:solu_dec:pf1}, we prove 
\bseq\label{eq:solu_dec:pf3}
\beq
 |\pa_{zz} \psis|  \les_{\e} \min( |\xx|^{\al} , \ang \xx^{-\d_{\e}-1} ).
\eeq 

Using \eqref{eq:solu_dec:pf1}, \eqref{eq:solu_dec:pf2}, \eqref{eq:solu_dec:pf3}, we obtain
\beq
  |\na \QQs - \cls \Id | \les_{\e} \ang \xx^{ -\d_{\e} }.
\eeq
\eseq

Since $\na X(0, \xx) = \Id$ for $|\xx| \geq 1$, estimating $\na X$ in \eqref{eq:solu_traj_dX} along the trajectory, using the above estimates and estimate \ref{eq:basic_traj}, we bound
\beq\label{eq:dX_far}
  |\na  X(s, \xx)| \les_{\e} e^{\cls s} .
\eeq

Given $\xx_1 \neq \xx_2$ with $|\xx_i| \geq 1$, using \eqref{eq:solu_traj}, we estimate 
the difference $| X(s, \xx_1) - X(s, \xx_2)|^2$:
\[
\bal
  & \tf{1}{2} \tf{d}{ds} | X(s, \xx_1) - X(s, \xx_2)|^2
   = \cls | X(s, \xx_1) - X(s, \xx_2)|^2  \\ 
  & \qquad +\big( (\QQs - \cls \xx) \cc X(s, \xx_1)  - (\QQs - \cls \xx) \cc X(s, \xx_2) \big) \cdot  ( X(s, \xx_1) - X(s, \xx_2) )  \\
\eal
\]

For any $|\xx| \geq 1$, since $|X(s, \xx) | \gtr_{\e} e^{\cls s}$ by \eqref{eq:traj_far}
using 
\eqref{eq:solu_dec:pf3}, we estimate 
$\ang {X(s, \xx)}^{-\d_{\e} } \les_{\e} e^{ -\d_{\e} \cls s}$ and 
\[
    \tf{1}{2} \tf{d}{ds} | X(s, \xx_1) - X(s, \xx_2)|^2
    \geq \cls | X(s, \xx_1) - X(s, \xx_2)|^2 - C e^{-\d_{\e} \cls s} | X(s, \xx_1) - X(s, \xx_2)|^2 .
\]
Integrating over $s$, we obtain
\beq\label{eq:traj_expand}
  | X(s, \xx_1) - X(s, \xx_2)| \gtr_{\e} e^{\cls s} |\xx_1 - \xx_2|,
  \quad \forall \, |\xx_i| \geq 1, \, s \geq 0.
\eeq

 For any $ \g \in [0, 1]$ and $g$, using the composition and \eqref{eq:dX_far}, and \eqref{eq:traj_far}, we have the following estimates for $\dot C^{\g_i}(\bar B_{1, 2})$ semi-norm 
\beq\label{eq:hol_compo_far}
  \| g(X(s, \cdot)) \|_{\dot C^{\g}(\bar B_{1, 2}) }
  \les \| g \|_{\dot C^{\g}( X(s, \bar B_{1,2}) )} \cdot  \| \na X(s, \cdot) \|_{L^{\infty}
(\bar B_{1, 2})}^{\g} 
\les    e^{\cls \g s} \cdot  \| g \|_{\dot C^{\g}(  c_1 e^{\cls s} \leq |\xx| \leq c_2 e^{\cls s} )} ,
\eeq
where  $X(s, A) = \{ \yy : \yy = X(s, \xx), \xx \in A\}$ denotes the image of $A$ under the flow map.

\vs{0.05in}

\paragraph{\bf Decay estimates of $\na \ws$}

Next, we estimate decay for $\na \ws$. 

Recall $\cS_{\cB}$ from \eqref{eq:traj_linf} and $\cBs = \cws - (1-\al) \pa_z \psis$. Using 
the bounds for $\na \pa_z \psis$ from  \eqref{eq:solu_dec:pf1}, \eqref{eq:solu_dec:pf3}, and \eqref{eq:hol_compo_far} with $\g = 1$,  we bound 
\[
  \|  \cS_{\cB}(s, \cdot) \|_{ \dot C^{0, 1}(\bar B_{ 1, 2 }) }
  \les \int_0^s \|  \pa_z \psis \|_{\dot C^{0, 1}( X(\tau , \bar B_{1,2}) ) } e^{\cls \tau } d \tau
  \les \int_0^s e^{- \cls \tau (\d_{\e} + 1) } e^{\cls \tau } d \tau \les_{\e } 1.
\]

Combining the above estimates and \eqref{eq:traj_linf:S1}, we bound 
\[
  \| e^{\cS_{\cB}(s, \cdot )} \|_{ \dot C^{0, 1}(\bar B_{ 1, 2 }) }
  +   \| e^{\cS_{\cB}(s, \cdot )} \|_{ L^{\infty}(\bar B_{ 1, 2 }) }
  \les (  1 +  \|  \cS_{\cB}(s, \cdot ) \|_{ \dot C^{0, 1}(\bar B_{ 1, 2 }) }  )
  \| e^{\cS_{\cB}(s, \cdot )} \|_{ L^{\infty}(\bar B_{ 1, 2 }) } \les e^{\cws s}.
\]

Applying the $C^1$ bound to $\ws(\xx)$ from \eqref{eq:solu_basic_reg} and using product rule to 
\eqref{eq:traj_linf}, we establish 
\[
  \| \ws(X(s, \cdot ))  \|_{ \dot C^{0, 1}(\bar B_{ 1, 2 }) } = \|  e^{\cS_{\cB}(s, \xx)} \ws(\xx) \|_{ \dot C^{0, 1}(\bar B_{ 1, 2 }) } \les e^{\cws s}.
\]

Using \eqref{eq:traj_expand}, for any $\xx_1\neq \xx_2 \in \bar B_{1,2}$, we obtain
\[
   \f{|\ws(X(s, \xx_1)) - \ws(X(s, \xx_2))|}{  | X(s,\xx_1) - X(s,\xx_2) |   }
  =    \f{|\ws(X(s, \xx_1)) - \ws(X(s, \xx_2))|}{  |\xx_1 - \xx_2|   }
  \cdot \f{  |\xx_1 - \xx_2|  }{| X(s,\xx_1) - X(s,\xx_2) | }
\les_{\e} e^{ \cws  s} \cdot e^{-\cls  s }.
\]

Since  $\ws \in C^{1,\al}(\bar B_R)$ for any $R > 0$ \eqref{eq:solu_basic_reg}, from the above estimates, for any $s \geq 0$,  we prove 
\[
 \| \na \ws \|_{L^{\infty}( (  X(s, B_{1, 2} ) ) ) } \les  \| \ws \|_{ \dot C^{0,1}(  X(s, B_{1, 2} ) ) } \les_{\e} e^{ (\cws - \cls) s} .
\]

Since for any $R > 1$, we can cover $\bar B_{R, 2 R}$ by $ X(s, B_{1, 2} ) $ with finite many
$s_i $ and $ e^{\cls s_i} \asymp_{\e} R$, from the above estimate, for any $R > 1$, we establish
\beq\label{eq:ws_dec2}
  \| \na \ws  \|_{ L^{\infty}(\bar B_{R, 2 R} )  } \les_{\e} R^{ (\cws - \cls) / \cls } = R^{-\als - 1}  .
\eeq

\subsubsection{High order decay estimates for $\ws$}

We consider the following sequences of regularity indices 
\[
 \g_i = \tf{1}{8} i, \quad 0\leq i\leq 8,
\quad n_{ \g} = 8
\]

Below, we use equation \eqref{eq:3D_tran:ws1} and decay estimates in \eqref{eq:schauder_decay} to obtain sharp decay estimates for $\| \pa_z \ws \|_{\dot C^{\g_i}}$. We consider $ 0\leq i \leq n_{\g}-1$.  We write $\dot C^0$ with $L^{\infty}$ and $\| g \|_{\dot C^{\th}( \S  )} = \| \na g \|_{\dot C^{1,\th-1}(\S) }$ for $\th \in ( 1, 2 )$.  Given the decay of $\pa_{zz} \psis$, we estimate $\pa_z \ws $ along trajectory:

\begin{lem}[\bf Decay estimates of $\pa_z \ws$]\label{lem:traj_decay}
Suppose that 
\beq\label{eq:ass:dec1}
   \|\pa_{zz} \psis \|_{ \dot C^{\g}(\bar B_{R, 2 R}) } \leq C(\e) R^{-\d_{\e}- 1 - \g},
   \quad \d_{\e} = \hal - \al - \epa \gtr \e .
\eeq
for any $R >1$. Then for any $R >1$, we have 
\beq\label{eq:traj_decay:ws}
    \| \pa_z \ws \|_{ \dot C^{\g}(\bar B_{R, 2 R}) } \les_{\e} R^{ -1 - \als - \g }.
\eeq
\end{lem}

We defer the proof to Section \ref{sec:traj_decay}. 

Now, we use Lemma \ref{lem:traj_decay} and Lemma \ref{lem:schauder_decay} to bound $\pa_z \ws$.

For $i=0$, since \eqref{eq:solu_dec:pf3} and \eqref{eq:ws_dec1} imply  the  assumption \eqref{eq:ass:dec1} with $\g = \g_0 = 0$. Applying Lemma \ref{lem:traj_decay}, we obtain 
$C^{\g_i}$ estimates for $\ws$ \eqref{eq:traj_decay:ws}.  Estimate \eqref{eq:traj_decay:ws} and \eqref{eq:ws_dec1}, \eqref{eq:ws_dec2}  implies assumption \eqref{eq:ass:schauder_decay} in Lemma \ref{lem:schauder_decay} with $\g = \g_i$. 
Thus, applying Lemma \ref{lem:schauder_decay} with $\g = \g_i$ 
and using $\g_{i+1} \leq \g_*$ with $\g_*$ chosen in Lemma \ref{lem:schauder_decay}, we obtain $C^{\g_{i+1}}$  estimates for $\pa_{zz} \psis$ \eqref{eq:schauder_decay},  which further implies the conditions 
\eqref{eq:ass:dec1} with $\g = \g_{i+1}$.

Applying the above arguments for $i=0,.., n_{\g}-1$,  for any $R>1$, we prove 
\[
   \| \na \pa_z \ws\|_{L^{\infty} (\bar B_{R, 2 R}) }
\les  \| \pa_z \ws \|_{ \dot C^{0, 1}(\bar B_{R, 2 R}) } \les R^{ -2 - \als  } \, .
\]

Using $\pa_z \ws \in C^{1,\al}$ from \eqref{eq:solu_high_reg1} and the above estimate for any $R>1$, we prove 
\beq\label{eq:ws_dec3:a}
  | \na \pa_z \ws(\xx) |  \les_{\e} \ang \xx^{-2- \als}.
\eeq

Next, we estimate $\pa_r \ws$. Taking $\pa_r$ in \eqref{eq:solu_al}, we obtain
\[
\QQQs^r \cdot \pa_r^2 \ws + \QQQs^z \pa_{r z} \ws  = \pa_r ( \cBs \ws ) - \pa_r \QQs \cdot \na \ws,
\quad \cBs = \cws - (1-\al) \pa_z \psis.
\]

Using \eqref{eq:ws_dec3:a}, \eqref{eq:ws_dec1}, \eqref{eq:ws_dec2} for $\ws$,
and \eqref{eq:solu_dec:pf1}, \eqref{eq:solu_dec:pf3} for $\psis$,  $ \QQQs^r = r q^r$ with $q^r \gtr \e^{-1}$ by \eqref{eq:3D_boot_q:1},  for any $|\xx| > 1$, we prove 
\beq\label{eq:ws_dec4}
| r \pa_{rr} \ws (\xx) | 
\les |  \QQQs^r \cdot \pa_r^2 \ws |
\les  \ang \xx^{-\als - 2 + 1}
+ |\na \ws| \cdot ( |\pa_r \QQs| + |\cBs | )
+ |\pa_{rz} \psis| \cdot |\ws|
  \les_{\e} \ang \xx^{-\als -1}.
\eeq

\paragraph{\bf Proof of Proposition \ref{prop:3D_profile}}
Using \eqref{eq:ws_dec1}, \eqref{eq:ws_dec2}, \eqref{eq:ws_dec3:a}, and \eqref{eq:ws_dec4}, we prove 
the decay estimates in \eqref{eq:3D_solu_decay:a}, \eqref{eq:3D_solu_decay:b} for $\ws$. 
Since $\pa_r \ws(r, z)$ is odd in $z$ and $\ws \in C^{1,\al}$, using \eqref{eq:solu_basic_reg}, 
we prove the estimate on $\pa_r \ws$ in  \eqref{eq:3D_solu_decay:a}.

 Using the estimates \eqref{eq:solu_dec:pf1}, \eqref{eq:solu_dec:pf3} for $\psis$, we prove the decay estimates \eqref{eq:3D_solu_decay:c}, \eqref{eq:3D_solu_decay:d} for $\psis$. 
Estimate \eqref{eq:3D_solu_decay:c2} follows from \eqref{eq:3D_solu_decay:d} and $\psis(r, 0) = 0, \pa_r \psis(r, 0) = 0$. Upon changing the absolute constants in \eqref{eq:3D_solu_decay}, we obtain the same constant $\upse$.

Since $\ws = \ww_{\al} + \wwb$ and $\nnn{\ww_{\al}} \leq \e^{\hk}$ by Theorem \ref{thm:3D_solu}, 
using \eqref{eq:traj_outgo} in Proposition \ref{prop:reg_solu}, we prove the outgoing properties \eqref{eq:outgo}.

Note that we have proved \eqref{eq:3D_solu_scale} in Section \ref{sec:3D_solu:scale}. 
Using \eqref{eq:solu_basic_reg} and \eqref{eq:solu_high_reg1}, we prove \eqref{eq:3D_solu_reg}.

\subsubsection{Proof of Lemma \ref{lem:traj_decay}}\label{sec:traj_decay}

Recall the equation  \eqref{eq:3D_tran:ws1}
\[
\bga
    \QQs \cdot \na \Ups = a(\xx) \cdot \Ups + f(\xx), \quad 
    \Ups(r, z) \teq  \tf{ \pa_z \ws }{q^r} , \\
  a(\xx) = \cws - \cls - (3 -\al) \pa_z \psis , \quad
    f(\xx)  
  =  ( \cws - (1-\al) \cls) \cdot \ws \cdot (q^r)^{-2}  \cdot \pa_{zz} \psis,
  \quad q_r = \cls - \pa_z \psis.
\ega
\]

Using \eqref{eq:3D_tran:ws1} and integrating $\Ups$ along trajectory, for any $s \geq 0$, we obtain
\beq\label{eq:solu_dec:pf4}
  \Ups(X(s, \xx)) 
  = e^{\cS(s, \xx)} \Ups(\xx)
  + \int_0^s e^{\cS(s, \xx) - \cS(\tau, \xx) } f( X(\tau, \xx) ) d \tau,
  \quad \cS(s, \xx) = \int_0^s a( X(\tau, \xx ) ) d \tau.
\eeq

Using \eqref{eq:solu_dec:pf1}, \eqref{eq:solu_dec:pf2}, \eqref{eq:solu_dec:pf3}, and \eqref{eq:ws_dec1}, and interpolation, for any $\s \in (0, 1]$, we bound 
\bseq\label{eq:solu_dec:pf5}
\beq
\bal
  & |a(\xx) - (\cws - \cls)| \les_{\e} \ang \xx^{-\d_{\e}},  
&    |\na a(\xx) | & \les_{\e} \ang \xx^{ -\d_{\e} - 1 }, 
&   \|  a  \|_{\dot C^{ \s }(\bar B_{R, 2 R}) } & \les \ang R^{-\d_{\e} - \s}, 
 \\
  & | f(\xx) |   \les_{\e} \ang \xx^{-\als - 1 - \d_{\e}} ,  & 
 q^r &  \asymp_{\e} 1,
&   \|  q^r  \|_{\dot C^{ \s }(\bar B_{R, 2 R}) } 
& \les_{\e} \ang R^{ -\d_{\e} - \s }, 
\eal
\eeq
where we have used $ \|  a  \|_{\dot C^{ \s }} = \|  a(\cdot) - (\cws - \cls)  \|_{\dot C^{ \s } }$. 

Using the assumption on decay \eqref{eq:ass:dec1} for $\pa_{zz} \psis$, 
interpolating \eqref{eq:ws_dec1}, \eqref{eq:ws_dec2} for $\ws$, and using the product rule, $q^r \gtr 1$, and \eqref{eq:solu_dec:pf3}, \eqref{eq:ws_dec1},  for any $R > 1$, we obtain
\beq
\bal
   \| f  \|_{\dot C^{ \g }(\bar B_{R, 2 R}) } 
 & \les_{\e}  
    \| \ws \|_{\dot C^{ \g }(\bar B_{R, 2 R}) }    \|  \pa_{zz} \psis \|_{ L^{\infty} (\bar B_{R, 2 R})  }
    +  \| \ws \|_{L^{\infty}(\bar B_{R, 2 R}) }    \|  \pa_{zz} \psis \|_{ \dot C^{ \g } (\bar B_{R, 2 R})  } \\
 &\quad  + \one_{\g >0}   \| \ws \|_{L^{\infty}(\bar B_{R, 2 R}) }    \|  q^r \|_{ \dot C^{ \g } (\bar B_{R, 2 R})  }
    \| \pa_{zz} \psis \|_{L^{\infty}(\bar B_{R, 2 R}) }  \les  R^{-\als -\g -\d_{\e}-1 }.
\eal
\eeq
\eseq

Heuristically,  for each $\g$-fractional derivatives, we gain a decay factor $R^{-\g}$.

Applying estimate \ref{eq:basic_traj} with $f=a(\xx) - (\cws - \cls) $, for $|\xx| \geq 1$,  we bound 
\beq\label{eq:solu_dec:pf6}
  |\cS(s, \xx) -  (\cws - \cls) s | \les_{\e} 1 ,\, \Rightarrow  \, e^{\cS(s, \xx) } \asymp_{\e} e^{  (\cws - \cls) s }.
\eeq

Recall the formula of $\cS$ from \eqref{eq:solu_dec:pf4} and denote $X(\tau, A) = \{ X(\tau, \xx) : \xx \in A\}$. Using \eqref{eq:traj_far},
we obtain $|\yy| \asymp_{\e} e^{\cls \tau}$ for any $\yy \in X(\tau, \bar B_{1,2})$. 
Using \eqref{eq:solu_dec:pf5}, \eqref{eq:hol_compo_far}, 
for $s_1 < s_2$, we bound $\cS$ as 
\[
\bal
& \one_{\g > 0} \|  \cS(s_2, \cdot) - \cS(s_1, \cdot) \|_{ \dot C^{\g}(\bar B_{1, 2}) }  
 \les  \one_{\g > 0} \int_{s_1}^{s_2} \| a \|_{ \dot C^{\g}( X(\tau, \bar B_{1, 2}) ) }  e^{\cls \g \tau } d \tau  \les \int_0^s e^{- \cls \tau ( \d_{\e}  + \g ) + \cls \g \tau} \les_{\e} 1 .
\eal
\]
Combining the above estimate and \eqref{eq:solu_dec:pf6}, we obtain
\[
\bal
  & \|  e^{\cS(s, \cdot) - \cS(\tau, \cdot)  } \|_{ L^{\infty}(\bar B_{1, 2}) } 
   + \|  e^{\cS(s, \cdot) - \cS(\tau, \cdot)  } \|_{ \dot C^{\g}(\bar B_{1, 2}) }  \\
 & \quad  \les     \|  e^{\cS(s, \cdot) - \cS(\tau, \cdot)  } \|_{ L^{\infty}(\bar B_{1, 2}) } 
  (1 + \one_{\g > 0} \|  \cS(s, \cdot) - \cS(\tau, \cdot)  \|_{ \dot C^{\g}(\bar B_{1, 2}) } )
\les e^{ (\cws - \cls ) (s-\tau)} .
\eal
\]

By covering the region $ \{ c_1 e^{\cls s} \leq |\xx| \leq c_2 e^{\cls s} \}$ 
in \eqref{eq:hol_compo_far} by finite many regions $ B_{ r_i e^{\cls s}, 2 r_i e^{\cls s} }$ with $r_i \asymp_{\e} 1$, and using 
\eqref{eq:solu_dec:pf5} with $R \asymp_{\e} e^{\cls s}$, \eqref{eq:hol_compo_far}, 
we bound $f$ as 
\[
\bal
   \|  f(X(s, \cdot)) \|_{ \dot C^{\g}(\bar B_{1, 2}) } 
  & \les e^{ - \cls s (\als + \g + \d_{\e} + 1 )} \cdot e^{\cls \g s}
   \les e^{ - \cls s (\als +  \d_{\e} + 1 )} ,\\
     \|  f(X(s, \cdot)) \|_{ L^{\infty}(\bar B_{1, 2}) } 
    &  \les e^{ - \cls s (\als  + \d_{\e} + 1)} .
\eal
\]

Applying product rule, the above estimates, and $\cls (\als + 1) =\cls - \cws$ by \eqref{eq:ws_dec1},
 for any $s \geq \tau \geq 0$, we bound 
\[
\bal
  \| e^{\cS(s, \cdot)} \Ups(\cdot) \|_{ \dot C^{\g}(\bar B_{1, 2}) } 
 & \les e^{(\cws - \cls) s}, \\
 \| e^{ \cS(s,\xx) -  \cS(\tau, \xx) } f( X(\tau, \xx) ) \|_{ \dot C^{\g}(\bar B_{1, 2}) } 
 & \les e^{(\cws - \cls) (s -\tau)}  \cdot e^{ - \cls \tau (\als  + \d_{\e} + 1 )}
 \les  e^{(\cws - \cls) s -\cls \tau  \d_{\e}}.
\eal
\]

Since $\d_{\e} > 0$ by \eqref{eq:solu_dec:pf1}, 
integrating the above estimates in \eqref{eq:solu_dec:pf4} over $\tau$, we establish
\[
    \| \Ups(X(s, \cdot)) \|_{ \dot C^{\g}(\bar B_{1, 2}) } 
    \les_{\e} e^{(\cws - \cls) s}( 1 + \int_0^{s} e^{-\cls \tau \d_{\e}} d \tau  )
    \les_{\e} e^{(\cws - \cls) s}.
\]

If $\g > 0$, using \eqref{eq:traj_expand}, for any $\xx_1\neq \xx_2 \in \bar B_{1,2}$, we obtain
\[
   \f{|\Ups(X(s, \xx_1)) - \Ups(X(s, \xx_2))|}{  | X(s,\xx_1) - X(s,\xx_2) |^{\g}   }
  =    \f{|\Ups(X(s, \xx_1)) - \Ups(X(s, \xx_2))|}{  |\xx_1 - \xx_2|^{\g}   }
  \cdot \f{  |\xx_1 - \xx_2|^{\g}  }{| X(s,\xx_1) - X(s,\xx_2) |^{\g} }
\les e^{(\cws - \cls) s} \cdot e^{-\cls \g s }.
\]

Thus, for any $s \geq 0$, we prove 
\[
  \| \Ups  \|_{ \dot C^{\g}(\bar X(s, B_{1, 2} ) ) } \les e^{ (\cws - \cls-\cls \g) s} .
\]

Since $X(s, \cdot)$ is a bijection, using \eqref{eq:traj_far}, 
for any $R > 1$, we can cover $\bar B_{R, 2 R}$ by 
$N_{\e}$ domains  $\bar X(s_i, B_{1, 2} ) $
with  $s_i  $ satisfying $ |X( s_i, 1)| \asymp_{\e}  e^{\cls s_i} \asymp_{\e}  R, i\leq N_{\e}$. Thus, we establish
\[
  \| \Ups  \|_{\dot C^{\g} ( \bar B_{R, 2 R} )} \les_{\e} R^{  (\cws - \cls-\cls \g) / \cls }
  = R^{ - 1 - \als -\g }.
\]

Recall $\Ups = \f{\pa_z \ws}{q^r}$ from \eqref{eq:3D_tran:ws1}. 
Using the H\"older estimate for $q^r$ in \eqref{eq:solu_dec:pf5}, the above estimate and product rule, we prove 
\[
  \| \pa_z \ws \|_{\dot C^{\g} ( \bar B_{R, 2 R} )} \les  R^{ - 1 - \als -\g }.
\]
for any $R>1$. We complete the proof of Lemma \ref{lem:traj_decay}.

\subsection{Proof of Theorem \ref{thm:self-similar}} \label{sec:proof_self_similar}

We are in a position to prove Theorem \ref{thm:self-similar} on the self-similar profile.

\begin{proof}[Proof of Theorem \ref{thm:self-similar} ]

We take \(\bar\e=\beps_{12}\) in Theorem \ref{thm:self-similar}, where
\(\beps_{12}\) is chosen in Proposition \ref{prop:3D_profile}.

\paragraph{\bf Proof of item (i), (ii)}
Since the self-similar equation \eqref{eq:profi_3D_omth} for $\omth = r^{\al} \www$ is equivalent to
\eqref{eq:profi_3D0}, \eqref{eq:profi_3D} for $\www$, 
using the profile $(\ws, \cls, \cws)$ constructed in Theorem \ref{thm:3D_solu} for \eqref{eq:profi_3D}, we construct the profile $(\omsth = r^{\al} \ws, \cls, \cws)$  for  \eqref{eq:profi_3D_omth}. 
From Proposition \ref{prop:3D_profile} and  \eqref{eq:3D_solu_decay}, we obtain 
\[
 \ws \in C^{1,\al}, \quad |\ws| \les \ang \xx^{ - \als}, \quad |\na \ws| \les \ang \xx^{-\als - 1}, 
 \quad \als > \al,
 \quad \Rightarrow
  \ \omsth = \ws r^{\al} \in C^{\al}(\R^3) .
 \]
Using \eqref{eq:3D_solu_scale} and $\e = \f13 - \al$, we estimate 
the rescaled scaling parameter 
$\cxs$ for \eqref{eq:SS_ansatz_1} and \eqref{eq:profi_3D_omth0}
\begin{align}\label{eq:thm_SS_pf0}
 \cxs & \teq - \f{ \cls}{ \cws + \al \cls}
 =  \f{1}{\als - \al } \in [ \f78 \e^{-1}, \f{9}{10} \e^{-1} ],  \\
\cxs & =  \f{1}{\als - \al }
 = \B( \f13 + \f{\e}{8} - \al + O(\e^{1 +  \mhk}) \B)^{-1}
 = ( \f{9 \e}{ 8} + O(\e^{1 +  \mhk}) )^{-1}
 =  \f{8}{9} \e^{-1} + 
  O(\e^{ -1+  \mhk}) . \notag
\end{align}

 Using the above estimate $\cls = \tf{64}{9} \e^{-1} + O(\e^{- \hk})
  \asymp \e^{-1}$  \eqref{eq:3D_solu_scale},  we obtain
$\cws + \al \cls  < 0$ and 
\beq\label{def:TTa}
\TTa \teq  - ( \cws + \al \cls)^{-1}=   \f{ \cxs } {\cls}
=  \f{   \f{8}{9} \e^{-1} }{ \f{64}{9} \e^{-1} }  + O( \e^{\mhk})
 =  \f18  + O( \e^{\mhk}) \asymp 1.
\eeq

Recall the parameter $\cpsia(\al) $ from \eqref{def:cpsi}. 
For any $\e = \f13 - \al \leq \beps_{ 12}$, using the relation \eqref{eq:profi_3D_rescale}
with $(c_l, c_{\om}, \cxx )$ replaced by $( \cls, \cws, \cxs )$, we construct the self-similar profile for \eqref{eq:profi_3D_omth0} and \eqref{eq:SS_ansatz_1} 
\beq\label{eq:thm_SS_pf2}
\omthss \teq \cpsia \TTa  \omsth  = \AAa r^{\al} \ws , 
\quad \AAa \teq  \cpsia \TTa .
\eeq
with scaling $\cxs$ satisfying \eqref{eq:thm_SS_pf0}. Using estimate \eqref{def:cpsi} for $\cpsia$  and  \eqref{def:TTa} for $\TTa$, we estimate $\AAa$:
\beq\label{def:AAa}
 \AAa =  \AAaa + O( \e^{\mhk}) \asymp 1,  \quad  \AAaa \teq \tf18 \barcp .
\eeq

Using the formula  \eqref{eq:SS_ansatz_1} and \eqref{eq:thm_SS_pf0}, we prove the self-similar blowup result \eqref{eq:thm_SS}, and estimates 
\eqref{eq:thm_profi_decay:a} and  \eqref{eq:thm_profi_decay:b}.
Since for any $\al \in (\f13- \beps_{12}, \f13)$ and any $\g \in [0, \al]$, we have $\omsth \in C^{\al}(\R^3 ) \subset C^{\g}(\R^3)$, we prove the existence of nontrivial $C^{\g}$ self-similar blowup solutions and profile
for any $\g \in (0, \f13)$.

\vspace{0.05in}
\paragraph{\bf Proof of item (ii)}
 Using \eqref{def:3d_wg} and 
\eqref{eq:Ja_hat}, we obtain 
\beq\label{eq:thm_SS_pf1}
\rhoo \gtr \ang {|\xx|}^{-\epa} |\cJaa(\xx)|^{-\kp} \gtr  \ang {|\xx|}^{-\epa} |\lgp x|^{-\kp},
\quad \rag \asymp  \tf{ \ang z^{\kag} }{ \ang {|\xx|}^{\kag} } \les 1 \, ,
\quad  |\xx| = |(r, z)|.
\eeq

Since $\ws = \wwb + \ww_{\al}$, 
using norms $\nnla{\cdot}, \nnn{\cdot}$ from \eqref{def:3d_norm}
and  Theorem \ref{thm:3D_solu},  we obtain
\[
\nlinf{ \ang z^{\hal + \kag }  \ang {|\xx|}^{-\kag -\epa} |\lgp \xx |^{-\kp}  (\ws - \wwb) }
\les \nlinf{ \ang z^{\hal} \rag \rhoo \ww_{\al}  }
\les \nnla{\ww_{\al}} \les  \nnn{\ww_{\al}}  \les \e^{1 - 5 \kp_1}.
\]

Since $|z| \leq |\xx|$, using estimates \eqref{eq:waa_error} for $ \waa - \ang x^{\b} \wa$ with $|z|^{-\bbb} \gtr 1$ \eqref{norm:Xc} when $|z| \leq 1$,
and $ | (\hal - \epa) - ( \alb - \b - \kp_1 \e) | \les \e^{2-\kp}$ by \eqref{ran:ep_all}, we obtain
\[
 \ang z^{\hal + \kag }  \ang {|\xx|}^{-\kag -\epa} |\lgp \xx |^{-\kp}  ( \waa(z) - \ang z^{\b} \wwwa(z))
\les  \ang z^{\hal - \epa } \cJaa(z)^{-\kp} |  ( \waa - \ang x^{\b} \wwwa |
\les \e \ang z^{ C \e^{2-\kp}}.
\]

Denote $h =  \ang {|\xx|}^{-\kag -\epa} |\lgp \xx |^{-\kp}$. Since $\hal + \b = \f13 + C \e^{2-\kp}$ by \eqref{eq:beta_est}, \eqref{eq:wa_upper_lower}, combining the above two estimates
and using $\wwb(r, z) = \waa(z)$ \eqref{def:3d_Om},
we prove 
\[
\bal
 \ang z^{\hal + \kag} h(x) | \ws(\xx) - \ang z^{\b} \wwwa(z)| 
& \les  \e^{1-5\kp_1}  \ang z^{ C \e^{2-\kp}}, \\
\Longrightarrow \quad  \ang z^{ \f13 + \kag } h(x) | \ang z^{-\b} \ws(\xx) - \wwwa(z)| 
& \les  \e^{1-5\kp_1}  \ang z^{ C \e^{2-\kp}} .
\eal
\]
Dividing $\ang z^{ C \e^{2-\kp}}$ on both sides and using $\lgp |\xx| \asymp \log (|\xx| + 2)$ 
by \eqref{def:lgp}, we prove \eqref{eq:thm_global_conv}.

For any $R> 1$, using the above estimates, 
$|\ang z^{\b}- 1| \les \b \lgp z ( \ang z^{\b} + 1)$, and $|\b| \les \e$, we prove 
\[
 |\ws - \wwwa|  \les  \ang  z^{\b} | \ang z^{-\b} \ws - \wwwa |
 + |(\ang z^{\b} - 1) \wwwa| \to 0  \, ,\quad \mbox{as} \ \e \to 0, 
\]
uniformly in $B_R$, and prove the first limit in \eqref{eq:thm_est_conv}. The convergence 
of $\omsth - \f18 r^{\al} \wwwa$ in \eqref{eq:thm_est_conv} follows from the above limit, $
\bar \Om^{\th}_{*, \al} = \TTa r^{\al}  \ws$, and $\TTa \to \f18 \barcp$ by \eqref{eq:thm_SS}. 

We complete the proof of Theorem \ref{thm:self-similar}.
\end{proof}

\section{Nonlinear finite codimension stability of 3D profiles and blowup}\label{sec:blowup}

In this section, we establish nonlinear finite codimension stability of the blowup profiles 
constructed in Theorem \ref{thm:3D_solu} and prove Theorem \ref{thm:main_blowup} on
asymptotically self-similar blowup in Section \ref{sec:main_blowup_pf}.

\subsection{Dynamic rescaling and linearization}

Following \cite{chen2019finite2,ChenHou2023a}, we consider the dynamic rescaling reformulation of \eqref{eq:Euler}. Let $ \om^{\th}$ be the angular vorticity in \eqref{eq:Euler} 
and $\psi^{\th}$ be the angular stream function in the elliptic equation \eqref{eq:Euler_b} or \eqref{eq:Euler2_psi}: 
 $\psi^{\th} =  r \BSa( r^{-\al} \om^{\th})$. It is easy to show that 
\footnote{
The convention $\CCw, \CCl$ are slightly different from those used in \cite{chen2019finite2,ChenHou2023a}. Moreover, we introduce $\cpsia$ for $\om^{\th}$.
}
\bseq\label{eq:rescal}
\beq
\bal
& \om^{\th}(   \xx,  t(s) ) = \CCw(s) \cpsia  \Om^{\th}( \CCl(s)^{-1} \xx, s),  & \quad  & 
\psi^{\th}(  \xx, t(s)) = \CCw(s) \CCl(s)^2     \Psi^{\th} ( \CCl(s)^{-1} \xx, s), \\
& \Psi^{\th} = r \BSa(  \cpsia r^{-\al} \Om^{\th})
 =  r \BS( r^{-\al} \Om^{\th}), 
& 
  \quad  &  U^r  = - \tf{1}{r} \pa_z ( r \Psi^{\th}),  \quad  U^z  = \tf{1}{r} \pa_r ( r \Psi^{\th}) \, ,
\eal
\eeq
is the solution to the dynamic rescaling equation
\beq\label{eq:rescal:b}
\bal
\pa_{s} \Om^{\th}  + ( c_{l}  \xx + \UU ) \cdot \na \Om^{\th}  &=  ( \comth-
 \tf{1}{r} \pa_z \Psi^{\th} ) \Om^{\th} \, .
\eal
\eeq
In \eqref{eq:rescal}, we \emph{only} multiply the constant $\cpsia$ in $\om^{\th}$ variable, 
not both $\om^{\th}, \psi^{\th}$, to obtain the Biot-Savart law $\Psi^{\th} =  r \BS( r^{-\al} \Om^{\th})$, as in \eqref{eq:profi_3D_omth}. 
The scaling relations are given by 
\beq\label{eq:rescal:c}
\bal
  \CCw(s) = \CCw(0) \exp\big(  - \int_0^{s} \comth (\tau )  d \tau \big), \  \CCl(s) = \exp\big( \int_0^{s} -c_{l}( \tau ) d \tau    \big) ,  \   t(s) = \int_0^{ \tau } \CCw^{-1}( \tau ) d \tau , 
\eal
\eeq
\eseq
for some $\CCw(0) > 0$. As in \eqref{eq:profi_3D_rescale}-\eqref{eq:profi_3D0} and \eqref{def:Q},  we introduce the 
variables 
\beq\label{eq:dyn_rela}
 \Om \teq r^{-\al} \Om^{\th}, \quad \Psi \teq r^{-1} \Psi^{\th}, 
 \quad
   \QQ = (c_l r - r \pa_z \Psi, \, c_l z + 2 \Psi  + r \pa_r \Psi )  , 
\quad  \com(\tau) \teq  \comth(\tau)  - \al  \cdot c_{l}(\tau), 
\eeq
where $\com$ can be seen as the scaling parameter for $\om$. Then we rewrite \eqref{eq:rescal:b} equivalently as 
\bseq\label{eq:dyn_full}
\beq
\pa_{\tau} \Om   + \QQ
 \cdot \na \Om   =  ( \com - (1-\al) \pa_z \Psi ) \Om ,  \quad \Psi = \BS(\Om) \, ,
\eeq
where $\BS$ is the Biot-Savart law between $\Psi$ and $\Om$ in \eqref{eq:Euler2_psi:BS}. 
We impose the normalization conditions \eqref{eq:normal_cond} on $c_l(\tau), c_{\om}(\tau)$ 
dynamically and choose $\CCw(0)$ for the initial rescaling size as
\beq\label{eq:dyn_full:b}
   c_l + 2 \Psi_z(0) = c_{\om} - (1 - \al) \Psi_z(0) \equiv 2 , 
   \quad  \CCw(0) = \TTa,
\eeq
\eseq
where $\TTa$ is defined in \eqref{def:TTa}. 
Then the above equation is equivalent to \eqref{eq:rescal:b} and Euler \eqref{eq:Euler2}.

\vs{0.1in}
\paragraph{\bf Exact profiles and linearization}

We have constructed the exact nontrivial steady state $(\ws, \cls, \cws)$ to \eqref{eq:dyn_full} 
in \eqref{eq:solu_al}, which satisfy the properties in Proposition \ref{prop:3D_profile}. 
We decompose the solution to \eqref{eq:dyn_full} into the profile and the perturbation $(\ww, 
\tcl, \tcw)$
\beq\label{eq:dyn_decomp}
 \Om =  \ws + \ww , \quad c_l = \cls + \tcl , \quad \com = \cws +  \tcw \, ,
\eeq
and decompose other variables similarly. Linearizing equation \eqref{eq:dyn_full} around the  profile $(\ws, \cls, \cws)$ and using the notations 
$\QQs ,\cBs$ in \eqref{eq:solu_al}, 
$\td \QQ$ in \eqref{def:Q}, and $\tcB$ in \eqref{eq:lin_2Db}, we derive 
\bseq\label{eq:dyn_lin}
\beq
  \pa_s \ww + (\QQs + \td \QQ(\ww)) \cdot \na \ww = 
  (\cBs + \tcB(\ww) )  \cdot \ww +  \tcB(\ww) \cdot \ws - \td \QQ(\ww) \cdot \na \ws.
\eeq

Following the derivations in \eqref{eq:lin_2D}, \eqref{eq:psio_res1}, and \eqref{eq:normal_cond}, we obtain
\beq
\bal
  \tcB(\ww) & = -(1-\al ) \psio_z, % ( \psi_z - \psi_z(0) ), 
  & \quad \td Q^r(\ww) & =  \td c_l r - r \pa_z \psi(\ww), \qquad  \td Q^z(\ww)  = 2 \psio + r \pa_r \psi(\ww),  \\
   \tcl & = - 2 \psi_z(\ww)(0), & \quad \tcw & = (1-\al) \psi_z(\ww)(0).
\eal
\eeq
\eseq

\subsection{Linear stability estimates via semigroups}

To establish linear stability, we use the semigroup method \cite{chen2024Euler,chen2024vorticity,bedrossian2026finite}.  The linear part of equation \eqref{eq:dyn_lin} without the nonlinear terms reads 
\beq\label{eq:dyn_lin2}
    \pa_s \ww = - \QQs \cdot \na \ww +  \cBs \cdot \ww +  \tcB(\ww) \cdot \ws - \td \QQ(\ww) \cdot \na \ws \teq \cA \ww.
\eeq
We decompose the linear operator as the local part $\cT$ and the nonlocal part $\cU$
\beq\label{def:cT}
\cA = \cT + \cU , \quad 
    \cT \ww \teq - \QQs \cdot \na \ww + \cBs \cdot \ww,
    \quad \cU \ww \teq  \tcB(\ww) \cdot \ws - \td \QQ(\ww) \cdot \na \ws .
\eeq

We construct the following weights to capture the outgoing effect:
\bseq\label{def:wg_s}
\beq
  \phi(\xx) = \phim \phi_1^{\NNN_2}  , 
  \quad \phim  =  |\xx|^{-\bbb} + |\xx|^{\aln} ,
  \quad \aln = \f{1}{2} ( \f13 + \als), 
    \quad \phi_{\e}(\xx) \teq  |\xx|^{\e^3} \ang \xx^{-2 \e^3} \, ,
\eeq
where $\als$ is the exponent defined in Proposition \ref{prop:3D_profile}. The weight $\phi_{\e}(\xx)\to 0$ as $\xx \to 0$ or $\xx \to \infty$.
Recall $ \hal$ from \eqref{ran:ep2}.
From the range  of $\e, \hal$  \eqref{ran:ep_all}
and of $\als$ \eqref{eq:3D_solu_scale}, we have 
\begin{gather}
 \phi_{\e} \leq 1 , \quad     \bbb + 4 \e^3 <  1 + \al,
 \quad \bbb > 3 \e^3 + 1.1, \quad \bbb + \aln \in [1.1, 1.6],  \\
 \quad \aln - 4 \e^3 , \ \aln, \  \aln + 4 \e^3 \in [ \f13 + \f{\e}{20}, \, \f13 + \f{\e}{12}]
 \subset [ \hal - \f{\e}{8}, \, \hal - \epa - \f{\e}{500}  ] ,
 \quad 
  \phi_{\e}^3 \phi \gtr_{\NNN_2} \ang \xx^{1/3} . \notag
 \end{gather}
We also introduce the weight $\phia$ for local well-posedness in Proposition \ref{prop:LWP}
\beq\label{def:phi0}
\phia =  |\xx|^{-1} + \ang \xx^{\aln} \, .
\eeq
\eseq

We use the weight $\phi_{\e}$ in the proof of compactness in the fixed point argument. See 
Section \ref{sec:dyn_comp}.

\begin{lem}\label{lem:tran_wg}

Let $\beps_{11}$ be as in Proposition \ref{prop:3D_profile}. 
There exists $\beps_{12} \in (0, \beps_{11}]$ such that,  
for any $ \e \in ( 0, \beps_{12} ] $, there exists a weight $\phi_1 \in C^{\infty}$, with $ \phi_1 \asymp_{\e} 1$, and a parameter $\NNN_2 = N(\e)> 0$, such that the following statements hold. 
We define  $\td \QQ(\eta)$ as in \eqref{eq:dyn_lin}.  Denote 
 \beq\label{def:lam_stab}
   \QQ_{\eta} \teq \QQs + \td \QQ(\eta), \quad 
   \lam =  \bbb - 1 = 0.2> 0.
  \eeq

 For any $\eta$ with $ \nlinf{ \phim \phi_{\e} \eta } \leq \e^3$,  $i \in \{r, z\}$, and $g = \phi$ or $\phi_{\e} \phi$, we have 
 \bseq
\begin{align}
  d_0(g, \eta)(\xx) & \teq \f{ \QQ_{\eta}  \cdot  \na g }{g}
  + \cBs 
  + |\xx| \cdot |\na \cBs|
  \leq - \lam , \label{eq:damp_3D_d0} \\
  d_{i}(g, \eta)(\xx) &  \teq  d_0(g, \eta)(\xx) + \f{ \QQ_{\eta} \cdot \xx }{  |\xx|^2 }
  -  \pa_i  \QQQs^i 
  + | \pa_z  \QQQs^r |
  + | \pa_r  \QQQs^z | 
  \leq -\lam . \label{eq:damp_3D_di}
\end{align}
\eseq
\end{lem}

In the following estimates, we focus on the main terms for $|\xx| \ll 1$ and $|\xx| \gg 1$. 
In the region $ |\xx| \asymp 1$, we use the crucial outgoing condition \eqref{eq:outgo} to generate 
\emph{arbitrary} strong damping terms.

\begin{proof}

Let $\xi$ be the radial coordinate for $\xx=(r, z)$. 
For some $R_{2,\e}, R_{3,\e} >0$ to be chosen, we design a radially symmetric weight 
\beq\label{def:wg_phi1}
\bal
   \phi_1(\xx)  & = 1, \quad \forall  \, |\xx| \leq R_{2,\e}  / 2, 
 &  \quad  \phi_1(\xx) &= 1/2, \quad \forall  \, |\xx| \geq 2 R_{3,\e}, \\
   \pa_{\xi} \phi_1(\xx) & \leq 0,  \quad  \forall \, \xx \in \R_+ \times \R,  
   & \quad  \pa_{\xi} \phi_1(\xx) & < 0, \quad \quad \forall \, |\xx| \in [R_{2,\e}, R_{3,\e} ].
\eal
\eeq

We determine the parameter and weight in the following order 
\[
  \beps_{12} \rsa R_{2,\e}, R_{3,\e} \rsa \phi_1 \rsa \NNN_2. 
\]

\paragraph{\bf Asymptotics of coefficients for $|\xx|\ll 1$ and $|\xx| \gg 1$ }

We first derive the asymptotics of coefficients $\QQ, \cB$.  Let $\psi  = \BS(\eta)$, $\psis = \BS(\ws)$ \eqref{eq:Euler2_psi}. Recall $\cBs$ from \eqref{eq:solu_al},  and $\td c_{\om}, \td c_l$ from \eqref{eq:dyn_lin}
\[
  \cBs 
   = \cws  - (1-\al ) \pa_z \psis ,
   \quad \td c_{\om} = (1 - \al) \psi_z(\eta)(0), \quad \td c_l = - 2 \psi_z(\eta)(0) .
\]

Using \eqref{eq:vel_iso:2nd} and $\nlinf{\phim \phi_{\e} \eta} \leq \e^3 $ we bound 
\bseq\label{eq:QB_asymp_1}
\beq\label{eq:QB_asymp_1:scale}
|\pa_z \psi(\eta)(0)| +  |\td c_{\om} | + |\td c_l | \les |\psi_z(\eta)(0) | \les \e^{-1} \nlinf{\phi \phi_{\e} \eta}
 \les \e^2.
\eeq

Note that $\cBs(0) = 2$ by \eqref{eq:normal_cond} and the weight $\phim, \phi_{\e}$ in  \eqref{def:wg_s} satisfies 
\[
  |\xx|^{\al} \les |\xx|^{-\bbb + \e^3} + |\xx|^{\aln -\e^3} \asymp  \phim \phi_{\e}, \quad 
   \aln - \e^3 \in [\hal -\e, \hal - \epa] . 
 \]
Using Proposition \ref{prop:iso_est} with $\bbu= \bbb - \e^3, \g = \aln - \e^3$ 
 and $\nlinf{\phim \phi_{\e} \eta} \leq \e^3$ and 
 \eqref{eq:3D_solu_decay},  we obtain
\beq\label{eq:QB_asymp_1:a}
\bga
  |\cBs  - 2|  
    \les  \upse  |\xx|^{1 + \al},   \quad 
 |\cBs  - \cws  |   \les \upse  \ang \xx^{ \al - \hal + \epa }, \\
    | \QQs  -  \cls  \xx |  \les  \e^{-1}  |\xx| \ang \xx^{\al- \hal + \epa } ,
   \quad  | \td \QQ(\eta)|  \les \e^{-1} |\xx| \cdot \nlinf{\phi \phi_{\e} \eta} 
  \les \e^2 |\xx|, \quad  |\QQ_{\eta} | \les \e^{-1} |\xx| ,
\ega
\eeq
where $\al-\hal + \epa < 0$ by \eqref{ran:ep_all} and  $\upse $ is the constant appearing in \eqref{eq:3D_solu_decay}. Using the outgoing condition for $\QQs$ \eqref{eq:outgo}, the above estimate for $\td \QQ(\eta)$, and by requiring $\e$ small enough, we bound 
\beq\label{eq:outgo_Q:recall}
  \QQ_{\eta}(\xx) \cdot \xx= ( \QQs + \td \QQ(\eta) ) \cdot \xx \geq ( \lam_Q - C \e^2) |\xx|^2 
  \geq \tf12 \lam_Q |\xx|^2.
\eeq

Using \eqref{eq:3D_solu_Q0} for $\pa_i \QQQs^i(0), i = r, z$ 
\eqref{eq:3D_solu_decay:d} for $\na (\psis - \pa_z \psis(0) \cdot z) $,  \eqref{eq:QB_asymp_1:scale},
and  \eqref{eq:QB_asymp_1:a}, we prove
\beq\label{eq:QB_asymp_1:b}
\bga
 | \QQ_{\eta}^r - q^r(0) r | +  | \QQ_{\eta}^z - q^z(0) z | \les \upse |\xx|^{2 + \al } 
 + \e^2 |\xx| , \\
  q^r(0) = \pa_r \QQ_{\eta}^r(0) %\geq C \e^{-1} - C_2 \e 
  \geq 4 - C \e^2 , \quad 
q^z(0) =  \pa_z \QQ_{\eta}^z(0) = 2 .
\ega
\eeq

\eseq

Applying \eqref{eq:3D_solu_decay} for $\psis$,  
we further obtain
\begin{align}\label{eq:QB_asymp_2}
| \pa_r \QQQs^z |  + |\pa_z \QQQs^r|   + |\xx| \cdot |\na \cBs| 
   \les  | r \na^2 \psis  | 
+ | \pa_r  \psis   | + |\xx| \cdot | \na \pa_z \psis | 
    & \les \upse \min( |\xx|^{1 + \al}, \ang \xx^{\al-\hal + \epa} ) ,  \notag \\
|\pa_r \QQQs^r - \pa_r \QQQs^r(0)| 
+ |\pa_z \QQQs^z - \pa_z \QQQs^z(0)| & \les \upse |\xx|^{\al + 1}, \notag \\
|\pa_r \QQQs^r - \cls | 
+ |\pa_z \QQQs^z - \cls| & \les \upse \ang \xx^{\al-\hal + \epa} .
\end{align}

\paragraph{\bf Decompose $d_0, d_i$}
Below, we simplify $ \tcB( \eta)$ as $\tcB$,  $d_0(g, \eta)(\xx)$ as $d_0(g)(\xx)$, and $d_i(g, \eta)(\xx)$ as $d_i(g)(\xx)$. 
Recall $\phi, \phi_{\e}$ from \eqref{def:wg_s}. 
Denote $\xi = |\xx|$. 
Since $\phi, \phi_{\e}$ is radially symmetric in $r, z$, 
for $g =  \phi \phi_{\e}^j
= \phim \phi_1^{\NNN_2} \phi_{\e}^j$ defined in \eqref{def:wg_s} with $j = 0, 1$, we expand 
\bseq\label{eq:damp_3D:d0_pf1}
\begin{align}
   \f{\QQ_{\eta} \cdot \xx }{|\xx|^2} & \cdot \f{ \xi \pa_{\xi} g}{g} 
  = \f{\QQ_{\eta} \cdot \xx }{|\xx|^2} \cdot \B( \f{ \xi \pa_{\xi}  \phim }{ \phim } +  
  \NNN_2 \f{ \xi \pa_{\xi} \phi_1 }{\phi_1} +  j \f{\xi \pa_{\xi} \phi_{\e} }{\phi_{\e}}  \B) \\
&  = \f{\QQ_{\eta} \cdot \xx }{|\xx|^2} \cdot \B( \f{    -\bbb \xi^{-\bbb} +\aln \xi^{\aln}  }{ \xi^{-\bbb} + \xi^{\aln}  } \B) +  
\f{\QQ_{\eta} \cdot \xx }{|\xx|^2} \cdot \NNN_2 \f{ \xi \pa_{\xi} \phi_1 }{\phi_1}  + \f{\QQ_{\eta} \cdot \xx }{|\xx|^2} \cdot  j \f{\xi \pa_{\xi} \phi_{\e} }{\phi_{\e}}  
\teq I_M + I_{\phi_1} + I_{\e} , \notag 
\end{align}
where $I_M$ is the main term from the weights. We then decompose $d_0$ \eqref{eq:damp_3D_d0} as 
\beq
\bal
  d_0(g)  
  =  \f{\QQ_{\eta} \cdot \xx }{|\xx|^2} \cdot \f{ \xi \pa_{\xi} g}{g} +\cBs 
  +  |\xx| \cdot |\na  \cBs| \teq I_M + I_{\phi_1} + I_{\e} + \cBs+   |\xx| \cdot |\na \cBs| \, .
\eal
\eeq
\eseq

For $d_i$ in \eqref{eq:damp_3D_di}, we denote 
\[
  I_{\cR} =    | \pa_z  \QQQs^r |
  + | \pa_r  \QQQs^z |  + |\xx| \cdot |\na  \cBs |.
\]

Using the above derivation of $d_0$, we derive
\beq\label{eq:damp_3D:di_pf1}
\bal
  d_i(g) & =  d_0(g) + \f{ \QQ_{\eta} \cdot \xx }{  |\xx|^2 }
  -  \pa_i  \QQQs^i  +   | \pa_z  \QQQs^r |
  + | \pa_r  \QQQs^z |
 \teq \big( I_{M} + \f{ \QQ_{\eta} \cdot \xx }{  |\xx|^2 } - \pa_i \QQs^i \big) + I_{\phi_1} + I_{\e} +\cBs+  I_{\cR}.
\eal 
\eeq

Recall $\phi_{\e}$ from \eqref{def:wg_s}. For $i\in \{0, 1\}$, using \eqref{eq:QB_asymp_1:a}, we obtain
\beq\label{eq:damp_Ie}
 |I_{\e}| \leq |\f{\QQ_{\eta} \cdot \xx }{|\xx|^2} | \cdot \B| \f{\xi \pa_{\xi} \phi_{\e} }{ \phi_{\e} } \B|
\les \e^{-1} \cdot  \B| \f{\xi \pa_{\xi} ( |\xi|^{\e^3} \ang \xi^{-2 \e^3 } )}{ |\xi|^{\e^3} \ang \xi^{-2 \e^3 }  }  \B| \les \e^2 ( \f{\xi \pa_{\xi} \xi }{\xi} + \f{\xi \pa_{\xi} \ang \xi^2}{\ang \xi^2} )
\les \e^2.
\eeq

For $I_{\phi_1}$, since $\phi_1$ is radially symmetric in $(r, z)$ and decreasing in $|\xx|$, using \eqref{eq:outgo_Q:recall}, we obtain
\beq\label{eq:damp_Iphi1}
  I_{\phi_1} \leq \f12 \lam_Q \cdot \NNN_2 \f{\xi \pa_{\xi} \phi_1}{\phi_1}  \leq 0.
\eeq

For $I_{\cR}$, using \eqref{eq:QB_asymp_2}, we bound 
\beq\label{eq:damp_IR}
   |I_{\cR} | \les \upse \min( |\xx|^{1 + \al}, \ang \xx^{\al -\hal + \epa} ) \, .
\eeq 

Next, we use \eqref{eq:QB_asymp_1}, \eqref{eq:QB_asymp_2} to derive the main terms for $d_0, d_i$ near $\xx= 0, \infty$. 

\paragraph{\bf Estimate near $0$}

Using \eqref{eq:QB_asymp_1} for $\cBs $, \eqref{eq:QB_asymp_2} for $\na \cBs$,  we estimate 
\[
  d_0 =  
  \f{\QQ_{\eta} \cdot \xx }{|\xx|^2}  \cdot ( -\bbb + O(\xi^{\aln + \bbb}) )  
   + 2  + I_{\phi_1} + I_{\e}  +  O( \upse |\xx|^{1 + \al} ) .
\]

Using \eqref{eq:QB_asymp_1:b}, we obtain 
\beq\label{eq:damp_d0:pf2}
  \f{\QQ_{\eta} \cdot \xx }{|\xx|^2} \geq \f{q^r(0) r^2 + 2 z^2}{ |\xx|^2 } - C \upse |\xx|^{1 + \al}
  - C \e^2
  \geq 2 - C \upse |\xx|^{1 + \al} - C \e^2,
  \quad  \B|  \f{\QQ_{\eta} \cdot \xx }{|\xx|^2} \B| \les \e^{-1}.
\eeq

Combining the above estimates,  using $\aln + \bbb \in (1,2)$, 
we bound $d_0$ from above 
\beq\label{eq:damp_d0:near0}
  d_0 \leq - 2 \bbb + 2 + C  \upse ( |\xx| + |\xx|^2) + C \e^2 + I_{\phi_1} + I_{\e}.
\eeq

For $d_i$ \eqref{eq:damp_3D:di_pf1}, similarly, using  \eqref{eq:QB_asymp_1}, 
we obtain
\[
  d_i =   \f{\QQ_{\eta} \cdot \xx }{|\xx|^2}  \cdot ( -\bbb + 1 + O(\xi^{\aln + \bbb}) )  
   + 2 - \pa_i \QQQs^i(0) + O( \upse ( |\xx| + |\xx|^2) )  + I_{\phi_1} + I_{\e} + I_{\cR} \, .
\]

Since $\pa_i \QQQs^i(0) \geq 2 - C \e^2$ \eqref{eq:QB_asymp_1:b}, $\bbb > 1$, 
and $\bbb + \aln < 2$ by \eqref{def:wg_s}, using \eqref{eq:damp_d0:pf2}, we prove 
\beq\label{eq:damp_di:near0}
  d_i \leq  2 (1 - \bbb) + C \e^2 + C \upse ( |\xx| + |\xx|^2)+ I_{\phi_1} + I_{\e} + I_{\cR}.
\eeq

\paragraph{\bf Estimate near $\infty$}

Applying \eqref{eq:QB_asymp_1} for $\QQ_{\eta},  \td c_l, \cBs $,
\eqref{eq:QB_asymp_2} for $\na  \cBs$, we estimate $d_0$ in \eqref{eq:damp_3D:d0_pf1}:
\[
\bal
  d_0  = % ( \cls + \td \cls +  )
   \f{\QQ_{\eta} \cdot \xx }{|\xx|^2}  \cdot ( \aln + O(\xi^{ -\bbb - \aln }) )
   + \cws  + O( \upse \ang \xx^{\al-\hal + \epa} )    + I_{\phi_1} + I_{\e}  .
  \eal
\]

Using estimates \eqref{eq:QB_asymp_1} for $\QQ_{\eta}$ 
and \eqref{eq:QB_asymp_1:scale}, we bound 
\[
   \f{\QQ_{\eta} \cdot \xx }{|\xx|^2}
   = \cls + \td c_l +  O( \upse \ang \xx^{\al-\hal + \epa })
   =  \cls +  O( \upse \ang \xx^{\al-\hal + \epa } + \e^2 ) ,
   \quad |\td c_{\om}| \les \e^2.
\]

Since $\xi = |\xx|, \bbb +\aln \in (1,2)$ by \eqref{def:wg_s},  combining the above estimates, we prove 
\beq\label{eq:damp_d0:nearinf}
  d_0 \leq  \cls  \aln +\cws +  C \upse  ( \ang \xx^{\al-\hal + \epa } + |\xx|^{-2}) + C \e^2
     +    I_{\phi_1} + I_{\e} .
\eeq

For $d_i$ \eqref{eq:damp_3D:di_pf1}, similarly, using  \eqref{eq:QB_asymp_1} 
and $|\pa_i \QQs^i - \cls | \les \upse \ang \xx^{\al-\hal + \epa}$ \eqref{eq:QB_asymp_2},  and the above bound for $\QQ_{\eta}$ and \eqref{eq:QB_asymp_1:scale} for $\td c_l, \td c_{\om}$, we obtain
\begin{align}\label{eq:damp_di:nearinf}
  d_i  & \leq  \f{\QQ_{\eta} \cdot \xx }{|\xx|^2}  \cdot ( \aln  + 1 + O( |\xx|^{-\aln - \bbb}) )  
  + \cws 
  - \cls  +  O( \upse \ang {\xx}^{\al-\hal + \epa} )  + I_{\phi_1} + I_{\e} + I_{\cR} \notag \\
& \leq \cls (\aln  + 1 )
+ \cws - \cls
+ C   \upse  ( \ang \xx^{\al-\hal + \epa} + |\xx|^{-2})  +C  \e^2
 + I_{\phi_1} + I_{\e} + I_{\cR} \notag \\
& \leq \cls \aln  + \cws 
+ C   \upse  ( \ang \xx^{\al-\hal + \epa} + |\xx|^{-2})+  C \e^2  + I_{\phi_1} + I_{\e} + I_{\cR} \, .
\end{align}

Combining estimates \eqref{eq:damp_d0:near0} and \eqref{eq:damp_d0:nearinf} for $d_0$,
and \eqref{eq:damp_di:near0} and \eqref{eq:damp_di:nearinf} for $d_i$, and using \eqref{eq:damp_Ie}
for $I_{\e}$, \eqref{eq:damp_IR} for $I_{\cR}$, we bound 
\beq\label{eq:damp_d0i:pf3}
\bal
 \max( d_0, d_i)  & \leq \min( -2 \bbb + 2  + \bar C \upse (|\xx| + |\xx|^2), \, \cls \aln + \cws + \bar C \upse 
  ( \ang \xx^{\al-\hal + \epa}  +|\xx|^{-2}  ) ) + \bar C \e^2 +I_{\phi_1}, \\
\eal
\eeq
for some absolute constant $\bar C$.  Here, we have combined the bound for $I_{\cR}$ in \eqref{eq:damp_IR} with 
$ C \upse (|\xx| + |\xx|^2)$ or $C \upse ( \ang \xx^{\al-\hal + \epa}  +|\xx|^{-2}  )$.

\vs{0.1in}
\paragraph{\bf Choosing $\phi_1, \NNN_{\e}$}

By definition of $\aln$ in \eqref{def:wg_s}, using Proposition \ref{prop:3D_profile} for $\als$ and $\cls$, we obtain 
\[
   \cls \aln + \cws = \cls( \aln - \als) = \f{1}{2} \cls ( \f{1}{3} - \als)
   = \f12 (\f{64}{9} \e^{-1} + O( \e^{-\hk})) \cdot ( - \f{\e}{8} + O( \e^{1 +\mhk})
   = -\f{4}{9}  + O( \e^{\mhk}).
\]

Recall $\bbb \in (1, 1.2]$  \eqref{norm:Xc} and $\lam = \bbb - 1 \leq 0.2$ by \eqref{def:lam_stab}.  We require $\e$ small enough so that 
\[
\bga
  -2 \bbb + 2 + \bar C \e^2 \leq \tf32( - \bbb + 1) =   - \tf32 \lam,
  \quad \cls \aln + \cws + \bar C \e^2 \leq - \tf13
  \leq - \tf32 \lam \, . 
\ega
\]

By continuity, there exists 
\[
  R_{2, \e} \ll 1 \ll R_{3, \e},
\]
such that 
\[
\bal
  -2 \bbb + 2  + \bar C \upse (|\xx| + |\xx|^2) & \leq -  \lam, \quad  \forall |\xx| \leq R_{2 ,\e}, \\
 \cls \aln + \cws +  \bar C \upse  ( \ang \xx^{\al-\hal + \epa}  +|\xx|^{-2}  )
 & \leq  -  \lam, \quad  \forall |\xx| \geq R_{3, \e}.
\eal 
\]
For $ |\xx| \in [R_{2,\e}, R_{3,\e}]$, we bound $d_0, d_i$ except 
for $I_{\phi_1}$-term  in \eqref{eq:damp_d0i:pf3}  by some constant \emph{independent} on $\phi_1$.  Combining the above estimates, we prove 
\beq
  \max(d_0, d_i) \leq -  \lam + (C_{ \eqref{eq:damp_pf:C1},\e} + 1)  \one_{ [R_{2,\e}, R_{3,\e}] }(|\xx| ) + I_{\phi_1} \, ,
  \label{eq:damp_pf:C1}
\eeq
for some  constant $C_{ \eqref{eq:damp_pf:C1},\e} > 0$ independent of $\phi_1$. We choose the above parameters $R_{2, \e}, R_{3,\e}$ for the weight $\phi_1$ \eqref{def:wg_phi1},
and construct $\phi_1$ with properties in  \eqref{def:wg_phi1}. From \eqref{def:wg_phi1}, we obtain
\beq
   \phi_1^{-1} \xi \pa_{\xi} \phi_1(\xx) \leq - C_{ \eqref{eq:damp_pf:C2}, \e} \one_{ [R_{2,\e}, R_{3, \e}] }(|\xx|).
   \label{eq:damp_pf:C2}
\eeq
for some $C_{ \eqref{eq:damp_pf:C2}, \e} > 0$.  We choose $\NNN_{2}$ as:
 $ 
 \NNN_{2} = \f{ 4( C_{ \eqref{eq:damp_pf:C1}, \e} + 1 ) }{ \lam_Q C_{ \eqref{eq:damp_pf:C2}, \e} } .
 $

Combining \eqref{eq:damp_IR}, \eqref{eq:damp_Iphi1} for $I_{\phi_1}$ and the above estimates, we prove 
\[
    \max(d_0, d_i)  \leq -  \lam_{\e} +( C_{ \eqref{eq:damp_pf:C1},\e} +1 )   \one_{ [R_{2,\e}, R_{3,\e}] }(|\xx| ) 
    - \tf12 \NNN_{2} \lam_Q C_{ \eqref{eq:damp_pf:C2},\e}   \one_{ [R_{2,\e}, R_{3, \e}] }(|\xx|) \leq -   \lam_{\e}.
\]
for any $\xx$. We prove \eqref{eq:damp_3D_d0}, \eqref{eq:damp_3D_di}.
\end{proof}

\subsubsection{Functional spaces}

Recall $\phi$ from \eqref{def:wg_s}. We consider the following functional space 
\bseq\label{def:cZ}
\beq
  \cZ \teq \{ f(r, z):  f \phi, \ \phi |\xx| \cdot \na f \in \cC_0(\R_+ \times \R) , 
  \ f(r, z) \mw{ \ is \ odd \ in  \ } z \} ,
\eeq
equipped with norm
\beq
  \| f \|_{\cZ} \teq  \nlinf{ \phi ( |f|^2 + |\xx|^2 |\na f|^2 )^{1/2} } , 
\eeq
where $\cC_0(\R_+ \times \R)$ denotes the space of continuous functions with 
\footnote{
To obtain strong continuity of the transport semigroup, it is important to consider functions that vanish for large \( |\xx| \) and for \( |\xx|\ll 1 \). 
See Lemma \ref{lem:tran_semi} and its proof. For a continuous function \(f\) that does not vanish at infinity, one need not have$\|f(\cdot+h)-f(\cdot)\|_{L^\infty}\to 0$ as  $h\to 0$. For example, $f(x) = \sin(x^2)$.
}
\beq\label{def:cC0}
  \cC_0 \teq \{ g:    \lim_{|\xx| \to \infty} g( \xx) = 0, \   \lim_{|\xx| \to 0 } g( \xx ) = 0 \} .
  \eeq
\eseq
We require vanishing near $\xx=0$ since $\phi$ is  singular  at $\xx=0$. Clearly, $\cZ$ is a Banach space. 

To apply spectral theory, we complexify the real Banach space $\cZ$ to the complex Banach space. 
\bseq
\beq
  \cZc \teq \{ f + i g : f, g \in \cZ \} \, .
\eeq
For any $z \in \cZc$, we define the associated norm
\beq
    \| z \|_{\cZc}  \teq  \sup\nolim_{\th \in \R} (  \cZZ{ \Re( e^{i \th} z ) }^2
    + \cZZ{ \Im( e^{ i \th } z ) }^2  )^{1/2} \, , 
\eeq
where $\Re(z), \Im(z)$ denote the real and imaginary part, respectively.  
\eseq

For any real bounded linear operator $\cL : \cZ \to \cZ$, we 
extend it to a complex bounded linear operator
$\cZc \to \cZc $: 
\beq\label{eq:op_lin}
  \cL (f + i g) = \cL f + i \cL g , \quad \forall f, \ g \in \cZ.
\eeq
It is not difficult to show that $\cL$  is bounded %complex linear operator 
and has the same operator norm
\footnote{
See \cite[Section 5.1.1]{BuDi2018} for some properties of complex Banach spaces and the complexification of a real Banach space.
}
\beq\label{eq:op_norm}
   \| \cL \|_{\cZc \to \cZc} = \| \cL \|_{\cZ \to \cZ}. 
\eeq

\subsubsection{Semigroup on the transport operator $\cT$}

We define the domain of $\cT$ as 
\beq\label{def:dom_cT}
    D(\cT) \teq \{ f \in \cZc: \QQs \cdot \na f \in \cZc \}.
\eeq

Clearly, \((\QQs\cdot\nabla):D(\cT)\subset\cZc\to\cZc\) is a closed operator. Moreover,
\(D(\cT)\) is dense in \(\cZc\), by a standard cutoff and mollification argument; see
Section \ref{sec:dyn_regu} for the cutoff argument. In the following two subsections, we establish the following results on the semigroup.

\begin{prop}\label{prop:semi}

Let $\lam$ be defined in \eqref{def:lam_stab}, and let $\cT,\cU$ be the operators in \eqref{def:cT}. The operator $ \cU$ is a bounded, linear, compact operator from $\cZc \to \cZc$. 
The operators $\cT , \cA = \cT + \cU  : D(\cT) \subset  \cZc \to \cZc$ generate strongly continuous semigroups
\[
   e^{\cT s} : \cZc \to \cZc, \quad e^{ \cA s } : \cZc \to \cZc , 
\]
with  the following decay estimates
\beq\label{eq:dissp_semi}
 ||  e^{s \cT} ||_{\cZc} \leq e^{-\lam s} , \quad 
 || e^{s \cA} ||_{\cZc} \leq  e^{C_{\e} s} .
\eeq
Moreover, we have the following spectral property of $\cT$
\beq\label{eq:spec}
 \{ z: \Re(z) > -\lam \} \subset \rho_{\mw{res}}(\cT) ,
\eeq
where $\rho_{\mw{res}}( \cL )$ denotes the resolvent set of a linear operator  $\cL$. 
\end{prop}

Let $\bar X$ be the characteristics associated with $\QQs$
\beq\label{eq:traj_s}
   \tf{d}{ds} \bar X(s, \xx) = \QQs( \XXb(s, \xx) ), \quad \XXb(0, \xx) = \xx. 
\eeq

Using estimate \eqref{eq:3D_solu_decay} and \eqref{eq:outgo}, for $s_1 < s_2$, we obtain
\bseq\label{eq:traj_ratio2}
\beq
\f{d}{ds} | \XXb(s, \xx)|^2 \leq C_{\e} |\XXb(s, \xx)|^2, \
\f{d}{ds} | \XXb(s, \xx)|^2 \geq  C | \XXb(s, \xx)|^2,
\  \Rightarrow \ 
1 \leq \f{ |\XXb(s_2, \xx )|}{ | \XXb(s_1, \xx)| } \leq e^{ C_{\e} (s_2 - s_1) }.
\eeq
Using the ODE for $\na \bar X$ similar to \eqref{eq:3D_ODE:X_C1} and \eqref{eq:3D_solu_decay}, we obtain
\beq\label{eq:traj_cont}
   |\na \bar X(s, \xx) - \Id| \les_{\e} s, \quad \forall |s| \leq 1.
\eeq 
\eseq

Using the regularity of $\QQs, \cBs$ from \eqref{eq:3D_solu_reg}
and \eqref{eq:3D_solu_decay} and method of characteristics , given initial data $\ww_0 \in \cZ$, 
we uniquely solve  the transport equation
\bseq\label{eq:tran_semi0}
\beq
   \pa_{s} \om + \QQs \cdot \na \om = 0,
\eeq
with solution operator
\beq
     (\MSS_{s} \ww_0)(\xx) \teq \ww(s, \xx ) = \ww_0( \bar X(-s, \xx ) ).
\eeq
\eseq

Below, we establish strong continuity and boundedness of the transport operator $\MSS_s$.

\begin{lem}[\bf Strong continuity]\label{lem:tran_semi}
For any $f \in \cZ$ and $|s| \leq 1$, we have $ f( \bar X( \cdot,s) ) \in \cZ $ and 
\beq\label{eq:tran_semi1}
  \lim_{s \to 0} \cZZ{f( \bar X( s, \cdot ) ) - f} =0,
  \quad  \cZZ{f(\bar X(s, \cdot))} \leq  C_{\e} \cZZ{f}.
\eeq
\end{lem}

\begin{proof}
It suffices to consider $|s| \leq 1 $.  Since $f \in \cZ$ implies that $f \in C^1$, and $ \QQs \in C^{1,\al}$, 
following Section \ref{sec:reg:traj_C1a}, we obtain $\XXb(s, \cdot ) \in C^{1,\al}$ and thus 
$f \cc \XXb(s, \cdot) \in C^1$. 

For $i ,j \in \{ r, z \}$, using \eqref{eq:traj_s} and a direct calculation, we yield 
\beq\label{eq:tran_cont0}
\bal
   \pa_i ( f \cc \bar X(s, \xx) )
   & = \pa_i \bar X(s, \xx) \cdot (\na f) \cc \bar X(s, \xx) \\
   & = ( \pa_i \bar X^j(s, \xx) - \d_{ij} ) \cdot (\pa_j f)  \cc \bar X(s, \xx)  +  (\pa_i f) 
   \cc \bar X(s, \xx) \teq I + II.
  \eal
\eeq

Since $|X(s, \xx)| \asymp_{\e} |\xx|$ by \eqref{eq:traj_ratio2} and the weight $\phi$ defined in \eqref{def:wg_s} satisfies 
\beq\label{eq:phi_double}
  \phi(\xx) \asymp_m \phi(\yy), \quad  \forall \, |\xx| / |\yy| \in [1/m, m],
\eeq
using \eqref{eq:traj_cont}, for $|s|\leq 1$, we obtain
\beq\label{eq:tran_cont1}
  |\phi(\xx) |\xx| \cdot I| \les  s \phi(\yy) | \yy|\cdot  | \pa_i f(\yy)| \,_{\yy = \bar X(s, \xx)  } 
  \les s \cZZ{f}.
\eeq

For $II$, for any $ h \in (0, 1), |s| \leq 1$, and $|\xx| \in [h, 1/h]$, from \eqref{eq:traj_ratio2}, and $\bar X(0, s) = 0$, we have 
\[
  C_{1,\e} h \leq |\bar X(s, \xx) | \leq C_{2,\e} h^{-1},
  \quad |\bar X(s, \xx) - \xx| \les  s |\xx| \les h^{-1} s.
\]

Since $\pa_i f( \yy )$ is uniformly continuous for $|\yy| \in [C_{1,\e} h, C_{2,\e} h^{-1}]$, 
using the above estimate, we obtain
\[
   \lim\nolim_{s \to 0} \sup\nolim_{|\xx| \in [h, h^{-1}]} \phi |\xx| \cdot |  (\pa_i f )\cc \bar X(s, \xx) - (\pa_i f )(\xx) | = 0.
\]

Since $|X(s, \xx)| \asymp_{\e} |\xx|$ by \eqref{eq:traj_ratio2}, 
using \eqref{eq:phi_double}, for any $|s| \leq 1$, we bound 
\[
\bal
  \sup_{|\xx| \notin [h, h^{-1}]} \phi(\xx) \cdot |\xx| \cdot |  (\pa_i f )\cc \bar X(s, \xx) - (\pa_i f )(\xx) |  
  \les \sup_{|\xx| \notin [ C_{\e} h, (C_{\e} h)^{-1} ] }  \phi(\xx) |\xx| |\na f(\xx) |
\eal
\]
with $C_{\e}$ uniformly in $h, s$.  From the definition of $\cZ$ \eqref{def:cZ}, we have  $\phi |\xx| \pa_i f(\xx) \to 0$ as $|\xx| \to 0$ and $|\xx| \to \infty$.
Combining the above two estimates and \eqref{eq:tran_cont1}, 
and taking $h\to 0$, we prove
\[
    \lim_{s \to 0} 
     \nlinf{ \phi |\xx| \na ( f(\bar X(s, \cdot )) -  f )} = 0.
\]

The proof of $\lim_{s \to 0} \nlinf{ \phi   ( f(\bar X(s, \cdot )) - f )}  = 0$ is similar and is easier, and thus it is omitted. We prove the first estimate in \eqref{eq:tran_semi1}.

Using $|\bar X(s, \xx)| \asymp_\e |\xx|$, $|\na \bar X| \les_{\e} 1$ for any $|s| \leq 1$,  \eqref{eq:tran_cont0}, and \eqref{eq:phi_double}, we prove 
\[
 \phi(\xx) ( |f(\bar X(s, \xx) )   | + |\xx| \cdot |\na f(\bar X(s, \xx)) | ) 
 \les \phi(\yy) ( |f(\yy)| + |\yy| \cdot | \na f(\yy)| ) \, |_{y =\bar X(s, \xx) }
 \les \cZZ{f}. 
\]
Taking supremum over $\xx$, we prove the second estimate in \eqref{eq:tran_semi1}.
\end{proof}

The uniqueness of solution to \eqref{eq:tran_semi0} implies that $\MSS_s$ form a semigroup. Lemma \ref{lem:tran_semi} implies that $\MSS_s$ is a bounded linear operator from $\cZ \to \cZ$, and is strongly continuous. Next, we verify that the generator of $\MSS_s$ is  $-\QQs \cdot \na$ with domain \eqref{def:dom_cT}. For $\ww_0 \in \cZ$, which is locally $C^1$, 
we define $\ww(s, \xx) = \MSS_s \ww_0 = \ww_0(\XXb(-s, \xx))$ as in \eqref{eq:tran_semi1}. 
Since $\ww_0$ is locally $C^1$, a direct calculation yields 
\beq\label{eq:tran_gen1}
  \pa_{s} \ww(s, \xx ) = -(\QQs \cdot \na \ww_0) \cc \bar X(-s, \xx )
\quad \Longrightarrow  \quad \lim_{s \to 0} \f{ \MSS_s \ww_0 - \ww_0 }{s} =  - (\QQs \cdot \na \ww_0) (\xx).
\eeq
The above limit denotes the pointwise limit in $\xx, s$. %rather than the limit in the space $\cZ$. 
Thus, for the limit to exist in $\cZ$,  a necessary condition is $\QQs \cdot \na \ww_0 \in \cZ \iff \ww_0 \in D(\cT)$, where $D(\cT)$ is defined in \eqref{def:dom_cT}. 

On the other hand, for $\ww_0 \in D(\cT)$, 
we obtain  $f = - \QQs \cdot \na \ww_0 \in \cZ$.  Using \eqref{eq:tran_gen1} and Lemma \ref{lem:tran_semi}
with $  f = - \QQs \cdot \na \ww_0$, we obtain
\[
  \lim_{s \to 0} \|  \pa_{s} \ww(s, \xx ) + (\QQs \cdot \na \ww_0)(\xx) \|_{\cZ} = 0,
  \quad \Longrightarrow \ \B| \B| \lim_{s \to 0}  \int \f{1}{s} \int_0^{s} \pa_{s} \ww  + 
   (\QQs \cdot \na \ww_0)(\xx)  \B|\B|_{\cZ} = 0.
\]
Thus, the limit in \eqref{eq:tran_gen1} holds, and we prove that $\MSS_s$ is a strongly continuous semigroup with generator $- \QQs \cdot \na$ and a domain $D(\cT)$ and conclude $\MSS_s = e^{\cT s}$. 

\subsubsection{Bounded and compact perturbation}

In this subsection, we prove that $\cBs  \ww$ is a bounded linear operator  $\cZ \to \cZ$  and 
$\cU \ww : \cZ \to \cZ$ is a compact operator.

Recall $\cBs = \cws - (1-\al) \pa_z \psis$ from \eqref{eq:solu_al}. Since $\cBs \in C^{1,\al}$,
for $\ww_0 \in \cZ$, using \eqref{eq:3D_solu_decay} and Leibniz rule, we obtain 
\beq\label{eq:stab_comp_cB}
\cZZ{ \cBs \ww }
\les \nlinf{\cBs + |\xx| \cdot |\na \cBs| \, } \cZZ{\ww}
\les_{\e} \cZZ{\ww}.
\eeq
Thus, the multiplication $\cM \ww \teq    \cBs \ww$,  $ \cM : \cZ  \to  \cZ$ is a bounded linear operator.

Let $\psi = \BS(\ww)$ \eqref{eq:Euler2_psi}.  Recall the linear operator $\cU$ from \eqref{def:cT} and $\td \QQ(\ww)$ from \eqref{eq:dyn_lin}
\beq\label{eq:stab_comp_cU}
\bal
 \td \QQ(\ww) &= ( (\tcl -\pa_z \psi) r, \, 2 \psio + r \pa_r \psi  ),  
 \qquad \tcB(\ww) =  - (1-\al) \psio_z , \\
       \cU (\ww) & \teq  \tcB(\ww) \cdot \ws - \td \QQ(\ww) \cdot \na \ws 
       =  - (1-\al) \psio_z \ws -
       (\tcl -\pa_z \psi) r \pa_r \ws 
       - ( 2 \psio + r \pa_r \psi ) \pa_z \ws.
\eal 
\eeq

Recall the weight and exponents from \eqref{def:wg_s}
\[
\phi \asymp_{\e} |\xx|^{-\bbb} + |\xx|^{\aln}, 
\quad \phi \phi_{\e}^2 \asymp_{\e} |\xx|^{-\bbb +2 \e^3} + \ang \xx^{ \aln - 2 \e^3 },
\quad \aln,  \aln - 2\e^3\in [ \tf13 +\tf{\e}{20}, \hal]
\subset [ \hal - \tf{\e}{8}, \hal] .
\]
Thus, the asymptotic of weights $\phi, \phi \phi_{\e}^2$ satisfy assumptions in \eqref{eq:iso_est:wg1}. Using Proposition \ref{prop:iso_est} with  weights $\phi, \phi \phi_{\e}^2$, 
the definition of $\cZ$-norm \eqref{def:cZ},
and  $\psi_z(0) = 2\al\JJ(\infty)$ \eqref{eq:vel_J}, we obtain
\bseq\label{eq:stab_comp:psi1}
\begin{align}
|\xx| \phi  \cdot \big(  | \na \psi_z| + r \ang \xx^{-1} |\pa_{rr} \psi| \big) 
+
\phi  \big(  | \tf{1}{z} \psio|  +   |\na \psio| \big) & \les_{\e}  |\xx|^{\al}  \cdot \nlinf{ \phi( |\ww| + |\xx| \cdot |\na \ww| ) }\les  \cZZ{\ww}, \label{eq:stab_comp:psi1:a} \\
   \phi \phi_{\e}^2 \cdot \big(  | \tf{1}{z} \psio|  +   |\na \psio| \big) & \les_{\e}  |\xx|^{\al} 
   \cdot \nlinf{  \phi_{\e}^2 \phi \ww }, \label{eq:stab_comp:psi1:b} \\
     | \tf{1}{z} \psio|  +   |\na \psio|  + |\td c_l| + |\td c_{\om}| +  |\psi_z(0) | & \les_{\e} \nlinf{   \phi_{\e}^2 \phi \ww } .
     \label{eq:stab_comp:psi1:c}
\end{align}
\eseq
where  we have used $\phi \phi_{\e}^2 \gtr_{\e} |\xx|^{\al}$ 
by \eqref{def:wg_s} in the last inequality. Since $|\ww| + |\na \ww| \les \cZZ{\ww}$ by \eqref{def:cZ}, estimates \eqref{eq:psi_reg:Om} in Lemma \ref{lem:psi_reg} hold for $ \ww \in \cZ$ with upper bound 
$1+ \nnn{\ww}$ in \eqref{eq:psi_reg:Om} replaced by  $\cZZ{\ww}$. Using the above estimate on $\psio$ and following the proof of Lemma \ref{lem:psi_reg}, we obtain
\[
  \|  \psi_z \|_{C^{1,\al}(B_R)} +  
  \|  r \psi \|_{C^{2,\al}(B_R)}  \les_{\e, R}  \cZZ{\ww},
   \quad |\psi_z(0)| \les_{\e} \cZZ{\ww}.
\]

Since $\psi = \KRF \ast (\ww \tr^{\al-1})$, applying Lemma \ref{lem:hol_singu} directly, we obtain
\[
  \| \na \psi \|_{\dot C^{\al}} \les \nlinf{\ww}  \les_{\e} \nlinf{ \phi_{\e}^2 \phi \ww} .
\]

Using $r \pa_r \ws, \pa_z \ws \in C^{1, \f18 }$ from \eqref{eq:3D_solu_reg}, 
and the above estimates for $ \psi $, we obtain
\bseq\label{eq:stab_comp_pf1}
\begin{align}
\| \tcB(\ww) \|_{C^{\al}(B_R) }  + \| \td \QQ(\ww) \|_{C^{\al}(B_R) } +   \| \cU \ww \|_{C^{\al}(B_R)} & \les_{\e, R} \nlinf{  \phi_{\e}^2 \phi \ww } \, , \label{eq:stab_comp_pf1:a} \\ 
    \|  \cU \ww \|_{C^{1,1/8}(B_R)}  
  & \les_{\e, R}  \cZZ{\ww} \, .  \label{eq:stab_comp_pf1:b}  
\end{align}
\eseq

To control the $C^{1,1/8}$ norm of $\cU \ww$, we require the weighted $W^{1,\infty}$ regularity $\ww \in \cZ$.

Using the decay estimates of the profile in \eqref{eq:3D_solu_decay} and estimates \eqref{eq:stab_comp:psi1:c}, we bound 
\bseq\label{eq:stab_comp_pf2}
\begin{align}
   | \tcB(\ww)(\xx) | +  |\xx|^{-1} \cdot |  \td \QQ(\ww)  | & \les_{\e} 
    |\na \psio| +|\xx|^{-1} \cdot |\psio  | + |\pa_z \psi(0)|
     \les_{\e} \nlinf{  \phi_{\e}^2 \phi \ww }  .
\label{eq:stab_comp_pf2:cB} 
\end{align}
For $\cU(\ww)$ in \eqref{eq:stab_comp_cU}, 
using  \eqref{eq:3D_solu_decay}, \eqref{eq:stab_comp:psi1:b},\eqref{eq:stab_comp:psi1:c},  
and $\phi_{\e}, \phi,  \al -\als + 2 \e^3 < \aln$ from \eqref{def:wg_s}, we have 
\begin{align}
\phi |\cU (\ww)|  & \les_{\e} \big(\phi_{\e}^{-2}   |\xx|^{\al}  ( |\ws| + |\xx| \cdot |\na \ws| ) 
+ \phi | r \pa_r \ws| \big) 
\nlinf{ \phi_{\e}^2 \phi \ww } \notag \\
& \les_{\e}  \min( |\xx|^{\al + 1 - \bbb}, \ang \xx^{ - \als + \aln } ) \nlinf{ \phi_{\e}^2 \phi \ww } . \label{eq:stab_comp_pf2:a} 
\end{align}
Similarly, using  \eqref{eq:stab_comp:psi1:a}, \eqref{eq:stab_comp:psi1:c}, 
and   \eqref{eq:3D_solu_decay}, we obtain the weighed  $W^{1,\infty}$ bound 
\begin{align}
  \phi ( |\cU (\ww)| + |\xx| \cdot |\na \cU (\ww) |  ) 
 & \les_{\e}  \min( |\xx|^{\al + 1 - \bbb}, \ang \xx^{ - \als + \aln } ) \cZZ{\ww}.
 \label{eq:stab_comp_pf2:b} 
\end{align}
\eseq

Since $ \bbb - 1 < \al$ and $\aln < \als$ by \eqref{def:wg_s}, the above upper bound vanishes 
for $|\xx| \to 0$ or $|\xx| \to \infty$. The above estimates imply 
\beq\label{eq:stab_comp_cU_est}
    \cU(\ww)  \in \cZ, \quad  \| \cU(\ww) \|_{\cZ} \les_{\e} \| \ww \|_{\cZ}. 
\eeq

Since $\cU(\ww)$ satisfies better local regularity than $\ww \in \cZ$ by \eqref{eq:stab_comp_pf1},
and vanishing estimates near $\xx=0$ and $\xx = \infty$ by \eqref{eq:stab_comp_pf2}, 
following the proof of compactness in Section \ref{sec:3D_comp}, we prove that 
$\cU(\ww)$ is a bounded, linear, compact operator from $\cZ \to \cZ$.

Using \eqref{eq:op_lin} and \eqref{eq:op_norm}, we extend the bounded linear operators $ \cM, \cU \colon \cZ \to \cZ$ to $ \cM, \cU \colon \cZ_{\C} \to\cZ_{\C}$. Using the Bounded Perturbation Theorem in~\cite[Theorem 1.3, Chapter III]{ElNa2000}, we establish that  
$\cT = - \QQs \cdot \na + \cM$ and $\cA = \cT + \cU$ generates a strongly continuous semigroup: 
\begin{align}\label{eq:semi_full}
  e^{s \cT } : \cZc \to \cZc, \quad  \|   e^{s \cT } \|_{\cZc \to \cZc}   
  \leq e^{ \bar C_{\e} s} , \quad 
  e^{s \cA } : \cZc \to \cZc, \quad  \|   e^{s \cA } \|_{\cZc \to \cZc}   
  \leq  e^{ \bar C_{\e} s} , 
\end{align}
for some constant $\bar C_{\e}$ depending on $C_{\e}$ in \eqref{eq:tran_semi1}, 
and the operator bounds of $\cM, \cU$ in \eqref{eq:stab_comp_cB} and \eqref{eq:stab_comp_cU_est}. 
We further estimate the growth bound for $\bar C_{\e}$ in next section.

\subsubsection{Decay estimates for $e^{ \cT s }$}\label{sec:decay_cT}

In this section, we estimate the growth bound for 
\beq\label{eq:cT_est1}
  \pa_s \ww = \cT \ww  \iff \pa_s \ww + \QQs \cdot \na \ww = \cBs \cdot \ww .
\eeq
In view of \eqref{eq:op_norm} and \eqref{eq:op_lin}, it suffices to consider real value $\ww_0$. 
Below, we consider 
$\om_0 \in D(\cT) \subset \cZ$ \eqref{def:dom_cT}. 
The decay estimate for $e^{\cT s} \om_0$ with general $\om_0$ follows from a density argument.
We estimate $\ww \phi$ and $\phi |\xx| \na \ww $. 
Using equations \eqref{eq:cT_est1}, for $i \in \{ r, z \}$, we derive 
\[
\bal
  (\pa_s + \QQs \cdot \na) (\ww \phi ) 
  & =  \big( \cBs + \f{\QQs \cdot \na \phi }{ \phi } \big) \cdot (\ww \phi) \teq I_0 \cdot (\ww \phi) ,  \\
  (\pa_s + \QQs \cdot \na) ( \phi  |\xx|  \pa_i \ww )
  & = \udb{  \B(  \cBs - \pa_i \QQQs^i + \f{\QQs \cdot \na ( \phi |\xx| )  }{ \phi |\xx| } \B) }_{:=I_i} ( \phi |\xx|  \pa_i \ww )  \udb{- \pa_i \QQQs^{\neq i}  \cdot  |\xx| \phi \pa_{\neq i} \ww  
 + |\xx|   \pa_i \cBs \cdot \phi \ww }_{:= II_i},
\eal
\]
where we use $ \pa_{\neq r }, \pa_{\neq z} $ to denote $\pa_z, \pa_r$, and similar notations for 
$\QQQs^{\neq i}$. We further derive the equations of  $ \phi^2 ( |\ww|^2 + |\xx|^2 \cdot |\na \ww|^2)$
\[
\bal
 \tf{1}{2}   (\pa_s + \QQs \cdot \na)  ( \phi^2 ( | \ww |^2 +|\xx|^2 |\na \ww|^2 ))
 & = I_0 \cdot |\ww \phi|^2
 + I_r  ( \phi |\xx| \pa_r \ww )^2 
 + I_z ( \phi |\xx| \pa_z \ww )^2  \\
 & \quad + II_r \cdot ( \phi |\xx| \pa_r \ww )
 +  II_z \cdot ( \phi |\xx| \pa_z \ww ).
\eal 
\]

 Using Cauchy-Schwarz inequality, we estimate the term $II_r, II_z$ perturbatively
\[
\bal
II_r \cdot ( \phi |\xx| \pa_r \ww )
 +  II_z \cdot ( \phi |\xx| \pa_z \ww ) 
& \leq   \tf12 |\xx|( |\pa_r \cBs | +  |\pa_z \cBs | ) \cdot ( |\ww \phi|^2 +|\phi |\xx| \pa_r \ww |^2 
 + | \phi |\xx| \pa_z \ww |^2   )  \\
 & \quad + ( |\pa_r \QQQs^z| + |\pa_z \QQQs^r| ) \cdot ( |\phi |\xx| \pa_r \ww |^2 
 + | \phi |\xx| \pa_z \ww |^2 ).
 \eal
\]

Thus using Lemma \ref{lem:tran_wg} with $\eta = 0$ on the damping coefficient, we obtain
\[
 \tf{1}{2} (\pa_s + \QQs \cdot \na) ( \phi^2 | \ww |^2 +\phi^2 |\xx|^2 |\na \ww|^2 )
  \leq - \lam ( \phi^2 | \ww |^2 +\phi^2 |\xx|^2 |\na \ww|^2 ).
\]

Using the flow map $\XXb(s, \xx )$, we establish the following estimates 
\beq\label{eq:cT_est2}
   \phi^2 ( | \ww |^2 + |\xx|^2 |\na \ww|^2 ) (\xx ,s)
  \leq e^{- 2 \lam s} \phi(\yy)^2 ( |\ww_0(\yy)|^2  + |\yy|^2 |\na \ww_0(\yy)|^2 ) |_{\yy= \XXb(-s, \xx)}.
\eeq

Recall the $\cZ$-norm from \eqref{def:cZ}. Since $\ww(\xx, s) = e^{\cT s} \ww_0$, taking supremum over $\xx$, we establish: %decay estimates:
\beq\label{eq:cT_est3}
\cZZ{ e^{\cT s} \ww_0 }
=  \|  \ww(s) \|_{\cZ}  \leq e^{-\lam s  } \| \ww_0 \|_{\cZ} , \quad \forall \ww_0 \in \cZ.
\eeq

Using the decay estimates \eqref{eq:cT_est3}, linearity, and \eqref{eq:op_norm},  we moreover obtain 
\beq\label{eq:dissp_semi0}
 \|  e^{ \cT s} \|_{\cZc \to \cZc} \leq e^{-\lam s} .
\eeq

The decay estimate of the semigroup~\eqref{eq:dissp_semi} implies the spectral property of $ \cT$ in \eqref{eq:spec} (see~\cite[Theorem 1.10, Chapter II]{ElNa2000}).
Combining estimates \eqref{eq:semi_full}, \eqref{eq:spec}, \eqref{eq:dissp_semi0}, we prove Proposition \ref{prop:semi}.

\subsubsection{Decay estimates of semigroups }\label{sec:semi_grow}

We follow \cite{chen2024vorticity,chen2024Euler} and  \cite{ElNa2000} to obtain decay estimates for the semigroup $e^{\cA t}$ with $\cA = \cT + \cU$ \eqref{def:cT}. 

For a semigroup $e^{\cA t} : \MXX \to \MXX$ on a Banach space $\MXX$, we denote by ${\sigma}(\cA)$ the spectrum of $\cA$: the set $ \{ z \in {\mathbb C}: z - \cA : D(\cA) \to \MXX  \mbox{ is not bijective}\}$. We introduce the  spectral bound ${\mathsf s}(\cA)$, the growth bound $\ssf_0(\cA)$, and the essential growth bound $\ssf_{ess}(\cA)$: %, defined by
\begin{align*}
{\mathsf s}(\cA) &\teq \sup \bigl\{ \Re (z) \colon z \in \sigma(\cA) \bigr\},  
\\
{\mathsf s}_0(\cA) &\teq \inf \bigl\{ \th \in {\mathbb R} \colon \mbox{ there exists } M_{\th} \geq 1,
 \mbox{ such that } \| e^{\cA t}\|_{\MXX \to \MXX} \leq M_{\th} e^{\th t} \mbox{ for all } t \geq 0 \bigr\},  \\
  \ssf_{ess}(\cA) &  \teq \inf\nolim_{t > 0} t^{-1} \log \| e^{\cA t} \|_{ess},
\end{align*}
where the norm $\| \cdot \|_{ess}$ is defined as 
$
 \| T \|_{ess} = \inf \{ \| T - K\|_{ \MXX \to \MXX} :   K\colon \MXX \to \MXX \mathrm{ \ is \ compact} 
    \}.$ 
   With this notation, we have the following result:
\begin{prop}[Corollary 2.11, Chapter IV,~\cite{ElNa2000}]\label{prop:growth_bd}
Let $e^{\cA t}$ be a strongly continuous semigroup generated by $\cA \colon D(\cA)\subset \MXX \to \MXX$, a closed operator. Then 
\[
\ssf_0 (\cA) = \max\{ \ssf_{ess}(\cA), {\mathsf s}(\cA) \}. 
\]
Moreover, for every $\th > \ssf_{ess}(\cA)$, the set $\s_c = \s(\cA) \cap \{ z \in {\mathbb C} \colon \Re (z) \geq \th  \} $ is finite, and the corresponding spectral projection has finite rank.
\end{prop}

We apply Proposition~\ref{prop:growth_bd} to $(\cA, \MXX) = (\cA, \cZc)$ with
$\cA = \cT + \cU$ defined in \eqref{def:cT}. 
From the decay estimate of $e^{\cT s}$ \eqref{eq:dissp_semi}, we obtain $\ssf_{ess}(\cT) \leq - \lam$. 
 Since $\cU$ is a compact operator from Proposition \ref{prop:semi} and $\ssf_{ess}$ is invariant under compact perturbations (see~\cite[Proposition 2.12, Chapter IV]{ElNa2000}),  we get 
\beq\label{eq:spec_L}
\ssf_{ess}(\cA) = \ssf_{ess}( \cT) \leq -\lam.
\eeq

\subsubsection{Hyperbolic decomposition}\label{sec:hyper}

From~\eqref{eq:spec_L} and Proposition~\ref{prop:growth_bd},  we have that 
\beq\label{eq:sig_unstab}
\s_{\sfu} \teq \s( \cA) \cap  \{ z\in {\mathbb C} : \Re(z) > - \tf{7}{10} \lam \} .
\eeq
only consists of finitely many eigenvalues of $\cA$ with finite multiplicity, and that $\s_{\sfu}$ is isolated from $\s(\cA) \backslash \s_{\sfu}$. 
\footnote{
  Since $\s_{\sfu}$ consists of finitely many eigenvalues of $\cA$, there exists $c > 0$ such that 
$\min_{z \in \s_{\sfu}} \Re (z) \geq - \f{7}{10} \lam + c > - \f{7}{10} \lam $. 
}
Applying %\cite[Theorem 2.1, Chapter XV, Part IV (Page 326)]{GoGoKa2013} 
\cite[Proposition 1.16, Chapter IV, Part I (Page 245)]{ElNa2000}
for the spectral projection, we can decompose $\cZc$ as a direct sum of the stable part $\cZs$ and unstable part $\cZu$ 
\footnote{
Note that the decompositions in \eqref{eq:dec_X} and \eqref{eq:dec_Xu} are 
spectral decompositions, \emph{not orthogonal decompositions}. 
}
\beq\label{eq:dec_X}
  \cZc = \cZu \oplus \cZs,
  \;\; 
  \s( \cA|_{\cZs} ) = \s(\cA) \backslash \s_{\sfu}
\subset \{ z: \Re(z) \leq - \tf{7}{10} \lam \}, 
  \;\; 
  \s( \cA|_{\cZu} ) = \s_{\sfu} .
\eeq

Moreover, using Proposition~\ref{prop:growth_bd} we can decompose the unstable part as 
\beq\label{eq:dec_Xu}
\cZu = \tts{\bigoplus}_{ z \in \s_{\sfu}} \ker( (z- \cA)^{\mu_z}), \quad \mu_z < \infty, \quad |\s_{\sfu}| < +\infty,
\quad  \dim( \cZu ) < \infty.
\eeq

Let $\Pi_u$ denote the \emph{spectral projection} from $\cZc$ onto $\cZu$ associated with
$\sigma_u$ in \eqref{eq:sig_unstab}, and set $\Pi_s = \Id - \Pi_u$. 
\footnote{
Note that $ \Pi_u$ is a \emph{spectral} projection rather than an orthogonal projection.
The operator norms as maps $\Pi_s : \cZc \to \cZs$, $\Pi_u : \cZc \to \cZu$ are finite but need not be less than $1$.
}
We have $\Pi_s + \Pi_u = \Id$ and $\Pi_u^2 = \Pi_u$. Since $\cZs$ is $\cA$ invariant, the restriction of the semigroup $e^{\cA s}$ on $\cZs$ is generated by $\cA|_{\cZs}$. 
\footnote{
This follows from the subspace semigroup theorem
\cite[Chapter II, Section 2.3]{ElNa2000}
and the properties of the part operator
\cite[Chapter IV, Part I, Lemma 1.15, p.~245]{ElNa2000}.
}
Thus, using \eqref{eq:spec_L} and \eqref{eq:sig_unstab}, we get
\[
\ssf_{ess}(\cA|_{\cZs} )  \leq \ssf_{ess}(\cA ) \leq -\lam, \quad  {\mathsf s}( \cA|_{\cZs} ) \leq - \tf{7}{10} \lam.
\]
From Proposition~\ref{prop:growth_bd} and the definition of $\ssf_0(\cA|_{ \cZs })$, we obtain
\beq\label{eq:decay_stab}
 \| e^{\cA s } \Pi_s f \|_{\cZc} \leq C_{\e} e^{ -  \f35 \lam s } \| \Pi_s f \|_{ \cZc }
 \leq C_{\e} e^{ -  \f35 \lam s } \|  f \|_{ \cZc } , \quad \forall \, f \in \cZc, \ s \geq 0 .
\eeq
We also note that since $\cZu $ has finite dimension, $\s(\cA|_{\cZu } ) = \s_{\sfu}$ \eqref{eq:sig_unstab}, the operator
$\cA|_{ \cZu }$ is bounded and may be represented as a matrix with eigenvalues $> - \f{7}{10} \lam$; taking $\f{7}{10}< \f{3}{4}$,  we obtain
\footnote{
Since $\cZu$ is finite-dimensional, polynomial factors from Jordan blocks 
of $\cA |_{\cZu}$ are absorbed into \(C_\varepsilon e^{\frac34\lambda s}\).
}
\beq\label{eq:decay_unstab}
 \| e^{-\cA s} \Pi_u f \|_{\cZc} \leq  
 C_{\e} e^{ \f{3}{4} \lam s} \| \Pi_u f \|_{\cZc}
 \leq  C_{\e} e^{ \f{3}{4} \lam s} \|  f \|_{\cZc} ,  \quad \forall \ f \in \cZc , \quad s \geq 0 .
\eeq

\subsection{Decomposition and regularization}
Below, we perform a few regularizations and decompose the equation \eqref{eq:dyn_lin} by 
following the general framework developed in \cite{chen2024Euler,chen2024vorticity,ChenHou2023a}.

\subsubsection{Regularization of the compact operator}\label{sec:dyn_regu}
Recall the cutoff function $\chi_1$ from \eqref{def:chi1}. 
We introduce a radially smooth function $\Xi(\yy)$ and $\Xi_c(\yy)$ 
in the whole space $\R^3$ (not  $\R_+ \times \R$) with 
\beq\label{def:Xi}
\Xi \in C_c^{\infty}(\R_3), 
\quad \supp(\Xi) \subset B(0, 1), \quad \Xi \geq 0, \quad  \int_{\R^3}  \Xi(\yy) d \yy = 1,
\quad \Xi_c(\yy) = c^{-3} \Xi( \yy / c ), 
\eeq
for any $c > 0$.  We consider the following 
regularization of the operator $\cU$ 
\beq\label{def:cU_h}
  (\cU_h \ww)(\xx) \teq \Xi_{h_2} \ast ( \chi_{2, h_1 }(\xx) \cdot (\cU \ww)(\xx) ),
  \quad \chi_{2, h_1}(\xx) \teq \chi_1( 4 h_1 \xx ) - \chi_1( \f{\xx}{2 h_1} ) , 
  \quad h_2 \leq h_1 \leq \tf14.
\eeq
We treat $ \chi_{2, h_1} \cdot \cU \ww$ as an axisymmetric function in $\R^3$ and apply the convolution in $\R^3$.  Since $\Xi(\yy)$ is radially symmetric, the function $\cU_h \ww$ is axisymmetric, and we evaluate the function at $\xx = (r, z)$. Since $\Xi_{h_2}$ has a support size at most $h_2$ and $h_2 \leq h_1$, by definition of $\chi_1$ \eqref{def:chi1}, we have 
\beq\label{eq:supp_cU_h}
  \supp \big(  \chi_{2, h_1}  \big) \subset \{ \xx : |\xx| \in [2 h_1, \tf{1}{2 h_1} ] \} , \quad \supp( \cU_h \ww ) \subset \{ \xx : |\xx| \in [ h_1, \tf{1}{h_1} ] \} ,
\eeq
and $\cU_h \ww(\xx) \in C_c^{\infty}$ both as a function in $\R^3$ 
and as a function on $(r, z) \in \R_+ \times \R$. Since 
\[
  |\na^k (\Xi_{h_2} \ast f)(\xx) | \les h_2^{-k} \sup\nolim_{ |\yy- \xx| \leq h_2 } |f(\yy)|, 
\]
using the definition of $\phi$ in \eqref{def:wg_s}, $\phi^{-1} \phi_{\e}^{-2} \les 1$, the \eqref{eq:supp_cU_h}, and \eqref{eq:stab_comp_pf1}, we obtain
\beq\label{eq:Uh_bound}
\bal
\cZZ{\cU_h \ww} 
+ \| \cU_h \ww \|_{C^{k, \al}}
& \les_{h_1}  \sup\nolim_{ |\yy| \in [h_1, 1 / h_1]  } |\na^{\leq k+1} \cU_h \ww(\yy)| \\
& \les_{h_1, h_2, k}  \sup\nolim_{ |\yy| \in [h_1, 1 / h_1]  }  |\cU \ww(\yy)| 
\les_{h_1, h_2 ,k, \e}   \nlinf{ \phi_{\e}^2 \phi \ww }.
\eal 
\eeq

Next, we estimate $ \nlinf{ \phi (\cU_h \ww - \cU \ww ) }$. Using definition  \eqref{def:cU_h}, 
 we decompose 
\[
   \phi |(\cU_h \ww - \cU \ww ) |
  \leq \phi | \Xi_{h_2} \ast  (\chi_{2, h_1} \cU \ww)  -  (\chi_{2, h_1} \cU \ww) |
  + \phi | (\chi_{2, h_1} \cU \ww) - \cU \ww | \teq  I + II.
\]

Since $\supp(I) \subset \{ |\xx| \in [h_1, h_1^{-1}] \}$, using estimate \eqref{eq:stab_comp_pf1:a} for $I$, we bound 
\footnote{
\label{foot:axi_reg}
Note that a $C^{\g}$ axisymmetric function $f $ with $\g \in (0,1)$ 
can be naturally extended to a $C^{\g}$ function in $\R^3$. See discussion 
at the beginning of Section  \ref{sec:qual_reg_stream}.
}
\[
\bal
  |I|  & \les_{h_1}  | \Xi_{h_2} \ast  (\chi_{2, h_1} \cU \ww)  -  (\chi_{2, h_1} \cU \ww) |  
  \les_{h_1} h_2^{\al / 2} \|   \chi_{2, h_1} \cU \ww \|_{C^{\al/2}} \\
 &  \les_{h_1} h_2^{\al/2} \| \cU \ww \|_{C^{\al/2}( B_{ 4 / h_1  } )}
\les_{\e, h_1} h_2^{\al/2} \nlinf{  \phi_{\e}^2 \phi \ww } .
\eal
\]

Since $\supp( 1- \chi_{2, h_1} ) \subset \{ |\xx|  \leq 4 h_1 \mw{ \ or \ } |\xx| \geq 1/ (4 h_1) \}$  by \eqref{eq:supp_cU_h}, using  estimate \eqref{eq:stab_comp_pf2:a}, we bound 
\[
  |II| \leq (1 - \chi_{2, h_1}) \cdot \phi |\cU \ww| 
   \les_{\e}  h_1^{ \min( \al + 1 - \bbb, \als - \aln ) } \nlinf{ \phi_{\e}^2 \phi  \ww }.
\]

Since $\bbb \in [1.01, 1.3)$ and $\al + 1 - \bbb > \als - \aln$, combining the above estimate, we prove 
\[
 \nlinf{ \phi (\cU_h \ww - \cU \ww ) }
  \leq (  C(\e)  h_1^{\aln- \als} +    C(\e , h_1) h_2^{ \al/2  } ) \nlinf{ \phi_{\e}^2 \phi \ww }.
 \]
for some constant $C(\e , h_1)$ depends on $\e, h_1$. 
Since $\als - \aln > 0$, we can first choose $h_1 = h_1(\e)$ small enough, and then  
$h_2 = h_2(h_1, \e)$ small enough to obtain
\beq\label{eq:cU_dif}
 \nlinf{ \phi (\cU_h \ww - \cU \ww ) }
  \leq 10^{-4} \lam \nlinf{ \phi_{\e}^2 \phi \ww } .
 \eeq

Since we fix $h_1, h_2$, in the following estimates, we treat $h_1, h_2$ as constant depending on $\e$.

\subsubsection{Decomposition}

We follow the splitting method developed in \cite{ChenHou2023a} to decompose  \eqref{eq:dyn_lin}.
Recall the equation for the perturbation from \eqref{eq:dyn_lin}
\beq\label{eq:dyn_recall}
    \pa_s  \ww + (\QQs + \td \QQ(\ww)) \cdot \na \ww =  ( \cBs + \tcB( \ww ) ) \cdot \ww +  \tcB(\ww) \cdot \ws - \td \QQ(\ww) \cdot \na \ws.
\eeq

Recall the operator $\cA, \cU$ from \eqref{def:cT}. We decompose the solution to \eqref{eq:dyn_lin} as
\[
  \ww = \ww_1 + \ww_2 ,
\]
and the equation \eqref{eq:dyn_lin} as 
\bseq\label{eq:split}
\begin{align}
  \pa_s \ww_1 + (\QQs + \td \QQ(\ww)) \cdot \na \ww_1 & =  \cBs \cdot \ww_1 
  + \cU(\ww_1) - \cU_h(\ww_1) + \tcB( \ww ) \cdot \ww,  \label{eq:split:a} \\
  \pa_s \ww_2 + \QQs \cdot \na \ww_2 & =  \cBs \cdot \ww_2 
  + \cU(\ww_2) + \cU_h(\ww_1)  \label{eq:split:b}.
\end{align}

It is not difficult to see that $\ww_1 + \ww_2$ solves \eqref{eq:dyn_lin}.  Using the operator $\cA$, we rewrite
\beq\label{eq:split:c}
  \pa_s \ww_2 = \cA \ww_2 + \cU_h(\ww_1).
\eeq
\eseq

Recall the decomposition from \eqref{eq:dec_X} and spectral 
projections $\Pi_s, \Pi_u$ of $\cZc$ onto $\cZs, \cZu$ introduced in \eqref{eq:dec_Xu}-\eqref{eq:decay_unstab}
To handle the unstable part, we construct $\ww_2$ as follows 
\footnote{
The idea of using the backward-in-time semigroup to construct the unstable manifold is classical in ODE theory, see, e.g. \cite[Section 9.2]{teschl2012ordinary}. 
}
\bseq\label{eq:w2_form}
\begin{align}
 \ww_{2}(s) & \teq \Re \int_{0}^s e^{\cA(s -\tau)} \Pi_s \cU_h \ww_1 (\tau) d \tau 
 - \Re \int_s^{\infty} e^{\cA(s - \tau)} \Pi_u \cU_h \ww_1(\tau) d \tau  \label{eq:w2_formb}  .
\end{align}
\eseq
We encode the above linear map from $\ww_1$ to $\ww_2$ as 
\beq\label{def:cH2}
  \cH_2 \ww_1 \teq \mbox{Right Side of }  \eqref{eq:w2_form}.
\eeq

Recall  $\lam$ from Lemma \ref{lem:tran_wg}. For later nonlinear estimates, we define the space $Y$ and ball
\beq\label{def:YY} 
    \| \ww \|_{Y} \teq \sup\nolim_{s \geq 0} e^{ \f45 \lam s } \nlinf{ \phi \phi_{\e}^2 \ww } ,
    \quad \bar B( R, Y ) \teq \{ \ww: \| \ww\|_Y \leq R \}.
\eeq

For operator $\cH_2$, we have the following estimates.

\begin{lem}\label{lem:semi_cH2}
Let $\lam$ be the parameter from Lemma \ref{lem:tran_wg}. For any $ f \in Y$ and $s \geq 0$, we have  
\[
\| \cH_2 f(s) \|_{\cZ} \les_{\e} e^{- \f35 \lam s} \| f \|_Y .
\]
\end{lem}

\begin{proof}
Since $f \in \cZ$ is real,  applying  \eqref{eq:decay_stab} and \eqref{eq:decay_unstab} to the formula \eqref{eq:w2_form}, we obtain
\[
\bal
   \| \cH_2 f(s) \|_{\cZ} 
   & = \B\|  \Re \int_{0}^s e^{\cA(s -\tau)} \Pi_s \cU_h f (\tau) d \tau -  \Re \int_s^{\infty} e^{\cA(s - \tau)} \Pi_u \cU_h f (\tau) d \tau  \B\|_{\cZ}  \\
  &  \les_{\e} \int_0^s e^{-\f35 \lam (s -\tau)} \| \cU_h f (\tau) \|_{\cZ} d \tau
   + \int_s^{\infty} e^{ - \f{3}{4} \lam (s -\tau)} \| \cU_h f (\tau) \|_{\cZ} d \tau .
\eal
\]

Applying the bound \eqref{eq:Uh_bound} on $\cU_h$ and the $Y$-norm defined in \eqref{def:cY}, we obtain
\[
  \| \cU_h f(\tau) \|_{\cZ} 
\les_{ \e} \nlinf{ \phi \phi_{\e}^2 f (\tau) } \les_{ \e} e^{-\f45 \lam \tau } \| f \|_Y . 
\]
Note that $h_1, h_2$ appearing in \eqref{eq:Uh_bound} have been chosen as constants depending on $\e$ in \eqref{eq:cU_dif}. Combining the above two estimates and using $\lam >0$, we prove 
\[
     \| \cH_2 f(s) \|_{\cZ}  
    \les_{h_1, h_2, \e} \B( \int_0^s e^{-\f35 \lam (s -\tau)} e^{-\f45 \lam \tau}   d \tau
   + \int_s^{\infty} e^{ -\f{3}{4} \lam (s -\tau)} e^{-\f45 \lam \tau}  d \tau \B) \| f \|_Y 
  \les_{h_1, h_2, \e} e^{- \f35 \lam s}\| f \|_Y .
\]
\end{proof}

\paragraph{\bf Regularization of unstable initial data}
Note that the initial data $\ww_{2}(0) = - \ww_{2,u}(0)$ depends on the ``future" of the solution, which we cannot prescribe. To construct initial data with sufficient regularity and compact support, we perform a regularization. Recall the mollifier $\Xi_{h_4}$ and cutoff function $\chi_{h_3}$ from \eqref{def:cU_h}
with parameters $h_1, h_2$ replaced by $h_3, h_4$. %To construct compact support initial data, 
Given any function $\ww_{1, \init} \in \cZ \cap C^{1,\al}$, we impose initial condition to \eqref{eq:split:a}
\bseq\label{eq:w1_init}
\beq
  \ww_1(0) = \ww_{1, \init} + \cH_{\init} \ww_2(0), 
  \quad \cH_{\init} f \teq  \Xi_{h_4} \ast (  \chi_{2, h_3}  f)   - f, 
  \quad h_4 \leq h_3, 
\eeq
for some $h_4 \leq h_3$ to be determined later. 
The above choice leads to initial data of \eqref{eq:dyn_recall}
\beq
  \ww(0) = \ww_1(0) + \ww_2(0) = \ww_{1,\init} +  \Xi_{h_4} \ast (  \chi_{2, h_1}  \ww_2(0) )  \,.
\eeq
\eseq

Consider $f \in \cZ$.  Since $\phi_{\e}$ \eqref{def:wg_s} is weaker than $1$ for $\xx$ near $0$ and $\infty$, $ \supp( \Xi_{h_4} \ast (  \chi_{2, h_3}  \ww_2(0) ) ) \in \{ \xx, |\xx| \in [h_3, h_3^{-1}] \} $ by \eqref{eq:supp_cU_h}, and since $\chi_{2, h_3} f$ can be extended to a $C^{1/2}$ function in $\R^3$ (see Footnote \ref{foot:axi_reg}), following the proof of \eqref{eq:cU_dif},  we have  
\begin{align}\label{eq:cH_init}
    \nlinf{ \phi \phi_{\e} \cH_{\init} f } 
& \leq \nlinf{ \phi \phi_{\e} \big( \Xi_{h_4} \ast (\chi_{2,h_3} f ) - (\chi_{2,h_3} f ) \big) 
+ \phi \phi_{\e} (\chi_{2,h_3} f - f )} \\
& \leq C(\e, h_3) h_4^{1/2} \| \chi_{2, h_3} f \|_{C^{\f12} ( B_{ 4/h_3  } ) }
+  C h_3^{\e^3} \nlinf{\phi f}
\leq ( C(\e, h_3) h_4^{1/2} + C h_3^{\e^3} ) \cZZ{f}. \notag
\end{align}

By definition of $\Xi_{h_4}$ \eqref{def:Xi} and $\chi_{2, h_3}$ and \eqref{def:cU_h}, 
and using estimates similar to \eqref{eq:Uh_bound}, we obtain
\beq\label{eq:w2_regu_init}
\bal
    \Xi_{h_4} \ast (  \chi_{2, h_3}  \ww_2(0) ) \in C_c^{\infty}(\R^3),
    \quad 
    \supp (     \Xi_{h_4} \ast (  \chi_{2, h_3}  \ww_2(0) ) ) & \in \{ \xx: |\xx| \in [h_3, h_3^{-1}]  , \\
     \|  \Xi_{h_4} \ast (  \chi_{2, h_3}  \ww_2(0) )  \|_{C^{k, \al}(\R^3) } 
   &  \les_{h_3, h_4, k ,\e}  \cZZ{ \ww_2(0) }.
   \eal
\eeq

\subsection{Nonlinear stability and Proof of Theorem \ref{thm:main_blowup}}
\label{sec:main_blowup_pf}

In this section, we reformulate the construction of a global solution to \eqref{eq:dyn_lin} 
as a fixed point problem and present the nonlinear stability results in Theorem \ref{thm:non_stab}. 
Then we use it to prove the main blowup result Theorem \ref{thm:main_blowup}. %We follow the general framework developed in \cite{chen2024Euler,chen2024vorticity,ChenHou2023a}. 

Our main nonlinear stability result is the following.

\begin{thm}\label{thm:non_stab}

Let $\beps_{10}, \beps_{12}, \lam$ be chosen in Theorem \ref{thm:3D_solu} and Lemma 
\ref{lem:tran_wg}. Fix any $ \e = \f13 - \al \in (0, \beps_{12}]$. There exists constants $h_4\leq h_3,
h_i = h_i(\e)$ and $\d_0(\e)$ small enough such that 
for any $\d \leq \d_0(\e) $, the following statement holds. 
For any $\ww_{1,\init} \in C^{1,\al}$ odd in $z$, with 
\beq\label{eq:ass_dyn_init}
\pa_z \ww_{1, \init}(0) = 0, \quad 
  \nlinf{ \phi \ww_{1, \init}} < \tf{1}{8} \d, 
  \quad  \ww_{1,\init} + \ws \in C_c^{1,\al},
\eeq
there exists an initial data $\ww_2(0) \in \cZ$, a global solution $\ww_2$ of \eqref{eq:split:b} given by~\eqref{eq:w2_form}, and a global solution $\ww_1$ to \eqref{eq:split:a} with initial data
$\ww_1(0) = \ww_{1,\init} + \cH_{\init} \ww_2(0)$, satisfying the bounds
\begin{align}\label{eq:solution:small}
\nlinf{ \phi_{\e} \phi \ww_1  }  & \leq  \d  e^{- \f{9}{10} \lam s }, \quad 
 \| \ww_2  \|_{ \cZ }  \les_{\e} \d   e^{- \f35 \lam s },
\end{align}
for all $s\geq 0$. We emphasize that the initial data $\ww_2(0)$ is constructed via \eqref{eq:w2_form} (simultaneously with $\ww_1$) to lie in a finite-dimensional subspace of $\cZ$. 

\end{thm}

\begin{remark}[\bf Set of initial data]\label{rem:init_data2}

From \eqref{eq:w2_form}, $ \ww_2(0)$ lies in the finite-dimensional subspace $ \Re( \cZu)$  of $\cZ$ 
with $\cZu$ defined in \eqref{eq:dec_Xu}. %Note that 
We can choose any $\ww_{1,\init}$ 
in an open set in some weighted $ C^{1, \al}$ space $X_1$ defined by \eqref{eq:ass_dyn_init} 
with $\ww_{1,\init} + \ws \in C_c^{1,\al}$. 
\footnote{
For $ \ww_{1, \init} \in C^{1,\al}$ odd in $z$, condition 
$ \nlinf{ \phi \ww_{1,\init}} < \infty$ implies $ \pa_z \ww_{1, \init} (0) = 0$. 
}
Moreover, $\Xi_{h_4} \ast (  \chi_{2, h_3}  \ww_2(0) ) , \ww_2(0) \in \Re( \cZu)$ forms a finite-dimensional subspace $X_2 \subset C_c^{\infty}(\R^3)$ by \eqref{eq:w2_regu_init}. From \eqref{eq:w1_init}, \eqref{eq:dyn_decomp}, the initial data for \eqref{eq:dyn_full} $\ws + \ww_{1,\init}  + \ww_{2}(0)$ forms a finite codimension set $X_1 \oplus X_2$ in $C_c^{1,\al}$.

\end{remark}

We defer the proof of Theorem \ref{thm:non_stab} to 
Sections \ref{sec:dyn_fix_pt}, \ref{sec:dyn_EE}, and \ref{sec:dyn_comp}

Based on 
Theorem \ref{thm:non_stab}, we are in a position to prove Theorem \ref{thm:main_blowup} 
on the blowup result.

\begin{proof}[Proof of Theorem \ref{thm:main_blowup}]

Due to the inclusion $C^{\al} \subset C^{\g}$ for any $\al \leq \g \in (0, 1)$,
we only need to prove Theorem \ref{thm:main_blowup} for $\e = \f13 - \al$ and $\delin$ sufficiently small.  Let $\beps_{12}$ be as in Lemma \ref{lem:tran_wg}. Taking $ \beps_{12}^{\pr} \in (0, \beps_{12}] $ small enough, for any $\e \in (0, \beps_{12}^{\pr} ]$, 
we estimate \eqref{def:TTa} as 
\beq\label{eq:ran_TT}
  \TTa^{-1}  = 8 + O( \e^{\mhk}) \in [7,10] .
\eeq

We choose $\al_0 =\f13 -  \beps_{12}^{\pr} $ in Theorem \ref{thm:main_blowup}. Let $\d_0(\e)$ be the constant determined in Theorem \ref{thm:3D_solu}.
We choose the parameter $\g$ in Theorem \ref{thm:main_blowup} as $\gamma = \al$, and require
$ \delin = \d \leq \d_0(\e)$. 

Below, since \(\e\) is fixed, we keep its dependence \emph{implicit} and the \(\d\)-dependence \emph{explicit}.

\vspace{0.05in}
\paragraph{\bf Initial data}
We choose initial perturbation $\ww_{1, \init}$ small enough with properties \eqref{eq:ass_dyn_init}. 
Recall the cutoff function $\chi_1$ from \eqref{def:chi1} and 
the mollifier $\Xi$ from \eqref{def:Xi}. We choose
\[
 \ww_{1,\init}(\xx) \teq 
 C_{\ell, L_1} \Xi_{ \ell } \ast (  \chi_1(\xx / L_1 )  \ws ) - \ws,
 \quad  C_{\ell, L_1} \teq  \f{ \pa_z \ws(0) }{ \pa_z  \big( \Xi_{ \ell } \ast (  \chi_1(\xx / L_1 )  \ws ) \big) (0)}, 
\]
with $L_1 > 10$ large enough and $\ell < 1$ small enough to be chosen. 
Since $\pa_z \ws(0) \neq 0$, by definition, $ \ww_{1,\init}(\xx)  + \ws $ 
is $C_c^{\infty}$ both as a function in $\R_+ \times \R$ and as an axisymmetric function in $\R^3$:
\bseq\label{eq:ww_init_reg}
\beq
 \ww_{1,\init}(\xx)  + \ws \in C_c^{\infty}(\R^3), 
 \quad \pa_z \ww_{1,\init}(0 ) = 0,
 \quad \pa_r \ww_{1,\init}(0 ) = 0,
 \quad \ww_{1,\init} \mbox{ \ odd  in \  }  z.
\eeq
Using the formula \eqref{eq:w1_init} for the initial data,
the $C_c^{\infty}$ regularity  \eqref{eq:ww_init_reg}, \eqref{eq:w2_regu_init}, we obtain
\beq
 \Om(0) = \ww(0) + \ws   = \ww_{1, \init} + \ws +  \Xi_{h_4} \ast (  \chi_{2, h_3}  \ww_2(0) ) 
 \in C_c^{\infty}(\R^3) .
\eeq
\eseq

Next, we show that $\ww_{1,\init}$ satisfies the smallness condition in \eqref{eq:ass_dyn_init}. We decompose 
\[
\Xi_{ \ell } \ast (  \chi_1(\xx / L_1 )  \ws ) - \ws 
 = \big( \Xi_{ \ell } \ast (  \chi_1(\xx / L_1 )  \ws ) - \chi_1(\xx/ L_1) \ws \big)
 +  ( \chi_1(\xx / L_1 ) - 1)  \ws \teq I(\xx) + II(\xx).
\]

Since $\ell < 1 < 10 < L_1$ and $\supp(I) \subset \bar B_{2 L_1}$, using \eqref{eq:3D_solu_reg}, 
for any $\g < \al$, we obtain 
\bseq\label{eq:est_win}
\begin{align}
 \| \pa_z I(\xx) \|_{C^{\g}( \bar B_1 )} & \les  \ell^{\al -\g } \| \pa_z ( \chi_1(\xx / L_1 )  \ws )  \|_{C^{\al}( \bar B_2 ) } \les_{\e ,\g} \ell^{\al -\g}, \label{eq:est_win:a} \\
   \| I(\xx) \|_{L^{\infty}( \R^3 )} & \les_{L_1} \ell \| \chi_1(\xx / L_1 )  \ws  \|_{
   C^{0, 1} ( \bar B_{2R}) } \les C(\e,  L_1) \ell \,  . \label{eq:est_win:b} 
\end{align}
\eseq
Note that $\pa_z ( \chi_1(\xx / L_1 )  \ws \in C^{\al}(\R^3), \chi_1(\xx / L_1 )  \ws  \in C^{0,1}(\R^3)$ as functions in $\R^3$. See Footnote \ref{foot:axi_reg}. Since $\pa_z \ws(0) \neq 0$ and $II(\xx) = 0$ for $|\xx| \leq L_1$, using the first estimate with $\g= 0$ and by requiring $ \ell \ll \e$,  we obtain 
\[
  |\pa_z (I + II)(0)| \les_{\e} \ell^{\al} ,
  \quad \Rightarrow  |C_{\ell, L_1} - 1| \leq C_{\e} \ell^{\al} \ll \tf12.
\]

We decompose $\ww_{1,\init}$ as 
\[
\ww_{1,\init} =  C_{\ell, L_1} ( I + II+ \ws) - \ws
=C_{\ell, L_1}  (I + II) +   ( C_{\ell, L_1} -1) \ws  
\]

For $|\xx| \leq 1$, since $II(\xx)=0$, using the above estimates for $I$, we bound 
\beq\label{eq:est_win2}
 \| \pa_z \ww_{1, \init} \|_{C^{\g}(\bar B_1)} \les \| \pa_z I \|_{C^{\g}(\bar B_1)}
 + |( C_{\ell, L_1} -1) |  \| \pa_z \ws \|_{C^{\g}(\bar B_1)}
 \les_{\e, \g} \ell^{\al-\g} + \ell^{\al} \les_{\e, \g} \ell^{\al -\g}.
\eeq

For general $\xx$ with $II(\xx)\neq 0$, using the decay estimate for $\ws$ in \eqref{eq:3D_solu_decay:a}, \eqref{eq:est_win:b}, and $\supp(\chi(x/ L_1)) \subset \bar B_{2 L_1}$, 
$\chi(x/L_1) - 1=0 $ for $|x| \leq L_1$, we bound 
\beq\label{eq:est_win3}
 |I| \les_{\e, L_1} \one_{|x| < 2L_1} \ell , 
 \quad |II| \les_{\e} \one_{|x|> L_1} \ang \xx^{-\als},
 \quad  | (C_{\ell, L_1} - 1) \ws | \les_{\e} \ell^{\al} \ang \xx^{-\als}.
\eeq

Since $\na \ww_{1,\init}(0) = 0$ by \eqref{eq:ww_init_reg} 
and $\ww(1,\init)(r, 0) = 0$, using \eqref{eq:est_win2}, \eqref{eq:est_win3}, we bound 
\[
 |\ww_{1,\init}(\xx)| \les_{\e, \g} \one_{|\xx| \leq 1} |z| \cdot |\xx|^{ \g}  \ell^{\al -\g}
 + \one_{|\xx| > 1} ( \one_{|x| < 2L_1} \ell  + \one_{|x|> L_1} \ang \xx^{-\als} 
 + \ell^{\al} \ang \xx^{-\als} ).
\]

Recall the weight $\phi , \phi_{\e}$ form \eqref{def:wg_s}. Since $\phi_{\e}\leq 1$, 
$\phi \les |\xx|^{-\bbb} + |\xx|^{\aln}$, $\aln < \als$ and $\bbb < 1 + \al$, choosing 
$\g = \bbb-1 \in (0, \al)$ and using the above estimates, we obtain
\[
 |\ww_{1,\init} \phi| \les_{\e} \ell^{\al-\g} + 
 \one_{|\xx|>1}
  ( \one_{|x| < 2L_1} \ell |\xx|^{\aln}
   + \one_{|x|> L_1} \ang \xx^{\aln-\als} 
 + \ell^{\al} \ang \xx^{\aln-\als} )
 \les_{\e} \ell^{\al+ 1 -\bbb} + \ell L_1^{\aln} + \ell^{\al} + L_1^{\aln-\als}. 
\]

By first choosing $L_1$ large enough and then $\ell$ small enough, we obtain $ \nlinf{\phi \ww_{1,\init} } < \tf18 \d$ 
with $\d$ chosen in Theorem \ref{thm:non_stab}. Thus, assumptions \eqref{eq:ass_dyn_init} are satisfied. Using Theorem \ref{thm:non_stab}, we construct a global solution $\ww_1, \ww_2$ to \eqref{eq:split}, and a global solution $\ww = \ww_1 + \ww_2$ to \eqref{eq:dyn_recall}. 
Since equation \eqref{eq:dyn_recall} is the linearized equation of \eqref{eq:dyn_full} around $\ws$, 
we obtain a global solution $\Om = \ww + \ws$ to \eqref{eq:dyn_full}.  
Using estimate  \eqref{eq:solution:small} in Theorem \ref{thm:non_stab}, and $\nlinf{ f \phi_{\e} \phi} \leq \| f \|_{\cZ}$ by \eqref{def:cZ}, we obtain
\bseq\label{eq:thm_blow_pf1}
\beq
\bga
  \nlinf{ \phi_{\e} \phi (\Om(s) - \ws) }  \les_{\e} \d e^{-3/5 \lam s}.
\ega
\eeq

Using estimate \eqref{eq:QB_asymp_1:scale} with $\eta = \ww$ and the above estimates, we obtain
\beq
 |\td c_l(\ww) |+ |\td c_{\om}(\ww)| \les \e^{-1} \nlinf{ \phi_{\e} \phi \ww} \les_{\e} 
\d 
 e^{-3/5 \lam s}.
\eeq
\eseq

\vs{0.05in}
\paragraph{\bf Modulation functions and blowup time}
Using  \eqref{eq:thm_blow_pf1}, we obtain the limiting values :
\bseq\label{eq:thm_blow_T}
\beq
 \mfL \teq \tts{\int}_0^{\infty} \tcl(\tau) d \tau  , 
 \quad 
 \mfM \teq \tts{\int}_0^{\infty}  (\tcw(\tau) + \al \tcl) d \tau ,
 \quad |\mfL|  + |\mfM| \les_{\e} \d \, .
\eeq

Using \eqref{eq:dyn_rela}, \eqref{eq:thm_SS_pf2}, and \eqref{eq:thm_SS_pf1}, we obtain
\beq
\bal
\comth &= c_l + \al \com = \cls + \al \cws +  \tcl + \al \tcw,   
\quad \bcomth  \teq \cls + \al \cws = - \TTa^{-1} \asymp - 1,  \\
|\bcomth|^{-1}  & = \TTa =  \f18 + O( \e^{\mhk}), % = - \f18 +  O( \e^{\mhk} ),
   \quad    \cxs = -\f{\cls}{ \bcomth}  \asymp \e^{-1}. 
\eal 
\eeq
\eseq

Using the decomposition \ref{eq:dyn_decomp} and estimate \eqref{eq:thm_blow_pf1}, \eqref{eq:thm_blow_T}, we obtain
\[
    | \int_0^{s} c_l(\tau) - \mfL  - s \cls | 
+    | \int_0^{s} \comth(\tau) - \mfM  - s \bcomth | 
\les \int_s^{\infty} |\tcl| + |\tcw| 
  \les_{\e} \d e^{-3/5 \lam s} .
\]

Recall $\CCw(0) = \TTa $ from \eqref{eq:dyn_full:b}. We estimate $\CCl, \CCw$ in \eqref{eq:rescal} as 
\[
\bal
& \CCl(s) = \exp( - \int_0^s  c_l(\tau) d \tau )
=e^{ - \mfL - s \cls  + O(  C_{\e} \d e^{-3/5 \lam s}) } ,
\quad \CCw(s) = \TTa e^{ - \mfM - s \bcomth  + O(  C_{\e} \d e^{-3/5 \lam s}) } .
\eal
\]

Since $ -\bcomth \asymp 1$, $|\bcomth^{-1} |= \TTa$, and $\d<1$, using \eqref{eq:rescal:b}
and the above estimates,  we estimate the blowup time $T = t(\infty)    =  \int_0^{\infty} \CCw^{-1}(\tau) d \tau$:
\begin{align}\label{eq:thm_blow_T2}
    |T -  1  |
  & \leq \B| \int_0^{\infty}(\CCw^{-1}(\tau) - \TTa^{-1} e^{s \bcomth}  ) d s\B|
  \les_{\e} \int_0^{\infty} (\mfM + \d e^{-3/5 \lam s} ) e^{s \bcomth} d s 
  \les_{\e}  \d, \notag \\
  |T - t(s)| &  = \int_s^{\infty} \CCw^{-1}(\tau) 
 = \TTa^{-1} \int_s^{\infty} e^{  \mfM + \tau \bcomth  + O( C_{\e} \d e^{-3/5 \lam \tau}) }
 = \TTa^{-1} \int_s^{\infty}
 e^{  \mfM + \tau \bcomth }  d \tau e^{  O( C_{\e} \d e^{-3/5 \lam s} )}  \notag \\ 
 & = \TTa^{-1} |\bcomth|^{-1}  e^{  \mfM + s \bcomth } e^{  O( C_{\e} \d e^{-3/5 \lam s} )} 
 =   e^{  \mfM + s \bcomth } e^{  O( C_{\e} \d e^{-3/5 \lam s} )} .
\end{align}

Combining the above scaling relation and using $\cxs = -\f{\cls}{ \bcomth}  \asymp \e^{-1}$ \eqref{eq:thm_blow_T}, we obtain
\beq\label{eq:scal_pf1}
\bal
\cE_l(s) & \teq \CCl(s) (T - t(s))^{ -\cxs} = e^{ -\mfL - \cxs \mfM } e^{ ( O ( C_{\e} \d e^{-3/5 \lam s} )}, \\
  \CCw(s) & = \TTa \cdot |T- t(s)|^{-1}  e^{  (O( C_{\e} \d e^{-3/5 \lam s} )}.
 \eal
\eeq

\paragraph{\bf Asymptotically self-similar blowup}
Note that the exponential error term converges to $1$ as $s \to \infty$. 
Recall the rescaling relation \eqref{eq:rescal},  \eqref{eq:thm_blow_pf1}
and the profile relations \eqref{eq:thm_SS_pf2}
\beq\label{eq:thm_blow_pf2}
  \om^{\th}( \CCl(s) \xx,  t(s) ) =  \CCw  \cpsia \Om^{\th}(\xx, s) ,
  \quad  \omthss = \TTa \cpsia \omsth 
  = \TTa \cpsia r^{\al} \ws.
\eeq

Recall $\phi_{\e} \phi \gtr_{\e} \ang \xx^{1/3} \gtr r^{\al}$ from \eqref{def:wg_s}. Using \eqref{eq:thm_blow_pf1} and
triangle inequality,  we obtain
\beq\label{eq:thm_blow_pf3}
\bal
\nlinf{ (\Om(s) - \ws ) \ang \xx^{1/3}  } &\les_{\e} \d e^{- \f35 \lam s},  \
\nlinf{ \Om^{\th}( s) -  \omsth }
 = \nlinf{ r^{\al} ( \Om(s) - \ws)}
  \les_{\e} \d e^{-\f35 \lam s} .
\eal 
\eeq
Since $\Om(0) \in C_c^{\infty}(\R^3)$ by \eqref{eq:thm_blow_pf1}, 
 the initial data for \eqref{eq:Euler} satisfies $\om_0^{\th} = \TTa \cpsia \Om_0^{\th} =\TTa \cpsia \Om r^{\al} \in C_c^{\al} $. 
We choose $g = \TTa \cpsia \Om(0) \in C_c^{\infty}$ in Theorem \ref{thm:main_blowup}.
 Estimate \eqref{eq:thm_blow_pf3} implies \eqref{eq:init_ass}.

Using $\CCl= \cE_l(s)  (T-t(s)  )^{\cxs} $ \eqref{eq:scal_pf1},  identities \eqref{eq:thm_blow_pf2},  \eqref{eq:thm_blow_pf3}, and $\omsth = r^{\al} \ws \in L^{\infty}$ by
\eqref{eq:3D_solu_decay}, $\cpsia \asymp 1$ \eqref{def:cpsi},  and $\d \leq 1$, we prove 
\[
\bal
  & \nlinf{ (T-t(s)) \om^{\th} ( \cE_l(s) \cdot (T-t(s)  )^{\cxs} \xx, t )  -   \omthss  }  \\
 & \qquad =\cpsia \nlinf{  (T-t(s))  
\CCw(s) \Om^{\th}(s) - \TTa  \omsth  }  
 =  \nlinf{ \TTa e^{  (O( C_\e \d e^{-3/5 \lam s} )}  \Om^{\th}(s) - \TTa  \omsth  }   \\
& \qquad \les_{\e}  \nlinf{    \Om^{\th}(s) -   \omsth }
+ \nlinf{  ( e^{  (O( C_\e \d e^{-3/5 \lam s} )} -1)   \omsth } 
  \les_{\e}   \d e^{-3/5 \lam s} 
 \les_{\e} \d (T - t(s))^{  \f{ 3 \lam}{ 5 |\bcomth |} }.
\eal
\]
Since $-\bcomth = \TTa^{-1} \in [7, 10]$, using $ \lam = 0.2$ from Lemma \ref{lem:tran_wg},
we obtain $ \f{ 3 \lam}{ 5 |\bcomth |} \geq 0.01$ and prove \eqref{eq:main_SS_blowup}. Since 
$\pa_z \QQs^{z}(0)\neq 0, \ws , \omthss \not \equiv 0$
by Proposition \ref{prop:3D_profile}, we obtain finite-time blowup at $T$ from 
the above estimate.  Using \eqref{eq:scal_pf1} and \eqref{eq:thm_blow_T}, we prove 
$ \cE_l(s) \asymp_{\e} 1, |\cE_l-1| \les_{\e} \d$. Using \eqref{eq:thm_blow_T2}, \eqref{eq:thm_blow_pf3} for $| \Om^{\th}_0 - r^{\al} \ws|$, we prove 
\eqref{eq:thm_blowup_time}, \eqref{eq:init_ass}. We complete the proof of  Theorem \ref{thm:main_blowup}. 
\end{proof}

\subsubsection{A fixed point problem}\label{sec:dyn_fix_pt}

In this section, we reformulate solving \eqref{eq:split} globally in self-similar time $s$ as a fixed point problem. We fix the initial data $\ww_{1,\init}$ as in Theorem \ref{thm:non_stab}. 

Recall the linear map $\cH_2$ from \eqref{def:cH2}, \eqref{eq:w2_form}.   Given an input $\hat \ww_1 \in Y$, we first define 
\bseq\label{eq:non_fix:V2}
\begin{equation}
\ww_2 = \cH_2 \hat \ww_1 .
\end{equation} 
 The above equation is equivalent to
\beq
    \pa_s \ww_2 + \QQs \cdot \na \ww_2  =  \cBs \cdot \ww_2 
  + \cU(\ww_2) + \cU_h(\hat \ww_1) ,\quad \ww_2(0) = \cH_2 \hat \ww_1(0).
\eeq
\eseq

 Second, we define $\ww_1 $ as the solution to the nonlinear equation ~\eqref{eq:split:a}, namely
\begin{subequations}\label{eq:non_fix}
\begin{align}\label{eq:non_fix:a}
  \pa_s \ww_1 + (\QQs + \td \QQ(\ww)) \cdot \na \ww_1 & =  \cBs \cdot \ww_1 
  + \cU(\ww_1) - \cU_h(\ww_1) + \tcB( \ww ) \cdot \ww -  \td \QQ(\ww) \cdot \na \ww_2 ,
\end{align}
with $\ww$ and initial data for $\ww_1$ given by 
\beq\label{eq:non_fix:init}
\ww = \ww_1 + \ww_2,\quad 
  \ww_{1}(0) = \ww_{1,\init} + \cH_{\init} \ww_2(0).
\eeq
\end{subequations}

We justify that there exists a global solution $\ww_1$ to \eqref{eq:non_fix:a} in Section 
\ref{sec:dyn_EE}.

Concatenating the above two steps  defines a map with input $\hat \ww_1$ and output the solution of \eqref{eq:non_fix}: 
\beq
\ww_1 
\stackrel{\eqref{eq:non_fix}}{=} \cH \hat \ww_1. 
\label{def:cH}
\eeq

The proof of Theorem~\ref{thm:non_stab} is broken down in two steps, according to Proposition~\ref{prop:onto} (which shows that the map $\cH$ maps the ball of radius $\d$ in $Y$ into itself), and Proposition~\ref{prop:dyn_compact}.

\begin{prop}\label{prop:onto}
There exists a positive $\d_0(\e) \ll 1$ and $h_4 \leq h_3$ depending on $\e$ such that for any $\d < \d_0(\e)$ and any 
$\hat \ww_1$ with $\| \hat \ww_1 \|_Y \leq \d$, the solutions $\ww_1= \cH \hat \ww_1$ and $\ww_2 = \cH_2 \hat \ww_1$ satisfy
\begin{align}\label{eq:dyn_onto}
 \| \phi \phi_{\e} \ww_1(s) \|  \leq \tf{1}{2} \d  e^{- \f{9}{10} \lam s }  ,
 \quad  \| \ww_2 \|_{\cZ} \les C(\e) e^{ - \f35 \lam  s } \d ,
 \quad \forall \, s \geq 0.
\end{align}
In particular, $   \| \cH \hat \ww_1 \|_{Y} = \| \ww_1 \|_Y \leq \tf{1}{2} \d $. 

\end{prop}

\begin{prop}\label{prop:dyn_compact}
 The map $\cH : \bar B(\d, Y) \to \bar B(\d / 2, Y) \subset \bar B(\d, Y)$ is continuous and 
$\cH(\bar B(\d, Y))$ is  compact with respect to the $Y$-norm.
\end{prop}

We defer the proof of Proposition \ref{prop:onto} to Section \ref{sec:dyn_EE}, 
and the proof of Proposition \ref{prop:dyn_compact} to Section \ref{sec:dyn_comp}.
Using these two propositions, we prove Theorem \ref{thm:non_stab}.

\begin{proof}[Proof of Theorem \ref{thm:non_stab}]
Clearly, $\bar B(\d, Y)$ defined in \eqref{def:YY} is a closed, convex, and nonempty set, with 
respect to $Y$-norm. Proposition \ref{prop:onto} proves that  $\cH: \bar B(\d, Y) \to \bar B(\d/2, Y)\subset \bar B(\d, Y)$, and Proposition \ref{prop:dyn_compact} proves that  $\cH: \bar B(\d, Y) \to \bar B(\d/2, Y)\subset \bar B(\d, Y)$ is continuous, and $\cH(  \bar B(\d, Y))$ is compact.
Applying Schauder fixed point theorem to the map $\cH$ and domain $\bar B(\d, Y)$, we construct a fixed point $\ww_1$ to \eqref{eq:non_fix}: $\cH \ww_1 = \ww_1 \in \bar B(\d, Y)$. Moreover, the solution $\ww_1$ and $\ww_2 = \cH_2 \ww_1$ solve \eqref{eq:split} for any $s \geq 0$.
The parameters  $ \d , h_3, h_4$ are determined in Proposition \ref{prop:onto}. Estimates in \eqref{eq:dyn_onto} imply estimates \eqref{eq:solution:small}. We prove Theorem \ref{thm:non_stab}.
\end{proof}

The key step in the construction is Proposition \ref{prop:onto}. We then prove continuity, compactness, and apply the Schauder fixed-point theorem, which is better suited to low-regularity data.

\subsection{Energy estimates}\label{sec:dyn_EE}

Firstly, we justify the solvability of $\ww_1, \ww_2$ in \eqref{eq:non_fix}. Applying Lemma \ref{lem:semi_cH2}, we prove the estimate of $\ww_2$ in Proposition \ref{prop:onto}
\beq\label{eq:EE_est:w2}
   \cZZ{ \ww_2 }  = \cZZ{ \cH_2 \hat \ww_1 } \les_{  \e } e^{- \f35 \lam s}\| \hat \ww_1 \|_Y \leq C(\e)  e^{- \f35 \lam s}\ \d.
\eeq

Combining the above estimate and \eqref{eq:cH_init},
 and by first choosing $h_3=h_3(\e)$ small enough and then $h_4 = h_4(h_3, \e)$ small enough, we obtain
\beq\label{eq:w2_init_bd}
  \nlinf{ \phi \phi_{\e} \cH_{\init} \ww_2(0) }
  \leq ( C(\e, h_3) h_4^{1/2} + C h_3^{\e^3} ) C(\e) 
  \cZZ{ \ww_2(0) }
  \leq \tf{1}{8}   \| \hat \ww_1 \|_Y \leq \tf18 \d.
\eeq

Introducing  $\Om = \ww_1 + \ww_2 + \ws$ and using \eqref{eq:non_fix}, 
we rewrite the equation of unknown $\ww_1$ in \eqref{eq:non_fix:a}
as that of $\Om$:
\beq\label{eq:non_fix:equiv}
\bal
   \pa_{s}  \Om + (\QQs + \td \QQ(\ww) \cdot \na \Om & = ( c_{\om} - (1-\al) \Psi_z ) \Om
  + \cU_h( \hat \ww_{ 1} + \ww_{2} + \ws -  \Om ) , \\
  \quad \Om(0) & = \ww_{1,\init} +\ws +  \Xi_{h_4} \ast (  \chi_{2, h_3} \ww_2(0)).
\eal
\eeq
The only difference between \eqref{eq:non_fix:equiv} and the original equation \eqref{eq:dyn_full} 
is the mismatch term $\cU_h( \hat \ww_1 - \ww_1)= \cU_h(\hat \ww_{ 1} + \ww_{2} + \ws -  \Om ) $. Using \eqref{eq:EE_est:w2}, \eqref{eq:ass_dyn_init}, and \eqref{eq:w2_regu_init}, we estimate the 
\emph{known} function $f$
\beq\label{eq:dyn_EE_force}
  f =  \hat \ww_{ 1} + \ww_{2} + \ws,
  \quad  \nlinf{ \phi \phi_{\e}^2 f }
  \les \d + 1 \les 1.
\eeq

Using Lemma \ref{lem:semi_cH2}, \eqref{eq:ass_dyn_init}, \eqref{eq:w2_regu_init}, we obtain $\Om(0) \in C_c^{1, \al}$. 
Since $\Om(0)$ is odd in $z$, we get $\Om(0) \in \CCs^{1,\al}$ for the $\CCs^{1,\al}$-norm defined 
in \eqref{def:cC1}. %and $\pa_z \Om(0) = \pa_z \ws(0) $.  
Using Proposition \ref{prop:LWP} with the above force $f$, \eqref{def:cC1},  and initial data $\Om(0)$, we construct a local-in-time $\CCs^{1,\al}$ solution to \eqref{eq:non_fix:equiv} for some $T >0$.
Moreover, using \eqref{eq:LWP_EE} on $\Om$, the $\ww_2$-equation \eqref{eq:non_fix:V2}, \eqref{eq:Uh_bound} for $\cU_h \hat \ww_1$, 
and estimates \eqref{eq:stab_comp_pf2}, we obtain
\beq\label{eq:boot_LWP1}
\bga
  \nlinf{ ( |\xx|^{-\al-1} +|\xx|^{\als} ) \pa_t \Om }
   \les_{\e} C(T, \| \Om_0 \|_{\CCs^{1,\al}} ), \\
    \nlinf{ \phi \pa_t \ww_2} 
  \les_{\e}  \cZZ{\ww_2} + \nlinf{ \phi \cU_h \hat \ww_1}
  \les \d + \nlinf{ \phi \phi_{\e}^2 \hat \ww_1 } \les \d.
  \ega
\eeq
Using $\ww_1 = \Om -\ww_2 - \ws$, we construct a local-in-time solution 
with
\beq\label{eq:LWP_w1}
   \nlinf{ \phi \pa_t \ww_1  } \les_{\e} C(T, \| \Om_0 \|_{\CCs^{1,\al}} ) + \d,
   \quad \ww_1(t) \in C^1(\R^3) \cap L^{\infty}(\phi),  \quad \forall t \in [0, T].
\eeq

\subsubsection{Energy estimates of $\ww_1$}

In this section, we perform energy estimates on $\ww_1$. 
Using the assumption \eqref{eq:ass_dyn_init} for $\ww_{\init}$, \eqref{eq:w2_init_bd}, \eqref{eq:non_fix:init}, and $\phi_{\e} \leq 1$ by \eqref{def:wg_s}, we obtain 
\beq\label{eq;w1_init_bd}
    \nlinf{ \phi \phi_{\e} \ww_1(0)} \leq 
        \nlinf{ \phi  \ww_1(0) } +   \nlinf{ \phi \phi_{\e} \cH_{\init} \ww_2(0) } \leq \tf18 \d + \tf18 \d \leq \tf14 \d.
\eeq

 First, we impose the following bootstrap assumption on $\ww_1$:
\beq\label{eq:dyn_boot1}
    \nlinf{ \phi \phi_{\e} \ww_1 } \leq  \d e^{- \f{9}{10} \lam \e} ,
\eeq
for $s \in [0, s_*]$, where $s_*$ is the maximal time where \eqref{eq:dyn_boot1} holds.
From \eqref{eq:LWP_w1} and \eqref{eq;w1_init_bd}, the above conditions hold with some $s_* > 0$. 
Moreover, by the blowup criterion \eqref{eq:LWP_blowup}, tbe estimates  \eqref{eq:EE_est:w2},
\eqref{eq:3D_solu_decay}  for $\ww_2$, and $\Om = \ww_1 + \ww_2 + \ws$, 
under the conditions \eqref{eq:dyn_boot1},  the local solution $\ww_1$ can be extended beyond $s_*$.
Below, we consider $s \in [0, s_*]$, and we shall improve the estimate \eqref{eq:dyn_boot1} to 
\beq\label{eq:dyn_boot1_impr}
      \nlinf{ \phi \phi_{\e} \ww_1 } \leq \tf12 \d e^{- \f{9}{10} \lam \e} , \quad \forall \  s \in [0, s_*].
\eeq
Then by a bootstrap argument, we yield $s_* = \infty$.

Recall $\phim, \phi$ from \eqref{def:wg_s}, $\phi_1$ constructed in Lemma \ref{lem:tran_wg},
and  $\cZ$-norm from \eqref{def:cZ}. Since $\phia \les_{\e} \phi $ and  $\ww= \ww_1 + \ww_2$, by requiring $\d$ small enough, using \eqref{eq:EE_est:w2} for $\ww_2$ and \eqref{eq:dyn_boot1} for $\ww_1$, we obtain
\beq\label{eq:damp_boot}
  \nlinf{ \phim \phi_{\e} \ww(s) } 
  \leq C(\e)  (  \nlinf{ \phi \phi_{\e} \ww_1 }
  +   \nlinf{ \phi \phi_{\e} \ww_2 } ) \leq C(\e) \d \leq \e^3, \quad \forall s \in [0, s_*].
\eeq

Multiplying \eqref{eq:non_fix:a} by $\phi_{\e} \phi$, we obtain 
\[
\bal
  ( \pa_s \ww_1 + (\QQs + \td \QQ(\ww)) ) \cdot \na  ( \ww_1 \phi_{\e} \phi)  & = 
  \B(   \cBs + \f{ (\QQs + \td \QQ(\ww)) \cdot \na (\phi_{\e} \phi) }{ 
\phi_{\e} \phi   } \B) \cdot \ww_1  \\
  & \quad + \phi \phi_{\e} \cdot ( \cU(\ww_1) - \cU_h(\ww_1) + \tcB( \ww ) \cdot \ww 
-  \td \QQ(\ww) \cdot \na \ww_2 )
  \teq I + II.
  \eal
\]

Since \eqref{eq:damp_boot} implies the assumption in Lemma \ref{lem:tran_wg} 
for $\eta = \ww$, 
applying Lemma \ref{lem:tran_wg} 
we obtain
\[
  d_0( \phi \phi_{\e}, \ww ) =
     \cBs + \f{ (\QQs + \td \QQ(\ww)) \cdot \na (\phi_{\e} \phi) }{ 
\phi_{\e} \phi   } \leq - \lam.
\]

We treat $II$ perturbatively. Using \eqref{eq:cU_dif}, the bootstrap bound \eqref{eq:dyn_boot1}, and $\phi_{\e} \leq 1$ by \eqref{def:wg_s}, we bound 
\[
   \nlinf{ \phi \phi_{\e} ( \cU(\ww_1) - \cU_h(\ww_1) ) } \leq 10^{-4} \lam \nlinf{ \phi \phi_{\e} \ww_1 }  \leq 10^{-4} \lam  \cdot \d e^{- \f{9}{10} \lam s}  \, .
\]

Recall $ |\phi \ww_2 | \leq \cZZ{\ww_2}$ by \eqref{def:cZ} and $\phi_{\e} \leq 1$. Applying \eqref{eq:stab_comp_pf2:cB} and $\ww = \ww_1 + \ww_2$, estimate \eqref{eq:EE_est:w2} for $\ww_2$, and \eqref{eq:dyn_boot1} for $\ww_1$, we bound 
\[
   \nlinf{ \phi \phi_{\e} \tcB( \ww ) \cdot \ww }
  \les_{\e} \nlinf{ \phi \phi_{\e}^2 \ww } \nlinf{ \phi \phi_{\e} \ww } 
\les_{\e} ( \d e^{-\f{9}{10} \lam s } + \d e^{- \f35 \lam s}  )^2 
\les_{\e} \d^2 e^{- \f{6}{5} \lam s}.
\]

Applying \eqref{eq:stab_comp_pf2:cB} and the definition of $\cZ$-norm \eqref{def:cZ}, we bound 
\[
\bal
  \nlinf{ \phi \phi_{\e} \td \QQ(\ww) \cdot \na  \ww_2 }
  & \les 
\nlinf{ \phi |\xx| \na \ww_2 }  \nlinf{\td \QQ(\ww)/ |\xx| }
\les_{\e}  \cZZ{\ww_2} \nlinf{ \phi \phi_{\e}^2 \ww } \\
& \les_{\e} ( \d e^{-\f{9}{10} \lam s } + \d  e^{-\f35 \lam s}  )   \cdot  \d  e^{-\f35 \lam s}
\les_{\e} \d^2 e^{-\f65 \lam s} \, .
\eal
\]

Applying the above estimates and following the estimates along the trajectory 
associated with $ \QQs + \td \QQ(\ww)$ in Section \ref{sec:decay_cT}, we prove 
\[
 \tf{1}{2} \cdot \tf{d}{d s} | \ww_1 \phi \phi_{\e}|^2 \cc X(s, \xx) 
\leq - \lam | \ww_1 \phi \phi_{\e}|^2 \cc X(s, \xx)
+ ( 10^{-4} \lam 
\cdot \d e^{- \f{9}{10} \lam s}  + C_{\e} \d^{2} e^{-\f65 \lam s} ) \cdot  | \ww_1 \phi \phi_{\e}| \cc X(s, \xx)  .
\]
Cancelling $| \ww_1 \phi \phi_{\e}| \cc X(s, \xx) $, applying Gronwall inequality, and taking supremum over $\xx$, we prove 
\[
  \nlinf{\ww_1(s) \phi \phi_{\e} } \leq 
  e^{-\lam s} \nlinf{ \ww_{1}(0) \phi \phi_{\e} }
  +  \int_0^s e^{- \lam (s -\tau) }
   ( 10^{-4} \lam  \cdot \d e^{- \f{9}{10} \lam \tau }
   + C_{\e} \d^{2} e^{-\f65 \lam \tau} ) d \tau .
\]

Using \eqref{eq;w1_init_bd} and the above estimate, we prove 
\[
    \nlinf{\ww_1(s) \phi \phi_{\e} } \leq \tf14 \d e^{-\lam s}  + 10^{-4} \lam \cdot ( \f{\lam}{10} )^{-1} \cdot e^{ \f{1}{10} \lam s } \cdot e^{-\lam s} \d + C_{\e} \d^{2} e^{-\lam s}
  \leq ( \f38 \d + C_{\e} \d^{2} ) e^{- \f{9}{10}\lam s}.
\]

By choosing $\d_0(\e)$ small enough, we obtain 
\[
  \tf38 \d + C_{\e} \d^{2}  \leq \tf12 \d , \quad \forall \ \d \leq \d_0(\e),
  \]
and prove \eqref{eq:dyn_boot1_impr}, which improves \eqref{eq:dyn_boot1}.  Thus, the bootstrap assumption \eqref{eq:dyn_boot1} hold for $s \in [0, \infty)$.

Combining the above estimate and \eqref{eq:EE_est:w2}, we prove 
estimate \eqref{eq:dyn_onto}. 
From the definition of $Y$-norm, estimate \eqref{eq:dyn_onto} implies 
$\| \ww_1 \|_Y \leq \f{1}{2} \d$ for  $\ww_1 = \cH \hat \ww_1$. We prove Proposition \ref{prop:onto}.

\subsection{Continuity and Compactness of the map}

In Proposition \ref{prop:onto}, we have proved $\cH : \bar B_{\d}(Y) \to \bar B_{\d/2}(Y)$.
In this section, we further prove that $\cH$ is compact and continuous. 

\subsubsection{Compactness}\label{sec:dyn_comp}

Consider a sequence of input $\hat \ww_i$ with  $ \| \hat \ww_{i} \|_{Y} \leq \d$. Denote $ \ww_{i, 1} = \cH \hat \ww_{i, 1}$ 
and $\ww_{i, 2} = \cH_2 \hat \ww_{i, 2}$. We shall prove that there exists a subsequence of 
$\ww_{i, 1}$ that convergences to some limit in $Y$. Using Proposition \ref{prop:onto}, we obtain
\beq\label{eq:dyn_comp_pf1}
   \nlinf{ \phi_{\e} \phi \ww_{i, 1} } \leq \tf{1}{2} \d e^{-\f{9}{10} \lam s },
   \quad \|\ww_{i, 1} \|_Y \leq \tf{1}{2 } \d ,
   \quad  \| \ww_{i, 2} \|_{\cZ} \les C(\e) e^{ - \f35 \lam  s } \d .
\eeq

Denote by $\Om_i$ 
\beq\label{eq:dyn_Omi}
   \Om_i = \ww_{i,1} + \ww_{i,2} + \ws, \quad \ww_i = \ww_{i, 1} + \ww_{i, 2}, 
\eeq
which solves \eqref{eq:non_fix:equiv}
\beq\label{eq:dyn_Omi_eqn}
  \pa_{s}  \Om_i + (\QQs + \td \QQ(\ww_i) \cdot \na \Om_i = ( c_{\om} - (1-\al) \Psi_z ) \Om_i
  + \cU_h( \hat \ww_{i, 1}  + \ww_{i, 2} + \ws - \Om_i ).
\eeq

Let $\phia$ be the weight in \eqref{def:cC1}. The above conditions and the estimates on $\ws$  in \eqref{eq:3D_solu_decay} implies 
\beq\label{eq:dyn_cont}
\nlinf{  \ang \xx^{-2 \e^3} \phia  \Om_i } 
\les_{\e}  \nlinf{ \phi \phi_{\e} \Om_i } \les_{\e} 1 ,\quad \forall s \geq 0.
\eeq

\paragraph{\bf Improved regularity of $\ww_{i, 1}$}

Recall the function $\ww_{1,\init}$ from \eqref{eq:ass_dyn_init}, which is fixed in the fixed point problem in \eqref{eq:non_fix}. For each pair $(\ww_{1,i}, \ww_{2, i})$, from the assumption \eqref{eq:non_fix:init} on $\ww_{1,\init}$, the estimate \eqref{eq:w2_regu_init} for $\ww_2$, 
and the definition of $\CCs^{1,\al}$ norm \eqref{def:cC1}, we obtain
\beq\label{eq:dyn_comp_pf2}
\bal
  \Om_i(0)  & = \ww_{i}(0) + \ws = \ww_{1, \init} + \ws +   \Xi_{h_4} \ast (  \chi_{2, h_3}  \ww_2(0))   \in C_c^{1,\al} ,  \\
       \| \Om_i(0) \|_{\CCs^{1,\al}} & \les C(\e , \d) + \| \ww_{1,\init}(0) + \ws \|_{\CCs^{1,\al}} 
   \les C(\e,  \d,  \ww_{1,\init} ). 
\eal 
\eeq

Since $\Om_i$ and the force $f = \hat \ww_{i,1} + \ww_{i, 2} + \ws$
are bounded  by \eqref{eq:dyn_cont}, \eqref{eq:dyn_EE_force}, for any $s_* \geq 0$, using the continuation criterion \eqref{eq:LWP_blowup} and estimate \eqref{eq:LWP_EE} for $\Om_i$, we obtain uniform estimates in $i$ 
\beq\label{eq:dyn_comp_pf3}
\bal
  \sup\nolim_{s \leq s_*}  \| \Om_i(s) \|_{C^{1,\al}} 
  +  \| \pa_s \Om_i(s) \|_{ L^{\infty} } 
  \leq   \sup\nolim_{s \leq s_*}  \| \Om_i(s) \|_{\CCs^{1,\al}} 
  & \leq  
   C(\e, \d,  s_*, \ww_{1,\init} ).
\eal
\eeq
Since $ \ws \in C^{1, \al}(\bar B_R)$ for any $R >0$ by \eqref{eq:3D_solu_decay}, \eqref{eq:3D_solu_reg}
and $ \Om_i = \ww_{i,1} + \ww_{i, 2} + \ws$, we obtain
\[
    \| \ww_{i, 1} + \ww_{i, 2}  \|_{ C^{1,\al}(B_R) }
    + 
    \| \pa_s (\ww_{i, 1 } + \ww_{i, 2 })  \|_{ L^{\infty}(B_R) }
     \les C(\e, \d,  s_*, \ww_{1,\init}, R),
    \quad  \forall \ R > 1 .
\]

Using the definition of $\cZ$ \eqref{def:cZ}, the equation of $\ww_{i, 2}$ \eqref{eq:split:b}, 
estimate similar to \eqref{eq:boot_LWP1}, we obtain that $\pa_s \ww_{i, 2}$ is continuous, and  
\[
  \| \ww_{i, 2}(s) \|_{C^1(B_R)} + \|  \pa_s \ww_{i, 2}(s) \|_{L^{\infty}(B_R)} \les 
 C(\e, \d,  s_*, \ww_{1,\init}, R).
\]
Combining the above two estimates, we obtain the \emph{qualitative improved regularity} of $\ww_{i, 1}$:
\beq\label{eq:w1_regu_impr}
    \| \ww_{i, 1}(s) \|_{C^1(B_R)} + \|  \pa_s \ww_{i, 1 }(s) \|_{L^{\infty}(B_R)} \les 
     C(\e, \d,  s_*, \ww_{1,\init}, R).
\eeq

\paragraph{\bf Compactness}
Now, for any $n\geq 1$, consider the region
\[
  \S_n = \{ (s, \xx) : s \in [0, n], |\xx| \in [1/n, n] \}.
\]

By \eqref{eq:w1_regu_impr}, $\ww_{i, 1}$ is equi-continuous in $\S_n$ for any $n \geq 1$. Using Arzela-Ascoli Theorem and a diagonalization argument, there exists 
a continuous function $\ww_{\infty}$ and a subsequence $ \ww_{m_i, 1}$ such that 
\[
  \lim\nolim_{i \to \infty} \| \ww_{m_i, 1} - \ww_{\infty} \|_{ L^{\infty}(\S_n ) }  = 0, \quad \forall \ n \geq 1,
\]
Recall the weight $\phi_{\e}, \phi$ from \eqref{def:wg_s}. The above estimates imply 
\[
  \lim\nolim_{i \to \infty}  \| \phi_{\e}^2 \phi e^{ \f{4}{5} \lam s } ( \ww_{m_i, 1} - \ww_{\infty} )
  \|_{ L^{\infty}(\S_n ) } = 0, \quad \forall \ n \geq 1.
\]

Moreover, the uniform convergence and \eqref{eq:dyn_comp_pf1} implies $\nlinf{ \phi_{\e} \phi \ww_{\infty} } \leq \tf{1}{2} \d e^{- \f{9}{10} \lam s}$. 
By definition of $\phi_{\e}$ in \eqref{def:wg_s} and the $Y$-norm in \eqref{def:YY}, 
we have 
\[
\bal
\| \phi_{\e}^2 \phi e^{ \f{4}{5} \lam s } ( \ww_{m_i, 1} - \ww_{\infty} )
  \|_{ L^{\infty}(\S_n^c ) }
 &  \les \sup\nolim_{ (s, \xx) \notin \S_n } \phi_{\e}(\xx) e^{- \f{1}{10} \lam s }
  \| \phi_{\e} \phi e^{ \f{9}{10} \lam s } ( | \ww_{m_i, 1} | + |\ww_{\infty}| )  \|_{L^{\infty}} 
  \les ( n^{- \e^3} + e^{ - \f{1}{10} \lam n } )  \d .
\eal
\]

Combining the above two estimates and \eqref{eq:dyn_comp_pf1}, we prove 
\[
  \| \ww_{m_i, 1} - \ww_{\infty} \|_{Y}  \to 0, \ \mbox{as \ } \ i \to \infty,
  \quad \| \ww_{\infty} \|_Y \leq \tf12 \d. 
\]
Thus, $\cH(\bar B(\d, Y))$ is compact.

\subsubsection{Closedness of $\cH$}

In this section, we prove that the map $\cH : \bar B(\d, Y) \to \bar B(\d/2, Y) $ is closed. 
Consider a sequence $\hat \ww_{i, 1} \in \bar B(\d, Y)$, $\ww_{i, 1} = \cH \hat \ww_{i, 1}, \ww_{i, 2} = \cH_2 \hat \ww_{i, 1}$ and suppose that they satisfy
\[
  \hat \ww_{i, 1} \to \hat \ww_{\infty},
  \quad \ww_{i, 1} \to \ww_{\infty, 1}  \, ,
\]
in the $Y$-norm \eqref{def:YY}. We shall prove $\ww_{\infty, 1} = \cH \hat \ww_{\infty}$. We define 
\beq\label{eq:close_pf1}
\bal
  \ww_{\infty, 2} =  \cH_2 \hat \ww_{\infty}, 
  \  \ww_i = \ww_{i, 1} + \ww_{i, 2} ,  
  \  \Om_i = \ww_i + \ws,
  \  \ww_{\infty} =  \ww_{\infty, 1} +  \ww_{\infty, 2} , \quad 
     \Om_{\infty} \teq \ww_{\infty} + \ws \, ,
\eal
\eeq
and 
\[
  \EEs_i \teq \| \ww_{i, 1} - \ww_{\infty, 1} \|_Y
  + \| \hat \ww_{i, 1} - \hat \ww_{\infty, 1} \|_Y
  = \sup\nolim_{s \geq 0} e^{ \f45 \lam s} \nlinf{\phi \phi_{\e}^2 (\ww_{i, 1} - \ww_{\infty, 1}) }
   + \sup\nolim_{s \geq 0} e^{ \f45 \lam s} \nlinf{\phi \phi_{\e}^2 ( \hat \ww_{i, 1} - \hat \ww_{\infty, 1}) }.
\]
The convergence implies $\EEs_i \to 0$ as $i\to \infty$.

Since  $\cH_2$ \eqref{def:cH2} and $\cH_{\init}$ \eqref{eq:w1_init} are linear,
using Lemma \ref{lem:semi_cH2} and estimate \eqref{eq:w2_init_bd}, we obtain
\beq\label{eq:dyn_close_pf0}
\bal
\cZZ{ ( \ww_{i, 2} - \ww_{\infty, 2} )(s) } &= \| \cH_2( \hat \ww_{i, 1} - \hat \ww_{\infty} )(s) \|_{\cZ} 
   \les_{\e} e^{ - \f35 \lam s } \|  \hat \ww_{i, 1} - \hat \ww_{\infty}  \|_Y
  \les_{\e} %e^{ - \f35 \lam s } 
   \EEs_i  \, , \\  
    \nlinf{ \phi \phi_{\e} \cH_{\init} ( \hat \ww_{i, 1} - \hat \ww_{\infty} )  }
 & \leq \tf{1}{8} \|  \hat \ww_{i, 1} - \hat \ww_{\infty}  \|_Y 
 \les  \EEs_i.
\eal
\eeq

From Proposition \ref{prop:onto}, $\ww_{i, 1}, \ww_{i,2}$ satisfy the estimates in \eqref{eq:dyn_comp_pf1}. From Section \ref{sec:dyn_comp}, $\Om_i$ satisfies estimates \eqref{eq:dyn_comp_pf2},
\eqref{eq:dyn_comp_pf3}. Note that these estimates are uniformly in $i$. 
Integrating \eqref{eq:dyn_Omi_eqn} in $s$ 
and using $\hat \ww_{i,1} + \ww_{i, 2} + \ws - \Om_i = \hat \ww_{i, 1} - \ww_{i, 1}$, we obtain %that $\Om_i$ solves 
\beq\label{eq:dyn_close_pf1}
    \Om_i(s) - \Om_i(0) =  \int_0^s \B( -( \QQs + \td \QQ(\ww_i) ) \cdot \na \Om_i + (\cBs + \tcB( \ww_i ) ) \Om_i
  + \cU_h( \hat \ww_{i, 1} - \ww_{i, 1} )  \B)(\tau) d \tau \, ,
\eeq
for any $s \geq 0$. Next, we aim to take $i \to \infty$.

 Using triangle inequality, the definitions of $\cZ$-norm in \eqref{def:cZ}, we obtain
\beq\label{eq:close_linf}
\nlinf{ \phi_{\e}^2 \phi ( \Om_i - \Om_{\infty}) }
  = \nlinf{ \phi_{\e}^2 \phi  ( \ww_i - \ww_{\infty}  )  } 
  \leq  \nlinf{ \phi_{\e}^2 \phi ( \ww_{i, 1} - \ww_{\infty,1} ) }  + 
  \cZZ{ \ww_{i, 2} - \ww_{\infty, 2} } \les_{\e} %e^{ - \f35 \lam s } 
   \EEs_i.
\eeq

Since $\phi_{\e}^2 \phi \gtr 1$ by \eqref{def:wg_s}, using the above convergence and the uniform bound
in \eqref{eq:dyn_comp_pf3}, we obtain
\beq\label{eq:Ominf_reg}
\Om_{\infty}(s) \in \CCs^{1,\al} \cap C^{1,\al},\quad   \| \Om_{\infty}(s) \|_{C^{1,\al}}
+ \| \Om_{\infty}(s) \|_{\CCs^{1,\al}}
\les_{\e, \d} C(s_*, \ww_{1,\init}), \quad \forall s \in [0, s_*].
\eeq

Using standard interpolation, $1\les_{\e} \phi \phi_{\e}^2$ by \eqref{def:wg_s}, and the uniform bound \eqref{eq:dyn_comp_pf3}, we obtain
\[
\bal
 \sup_{s \leq s_*} \| \na (\Om_i  - \Om_{\infty} )(s) \|_{ L^{\infty}(\bar B_R) }
 \les  \sup_{s \leq s_*}  \| (\Om_i  - \Om_{\infty} )(s) \|_{ L^{\infty}(\bar B_{2R}) }^{ \f{1}{1+\al} }
  \| (\Om_i  - \Om_{\infty} )(s) \|_{ C^{1,\al}(\bar B_{2R} ) }^{ \f{\al}{1 + \al} }
  \les  \EEs_i^{ \f{1}{1+\al} } C(\e, s_*, \d, \ww_{1,\init} ) \, .
\eal
\]
The above \emph{qualitative} bound \(C(\e,s_*,\d,\ww_{1,\init})\) in $s$ suffices for our purpose. Since $\td \QQ$ and $\cU, \cU_h$ are linear, using the above convergence, estimate \eqref{eq:stab_comp_pf1:a} for $\td \QQ, \tcB$, and \eqref{eq:Uh_bound} for $\cU_h$, we obtain
\[
\bal
  \| \td \QQ(\ww_{i} - \ww_{\infty}) \|_{C^{\al}(B_R)} 
  +   \| \tcB(\ww_{i} - \ww_{\infty}) \|_{C^{\al/2}(B_R)} 
 &  \les_{\e, R} \nlinf{  \phi_{\e}^2 \phi (\ww_i -\ww_{\infty}) } \les 
 \EEs_i, \\
  \| \cU_h(  \hat \ww_{i, 1} - \ww_{i, 1} )
  - \cU_h(  \hat \ww_{\infty, 1} - \ww_{\infty, 1} ) \|_{C^{\al/2}(B_R)} 
  & \les \nlinf{ \phi_{\e}^2 \phi ( \hat \ww_{i, 1} -\hat \ww_{\infty, 1} )}
  + \nlinf{ \phi_{\e}^2 \phi (  \ww_{i, 1} - \ww_{\infty, 1} )} \les 
   \EEs_i.
\eal
\]

Moreover,  from  \eqref{eq:stab_comp_pf1:a}, $\td \QQ(\ww_i) , \QQs, \cBs, \tcB(\ww_i)$ are bounded for $\xx \in B_R$  by uniformly in $i$ and $s$.

We fix $R> 0$ and $s_* > 0$. For any $\xx \in \bar B_R, s \in [0, s_*]$, 
using the above uniform convergence in $[0, s_*] \times \bar B_R$,  we obtain  the uniform convergence 
as $i\to \infty$: 
\[
\bal
  (\QQs + \td \QQ(\ww_i) \cdot \na \Om_i & \to  (\QQs + \td \QQ(\ww_{\infty}) \cdot \na \Om_{\infty}  , & \quad 
(\cBs + \tcB( \ww_i ) ) \Om_i  & \to (\cBs + \tcB( \ww_{\infty} ) ) \Om_{\infty}, \\
   \cU_h( \hat \ww_{i, 1} - \ww_{i, 1} ) & \to
      \cU_h( \hat \ww_{\infty, 1} - \ww_{\infty, 1} ),
      & \quad  \Om_i &  \to \Om_{\infty}  .
  \eal
\]
Thus, taking $i \to \infty$ in \eqref{eq:dyn_close_pf1}, 
for $s \in [0, s_*] , \xx \in \bar B_R$, we prove 
\[
      \Om_{\infty}(s) - \Om_{\infty}(0) =  \int_0^s \B( -(\QQs + \td \QQ(\ww_{\infty}) \cdot \na \Om_{\infty} + (\cBs + \tcB( \ww_{\infty} ) ) \Om_{\infty}
  + \cU_h( \hat \ww_{\infty, 1} - \ww_{\infty, 1} )  \B)(\tau) d \tau  \, .
\]
 Since $s_*, R>0$ are arbitrary and the integrand is continuous in $\bar B_R$, taking $s$-derivative, using the regularity on $\Om_{\infty}$ \eqref{eq:Ominf_reg}, we obtain that $\Om_{\infty}$ is a $\CCs^{1,\alpha}$ solution to 
\[
        \pa_s \Om_{\infty}(s) + (\QQs + \td \QQ(\ww_{\infty}) \cdot \na \Om_{\infty} = (\cBs + \tcB( \ww_{\infty} ) ) \Om_{\infty}
  + \cU_h( \hat \ww_{\infty, 1} + \ww_{\infty, 2} + \ws  - \Om_{\infty} ) .
\]
where we use \eqref{eq:close_pf1} to rewrite $\hat \ww_{\infty, 1} + \ww_{\infty, 2} + \ws  - \Om_{\infty} = \hat \ww_{\infty, 1} - \ww_{\infty,1}$. 

Finally, we prove that $\Om_{\infty} = \ww_{\infty,1} + \ww_{\infty,2}+ \ws $ and $\eta_{\infty} \teq \cH \hat \ww_{\infty} + \ww_{\infty,2}+ \ws$ 
are the same.  Recall the initial data from \eqref{eq:dyn_comp_pf2}
\beq\label{eq:close_pf2}
  \Om_i(0) = \ww_i(0) + \ws = \ww_{1, \init} + \ws +   \Xi_{h_4} \ast (  \chi_{2, h_3}  \ww_{i,2}(0))  .
\eeq

Applying estimates \eqref{eq:dyn_close_pf0} and \eqref{eq:w2_regu_init} 
to the right hand side 
and the convergence \eqref{eq:close_linf} to the left hand side, and taking $i \to \infty$, we prove
\[
\Om_{\infty}(0) =  \ww_{\infty}(0) + \ws  = \ww_{1, \init} + \ws +   \Xi_{h_4} \ast (  \chi_{2, h_3}  \ww_{\infty,2}(0)).
\]

From the definition of the $\cH$-map in \eqref{eq:non_fix}, 
given the input $\hat \ww_{\infty, 1}$, we construct $\ww_{\infty, 2} = \cH_2 
\hat \ww_{\infty, 1}$ as in \eqref{eq:non_fix:V2}. 
Then from Proposition \ref{prop:LWP},  $\eta_{\infty} = \cH \hat \ww_{\infty,1}
+ \ww_{\infty,2} + \ws$ is the unique $\CCs^{1,\al}$ solution to \eqref{eq:non_fix:equiv} 
from initial data $ \ww_{1, \init} + \ws +   \Xi_{h_4} \ast (  \chi_{2, h_3}  \ww_{\infty,2}(0)) $.
Since $\Om_{\infty} \in \CCs^{1,\al}$ constructed above solves the same equation with the same initial data, by uniqueness in Proposition \ref{prop:LWP}, we establish
\[
  \Om_{\infty} = \eta_{\infty}  \  \Longrightarrow  \ \ww_{\infty, 1} = \cH \hat \ww_{\infty, 1}.
\]
We prove the closedness of the map $\cH$. Since $\cH$ is compact and closed, following the proof of Theorem \ref{thm:3D_contin}, we prove that $\cH: \bar B(\d, Y) \to \bar B(\d/2, Y)$ is continuous. 
Combining the proof of compactness in Section \ref{sec:dyn_comp} and the above proof, 
we prove Proposition \ref{prop:dyn_compact}.

\vs{0.1in}

\paragraph{\bf Acknowledgments.}
JC was partially supported by NSF Grant DMS--2408098. 
He would like to thank Dr.~Shumao Zhang for insightful discussions on the numerical results in \cite{hou2024potential}.

\appendix

\section{Basic functional inequalities and estimates}\label{app:basic}

In this Appendix, we develop several basic functional inequalities and estimates.

 Recall the norms from \eqref{def:3d_wg}, \eqref{def:3d_norm}.  We control $\linfa$-norm using $\linfza, \wwra$ norms:

\begin{lem}\label{lem:linf_R3}  
Let $\rag$ be as in \eqref{def:rag} with parameter $\kag> 0$.
Suppose that weight $\rhoo$ satisfies $\rhoo(\yy) \leq C \rhoo(\xx) $ for any $|\yy| \geq |\xx|$ 
for some absolute constant $C$ independent of $\xx, \yy$.  Let $\nnla{\cdot}, \nnr{\cdot}$ be the norms associated with $\rhoo$ via \eqref{def:3d_wg}, \eqref{def:3d_norm}
For any function $f$, we have 
\bseq\label{eq:linf_R3}
\begin{align}
  |\rhoo  \rag f(r,z) | & \les \kag^{-1}  ( |\rhoo f|(0 , z )| + \max_{\tr \leq r}
   | \rhoo \rag \cdot \ang {\tr, z} \pa_r f(\tr, z) | ), \label{eq:linf_R3:a}   \\
  \nnla{f} & \les \kag^{-1} ( \nnr{f} + \nna{f}) ,
  \label{eq:linf_R3:b} 
\end{align}
\eseq
with absolute constant independent of $ \rhoo,  f, \rag$.

\end{lem}

\begin{proof}

We first prove \eqref{eq:linf_R3:a}. We have
\[
  |f(r, z) \rhoo(\xx)| \leq ( |f(0, z)| + \int_0^r  | \pa_r f(s, z) d s | ) \rhoo(\xx) \teq I_1 + I_2. 
\]
Recall $\rhoo, \rhor $ from \eqref{def:3d_wg:a}. For any $s \in [0, r]$, we have 
$|(s,z) | \leq |\xx|$ and $\rhoo(r, z) \les \rhoo(s, z)$. It follows 
\beq\label{eq:linf_R3:pf1}
\bal
  I_1 &\leq |f(0,z) \rhoo(0,z) | ,  \\ 
  I_2 & \leq \int_0^r \f{ \rhoo(\xx)}{ \rhoo(s, z)}  \rag(s, z)^{-1} \la s, z\ra^{-1} d s 
  \cdot \max_{\tr \leq r} |(\rhoo \rag \ang \xx \pa_r f )(\tr, z)| \\
& \leq \int_0^r    \rag(s, z)^{-1} 
 \la s, z\ra^{-1} d s \cdot \max_{\tr \leq r} |( \rhor \pa_r f )(\tr, z)|
 \teq I_3 \cdot  \max_{\tr \leq r} |( \rhor \pa_r f )(\tr, z)|.
\eal
\eeq

Using the definition of $\rag$  \eqref{def:3d_wg},  $\kag \in (0, \f{1}{100} )$, 
and $r \leq |\xx|$, we obtain
\beq\label{eq:linf_R3:pf2}
\bal
I_3 & \les \ang z^{-\kag} \int_0^r \la s, z \ra^{\kag-1}  d s  
  \les \ang z^{-\kag}  \int_0^r \ang s^{\kag-1} d s    \\
  &  \les  \kag^{-1}  \ang z^{-\kag} \min(r, r^{\kag} )  \les \kag^{-1}  \rag^{-1}(r, z)
  \min(|\xx|, 1) .
\eal
\eeq

Combining the above estimates for $I_1, I_2, I_3$, we prove \eqref{eq:linf_R3:a}. 

Applying estimate \eqref{eq:linf_R3:pf1}, \eqref{eq:linf_R3:pf2} with $f$ replaced by $ \ang z^{\hal} f$, and using the  norms \eqref{def:3d_norm}, we obtain
\beq\label{eq:linf_R3:pf3}
\bal
  | \rhoo  \ang z^{\hal} f(r,z) | 
& \les
 \big( \rhoo \ang z^{\hal} |f| \big)(0, z) + \kag^{-1} \rag^{-1} \min(|\xx|, 1)
 \cdot  \|  \rhor  \ang z^{\hal} \pa_r f \|_{L^{\infty}} \\
&  \les \kag^{-1}   ( \min( |z|, 1) \nna{f} + \rag^{-1} \min(|\xx|, 1)  \nnr{f} ) \\
& \les \kag^{-1} (1 + \rag^{-1} )  \min(|\xx|, 1) ( \nnr{f}  + \nna{f} ).
 \eal
\eeq
Using $\rag \les 1$, definition of $\nnla{\cdot}$-norm, and multiplying  \eqref{eq:linf_R3:pf3}
by $\rag$, we prove \eqref{eq:linf_R3:b}.
\end{proof}

We have the following H\"older estimates for the Poisson equation with a singular source.

\begin{lem}\label{lem:hol_singu}
Consider $n \geq 3$.  Denote $x_h = (x_1, .. , x_{n-1} )$. 
 Suppose that $|f(x)| \leq M |x_h|^{\al-1} $ for some $\al \in (0, 1)$  and $f$ has a compact support. 
Consider  $ \psi = c_n \int |x-y|^{2- n} f(y)  d y$. We have 
\[ 
  \| \na \psi \|_{\dot C^{\al}} \les  M .
\]
where  $\dot C^{\al}$ is the standard H\"older semi-norm.
\end{lem}

\begin{proof}

Taking $\pa_{x_i}$,  we obtain
\[
  \pa_i \psi = \int K(x- y) f(y) d y, \quad K(z) = C_n \f{z_i}{|z|^{n}}.
\]

Given arbitrary $x_1, x_2$. Denote $a = x_1 - x_2$ and we decompose the integral as 
\[
\bal
 & \pa_i \psi( x_1  ) - \pa_i \psi(x_2)
   = \int K(x_1 - y) f(y) - \int K(x_2 - y) f(y) \\
 &   = \int_{|y- x_1| < 2 a} K(x_1 - y) f(y)
  - \int_{|y-x_1| < 2 a}  K(x_2 - y) f(y)
  + \int_{ |y-x_1| > 2 a  } (K(x_1 - y) - K(x_2 - y) ) f(y) d y \\
  & \teq I_1 + I_2 + I_3.
\eal
\]

Since $|x_1 - x_2| = a$, for $|y-x_1| \leq 2 a$, we get $|y-x_2| \leq  4 a$. 
For $j=1,2$, we estimate $I_j$ in the same way and denote $z = x_j$. Using $|f(y)| \les |y_h|^{\al-1}$, and changing $ z = a \td z,  y = a \td y$, 
 we estimate 
\[
  |I_j| \les \int_{ |y-z| \leq 4 a } |y- z|^{1-n}  \cdot |y_h|^{\al-1} d y 
  \les a^{\al} \int_{ |\td y- \td z| \leq 4  } |\td y- \td z|^{1-n}  \cdot |\td y_h|^{\al-1} d y  \,.
\]

Denote $s = \td y - \td z \in \R^n$.
Using $|s| \asymp |s_h| + |s_n|, n \geq 3$, and integrating over $s_n$,  we obtain
\[
  |I_j| 
    \les a^{\al} \int_{ |s_h| \leq 4 } \int_{|s_n| \leq 4}  ( |s_h| +|s_n| )^{1-n} |\td y_h|^{\al-1} d \td y_h d s_n
    \les a^{\al} \int_{ |s_h| \leq 4 }   |s_h|^{2-n} |\td y_h|^{\al-1} d \td y_h  \, . 
  \]

Using Young's inequality, we decompose the integrand as
\[
 M= \one_{|s_h| \leq 4 }   |s_h|^{2-n} |\td y_h|^{\al-1}
  \les   \one_{|s_h| \leq 4 } ( |s_h|^{\al+1 -n} + |\td y_h|^{\al+1 - n} 
(\one_{|\td y_h| \leq 1} +(\one_{|\td y_h| \geq 1}  )
   ) .
\]

Since  $\al+1 - n> 1-n$, and each part in $M$ is $L^1(\R^{n-1})$-integrable, we obtain
 \[
   |I_j| \les a^{\al}, \quad j = 1, 2.
 \]

For $I_3$, since $|y-x_1| \geq 2 a$, we obtain 
\[
|y-x_1| \asymp |y-x_2| ,
\quad  |K(y-x_1) - K(y-x_2)| \les |x_1 - x_2| \cdot |y-x_1|^{-n}. 
\]
Using $|f(y)| \les |y_h|^{\al-1}$ and then changing $y = a \td y, z = x_1 = a \td z$, we estimate $I_3$ as 
\[
  |I_3| \les a \int_{ |y-x_1| > 2 a } |y-x_1|^{-n} |y_h|^{\al-1} d y
  \les a^{\al} \int_{ |\td y - \td z| > 2  } |\td y - \td z|^{-n} |\td y_h|^{\al-1} d \td y  \, .
\]

Denote $s = \td y - \td z $. By first integrating $s_n$, we estimate
\[
  |I_3| \les a^{\al} \int_{ |s_h| \geq 0 }
  \int_{ |s_n| \geq \max(0, 2 - |s_h|) }  (s_h + s_n)^{-n} |\td y_h|^{\al-1} d \td y_h d s 
  \les  a^{\al} \int ( \max( s_h, 2 ) )^{1 - n} 
   |\td y_h|^{\al-1} d \td y_h .
\]

Since $\al-1 > -1, \al-n < -(n-1)$, it is easy to obtain that the integral 
in $\R^{n-1}$ is finite. We prove $|I_3| \les a^{\al}$.
Combining the above estimates of $I_j$, we complete the proof.
\end{proof}

\begin{lem}\label{lem:r_hol}
For any $x , y \neq 0, a\in \R$ and $b \in (0, 1]$, we have
\[
 |  |x|^{a} - |y|^a |  \les_{a, b} |x-y|^{b} ( |x|^{a-b} + |y|^{a-b}).
\]

\end{lem}

\begin{proof}
Using mean-value theorem and $| \, |x| - |y| \,| \leq |x-y| \leq |x| + |y|$, we prove
\[
\bal
 |  |x|^{a} - |y|^a | 
& \les_{a} \min( |x-y| \cdot ( |x|^{a-1} + |y|^{a-1} ), |x|^{a} + |y|^a)
\les_a  \min( |x-y|, |x | + |y| ) \cdot ( |x|^{a-1} + |y|^{a-1} ) \\
& \les_{a,b} |x-y|^b ( |x | + |y| )^{1-b} \cdot ( |x|^{a-1} + |y|^{a-1} ) 
\les_{a,b} |x-y|^b  ( |x|^{a-b} + |y|^{a-b} ) .
\eal
\]

\end{proof}

We define the annulus region. 
\[
B_{a, b} \teq \{ \yy : |\yy| \in (a, b) \},
\quad \bar B_{a, b}  \teq \{ \yy : |\yy| \in [ a, b ] \},
\quad B_R \teq \{ \yy: |\yy| < R \}. 
\]

\begin{lem}\label{lem:axi_reg}
Let $\chi_1$ be the cutoff function in \eqref{def:chi1}. Suppose that 
$\al, \g \in [0, 1), \g + \al > 1$, and the axisymmetric function $\Om(r, z)$ satisfies $\pa_z \Om(r, z)\in C^{ \g}(\bar B_R), \Om \in C^1(\bar B_R)$ 
for any $R \geq 1$. Denote 
\[
f(r, z)  \teq \pa_z \big( \Om(r, z) -  \Om(0, z) \cdot \chi_1( \tf{8 r^2}{\ang z^2} ) \big)  
\]
For any  $R \geq 1$, we have 
\bseq\label{eq:axi_reg}
\begin{align} 
        \|  f  r^{\al-1} \|_{C^{\al + \g-1}(\bar B_R)} & \les_R \| \pa_z \Om \|_{ C^{\g}(\bar B_{2R}) }
+ \|  \Om(0,\cdot)  \|_{ C^1(\bar B_{2R}) } ,
        \label{eq:axi_reg:a}  \\
  \| f r^{\al-1}  \|_{C^{\al + \g-1}(\bar B_{1.5R, 3.5R})}  & \les_R \| \pa_z \Om \|_{ C^{\g}(\bar B_{R, 4R}) }
        +  \| \Om(0,\cdot)  \|_{ C^{1}(\bar B_{R, 4 R}) } . 
        \label{eq:axi_reg:b}
\end{align}
\eseq

\end{lem}

\begin{proof}

We fix $R  \geq 0.1$.  Below, we focus on the proof for \eqref{eq:axi_reg:b}. 
The implicit constants in the following derivation may depend on $\g, \al, R$, we drop it to simplify notations.  Denote
\[
   m_R = \| \pa_z \Om \|_{ C^{\g}(\bar B_{R, 4R}) }
        +  \| \pa_z^{\leq 1} \Om(0, z)  \|_{ L^{\infty}(\bar B_{R, 4 R}) } ,
        \quad \xx = (r, z).
\]

A direct calculation yields 
\[
\bal
f(r,z) &= ( \pa_z \Om(r, z) - \pa_z \Om(0, z ) )   \cdot \chi_1( \tf{8 r^2}{\ang z^2} )  
+ \pa_z \Om(r, z) \cdot ( 1 - \chi_1( \tf{8 r^2}{\ang z^2} ) )  
- 8\pa_z \tf{r^2}{\ang z^2} \cdot \chi_1^{\pr} ( \tf{ 8 r^2}{\ang z^2} ) \cdot \Om(0, z)  \\
& \teq P_1 + P_2 + P_3.
\eal
\]

From the definition \eqref{def:chi1}, $\chi_1(y)$ is supported in $\{ y: |y|\leq 2 \}$.
For $|\xx | \geq R \geq 1$, in the support of $\chi_1( \f{8 r^2}{ \ang z^2}), \chi_1^{\pr} ( \tf{8 r^2}{\ang z^2} ) $, we obtain $8 r^2 \leq 2  \ang z^2$, which implies 
\[
  (1 + \tf{1}{4} ) |z|^2 \geq |z|^2 +  \tf{1}{4} \ang z^2 - \tf{1}{4} \geq |\xx|^2 - \tf{1}{4} \geq \tf{3}{4 } |\xx|^2 , \ \Rightarrow \ |z| > \tf{3}{4} |\xx|.
\] 

Thus, for $\xx \in \bar B_{1.4 R, 3.6 R} \cap \supp(\chi_1 ( \f{8 r^2}{\ang z^2} ) )$,  we obtain $|(0, z)| \in [1.05 R,  3.6 R]$. It follows
\beq\label{eq:hol_axi_pf1}
 \| f(r, z) \|_{ C^{\g}(\bar B_{1.4 R, 3.6 R}) }  \les m_R .
\eeq

For $P_2, P_3$, since $\chi_1^{\pr}(y), 1- \chi_1(y)$ are supported in $\{ y : |y| \geq 1 \} $, 
for any $(r, z) \in B_{4R}$, we have
\[
  | 1 - \chi_1( \tf{8 r^2}{\ang z^2} ) | + | \chi_1^{\pr} ( \tf{ 8 r^2}{\ang z^2} ) |  \les \tf{r^2}{\ang z^2} \les r^2,  \quad 
\Rightarrow \quad  \  |P_2| + |P_3|  \les r^2 m_R \les_R r^{\g} m_R.
\]

For any $\xx = (r , z) , \td \xx \in \bar B_{1.5R, 3.5R}$, since $R \gtr 1$
, using the above two estimates and $\pa_z \Om \in C^{\g}(\bar B_{R, 4 R})$, we obtain
\beq\label{eq:hol_bd}
\bal
  |f(r, z) | \les r^{  \g }  m_R ,
  \quad |f( \xx ) - f(\td \xx)| 
  \les |\xx- \td \xx |^{\g} \| f \|_{ C^{\g}(\bar B_{1.4R, 3.6 R} ) } 
  \les |\xx- \td \xx|^{\g} m_R.
\eal
\eeq

Suppose that $\td \xx = (\tr, \tz), \xx = (r, z)$ with $0 < \tr \leq r$ and $\td \xx \neq \xx$. The case $\tr = 0$ can be obtained using continuity since $f(\tr, z) \tr^{\al-1} |_{\tr =0}  = 0$. We estimate 
\[
|f(\xx) r^{\al-1}- f(\td \xx) \tr^{\al-1}|
\leq  | f( \td \xx ) (r^{\al-1} -\tr^{\al-1} )|
+ | f(\xx) -f( \td \xx )| r^{\al-1} = I_1 + I_2.
\]
For $I_1$, using  the pointwise bound for $f(\td \xx)$, Lemma \ref{lem:r_hol} with $a = \al-1$ and $b = \al + \g - 1 \in (0, 1]$, and $0< \tr \leq r$, we bound 
\[
   |I_1| \les \tr^{\g } m_R \cdot  | r^{\al-1} -\tr^{\al-1} | 
\les \tr^{\g } m_R \cdot ( r^{ - \g} + \tr^{  - \g })  |r-\tr|^{\al+ \g -1}
\les m_R |r-\tr|^{\al+ \g - 1}. 
\]

For $I_2$, interpolating two estimates in \eqref{eq:hol_bd} for $f$, and using $\al + \g-1 \in (0,\g)$ and $\td r \leq r$, we get 
\[
\bal
|I_2 | & \les \min( |\xx - \td \xx|^{\g}, r^{\g} + \td r^{\g} )   r^{\al-1} m_R
\les |\xx-\td \xx|^{\g} \min(1, r |\xx - \td \xx|^{-1} )^{\g}  r^{\al-1} m_R \\
& \les 
 |\xx-\td \xx|^{\g} ( r |\xx - \td \xx|^{-1} )^{1-\al} \cdot  r^{\al-1} m_R 
\les |\xx-\td \xx|^{\al + \g-1} m_R.
\eal
\] 
Combining the above estimates for $I_i$, we prove \eqref{eq:axi_reg:b}. 

The proof for \eqref{eq:axi_reg:a} is similar and simpler, and thus it is omitted.
\end{proof}

\subsection{$\e$-interpolation estimates}

Below, we develop several $\e$-interpolation estimates.

\begin{lem}\label{lem:interp1}
Let $\e = \f13 - \al$. Recall $ \kp_1, \epa$ from \eqref{def:kp} and \eqref{ran:ep_all}. 
Suppose that the parameter $\hal, \hau, \b$ satisfies the range in \eqref{ran:ep_all} and \eqref{def:kp}
\beq\label{ran:ep_rec2}
 \e - \b, \  \hal - \al   , \  \hau - \al  \in [ \tf{9}{8} \e - C \e^{ 2 -\kp }, \, \tf{9}{8} \e + C \e^{ 2 -\kp } ]
     \subset [ \tf{8.9}{8} \e,  \ \tf{9.1}{8} \e ] ,   \quad 
     \epa \in [0, \kp_1 \e],  \ \kp_1 = \tf{1-\kp}{1000} \in (0, 10^{-3} ) . 
\eeq
There exists an absolute constant $ \beps_{\mw{pow}} \in (0, \f13)$ such that for any  $\e \in (0, \beps_{\mw{pow}}], a >0$, we have
\begin{align}
  \ang \xx^{\al -\hal + \epa} & \leq  a^{ -\kp_1 } \ang \xx^{\al -\hau } 
  +  a , \label{eq:interp:pow} \\
  \ang \xx^{\al-\hal} +  \ang \xx^{\al -\hal + \e^2} 
 + \ang \xx^{-\e + \b + \e^2}
 &  \leq ( 1 + C \e^{1/6} ) \ang \xx^{\al-\hau} + \e^4 , \label{eq:interp:pow2} \\
    \cJaa(\xx) \cdot \e +  \ang \xx^{ -  10 \e  } & \gtr 1 \, , \label{eq:interp:J}
\end{align}
with implicit constant independent of $\e$.
\end{lem}

\begin{proof}

 Since $\kp_1 \in (0, 1)$,  using Young's inequality we obtain
\[
  a^{ - \kp_1} \ang \xx^{\al - \hau} + a 
  \geq (   a^{ - \kp_1} \ang \xx^{\al - \hau} )^{ \f{ 1 }{ 1 + \kp_1 } } \cdot a^{ \f{\kp_1}{1 + \kp_1} }
  \geq \ang \xx^{ \f{\al -\hau}{1 + \kp_1} }.
\]

Using \eqref{ran:ep_rec2}, $\epa \leq \kp_1 \e$, and $\kp_1 < \f{1}{100}$ by \eqref{ran:ep_rec2}, for $\e$ small enough, we prove 
\[
    \f{\al-\hau}{\kp_1 +1}  - (\al - \hal + \epa)
    \geq O( \e^{ 2 -\kp } ) - \f{9 \e}{8 (\kp_1 +1)} + \f{9 \e}{8} -  \kp_1 \e 
    = O(\e^{2 -\kp}) + \kp_1 \e (  \f{9}{ 8(\kp_1 + 1) } - 1 ) > 0.
\]
It follows \eqref{eq:interp:pow}.

To prove \eqref{eq:interp:pow2}, from \eqref{ran:ep_rec2},  the powers among $\al-\hal, \al -\hau, -\e + \b + \e^2$ 
differs by $O(\e^{2 -\kp})$. 

If $\ang \xx^{\e^{1 + \mhk}} \leq 2$, 
we yield  $\ang \xx^{  \e^{2 - \kp} } \leq 2^{\e^{ \mhk } } \leq 1 + C \e^{\mhk}$. Thus, estimate \eqref{eq:interp:pow2} holds trivially. 
If $\ang \xx^{\e^{1 + \mhk}} \geq 2$, using $\al-\hal,  \al -\hal + \e^2, -\e + \b + \e^2 \leq -\e$ 
for $\e$ small enough by \eqref{ran:ep_rec2}, we prove
\[
  \ang \xx^{\al-\hal + \e^2}
  +  \ang \xx^{\al-\hal }
  +   \ang \xx^{-\e + \b + \e^2 }
   \les \ang \xx^{ - \e} 
  \leq ( \ang \xx^{\e^{1 + \mhk}} )^{- \e^{-\mhk}} \leq 2^{-\e^{-\mhk}} \les \e^4.
\]

For \eqref{eq:interp:J}, recall $\cJaa \asymp \min( \lgp x , \f{1}{\e} )$ from \eqref{eq:Ja_hat}. 
If $\lgp |\xx| \geq \f{1}{\e}$, we obtain \eqref{eq:interp:J}. If $\lgp |\xx| \leq \e^{-1}$, 
we obtain $ \log \ang \xx \les \e^{-1}, \ang \xx^{-\e} \gtr 1$, and $\ang \xx^{- 10 \e} \gtr 1$.
We prove \eqref{eq:interp:J}.
\end{proof}

\begin{lem}\label{lem:interp2}
Let $\xx = (r, z)$. Suppose $k, \al$ satisfies $\al - k \in (-1, 0), \al \in (0, 1/3), \e =\alb - \al $,
and $\e \leq \bar c_1$, where $\bar c_1$ is the constant in Lemma \ref{lem:interp1}. For any $\ell > 0$, we obtain 
\beq\label{eq:est_interp2}
    \min( \one_{r > |z| }\ang z^{\al- k},  |\xx|^{\al + 1} ) 
    \les 
 \f{r}{\ang \xx} \big( \ell^{-1} \ang \xx^{-10\e}  +     \e \cJaa(\xx)  \cdot \ang z^{-2} )
  + \e \cJaa(\xx) \cdot \f{|z|}{\ang z} \ang z^{\al-k} 
  + \ell ,
\eeq
with implicit constants independent of $\ell, \e, a, \al, k$. 
The $\cJaa(\xx)$ evaluates at $\xx$ rather than $z$.

\end{lem}

\begin{proof}

Denote by $\mw{LHS}, \mw{RHS}$ the left hand side, and the right hand side of the inequality, 
respectively.

If $r \leq |z|$, the estimate is trivial since $\mw{LHS}_{\eqref{eq:est_interp2} } = 0$. 

Next, we consider $|z| \leq r$. If $|\xx| \leq 1$, since $r \asymp |\xx|$, using Young's inequality for $\mw{RHS}$, we obtain
\[
 \mw{RHS}_{\eqref{eq:est_interp2} } 
\geq \ell + r \ell^{-1} 
\gtr r^{1/2} 
  \gtr |\xx|^{1/2} \gtr |\xx|^{1 + \al} \gtr \mw{LHS}_{\eqref{eq:est_interp2} }.
\]

If  $|\xx| \geq 1$, we obtain $ r \gtr \ang \xx$.  Since 
\[
\ang z^{-2} + |z|\ang z^{\al-k-1}  
\gtr  \ang z^{\al-k} \one_{|z|\leq 1} + |z| \ang z^{\al-k-1}   \gtr \ang z^{\al-k}  , 
\]
and $1 \geq \ang z^{\al -k}$,  we obtain
\[
  \mw{RHS}_{\eqref{eq:est_interp2} }
  \gtr \e \cJaa(\xx) \cdot \ang z^{\al - k} + \ang \xx^{-5\e} 
  \geq (\e \cJaa(\xx) + \ang \xx^{-5 \e} ) \ang z^{\al - k} . 
\]
Using \eqref{eq:interp:J}, we prove $  \mw{RHS}_{\eqref{eq:est_interp2} } \gtr \ang z^{\al-k}$ and  the desired estimates.
\end{proof}

Recall $\cJaa \asymp \min( \lgp x, \e^{-1} )$ from \eqref{eq:Ja_hat}. 
The following lemma compares $\cJaa$ and $\lgp x$:

\begin{lem}\label{lem:lgx_pow}

For any $k>0, \ell > 0$ and $\e < 1$, we have %$\cJaa \les \lgp x$ and
 \[
   \ang x^{-k \e} |\lgp x|^\ell \les_{k,\ell} |\cJaa(x)|^\ell  \les  \e^{-\ell}.
 \]
\end{lem}

Note that the constant is independent of $\e$. The proof follows easily from $\ang x^{-k \e} = e^{ -k \e \log \ang x }, \lgp x \asymp \log (1 + \ang x)$ and discussing $ \e \log \ang x < 1 $ and $ \e \log \ang x > 1$.

We have the following estimates for the $\cJ$-integral.

\begin{lem}\label{lem:log_ineq_J}
For $ k > -1$ and $c > 0$, we obtain
\begin{align}\label{eq:log_ineq_J:a}
 \int_0^x  \ang y^{ - c \e - 1 }  \cJaa^{k}(y) d  y
 \les_{c, k} \min( |x|, \cJaa^{k+1}(x) ) .
\end{align}

\begin{comment}

\bseq\label{eq:log_ineq_log}

For any $ c>0$ and $k>-1$, we obtain
\beq\label{eq:log_ineq_log:a}
  \int_0^{x} \ang z^{- c \e - 1} |\lgp z|^{k} d z 
  \les_{c, k} \min( |x| , \cJaa^{k+1}(x) ).
\eeq
For $k = -1$ and any $a > 0$, we obtain
\beq\label{eq:log_ineq_log:b}
 \int_0^{x} \ang z^{- c \e - 1} |\lgp z|^{-1} d z 
  \les_{c, a} \min( |x| , \cJaa^{a}(x) ).
\eeq

For $ k < -1 $ and any $c \geq 0$, we obtain 
\beq\label{eq:log_ineq_log:c}
   \int_0^{x} \ang z^{- c \e - 1} |\lgp z|^{k} d z \les_k \min(|x|, 1).
\eeq

For $k > -1$, we obtain 
\beq\label{eq:log_ineq_log:d}
   \int_0^{x} \ang z^{ - 1} |\lgp z|^{k} d z \les_k \min( x,   |\lgp x|^{k+1}).
\eeq

\eseq

\end{comment}

\end{lem}

\begin{proof}
The first bound in \eqref{eq:log_ineq_J:a} by $C_k |x|$ is trivial.  We define
$H(x) = 1 + \int_0^x \ang y^{-c \e - 1} d y$ and obtain
\[
H(0) = 1, \quad  H(x) \asymp_c \min( \e^{-1}, \lgp x ),
  \quad \pa_x H(x) = \ang x^{- c \e - 1} > 0 .
\]
Since $\cJaa \asymp \min( \e^{-1}, \lgp x) \asymp_c H$ \eqref{eq:Ja_hat}, and $k+1>0$,   we bound 
\[
\bal
  I & = \int_0^x \ang y^{-c \e - 1} \cJaa^k(y) d y 
\les_{c, k} \int_0^x \pa_y H(y) \cdot H(y)^k d y   \les_{c, k} 1 +  H(x)^{k+1}
\les \cJaa^{k+1}(x).
\eal
\]
We prove  \eqref{eq:log_ineq_J:a}.  
\end{proof}

\section{Local well-posedness}\label{app:LWP}

In this section, we prove local well-posedness for
\beq\label{eq:LWP}
\bal
   \pa_{s}  \Om + \QQ(\Om) \cdot \na \Om & = \cB(\Om) \cdot \Om
  + \sss \cdot  \cU_h ( f  -   \Om  ), \quad \sss \in \{ 0, 1\} , \\
\eal
\eeq
on the domain $(r, z) \in \R_+ \times \R$, where $\Psi = \BS(\Om)$ is defined in \eqref{eq:5D_BS}, 
$\cU_h$ is defined in \eqref{def:cU_h}, and
\[
  \QQ(\Om) = ( c_l r - r \pa_z \Psi, \,  c_l z + 2 \Psi + r \pa_r \Psi  ), 
  \quad  c_l(\Om) = 2 - 2 \Psi_z(0),  \quad 
  \cB(\Om) =  2 - (1-\al) (\Psi_z - \Psi_z(0)) .
\]

Below, we work on the domain $(r,z)\in \R_+\times\R$ rather than full $\R^3$.
Recall the weight $\phia$ from \eqref{def:phi0}. We define the weighted $C^1$ and $C^{1,\al}$ norm 
in $(r, z) \in \R_+ \times \R$:
\bseq\label{def:cC1}
\beq
\phia =  |\xx|^{-1} + \ang \xx^{\aln}, \quad 
   \| g \|_{\CCs^{1}} \teq 
    \nlinf{ \phia g }
   + \nlinf{  \phia  |\xx|  \cdot \na g },
   \quad 
   \| g \|_{\CCs^{1,\al}} \teq    \| g \|_{\CCs^{1}}  + \| \na g \|_{\dot C^{\al}} ,
\eeq
where $\aln$ is defined in \eqref{def:cZ} with $\aln > \al$.  Clearly, by definition, for any $ g(r, z) \in C_c^{1, \al}$ that is odd in $z$, we obtain 
$g \in \CCs^{1,\al}$. From the definition of $\cZ$-norm in \eqref{def:cZ}, we obtain
\beq
  \| g \|_{\CCs^1} \les \| g \|_{\cZ} ,
  \quad   \| g \|_{\CCs^{1,\al}} \les \| g \|_{\cZ} + \| g \|_{C^{1,\al}},
  \quad 
  \| g \|_{\CCs^{1,\al}} \les (1 + |\supp(g) | )^2 \cdot \| g \|_{C^{1,\al}} ,
\eeq
where $|\supp(g) | \teq \inf \{ R : \supp(g) \in B_R \}$.
\eseq

 We have the following well-posedness results.

\begin{prop}\label{prop:LWP}
Let $\sss \in \{ 0, 1 \}$ and $\cU_h$ be defined in \eqref{def:cU_h}.  Suppose that $\Om(0) \in C_c^{1,\al},   \ang \xx^{-2 \e^3} \phia f \in L^{\infty}_{s, \xx }$. There exists 
$T = T(  \| \ang \xx^{-2 \e^3} \phia  f \|_{L_{s, x}^{\infty} } , \| \Om(0) \|_{\CCs^{1,\al}} ) > 0$ such that \eqref{eq:LWP} admits a unique  $\CCs^{1,\al}$ solution for $s \in [0, T]$.
The solution can be extended beyond $T_1$ for $T_1 \geq T$ if 
\beq\label{eq:LWP_blowup}
 M(T_1) \teq \tts{ \int_0^{T_1} } \nlinf{   \ang \xx^{-2 \e^3}  \phia  \Om(s)  } d s < \infty.
\eeq
Moreover, we have the following estimates 
\begin{align}\label{eq:LWP_EE}
M(T) & \leq 1 , \\
 \sup\nolim_{s \leq T_1}  \| \Om(s) \|_{\CCs^{1,\al}} 
 +      \nlinf{ ( |\xx|^{-\al-1} + |\xx|^{\aln} ) \pa_s \Om(s) }
 & \leq C( M(T_1), \| \Om_0 \|_{\CCs^{1,\al}} ,  \|  \ang \xx^{-2 \e^3} \phia f \|_{L_{s, x}^{\infty} } ). \notag 
\end{align}

\end{prop}

Taking $\sss = 0$, we obtain the local well-posedness result for the dynamic equation \eqref{eq:dyn_full}.
We do not seek a strong blowup criterion.
We solve \eqref{eq:LWP} using an iterative scheme similar to that in 
\cite{majda2002vorticity}.

\vs{0.05in}
\paragraph{\bf Estimates on $\Psi$}
Below, the implicit constants may depend on $\e = \f13 - \al $. 
We omit the dependence on $\e$ since it does not affect the proof. We have the following estimates on $\Psi = C \KRF \ast ( \Om r^{\al-1} )$ from Lemma \ref{lem:hol_singu} and Proposition \ref{prop:iso_est}
with $\bbu = 1, \g = \aln - 2 \e^3 \in [\hal - \e, \hal]$ by \eqref{def:wg_s}:
\bseq\label{eq:LWP_psi_est1}
\begin{gather}
 \tf{1}{|\xx|} | \QQ(\Om)| + |\cB(\Om)| + | \na \Psi | + |\tf{1}{z} \Psi | \les 1 + \nlinf{
  \ang \xx^{-2 \e^3} \phia \Om } ,  
\label{eq:LWP_psi_est1:a}
   \\
 \quad | \cB(\Om) - 2 | + \tf{1}{|\xx|}   |\QQ^z(\Om) - 2 z| 
 \les |\xx|^{\al} \nlinf{  \ang \xx^{-2 \e^3} \phia  \Om } ,  \label{eq:LWP_psi_est1:b} \\ 
\| \na \cB(\Om) \|_{ C^{\al} }
\les \nlinf{ \na \Psi_z } 
+ \| \na \Psi_z \|_{\dot C^{\al}} 
\les    \nlinf{ \pa_z \Om }  + \nlinf{ \phia ( |\Om| + |\xx| \cdot |\na \Om| ) }
\les   \| \Om \|_{\CCs^{1}}, 
\label{eq:LWP_psi_est1:c}
 \\
|\na \QQ(\Om)| \les 1 + |\na \Psi(\Om)|  + |\na_{r, z}^2 ( r \Psi(\Om) )| 
\les  1+ \nlinf{ \phia ( |\Om| + |\xx| \cdot |\na \Om| ) }  \les 1+  \| \Om \|_{\CCs^1}.
\label{eq:LWP_psi_est1:d}
\end{gather}
\eseq

Using \eqref{eq:Uh_bound} and \eqref{eq:supp_cU_h}, for any $F$, we bound 
\beq\label{eq:LWP_Uh}
\nlinf{ \phia \cU_h F }   +  \| \cU_h F \|_{C^{1,\al}} \les  
 \nlinf{  \ang \xx^{-2 \e^3}  \phia  F} ,
 \quad \supp(\cU_h F) \subset \{ |\xx| \in [h_1, h_1^{-1}] \} .
\eeq

Using product rule, we obtain
\beq\label{eq:LWP_pf1_Cal}
  \| \Om r^{\al} \|_{C^{\al}}
  \les \| \Om \ang \xx^{\al} \|_{C^{\al}}  \| r^{\al} \ang \xx^{-\al} \|_{C^{\al}}
  \les \| \Om \ang \xx^{\al} \|_{C^1} \les \| \Om \|_{\CCs^{1}}.
\eeq

Since $\D_{\R^3}( r \Psi ) =  \Om r^{\al}$, 
using  H\"older estimates for the Calderon Zygmund operator 
$\na^2(r \Psi) = \na^2 (-\D)^{-1} (\Om r^{\al}) $ \cite[Lemma 4.6]{majda2002vorticity}, 
Lemma \ref{lem:hol_singu}, and \eqref{eq:LWP_pf1_Cal}, we obtain
\beq\label{eq:LWP_psi_Cal}
\| \na \QQ \|_{\dot C^{\al}}
\les \| \na \Psi \|_{\dot C^{\al}}
+ 
 \| \na^2 (r \Psi) \|_{\dot C^{\al} } \les 
\nlinf{ \Om } + \| \Om r^{\al} \|_{\dot C^{\al}}
\les \| \Om \|_{\CCs^{1}}.
\eeq

\vs{0.05in}

\paragraph{\bf Iterative scheme}

Let  $\Om_0 \equiv \Om(0)$. We construct approximate solution to \eqref{eq:LWP} as 
\beq\label{eq:LWP_iter1}
  \pa_s \Om_{n+1} + \QQ(\Om_n) \cdot \na \Om_{n+1} = \cB(\Om_n) \Om_{n+1} 
+ \sss \cdot \cU_h( f-  \Om_n ) ,
\quad  \Om_{n+1}(0) = \Om_0.
\eeq

We introduce the characteristic $X_n$ associated with $\QQ(\Om_n)$ and coefficient $S_n$
\[
  \f{d}{d s} X_{n}(s, \xx ) = \QQ(\Om_n(s))( X_n(s, \xx) ),
  \quad 
\cS_n(s, \xx) \teq \int_0^s \cB(\Om_n(\tau)) \cc X_n( \xx  , \tau) d \tau.
\]

Using method of characteristics and following Section \ref{sec:reg:eta_C1a}, we solve $\Om_{n+1}$ 
\beq\label{eq:LWP_solu_form}
  \Om_{n+1}( s, X_n(s, \xx ) ) = e^{\cS_n(s, \xx)} \Om_0(\xx)
  +\sss \int_0^s e^{\cS_n(s, \xx) - \cS_n(\tau, \xx)}    \cU_n (f-\Om_n ) \cc X_n(\tau, \xx) d \tau .
\eeq

\paragraph{\bf Energy estimates}
Denote 
\beq\label{def:Lam_n}
  \Lam_n(s) = \int_0^s ( 1+  \| \Om(s) \|_{\CCs^1} )  d s. 
\eeq

Performing weighted $C^1$ estimate on $\Om_{n+1}$
by following the argument in Section \ref{sec:decay_cT}, and using 
\beq\label{eq:LWP_damp1}
 \phia^{-1} | \QQ_n \cdot \na \phia |  \les 
 |\QQ_n | \cdot |\xx|^{-1}
\les  \nlinf{  \ang \xx^{-2 \e^3} \phia \Om_n } \, , 
\eeq
and \eqref{eq:LWP_psi_est1} on $\cB(\Om_n), |\na \QQ_n|$, we obtain
\beq\label{eq:LWP_EE0}
\bal
\tf{d}{d s}
   & \phia ( |\Om_{n+1}|^2 + |\xx|^2 \cdot |\na \Om_{n+1}|^2 )^{1/2} \cc X_n(s, \xx) \\
  & \leq  C ( \nlinf{ \ang \xx^{-2 \e^3} \phia \Om_n } +  \nlinf{ \na \QQ_n} )
  \cdot 
 \phia ( |\Om_{n+1}|^2 + |\xx|^2 \cdot |\na \Om_{n+1}|^2 )^{1/2} \cc X_n(s, \xx) \\
 & \quad  +  \nlinf{ \ang \xx^{-2 \e^3} \phia  f } + \nlinf{ \phia \Om_n }.
 \eal
\eeq

Applying Gronwall's inequality and then taking supremum over $\xx$, we obtain
\beq\label{eq:LWP_C1}
  \| \Om_{n+1}(s) \|_{\CCs^1}
  \leq e^{ \bar C   \Lam_n(s)  } \| \Om_0 \|_{\CCs^1}
  + \int_0^s  e^{ \bar C (  \Lam_n(s)  - \Lam_n(\tau) ) } 
   (  \nlinf{ \ang \xx^{-2 \e^3} \phia  f(\tau) } + \nlinf{ \phia \Om_n(\tau) } ) d \tau.
\eeq
where $\bar C$ is some absolute constant which is fixed for later derivations. We consider the following uniform bound 
\bseq\label{eq:LWP_unif}
\beq\label{eq:LWP_unif:a}
   \sup\nolim_{s \leq T} \| \Om_{n}(s) \|_{\CCs^1} \leq 4 \| \Om_0 \|_{\CCs^1}, \quad n \geq 0,
\eeq
which holds for $n=0$. We choose
\beq\label{eq:LWP_T}
  T =  \f{ \min(1,\| \Om_0 \|_{\CCs^1} )  } {8 ( \bar C +1) ( | \Om_0 \|_{\CCs^1} + 
\| \ang \xx^{-2 \e^3}  \phia  f \|_{L_{s, x}^{\infty} } + 1)  }.
\eeq
If \eqref{eq:LWP_unif:a} holds for $n$, 
using \eqref{def:Lam_n}, and $\nlinf{\phia \Om_n} \leq   \| \Om_{n} \|_{\CCs^1}
\leq 4  \| \Om_0 \|_{\CCs^1} $ by \eqref{def:cC1}, we get
\beq
\bar C \Lam_n(s)  \leq \tf18, \quad  e^{\bar C \Lam_n(s)} \leq e^{1/8} \leq  2,
\quad \forall \ s \in [0, T].
\eeq
\eseq
Thus, using \eqref{eq:LWP_C1} and induction on $n \geq 0$, we prove \eqref{eq:LWP_unif:a}.

\vs{0.1in}
\paragraph{\bf Higher regularity}

To estimate $ \|\Om_{n+1}(s) \|_{\dot C^{1,\al}}$, we use \eqref{eq:LWP_solu_form}. Below, we consider $s \in [0, T]$.
Following \eqref{eq:3D_ODE:X_C1} and \cite[Section 1.3]{majda2002vorticity}, we have the following 
identities  for $\na X_n$ and its determinant:
\[
\bal
\tf{d}{d s} (\na X_n)  \cc X_n(s,\xx) &=  (\na \QQ_n(s) \cdot \na X_n) \cc X_n(\xx, s) ,
& \quad  \na X_n( 0, \xx) &= \Id, \\
  \tf{d}{d s} \det (\na X_n) & = (\na_{r, z} \cdot \QQ) \det (\na X_n), 
  & \quad  \det (\na X_n) |_{s = 0} & = 1 . 
\eal
\]
Using the uniform estimate \eqref{eq:LWP_unif} and the bound on $\na \QQ_n$ by \eqref{eq:LWP_psi_est1:d}
and \eqref{eq:LWP_psi_Cal}, we bound 
\[
1  \les \det (\na X_n)(s, \xx ) \les 1, \quad 1 \les   |\na X_n(s, \xx )| \les 1, \quad \forall \ s \in [0, T].
\]

Thus, the flow map $X_n(s, \cdot )$ is invertible.  Performing $\dot C^{\al}$ estimate on $\na X_n$, using the above bound, 
\eqref{eq:LWP_psi_est1:d} for $|\na \QQ(\Om_n)|$,  \eqref{eq:LWP_psi_Cal},
\eqref{eq:LWP_unif}  for $ \| \na \QQ(\Om_n) \|_{\dot C^{\al} } $,
and Gronwall's inequality,  we obtain
\[
  \| \na X_n \|_{\dot C^{\al}} \les C  + \tts{\int_0^s}  \| \na \QQ(\Om_n(\tau)) \|_{\dot C^{\al} }  d \tau
  \les C + T \| \Om_0 \|_{\CCs^1} \les 1.
\]
where the implicit constants can depend on $\| \Om_0 \|_{\CCs^1}$.

Applying the above estimates and \cite[Lemma 4.2]{majda2002vorticity} on the flow map $X$,
and \cite[Lemma 4.3]{majda2002vorticity} on the composition $g \cc X_n, g \cc X_n^{-1}$, we obtain
\beq\label{eq:hol_compose}
\| g \cc X_n(s, \cdot) \|_{C^{1,\al}}  \les \| g \|_{C^{1,\al}} , 
\quad  \| g \cc X_n^{-1}(s, \cdot) \|_{C^{1,\al}}  \les \| g \|_{C^{1,\al}} ,
\quad s \in [0, T].
\eeq

As a result, applying the $C^{1,\al}$ bound on $\cB(\Om), f, \cU_h \Om$ from 
\eqref{eq:LWP_psi_est1:c}, \eqref{eq:LWP_Uh}, and \eqref{eq:LWP_unif}, we bound 
\[
\bal
  & \| \cB( \Om_n(s) ) \cc X_n(s, \cdot) \|_{C^{1,\al}} 
  +   \| \cU_h( f(s) - \Om_n(s)  )  \cc X_n(s, \cdot) \|_{C^{1,\al}}  
  \les \| \Om_0 \|_{\CCs^1} 
+  \| \ang \xx^{-2 \e^3}  \phia f \|_{L_{s, x}^{\infty} } \, .
\eal
\]

Applying the above estimate and \eqref{eq:LWP_unif} to \eqref{eq:LWP_solu_form}, the product rule, 
and \eqref{eq:hol_compose}, we obtain %$C^{1,\al}$ bound
\beq\label{eq:LWP_unifb}
   \| \Om_{n+1}(s) \cc X_n(s, \cdot) \|_{C^{1,\al}} \les  C(T, \| \Om_0 \|_{C^{1,\al}} ,  
   \| \Om_0 \|_{\CCs^1} )  
     \  \Rightarrow  \   \| \Om_{n+1}(s) \|_{C^{1,\al}}  \les   C( T,   \| \Om_0 \|_{\CCs^{1,\al}} ) ,
\eeq
uniformly for $s \in [0, T]$ with $T$ chosen in \eqref{eq:LWP_T}. 

\vs{0.1in}

\paragraph{\bf Convergence}
Using \eqref{eq:LWP_iter1}, we derive the equation for $\Om_{n+1} - \Om_n$
\[
\bal
 & \pa_s ( \Om_{n+1} - \Om_n )
  + \QQ(\Om_n) \cdot \na ( \Om_{n+1} - \Om_n )
  + \QQ(\Om_n - \Om_{n-1}) \cdot \na \Om_n \\
 & \qquad \quad  = \cB(\Om_n) ( \Om_{n+1} - \Om_n ) + \cB(\Om_n - \Om_{n-1}) \Om_n
  - \sss  \cdot \cU_h( \Om_n -\Om_{n-1} ). 
\eal
\]

Estimating $ \phia (\Om_{n+1} - \Om_n)$ \eqref{eq:LWP_iter1} along the characteristics $X_n(s, \cdot)$
and using the bound in \eqref{eq:LWP_psi_est1:a}, and uniform bound \eqref{eq:LWP_iter1}, we obtain
\[
   \tf{d}{d s}  |\phia (\Om_{n+1} - \Om_n) \cc X_n(s, \xx)| 
   \les ( 1+ \| \Om_0 \|_{\CCs^1} ) \cdot ( |\phia (\Om_{n+1} - \Om_n) \cc X_n(s, \xx)| 
   +  |\phia (\Om_{n} - \Om_{n-1}) \cc X_n(s, \xx)| ).
\]

Since $\Om_{n+1}(0, \xx) = \Om_n(s, \xx) = \Om_0$, by taking $T$ small, if necessary, 
using Gronwall's inequality and taking supremum over $\xx$, we prove 
\[
 \sup\nolim_{s \in [0, T] } \nlinf{ \phia (\Om_{n+1} - \Om_n)  }
  \leq \tf{1}{2} \sup\nolim_{s \in [0, T]}  \nlinf{ \phia (\Om_{n} - \Om_{n-1})  }.
\]

Thus, $\Om_{n}$ converges to some $\Om$ uniformly in $L^{\infty}( [0, T], L^{\infty}( \phia ))$. By the uniform estimate \eqref{eq:LWP_unif} and \eqref{eq:LWP_unifb}, we obtain 
strong convergence of $\Om_n - \Om$ in% $L^{\infty}( [0, T], \CCs^1)$-norm and 
$L^{\infty}( [0, T], C^{1, \g})$-norm with $\g < \al$. Moreover,
$\Om \in L^{\infty}( [0, T], \CCs^{1, \al})$ and it satisfies estimates \eqref{eq:LWP_unif} and \eqref{eq:LWP_unifb}.

Taking $n \to \infty$ in \eqref{eq:LWP_iter1} and using the strong convergence, we prove that $\Om$ is a $L^{\infty}( [0, T], \CCs^{1, \al})$ solution to \eqref{eq:LWP}. The uniqueness of $\CCs^{1,\al}$ solution follows from estimates similar to those established above for $\Om_{n+1} - \Om_n$.
The bound $M(T) \leq 1$ in \eqref{eq:LWP_EE} follows from \eqref{eq:LWP_unif} with $\Om_n = \Om$.

\vs{0.05in}

\paragraph{\bf Continuation criterion}

 Since $\D_{\R^3}(r \Psi) = \Om r^{\al}$, 
 by estimating the singular integral in the region 
$|\xx - \yy| \geq 1, |\xx - \yy| \in [\d, 1]$, and $|\xx-\yy| \leq \d$, 
using the $L^p$, $L^{\infty}$, and $C^{\al}$-norm, respectively, and then optimizing $\d$, 
 we obtain the log-extrapolation
\[
 \nlinf{ \na^2 (r \Psi ) } \les_p  \log ( 2+ \| \Om r^{\al} \|_{C^{\al} }   )    \nlinf{ r^{\al} \Om }
  +\| r^{\al} \Om \|_{L^p} ,\quad \forall p \in [1, \infty).
\]
Since $ \ang \xx^{-2 \e^3} \phia  \gtr \ang \xx^{\aln -  2 \e^3}$ and $\aln - 2 \e^3 > \al$, we obtain $\| r^{\al} \Om \|_{L^p} \les \nlinf{ 
\ang \xx^{-2 \e^3} \phia \Om }$ for 
$p \in [p_0 , \infty]$ with  $p_0$ large enough. Using  \eqref{eq:LWP_psi_est1} and  \eqref{eq:LWP_pf1_Cal}, we obtain
\beq\label{eq:LWP_log}
 \nlinf{ \na \QQ }
 \les  1 + |\na \Psi| + | \na^2 (r \Psi ) | 
    \les 1+    \log ( 2+  \| \Om \|_{\CCs^{1}}  ) \cdot \nlinf{  \ang \xx^{-2 \e^3}  \phia \Om } .
\eeq

Using the energy estimates \eqref{eq:LWP_EE0} with $\Om_n = \Om_{n+1} = \Om$, we prove that
$ M(s) = \int_0^s \nlinf{
\ang \xx^{-2 \e^3} \phia \Om }$ controls the regularity $\| \Om \|_{\CCs^1}$. 
From the estimate in \eqref{eq:LWP_unifb}, $M(s)$ further controls the $\CCs^{1,\al}$-regularity.
Thus, we prove the continuation criterion \eqref{eq:LWP_blowup} and the estimate 
on $\| \Om \|_{\CCs^{1,\al}}$ in \eqref{eq:LWP_EE}.

\vs{0.05in}

\paragraph{\bf Estimate of $\pa_s \Om$ in \eqref{eq:LWP_EE}} 

Since $\Om \in C^{1,\al}$ and is odd in $z$, we get $\pa_r \Om(0) = 0$ and 
\[
 |\pa_r \Om| + |z \pa_z \Om - \Om| \les \min( |\xx|^{\al}, 1) \| \Om \|_{\CCs^{1,\al}}.
\]

Using estimate \eqref{eq:LWP_psi_est1:a}, \eqref{eq:LWP_psi_est1:b} on $\QQ, \cB(\Om)$, and the above estimates, we obtain
\[
\bal
  I \teq \phia | - \QQ \cdot \na \Om + \cB(\Om) \Om   |
& \les \phia \B( | - (\QQ^r, \QQ^z - 2 z) \cdot \na \Om 
+ ( \cB(\Om) - 2 ) \Om |  + | 2 \Om - 2 z \pa_z \Om| \B) \\ 
& \les \min( |\xx|^{\al} , 1) ( \| \Om \|_{\CCs^{1,\al}}^2  + \| \Om \|_{\CCs^{1,\al}}).
\eal
\]

Using the equation \eqref{eq:LWP}, the above estimate, and \eqref{eq:LWP_Uh}, for $s \in [0, T]$, we obtain
\[
  |\phia \pa_s \Om | \leq 
    I + 
   \phia | \cU_h( f-  \Om)|
\les \min(|\xx|^{\al}, 1) ( \| \Om \|_{\CCs^{1,\al} }^2 + 
 \| \Om \|_{\CCs^{1,\al}} +  \|  \ang \xx^{-2 \e^3}  \phia  f \|_{L_{s, x}^{\infty} } ).
\]
Using the estimate for $ \| \Om \|_{\CCs^{1,\al}}$ in \eqref{eq:LWP_EE} proved in the 
Paragraph \textit{Continuation criterion}, using $\phia$ in \eqref{def:cC1}, and dividing the above estimates by $\min( |\xx|^{\al},1)$, we prove the estimate for $\pa_s \Om$ in \eqref{eq:LWP_EE}. We complete the proof of Proposition \ref{prop:LWP}.

\bibliographystyle{plain}
\bibliography{selfsimilar}

\end{document}